\documentclass[a4paper,11pt]{book}

\title{Global nonlinear stability of Minkowski space \\ for spacelike-characteristic initial data}
\author{Olivier Graf}
\date{\today}

\usepackage{fontenc}
\usepackage[english]{babel}
\usepackage{appendix}
\usepackage{setspace} 
\usepackage{graphicx}
\usepackage[dvipsnames]{xcolor}
\usepackage{caption}
\usepackage{subcaption}
\usepackage{booktabs}
\usepackage{microtype}
\usepackage{url}
\usepackage[final]{pdfpages}
\usepackage[inner=20mm,outer=20mm,top=20mm,bottom=30mm]{geometry} 
\usepackage{fancyhdr}

\fancypagestyle{plain}{
	\fancyhf{}

	\fancyfoot[C]{\thepage}}
\fancypagestyle{addpagenumbersforpdfimports}{
	\fancyhead{}
	
	\fancyfoot{}
	\fancyfoot[C]{\thepage}
}

\usepackage{listings}
\usepackage{hyperref}
\hypersetup{pdfborder={0 0 0},
	colorlinks=true,
	linkcolor=black,
	citecolor=black,
	urlcolor=black}
\makeatletter
\def\cleardoublepage{\clearpage\if@twoside \ifodd\c@page\else
    \if@twocolumn\hbox{}\newpage\fi\fi\fi
  }
\makeatother \clearpage{\pagestyle{plain}\cleardoublepage
}

\usepackage{color}
\usepackage{tikz}
\usepackage[explicit]{titlesec}
\titleformat{\chapter}[hang]{\normalfont\bfseries\huge}{\thechapter}{0.5em}{#1}

\titlespacing*{\chapter}{0pt}{0pt}{*0} 
\titlespacing*{\section}{0pt}{13.2pt}{*0}  
\titlespacing*{\subsection}{0pt}{13.2pt}{*0}
\titlespacing*{\subsubsection}{0pt}{13.2pt}{*0}

\usepackage{amsmath}
\usepackage{amsthm}
\usepackage{amsfonts}
\usepackage{amssymb}
\usepackage{euscript}
\usepackage{fancyhdr}
\usepackage{mathtools}
\usepackage{pifont}
\usepackage{mathrsfs}

\usepackage{enumerate}
\usepackage{enumitem}

\usepackage[nottoc, notlof, notlot]{tocbibind}

\newtheorem{theorem}{Theorem}[chapter]
\newtheorem{lemma}[theorem]{Lemma}
\newtheorem{proposition}[theorem]{Proposition}
\newtheorem{corollary}[theorem]{Corollary}
\newtheorem{definition}[theorem]{Definition}
\newtheorem{remark}[theorem]{Remark}

\newtheorem{BA}[theorem]{Bootstrap Assumption}

\newtheorem*{conjecture*}{Conjecture}
\newtheorem*{question*}{Question}

\numberwithin{equation}{chapter}

\makeatletter
\newcommand\footnoteref[1]{\protected@xdef\@thefnmark{\ref{#1}}\@footnotemark}
\makeatother

\usepackage{chngcntr}
\usepackage{etoolbox}

\AtBeginEnvironment{subappendices}{%
\chapter*{Appendices}
\addcontentsline{toc}{chapter}{Appendices}
\counterwithin{figure}{section}
\counterwithin{table}{section}
}

\AtEndEnvironment{subappendices}{%
\counterwithout{figure}{section}
\counterwithout{table}{section}
\setcounter{chapter}{0}
}
\AtEndEnvironment{appendices}{%
\counterwithout{figure}{section}
\counterwithout{table}{section}
\setcounter{chapter}{0}
}

\hypersetup{
colorlinks=true,
linkcolor=MidnightBlue,
citecolor=MidnightBlue,
linktocpage=true
}

\def\a{{\alpha}}
\def\al{\alpha}

\def\be{{\beta}}
\def\ga{\gamma}
\def\Ga{\Gamma}
\def\de{\delta}

\def\ep{\epsilon}

\def\la{\lambda}

\def\si{\sigma}
\def\Si{\Sigma}
\def\om{\omega}
\def\Om{\Omega}
\def\th{\theta}
\def\Th{\Theta}
\def\ze{\zeta}

\def\CC{{\mathcal C}}
\def\MM{{\mathcal M}}
\def\NN{{\mathcal N}}

\def\FF{{\mathcal F}}
\def\EE{{\mathcal E}}
\def\HH{{\mathcal H}}
\def\LL{{\mathcal L}}

\def\TT{{\mathcal T}}

\def\OO{{\mathcal O}}

\def\NN{{\mathcal N}}

\def\UU{{\mathcal U}}

\def\Lie{{\mathcal L}}

\def\DD{{\mathcal D}}

\def\RR{{\mathcal R}}
\def\QQ{{\mathcal Q}}

\def\VV{{\mathcal V}}

\def\Lie{{\mathcal L}}

\def\D{{\bf D}}

\def\M{{\bf M}}

\def\R{{\bf R}}

\def\g{{\bf g}}

\def\SSS{{\mathbb S}}
\def\RRR{{\mathbb R}}

\def\Yb{\underline{Y}}

\def\Kb{{\,\underline K}}
\def\chih{{\hat \chi}}
\def\chib{{\underline \chi}}
\def\chibh{{\underline{\chih}}}

\def\etab{{\underline \eta}}
\def\omb{{\underline{\om}}}

\def\alphab{{\underline{\alpha}}}
\def\betab{{\underline{\beta}}}
\def\xib{{\underline{\xi}}}

\newcommand{\thh}{\hat \theta}

\newcommand{\Divd}{\Div \mkern-17mu /\ \,}
\newcommand{\Curld}{\Curl \mkern-17mu /\ \,}
\newcommand{\Nd}{\nabla \mkern-13mu /\ }

\newcommand{\Ld}{\triangle \mkern-12mu /\ }
\newcommand{\iin}{\in \mkern-16mu /\ \mkern-5mu}
\newcommand{\gd}{{g \mkern-8mu /\ \mkern-5mu }}

\newcommand{\cDd}{\mathcal{D} \mkern-11mu /\  \mkern-3mu}
\newcommand{\tr}{{\mathrm{tr}}}
\newcommand{\nab}{\nabla}
\providecommand{\pr}{{\partial}}
\providecommand{\varep}{{\varepsilon}}
\providecommand{\les}{\lesssim}
\providecommand{\half}{\frac{1}{2}}

\providecommand{\otimesh}{{\hat{\otimes}}}
\providecommand{\XE}{{\mathbf{X}^\ext}}

\providecommand{\Ga}{{\Gamma}}
\providecommand{\Fslash}{{F\mkern-9mu/\ \mkern-7mu}}
\providecommand{\Gslash}{{G\mkern-9mu/\ \mkern-7mu}}
\providecommand{\Eslash}{{E\mkern-9mu/\ \mkern-7mu}}

\providecommand{\psl}{{p \mkern-9mu /\ \mkern-7mu}}
\providecommand{\Omt}{{\widetilde{\Om}}}
\providecommand{\elt}{{\tilde{e}_4}}
\providecommand{\elbt}{{\tilde{e}_3}}
\providecommand{\ut}{{\tilde{u}}}
\providecommand{\ubt}{{\tilde{\ub}}}

\providecommand{\chit}{{\widetilde{\chi}}}
\providecommand{\chibt}{{\widetilde{\chib}}}
\providecommand{\etat}{{\tilde{\eta}}}
\providecommand{\zet}{{\tilde{\ze}}}
\providecommand{\etabt}{{\tilde{\etab}}}

\providecommand{\omt}{{\tilde{\om}}}
\providecommand{\ombt}{{\tilde{\omb}}}
\providecommand{\xit}{{\tilde{\xi}}}
\providecommand{\xibt}{{\tilde{\xib}}}
\providecommand{\chiht}{{\widehat{\widetilde{\chi}}}}
\providecommand{\chibht}{{\widehat{\widetilde{\underline{\chi}}}}}
\providecommand{\trchit}{{\tr \chit}}

\providecommand{\trchibt}{{\tr \chibt}}
\providecommand{\alphat}{{\tilde{\alpha}}}
\providecommand{\betat}{{\tilde{\beta}}}
\providecommand{\rhot}{{\tilde{\rho}}}
\providecommand{\sigmat}{{\tilde{\sigma}}}
\providecommand{\betabt}{{\tilde{\betab}}}
\providecommand{\alphabt}{{\tilde{\alphab}}}
\providecommand{\Stt}{{\widetilde{S}}}
\providecommand{\rhoot}{{\overline{\tilde{\rho}}}}
\providecommand{\sigmaot}{{\overline{\tilde{\sigma}}}}

\providecommand{\tr}{\mathrm{tr}}
\providecommand{\trch}{\tr \chi}
\providecommand{\trchb}{\tr \chib}
\providecommand{\trchi}{\trch}
\providecommand{\trchib}{\trchb}
\providecommand{\RRb}{{\underline{\mathcal{R}}}}

\providecommand{\norm}[1]{\le \Vert#1\ri \Vert}

\providecommand{\Div}{{\mathrm{div}\,}}
\providecommand{\Curl}{{\mathrm{curl}\,}}
\providecommand{\DIV}{{\text{\bf div}}}

\providecommand{\dual}{{^*}} 
\providecommand{\RRRic}{\mathrm{Ric}}

\providecommand{\ub}{{\underline{u}}}

\providecommand{\varth}{{\vartheta}}
\providecommand{\le}{\left}
\providecommand{\ri}{\right}

\providecommand{\el}{{e_4}}
\providecommand{\elb}{{e_3}}
\providecommand{\ea}{{e_a}}
\providecommand{\eb}{{e_b}}
\providecommand{\ec}{{e_c}}

\providecommand{\cc}{{\tau}}
\providecommand{\yy}{{\mathfrak{y}}}
\providecommand{\d}{{\mathrm{d}}}
\providecommand{\fb}{{\underline{f}}}
\providecommand{\Err}{{\mathrm{Err}}}
\providecommand{\lot}{l.o.t.}
\providecommand{\quar}{\frac{1}{4}}
\providecommand{\beb}{{\betab}}
\providecommand{\alb}{{\alphab}}
\providecommand{\dg}{{\dagger}}
\providecommand{\ddg}{{\ddagger}}
\providecommand{\Lieh}{{\hat \Lie}}
\providecommand{\trchibo}{{\overline{\trchib}}}
\providecommand{\trchio}{{\overline{\trchi}}}

\providecommand{\nulld}{{\text{null}}}

\providecommand{\CCba}{{\CCb^\ast}}

\providecommand{\uba}{{\ub^\ast}}
\providecommand{\bott}{{\mathrm{bot}}} 
\providecommand{\topp}{{\mathrm{top}}} 
\providecommand{\Dd}{{\cDd}}
\providecommand{\Pb}{{\underline{P}}}
\providecommand{\Mb}{{\underline{M}}}
\providecommand{\Nb}{{\underline{N}}}
\providecommand{\Qb}{{\underline{Q}}}
\providecommand{\RRf}{{\frak{R}}}
\providecommand{\RRfb}{{\underline{\frak{R}}}}
\providecommand{\LLc}{{\LL_{\mathrm{con}}}}
\providecommand{\LLcext}{{\LL^\ext_{\mathrm{con}}}}
\providecommand{\LLbint}{{\LL^\intr_{\mathrm{bot}}}}

\providecommand{\LLb}{{\LL_\bott}}
\providecommand{\muo}{{\overline{\mu}}}

\providecommand{\ext}{{\mathrm{ext}}}
\providecommand{\intr}{{\mathrm{int}}}

\providecommand{\Lied}{\mathcal{L} \mkern-9mu/\ \mkern-7mu}
\providecommand{\Liedh}{{\hat{\Lied}}}
\providecommand{\pih}{{\hat{\pi}}}
\providecommand{\ibf}{{\mathbf{i}}}
\providecommand{\jbf}{{\mathbf{j}}}
\providecommand{\mbf}{{\mathbf{m}}}
\providecommand{\mbbf}{{\underline{\mathbf{m}}}}
\providecommand{\nbf}{{\mathbf{n}}}
\providecommand{\nbbf}{{\underline{\mathbf{n}}}}
\providecommand{\DDf}{{\mathfrak{D}}}

\providecommand{\Gac}{{\check{\Ga}}}

\providecommand{\qq}{q}

\providecommand{\Nf}{{\overline{N}}} 

\providecommand{\kt}{{k}}

\providecommand{\TX}{{\mathbf{T}}}
\providecommand{\SX}{{\mathbf{S}}}
\providecommand{\KX}{{\mathbf{K}}}
\providecommand{\TI}{{\TX^\intr}}
\providecommand{\SI}{{\SX^\intr}}
\providecommand{\KI}{{\KX^\intr}}
\providecommand{\TE}{{\TX^\ext}}
\providecommand{\SE}{{\SX^\ext}}
\providecommand{\KE}{{\KX^\ext}}
\providecommand{\nt}{{n}}

\providecommand{\Tf}{{\overline{T}}}

\providecommand{\Et}{{E}}
\providecommand{\Ht}{{H}}

\providecommand{\nut}{{\nu}}

\providecommand{\Ntf}{{\Nf}}

\providecommand{\delt}{{\de}}
\providecommand{\ept}{{\ep}}
\providecommand{\kapt}{{\kappa}}

\providecommand{\trth}{{\tr\th}}

\providecommand{\sigmac}{{\check{\sigma}}}

\providecommand{\Lambdab}{{\underline{\Lambda}}}
\providecommand{\Kbb}{{\underline{K}}}
\providecommand{\Xib}{{\underline{\Xi}}}
\providecommand{\Ib}{{\underline{I}}}
\providecommand{\Thetab}{{\underline{\Theta}}}
\providecommand{\OOI}{{\mathbf{O}^\intr}}
\providecommand{\OOE}{{\mathbf{O}^\ext}}

\providecommand{\XI}{{\mathbf{X}^\intr}}
\providecommand{\tast}{{t^\ast}}
\providecommand{\OOO}{{\mathbf{O}}}

\providecommand{\rhoo}{{\overline{\rho}}}
\providecommand{\sigmao}{{\overline{\sigma}}}

\providecommand{\OOEi}{{^{(\ell)}\OOE}}
\providecommand{\YEi}{{^{(i)}Y}}

\providecommand{\HEi}{{^{(i)}H}}
\providecommand{\QQb}{{\underline{\QQ}}}

\providecommand{\ZZt}{{\tilde{Z}}}
\providecommand{\ZZterr}{{\ZZt^\Err}}
\providecommand{\Xf}{{\overline{X}}}
\providecommand{\zdot}{{\dot z}}

\providecommand{\Psibar}{{\underline{\Psi}}}
\providecommand{\vdot}{{\dot v}}
\providecommand{\Id}{{\mathrm{Id}}}

\providecommand{\mubt}{{\tilde{\mu}}}
\providecommand{\nubt}{{\tilde{\nu}}}
\providecommand{\abt}{{\tilde{a}}}
\providecommand{\bbt}{{\tilde{b}}}
\providecommand{\cbt}{{\tilde{c}}}
\providecommand{\prt}{{\tilde{\pr}}}

\providecommand{\etadt}{{\eta \mkern-8mu /\ \mkern-5mu }}
\providecommand{\OS}{{\mathscr{O}}}
\providecommand{\pdot}{{\dot p}}

\providecommand{\Nds}{{\dual\Nd}}
\providecommand{\io}{{\iota}}
\providecommand{\ombs}{{\underline{\om}_\sigma}}
\providecommand{\ombr}{{\underline{\om}_\rho}}
\providecommand{\ombo}{{\overline{\omb}}}

\providecommand{\Pso}{{\varsigma}}
\providecommand{\HHt}{{\tilde{H}^{1/2}}}
\providecommand{\iob}{{\underline{\io}}}
\providecommand{\ombd}{{\dual\omb}}
\providecommand{\OOfgoodext}{{\mathfrak{O}^{\ext,\mathfrak{g}}_{\leq 1}}}
\providecommand{\OOfbadext}{{\mathfrak{O}^{\ext,\mathfrak{b}}_{\leq 1}}}
\providecommand{\OOgoodext}{{\OO^{\ext,\mathfrak{g}}_{\leq 2,\ga}}}
\providecommand{\OObadext}{{\OO^{\ext,\mathfrak{b}}_{\leq 2,\ga}}}
\providecommand{\OOfext}{{\mathfrak{O}^{\ext}_{\leq 1}}}
\providecommand{\OOext}{{\OO^{\ext}_{\leq 2,\ga}}}
\providecommand{\OOfexti}{{\mathfrak{O}^{\ext}_{\leq 1}}}
\providecommand{\OOfextii}{{\mathfrak{O}^{\ext}_{\leq 2}}}
\providecommand{\OOexti}{{\OO^{\ext}_{\leq 2,\ga}}}
\providecommand{\OOextii}{{\OO^\ext_{\leq 3,\ga}}}
\providecommand{\OOof}{{\overline{\OO}^{\ext}_{\leq 2,\ga}}}
\providecommand{\OOofb}{{\overline{\mathfrak{O}}^\ext_{\leq 1}}}

\providecommand{\yyo}{{\overline{\yy}}}
\providecommand{\OOastgood}{{\OO^{\ast, \mathfrak{g}}_{\leq 3}}}
\providecommand{\OOastbad}{{\OO^{\ast,\mathfrak{b}}_{\leq 3}}}
\providecommand{\OOastbadbad}{{\OO^{\ast,\mathfrak{b}}_{\leq 2+}}}
\providecommand{\OOast}{{\OO^\ast_{\leq 3}}}
\providecommand{\OOfastgood}{{\mathfrak{O}^{\ast, \mathfrak{g}}_{\leq 2}}}
\providecommand{\OOfastbad}{{\mathfrak{O}^{\ast,\mathfrak{b}}_{\leq 2}}}
\providecommand{\OOfastbadbad}{{\mathfrak{O}^{\ast,\mathfrak{b}}_{\leq 1+}}}
\providecommand{\OOfast}{{\mathfrak{O}^\ast_{\leq 2}}}
\providecommand{\OOoast}{{\overline{\mathfrak{O}}^\ast_{\leq 2}}}
\providecommand{\RRast}{{\RR^{\ast}_{\leq 2}}}
\providecommand{\RRfast}{{\mathfrak{R}^\ast_{\leq 1}}}
\providecommand{\RRfoast}{{\overline{\mathfrak{R}}^\ast_{\leq 2}}}
\providecommand{\RRext}{{\RR^\ext_{\leq 2,\ga}}}
\providecommand{\RRfext}{{\mathfrak{R}^\ext_{\leq 1}}}
\providecommand{\RRfoext}{{\overline{\mathfrak{R}}^\ext_{\leq 2}}}

\providecommand{\Sit}{{\widetilde{\Si}}}
\providecommand{\Sitt}{{\overline{\Si}}}
\providecommand{\POE}{{\Psi}}
\providecommand{\FFF}{{\mathrm{F}}}

\providecommand{\gao}{{\ga_0}}
\providecommand{\CCbb}{{\overline{\CCb}}}
\providecommand{\rast}{{r^\ast}}
\providecommand{\too}{{t^\circ}}

\providecommand{\etabold}{{\boldsymbol{\eta}}}

\providecommand{\ft}{{\tilde{f}}}
\providecommand{\fbt}{{\tilde{\fb}}}
\providecommand{\lat}{{\tilde{\la}}}
\providecommand{\rt}{{\tilde{r}}}
\providecommand{\LLbext}{{\LL_\bott^\ext}}
\providecommand{\EEE}{{\mathrm{E}}}
\providecommand{\CCt}{{\widetilde{\mathcal{C}}}}

\providecommand{\Ndt}{{\boldsymbol{\Nd}}}
\providecommand{\Ndtt}{{\widetilde{\Ndt}}}
\providecommand{\DDD}{{\mathbb{D}}}
\providecommand{\doto}{{\dot{\o}}}
\providecommand{\CCb}{{\underline{\CC}}}
\providecommand{\MMt}{{\widetilde{\mathcal{M}}}}
\providecommand{\scri}{{\mathscr{I}}}
\providecommand{\GGG}{{\mathrm{G}}}
\providecommand{\BBB}{{\mathrm{B}}}

\renewcommand{\Ga}{{\Gamma}}
\renewcommand{\Fslash}{{F\mkern-9mu/\ \mkern-7mu}}
\renewcommand{\Gslash}{{G\mkern-9mu/\ \mkern-7mu}}
\renewcommand{\Eslash}{{E\mkern-9mu/\ \mkern-7mu}}

\renewcommand{\psl}{{p \mkern-9mu /\ \mkern-7mu}}
\renewcommand{\Omt}{{\widetilde{\Om}}}
\renewcommand{\elt}{{\tilde{e}_4}}
\renewcommand{\elbt}{{\tilde{e}_3}}
\renewcommand{\ut}{{\tilde{u}}}
\renewcommand{\ubt}{{\tilde{\ub}}}

\renewcommand{\chit}{{\widetilde{\chi}}}
\renewcommand{\chibt}{{\widetilde{\chib}}}
\renewcommand{\etat}{{\tilde{\eta}}}
\renewcommand{\zet}{{\tilde{\ze}}}
\renewcommand{\etabt}{{\tilde{\etab}}}

\renewcommand{\omt}{{\tilde{\om}}}
\renewcommand{\ombt}{{\tilde{\omb}}}
\renewcommand{\xit}{{\tilde{\xi}}}
\renewcommand{\xibt}{{\tilde{\xib}}}
\renewcommand{\chiht}{{\widehat{\widetilde{\chi}}}}
\renewcommand{\chibht}{{\widehat{\widetilde{\underline{\chi}}}}}
\renewcommand{\trchit}{{\tr \chit}}

\renewcommand{\trchibt}{{\tr \chibt}}
\renewcommand{\alphat}{{\tilde{\alpha}}}
\renewcommand{\betat}{{\tilde{\beta}}}
\renewcommand{\rhot}{{\tilde{\rho}}}
\renewcommand{\sigmat}{{\tilde{\sigma}}}
\renewcommand{\betabt}{{\tilde{\betab}}}
\renewcommand{\alphabt}{{\tilde{\alphab}}}
\renewcommand{\Stt}{{\widetilde{S}}}
\renewcommand{\rhoot}{{\overline{\tilde{\rho}}}}
\renewcommand{\sigmaot}{{\overline{\tilde{\sigma}}}}

\renewcommand{\tr}{\mathrm{tr}}
\renewcommand{\trch}{\tr \chi}
\renewcommand{\trchb}{\tr \chib}
\renewcommand{\trchi}{\trch}
\renewcommand{\trchib}{\trchb}
\renewcommand{\RRb}{{\underline{\mathcal{R}}}}

\renewcommand{\norm}[1]{\le \Vert#1\ri \Vert}

\renewcommand{\dual}{{^*}} 
\renewcommand{\RRRic}{\mathrm{Ric}}

\renewcommand{\ub}{{\underline{u}}}

\renewcommand{\varth}{{\vartheta}}
\renewcommand{\le}{\left}
\renewcommand{\ri}{\right}

\renewcommand{\el}{{e_4}}
\renewcommand{\elb}{{e_3}}
\renewcommand{\ea}{{e_a}}
\renewcommand{\eb}{{e_b}}
\renewcommand{\ec}{{e_c}}

\renewcommand{\cc}{{\tau}}
\renewcommand{\yy}{{\mathfrak{y}}}
\renewcommand{\d}{{\mathrm{d}}}
\renewcommand{\fb}{{\underline{f}}}
\renewcommand{\Err}{{\mathrm{Err}}}
\renewcommand{\lot}{l.o.t.}
\renewcommand{\quar}{\frac{1}{4}}
\renewcommand{\beb}{{\betab}}
\renewcommand{\alb}{{\alphab}}
\renewcommand{\dg}{{\dagger}}
\renewcommand{\ddg}{{\ddagger}}
\renewcommand{\Lieh}{{\hat \Lie}}
\renewcommand{\trchibo}{{\overline{\trchib}}}
\renewcommand{\trchio}{{\overline{\trchi}}}

\renewcommand{\nulld}{{\text{null}}}

\renewcommand{\CCba}{{\CCb^\ast}}

\renewcommand{\uba}{{\ub^\ast}}
\renewcommand{\bott}{{\mathrm{bot}}} 
\renewcommand{\topp}{{\mathrm{top}}} 
\renewcommand{\Dd}{{\cDd}}
\renewcommand{\Pb}{{\underline{P}}}
\renewcommand{\Mb}{{\underline{M}}}
\renewcommand{\Nb}{{\underline{N}}}
\renewcommand{\Qb}{{\underline{Q}}}
\renewcommand{\RRf}{{\frak{R}}}
\renewcommand{\RRfb}{{\underline{\frak{R}}}}
\renewcommand{\LLc}{{\LL_{\mathrm{con}}}}
\renewcommand{\LLcext}{{\LL^\ext_{\mathrm{con}}}}
\renewcommand{\LLbint}{{\LL^\intr_{\mathrm{bot}}}}

\renewcommand{\LLb}{{\LL_\bott}}
\renewcommand{\muo}{{\overline{\mu}}}

\renewcommand{\ext}{{\mathrm{ext}}}
\renewcommand{\intr}{{\mathrm{int}}}

\renewcommand{\Lied}{\mathcal{L} \mkern-9mu/\ \mkern-7mu}
\renewcommand{\Liedh}{{\hat{\Lied}}}
\renewcommand{\pih}{{\hat{\pi}}}
\renewcommand{\ibf}{{\mathbf{i}}}
\renewcommand{\jbf}{{\mathbf{j}}}
\renewcommand{\mbf}{{\mathbf{m}}}
\renewcommand{\mbbf}{{\underline{\mathbf{m}}}}
\renewcommand{\nbf}{{\mathbf{n}}}
\renewcommand{\nbbf}{{\underline{\mathbf{n}}}}
\renewcommand{\DDf}{{\mathfrak{D}}}

\renewcommand{\Gac}{{\check{\Ga}}}

\renewcommand{\qq}{q}

\renewcommand{\Nf}{{\overline{N}}} 

\renewcommand{\kt}{{k}}

\renewcommand{\TX}{{\mathbf{T}}}
\renewcommand{\SX}{{\mathbf{S}}}
\renewcommand{\KX}{{\mathbf{K}}}
\renewcommand{\TI}{{\TX^\intr}}
\renewcommand{\SI}{{\SX^\intr}}
\renewcommand{\KI}{{\KX^\intr}}
\renewcommand{\TE}{{\TX^\ext}}
\renewcommand{\SE}{{\SX^\ext}}
\renewcommand{\KE}{{\KX^\ext}}
\renewcommand{\nt}{{n}}

\renewcommand{\Tf}{{\overline{T}}}

\renewcommand{\Et}{{E}}
\renewcommand{\Ht}{{H}}

\renewcommand{\nut}{{\nu}}

\renewcommand{\Ntf}{{\Nf}}

\renewcommand{\delt}{{\de}}
\renewcommand{\ept}{{\ep}}
\renewcommand{\kapt}{{\kappa}}

\renewcommand{\trth}{{\tr\th}}

\renewcommand{\sigmac}{{\check{\sigma}}}

\renewcommand{\Lambdab}{{\underline{\Lambda}}}
\renewcommand{\Kbb}{{\underline{K}}}
\renewcommand{\Xib}{{\underline{\Xi}}}
\renewcommand{\Ib}{{\underline{I}}}
\renewcommand{\Thetab}{{\underline{\Theta}}}
\renewcommand{\OOI}{{\mathbf{O}^\intr}}
\renewcommand{\OOE}{{\mathbf{O}^\ext}}

\renewcommand{\XI}{{\mathbf{X}^\intr}}
\renewcommand{\tast}{{t^\ast}}
\renewcommand{\OOO}{{\mathbf{O}}}

\renewcommand{\rhoo}{{\overline{\rho}}}
\renewcommand{\sigmao}{{\overline{\sigma}}}

\renewcommand{\OOEi}{{^{(\ell)}\OOE}}
\renewcommand{\YEi}{{^{(i)}Y}}

\renewcommand{\HEi}{{^{(i)}H}}
\renewcommand{\QQb}{{\underline{\QQ}}}

\renewcommand{\ZZt}{{\tilde{Z}}}
\renewcommand{\ZZterr}{{\ZZt^\Err}}
\renewcommand{\Xf}{{\overline{X}}}
\renewcommand{\zdot}{{\dot z}}

\renewcommand{\Psibar}{{\underline{\Psi}}}
\renewcommand{\vdot}{{\dot v}}
\renewcommand{\Id}{{\mathrm{Id}}}

\renewcommand{\mubt}{{\tilde{\mu}}}
\renewcommand{\nubt}{{\tilde{\nu}}}
\renewcommand{\abt}{{\tilde{a}}}
\renewcommand{\bbt}{{\tilde{b}}}
\renewcommand{\cbt}{{\tilde{c}}}
\renewcommand{\prt}{{\tilde{\pr}}}

\renewcommand{\etadt}{{\eta \mkern-8mu /\ \mkern-5mu }}
\renewcommand{\OS}{{\mathscr{O}}}
\renewcommand{\pdot}{{\dot p}}

\renewcommand{\Nds}{{\dual\Nd}}
\renewcommand{\io}{{\iota}}
\renewcommand{\ombs}{{\underline{\om}_\sigma}}
\renewcommand{\ombr}{{\underline{\om}_\rho}}
\renewcommand{\ombo}{{\overline{\omb}}}

\renewcommand{\Pso}{{\varsigma}}
\renewcommand{\HHt}{{\tilde{H}^{1/2}}}
\renewcommand{\iob}{{\underline{\io}}}
\renewcommand{\ombd}{{\dual\omb}}
\renewcommand{\OOfgoodext}{{\mathfrak{O}^{\ext,\mathfrak{g}}_{\leq 1}}}
\renewcommand{\OOfbadext}{{\mathfrak{O}^{\ext,\mathfrak{b}}_{\leq 1}}}
\renewcommand{\OOgoodext}{{\OO^{\ext,\mathfrak{g}}_{\leq 2,\ga}}}
\renewcommand{\OObadext}{{\OO^{\ext,\mathfrak{b}}_{\leq 2,\ga}}}
\renewcommand{\OOfext}{{\mathfrak{O}^{\ext}_{\leq 1}}}
\renewcommand{\OOext}{{\OO^{\ext}_{\leq 2,\ga}}}
\renewcommand{\OOfexti}{{\mathfrak{O}^{\ext}_{\leq 1}}}
\renewcommand{\OOfextii}{{\mathfrak{O}^{\ext}_{\leq 2}}}
\renewcommand{\OOexti}{{\OO^{\ext}_{\leq 2,\ga}}}
\renewcommand{\OOextii}{{\OO^\ext_{\leq 3,\ga}}}
\renewcommand{\OOof}{{\overline{\OO}^{\ext}_{\leq 2,\ga}}}
\renewcommand{\OOofb}{{\overline{\mathfrak{O}}^\ext_{\leq 1}}}

\renewcommand{\yyo}{{\overline{\yy}}}
\renewcommand{\OOastgood}{{\OO^{\ast, \mathfrak{g}}_{\leq 3}}}
\renewcommand{\OOastbad}{{\OO^{\ast,\mathfrak{b}}_{\leq 3}}}
\renewcommand{\OOastbadbad}{{\OO^{\ast,\mathfrak{b}}_{\leq 2+}}}
\renewcommand{\OOast}{{\OO^\ast_{\leq 3}}}
\renewcommand{\OOfastgood}{{\mathfrak{O}^{\ast, \mathfrak{g}}_{\leq 2}}}
\renewcommand{\OOfastbad}{{\mathfrak{O}^{\ast,\mathfrak{b}}_{\leq 2}}}
\renewcommand{\OOfastbadbad}{{\mathfrak{O}^{\ast,\mathfrak{b}}_{\leq 1+}}}
\renewcommand{\OOfast}{{\mathfrak{O}^\ast_{\leq 2}}}
\renewcommand{\OOoast}{{\overline{\mathfrak{O}}^\ast_{\leq 2}}}
\renewcommand{\RRast}{{\RR^{\ast}_{\leq 2}}}
\renewcommand{\RRfast}{{\mathfrak{R}^\ast_{\leq 1}}}
\renewcommand{\RRfoast}{{\overline{\mathfrak{R}}^\ast_{\leq 2}}}
\renewcommand{\RRext}{{\RR^\ext_{\leq 2,\ga}}}
\renewcommand{\RRfext}{{\mathfrak{R}^\ext_{\leq 1}}}
\renewcommand{\RRfoext}{{\overline{\mathfrak{R}}^\ext_{\leq 2}}}

\renewcommand{\Sit}{{\widetilde{\Si}}}
\renewcommand{\Sitt}{{\overline{\Si}}}
\renewcommand{\POE}{{\Psi}}
\renewcommand{\FFF}{{\mathrm{F}}}

\renewcommand{\gao}{{\ga_0}}
\renewcommand{\CCbb}{{\overline{\CCb}}}
\renewcommand{\rast}{{r^\ast}}
\renewcommand{\too}{{t^\circ}}

\renewcommand{\etabold}{{\boldsymbol{\eta}}}

\renewcommand{\ft}{{\tilde{f}}}
\renewcommand{\fbt}{{\tilde{\fb}}}
\renewcommand{\lat}{{\tilde{\la}}}
\renewcommand{\rt}{{\tilde{r}}}
\renewcommand{\LLbext}{{\LL_\bott^\ext}}
\renewcommand{\EEE}{{\mathrm{E}}}
\renewcommand{\CCt}{{\widetilde{\mathcal{C}}}}

\renewcommand{\Ndt}{{\boldsymbol{\Nd}}}
\renewcommand{\Ndtt}{{\widetilde{\Ndt}}}
\renewcommand{\DDD}{{\mathbb{D}}}
\renewcommand{\doto}{{\dot{\o}}}
\renewcommand{\CCb}{{\underline{\CC}}}
\renewcommand{\MMt}{{\widetilde{\mathcal{M}}}}
\renewcommand{\scri}{{\mathscr{I}}}
\renewcommand{\GGG}{{\mathrm{G}}}
\renewcommand{\BBB}{{\mathrm{B}}}
\providecommand{\DIV}{{\boldsymbol{\mathrm{div}}}}
\newcommand{\cco}{{\cc_0}}
\providecommand{\Hrot}{{M}}
\providecommand{\HEi}{{^{(\ell)}M}}
\renewcommand{\HEi}{{^{(\ell)}M}}
\renewcommand{\Ndt}{{\mathfrak{d}}}
\newcommand{\vb}{{\underline{v}}}

\begin{document}

\frontmatter

\makeatletter
\begin{titlepage}
  \vspace*{\fill}
  \begin{center}
    \LARGE
    \@title \\
    \vspace{1.5cm}
    \Large
    \@author\\
    \vspace{1cm}
    \@date \\
    \vspace{2cm}
    \normalsize
    \textbf{Abstract}\\

    \begin{quotation}
      In this paper, we prove the global nonlinear stability of Minkowski space in the context of the spacelike-characteristic Cauchy problem for Einstein vacuum equations. Spacelike-characteristic initial data are posed on a compact 3-disk and on the future complete null hypersurface emanating from its boundary. Our result extends the seminal stability result for Minkowski space proved by Christodoulou and Klainerman for which initial data are prescribed on a spacelike 3-plane.

      The proof relies on the classical vectorfield method and bootstrapping argument from Christodoulou-Klainerman. The main novelty is the introduction and control of new geometric constructions adapted to the spacelike-characteristic setting. In particular, it features null cones with prescribed vertices, spacelike maximal hypersurfaces with prescribed boundaries and global harmonic coordinates on Riemannian 3-disks. 
    \end{quotation}

  \end{center}
  \vspace*{\fill}
\end{titlepage}
\makeatother


\tableofcontents



\mainmatter


\chapter{Introduction}
\section{Einstein equations and the stability of Minkowski space}\label{sec:IntroStabMinkCK}
A Lorentzian $4$-dimensional manifold $(\MM,\g)$ is called a \emph{vacuum spacetime} if it solves the Einstein vacuum equations
\begin{align}\label{eq:EVEintroSTAB}
\RRRic(\g) = 0,
\end{align}
where $\RRRic(\g)$ denotes the Ricci tensor of the Lorentzian metric $\g$. For the metric components  $\g_{\mu\nu}$ in general coordinates, equation~\eqref{eq:EVEintroSTAB} writes as a system of non-linear coupled partial differential equations of order $2$ for $\g_{\mu\nu}$. In so-called \emph{wave coordinates}, it can be shown that \eqref{eq:EVEintroSTAB} is a system of nonlinear wave equations. It therefore admits a well-posed initial value formulation (or \emph{Cauchy problem}).\\

Cauchy data for equations~\eqref{eq:EVEintroSTAB} are classically described by a triplet $(\Si,g,k)$ such that 
\begin{itemize}
\item $(\Si,g)$ is a $3$-dimensional Riemannian manifold,
\item $k$ is a symmetric covariant $2$-tensor on $\Si$,
\item $(g,k)$ satisfy so-called constraint equations on $\Si$.\footnote{See~\cite[equations (10.2.28), (10.2.30)]{Wal84}.}
\end{itemize}

The seminal well-posedness results for the Cauchy problem obtained in~\cite{Fou52,Cho.Ger69} ensure that for any smooth Cauchy data, there exists a unique smooth \emph{maximal globally hyperbolic development} $(\MM,\g)$ solution of Einstein equations~\eqref{eq:EVEintroSTAB} such that $\Si\subset \MM$ and $g,k$ are respectively the first and second fundamental forms of $\Si$ in $\MM$. We refer to~\cite[Chapter 10]{Wal84} for definitions and further discussions on the Cauchy problem in general relativity. 

\begin{remark}
  Here and in the rest of this paper, a \emph{smooth} or $\mathscr{C}^\infty$ manifold admits by definition an atlas of charts such that all coordinates changes are $\mathscr{C}^\infty$ with respect to the standard $\mathscr{C}^\infty$-topology for functions of $\RRR^n$. As all manifolds we consider will be smooth submanifolds of a fixed smooth $4$-dimensional manifold $\MM$ and as all vector bundles we consider will be constructed upon $\mathrm{T}\MM$ and $\mathrm{T}^\ast\MM$, we shall assume that such an atlas is fixed on $\MM$, which then canonically determines the $\mathscr{C}^k$-topology for all tensors on all smooth submanifolds of $\MM$ in this paper.
\end{remark}

The prime example of a vacuum spacetime is \emph{Minkowski space}
\begin{align*}
  \MM & = \RRR^4, & \g & = -(\d x^0)^2+ (\d x^1)^2 + (\d x^2)^2 + (\d x^3)^2 =: \etabold,
\end{align*}
for which Cauchy data are given by
\begin{align*}
  \Si & = \RRR^3, & g & = (\d x^1)^2 + (\d x^2)^2 + (\d x^3)^2 =: \de, & k & = 0.
\end{align*}



In the breakthrough~\cite{Chr.Kla93}, it is shown that Minkowski space is an asymptotically stable solution of Einstein vacuum equations in the following sense (see also Theorem~\ref{thm:StabCKintroSTAB} for a precised version).
\begin{theorem}[Stability of Minkowski space~\cite{Chr.Kla93}, version 1]\label{thm:StabCKintroSTABrough}
  For Cauchy data $(\Si,g,k)$ such that
  \begin{itemize}
  \item $\Si$ is diffeomorphic to $\RRR^3$,
  \item $(\Si,g,k)$ is asymptotically flat (\emph{i.e.} tends to Minkowski initial data $(\RRR^3,\de,0)$ when $r\to\infty$)
  \item $(g,k)$ are close to Minkowski initial data $(\de,0)$ measured in an (weighted) $L^2$-sense,
  \end{itemize}
  then its maximal globally hyperbolic development $(\MM,\g)$ is \emph{geodesically complete} and admits global time and optical functions $t$ and $u$ such that, measured in these coordinates, $\g$ is bounded and decays towards $\etabold$.   
\end{theorem}

A localised version of Theorem~\ref{thm:StabCKintroSTABrough} was proved for initial data posed on the exterior of a $3$-disk.
\begin{theorem}[Exterior stability of Minkowski space~\cite{Kla.Nic03}, rough version]\label{thm:StabKNintroSTABrough}
  For Cauchy data $(\Si,g,k)$ such that
  \begin{itemize}
  \item $\Si$ is diffeomorphic to $\RRR^3\setminus\DDD$ where $\DDD$ denotes the disk of $\RRR^3$,
  \item the same asymptotic flatness and closeness to Minkowski space assumptions as in Theorem~\ref{thm:StabCKintroSTABrough} hold,
  \end{itemize}
  then, the maximal globally hyperbolic development $(\MM,\g)$ admits global optical functions $u,\ub$ such that, measured in these coordinates, $\g$ is bounded and decays towards $\etabold$.
\end{theorem}

\begin{remark}
  In the proof of Theorem~\ref{thm:StabCKintroSTABrough} in~\cite{Chr.Kla93}, the topology assumption $\Si \simeq \RRR^3$ is crucially used to define a \emph{global time function} such that its level sets are \emph{maximal hypersurfaces} with prescribed asymptotic conditions when $r\to\infty$. The main novelty in the proof~\cite{Kla.Nic03} of Theorem~\ref{thm:StabKNintroSTABrough} is the definition of a \emph{double-null foliation} by the level sets of two optical functions $u,\ub$. It replaces the global time function and enables a localisation of the global nonlinear stability proof~\cite{Chr.Kla93} to the exterior of a disk. See definitions and discussions in Section~\ref{sec:preciseSTABintro} and see also discussions in~\cite[Section 2]{Kla.Nic03}.  
\end{remark}

\begin{figure}[!h]
  \centering
  \includegraphics[width=7cm]{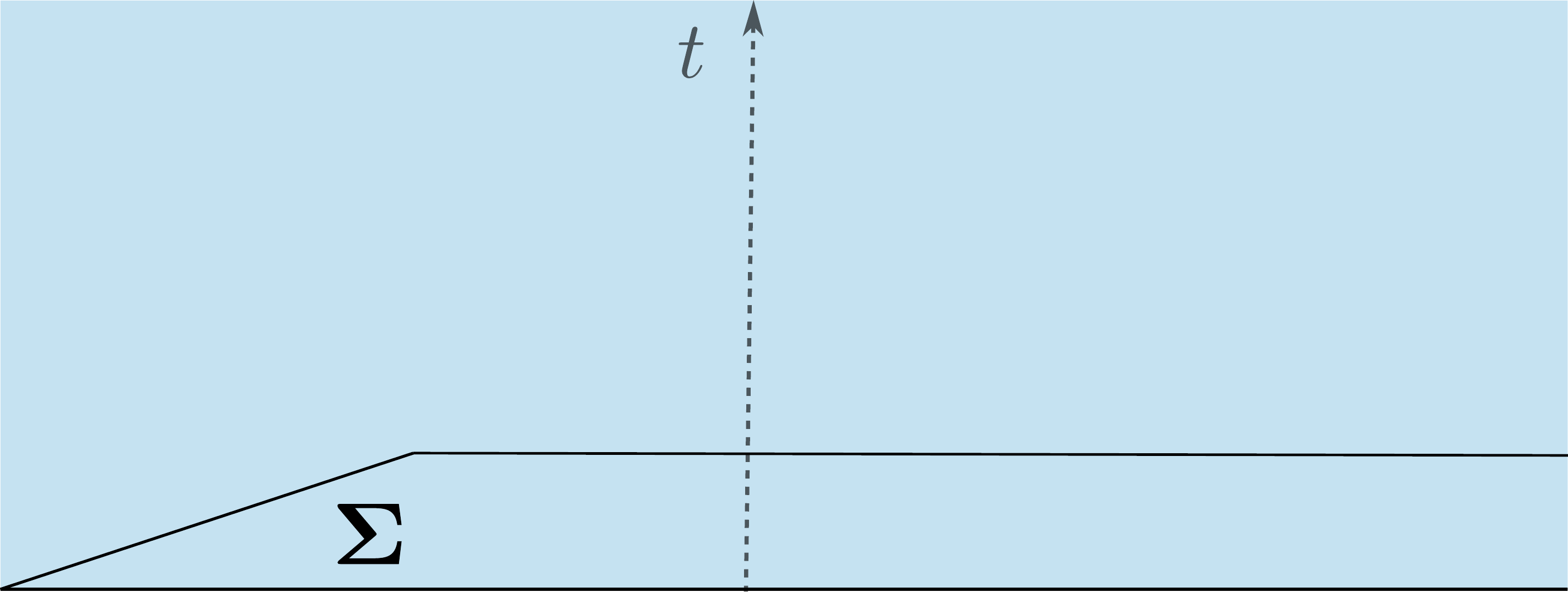}
  \vspace{0.5em}
  \includegraphics[width=7cm]{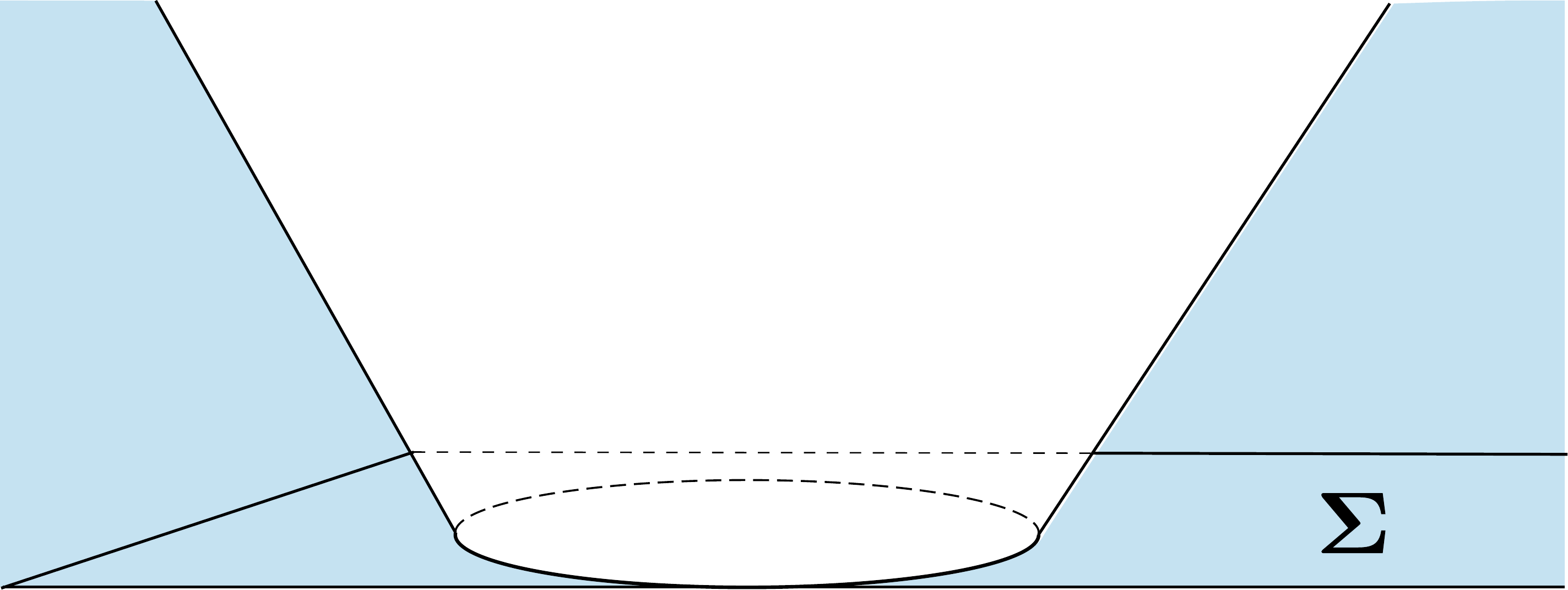}
  \caption{The stability of Minkowski space of Theorems~\ref{thm:StabCKintroSTABrough},~\ref{thm:StabKNintroSTABrough} proved in~\cite{Chr.Kla93} and~\cite{Kla.Nic03}.}
  \label{fig:Stabroughintro}
\end{figure}

In view of Theorems~\ref{thm:StabCKintroSTABrough} and~\ref{thm:StabKNintroSTABrough}, we have the following natural question.

\begin{question*}
Can we complete the result of Theorem~\ref{thm:StabKNintroSTABrough} to re-obtain the result of Theorem~\ref{thm:StabCKintroSTABrough}? In other words, can we prove the global nonlinear stability of Minkowski space for initial data posed on a spacelike disk and an outgoing null hypersurface?
\end{question*}

\begin{figure}[!h]
  \centering
    \includegraphics[width=8cm]{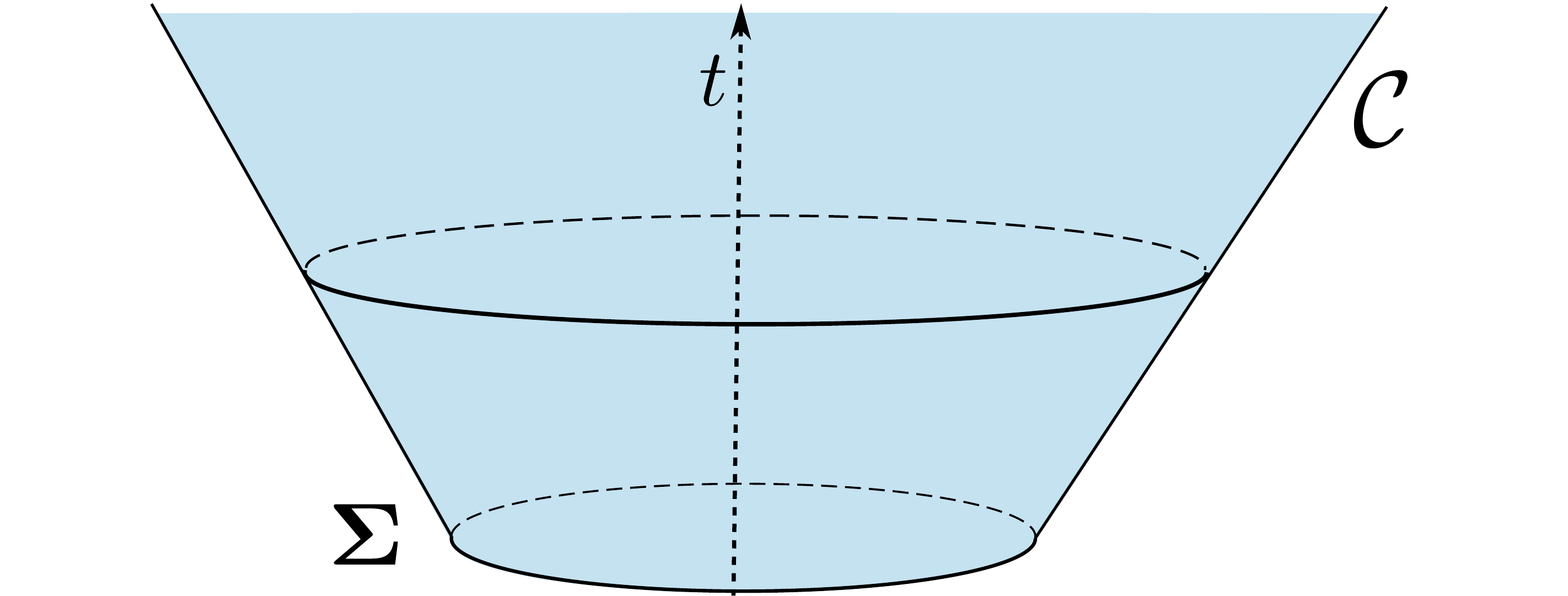}
  \caption{The stability of Minkowski space of Theorem~\ref{thm:mainv1rough} proved in the present paper.}
  \label{fig:StabroughintroMAIN}
\end{figure}

This paper is dedicated to the proof of the following theorem, which provides a positive answer to this question. See Theorem~\ref{thm:mainv1} for a more precise version, and Theorem~\ref{thm:mainv2} for the exact result proved in this paper.
\begin{theorem}[Main theorem, rough version]\label{thm:mainv1rough}
  Let initial data for Einstein equations~\eqref{eq:EVEintroSTAB} be given on
  \begin{itemize}
  \item an initial spacelike hypersurface $\Si$ diffeomorphic to $\DDD$ the unit disk of $\RRR^3$,
  \item an initial outgoing null hypersurface $\CC$ emanating from $\pr\Si$.
  \end{itemize}
  Assume that
  \begin{itemize}
  \item $\CC$ is future geodesically complete,
  \item the initial data are close to the corresponding Minkowski data on $\Si$ and $\CC$ consistently with the boundedness and decay assumptions and results of~\cite{Chr.Kla93} and~\cite{Kla.Nic03}.
  \end{itemize}
  Then the maximal globally hyperbolic development $(\MM,g)$ is future causally geodesically complete and admits global time and optical functions $t$ and $u$ such that, measured with respect to these coordinates, $\g$ is bounded and decays towards $\etabold$ when $t\to +\infty$.
\end{theorem}

\begin{remark}
  The so-called \emph{spacelike-characteristic} initial data of Theorem~\ref{thm:mainv1rough} have to satisfy constraints. In this paper, we do not discuss the prescription of such initial data. See~\cite{Chr.Pae12} or~\cite{Cac.Nic05} for discussions.
\end{remark}

\section{The \emph{global nonlinear stability of Minkowski space}}\label{sec:preciseSTABintro}
In this section, we give a more detailed statement of Theorem~\ref{thm:StabCKintroSTABrough} and highlight key features of its proof in~\cite{Chr.Kla93}. We first start with preliminary definitions.\\

For a fixed spacetime $(\MM,\g)$, we consider foliations of $\MM$ by the level sets $\Si_t$ a time function $t$. We write $T$ the future-pointing unit normal to $\Si_t$, and we recall that the second fundamental form $k$ and the \emph{time lapse} $n$ of $\Si_t\subset\MM$ are defined by
\begin{align*}
  k(X,Y) & := -\g(\D_XT,Y), & n^{-2} & := -\g(\D t,\D t),
\end{align*}
where $X,Y\in\mathrm{T}\Si_t$ and $\D$ denotes the spacetime covariant derivative. In this paper, we will consider \emph{maximal hypersurfaces} $\Si_t$, \emph{i.e.} such that
\begin{align*}
  \tr_gk & = 0,
\end{align*}
where $\tr_g$ is the trace with respect to the induced metric $g$ on $\Si$.\\

We also consider foliations of $\MM$ by $2$-spheres $S_{u,\ub}$, intersections of the level sets two functions $u,\ub$. A \emph{null pair} $(\elb,\el)$ adapted to $S_{u,\ub}$ is a pair of vectorfields orthogonal to the $2$-spheres $S_{u,\ub}$ which satisfy
\begin{align*}
  \g(\elb,\elb) & = \g(\el,\el) = 0, & \g(\elb,\el) & = -2.
\end{align*}
To a null pair, we associate \emph{null connection coefficients}, which are the $S_{u,\ub}$-tangent tensors such that
\begin{align}
  \label{eq:defnulldecompSTAB}
  \begin{aligned}
    \chi(X,Y) & := \g(\D_X\el,Y), & \xi(X)  & := \half\g(\D_4\el,X), & \eta(X) & := \half\g(\D_3\el,X),\\    \ze(X) & := \half \g(\D_X\el,\elb), & \om  & := \frac{1}{4} \g(\D_4\el,\elb), \\
     \chib(X,Y) & := \g(\D_X\elb,Y), & \xib(X) & := \half\g(\D_3\elb,X), & \etab(X) & := \half\g(\D_4\elb,X),\\  
     \ze(X) & := -\half \g(\D_X\elb,\el), & \omb & := \frac{1}{4} \g(\D_3\elb,\el),
  \end{aligned}
\end{align}
where $X,Y \in \text{T}S_{u,\ub}$, and the \emph{null curvature components}, which are the $S_{u,\ub}$-tangent tensors such that
\begin{align}
  \begin{aligned}
  \alpha(X,Y) & := \R(\el,X,\el,Y), & \beta(X) & := \half \R(X,\el,\elb,\el), & \rho & := \frac{1}{4} \R(\elb,\el,\elb,\el), \\
  \alphab(X,Y) & := \R(\elb,X,\elb,Y), & \betab(X) & := \half \R(X,\elb,\elb,\el), & \sigma & := \frac{1}{4}\dual\R(\elb,\el,\elb,\el),
  \end{aligned}
\end{align}
where $X,Y \in \text{T}S_{u,\ub}$ and where $\R$ is the spacetime curvature tensor and $\dual\R$ denotes its Hodge dual.\\

We have the following precised version of the stability result of Theorem~\ref{thm:StabCKintroSTABrough}.
\begin{theorem}[Stability of Minkoswki space~\cite{Chr.Kla93}, version 2]\label{thm:StabCKintroSTAB}
  Let $(\Si,g,k)$ be Cauchy data such that:
  \begin{subequations}\label{est:StabCKintroassumptions}
  \begin{itemize}
  \item $\Si$ is maximal, diffeomorphic to $\RRR^3$,
  \item $\Si$ is asymptotically flat, \emph{i.e.} there exists coordinates $(x^1,x^2,x^3)$ in a neighbourhood of infinity such that
    \begin{align}\label{est:AFStabCKintro}
      (r\pr)^{\leq 4}\le(g_{ij} - \le(1+\frac{2M}{r}\ri)\de_{ij}\ri) = O(r^{-3/2}), 
    \end{align}
    when $r\to\infty$ and where here $r := \sqrt{\sum_{i=1}^3(x^i)^2}$ and $M\geq 0$, and we have the following sup-norm bound for the curvature of $(\Si,g)$
    \begin{align}
      \norm{(1+d)^3\RRRic}_{L^\infty(\Si)} & \leq \varep,\label{est:curvfluxStabCKintro1}
    \end{align}
    where $\RRRic$ denotes the Ricci tensor of the metric $g$ and $d$ denotes the geodesic distance to a fixed point of $\Si$,
  \item the following bounds hold for curvature $L^2$-fluxes through $\Si$
    \begin{align}
      \norm{(1+d)\le((1+d)\nab\ri)^{\leq 3}k}_{L^2(\Si)} + \norm{(1+d)^3 \le((1+d)\nab\ri)^{\leq 1}\mathrm{B}}_{L^2(\Si)} & \leq \varep,  \label{est:curvfluxStabCKintro2}
    \end{align}
     where $\mathrm{B} := \Curl \le(\RRRic-\frac{1}{3}\mathrm{R}g\ri)$.
  \end{itemize}
  \end{subequations}
  Then, there exists $\varep_0>0$, such that if $\varep < \varep_0$, the following holds for the maximal globally hyperbolic development $(\MM,\g)$ of $(\Si,g,k)$.
  \begin{subequations}\label{est:StabCKintroresult}
  \begin{itemize}
  \item $(\MM,\g)$ is geodesically complete.
  \item There exists a global time function $t$ on $\MM$ ranging from $-\infty$ to $+\infty$ which foliates $\MM$ by maximal spacelike hypersurfaces $\Si_t$ such that $\Si_0=\Si$.
  \item There exists a future exterior region\footnote{The future exterior region is of the type $\{u\leq c t\}$ with $c<1$.} $\MM^\ext$ foliated by outgoing null hypersurfaces $\CC_u$ level sets of a global optical function $u$ ranging from $-\infty$ to $+\infty$ on $\MM^\ext$, and a past exterior region with symmetric constructions.
  \item We have the following decay in the interior region $\MM^\intr := \MM\setminus\MM^\ext$ of the spacetime curvature tensor $\R$ of $\g$\footnote{Here the norm of $\R$ is taken with respect to an orthonormal frame adapted to the maximal hypersurfaces $\Si_t$. See precise definitions in Section~\ref{sec:normsBA}.}
    \begin{align}\label{est:intdecayStabCKintro}
      |\R| & \les \varep t^{-7/2}.
    \end{align}
  \item We have the following differentiated decay in the exterior region $\MM^\ext$ of the spacetime curvature tensor $\R$ according to its null decomposition\footnote{Here the null pair is adapted to the $2$-spheres intersections of the maximal-null foliation of $\MM^\ext$ by $\Si_t$ and $\CC_u$. See~\cite[Introduction]{Chr.Kla93}}
    \begin{align}\label{est:extdecayStabCKintro}
      \begin{aligned}
      |\al(\R)| & \les \varep r^{-7/2}, & |\be(\R)| & \les \varep r^{-7/2}, & |\rho(\R)| & \les \varep r^{-3},\\
      |\alb(\R)| & \les \varep r^{-1}u^{-5/2}, & |\beb(\R)| & \les \varep r^{-2}u^{-3/2}, & |\sigma(\R)| & \les \varep r^{-3}u^{-1/2},
      \end{aligned}
    \end{align}
    where here $r := t-u$.
  \item The induced metric and connection coefficients adapted to the maximal foliation $\Si_t$ and maximal-null foliation $\Si_t$ and $\CC_u$ satisfy decay statements consistent with~\eqref{est:intdecayStabCKintro} and~\eqref{est:extdecayStabCKintro}.
  \item The spacetime $(\MM,\g)$ admits a past/future timelike, past/future null and spacelike infinities $i^-,i^+$, $\scri^-,\scri^+$ and $i_0$ on which one can make sense of asymptotic quantities and their evolution equations.
  \end{itemize}
  \end{subequations}
\end{theorem}
\paragraph{Remarks on Theorem~\ref{thm:StabCKintroSTAB}}
\begin{enumerate}[ref=\thetheorem\alph*,label=\thetheorem\alph*]
\item\label{item:vectorfieldmethod} The proof of Theorem~\ref{thm:StabCKintroSTAB} in~\cite{Chr.Kla93} is based on the \emph{vectorfield method}, which proceeds in the following two steps:
  \begin{enumerate}
  \item[Step 1] Einstein equations~\eqref{eq:EVEintroSTAB} induce the following \emph{Bianchi equations}
    \begin{align}\label{eq:BianchiSTAB}
      \D^\al\R_{\al\be\ga\de} & = 0.
    \end{align}
    Multiplying and commuting~\eqref{eq:BianchiSTAB} with a set of approximate conformal Killing vectorfields, wave-type energy estimates are obtained. These estimates hold provided that the nonlinear error terms produced by the conformal Killing approximations are controlled. 
  \item[Step 2] From the boundedness of $L^2$ fluxes for (derivatives of) $\R$ through the hypersurfaces $\Si_t$ and $\CC_u$ resulting from the energy estimates of Step 1, one deduces decay estimates for $\R$ using Klainerman-Sobolev embeddings, as well as boundedness and decay estimates for the induced metric and connection coefficients associated to the maximal-null foliation. This is done using \emph{structure equations} which schematically read
    \begin{align*}
      \nab \Ga = \R + \nab\Ga + \Ga\cdot\Ga.
    \end{align*}
    Here $\nab$ are derivatives in the $\elb,\el$ and tangential directions, $\Ga$ are connection coefficients as defined in~\eqref{eq:defnulldecompSTAB}, the terms $\R,\nab\Ga$ on the right-hand side are treated as linear source terms and the terms $\Ga\cdot\Ga$ as nonlinear error terms.
  \end{enumerate}
  The crux of the proof of Theorem~\ref{thm:StabCKintroSTAB} in~\cite{Chr.Kla93} is the control of the nonlinear error terms of Step 1. Since the approximate conformal Killing vectorfields are constructed upon the geometric time and optical functions, these nonlinear error terms can be expressed in terms of the connection coefficients $\Ga$ for the foliations $\Si_t$ and $\CC_u$, and their control thus crucially relies on the decay estimates obtained in Step 2.
\item Assumption~\eqref{est:curvfluxStabCKintro2}, together with the maximal assumption on $\Si$ guarantee that the boundary fluxes on $\Si$ arising when performing energy estimates for Bianchi equations commuted with a set of approximate Killing vectorfields are controlled by $\varep$.\footnote{Estimates~\eqref{est:curvfluxStabCKintro2} do not bound the boundary flux for the uncommuted Bianchi equations, which would in fact cause the ADM mass of $\Si$ to vanish. The control of the ADM mass is obtained from~\eqref{est:AFStabCKintro} and~\eqref{est:curvfluxStabCKintro1}.}
\item\label{item:restframeStabCKintro} The global time function $t$ is constructed by imposing that its level sets $\Si_t$ are maximal hypersurfaces of $\MM$, that $\Si = \{t=0\}$ and that $n \to 1$ when $r\to\infty$. 
  These last conditions are equivalent to the choice of boundary for $\Si_t$ at infinity. It physically corresponds to considering a \emph{centre-of-mass frame} for the system (see the discussion in the introduction of~\cite{Chr.Kla93}).
\item The vectorfield method of~\cite{Chr.Kla93} is wrapped in an elaborate \emph{bootstrap argument}. One of the main challenge is to define the geometric constructions within this bootstrap argument to obtain sufficient decay rates for the associated metric and connection coefficients (see Item~\ref{item:vectorfieldmethod}). In~\cite{Chr.Kla93} the global optical function $u$ is constructed by an \emph{initialisation on the last slice}, \emph{i.e.} by imposing that the outgoing null hypersurfaces $\CC_u$ are backwards emanating from the $2$-spheres level sets of a \emph{canonical foliation} on a last slice $\Si^\ast$ corresponding to the future boundary of the bootstrap region. 
\item A first stability result for initial data with stronger decay assumptions was obtained in~\cite{Fri83}. A global stability result has been obtained in~\cite{Bie10,Bie.Zip09} using the same general techniques as in~\cite{Chr.Kla93} but under relaxed assumptions for both the regularity and decay of the initial data. A global stability result for Minkowski space has also been obtained using wave coordinates, see~\cite{Lin.Rod10}. See also~\cite{Hin.Vas20} for an alternative proof.
\end{enumerate}

\section{Main theorem}
This section is dedicated to the following precised version of Theorem~\ref{thm:mainv1rough} which is the main result of this paper. We also refer the reader to Theorem~\ref{thm:mainv2} for the detailed assumptions and conclusions. 
\begin{theorem}[Main theorem, more precise version]\label{thm:mainv1}
  Let $(\Sit_1,\CCt_0)$ be smooth spacelike-characteristic initial data, such that
  \begin{subequations}\label{est:assthmmainv1STAB}
  \begin{itemize}
  \item we have the following curvature fluxes bounds through $\Sit_1$
    \begin{align}
      \int_{\Sit_1}\le|\D^{\leq 2}\R\ri|^2 & \leq \varep^2,
    \end{align}
    together with consistent bounds for a Cartesian coordinates system $(x^i)$ on $\Sit_1$,
  \item the null hypersurface $\CCt_0$ is future geodesically complete, foliated by the $2$-spheres of a geodesic foliation $(S'_s)_{1\leq s < +\infty}$, and for the associated geodesic null pair, 
    we have the following curvature fluxes bounds
    \begin{align}
      & \int_{1}^{\infty}\int_{S'_s} \bigg(\le|\Ndt^{\leq 2}\beb\ri|^2 + \le|s\Ndt^{\leq 2}\rho\ri|^2 + \le|s\Ndt^{\leq 2}\sigma\ri|^2 + \le|s^2\Ndt^{\leq 2}\be\ri|^2 + \le|s^2\Ndt^{\leq 2}\al\ri|^2 \bigg) \d s \leq \varep^2
    \end{align}
    where $\Ndt\in\le\{(s\Nd), (s\Nd_4), \Nd_3\ri\}$, 
    together with consistent bounds for the metric and connection coefficients.
  \end{itemize}
  \end{subequations}
  There exists $\varep_0>0$ such that if $\varep < \varep_0$, the following holds for the future maximal globally hyperbolic development $(\MMt,\g)$ of $(\Sit_1,\CCt_0)$.
  \begin{itemize}
  \item The spacetime $(\MMt,\g)$ is future geodesically complete.
  \item The spacetime $(\MMt,\g)$ is covered by an interior and an exterior region $\MM^\intr$ and $\MM^\ext$, intersecting at a \emph{timelike transition hypersurface} $\TT = \MM^\intr\cap\MM^\ext$.
  \item There exists a global time function $t$ on $\MM^\intr$ ranging up to $+\infty$ foliating $\MM^\intr$ by spacelike maximal hypersurface $\Si_t$.
  \item There exists a global optical function $u$ on $\MM^\ext$ ranging up to $+\infty$ foliating $\MM^\ext$ by outgoing null hypersurfaces $\CC_u$. There exists a global function $\ub$ on $\MM^\ext$ which is a geodesic affine parameter on $\CC_u$, foliating $\CC_u$ by $2$-spheres $S_{u,\ub}$. Moreover, on the transition hypersurface $\TT$, we have
    \begin{align*}
      u & = \cc \ub, & t & = \half (u+\ub),
    \end{align*}
    where $0<\cc<1$ is a fixed parameter.
  \item We have the following decay bounds in $\MM^\intr$
    \begin{align*}
      |\R| & \les \varep t^{-7/2},
    \end{align*}
    together with consistent bounds for the metric and connection coefficients.
  \item We have the following decay bounds in $\MM^\ext$
    \begin{align*}
      \begin{aligned}
        |\al| & \les \varep \ub^{-7/2}, & |\be| & \les \varep \ub^{-7/2}, & |\rho| & \les \varep \ub^{-3}u^{-1/2},\\
        |\alb| & \les \varep \ub^{-1}u^{-5/2}, & |\beb| & \les \varep \ub^{-2}u^{-3/2}, & |\sigma| & \les \varep \ub^{-3}u^{-1/2},
      \end{aligned}
    \end{align*}
    together with consistent bounds for the metric and connection coefficients.
  \item The spacetime $(\MMt,\g)$ admits a future timelike and future null infinity $i^+$ and $\scri^+$. The future null infinity $\scri^+$ is future geodesically complete, admits well-defined notions of Bondi mass and angular momentum for which we obtain Bondi mass loss formula and angular momentum evolution equation along $\scri^+$, and which tend to $0$ at future timelike infinity $i^+$.
  \end{itemize}
\end{theorem}
\begin{figure}[h!]
  \centering
  \includegraphics[width=0.8\linewidth]{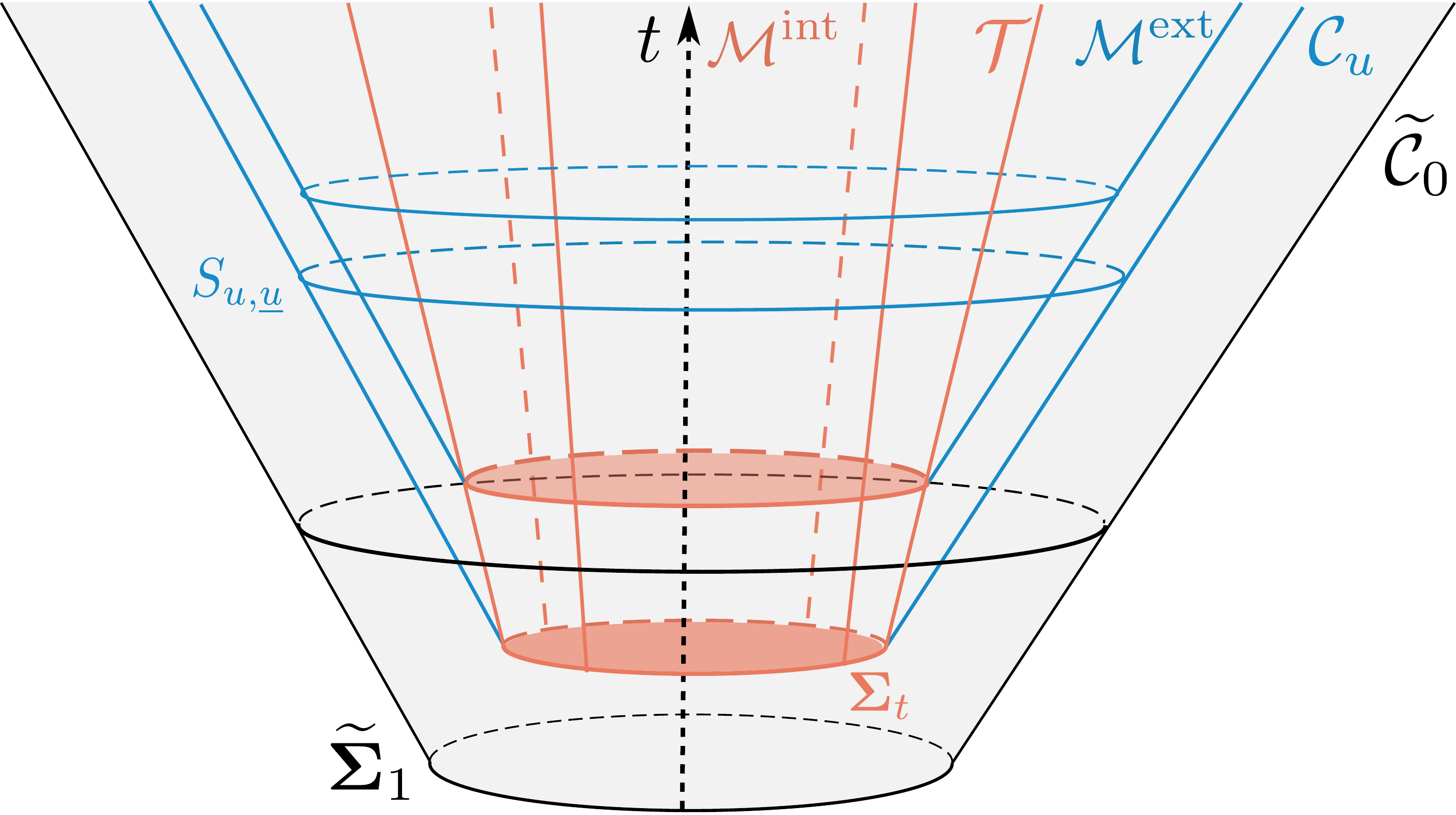}
  \caption{The global nonlinear stability of Minkowski space for characteristic data.}
  \label{fig:stabMinkSTABthmmainv1}
\end{figure}
\paragraph{Remarks on Theorem~\ref{thm:mainv1}}
\begin{enumerate}[ref=\thetheorem\alph*,label=\thetheorem\alph*]
\item The initial data assumptions~\eqref{est:assthmmainv1STAB} match what can be obtained for an outgoing null hypersurface in~\cite{Chr.Kla93} or~\cite{Kla.Nic03}. Therefore, Theorem~\ref{thm:mainv1} provides a stability result for the complementary region to the exterior region considered in~\cite{Kla.Nic03}. Together, they amount to a stability result for initial data posed on a spacelike hypersurface.

\item Theorem~\ref{thm:mainv1} was conjectured to hold true in~\cite{Kla.Nic03,Tay17} and its conclusions were used in~\cite{Tay17}.
\item The basic scheme of proof of Theorem~\ref{thm:mainv1} is a vectorfield method wrapped in a bootstrap argument as in~\cite{Chr.Kla93}. The main novelty is the introduction and control of new geometric constructions
  , which provide suitable spacetime decompositions to run these arguments.
\item Our constructions display the following new crucial geometric features.
  \begin{itemize}
  \item They virtually emanate from the future infinity of a (timelike) central axis. This guarantees optimal decay rates. It replaces an asymptotically flat spacelike infinity which plays a similar crucial role in~\cite{Chr.Kla93,Kla.Nic03}.
  \item In the interior region, our constructions are build on \emph{spacelike maximal hypersurfaces with prescribed boundaries} and \emph{global harmonic coordinates}. This makes any reference to null decompositions and spherical foliations -- which degenerate at the central axis -- disappear in that region. 
  \end{itemize}
\item In the proof of Theorem~\ref{thm:mainv1}, we match discontinuous gauge choices across the timelike interface $\TT$ without using the gluing procedure of~\cite{Chr.Kla93,Kla.Sze17}. Our matching features a \emph{mean value argument} which compensates regularity losses at the timelike interface. We believe that this new treatment gains in concision and clarity.
\item In the appendix to this paper, we also provide new optimal estimates and control for harmonic coordinates on a $3$-dimensional Riemannian manifold, only based on elementary energy and Bochner estimates (see Theorem~\ref{thm:globharmonics} and Appendix~\ref{sec:globharmo}). We moreover give a full statement and proof for general limits of the metric and connection coefficients and their derivatives in all directions at the vertex of a (general family of) null cones (see Theorem~\ref{thm:vertex} and Appendix~\ref{app:vertexlimits}). 
\end{enumerate}
  

In the next section, we give an overview of the proof of Theorem~\ref{thm:mainv1}. We postpone discussions/comparisons with other results to Section~\ref{sec:discSTAB}.

\section{Overview of the proof of Theorem~\ref{thm:mainv1}}\label{sec:proofStabMinkintro}

The proof of the global nonlinear stability result of Theorem~\ref{thm:mainv1} goes by a standard continuity argument on the maximal parameter $\uba$ such that the smooth maximal globally hyperbolic development $(\MMt,\g)$ admits a subregion $\MM_\uba$ (of size $\uba$) which we geometrically describe next in Section~\ref{sec:geomsetupstabintro}. In this overview, we will focus on the geometric setup of $\MM_\uba$ and the obtention of bounds for the curvature and the geometric structures of $\MM_\uba$. We refer the reader to Section~\ref{sec:mainresult} for the full setup and conclusion of the bootstrap argument and for its consequences in the limit $\uba \to + \infty$, from which the conclusions of Theorem~\ref{thm:mainv1} follow. 

\subsection{Geometric setup of the bootstrap region $\MM_{\protect\uba}$}\label{sec:geomsetupstabintro}

Let $O = \{x^i = 0\}$ be the centre of the initial spacelike hypersurface $\Sit_1$ for the Cartesian coordinates $x^i$ given as assumption on $\Sit_1$. We define the \emph{central axis} $\o \subset \MMt$ to be the timelike geodesic parametrised by arc-length such that $\o(1) = O$ and $\doto(1)$ is future-pointing and normal to $\Sit_1$ at $O$.\\

For a fixed parameter $\uba\geq 1$, we define the \emph{last cone} $\CCba$ to be the ingoing null cone backwards emanating from the point $\o(\uba)$. The cone $\CCba$ is foliated by the $2$-spheres $S_{u,\uba}$ of a \emph{canonical foliation} with parameter $u$ which ranges from $u|_{\o(\uba)} = \uba$ to $0$.\\

\begin{figure}[h!]
  \centering
  \includegraphics[height=9cm]{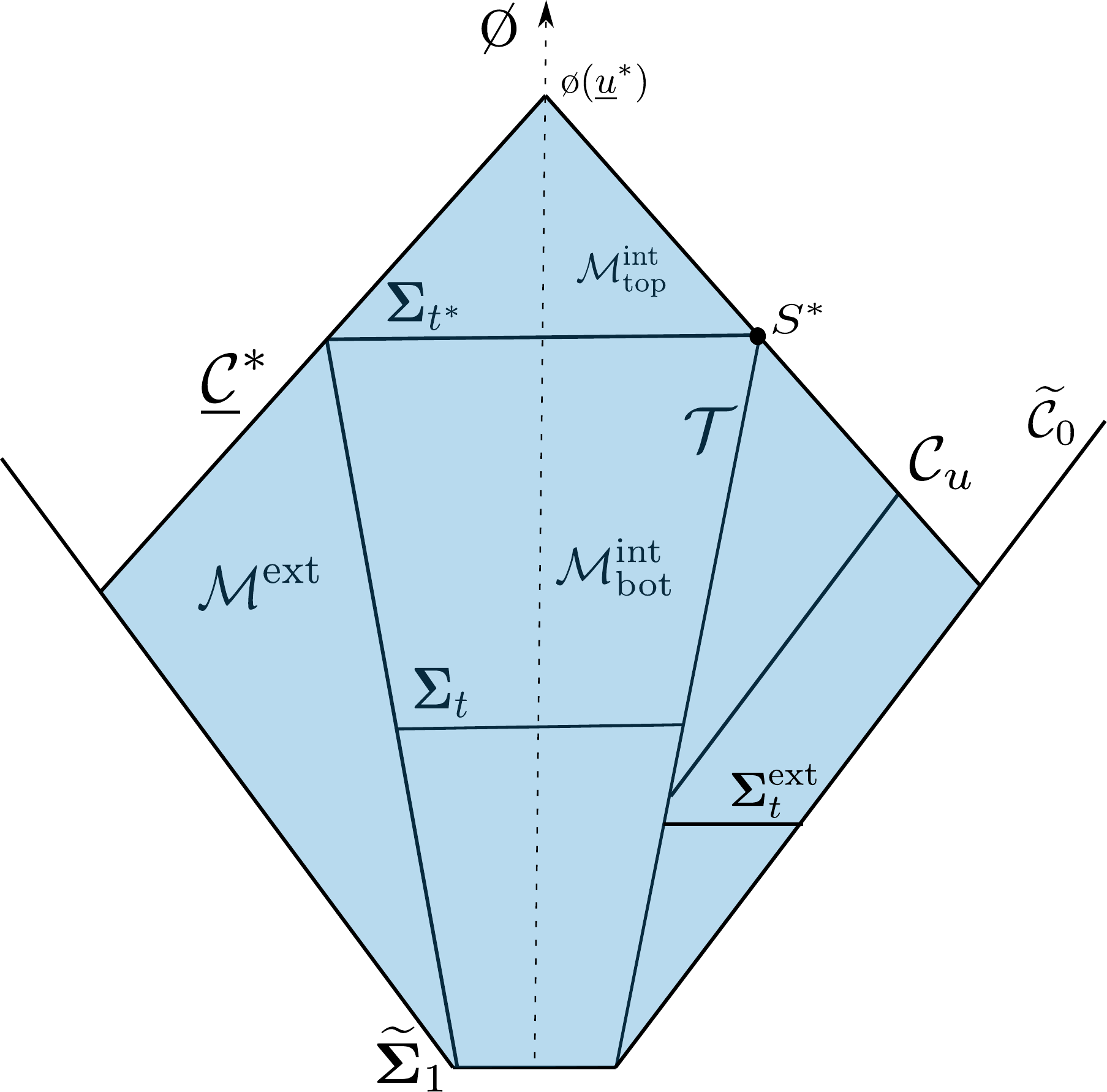}
  \caption{The bootstrap region $\MM_\uba$}
  \label{fig:MMubaintro}
\end{figure}

We define $\CC_u$ to be the outgoing null hypersurfaces backwards emanating from the $2$-spheres $S_{u,\uba}$, and we denote by $u$ the associated optical function. We foliate the hypersurfaces $\CC_u$ by the $2$-spheres $S_{u,\ub}$ of the \emph{geodesic foliation} with parameter $\ub$ ranging from $\ub|_{\CCba} = \uba$ to $\cc^{-1}u$, where $0 < \cc < 1$ is a (suitably chosen) parameter. We define the following \emph{exterior region} (see Figure~\ref{fig:MMubaintro} for a graphic representation)
\begin{align*}
  \MM^\ext & := \bigcup_{\substack{0~\leq~u~\leq~\cc\uba \\ \cc^{-1}u~\leq~\ub~\leq~\uba}} S_{u,\ub},
\intertext{the following (timelike) \emph{transition hypersurface}}\\
  \TT & := \bigcup_{\substack{0~\leq~u~\leq~\cc\uba \\ \ub = \cc^{-1}u}} S_{u,\ub},\\
\intertext{and the following \emph{last sphere}} \\
  S^\ast & := S_{\cc\uba,\uba} = \CCba\cap\TT.
\end{align*}

Let $\Si_t$ be the \emph{maximal hypersurfaces} coinciding on the transition hypersurface with the $2$-spheres $S_{u,\ub}$, \emph{i.e.} such that the associated maximal time function $t$ satisfies
\begin{align*}
  t|_{\TT} & := \half(u+\ub).
\end{align*}
We define the following \emph{bottom interior region}
\begin{align*}
  \MM^\intr_\bott & := \bigcup_{1\leq t \leq \tast} \Si_t
\end{align*}
where $\tast := \half(1+\cc)\uba$, and we call $\Si_\tast$ the \emph{last maximal hypersurface} (note that the boundary of the last maximal hypersurface is the last sphere, \emph{i.e.} $\pr\Si_\tast = S^\ast$). We further define the \emph{top interior region} $\MM^\intr_{\mathrm{top}}$ to be the domain of dependence of $\Si_\tast$. We are now able to define the full spacetime bootstrap region $\MM_\uba$
\begin{align*}
  \MM_\uba & := \MM^\ext \cup \MM^\intr,\\
  \intertext{where}\\
  \MM^\intr & := \MM^\intr_\bott \cup \MM^\intr_{\mathrm{top}}.
\end{align*}

As it will play a key role in the control of the interior region, we also define \emph{global harmonic Cartesian coordinates} $(x^i)$ on the last slice $\Si_\tast$ to be functions such that on $\Si_\tast$
\begin{align*}
  \Delta x^i & = 0,
\end{align*}
and such that Dirichlet boundary conditions for $x^i$ at $\pr\Si_\tast= S^\ast$ are fixed \emph{via} (a suitable class of) conformal isomorphism of $S^\ast$ to the Euclidean $2$-sphere. 


\paragraph{Remarks}
\begin{itemize}
\item The central axis and its last point in the future direction is the starting point of our \emph{initialisation from timelike infinity} procedure. 
\item The last cone $\CCba$ and the last maximal hypersurface $\Si_\tast$ play a similar role as the last maximal hypersurface in~\cite{Chr.Kla93} or the last cone in~\cite{Kla.Nic03}.
\item\label{item:canonicalregfeatures} The canonical foliation on $\CCba$ must provide sufficient regularity for its transverse geometric quantities (\emph{i.e.} $\trchi$ and $\chih$) since this determines the regularity of the transversely emanating foliation of null cones $\CC_u$. The canonical foliations on the last slices in~\cite{Chr.Kla93,Kla.Nic03} are build to satisfy the same transverse regularity features. See also the discussions in~\cite{Czi.Gra19,Czi.Gra19a}.
\item Considering the geodesic foliation on the outgoing cones $\CC_u$ is geometrically simpler than constructing a double null foliation which requires a second transverse hypersurface to be initialised. The counterpart is additional difficulty in the analysis of the null structure equations to avoid regularity loss. See a similar foliation choice and difficulties in~\cite{Kla.Sze17}.
\item\label{item:CCudegvertex} Since they are not required to emanate from one point, but are defined from $2$-spheres on $\CCba$, the outgoing null hypersurfaces $\CC_u$ may degenerate before reaching their potential vertex $\ub \to u$. We thus restrict their definition to the exterior region $\cc^{-1}u \leq \ub$, with $0<\cc<1$, on which they remain regular.
\item\label{item:noinitlayerintro} In this overview, we shall assume that the initial hypersurfaces $\Sit_1$ and $\CCt_0$ match the geometric constructions described above, \emph{i.e.} $\Sit_1 = \Si_1$ and $\CCt_0 = \CC_0$. In general this does not hold. See Section~\ref{sec:initlayer} for the \emph{initial layer} existence and comparisons arguments. In particular, in Sections~\ref{sec:definition}--\ref{sec:initlayer}, we will only assume that $u$ ranges up to $1$ (and not to $0$) and $t$ ranges up to $(1+\cc^{-1})/2$ (and not to $1$).
\item That the harmonic functions $x^i$ on the last slice $\Si_\tast$ form a global coordinate system on $\Si_\tast$ is a result, obtained as a consequence of estimates for the functions $x^i$. These estimates are obtained using energy and Bochner estimates for the above defined Dirichlet problem on $\Si_\tast$. This (only) involves intrinsic quantities (the Ricci curvature and fundamental forms of the boundary) and basic functional estimates on $\Si_\tast$. See Theorem~\ref{thm:globharmonics}.  
\end{itemize}

\subsection{Global energy estimates}\label{sec:globenerintro}
The final goal is to obtain bounds on the spacetime metric evaluated on the above defined geometric structures of $\MM_\uba$ and decay estimates to the corresponding Minkowskian quantities, in terms of the initial data on $\Si_1\cup\CC_0$. To that end, the spacetime curvature $\R$ is the key dynamical object since it satisfies Bianchi equations~\eqref{eq:BianchiSTAB} for which wave-type energy estimates can be obtained. We use the \emph{vectorfield method} (see Item~\ref{item:vectorfieldmethod}) for these equations in the region $\MM_\uba$.\\

We first recall tools developed in~\cite{Chr.Kla93} to perform energy estimates for Bianchi equations.\\

First, commuting Bianchi equation~\eqref{eq:BianchiSTAB} with a vectorfield $\mathbf{X}$ gives schematically
\begin{align}\label{eq:BianchicommutedSTAB}
  \D^\al\le(\Lieh_{\mathbf{X}}\R\ri)_{\al\be\ga\de} & = \D^{\leq 1}\le(^{(\mathbf{X})}\pih\ri) \cdot \D^{\leq 1}\R, 
\end{align}
where $\Lieh$ is a normalised derivative in the $\mathbf{X}$ direction, and where $^{(\mathbf{X})}\pih$ is the following traceless part of the \emph{deformation tensor} of $\mathbf{X}$
\begin{align*}
  ^{(\mathbf{X})}\pih_{\mu\nu} & := \D_\mu\mathbf{X}_\nu + \D_\nu\mathbf{X}_\mu - \frac{1}{2}\le(\D^\al\mathbf{X}_\al\ri)\g_{\mu\nu}. 
\end{align*}
\begin{remark}
  For a \emph{conformal Killing vectorfield} $\mathbf{X}$, we have $^{(\mathbf{X})}\pih = 0$. Thus from~\eqref{eq:BianchicommutedSTAB} we obtain that commuting Bianchi equation by $\mathbf{X}$ gives
  \begin{align*}
    \D^\al\le(\Lieh_{\mathbf{X}}\R\ri)_{\al\be\ga\de} & = 0.
  \end{align*}
  For \emph{approximate} conformal Killing vectorfields, the source terms of the commuted Bianchi equations are \emph{nonlinear error terms}.
\end{remark}

Second, we introduce the \emph{Bel-Robinson tensor} $Q(\mathbf{W})$ associated to a $4$-tensor $W$ by
\begin{align*}
  Q(\mathbf{W})_{\a \be \ga \de}:= \mathbf{W}_{\a \nu \ga \mu} \mathbf{W}_{\be \,\,\, \de}^{\,\,\,\nu \,\,\, \mu} + {}^\ast \mathbf{W}_{\a \nu \ga \mu} {}^\ast \mathbf{W}_{\be \,\,\, \de}^{\,\,\,\nu \,\,\, \mu},
\end{align*}
where $\dual$ denotes the Hodge dual of $\mathbf{W}$. When $\mathbf{W}=\R$, we have the following consequence of the Bianchi equations~\eqref{eq:BianchiSTAB}
\begin{align*}
  \D^\al Q(\R)_{\al\be\ga\de} & = 0.
\end{align*}
Multiplying and commuting the Bel-Robinson tensor by vectorfields $\mathbf{X}_1,\mathbf{X}_2,\mathbf{Y}_1,\mathbf{Y}_2,\mathbf{Y}_3$, we similarly obtain from~(\ref{eq:BianchicommutedSTAB}) the following schematic formula
\begin{align}\label{eq:contractedcommutedBelSTAB}
  \DIV \le(Q\le(\Lieh_{\mathbf{X}_1}\Lieh_{\mathbf{X}_2}\R\ri)(\mathbf{Y}_1,\mathbf{Y}_2,\mathbf{Y}_3)\ri) & = \D^{\leq 2}\pih \cdot \D^{\leq 2}\R,  
\end{align}
where $\DIV$ is the spacetime divergence for $1$-tensor and where $\pih$ denotes the traceless deformation tensors of $\mathbf{X}_1,\mathbf{X}_2,\mathbf{Y}_1,\mathbf{Y}_2,\mathbf{Y}_3$.
\begin{remark}
  For exact conformal Killing vectorfields, formula~\eqref{eq:contractedcommutedBelSTAB} produces a spacetime divergence-free vectorfield
  \begin{align*}
    Q\le(\Lieh_{\mathbf{X}_1}\Lieh_{\mathbf{X}_2}\R\ri)(\mathbf{Y}_1,\mathbf{Y}_2,\mathbf{Y}_3).
  \end{align*}
  An application of Stokes theorem yields an exact identity of boundary fluxes for the above vectorfield. Provided that the multiplier vectorfields $\mathbf{Y}_1,\mathbf{Y}_2,\mathbf{Y}_3$ are suitably chosen, this yields an \emph{energy estimate} in $\MM$.
\end{remark}
\begin{remark}
  For approximate conformal Killing vectorfields, the same procedure using formula~\eqref{eq:contractedcommutedBelSTAB} produces an energy estimate with a spacetime integral of nonlinear error terms.
\end{remark}

Using the above tools, we can now perform \emph{global energy estimates} in $\MM^\ext\cup\MM^\intr_\bott$. These estimates are obtained using the following set of contracted and commuted Bel-Robinson tensors
\begin{align}\label{eq:defBelRobintro}
  \begin{aligned}
  Q\le(\Lieh_\TX\R\ri)(\KX,\KX,\KX), \quad Q\le(\Lieh_{\OOO}\R\ri)(\KX,\KX,\TX), \\
  Q\le(\Lieh_\OOO\Lieh_\OOO\R\ri)(\KX,\KX,\TX), \quad Q\le(\Lieh_\SX\Lieh_\TX\R\ri)(\KX,\KX,\KX), \quad Q\le(\Lieh_\OOO\Lieh_\TX\R\ri)(\KX,\KX,\TX),
  \end{aligned}
\end{align}
where 
\begin{itemize}
\item $\TX$ is an approximation for the \emph{time translation Killing vectorfield} $\pr_t$ of Minkowski space, 
\item $\SX$ is an approximation for the \emph{scaling conformal Killing vectorfield} $t\pr_t + r\pr_r$ of Minkowski space, 
\item $\KX$ is an approximation for the \emph{Morawetz conformal Killing vectorfield} $(t^2+r^2)\pr_t + 2tr\pr_r$ of Minkowski space,
\item $\OOO$ are approximations for the three \emph{rotation Killing vectorfields} $x^1\pr_2-x^2\pr_1$, $x^2\pr_3-x^3\pr_2$ and $x^3\pr_1-x^1\pr_3$ of Minkowski space.
\end{itemize}
We assume for the moment that these vectorfields are given and postpone their respective definitions in the bottom interior and exterior region to Section~\ref{sec:defKillintro}.\\

Applying Stokes theorem and formula~(\ref{eq:contractedcommutedBelSTAB}) to the set of vectorfields~\eqref{eq:defBelRobintro} simultaneously in the bottom interior region $\MM^\intr_\bott$ and the exterior region $\MM^\ext$, we obtain the following energy estimates (see Figure~\ref{fig:MMubaintro} for a graphic representation of these hypersurfaces)
\begin{align*}
  \int_{\Si_t} + \int_{\CC_u} + \int_{\Si_t^\ext} + \int_{\CCba\cap\MM^\ext} & \les \int_{\Si_1} + \int_{\CC_0} + ~\EE^\TT + \EE^\intr + \EE^\ext,
\end{align*}
where
\begin{itemize}
\item the integrands are the contracted and commuted Bel-Robinson tensors~\eqref{eq:defBelRobintro},
\item for all $1\leq t \leq \tast$ and all $0 \leq u \leq \cc\uba$ (we recall that)
  \begin{itemize}
  \item $\Si_t$ are the maximal hypersurfaces of $\MM^\intr$,
  \item $\CC_u$ are the outgoing null hypersurface of $\MM^\ext$,
  \item $\Si_t^\ext$ are the level sets of the time function $\half(u+\ub)$ in $\MM^\ext$,
  \item $\CCba\cap\MM^\ext$ is the exterior part of the last cone $\CCba$,
  \item $\Si_1$ and $\CC_0$ are the initial hypersurfaces for which we have bounds for the curvature fluxes (see the assumptions of Theorem~\ref{thm:mainv1}),
  \end{itemize}
\item the nonlinear error term $\EE^\TT$ is the difference of boundary fluxes on the timelike transition hypersurface $\TT$ for Stokes formula applied in $\MM^\intr_\bott$ and $\MM^\ext$,
\item the nonlinear error terms $\EE^\intr$ and $\EE^\ext$ are spacetime integrals over $\MM^\intr_\bott$ and $\MM^\ext$ respectively, involving (two derivatives of) the spacetime curvature tensor $\R$ and (two derivatives of) the deformations tensors $\pih$ for the approximate conformal Killing vectorfields.  
\end{itemize}

The control of the error terms $\EE^\intr$ and $\EE^\ext$ is the crux of the analysis and is obtained provided that sufficient \emph{decay} can be obtained for (null decompositions of) $\R$ and (null decompositions of) $\pih$, as well as sufficient \emph{regularity} can be obtained for $\pih$.\\

The control of the interface error term $\EE^\TT$ is obtained provided that the difference of the corresponding interior/exterior approximate conformal Killing vectorfields at the interface $\TT$ can be controlled with sufficient decay and regularity. It also requires that the spacetime curvature tensor $\R$ has optimal regularity on $\TT$, \emph{i.e.} that $\D^2\R \in L^2(\TT)$. Since the hypersurface $\TT$ is timelike, this cannot be obtained from bounds on energy estimates boundary fluxes on $\TT$. We thus rely on a \emph{mean value argument} which selects a suitable transition parameter $\cc$/transition hypersurface $\TT$ on which such a control holds.\\

From these controls, we obtain
\begin{align*}
  \int_{\Si_t} + \int_{\CC_u} + \int_{\Si_t^\ext} + \int_{\CCba\cap\MM^\ext} & \les \varep^2 + (D\varep)^3 \les \varep^2.
\end{align*}

\paragraph{Remarks}
\begin{itemize}
\item The contracted and commuted Bel-Robinson tensors~(\ref{eq:defBelRobintro}) are identical to the ones used in~\cite{Chr.Kla93,Kla.Nic03}.
\item The decay rates obtained in this paper are similar, but slightly different to the decay rates of~\cite{Chr.Kla93,Kla.Nic03} due to the difference of geometric constructions.
\item In Section~\ref{sec:globener}, we use systematically that the tensors $\R$ and $\pih$ are respectively controlled with the following regularity in $\MM^\ext$
  \begin{align*}
    \R & \in L^\infty(\MM^\ext), & \D\R & \in L^\infty_{u,\ub} L^4(S_{u,\ub}), & \D^2\R & \in L^2(\MM^\ext),\\
    \pih & \in L^\infty(\MM^\ext), & \D\pih & \in L^\infty_{u,\ub} L^4(S_{u,\ub}), & \D^2\pih & \in L^2(\MM^\ext),
  \end{align*}
 and with the following regularity in $\MM^\intr_\bott$
 \begin{align*}
   \R & \in L^\infty(\MM^\intr_\bott), & \D\R & \in L^\infty_{t} L^6(\Si_t), & \D^2\R & \in L^2(\MM^\intr_\bott),\\
   \pih & \in L^\infty(\MM^\intr_\bott), & \D\pih & \in L^\infty_{t} L^6(\Si_t), & \D^2\pih & \in L^2(\MM^\intr_\bott),
 \end{align*}
 which is in the spirit of the (bulk) Morawetz estimates of the $r^p$-method (see for example~\cite{Hol10}) and simplifies the analysis of~\cite{Chr.Kla93,Kla.Nic03}.
\item The main new feature of the global energy estimates in the region $\MM^\ext\cup\MM^\intr$ is to avoid the gluing procedure of~\cite{Chr.Kla93,Kla.Sze17} for functions/vectorfields across the interface $\TT$. This is replaced by a treatment of transition nonlinear error terms $\EE^\TT$ which arise from the discontinuity of gauge choices. We believe that this new procedure gains in clarity.
\end{itemize}

\begin{figure}[!h]
  \centering
  \includegraphics[width=15cm]{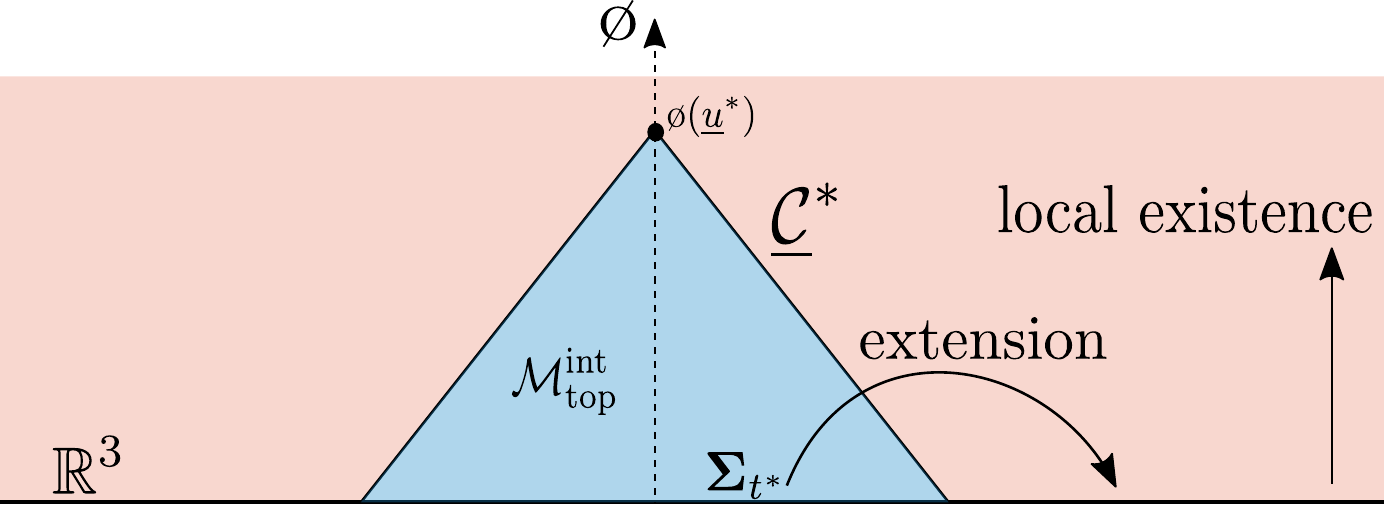}
  \caption{Energy estimates in (the rescaled) $\MM^\intr_\topp$}
  \label{fig:energestMMintrtopp}
\end{figure}
To obtain a control for the geometry of the cone $\CCba$, we also need a control of the boundary fluxes through $\CCba\cap \MM^\intr_\topp$. To this end, we perform a \emph{$\tast$-rescaling} of the spacetime region $\MM^\intr_\topp$. The metric, connection and curvature components scale homogeneously in that region, which thus reduces to a size-$1$ region. Then, we \emph{extend} the data on the last slice $\Si_\tast$ to $\RRR^3$ using the result of~\cite{Czi18} and apply a \emph{local (time-$1$) existence result for the extended (small) data}. Then, we can perform energy estimates in $\MM^\intr_{\mathrm{top}}$ -- which we recall is the future domain of dependence of $\Si_\tast$ --, using \emph{Cartesian} approximate Killing vectorfields (\emph{i.e.} the approximate Minkowskian time and space translations $\pr_\mu$). The error terms are controlled using the local existence result. From a comparison argument on $\CCba\cap\MM^\intr$, we obtain control for the (null decomposition of the) curvature on $\CCba\cap\MM^\intr$, that is, an estimate of the type
\begin{align*}
  \int_{\CCba\cap\MM^\intr} & \les \int_{\Si_\tast} + \; \text{nonlinear error terms} \; \les \varep^2 + (D\varep)^3 \les \varep^2,
\end{align*}
where we used that the energy boundary fluxes through the last slice $\Si_\tast$ were controlled by the previous energy estimates in $\MM^\intr_\bott\cap\MM^\ext$.

\paragraph{Remarks}
\begin{itemize}
\item Spacetime regions as $\MM^\intr_\topp$ with conical degeneracies and outside of spherical symmetry are not treated as such in other works. Our procedure is in the spirit of~\cite{Czi.Gra19,Czi.Gra19a}.  
\item Energy estimates using the set of (spherical) vectorfields $\TX,\SX,\KX,\OOO$ do not provide optimal bounds and is the cause of degeneracies -- even for the wave equation in Minkowski space -- in the region $\MM^\intr_{\mathrm{top}}$.\footnote{At the vertex, $|\OOO| \to 0$ and $\TX,\SX,\KX$ become colinear.}
\item Classical Sobolev estimates degenerate at the vertex of the cone, which without the extension to a larger spacetime (artificially) complicates the control of the nonlinear terms at the vertex.
\item The extension result for $\Si_\tast$ uses an \emph{optimal control} for the fundamental forms of $\Si_\tast$ which can be obtained using global harmonic Cartesian coordinates and the new optimal control for these coordinates established in this paper (see Theorem~\ref{thm:globharmonics} and Appendix~\ref{sec:globharmo}).
\end{itemize}


\subsection{Curvature, connection and metric control}\label{sec:curvconnmetricintro}
From the boundary fluxes controls obtained in Section~\ref{sec:globenerintro}, control for the spacetime curvature follows similarly as in~\cite{Chr.Kla93,Kla.Nic03}.\\

From the control of the spacetime curvature and using the structure equations, one deduces a control for the connection and metric by
\begin{itemize}
\item integration of transport equations from the vertex $\o(\uba)$ of $\CCba$ to $\CCba$ and elliptic estimates on the $2$-spheres of the canonical foliation on $\CCba$,
\item integration of transport equations from $\CCba$ to $\MM^\ext$ and $\TT$, and elliptic estimates on the $2$-spheres of the geodesic foliation in $\MM^\ext$,
\item elliptic estimates on the maximal slices $\Si_t$ using mixed (implicit) Dirichlet-Neumann boundary values on $\pr\Si_t \subset \TT$.
\end{itemize}

\paragraph{Remarks}
\begin{itemize}
\item We use that the spacetime is assumed to be \emph{smooth} and in consequence, we establish general limits for (all derivatives of) the  metric and connection coefficients at the vertex of a cone.\footnote{In fact, we even establish limits in the more general setting of a foliation of cones emanating from a central axis.} These limits are used as initial data to integrate the transport equations on $\CCba$.
\item The main challenge in the succession of geometric constructions of $\MM_\uba$ which are build upon one another, is to avoid the addition of \emph{loss of regularity} for the connection coefficients. The regularity must eventually be sufficient for the control of the nonlinear error terms in the energy estimates of Section~\ref{sec:globenerintro}.
\end{itemize}

\subsection{Approximate conformal Killing vectorfields}\label{sec:defKillintro}
In the exterior region $\MM^\ext$, the functions $u,\ub$ are associated to a null pair $(\elb,\el)$, such that
\begin{align*}
  \el & = - (\D u)^\sharp, & \g(\elb,\elb) & = 0, & \g(\el,\elb) & = -2,
\end{align*}
and such that $\elb,\el$ are orthogonal to the $2$-spheres $S_{u,\ub}$ of the null-geodesic foliation. Here $(\D u)^\sharp$ is the spacetime gradient of $u$. These vectorfields serve as approximations for the Minkowskian $\pr_t-\pr_r$ and $\pr_t+\pr_r$. Upon these we thus define the following approximate exterior conformal Killing vectorfields
\begin{align*}
  \TE & := \half(\elb+\el), & \SE & :=\half(u\elb+\ub\el), & \KE & := \half(u^2\elb+\ub^2\el),
\end{align*}
and we postpone the definition of the exterior rotation vectorfields $\OOE$ to the end of this section.\\

In the bottom interior region $\MM^\intr_\bott$, we define $\TI$ to be the future-pointing unit normal to the maximal hypersurfaces $\Si_t$. To obtain further definitions for $\SI$ and $\KI$ we need an approximation for the Minkowskian vectorfield $r\pr_r$. This can be obtained by defining the vectorfield $\XI$ on the last slice $\Si_\tast$ using the global harmonic Cartesian coordinates 
\begin{align*}
  \XI & := \sum_{i=1}^3 x^i \nab x^i,
\end{align*}
and by extending it on $\MM^\intr_\bott$ by parallel transport along the flow of $t$. 
Using $\XI$, we have the following definitions for $\SI$ and $\KI$
\begin{align*}
  \SI & := t \TI + \XI, & \KI & := \le(t^2+\g(\XI,\XI)\ri)\TI + 2t \XI.
\end{align*}
Similarly, we can define the interior rotation vectorfields on the last slice $\Si_\tast$ by
\begin{align*}
  ^{(1)}\OOI & := x^2\nab x^3 - x^3\nab x^2, & ^{(2)}\OOI & := x^3\nab x^1 - x^1\nab x^3, & ^{(3)}\OOI & := x^1\nab x^2 - x^2\nab x^1,
\end{align*}
and extend this definition on $\MM^\intr_\bott$ by parallel transport along the flow of $t$.\\

To match with the definitions of the interior rotations, the definition of the exterior rotation vectorfields is initialised at $S^\ast = \pr\Si_\tast$ using the harmonic coordinates of $\Si_\tast$ by
\begin{align*}
  ^{(1)}\OOE & := x^2\Nd x^3 - x^3\Nd x^2, & ^{(2)}\OOE & := x^3\Nd x^1 - x^1\Nd x^3, & ^{(3)}\OOE & := x^1\Nd x^2 - x^2\Nd x^1.
\end{align*}
We then define them by Lie transport along $\CCba\cap\MM^\ext$
\begin{align*}
  \le[^{(\ell)}\OOE,\elb\ri] & = 0,
\end{align*}
and in $\MM^\ext$
\begin{align*}
  \le[^{(\ell)}\OOE,\el\ri] & = 0,
\end{align*}
for all $\ell=1,2,3$.\\

The control of the deformation tensors of
\begin{itemize}
\item $\TE,\SE,\KE$ is obtained from the control of the null connection coefficients associated to the null pair $(\elb,\el)$ in $\MM^\ext$,
\item $\OOE$ is obtained by estimates for the harmonic coordinates on $S^\ast$ and integration of transport equations from $S^\ast$ to $\CCba\cap\MM^\ext$ and to $\MM^\ext$,
\item $\TI$ is obtained from the control of the maximal connection coefficients in $\MM^\intr_\bott$,
\item $\XI,\SI,\KI,\OOI$ is obtained from estimates on $\Si_\tast$, using the bounds on the harmonic coordinates $(x^i)$, and by integration in $t$.
\end{itemize}

The control at the interface $\TT$ of the difference of vectorfields
\begin{itemize}
\item $\TE-\TI$ is obtained by a control of the \emph{slope} between the maximal hypersurfaces and the boundary $\TT$ (see a similar result in~\cite{Czi.Gra19a}),
\item $\XE-\XI$\footnote{where $\XE := \quar (\ub-u) (\el-\elb)$.} (and subsequently $\SE-\SI$ and $\KE-\KI$) is obtained by a control on $S^\ast$ using the harmonic coordinates and the control of the slope, and by integration along $\TT$, using that $\D^2\mathbf{X} \simeq 0$,
\item $\OOE-\OOI$ is obtained by a control on $S^\ast$ using the harmonic coordinates control, and by integration along $\TT$, using that $\D^2\OOO \simeq 0$.
\end{itemize}

\paragraph{Remarks}
\begin{itemize}
\item The rotation vectorfields in the interior region are only defined to extend the rotation vectorfields of the exterior region. They are not used to estimate the curvature in $\MM^\intr_\bott$.
\item The rotation vectorfields in the exterior region are used to estimate the tangential derivatives of the curvature. Thus, they have to be tangent to the $2$-spheres of the canonical and geodesic foliation on $\CCba$ and $\MM^\ext$ respectively. This motivates their definition by Lie transport. 
\end{itemize}

\section{Comparison to previous works}\label{sec:discSTAB}
In this section, we discuss the strategy and techniques of proof of Theorem~\ref{thm:mainv1} and compare them to other works.
\begin{enumerate}
\item\label{item:semiglobal} In the literature, the characteristic Cauchy problem outside of spherical symmetry is rather studied under a local (\emph{e.g.} \cite{Ren90,Cho.Chr.Mar11,Czi.Gra19a}) or semi-global perspective -- \emph{i.e.} in a size $1$ region from the initial null hypersurface -- as in~\cite{Luk12,Li.Zhu18} or~\cite{Cac.Nic05,Cac.Nic10}.
\item\label{item:modulationtimelikeinfinity} To obtain the \emph{full} global result of Theorem~\ref{thm:mainv1}, the spacetime geometric constructions are \emph{initialised at timelike infinity}, which corresponds in this paper to starting our construction at the last point $\o(\uba)$ on the central axis $\o$. This differs from:
  \begin{itemize}
  \item the semi-global existence results of Item~\ref{item:semiglobal} where the geometric constructions are initialised from the null initial data hypersurface,
  \item the exterior region stability result of~\cite{Kla.Nic03}, where the double null foliation is constructed by imposing that:
    \begin{itemize}
    \item the ingoing null hypersurfaces $\CCb_\ub$ emanate from the $2$-spheres of a canonical foliation of the initial spacelike hypersurface $\Si$,
    \item the outgoing null hypersurfaces $\CC_u$ are backwards emanating from the $2$-spheres of a canonical foliation on the last ingoing null hypersurface $\CCb^\ast$,
    \end{itemize}
    and which thus rather uses \emph{spacelike infinity} to initialise the constructions,
  \item the original stability result of~\cite{Chr.Kla93} where a global maximal time function is given, based on the existence of an asymptotically flat spacelike infinity,
  \item the stability of Schwarzschild spacetime established in~\cite{Kla.Sze17} (see also constructions for Kerr spacetime in~\cite{Kla.Sze19a,Kla.Sze19}), where the initialisation is performed from a \emph{last sphere} corresponding to timelike infinity, but which tracks \emph{the} central axis of Schwarzschild/Kerr spacetime. This is done by the construction of the last sphere as \emph{the intrinsic} GCM sphere (see~\cite[Section 7]{Kla.Sze19} or~\cite{Kla.Sze17}) and uses that the mass of the spacetime is \emph{non vanishing}. In the case of the perturbations of Minkowski space of Theorem~\ref{thm:mainv1}, there is no such intrinsic choice since there is no canonical central axis (the mass is not non-vanishing). We rather \emph{prescribe} a central axis, at the future timelike infinity of which we perform our initialisation.
  \end{itemize}
\item\label{item:interiorregion} We treat an \emph{interior region} enclosing the central axis, where constructions related to spherical coordinates \emph{degenerate} (\emph{i.e.} null pairs, null decompositions, etc.). This again differs from the semi-global results of Item~\ref{item:semiglobal}, as well as from the exterior region global stability result of~\cite{Kla.Nic03} and from the global stability of Schwarzschild spacetime\footnote{We recall that in the case of Schwarzschild spacetime, the axis $r=0$ is the Schwarzschild (spacelike) singularity and lies inside the Schwarzschild black hole.} from~\cite{Kla.Sze17}, which do not treat such an interior region and rely on null projected equations. In particular, the $r^p$-methods developed in~\cite{Hol10} for Bianchi equations projected on null pairs (see also the seminal $r^p$-method of~\cite{Daf.Rod10}) are unsuited for such a region. 

  In~\cite{Chr.Kla93}, the interior region is treated using the global maximal time function constructed from spacelike infinity and an interior optical function initialised at the central axis.

  In the case of Theorem~\ref{thm:mainv1}, no such global time function is available, and we rely instead on a construction of maximal hypersurfaces in the interior region, by prescribing their boundaries on the transition hypersurface, which is the timelike boundary between the interior and the exterior region. This is close in spirit to the proof of the spacelike-characteristic bounded $L^2$ curvature result from~\cite{Czi.Gra19a}.

  We moreover get rid of the interior optical function of~\cite{Chr.Kla93}, and rather use (transported) \emph{Cartesian harmonic coordinates}, which virtually makes any reference to the central axis disappear in the analysis of that region. See Section~\ref{sec:proofStabMinkintro}.
\item\label{item:DDM} Matching the (Cartesian) setting of the interior region to the (spherical) setting of the exterior region, is \emph{not performed} following functions or frames gluing procedures as in~\cite{Chr.Kla93,Kla.Sze17}, but is rather done by integration by parts, allowing discontinuities for the gauge choices at the interfaces. This procedure features a \emph{mean value argument} to avoid regularity losses at the timelike interface. We believe that this treatment -- although not fundamentally different -- gains in clarity with respect to previous works.
\item\label{item:GeomCK} At the core of the proof of Theorem~\ref{thm:mainv1} are the global energy estimates obtained by performing simultaneous energy estimates in the interior and exterior regions. Because of the different Cartesian/spherical setting used in each region, we choose to rely on the fully geometric framework (the spacetime Bel-Robinson tensors, the interior/exterior approximate conformal Killing vectorfields) of~\cite{Chr.Kla93} to match these estimates across the timelike transition hypersurface.
\end{enumerate}

\section{Organisation of the paper}
 We outline the structure of this paper.
\begin{itemize}
\item The geometric set up, definitions and formulas are collected in Section~\ref{sec:definition}. 
\item Section~\ref{sec:normsBA} is dedicated to collecting the definitions of the norms and bootstrap assumptions used in Sections~\ref{sec:globener}~--~\ref{sec:initlayer}.
\item Section~\ref{sec:mainresult} is dedicated to the statement of the main theorem as well as auxiliary theorems (global existence of harmonic coordinates, vertex/axis limits, etc.). We also set up and prove the bootstrap argument from which the conclusions of the main theorem follow. This relies on the improvement of the bootstrap assumptions, which is the core of this paper and is obtained in Sections~\ref{sec:globener}~--~\ref{sec:initlayer}.
\item In Section~\ref{sec:globener}, we perform the global energy estimates which are split into simultaneous energy estimates in the interior and exterior region. That is, we analyse and control the nonlinear error terms at the timelike transition interface, in the interior and in the exterior region.
\item In Sections~\ref{sec:curvest} and~\ref{sec:planehypcurvest}, we deduce from the bounds for the boundary energy fluxes of Section~\ref{sec:globener}, bounds for the curvature in the exterior and bottom interior region.
\item In Section~\ref{sec:remainingcurvestfinal}, we perform local energy estimates to obtain the remaining curvature bounds. In particular in the interior top region, we use an extension and local existence argument to obtain curvature bounds on the last cone.
\item In Sections~\ref{sec:connestCCba} and~\ref{sec:connest}, we obtain bounds for the null connection and rotation coefficients on the last cone and in the exterior region.
\item In Section~\ref{sec:planehypconnest}, we obtain bounds for the interior connection coefficients and approximate interior Killing vectorfields. We also obtain the bounds for the difference of Killing fields at the interface.
\item Section~\ref{sec:initlayer} is dedicated to the initial layer comparisons arguments.
\item Appendices~\ref{sec:globharmo} and~\ref{app:vertexlimits} are respectively dedicated to the obtention of global harmonic coordinates and vertex/axis limits, which are the main, independent and general results of the auxiliary theorems stated in Section~\ref{sec:mainresult}.
\item Appendices~\ref{app:canlocalex} and~\ref{app:KlSobH12} are dedicated to the obtention of auxiliary local existence and functional results.
\end{itemize}

  






\section{Acknowledgements} The author is very grateful to J\'er\'emie Szeftel for his support and encouragements. He also thanks Grigorios Fournodavlos for interesting discussions around the results of Theorem~\ref{thm:globharmonics} and Arick Shao for a detailed set of corrections and remarks on a first version of this paper. This work has been partially funded by the ERC grant ERC-2016-CoG 725589 EPGR and partially funded by the Deutsche Forschungsgemeinschaft (DFG, German Research Foundation) under Germany's Excellence Strategy EXC 2044 –390685587, Mathematics Münster: Dynamics–Geometry–Structure.

\chapter{Definitions \& formulas}\label{sec:definition}

Let $(\MMt,\g)$ be a smooth vacuum spacetime. Let $\Sit_1\subset\MMt$ be a smooth spacelike hypersurface diffeomorphic to the unit disk of $\RRR^3$ and $\CCt_0\subset\MMt$ the outgoing null hypersurface emanating from $\pr\Sit_1$. We call $(\Sit_1,\CCt_0)$ \emph{spacelike-characteristic initial data} and we assume that $\MMt$ coincides with the future maximal globally hyperbolic development of $\Sit_1\cup\CCt_0$.\\

Let $O$ be a fixed point of $\Sit_1$, which we call the \emph{centre of $\Sit_1$}.\footnote{The centre $O$ will be defined in Section~\ref{sec:mainresult} to be the point $x^i=0$, where $(x^i)$ are coordinates on $\Sit_1$ such that $\Sit_1$ is close to the unit Euclidean disk of $\RRR^3$.} We note $\o \subset \MMt$ the timelike geodesic emanating from $O$ orthogonally to $\Sit_1$ parametrised by arc-length and such that $\o(1) = O$. We call $\o$ the \emph{central axis} of the spacetime $\MMt$.

\section{The bootstrap region $\MM_{\protect\uba}$}\label{sec:defcannull}
In this section, we define the spacetime region $\MM_\uba\subset\MMt$ involved in our bootstrap argument, as well as its foliation by geometric hypersurfaces and its decomposition into subdomains.\footnote{The existence and regularity of the spacetime region $\MM_\uba\subset\MMt$ will be an assumption of the bootstrap argument of Section~\ref{sec:mainresult}.} \\

\begin{figure}[!h]
  \centering
  \includegraphics[height=12cm]{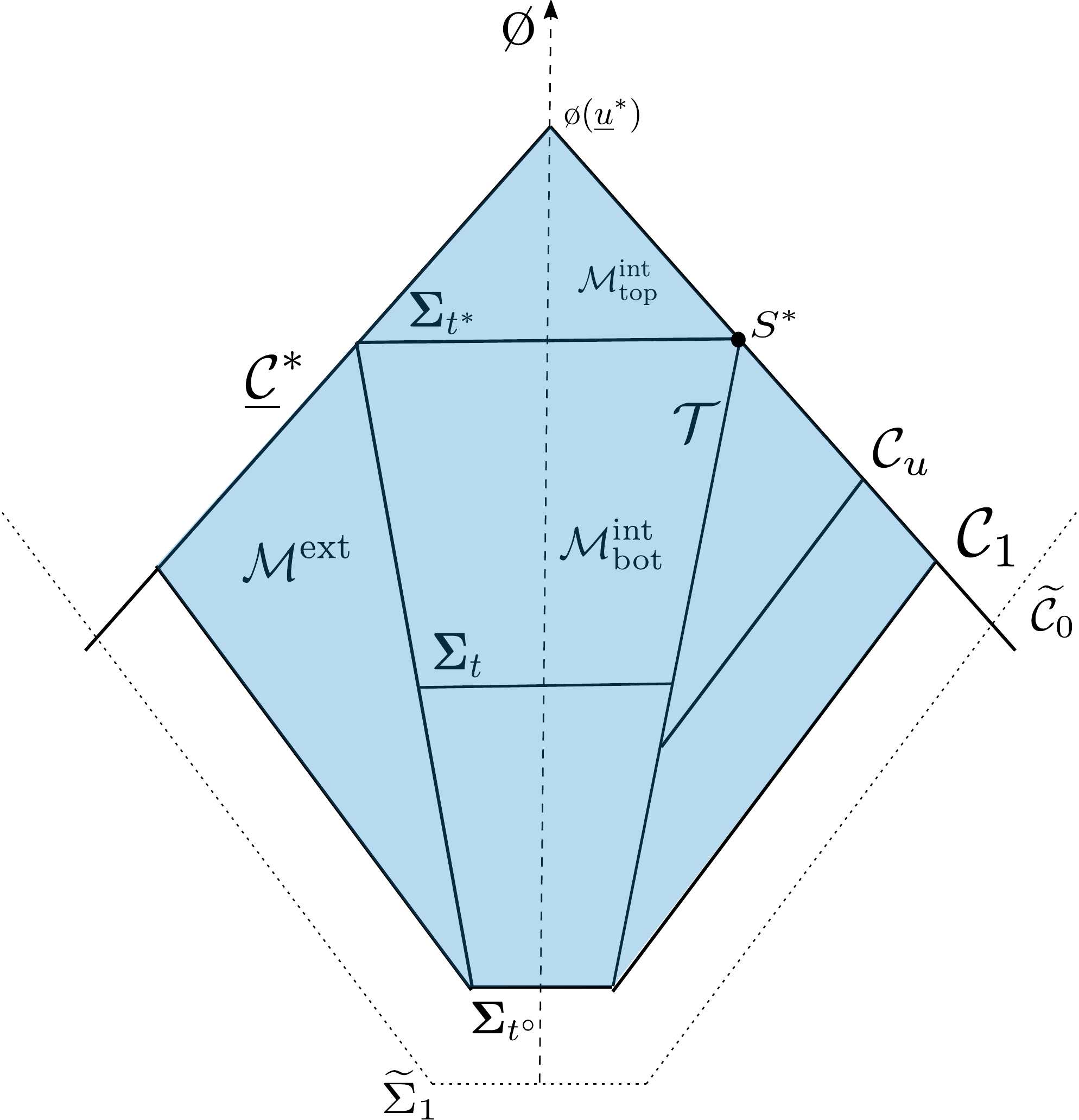}
  \caption{The spacetime region $\MM=\MM_\uba$.}
  \label{fig:MMuba}
\end{figure}

The constructions of this section are graphically summarised in Figure~\ref{fig:MMuba}.


\paragraph{The last cone $\CCba$ and the optical function $u$.}
Let $\uba \in[1,\infty)$. We define $\CCba$ to be the incoming null cone with vertex $\o(\uba)$. 
Let $u$ be a scalar function on $\CCba$. We assume that the level sets of $u$ define a foliation of $\CCba$ by regular $2$-spheres $S_{u,\uba}$, which will further be defined to be the $2$-spheres of the so-called \emph{canonical foliation} (see Definition~\ref{def:can}). We consider $\CC_u$ the outgoing null hypersurface backward emanating from the $2$-sphere $S_{u,\uba}\subset\CCba$. We define $u$ to be the associated optical function, and $L$ to be the (normalised) spacetime gradient of $u$, \emph{i.e.} $L := -(\D u)^{\sharp}$.

\paragraph{The geodesic parameter $\ub$.}
We define $\ub$ to be the following null geodesic affine parameter of the hypersurfaces $\CC_u$
\begin{align*}
  L(\ub) & = 2,\\
  \ub|_{\CCba} & = \uba.
\end{align*}
We assume that its level sets 
define a regular foliation of hypersurfaces transverse to the null hypersurfaces $\CC_u$ and that their intersections $S_{u,\ub}$ define spacelike regular $2$-spheres.\\

Let $\Yb := - (\D \ub)^\sharp$ be the (normalised) spacetime gradient of $\ub$. By definition of $\ub$, we have
\begin{align*}
  \g(\Yb,L) & = - 2.
\end{align*}
We define the \emph{optical defect} of the foliation $S_{u,\ub}$ to be the scalar function $\yy$ given by
\begin{align*}
  \g(\Yb,\Yb) =: -2\yy.
\end{align*}
We define the \emph{null pair  $(\elb,\el)$ associated to the foliation $S_{u,\ub}$} to be the pair of null vectorfields orthogonal to the $2$-spheres $S_{u,\ub}$ such that
\begin{align*}
  \el = L, \quad \text{and} \quad \g(\elb,\el) = -2.
\end{align*}
From the above definition, we have the following relations
\begin{align}\label{eq:elu}\begin{aligned}
    \elb(\ub) & = \yy, & \elb(u) & = 2, \\
    \el(\ub) & = 2, & \el(u) & = 0,
  \end{aligned}
\end{align}
and
\begin{align}\label{eq:relYelelb}
  \Yb & = \elb + \half\yy\el.
\end{align}

\paragraph{The canonical foliation on $\CCba$.}
We define $\chi,\xi,\eta,\ze,\om$ and $\chib,\xib,\etab,\ze,\omb$ to be the null decomposition of the connection coefficients associated to the null pair $(\elb,\el)$, and $\al,\be,\rho,\sigma,\beb,\alb$ to be the null decomposition of the spacetime curvature tensor with respect to the null pair $(\elb,\el)$ (see Section~\ref{sec:nulldecomp} for definitions).\\

With respect to these null decompositions, we have the following definition, which determines the function $u$ on $\CCba$. 
\begin{definition}[Canonical foliation on $\CCba$]\label{def:can}
  The scalar function $u$ is said to be \emph{canonical} on $\CCba$ if on each $2$-sphere $S_{u,\uba} \subset \CCba$, the following condition holds
  \begin{align}\label{eq:can}
    \begin{aligned}
      \Divd\eta + \rho & = \rhoo,\\
      \ombo & = 0, 
    \end{aligned}
  \end{align}
  where we refer to Section~\ref{sec:nulldecomp} for definitions, and if at the vertex $\o(\uba)$, the function $u$ is normalised by
  \begin{align*}
    u|_{\o(\uba)} & = \uba, & \g(\elb,\doto)|_{\o(\uba)} & = -1.
    \end{align*}
\end{definition}
\begin{remark}
  The null pair $(\elb,\el)$ is constructed to be adapted to the canonical foliation on $\CCba$ (\emph{e.g.} $\elb(u)=2$). It can be related to a geodesic null pair \emph{via} a null lapse. This will only be introduce in the proof of the existence of the canonical foliation in Appendix~\ref{app:canlocalex}. Other arguments do not require the existence or control of a background geodesic foliation.
\end{remark}
\begin{remark}
  The canonical foliation of Definition~\ref{def:can} does not coincide with the canonical foliations of~\cite{Kla.Nic03,Nic04,Czi.Gra19}, since we replaced for the conciseness of the argument the condition on the mean value of the null lapse in that papers by the condition $\ombo = 0$. This is purely to ease the notations and does not change anything to the motivations, its local or global existence on $\CCba$. See Theorem~\ref{thm:canonical}.
\end{remark}

\paragraph{The exterior region $\MM^\ext_{\protect\uba}$.} 
Let $0<\cco<1$ be a constant sufficiently close to $1$, which we call the \emph{transition constant}. Its value is determined in Section~\ref{sec:controlRRintr2} (see also Section~\ref{sec:constantsofthepaper} for a recapitulation of the constants of this paper and their dependencies).\\

Let $\cc$ be a \emph{transition parameter}, such that
\begin{align*}
   \cco \leq \cc \leq \half(1+\cco).
\end{align*}

For a fixed transition parameter, we define the \emph{exterior region} $^{(\cc)}\MM^\ext_\uba$ to be
\begin{align*}
  ^{(\cc)}\MM^\ext_\uba & := \le\{1 \leq u \leq \cc\ub\ri\},
\end{align*}
and we define the \emph{timelike transition hypersurface} $^{(\cc)}\TT$ to be
\begin{align*}
  ^{(\cc)}\TT & := \{u = \cc\ub\}.
\end{align*}

\begin{remark}\label{rem:transitionparameter}
  The freedom on the transition parameter $\cc$ is used to perform a \emph{mean value argument} (see Section~\ref{sec:meanvalue}). This mean value argument provides optimal regularity for the spacetime curvature tensor on \emph{a well chosen} timelike interface $^{(\cc)}\TT$. This enables a control of the nonlinear transition error terms. See Remark~\ref{rem:optregvsoptdecay} and Section~\ref{sec:estEETT12}.
\end{remark}

\paragraph{The interior region $\MM^\intr_{\protect\uba}$.}
Let 
\begin{align*}
  ^{(\cc)}\too & := \le(1+\cc^{-1}\ri)/2, & ^{(\cc)}\tast & := \le(1+\cc\ri)\uba/2.
\end{align*}
For all $\too \leq t \leq \tast$, we define $^{(\cc)}\Si_t$ to be the \emph{maximal hypersurfaces}, \emph{i.e.} such that
\begin{align*}
  \tr_g\kt & = 0,
\end{align*} 
where $g$ and $k$ are the first and second fundamental form of $^{(\cc)}\Si_t$ (see Section~\ref{sec:decSiHH} for definitions), with prescribed boundary
\begin{align*}
  \pr^{(\cc)}\Si_t & = S_{u= 2t\cc/(1+\cc), \ub = 2t/(1+\cc)} \subset ^{(\cc)}\TT,
\end{align*}
\emph{i.e.}
\begin{align}\label{eq:deftTT}
  ^{(\cc)}t\,|_{^{(\cc)}\TT} & = \half(u+\ub).
\end{align}

We define the \emph{bottom interior region} $^{(\cc)}\MM^\intr_{\bott,\uba}$ to be
\begin{align*}
  ^{(\cc)}\MM^\intr_{\bott,\uba} & := \le\{{^{(\cc)}\too} \leq {^{(\cc)}t} \leq {^{(\cc)}\tast}\ri\}.
\end{align*}

We call $^{(\cc)}\Si_\tast$ the \emph{last maximal slice} and we note its boundary $^{(\cc)}S^\ast := \pr^{(\cc)}\Si_\tast = S_{\cc\uba,\uba} \subset\CCba$. We define the \emph{top interior region} $^{(\cc)}\MM^\intr_{\topp,\uba}$ to be the domain of dependence of the last maximal hypersurface $^{(\cc)}\Si_\tast$.
\begin{remark}
  Since $\pr^{(\cc)}\Si_\tast \subset \CCba$, the future boundary of $^{(\cc)}\MM^\intr_{\topp,\uba}$ is contained in the last cone $\CCba$.
\end{remark}
We define the \emph{interior region} $^{(\cc)}\MM^\intr_\uba$ and the spacetime region $^{(\cc)}\MM_\uba$ by
\begin{align*}
  ^{(\cc)}\MM^\intr_\uba & := ^{(\cc)}\MM^\intr_{\bott,\uba}\cup{^{(\cc)}\MM}^\intr_{\topp,\uba},\\
  ^{(\cc)}\MM_\uba & := ^{(\cc)}\MM^\ext_\uba \cup {^{(\cc)}\MM}^\intr_\uba.
\end{align*}

From now on, we shall drop the labels $\uba,\cc$ in the above geometric constructions, unless that precision is relevant.






 

\section{Maximal hypersurfaces $\Si_t$ of $\MM^\intr_\bott$}

\subsection{Fundamental forms}\label{sec:decSiHH}
The \emph{first fundamental form} is defined to be the induced metric $g$ by $\g$ on $\Si_t$.\\

We define the \emph{time lapse} $\nt$ of the maximal hypersurfaces $\Si_t$ by
\begin{align}\label{eq:defnt}
  \nt^{-2} & := -\g(\D t,\D t).
\end{align}
We define $\Tf$ to be the future-oriented unit normal to $\Si_t$. It satisfies the following relations
\begin{align}\label{eq:relTnt}
  \Tf & = -\nt(\D t)^\sharp, & \D_\Tf\Tf & = \nt^{-1}\nab \nt, & \Tf(t) & = \nt^{-1},
\end{align}
where $\nab$ denotes the induced covariant derivative on $\Si_t$.\\

We define the \emph{second fundamental form} $\kt$ of $\Si_t$ to be
\begin{align}\label{eq:defkt}
  \kt(X,Y) & := -\g(\D_X\Tf,Y), 
\end{align}
for all $X,Y$ in $T\Si_t$. From the definitions of Section~\ref{sec:defcannull}, we recall that we have
\begin{align}\label{eq:defmax}
  \tr_g\kt & = 0.
\end{align}




\subsection{Electric-magnetic decompositions}\label{sec:defelecmag}
The \emph{electric-magnetic} decomposition of the curvature tensor $\R$ consists of the $\Si_t$ tangent tensors $\Et$ and $\Ht$ defined by
\begin{align}\label{eq:defEH}
  \Et(X,Y) & := \R(\Tf,X,\Tf,Y), & \Ht(X,Y) & := \dual\R(\Tf,X,\Tf,Y),
\end{align}
for all $X,Y$ in $T\Si_t$.\\

Using the maximal condition~\eqref{eq:defmax}, we have the following Gauss equation for the Ricci curvature tensor of $\Si_t$
\begin{align}\label{eq:GaussSit}
  \RRRic_{ij} & = \Et_{ij} + \kt_{ia}\kt^{a}_j.
\end{align}

Using the maximal condition~\eqref{eq:defmax}, we have the following Hodge-type system for $\kt$
\begin{align}\label{eq:Hodgekt}
  \Div \kt & = 0, & \Curl \kt & = \Ht, & \tr\kt & = 0,
\end{align}
and Laplace-type equation for $\nt$
\begin{align}\label{eq:Lapnt}
  \Delta \nt & = \nt |\kt|^2,
\end{align}
see~\cite[p. 6]{Chr.Kla93}.\\

One also has the following equation for $\Lieh_\Tf\kt$ (see~\cite[p. 6]{Chr.Kla93})
\begin{align}\label{eq:LiehTfkt}
  n\Lieh_\Tf \kt & = -\nab^2\nt + \nt\le(\Et - \kt \cdot \kt\ri),
\end{align}
where for a symmetric $\Si_t$-tangent $2$-tensor $F$
\begin{align*}
  (F\cdot F)_{ij} & := F_{il}F^{l}_j. 
\end{align*}

The tensors $\Et$ and $\Ht$ satisfy the following Maxwell-type system of equations (see~\cite[pp. 144-146]{Chr.Kla93})
\begin{align}\label{eq:Maxwellt}
  \begin{aligned}
    \Div \Et & = \kt\wedge \Ht, \\
    \Div \Ht & = -\kt\wedge \Et,\\
    -\Lieh_\Tf \Ht +\Curl \Et & = -\nt^{-1}\nab \nt \wedge \Et -\half \kt \times \Ht,\\
    \Lieh_\Tf \Et + \Curl \Ht & = -\nt^{-1}\nab \nt\wedge \Ht + \half \kt\times \Et,
  \end{aligned}
\end{align}
where for a $\Si_t$-tangent symmetric traceless $2$-tensor $U$, $\Lieh_\Tf U_{ij} := \Lie_\Tf U_{ij} +\frac{2}{3}(k\cdot U)g_{ij}$, and where for symmetric $\Si_t$-tangent $2$-tensors $F,G$ and a $\Si_t$-tangent $1$-tensor $X$, we have
\begin{align*}
  (F\wedge G)_{i} & := \in_i^{mn}F_m^lG_{ln},\\
  (X\wedge F)_{ij} & := \in_{i}^{mn}X_mF_{nj}+\in_j^{mn}X_mF_{in},\\
  (F\times G)_{ij} & := \in_i^{ab}\in_j^{cd}F_{ac}G_{bd}+\frac{1}{3}(F_{ab}G^{ab})g_{ij}.
\end{align*}

\section{Null decompositions}\label{sec:nulldecomp}
In this section, we recall the general null decomposition of the connection and curvature as defined in~\cite{Chr.Kla93}.\\

Let $S_{u,\ub}$ be a (local) foliation of $\MM$ by spacelike $2$-spheres. We note $\gd$ the Riemannian metric induced by $\g$ on $S_{u,\ub}$, and we note $r$ its area radius $4\pi r^2:=|S|$. A \emph{null pair} $(\el,\elb)$ is a pair of vectorfields on $\MM$ orthogonal to the $2$-spheres $S_{u,\ub}$ such that
\begin{align*}
  \g(\el,\elb) = -2, \quad \g(\el,\el) = \g(\elb,\elb) = 0.
\end{align*}

We define the connection coefficients relative to a null pair $(\el,\elb)$ to be the $S$-tangent tensors defined by 
\begin{align*}
  \chi(X,Y) & := \g(\D_X\el,Y), & \xi(X)  & := \half\g(\D_4\el,X), & \eta(X) & := \half\g(\D_3\el,X),\\   \ze(X) & := \half \g(\D_X\el,\elb), & \om  & := \frac{1}{4} \g(\D_4\el,\elb),    
\end{align*}
and
\begin{align*}
    \chib(X,Y) & := \g(\D_X\elb,Y), & \xib(X) & := \half\g(\D_3\elb,X), & \etab(X) & := \half\g(\D_4\elb,X),\\  
  \ze(X) & := -\half \g(\D_X\elb,\el), & \omb & := \frac{1}{4} \g(\D_3\elb,\el),
\end{align*}
where $X,Y \in \text{T}S$.

\begin{remark}
  These correspond to the same definitions as in~(\cite[p. 147]{Chr.Kla93}) with
  \begin{align*}
    \chi = H,~ \xi = Y,~ \eta = Z,~ \ze = V,~\om = \Om, 
  \end{align*}
  and
  \begin{align*}
    \chib = \underline{H}, ~ \xib = \underline{Y},~ \etab = \underline{Z},~ \ze = V, ~ \omb = \underline{\Om}. 
  \end{align*}
\end{remark}

For a (local) orthonormal frame $(\ea)_{a=1,2}$ on $\text{T}S$, we have the following relations for the covariant derivatives of the orthonormal null frame $(\el,\elb,e_a)$ (see~\cite[p. 147]{Chr.Kla93})
\begin{align}\label{eq:Riccirel}
  \begin{aligned}
    \D_a\elb & = \chib_{ab}e_b + \ze_a\elb, & \D_a\el & = \chi_{ab}\eb-\ze_a\el,\\
    \D_3\elb & = 2\xib_a\ea-2\omb\elb, & \D_3\el & = 2\eta_a\ea + 2\omb\el,\\
    \D_4\elb & = 2\etab_a\ea+2\om\elb, & \D_4\el & = 2\xi_a\ea-2\om\el, \\
    \D_3\ea & = \Nd_3\ea + \eta_a\elb + \xib_a\el, & \D_4\ea & = \Nd_4\ea + \xi_a\elb + \etab_a\el,\\
    \D_a\eb & = \Nd_a\eb + \half\chi_{ab}\elb+\half\chib_{ab}\el.
  \end{aligned}
\end{align}

We define the curvature components relative to a null pair $(\el,\elb)$ to be the $S$-tangent tensors defined by
\begin{align*}
  \alpha(X,Y) & := \R(\el,X,\el,Y), & \beta(X) & := \half \R(X,\el,\elb,\el), & \rho & := \frac{1}{4} \R(\elb,\el,\elb,\el), \\
  \alphab(X,Y) & := \R(\elb,X,\elb,Y), & \betab(X) & := \half \R(X,\elb,\elb,\el), & \sigma & := \frac{1}{4}\dual\R(\elb,\el,\elb,\el), 
\end{align*}
where $X,Y \in \text{T}S$ and where $\dual\R$ denotes the Hodge dual of $\R$.\\

We have the following null structure relating the null connection coefficients and the null curvature components (see~\cite[pp. 168-170]{Chr.Kla93}).\\

\begin{subequations}\label{eq:null}
  We have the following transport equations along $\elb$ and $\el$
  \begin{align}
    \Lied_{\elb}\gd & = 2\chib, & \Lied_\el\gd & = 2\chi,\label{eq:elbgd}
  \end{align}
  and
  \begin{align}
    \Nd_3\chibh + \trchib\chibh & = \Nd\otimesh\xib - 2\omb\chibh + \le(\eta+\etab-2\ze\ri)\otimesh\xib - \alb, \label{eq:Nd3chibh}\\
    \Nd_3\trchib +\half(\trchib)^2 & = 2\Divd\xib-2\omb\trchib + 2\xib\cdot\le(\eta+\etab-2\ze\ri) -|\chibh|^2,  \label{eq:Nd3trchib}\\
    \Nd_3\ze & = -2\Nd\omb - \chib\cdot(\ze+\eta) +2\omb(\ze-\eta) +\chi\cdot\xib +2\om\xib -\beb,\label{eq:Nd3ze}\\
    \Nd_3\chih + \half\trchib \chih & = \Nd\otimesh\eta +2\omb\chih -\half\trchi\chibh + \xib\otimes\xi + \eta\otimesh\eta,\label{eq:Nd3chih}\\
    \Nd_3\trchi + \half\trchib\trchi & = 2\Divd\eta +2\omb\trchi-\chih\cdot\chibh +2(\xi\cdot\xib +|\eta|^2) +2\rho,\label{eq:Nd3trchi}\\
    \Nd_3\xi - \Nd_4\eta & = 4\omb\xi + \chi\cdot\le(\eta-\etab\ri) + \beta,\label{eq:Nd3xiOLD}\\
    \Nd_3\etab -\Nd_4\xib & = - 4\om\xib - \chib\cdot(\etab-\eta) +\betab,\label{eq:Nd3etab}\\
    \Nd_3\om +\Nd_4\omb & = \xi\cdot\xib + \ze\cdot(\eta-\etab) - \eta\cdot\etab + 4\om\omb+\rho.\label{eq:Nd3om}\\
    \Nd_4\chih + \trchi\chih & = \Nd\otimesh\xi - 2\om\chih + \le(\eta+\etab+2\ze\ri)\otimesh\xi - \al, \label{eq:Nd4chih}\\
    \Nd_4\trchi +\half(\trchi)^2 & = 2\Divd\xi-2\om\trchi + 2\xi\cdot\le(\eta+\etab+2\ze\ri) -|\chih|^2, \label{eq:Nd4trchi}\\
    \Nd_4\ze & = 2\Nd\om + \chi\cdot(-\ze+\etab) +2\om(\ze+\etab) -\chib\cdot\xi -2\omb\xi -\beta,\label{eq:Nd4ze}\\
    \Nd_4\chibh + \half\trchi \chibh & = \Nd\otimesh\etab +2\om\chibh -\half\trchib\chih + \xi\otimes\xib + \etab\otimesh\etab,\label{eq:Nd4chibh}\\
    \Nd_4\trchib + \half\trchi\trchib & = 2\Divd\etab +2\om\trchib-\chibh\cdot\chih +2(\xib\cdot\xi +|\etab|^2) +2\rho.\label{eq:Nd4trchib}
  \end{align}
  We have the following elliptic equations on the $2$-spheres
  \begin{align}
    \Curld\etab & = - \Curld\eta = \half \chih\wedge\chibh - \xi\wedge\xib - \sigma,\label{eq:Curlze}\\
    \Divd\chibh & = \Nd\trchib + \chib\cdot\ze - \trchib\ze + \betab,\label{eq:Divdchibh}\\
    \Divd\chih & = \Nd\trchi -\chi\cdot\ze + \trchi\ze - \beta,\label{eq:Divdchih}\\
    \Curld\xi & = \xi\wedge\le(\eta+\etab+2\ze\ri),\label{eq:Curldxi}\\
    \Curld\xib & = \xib\wedge\le(\eta+\etab-2\ze\ri),\label{eq:Curldxib}\\
    K & = -\quar \trchi\trchib + \half \chih\cdot\chibh - \rho.\label{eq:Gauss}
  \end{align}
\end{subequations}

\begin{subequations}\label{eq:bianchi}
  We have the following Bianchi equations relating the null connection coefficients and the null curvature components (see~\cite[p. 161]{Chr.Kla93}).
  \begin{align}
    \Nd_3\al +\half\trchib \al & = \Nd\otimesh\beta + 2\omb\al - 3(\chih\rho+\dual\chih\sigma) + (\ze+4\eta)\otimesh\be,\label{eq:Nd3al}\\
    \Nd_4\beta +2\trchi\be & = \Divd\al -2\om\be+(2\ze+\etab)\cdot\al +3(\xi\rho+\dual\xi\sigma),\label{eq:Nd4be}\\
    \Nd_3\be +\trchib\be & = \Nd\rho + \dual\Nd\sigma + 2\omb\be+\xib\cdot\alpha+3(\eta\rho+\dual\eta\sigma),\label{eq:Nd3be}\\
    \Nd_4\rho + \frac{3}{2}\trchi\rho & = \Divd\be -\half \chibh\cdot\alpha + \ze\cdot\beta +2(\etab\cdot\be-\xi\cdot\betab),\label{eq:Nd4rho}\\
    \Nd_3\rho + \frac{3}{2}\trchib\rho & = - \Divd\beb -\half\chih\cdot\alphab + \ze\cdot\beb+2(\xib\cdot\be-\eta\cdot\beb),\label{eq:Nd3rho}\\
    \Nd_4\sigma +\frac{3}{2}\trchi\sigma & = - \Curld\be + \half\chibh\cdot\dual\alpha-\ze\cdot\dual\be-2(\etab\cdot\dual\be+2\xi\cdot\dual\beb),\label{eq:Nd4sigma}\\
   \Nd_3\sigma + \frac{3}{2}\trchib \sigma & = -\Curld\beb -\half\chih\cdot\dual\alb+\ze\cdot\dual\beb-2(\etab\cdot\dual\be+\eta\cdot\dual\beb),\label{eq:Nd3sigma}\\
    \Nd_4\beb + \trchi \beb & = -\Nd\rho +\dual\Nd\sigma + 2\chibh\cdot\be + 2\om\beb-\xi\cdot\alb-3(\etab\rho-\dual\etab\sigma),\label{eq:Nd4beb}\\
    \Nd_3\beb +2\trchib\beb & = -\Divd\alb -2\omb\beb -(-2\ze+\eta)\cdot\alb+3(-\xib\rho+\dual\xib\sigma),\label{eq:Nd3beb}\\
    \Nd_4\alb +\half \trchi \alb & = -\Nd\otimesh\beb + 4\om\alb -3(\chibh\rho-\dual\chibh\sigma) +(\ze-4\etab)\otimesh\beb.\label{eq:Nd4alb}
  \end{align}
\end{subequations}

We have the following commutations formulas for the covariant derivatives $\Nd,\Nd_4,\Nd_3$ (see~\cite[p. 159]{Chr.Kla93})
\begin{subequations}\label{eq:comm}
  \begin{align}
    [\Nd_4,\Nd] F & = -\half\trchi \Nd F - \chih\cdot\Nd F + \xi \Nd_3F + (\etab+\ze)\Nd_4F + \EEE(\Nd_4,\Nd)\cdot F,\label{eq:commNd4Nd}\\
    [\Nd_3,\Nd] F & = -\half\trchib \Nd F - \chibh\cdot\Nd F+ \xib\Nd_4F + (\eta-\ze)\Nd_3F +\EEE(\Nd_3,\Nd)\cdot F,\label{eq:commNd3Nd}\\
    [\Nd_3,\Nd_4] F & = 2\om\Nd_3F+2\omb\Nd_4F + (\eta-\etab)\cdot\Nd F + \EEE(\Nd_3,\Nd_4)\cdot F,\label{eq:commNd3Nd4}
  \end{align}  
  where $F$ is a $S$-tangent $k$-tensor and where the tensors $\EEE$ are given by
  \begin{align*}
    \le(\EEE(\Nd_4,\Nd) \cdot F\ri)_{a a_1\cdots a_k} & := \sum_{i=1}^k \le(\chib_{a_i a}\xi_b-\chib_{ab}\xi_{a_i} + \chi_{a_i a}\etab_b-\chi_{ab}\etab_{a_i} + \ep_{a_ib}\dual\be_a\ri) F_{a_1 \cdots b \cdots a_k} ,\\
    \le(\EEE(\Nd_3,\Nd) \cdot F\ri)_{a a_1\cdots a_k} & := \sum_{i=1}^k \le(\chi_{a_i a}\xib_b-\chi_{ab}\xib_{a_i} + \chib_{a_i a}\eta_b-\chib_{ab}\eta_{a_i} - \ep_{a_ib}\dual\beb_a\ri) F_{a_1 \cdots b \cdots a_k} ,\\
    \le(\EEE(\Nd_3,\Nd_4)\cdot F\ri)_{a_1 \cdots a_k} & := 2\sum_{i=1}^k \le(\xi_{a_i}\xib_{b}-\xi_b\xib_{a_i} + \etab_{a_i}\eta_b-\etab_{b}\eta_{a_i}+\ep_{a_i b}\sigma\ri)F_{a_1 \cdots b \cdots a_k}.
  \end{align*}
\end{subequations}


\section{Transition relations on the boundary $\TT$}
In this section, we derive boundary relations between the maximal hypersurface decompositions of Sections~\ref{sec:decSiHH} and~\ref{sec:defelecmag} and the null decompositions of Section~\ref{sec:nulldecomp}.\\

There exists $\nut>0$ such that on $\TT$, we have
\begin{align}\label{eq:defnutnur}
  \Tf & := \half\le(\nut^{-1}\elb + \nut\el\ri),
\end{align}
and the outward-pointing unit normal to the boundary $\pr\Si_t$ writes
\begin{align}\label{eq:NTT}
  \Ntf & := \half\le(-\nut^{-1}\elb+\nut\el\ri).
\end{align}

We define the boundary decompositions of $\kt$ by
\begin{align}\label{eq:defbdydecomp}
  \begin{aligned}
  \delt & := \kt_{\Ntf\,\Ntf}, & \ept_{a} & := \kt_{\Ntf a}, & \kapt_{ab} & := \kt_{ab}.
  \end{aligned}
\end{align}

\begin{lemma}\label{lem:transkTT}
  We have
\begin{align}\label{eq:transkTTt}
  \begin{aligned}
    \kapt_{ab} & = -\half\nut^{-1}\chib_{ab} - \half\nut\chi_{ab}, \\
    \ept_{a} & = -\Nd_a(\log\nut) +\ze_a,\\
    \delt & = \half \nut^{-1}\D_3(\log\nut) -\half\nut\D_4(\log\nut) + \nut^{-1}\omb + \nut\om.
  \end{aligned}
\end{align}
\end{lemma}

\begin{proof}
  The first relations for $\kapt$ immediately follow from~\eqref{eq:Riccirel} and~\eqref{eq:defkt}. Using \eqref{eq:NTT}, we have
  \begin{align*}
    \kt_{\Ntf a} & = -\quar\g\le(\D_{a}\le(\nut^{-1}\elb+\nut\el\ri), -\nut^{-1}\elb+\nut\el\ri) \\
              & = - \Nd_a(\log \nut) +\ze_a,
  \end{align*}
  and
  \begin{align*}
    \kt_{\Ntf\,\Ntf} & = -\frac{1}{8}\g\le(\D_{-\nut^{-1}\elb+\nut\el}\le(\nut^{-1}\elb+\nut\el\ri),-\nut^{-1}\elb+\nut\el\ri) \\
               & = \half \nut^{-1} \D_3(\log \nut) - \half \nut \D_4(\log \nut) \\
               & \quad + \frac{1}{8}\nut^{-1}\g(\D_{\elb}(\elb+\el),-\elb+\el) \\
                   & \quad -\frac{1}{8}\nut\g(\D_{\el}(\elb+\el),-\elb+\el) \\
               & = \half \nut^{-1} \D_3(\log \nut) - \half \nut \D_4(\log \nut) + \nut^{-1}\omb + \nut \om.
  \end{align*}
 This finishes the proof of the lemma.
\end{proof}

\begin{lemma}\label{lem:relntnut}
  We have at $\TT$
  \begin{align}\label{eq:relntnut}
    \nt & = \frac{\cc}{\cc+1}\le(\nut+\le(\cc^{-1}-\half\yy\ri)\nut^{-1}\ri).
  \end{align}
\end{lemma}
\begin{proof}
Using relations~(\ref{eq:elu}), the future-pointing unit normal $Z$ to $S_{u,\ub}$ in $\TT$ writes
\begin{align*}
  Z & = \half \le(a^{-1/2}\elb +a^{1/2}\el\ri),
\end{align*}
with
\begin{align*}
  a & := \cc^{-1} -\half \yy.
\end{align*}

Using that we have on $\TT$
\begin{align*}
  Z(t) & = \half Z(u+\ub),
\end{align*}
we therefore deduce at $\TT$ that
\begin{align}\label{eq:Z(t)}
  \begin{aligned}
    a^{-1/2}\elb(t) + a^{1/2}\el(t) & = \half(a^{-1/2}\elb+a^{1/2}\el)(u+\ub) \\
    & = a^{-1/2} +a^{1/2} + \half a^{-1/2}\yy.
  \end{aligned}
\end{align}

Moreover from the relations~(\ref{eq:relTnt}) and~definitions~(\ref{eq:defnutnur}), (\ref{eq:NTT}), we have at $\TT$
\begin{align*}
  \begin{aligned}
    \half \nut^{-1} \elb(t) + \half \nut \el(t) & = \nt^{-1},\\
    -\half\nut^{-1}\elb(t) + \half\nut\el(t) & = 0,
  \end{aligned}
\end{align*}
which gives
\begin{align}
  \begin{aligned}
    \elb(t) & = \nut\nt^{-1}, & \el(t) & = \nut^{-1}\nt^{-1}. 
  \end{aligned}
\end{align}
Plugging this into~\eqref{eq:Z(t)} gives
\begin{align*}
  \nt & = \frac{a^{-1/2}\nut+a^{1/2}\nut^{-1}}{a^{-1/2}+a^{1/2}+ \half a^{-1/2} \yy},
\end{align*}
which, after rewriting, using the expression for $a$, gives the desired formula.
\end{proof}

\section{Uniformisation of the sphere $S^\ast$ and harmonic Cartesian coordinates on $\Si_\tast$}\label{sec:defUnifandharmo}
We recall the following definition from Section~\ref{sec:defcannull}
\begin{align*}
  S^\ast & := \pr\Si_\tast = S_{u=\cc\uba,\uba},
\end{align*}
and we note $(r^\ast)^2 := \frac{1}{4\pi}|S^\ast|$ its area radius.\\

A \emph{conformal isomorphism} between $S^\ast$ and the Euclidean unit $2$-sphere $\SSS$ is a diffeomorphism $\Phi~:~S^\ast \to \SSS$ such that there exists a \emph{conformal factor} $\phi>0$ on $\SSS$ satisfying
\begin{align*}
  \Phi_{\sharp}\gd_{S^\ast} & = (r^\ast)^2\phi^2\gd_{\SSS},
\end{align*}
where $\Phi_\sharp\gd_{S^\ast}$ is the push-forward of the metric $\gd$ by $\Phi$.\\

To a fixed conformal isomorphism of $S^\ast$, we associate the (normalised) \emph{Cartesian coordinates} $\le(\frac{x^i}{r^\ast}\ri)_{i=1\ldots 3}$ on $S^\ast$ to be the pull-back by $\Phi$ of the standard Cartesian coordinates on $\SSS$.\\

We say that the conformal isomorphism $\Phi$ is \emph{centred} if the functions $x^i$ satisfy the following conditions on $S^\ast$
\begin{align*}
  \int_{S^\ast} x^i & = 0, & i & = 1, \cdots, 3.
\end{align*}

Using these coordinates, we further define the associated \emph{harmonic Cartesian coordinates} on $\Si_\tast$ to be the solution of the following Dirichlet problem on $\Si_\tast$
\begin{align}\label{def:harmocartcoords}
  \begin{aligned}
    \Delta_g x^i & = 0,\\
    x^i|_{S^\ast = \pr\Si_\tast} & = x^i,
  \end{aligned}
\end{align}
for all $i=1,2,3$.

\begin{remark}
  Further constructions and bounds in this paper will hold for all centred conformal isomorphisms. In this paper, we will use the existence, uniqueness and control of all the centred conformal isomorphism established in~\cite[Theorem 3.1]{Kla.Sze19}.
\end{remark}

\section{Approximate interior Killing fields $\TI$, $\SI$, $\KI$ and $\OOI$ in $\MM^\intr_\bott$}\label{sec:defKillingint}
We first define the \emph{approximate Killing time translation vectorfield} $\TI$ on $\MM^\intr_\bott$ by
\begin{align*}
  \TI & := \Tf.
\end{align*}

We define $\XI$ on $\Si_\tast$ to be the $\Si_\tast$-tangent vectorfield given by
\begin{align}\label{eq:defXISitast}
  \XI & := \sum_{i=1}^3 x^i \nab x^i,
\end{align}
where $\nab$ is the induced gradient on $\Si_\tast$ and $x^i$ are the harmonic Cartesian coordinates defined in Section~\ref{sec:defUnifandharmo}.\\

We define $\XI$ on $\MM^\intr_\bott$ by parallel transport, \emph{i.e.}
\begin{align}\label{eq:defXI}
  \begin{aligned}
  \D_\Tf\XI & = 0,
  \end{aligned}
\end{align}

We define the \emph{approximate conformal Killing scaling vectorfield} $\SI$ on $\MM^\intr_\bott$ by
\begin{align}\label{eq:defSI}
  \SI & := t\TI + \XI.
\end{align} 

We define the \emph{approximate conformal Killing Morawetz vectorfield} $\KI$ on $\MM^\intr_\bott$ by
\begin{align}\label{eq:defKI}
  \begin{aligned}
    \KI & := \le(t^2 +\g(\XI,\XI)\ri)\TI + 2t\XI.
  \end{aligned}
\end{align}

We define the \emph{approximate Killing rotation vectorfields} $\OOI$ on $\Si_\tast$ by
\begin{align}\label{eq:defOOISitast}
  ^{(\ell)}\OOI & := \in_{\ell i j} x^i\nab x^j,
\end{align}
for $\ell = 1,2,3$, and we extend them on $\MM^\intr_\bott$ by parallel transport along $\Tf$, \emph{i.e.}
\begin{align}\label{eq:defOOI}
  \D_\Tf {^{(\ell)}\OOI} & = 0.
\end{align}

We recall that the spacetime \emph{deformation tensor} of a spacetime vectorfield $X$ is given by
\begin{align*}
  ^{(X)}\pi_{\mu\nu} & := \D_\mu X_\nu + \D_\nu X_\mu,
\end{align*}
and that we note $\pih$ its traceless part, \emph{i.e.} $\pih := \pi - \quar \tr\pi \g$.\\

Using the maximal condition~(\ref{eq:defmax}), we have
\begin{align}\label{eq:trpiTI}
  \tr\le(^{(\TI)}\pi\ri) & = 0,
\end{align}
and using relations~\eqref{eq:relTnt}, we have
\begin{align}\label{eq:piTI}
  ^{(\TI)}\pi_{\Tf\Tf} & = 0, & ^{(\TI)}\pi_{\Tf i} & = \nt^{-1}\nab_i\nt, & ^{(\TI)}\pi_{ij} & = -2\kt_{ij}.
\end{align}

\section{Exterior hypersurfaces $\Si_t^\ext$ of $\MM^\ext$}\label{sec:defSitext}
For all $\too \leq t \leq \tast$, we define the exterior spacelike hypersurfaces $\Si^\ext$ by
\begin{align*}
  \Si^\ext_t & := \le\{u+\ub = 2t\ri\}\cap\MM^\ext.
\end{align*}
From~\eqref{eq:elu}, one has the following relations for the future-pointing unit normal $\Tf^\ext$ to $\Si^\ext_t$ and the outward-pointing unit normal $\Nf^\ext$ to the $2$-spheres $S_{u,\ub}\subset\Si^\ext_t$ in $\Si^\ext_t$ 
\begin{align}
  \label{eq:relTfNfextelelb}
  \begin{aligned}
    \Tf^\ext & = \half\le(1+\half\yy\ri)^{-1/2}\elb + \half \le(1+\half\yy\ri)^{1/2}\el,\\
    \Nf^\ext & = -\half\le(1+\half\yy\ri)^{-1/2}\elb + \half \le(1+\half\yy\ri)^{1/2}\el.
  \end{aligned}
\end{align}
Moreover, we have the following definitions and expression of the time lapse $n^\ext$
\begin{align}
  \label{eq:timelapseext}
  \begin{aligned}
    n^\ext & := \le(-\g\le(\D\le(\half(u+\ub)\ri), \D\le(\half(u+\ub)\ri)\ri)\ri)^{-1/2}\\
    & = \le(1+\half\yy\ri)^{-1/2}.
  \end{aligned}
\end{align}

\section{Commutation relations for integrals and averages on $S_{u,\protect\ub}$}

We have the following commutation relation between $\elb,\el$-derivatives and the integral on a $2$-sphere $S_{u,\ub}$
\begin{lemma}
  For all scalar function $\phi$, we have
  \begin{align}\label{eq:commelbint}
    \begin{aligned}
      \elb\le(\int_S \phi\ri) & = \int_S \le(\elb(\phi) + \trchib \phi\ri) &&& \text{on $\CCba$},\\
      \el\le(\int_S \phi\ri) & = \int_S \le(\el(\phi) + \trchi \phi\ri) &&& \text{on $\MM^\ext$}.
    \end{aligned}
  \end{align}
\end{lemma}
\begin{proof}
  Using~\eqref{eq:elu} and defining appropriate transported spherical coordinates, we have $\el = 2\pr_\ub$ in $\MM^\ext$, which thus commutes with $\int_{S_{u,\ub}}$. Using~\eqref{eq:elbgd}, the desired formula in $\MM^\ext$ follows.\\

  Since $\CCba$ is null, $\Yb$ is null and $\yy = 0$. Using~\eqref{eq:elu}, and defining appropriate transported spherical coordinates, we have $\elb = 2\pr_u$ on $\CCba$. Using~\eqref{eq:elbgd}, the desired formula follows on $\CCba$. This finishes the proof of the lemma.
\end{proof}

Using~\eqref{eq:commelbint}, we also have the following commutation relation between $\elb,\el$-derivatives and the mean value on a $2$-sphere $S_{u,\ub}$
\begin{align}\label{eq:commelbov}
  \begin{aligned}
    \elb\le(\overline{\phi}\ri) & = \overline{\elb(\phi)} + \overline{\le(\trchib-\trchibo\ri)\phi}& \text{on $\CCba$},\\
    \el\le(\overline{\phi}\ri) & = \overline{\el(\phi)} + \overline{\le(\trchi-\trchio\ri)\phi} & \text{on $\MM^\ext$},
    \end{aligned}
\end{align}
where $\phi$ is a scalar function and $\overline{\phi} := (4\pi r^2)^{-1} \int_S\phi$.\\

Using~\eqref{eq:commelbint}, we also have
\begin{align}\label{eq:elbr}
  \elb(r) & = \half r \trchibo & \text{on $\CCba$},\\
  \el(r) & = \half r\trchio, & \text{on $\MM^\ext$}.
\end{align}

\section{Null decomposition of the geodesic-null foliation in $\MM^\ext$}\label{sec:nulldecompMMext}
In this section, we derive additional equations to the general null decompositions defined in Section~\ref{sec:nulldecomp}, in the case of the null pair $(\elb,\el)$ of the geodesic-null foliation defined in Section~\ref{sec:defcannull}.\\

We have the following relations.
\begin{lemma}\label{lem:relgeodnull}
 In $\MM^\ext$ the following relations hold
 \begin{align}
   \xi & = 0,\label{eq:xib0}\\
   \om & = 0, \label{eq:Nd4Om}\\
   0 & = \eta -\ze,\label{eq:NdOmbis}\\
   0 & = \etab + \ze,\label{eq:NdOm}\\
   \Nd\le(\yy\ri) & = -2\xib,\label{eq:Ndyy}\\
   \Nd_4\le(\yy\ri) & = -4\omb.\label{eq:Nd4yy}
 \end{align}
\end{lemma}
\begin{proof}
  Identities~\eqref{eq:xib0} and~\eqref{eq:Nd4Om} are a consequence of the geodesic equation $\D_LL = 0$ and the definition $\el = L$.\\
  
  From the definition of $\eta$ and $\ze$, we have
  \begin{align*}
    \eta_a & = \half \g(\D_3\el,\ea) = - \half \g\le(\el,[\elb,\ea]\ri) - \half \g(\el,\D_a\elb) = - \half \g\le(\el,[\elb,\ea]\ri) + \ze_a. 
  \end{align*}
  From the definition of $\el$ and the relations~\eqref{eq:elu}, we have
  \begin{align*}
    -\half\g\le(\el,[\elb,\ea]\ri) & = \half[\elb,\ea](u) = -\ea\le(1\ri) = 0,
  \end{align*}
  and we deduce~\eqref{eq:NdOmbis}.\\

  From the definition of $\etab$ and $\ze$, we have
  \begin{align*}
    \etab_a & = \half \g(\D_4\elb,\ea) = - \half \g\le(\elb,[\el,\ea]\ri) - \ze_a.
  \end{align*}

  From the definition of $\elb$ and the relations~(\ref{eq:relYelelb}),~\eqref{eq:elu}, we have
  \begin{align*}
    \g\le(\elb,[\el,\ea]\ri) & = \g(\Yb,[\el,\ea]) - \half\yy\g(\el,[\el,\ea]),\\
                             & = - [\el,\ea](\ub) + \half\yy [\el,\ea](u) \\
                             & = 0,
  \end{align*}
  and we deduce~\eqref{eq:NdOm}.\\
  
  From the definition of $\xib$, we have
  \begin{align*}
    \xib_a & = \half\g(\D_3\elb,\ea) = - \half\g(\elb,[\elb,\ea]).
  \end{align*}
  Using the definition of $\elb$ and $\Yb$, we infer
  \begin{align*}
    \xib_a & = - \half \g\le(\Yb-\half\yy\el,[\elb,\ea]\ri) \\
          & = \half [\elb,\ea](\ub) - \quar\yy [\elb,\ea](u) \\
          & = - \half \ea\le(\yy\ri),
  \end{align*}
  and we deduce the identity~\eqref{eq:Ndyy}. \\

  From the definition of $\omb$, the definition of $\el$ and the relations~\eqref{eq:elu}, we have
  \begin{align*}
    \omb & = \quar \g(\D_3\elb,\el) \\
        & = -\quar \g(\elb,[\elb,\el]) \\
        & = -\quar\g\le(\Yb-\half\yy\el,[\elb,\el]\ri) \\
        & = \quar[\elb,\el](\ub) - \frac{1}{8} \yy [\elb,\el](u) \\
        & = -\quar \el(\yy),
  \end{align*}
  and we deduce the identity~\eqref{eq:Nd4yy}. This finishes the proof of the lemma.
\end{proof}

\subsection{Averages and renormalisations}
\begin{lemma}[Average of $\rho$ and $\sigma$]
  The following relations hold in $\MM^\ext$
  \begin{align}
    \Nd_4\rhoo + \frac{3}{2}\trchio\rhoo & = \Err\le(\Nd_4,\rhoo\ri),\label{eq:Nd4rhoo}\\
    \sigmao & = \half \overline{\chih\wedge\chibh},\label{eq:sigmao}
  \end{align}
  where
  \begin{align*}
    \Err\le(\Nd_4,\rhoo\ri) & := -\half \overline{\chibh\cdot\al} - \overline{\ze\cdot\be} -\half \overline{(\trchi-\trchio)(\rho-\rhoo)}.
  \end{align*}
\end{lemma}
\begin{proof}
  Equation~\eqref{eq:Nd4rhoo} follows from taking the average in equation~\eqref{eq:Nd4rho} and using commutation formula~\eqref{eq:commelbov}. Equation~\eqref{eq:sigmao} follows from taking the average in equation~(\ref{eq:Curlze}).
\end{proof}

\begin{lemma}[Average and renormalisation of $\trchi$]
  We have the following transport equations in the $\el$ direction
  \begin{align}
    \Nd_4\le(\trchi-\trchio\ri) + \trchio \le(\trchi-\trchio\ri) & = \Err\le(\Nd_4,\trchi-\trchio\ri),\label{eq:Nd4trchitrchio}\\
    \Nd_4\le(\trchio-\frac{2}{r}\ri) + \half \trchio\le(\trchio-\frac{2}{r}\ri) & = \Err\le(\Nd_4,\trchio-\frac{2}{r}\ri),\label{eq:Nd4trchio}
  \end{align}
  where
  \begin{align*}
    \Err\le(\Nd_4,\trchi-\trchio\ri) & := - |\chih|^2 - \half (\trchi-\trchio)^2 + \overline{|\chih|^2} - \half\overline{(\trchi-\trchio)^2},\\
    \Err\le(\Nd_4,\trchio-\frac{2}{r}\ri) & := - \overline{|\chih|^2} + \half\overline{(\trchi-\trchio)^2}.
  \end{align*}
\end{lemma}
\begin{proof}
  Rewriting~\eqref{eq:Nd4trchi}, we have
  \begin{align}\label{eq:Nd4trchiNEW}
    \Nd_4\trchi + \half (\trchi)^2 & = - |\chih|^2.
  \end{align}
  
  Taking the mean value in~\eqref{eq:Nd4trchiNEW}, using commutation formula~(\ref{eq:commelbov}), we obtain
  \begin{align}\label{eq:Nd4trchioNEW}
    \begin{aligned}
      \Nd_4\trchio + \half \trchio^2 & = - \overline{|\chih|^2} + \half\le(\trchio^2-\overline{(\trchi)^2}\ri) + \overline{(\trchi-\trchio)\trchi} \\
      & = - \overline{|\chih|^2} + \half\overline{(\trchi-\trchio)^2}.
    \end{aligned}
  \end{align}

Combining equations~\eqref{eq:Nd4trchiNEW} and~\eqref{eq:Nd4trchioNEW} we obtain~\eqref{eq:Nd4trchitrchio}. Using equation~\eqref{eq:Nd4trchioNEW} and relation~\eqref{eq:elbr} we obtain~\eqref{eq:Nd4trchio} and it finishes the proof of the lemma.
\end{proof}

\begin{lemma}[Average and renormalisation of $\trchib$]
  We have the following equations in the $\el$ direction
  \begin{align}
    \Nd_4(\trchib-\trchibo) +\half \trchio (\trchib-\trchibo) & = -\half\trchib(\trchi-\trchio) + 2\Divd\ze +2(\rho-\rhoo)\label{eq:Nd4trchibtrchibo} \\
    & \quad + \Err\le(\Nd_4,\trchib-\trchibo\ri)\nonumber\\
    \Nd_4\le(\trchibo+\frac{2}{r}\ri) + \half\trchio \le(\trchibo+\frac{2}{r}\ri) & = 2\rhoo + \Err\le(\Nd_4,\trchibo+\frac{2}{r}\ri), \label{eq:Nd4trchibo}
  \end{align}
  where
  \begin{align*}
    \Err\le(\Nd_4,\trchib-\trchibo\ri) & := -\chih\cdot\chibh +2|\ze|^2 + \overline{\chih\cdot\chibh} - 2 \overline{|\ze|^2} -\half \overline{(\trchi-\trchio)(\trchib-\trchibo)},\\
    \Err\le(\Nd_4,\trchibo+\frac{2}{r}\ri) & := -\overline{\chih\cdot\chibh} +2\overline{|\ze|^2} +\half \overline{(\trchi-\trchio)(\trchib-\trchibo)}.
  \end{align*}
\end{lemma}
\begin{proof}
  Equation~\eqref{eq:Nd4trchib} together with the relations of Lemma~\ref{lem:relgeodnull} rewrites
  \begin{align}
    \label{eq:Nd4trchibNEW}
    \Nd_4\trchib +\half\trchi\trchib & = 2\Divd\ze -\chih\cdot\chibh + 2|\ze|^2 + 2\rho.
  \end{align}

  Taking the average in~\eqref{eq:Nd4trchibNEW}, using commutation formula~(\ref{eq:commelbov}) gives
  \begin{align}
    \label{eq:Nd4trchiboNEW}
    \Nd_4\trchibo + \half\trchio\;\trchibo & = -\overline{\chih\cdot\chibh} +2\overline{|\ze|^2} +2\rhoo +\half \overline{(\trchi-\trchio)(\trchib-\trchibo)}.
  \end{align}

  Combining~\eqref{eq:Nd4trchibNEW} and~\eqref{eq:Nd4trchiboNEW} directly gives~\eqref{eq:Nd4trchibtrchibo}. Combining~\eqref{eq:Nd4trchiboNEW} and formula~\eqref{eq:elbr} gives~\eqref{eq:Nd4trchibo} and this finishes the proof of the lemma.
\end{proof}

\begin{lemma}[Average and renormalisation of $\omb$]
  We have the following equations in the $\el$ direction
  \begin{align}
    \Nd_4(\omb-\ombo) & = \rho-\rhoo +\Err(\Nd_4,\omb-\ombo),\label{eq:Nd4ombombo} \\
    \Nd_4\ombo & = \rhoo + \Err(\Nd_4,\ombo),\label{eq:Nd4ombo}
  \end{align}
  where
  \begin{align*}
    \Err(\Nd_4,\omb-\ombo) & := 3|\ze|^2-3\overline{|\ze|^2} - \overline{(\trchi-\trchio)(\omb-\ombo)},\\
    \Err(\Nd_4,\ombo) & := 3\overline{|\ze|^2}+ \overline{(\trchi-\trchio)(\omb-\ombo)}.
  \end{align*}
\end{lemma}
\begin{proof}
  Using the result of Lemma~\ref{lem:relgeodnull}, equation~\eqref{eq:Nd3om} rewrites
  \begin{align*}
    \Nd_4\omb = 3|\ze|^2 + \rho.
  \end{align*}
  The result of the lemma then follows from taking the average in the above equation. The details are left to the reader.
\end{proof}

\subsection{The null coefficient $\ze$ and renormalisations}
\paragraph{The mass aspect function $\mu$}
Equation~\eqref{eq:Nd4ze} together with the relations of Lemma~\ref{lem:relgeodnull} rewrites
\begin{align}\label{eq:Nd4zeNEW}
  \Nd_4\ze + \trchi\ze & = -2\chih\cdot\ze -\beta.
\end{align}

We define the mass aspect function $\mu$ by
\begin{align}\label{eq:defmu}
  \mu & := \Divd\ze + \rho
\end{align}


Using the relations of Lemma~\ref{lem:relgeodnull} and~\eqref{eq:defmu} and~(\ref{eq:Curlze}), we have
\begin{align}\label{eq:Dd1etab}
  \begin{aligned}
    \Dd_1(\ze) & = \le(-\rho +\mu,\sigmac\ri), \\
    & = \le(-\rho+\rhoo + \mu-\muo, \sigmac-\overline{\sigmac}\ri).
    \end{aligned}
\end{align}
where for a $S$-tangent $1$-tensor $U$, we have $\Dd_1U := \le(\Divd U,\Curld U\ri)$, and where 
\begin{align*}
  \sigmac & := \sigma -\half \chih\wedge\chibh.
\end{align*}

We have the following transport equation in the $\el$ direction for $\mu$.
\begin{lemma}
  We have
  \begin{align}\label{eq:Nd4mu}
    \Nd_4\mu + \frac{3}{2}\trchi \mu & = \Err(\Nd_4,\mu),
  \end{align}
  where
  \begin{align*}
    \Err(\Nd_4,\mu) & := -2\Nd\chih\cdot\ze -3\chih\cdot\Nd\ze - \chi\cdot\ze\cdot\ze + \trchi|\ze|^2 - 2\ze\cdot\be -\Nd\trchi\cdot\ze   -\half\chibh\cdot\alpha.
  \end{align*}
\end{lemma}
\begin{proof}
  Commuting the transport equation~\eqref{eq:Nd4zeNEW} with $\Divd$, using commutation formula~\eqref{eq:commNd4Nd}, we have
  \begin{align}\label{eq:Nd4divdze}
    \Nd_4\Divd\ze + \frac{3}{2}\trchi\Divd\ze & = -\Divd\beta +\Err_1,
  \end{align}
  where
  \begin{align*}
    \Err_1 & := \le([\Nd_4,\Divd]\ze +\half\trchi\Divd\ze\ri) -\Nd\trchi\cdot\ze -2\Nd\chih\cdot\ze -2\chih\cdot\Nd\ze \\
           & = \le(-\chih\cdot\Nd\ze - \chi\cdot\ze\cdot\ze + \trchi|\ze|^2 - \ze\cdot\be\ri) -\Nd\trchi\cdot\ze -2\Nd\chih\cdot\ze -2\chih\cdot\Nd\ze.
  \end{align*}

  From Bianchi identity~\eqref{eq:Nd4rho}, we have
  \begin{align}\label{eq:Nd4rhobis}
    \Nd_4\rho + \frac{3}{2}\rho & = \Divd\be +\Err_2,
  \end{align}
  where
  \begin{align*}
    \Err_2 & := -\half\chibh\cdot\alpha-\ze\be,
  \end{align*}
  and summing~\eqref{eq:Nd4divdze} and~\eqref{eq:Nd4rhobis} gives the desired formula~\eqref{eq:Nd4mu}.
\end{proof}

We have the following transport equation for $\mu-\muo$.
\begin{lemma}[Average and renormalisation of $\mu$]
  We have
  \begin{align}\label{eq:muorhoo}
    \muo = \rhoo,
  \end{align}
  and
  \begin{align}\label{eq:Nd4mumuo}
    \Nd_4(\mu-\muo) + \frac{3}{2}\trchio(\mu-\muo) & = \Err\le(\Nd_4,\mu-\muo\ri),
  \end{align}
  where
  \begin{align*}
    \Err\le(\Nd_4,\mu-\muo\ri) & := \Err(\Nd_4,\mu) - \overline{\Err(\Nd_4,\mu)} +\half \overline{(\trchi-\trchio)(\mu-\muo)} + 3/2 (\trchi-\trchio)\mu.
  \end{align*}
\end{lemma}
\begin{proof}
  Equation~\eqref{eq:muorhoo} follows from taking the average in the definition~(\ref{eq:defmu}) of $\mu$.\\
  Taking the average in equation~\eqref{eq:Nd4mu} using formula~\eqref{eq:commelbov}, we have
  \begin{align}\label{eq:Nd4muo}
    \Nd_4\muo + \frac{3}{2}\trchio\muo & = -\half\overline{(\trchi-\trchio)(\mu-\muo)} + \overline{\Err\le(\Nd_4,\mu\ri)}.
  \end{align}
  Combining~\eqref{eq:Nd4mu} and~\eqref{eq:Nd4muo} then gives the desired~\eqref{eq:Nd4mumuo}.
\end{proof}

\paragraph{Renormalisation of $\Nd_3\ze$}
We define the auxiliary coefficients $\ombr$ and $\ombs$ to be the solution of the following transport equations in the $\el$ direction
\begin{align}\label{eq:Nd4ombrs}
  \begin{aligned}
    \Nd_4(r\ombr) & = r(\rho-\rhoo), & \Nd_4(r\ombs) & = r(\sigma-\sigmao), \\
    \ombr|_{\CCba} & = 0, & \ombs|_{\CCba} & = 0.
  \end{aligned}
\end{align}

We define $\io$ to be the following $S$-tangent tensor 
\begin{align}\label{eq:defio}
  \io & := \Nd\ombr - \Nds\ombs + \beb = \Dd^\ast_1(\ombr,\ombs) + \beb.
\end{align}

The tensor $\io$ satisfies the following transport equation.
\begin{lemma}
  We have
  \begin{align}\label{eq:Nd4io}
    \Nd_4\le(r^2\io\ri) & = \Err\le(\Nd_4,\io\ri),
  \end{align}
  where
  \begin{align*}
    \Err\le(\Nd_4,\io\ri) & := -\half\le(\trchi-\trchio\ri)r\Nd(r\ombr) - \chih\cdot (r\Nd)(r\ombr) \\
                          & \quad +\half\le(\trchi-\trchio\ri)(r\Nds)(r\ombs) - \chih\cdot (r\Nds)(r\ombs) \\
                          & \quad - r^2(\trchi-\trchio)\beb -2r^2\chibh\cdot\be-3r^2\le(\ze\rho-\dual\ze\sigma\ri).
  \end{align*}
\end{lemma}
\begin{proof}
  Commuting the transport equations~\eqref{eq:Nd4ombrs} by respectively $r\Nd$ and $r\Nds$, using commutation formula~(\ref{eq:commNd4Nd}), we have
  \begin{align}\label{eq:rNdNd3omoms}
    \begin{aligned}
      \Nd_4(r^2\Nd\ombr) & = r^2\Nd\rho + [\Nd_4,r\Nd](r\ombr) \\
      & = r^2\Nd\rho -\half\le(\trchi-\trchio\ri)(r\Nd)(r\ombr) - \chih\cdot (r\Nd)(r\ombr), \\
      \Nd_4(r\Nds\ombs) & = r^2\Nds\sigma + \dual\le([\Nd_4,r\Nd](r\omb)\ri) \\
      & = r^2\Nds\sigma -\half\le(\trchi-\trchio\ri)r\Nds(r\ombs) - \dual\le(\chih\cdot (r\Nd)(r\ombs)\ri) \\
      & = r^2\Nds\sigma -\half\le(\trchi-\trchio\ri)r\Nds(r\ombs) + \chih\cdot (r\Nds)(r\ombs),
    \end{aligned}
  \end{align}
  where we used standard Hodge dual computations. From Bianchi identity~\eqref{eq:Nd4beb} and~\eqref{eq:elbr} we have
  \begin{align}\label{eq:rNdrhosigma}
    \begin{aligned}
      -r^2\Nd\rho + r^2\Nds\sigma & = r^2\Nd_4\beb+r^2\trchi\beb + 2r^2\chibh\cdot\be+3r^2\le(\ze\rho-\dual\ze\sigma\ri) \\
      & = \Nd_4(r^2\beb) + r^2(\trchi-\trchio)\beb + 2r^2\chibh\cdot\be+3r^2\le(\ze\rho-\dual\ze\sigma\ri).
    \end{aligned}
  \end{align}
  
  Combining~\eqref{eq:rNdNd3omoms} and~\eqref{eq:rNdrhosigma}, then directly gives the desired formula~\eqref{eq:Nd4io}.  
\end{proof}

We define $\Pso$ to be the following renormalisation of $\Nd_3\ze$
\begin{align}\label{eq:defPso}
  \Pso & := \Nd_3\ze + \Nd\ombr + \Nd\ombs.
\end{align}

We have the following transport equation for $\Pso$ in the $\el$ direction.
\begin{lemma}
  We have
  \begin{align}\label{eq:Nd4Pso}
    \Nd_4\le(r^2\Pso\ri) & = r^2\trchib \be + \half r^2\trchi\trchib\ze + \Err\le(\Nd_4,\Pso\ri),
  \end{align}
  where
  \begin{align*}
    \Err\le(\Nd_4,\Pso\ri) & := \le(r^2\trchio-r^2\trchi\ri)\Nd_3\ze \\
                           & \quad -2r^2\omb\be-r^2\xib\cdot\al -3r^2(\ze\rho+\dual\ze\sigma) \\
      & \quad  -2r^2(\Divd\ze)\ze-2r^2\omb\trchi\ze+r^2\chih\cdot\chibh\ze -2r^2|\ze|^2\ze -2r^2\rho\ze \\
                           & \quad -2r^2\Nd_3(\chih\cdot\ze) -2r^2\omb\Nd_4\ze -2r^2\ze\cdot\Nd\ze - 2r^2\dual\ze\sigma \\
                           & \quad +\half(\trchio-\trchi)\Nd(r\ombr) -\chih\cdot\Nd(r\ombr) \\
           & \quad + \half(\trchio-\trchi)\Nds(r\ombs) + \chih\cdot\Nds(r\ombs).
  \end{align*}
\end{lemma}
\begin{proof}
  Using formula~\eqref{eq:elbr}, commuting the transport equation~\eqref{eq:Nd4zeNEW} with $\Nd_3$, using commutation formula~(\ref{eq:commNd3Nd4}), Bianchi equation~\eqref{eq:Nd3be} for $\Nd_3\be$ and the equation~\eqref{eq:Nd3trchi} for $\Nd_3\trchi$ together with the relations from Lemma~\ref{lem:relgeodnull} gives
  \begin{align}\label{eq:Nd4Nd3zeproof}
    \begin{aligned}
      \Nd_4\le(r^2\Nd_3\ze\ri) & = r^2\le(\Nd_4\Nd_3\ze + \trchi \Nd_3\ze\ri) + \Err_1 \\
      & = -r^2\Nd_3\be - r^2\Nd_3(\trchi)\ze -2r^2\Nd_3(\chih\cdot\ze) + r^2[\Nd_4,\Nd_3]\ze + \Err_1 \\
      & = r^2\trchib\be -r^2\Nd\rho-r^2\Nds\sigma + \half r^2\trchi\trchib\ze + \Err_1 + \Err_2,
    \end{aligned}
  \end{align}
  where
  \begin{align*}
    \begin{aligned}
      \Err_1 & := \le(r^2\trchio-r^2\trchi\ri)\Nd_3\ze, \\ \\
      \Err_2 & := -2r^2\omb\be-r^2\xib\cdot\al -3r^2(\ze\rho+\dual\ze\sigma) \\
      & \quad  -2r^2(\Divd\ze)\ze-2r^2\omb\trchi\ze+r^2\chih\cdot\chibh\ze -2r^2|\ze|^2\ze -2r^2\rho\ze \\
      & \quad -2r^2\Nd_3(\chih\cdot\ze) -2r^2\omb\Nd_4\ze -2r^2\ze\cdot\Nd\ze - 2r^2\dual\ze\sigma.
    \end{aligned}
  \end{align*}

  From the transport equation~(\ref{eq:Nd4ombrs}) for the auxiliary coefficients $\ombr$ and $\ombs$, using commutation formula~\eqref{eq:commNd4Nd}, we have
  \begin{align}\label{eq:NdrNdsig2}
    \begin{aligned}
    -r^2\Nd\rho-r^2\Nds\sigma & = -r\Nd\Nd_4(r\ombr) - r\Nds\Nd_4(r\ombs) \\
    & = -\Nd_4(r^2\Nd\ombr)-\Nd_4(r^2\Nd\ombs) + \Err_3,
    \end{aligned}
  \end{align}
  where
  \begin{align*}
    \Err_3 & := -[r\Nd,\Nd_4](r\ombr) -[r\Nds,\Nd_4]\ombs \\
           & = \half(\trchio-\trchi)\Nd(r\ombr) -\chih\cdot\Nd(r\ombr) \\
           & \quad + \half(\trchio-\trchi)\Nds(r\ombs) + \chih\cdot\Nds(r\ombs).
  \end{align*}

  Combining~\eqref{eq:Nd4Nd3zeproof} and~\eqref{eq:NdrNdsig2} gives the desired formula.
\end{proof}

\subsection{The null coefficient $\protect\omb$ and renormalisations}
Using equation~(\ref{eq:Nd3om}) and the relations from Lemma~\ref{lem:relgeodnull}, we have
\begin{align}\label{eq:Nd4ombNEW}
  \Nd_4\omb = 3|\ze|^2 + \rho. 
\end{align}

Let $\ombd$ be the \emph{ad hoc} dual of $\omb$, which we define as the solution of the following transport equation
\begin{align}\label{eq:Nd4oms}
  \begin{aligned}
  \Nd_4\ombd & := \sigma,\\
  \ombd|_{\CCba} & := 0.
  \end{aligned}
\end{align}

We define $\iob$ to be the $S$-tangent tensor given by
\begin{align}
  \label{eq:defiob}
  \iob & := \Nd\omb - \Nds\ombd +\beb.
\end{align}

The tensor $\iob$ satisfies the following transport equation.
\begin{lemma}
  We have
  \begin{align}\label{eq:Nd4iob}
    \Nd_4\le(r\iob\ri) & = -\half r\trchio\beb + \Err\le(\Nd_4,\iob\ri),
  \end{align}
  where
  \begin{align*}
    \Err\le(\Nd_4,\iob\ri) & := 6\ze\cdot(r\Nd)\ze \\
                          & \quad -\half\le(\trchi-\trchio\ri)r\Nd\omb - \chih\cdot (r\Nd)\omb \\
                          & \quad +\half\le(\trchi-\trchio\ri)(r\Nds)\ombd - \chih\cdot (r\Nds)\ombd \\
                          & \quad - r(\trchi-\trchio)\beb -2r\chibh\cdot\be-3r\le(\ze\rho-\dual\ze\sigma\ri).
  \end{align*}
\end{lemma}
\begin{proof}
  Commuting the transport equations~\eqref{eq:Nd4ombNEW} and~\eqref{eq:Nd4oms} by respectively $r\Nd$ and $r\Nds$, using commutation formula~(\ref{eq:commNd4Nd}), we have
  \begin{align}\label{eq:rNdNd3ombombs}
    \begin{aligned}
      \Nd_4(r\Nd\omb) & = 6\ze\cdot(r\Nd)\ze + (r\Nd)\rho + [\Nd_4,r\Nd]\omb \\
      & = 6\ze\cdot(r\Nd)\ze + (r\Nd)\rho \\
      & \quad -\half\le(\trchi-\trchio\ri)(r\Nd)\omb - \chih\cdot (r\Nd)\omb, \\
      \Nd_4(r\Nds\ombd) & = (r\Nds)\sigma + \dual\le([\Nd_4,r\Nd]\ombd\ri) \\
      & = (r\Nds)\sigma -\half\le(\trchi-\trchio\ri)r\Nds\ombd - \dual\le(\chih\cdot (r\Nd)\ombd\ri) \\
      & = (r\Nds)\sigma -\half\le(\trchi-\trchio\ri)r\Nds\ombd + \chih\cdot (r\Nds)\ombd,
    \end{aligned}
  \end{align}
  where we used standard Hodge dual computations. From Bianchi identity~\eqref{eq:Nd4beb} and~\eqref{eq:elbr} we have
  \begin{align}\label{eq:rNdrhosigmabis}
    \begin{aligned}
      -(r\Nd)\rho + (r\Nds)\sigma & = r\Nd_4\beb+r\trchib\beb + 2r\chibh\cdot\be+3r\le(\ze\rho-\dual\ze\sigma\ri) \\
      & = \Nd_4(r\beb) +\half r\trchio\beb + r(\trchi-\trchio)\beb + 2r\chibh\cdot\be+3r\le(\ze\rho-\dual\ze\sigma\ri).
    \end{aligned}
  \end{align}
  
  Combining~\eqref{eq:rNdNd3ombombs} and~\eqref{eq:rNdrhosigmabis}, then directly gives the desired formula~\eqref{eq:Nd4iob}.  
\end{proof}


\section{Null decomposition of the canonical foliation in $\protect\CCba$}\label{sec:nulldecCCba}

Additionally to the relations of Section~\ref{sec:nulldecompMMext}, we have the following relations on $\CCba$.
\begin{lemma}\label{lem:relCCba}
  On $\CCba$, the following relations hold
  \begin{align}
    \yy & = 0,\label{eq:yyCCba} \\
    \xib & = 0, \label{eq:xibCCba}
  \end{align}
\end{lemma}
\begin{proof}
  Since $\CCba$ is null, the gradient $\Yb$ of $\ub$ is null and $\yy = 0$ on $\CCba$. This proves~\eqref{eq:yyCCba}.\\

  We also deduce on $\CCba$ that
  \begin{align*}
    \xib_a & = \half \g(\D_3\elb,\ea) = \half \g(\elb,[\ea,\elb]) = -\half\g(\D\ub,[\ea,\elb]) = -\half [\ea,\elb](\ub) = 0. 
  \end{align*}
  This finishes the proof of the lemma.
\end{proof}

\subsection{Averages and renormalisations}
\begin{lemma}[Average of $\rho$ and $\sigma$]
  The following relations hold on $\CCba$
  \begin{align}
    \Nd_3\rhoo + \frac{3}{2}\trchibo\rhoo & = \Err\le(\Nd_3,\rhoo\ri), \label{eq:Nd3rhoo}\\
    \sigmao & = \half \overline{\chih\wedge\chibh},
  \end{align}
  where
  \begin{align*}
    \Err\le(\Nd_3,\rhoo\ri) & := -\overline{\half\chih\cdot\alb} - \overline{\ze\cdot\beb} -\half \overline{\le(\trchib-\trchibo\ri)\le(\rho-\rhoo\ri)}. 
  \end{align*}
\end{lemma}
\begin{proof}
  Equation~\eqref{eq:Nd3rhoo} follows from taking the average in Bianchi equation~\eqref{eq:Nd3rho}, using the relations of Lemmas~\ref{lem:relgeodnull} and \ref{lem:relCCba}, and using formula~\eqref{eq:commelbov}.
\end{proof}

\begin{lemma}[Average and renormalisation of $\trchib$]
  We have the following transport equation along $\CCba$
  \begin{align}
    \Nd_3\le(\trchib-\trchibo\ri) + \trchibo \le(\trchib-\trchibo\ri) & = -2(\omb-\ombo)\trchib + \Err\le(\Nd_3,\trchib-\trchibo\ri),\label{eq:Nd3trchibtrchibo}\\
    \Nd_3\le(\trchibo+\frac{2}{r}\ri) +\half \trchibo \le(\trchibo+\frac{2}{r}\ri) & = \Err\le(\Nd_3,\trchibo+\frac{2}{r}\ri),\label{eq:Nd3trchibo}
  \end{align}
  where\footnote{On $\CCba$ we have by definition of the canonical foliation $\ombo = 0$. For completeness, we keep track of the factor $\ombo$ in the formulas of this section.}
  \begin{align*}
    \Err\le(\Nd_3,\trchib-\trchibo\ri) & := -|\chibh|^2 + \overline{|\chibh|^2} +2\overline{(\omb-\ombo)(\trchib-\trchibo)}-\half \le(\trchib-\trchibo\ri)^2-\half\overline{(\trchib-\trchibo)^2},\\
    \Err\le(\Nd_3,\trchib+\frac{2}{r}\ri) & := -2\overline{(\omb-\ombo)(\trchib-\trchibo)} - \overline{|\chibh|^2} + \half \overline{(\trchib-\trchibo)^2}.
  \end{align*}
\end{lemma}
\begin{proof}
  Using the relations from Lemmas~\ref{lem:relgeodnull} and~\ref{lem:relCCba}, equation~\eqref{eq:Nd3trchib} rewrites
  \begin{align}\label{eq:Nd3trchibNEW}
    \Nd_3\trchib + \half(\trchib)^2 & = -2(\omb-\ombo)\trchib - |\chibh|^2,
  \end{align}
  where we used that $\ombo = 0$ for the canonical foliation on $\CCba$ (see Definition~\ref{def:can}).
  Taking the average in~\eqref{eq:Nd3trchibNEW}, using formula~\eqref{eq:commelbov}, we have
  \begin{align}\label{eq:Nd3trchiboo}
    \Nd_3\trchibo + \half (\trchibo)^2 & = -2\overline{(\omb-\ombo)(\trchib-\trchibo)} - \overline{|\chibh|^2} + \half \overline{(\trchib-\trchibo)^2}.
  \end{align}

  Combining~\eqref{eq:Nd3trchibNEW} and~\eqref{eq:Nd3trchiboo}, we obtain~\eqref{eq:Nd3trchibtrchibo}. Combining~\eqref{eq:Nd3trchiboo} and formula~\eqref{eq:elbr} gives~\eqref{eq:Nd3trchibo}.  
\end{proof}

\begin{lemma}[Average and renormalisation of $\trchi$]
  We have the following transport equation along $\CCba$
  \begin{align}
    \Nd_3\le(\trchi-\trchio\ri) +\half\trchibo(\trchi-\trchio) & = 2(\omb-\ombo)\trchi + \Err\le(\Nd_3,\trchi-\trchio\ri),\label{eq:Nd3trchitrchio}\\
    \Nd_3\le(\trchio-\frac{2}{r}\ri) + \half\trchibo\le(\trchio-\frac{2}{r}\ri)& = 2\rhoo + \Err\le(\Nd_3,\trchio-\frac{2}{r}\ri),\label{eq:Nd3trchio}
  \end{align}
  where
  \begin{align*}
    \Err\le(\Nd_3,\trchi-\trchio\ri) & := -\chih\cdot\chibh+2|\ze|^2 + \overline{\chih\cdot\chibh} -2\overline{|\ze|^2} -2\overline{(\omb-\ombo)(\trchi-\trchio)} \\
    & \quad - \half(\trchib-\trchibo)(\trchi-\trchio) - \half\overline{(\trchib-\trchibo)(\trchi-\trchio)},\\
    \Err\le(\Nd_3,\trchio-\frac{2}{r}\ri) & := 2\overline{(\omb-\ombo)(\trchi-\trchio)} - \overline{\chih\cdot\chibh} + 2\overline{|\ze|^2} + \half\overline{(\trchib-\trchibo)(\trchi-\trchio)}. 
  \end{align*}
\end{lemma}
\begin{proof}
  Rewriting equation~\eqref{eq:Nd3trchi}, using the relations of Lemmas~\ref{lem:relgeodnull} and~\ref{lem:relCCba}, and the Definition~\ref{def:can} of the canonical foliation on $\CCba$, we have
  \begin{align}
    \label{eq:Nd3trchiNEW}
    \Nd_3\trchi +\half\trchib\trchi & = 2\rhoo + 2(\omb-\ombo)\trchi -\chih\cdot\chibh+2|\ze|^2.
  \end{align}
  Taking the average in~\eqref{eq:Nd3trchiNEW} gives
  \begin{align}
    \label{eq:Nd3trchioo}
    \Nd_3\trchio + \half\trchibo\trchio & = 2\rhoo + 2\overline{(\omb-\ombo)(\trchi-\trchio)} - \overline{\chih\cdot\chibh} + 2\overline{|\ze|^2} + \half\overline{(\trchib-\trchibo)(\trchi-\trchio)}.
  \end{align}
  Combining~\eqref{eq:Nd3trchiNEW} and~\eqref{eq:Nd3trchioo} gives~\eqref{eq:Nd3trchitrchio}. Combining~\eqref{eq:Nd3trchioo} and formula~\eqref{eq:elbr} gives~\eqref{eq:Nd3trchio}.
\end{proof}

\section{Approximate exterior Killing vectorfields $\TE,\SE,\KE$}\label{sec:defKillingext}
Associated to the foliation $S_{u,\ub}$ and the null pair $(\el,\elb)$, we define the following \emph{approximate conformal Killing fields} in the exterior region $\MM^\ext$
\begin{align}\label{eq:defTESEKE}
  \begin{aligned}
  \TE & := \half\le(\el +\elb\ri), & \SE & := \half\le(\ub \el + u \elb\ri), \\
  \KE & := \half\le(\ub^2\el + u^2\elb\ri).
  \end{aligned}
\end{align}

We remark that we have
\begin{align}\label{eq:relKESETE}
  \KE & = (u+\ub)\SE - \quar (\ub u)\TE.
\end{align}

\section{Approximate exterior Killing rotations $\OOE$}\label{sec:defrotext}
\subsection{Approximate exterior Killing rotations on $\protect\CCba$}




For a fixed centred conformal isomorphism $\Phi$ of $S^\ast$ and for the associated Cartesian functions $x^i$ (see the definitions of Section~\ref{sec:defUnifandharmo}), we define the \emph{approximate Killing exterior rotations} on $S^\ast$ by
\begin{align}
  \label{eq:defOOESast}
  ^{(1)}\OOE & := x^2 \Nd x^3 - x^3\Nd x^2, & ^{(2)}\OOE & := x^3 \Nd x^1 - x^1\Nd x^3, &  ^{(3)}\OOE & := x^1 \Nd x^2 - x^2\Nd x^1,
\end{align}
where $x^i$ are the Cartesian functions on $S^\ast$ defined in Section~\ref{sec:defUnifandharmo}.\\

We extend the rotations $\OOE$ by Lie transport along $\CCba$, \emph{i.e.}
\begin{align}\label{eq:Lie3OOE}
  \le[\elb, ^{(\ell)}\OOE\ri] = 0,
\end{align}
for all $\ell=1,2,3$.

\begin{lemma}\label{lem:OOEitangent}
  The vectorfields $\OOEi$ are tangent to the $2$-spheres $S_{u,\ub}$ of the canonical foliation on $\CCba$.
\end{lemma}
\begin{proof}
  From the definition~\eqref{eq:Lie3OOE} and relations~\eqref{eq:elu}, we have
  \begin{align*}
    \elb\le(\OOEi(\ub)\ri) & = \OOEi\le(\elb(\ub)\ri) = 0, && \text{and} && \elb\le(\OOEi(u)\ri) = \OOEi\le(\elb(u)\ri) = 0.
  \end{align*}
  From the definition of $\OOEi$ on $S^\ast$, we have $\OOEi(u)|_{S^\ast} = \OOEi(\ub)|_{S^\ast} = 0$, and the result follows.
\end{proof}

\begin{remark}
  The exterior rotation vectorfields have to be \emph{exactly tangent} to the $2$-spheres of the foliation $S_{u,\ub}$ since they are used to estimate the tangential derivatives of the curvature. This is the reason for the definition of $\OOE$ by Lie transport. In contrast, the interior rotation $\OOI$ are not used to estimate the curvature and thus can simply be defined by parallel transport. See Section~\ref{sec:defKillingint}. 
\end{remark}





In the rest of this section, we simply call $\OOO$ a rotation vectorfield $\OOEi$.\\

Projecting equation~\eqref{eq:Lie3OOE} on $S$ provides in particular
\begin{align}\label{eq:Nd3OOE}
  \Nd_3\OOO_a & = \chib_{ab}\OOO^b.
\end{align}

Commuting equation~\eqref{eq:Nd3OOE} with $\Nd$ using commutation formula~\eqref{eq:commNd3Nd}, we further have the following transport equations,
\begin{align}\label{eq:Nd3NdO}
  \begin{aligned}
    \Nd_3\Nd_a \OOO_b  & = \half\Nd_a(\trchib) \OOO_b + \Nd_a\chibh_{bc} \OOO^c + \chibh_{bc}\Nd_a\OOO^c  -(\chibh\cdot \Nd)_a \OOO_b  \\
    & \quad + (\ze\cdot \OOO)\chib_{ab} -(\chib\cdot \OOO)_a\ze_b - \dual\beb_a\dual \OOO_b,
  \end{aligned}
\end{align}
and the symmetrised version
\begin{align}\label{eq:Nd3H}
  \begin{aligned}
    \Nd_3\Hrot_{ab} & = \le(\Nd(\trchib)\otimes \OOO\ri)_{ab} + \le((\Nd\otimes\chibh)\cdot \OOO\ri)_{ab} + \chibh_{bc}\Nd_a\OOO^c \\
    & \quad + \chibh_{ac}\Nd_b\OOO^c - \le((\chibh\cdot \Nd) \otimes \OOO\ri)_{ab} + 2(\ze\cdot \OOO)\chib_{ab} \\
      &\quad - \le(\chib\cdot \OOO \otimes \ze\ri)_{ab} - (\dual\beb\otimes\dual \OOO)_{ab},
    \end{aligned}
\end{align}
where
\begin{align}\label{eq:defH}
  \Hrot_{ab} := \Nd_a\OOO_b + \Nd_b\OOO_a.
\end{align}

We also define the following $S$-tangent $3$-tensor $\POE$\footnote{The definition of $\POE$ is motivated by the fact that $\POE = 0$ in the Minkowskian case. See Lemma~\ref{lem:D2Nd2X}.}
\begin{align}
  \label{eq:defPOE}
  \POE_{abc} & := \Nd^2_{a,b}\OOO_c -r^{-2}\OOO_b\gd_{ac} +r^{-2}\OOO_c\gd_{ab}.
\end{align}

\begin{lemma}
  We have the following transport equation for $r\POE$
  \begin{align}\label{eq:Nd3POE}
    \begin{aligned}
      \Nd_3\le(r\POE\ri)_{abc} & = \half(r\Nd)_a\Nd_b(\trchib) \OOO_c + (r\Nd)_a\Nd_b\chibh_{cd} \OOO_d + (r\Nd)_a\chibh_{cd}\Nd_b\OOO_d  -(r\Nd)_a\chibh_{bd}\Nd_d\OOO_c\\
      & \quad + (r\Nd)_a\ze_d\OOO_d\chib_{bc} -(r\Nd)_a\chib_{bd}\OOO_d\ze_c - (r\Nd)_a\dual\beb_b\dual \OOO_c + \half\Nd_b(\trchib) (r\Nd)_a\OOO_c \\
      & \quad + \Nd_b\chibh_{cd} (r\Nd)_a\OOO_d + \chibh_{cd}(r\Nd)_a\Nd_b\OOO_d  -\chibh_{bd}(r\Nd)_a\Nd_d\OOO_c + \ze_d(r\Nd)_a\OOO_d\chib_{bc} \\
      & \quad -\chib_{bd}(r\Nd)_a\OOO_d\ze_b - \dual\beb_b(r\Nd)_a\dual \OOO_c + \ze_d\OOO_d (r\Nd)_a\chib_{bc}-\chib_{bd}\OOO_d(r\Nd)_a\ze_b \\
      & \quad + \half (\trchibo-\trchib)(r\Nd)_a\Nd_b\OOO_c -r\chibh_{ad}\Nd_d\Nd_b\OOO_c \\
      & \quad + \le(\chib_{ba}\ze_d-\chib_{da}\ze_b-\iin_{bd}\dual\beb_a\ri)\Nd_d\OOO_c \\
      & \quad + \le(\chib_{ca}\ze_d-\chib_{da}\ze_c-\iin_{cd}\dual\beb_a\ri)\Nd_b\OOO_d \\
      & \quad -\half r^{-1}\le(\trchib-\trchibo\ri)\OOO_b\gd_{ac} - r^{-1}\chibh_{bd}\OOO_d\gd_{ac} \\
      & \quad +\half r^{-1}\le(\trchib-\trchibo\ri)\OOO_c\gd_{ab} + r^{-1}\chibh_{cd}\OOO_d\gd_{ab}.
    \end{aligned}
  \end{align}
\end{lemma}
\begin{proof}
  Using commutation formula~\eqref{eq:commNd3Nd} and the relations of Lemmas~\ref{lem:relgeodnull} and~\ref{lem:relCCba}, we have
  \begin{align}\label{eq:Nd3Nd2OPOE}
    \begin{aligned}
      \Nd_3r\Nd_a\Nd_b\OOO_c & = (r\Nd)_a\Nd_3\Nd_b\OOO_c + \half (\trchibo-\trchib)(r\Nd)_a\Nd_b\OOO_c -r\chibh_{ad}\Nd_d\Nd_b\OOO_c \\
      & \quad + \le(\chib_{ba}\ze_d-\chib_{da}\ze_b-\iin_{bd}\dual\beb_a\ri)\Nd_d\OOO_c + \le(\chib_{ca}\ze_d-\chib_{da}\ze_c-\iin_{cd}\dual\beb_a\ri)\Nd_b\OOO_d.
    \end{aligned}
  \end{align}

  Using equation~\eqref{eq:Nd3NdO}, we have
  \begin{align}\label{eq:NdaNd3NdOOOPOE}
    \begin{aligned}
      & (r\Nd)_a\Nd_3\Nd_b\OOO_c \\
       = & \half (r\Nd)_a\Nd_b(\trchib) \OOO_c + (r\Nd)_a\Nd_b\chibh_{cd} \OOO_d + (r\Nd)_a\chibh_{cd}\Nd_b\OOO_d  -(r\Nd)_a\chibh_{bd}\Nd_d\OOO_c\\
      & + (r\Nd)_a\ze_d\OOO_d\chib_{bc} -(r\Nd)_a\chib_{bd}\OOO_d\ze_c - (r\Nd)_a\dual\beb_b\dual \OOO_c \\
      &  + \half \Nd_b(\trchib) (r\Nd)_a\OOO_c + \Nd_b\chibh_{cd} (r\Nd)_a\OOO_d + \chibh_{cd}(r\Nd)_a\Nd_b\OOO_d  -\chibh_{bd}(r\Nd)_a\Nd_d\OOO_c  \\
      &  + \ze_d(r\Nd)_a\OOO_d\chib_{bc} -\chib_{bd}(r\Nd)_a\OOO_d\ze_b - \dual\beb_b(r\Nd)_a\dual \OOO_c \\
      & + \ze_d\OOO_d (r\Nd)_a\chib_{bc}-\chib_{bd}\OOO_d(r\Nd)_a\ze_b.
    \end{aligned}
  \end{align}

  Using formula~\eqref{eq:elbr}, we have
  \begin{align}\label{eq:Nd3rOPOE}
    \begin{aligned}
      -\Nd_3r^{-1}\OOO_b\gd_{ac} & = -\half r^{-1}\le(\trchib-\trchibo\ri)\OOO_b\gd_{ac} - r^{-1}\chibh_{bd}\OOO_d\gd_{ac},\\
      \Nd_3r^{-1}\OOO_c\gd_{ab} & = \half r^{-1}\le(\trchib-\trchibo\ri)\OOO_c\gd_{ab} + r^{-1}\chibh_{cd}\OOO_d\gd_{ab}.
    \end{aligned}
  \end{align}

  Combining~\eqref{eq:Nd3Nd2OPOE},~\eqref{eq:NdaNd3NdOOOPOE} and~\eqref{eq:Nd3rOPOE} gives the desired formula and finishes the proof of the lemma.
\end{proof}

\subsection{Approximate exterior Killing rotations in $\MM^\ext$}
We define the vectorfields $\OOE$ in $\MM^\ext$ by Lie transport from $\CCba$ along $\el$, \emph{i.e.}
\begin{align}\label{eq:Lie4OOE}
  \begin{aligned}
  \le[\el,{^{(\ell)}\OOE}\ri] & = 0,
  \end{aligned}
\end{align}
Arguing as in Lemma~\ref{lem:OOEitangent}, the vectorfields $\OOEi$ are tangent to the $2$-spheres $S_{u,\ub}$.\\

In the rest of this section, we simply call $\OOO$ a rotation vectorfield $\OOEi$.\\

From~\eqref{eq:Lie4OOE}, Lemma~\ref{lem:OOEitangent} and the relations~\eqref{eq:Riccirel}, we have the following relations
\begin{align}\label{eq:relOOEnull}
  \begin{aligned}
    \g\le(\D_4\OOO,\el \ri) & = 0,\\
    \g\le(\D_4\OOO,\elb \ri) & = 2\ze_a\OOO^a,\\
    \g\le(\D_4\OOO, \ea \ri) & = \chi_{ab}\OOO^b,\\
    \g\le(\D_3\OOO, \elb \ri) & = -2\xib_a\OOO^a,\\
    \g\le(\D_3\OOO, \el\ri) & = -2\ze_a\OOO^a,\\
    \g\le(\D_a\OOO,\el \ri) & = - \chi_{ab}\OOO^b,\\
    \g\le(\D_a\OOO,\elb\ri) & = -\chib_{ab}\OOO^b,\\
    \g\le(\D_3\OOO,\ea\ri) & = \chib_{ab}\OOO^b + Y_a,
  \end{aligned}
\end{align}
where
\begin{align}\label{eq:defY}
  Y_a & := \g\le(\D_3\OOO,\ea\ri) - \chib_{ab}\OOO^b.
\end{align}

\begin{remark}
  From the definition of $\OOE$ on $\CCba$ and the definition~\eqref{eq:defY} of $Y$, we have
  \begin{align}\label{eq:YCCba}
    Y|_{\CCba} & = 0.
  \end{align}
\end{remark}

We deduce the following expression for the deformation tensor $\pi$ of $\OOO$
\begin{align}\label{eq:relpiOOEnull}
  \begin{aligned}
    ^{(\OOO)}\pi_{44} & = 0,\\
    ^{(\OOO)}\pi_{34} & = 0,\\
    ^{(\OOO)}\pi_{33} & = -4\xib_a\OOO^a,\\
    ^{(\OOO)}\pi_{4a} & = 0, \\
    ^{(\OOO)}\pi_{3a} & = Y_a,\\
    ^{(\OOO)}\pi_{ab} & =: \Hrot_{ab},
  \end{aligned}
\end{align}
and
\begin{align*}
  \tr ^{(\OOO)}\pi & = \tr \Hrot. 
\end{align*}

From equation~\eqref{eq:Lie4OOE}, we have in particular the following transport equation for $\OOO$ in the $\el$ direction
\begin{align}\label{eq:Nd4OOE}
  \Nd_4\OOO_a & = \chi_{ab}\OOO^b.
\end{align}

We further have the following transport equations for $\Nd \OOO$, $Y$, $\Hrot$ and $\POE$ in the $\el$ direction.
\begin{lemma}
  We have
  \begin{align}\label{eq:Nd4NdO}
    \begin{aligned}
    \Nd_4\Nd_a \OOO_b  & = \half\Nd_a(\trchi) \OOO_b + \Nd_a\chih_{bc} \OOO^c + \chih_{bc}\Nd_a\OOO^c  -(\chih\cdot \Nd)_a \OOO_b  \\
    & \quad - (\ze\cdot \OOO)\chi_{ab} +(\chi\cdot \OOO)_a\ze_b + \dual\be_a\dual \OOO_b.
  \end{aligned}
  \end{align}
  
  We have
  \begin{align}\label{eq:Nd4YOOE}
    \begin{aligned}
      \Nd_4Y- \half\trchi Y & = 2(\Divd\ze)\OOO+ \le(2\Nd\otimesh\ze + \chih\cdot\chibh - \chibh\cdot\chih\ri)\cdot\OOO \\
      & \quad + \chih\cdot Y-4\ze\cdot\Nd\OOO + 2\sigma\dual\OOO + 2(\xib\cdot\OOO)\ze,
    \end{aligned}
  \end{align}
  We have
  \begin{align}\label{eq:Nd4H}
    \begin{aligned}
      \Nd_4\Hrot_{ab} & = \le(\Nd(\trchi)\otimes \OOO\ri)_{ab} + \le((\Nd\otimes\chih)\cdot \OOO\ri)_{ab} + \chih_{bc}\Nd_a\OOO^c \\
      & \quad + \chih_{ac}\Nd_b\OOO^c - \le((\chih\cdot \Nd) \otimes \OOO\ri)_{ab} - 2(\ze\cdot \OOO)\chi_{ab} \\
      &\quad + \le(\chi\cdot \OOO \otimes \ze\ri)_{ab} + (\dual\be\otimes\dual \OOO)_{ab}.
    \end{aligned}
  \end{align}
  We have
  \begin{align}
    \label{eq:Nd4POE}
    \begin{aligned}
      \Nd_4\le(r\POE\ri)_{abc} & = \half(r\Nd)_a\Nd_b(\trchi) \OOO_c + (r\Nd)_a\Nd_b\chih_{cd} \OOO_d + (r\Nd)_a\chih_{cd}\Nd_b\OOO_d  -(r\Nd)_a\chih_{bd}\Nd_d\OOO_c\\
      & \quad - (r\Nd)_a\ze_d\OOO_d\chi_{bc} +(r\Nd)_a\chi_{bd}\OOO_d\ze_c + (r\Nd)_a\dual\be_b\dual \OOO_c + \half\Nd_b(\trchi) (r\Nd)_a\OOO_c \\
      & \quad + \Nd_b\chih_{cd} (r\Nd)_a\OOO_d + \chih_{cd}(r\Nd)_a\Nd_b\OOO_d  -\chih_{bd}(r\Nd)_a\Nd_d\OOO_c  -\ze_d(r\Nd)_a\OOO_d\chi_{bc} \\
      & \quad +\chi_{bd}(r\Nd)_a\OOO_d\ze_b + \dual\be_b(r\Nd)_a\dual \OOO_c  -\ze_d\OOO_d (r\Nd)_a\chi_{bc}+\chi_{bd}\OOO_d(r\Nd)_a\ze_b \\
      & \quad + \half (\trchio-\trchi)(r\Nd)_a\Nd_b\OOO_c -r\chih_{ad}\Nd_d\Nd_b\OOO_c \\
      & \quad + \le(-\chi_{ba}\ze_d+\chi_{da}\ze_b+\iin_{bd}\dual\be_a\ri)\Nd_d\OOO_c \\
      & \quad + \le(-\chi_{ca}\ze_d+\chi_{da}\ze_c+\iin_{cd}\dual\be_a\ri)\Nd_b\OOO_d \\
      & \quad -\half r^{-1}\le(\trchi-\trchio\ri)\OOO_b\gd_{ac} - r^{-1}\chih_{bd}\OOO_d\gd_{ac} \\
      & \quad +\half r^{-1}\le(\trchi-\trchio\ri)\OOO_c\gd_{ab} + r^{-1}\chih_{cd}\OOO_d\gd_{ab}.
    \end{aligned}
  \end{align}
\end{lemma}
\begin{proof}
  Equations~\eqref{eq:Nd4NdO},~\eqref{eq:Nd4H} and~\eqref{eq:Nd4POE} follow by duality from respectively~\eqref{eq:Nd3NdO},~\eqref{eq:Nd3H} and~\eqref{eq:Nd3POE}.\\ 

  Using relations~\eqref{eq:relOOEnull}, we have
  \begin{align*}
    \D_3\OOO & = (\xib \cdot \OOO) \el + (\ze\cdot \OOO)\elb +(\chib\cdot \OOO) + Y.
  \end{align*}
  Using equation~\eqref{eq:Lie4OOE}, we have
  \begin{align*}
    \D_4(\D_3\OOO) & = \D_3(\D_4\OOO) + \D_{[\el,\elb]}\OOO + \R(\el,\elb,\eb,\ea)\OOO^a\eb \\
                   & =\D_3(\chi\cdot \OOO) + \D_{[\el,\elb]}\OOO + \R(\el,\elb,\eb,\ea)\OOO^a\eb.
  \end{align*}
  Using the above two equations and projecting on $\text{T}S$, this gives
  \begin{align*}
    \Nd_4Y & = \Nd_3(\chi\cdot \OOO) + \Nd_{[\el,\elb]}\OOO - 2\sigma\dual \OOO +2(\xib \cdot \OOO) \ze -\Nd_4(\chib\cdot \OOO).
  \end{align*}
  Using relations~\eqref{eq:Riccirel}, we have
  \begin{align*}
    \Nd_{[\el,\elb]}\OOO & = \Nd_{-4\ze_a\ea -2\omb\el}\OOO \\
                      & = -4\ze\cdot\Nd \OOO -2\omb\Nd_4\OOO \\
                      & = -4\ze\cdot\Nd \OOO -2\omb(\chi\cdot \OOO),
  \end{align*}
  and, using~\eqref{eq:relOOEnull} again, we obtain
  \begin{align}\label{eq:Nd4YOOE1}
    \begin{aligned}
      \Nd_4Y -\half \trchi Y & = \le(\Nd_3\chi-\Nd_4\chib\ri)\cdot \OOO + \le(\chi\cdot\chib-\chib\cdot\chi\ri) \cdot \OOO + \chih\cdot Y \\
      & \quad - 4\ze\cdot\Nd \OOO -2\omb\chi\cdot \OOO -2\sigma\dual \OOO +2(\xib\cdot \OOO)\ze.
    \end{aligned}
  \end{align}
  Using equations~\eqref{eq:Nd3trchi} and~\eqref{eq:Nd4trchib} together with the relations of Lemma~\ref{lem:relgeodnull}, we have
  \begin{align*}
    \half\le(\Nd_3\trchi-\Nd_4\trchib\ri) & = 2\Divd\ze + \omb\trchi. 
  \end{align*}
  Using equations~\eqref{eq:Nd3chih} and~\eqref{eq:Nd4chibh} together with the relations of Lemma~\ref{lem:relgeodnull}, we have
  \begin{align*}
    \Nd_3\chih-\Nd_4\chibh & = 2\Nd\otimesh\ze + \omb\chih. 
  \end{align*}
  Additionally, we also check that
  \begin{align*}
    \le(\chi\cdot\chib - \chib\cdot\chi\ri)_{ab} & = \chi_{ac}\chib_{cb}-\chib_{ac}\chi_{cb} \\
                                           & = \chih_{ac}\chibh_{cb} - \chibh_{ac}\chih_{cb} \\
                                           & = \le(\chih\cdot\chibh-\chibh\cdot\chih\ri)_{ab}.
  \end{align*}

  Using the above formulas, equation~\eqref{eq:Nd4YOOE1} directly rewrites as~\eqref{eq:Nd4YOOE} as desired.
\end{proof}

We record the following lemma, which motivates the definition of $\POE$.\footnote{Using that $\D^2\OOO= 0$ in the Euclidean space, one can deduce from the formula of Lemma~\ref{lem:D2Nd2X} that $\POE =0$ in the Euclidean case (see Lemma~\ref{lem:preciseOOESast}). Reciprocally, from the control of the rotation coefficients $\POE,\Hrot,Y \simeq 0$ defined in this section, one can deduce from this formula a control of $\D^2\OOE \simeq 0$ (see Lemma~\ref{lem:DOOEOOITTL4}).}
\begin{lemma}\label{lem:D2Nd2X}
  Let $X$ be an $S$-tangent vectorfield. We have the following formula
  \begin{align}\label{eq:D2Nd2X}
    \begin{aligned}
      \D^2_{a,b}X & = \Nd^2_{a,b}X-r^{-2}X_b\ea + \half r^{-1}\gd_{ab}\le(\D_4-\D_3\ri)X \\
      & \quad + \half r^{-1}\le(\Nd_aX_b+\Nd_bX_a\ri)(\elb-\el) + \le(\EEE\le(\D^2,\Nd^2\ri)\cdot X\ri)_{ab},
    \end{aligned}
  \end{align}
  where
  \begin{align*}
    \le(\EEE\le(\D^2,\Nd^2\ri)\cdot X\ri)_{ab} & :=  \half X_c\Nd_a\le(\chi-r^{-1}\gd\ri)_{bc} \elb + \half X_c\Nd_a\le(\chib+r^{-1}\gd\ri)_{bc}\el \\
                                               & \quad + \half \Nd_bX_c\le(\chi_{ca}-r^{-1}\gd_{ca}\ri)\elb + \half \Nd_bX_c\le(\chib_{ca}+r^{-1}\gd_{ca}\ri)\el \\
                                               & \quad + \half r^{-1}X_b\le((\chib_{ac}+r^{-1}\gd_{ac})\ec - (\chi_{ac}-r^{-1}\gd_{ac})\ec + \ze_a\elb + \ze_a\el\ri) \\
                                               & \quad -\half (\chi_{ab}-r^{-1}\gd_{ab})\D_3X - \half (\chib_{ab}+r^{-1}\gd_{ab})\D_4X.
  \end{align*}
\end{lemma}
\begin{proof}
  We first start by recording the following formula, which follows from the relations~\eqref{eq:Riccirel}
  \begin{align}\label{eq:DNdX}
    \begin{aligned}
      \D_bX & = \Nd_bX +\half X_c\chi_{cb}\elb + \half X_c \chib_{cb}\el \\
      & = \Nd_bX + \half r^{-1} X_b (\elb-\el) + \le(\EEE(\D,\Nd)\cdot X\ri)_b,
    \end{aligned}
  \end{align}
  where
  \begin{align*}
    \le(\EEE(\D,\Nd)\cdot X\ri)_{b} & := \half X_c \le(\chi_{cb} - r^{-1}\gd_{cb}\ri)\elb + \half X_c \le(\chib_{cb}+r^{-1}\gd_{cb}\ri)\el.
  \end{align*}
  
  Using formula~\eqref{eq:DNdX} and the relations~\eqref{eq:Riccirel}, we have
  \begin{align}\label{eq:D2X1}
    \begin{aligned}
      \D^2_{a,b}X & = \D_a(\D_bX) - \D_{\D_\ea \eb}X \\
      & = \D_a\le(\Nd_bX + \half r^{-1}X_b(\elb-\el) + \le(\EEE(\D,\Nd)\cdot X\ri)_b\ri) \\
      & \quad -\D_{\Nd_\ea \eb}X - \half\chi_{ab}\D_3X - \half \chib_{ab}\D_4X.
    \end{aligned}
  \end{align}

  Using formula~\eqref{eq:DNdX}, we have the following equations
  \begin{align}\label{eq:D2X1annex}
    \begin{aligned}
      \D_a\le(\Nd_bX\ri) & = \Nd_a(\Nd_bX) + \half r^{-1}\Nd_bX_a\le(\elb-\el\ri) + \le(\EEE(\D,\Nd)\cdot\Nd_b X\ri)_a \\
      & = \Nd^{2}_{a,b}X + \Nd_{\Nd_a\eb}X + \half r^{-1}\Nd_bX_a\le(\elb-\el\ri) + \le(\EEE(\D,\Nd)\cdot\Nd_b X\ri)_a,
    \end{aligned}
  \end{align}
  and
  \begin{align}\label{eq:D2X1annexbis}
    \begin{aligned}
      \D_a\le(\half r^{-1}X_b(\elb-\el)\ri) & = \half r^{-1}\Nd_aX_b(\elb-\el) +\half r^{-1}X_c(\Nd_a\eb)^c(\elb-\el) \\
      & \quad + \half r^{-1}X_b \D_a(\elb-\el).
    \end{aligned}
  \end{align}

  Combining the terms of~\eqref{eq:D2X1},~\eqref{eq:D2X1annex} and~\eqref{eq:D2X1annexbis} containing $\Nd_\ea\eb$, we obtain
  \begin{align}\label{eq:D2X1annexbisbis}
    \begin{aligned}
    & \Nd_{\Nd_a\eb}X +\half r^{-1}X_c(\Nd_a\eb)^c(\elb-\el) -\D_{\Nd_\ea \eb}X \\
    = & \Nd_{\Nd_a\eb}X +\half r^{-1}X_c(\Nd_a\eb)^c(\elb-\el) \\
    & - \le(\Nd_{\Nd_\ea\eb}X +\half r^{-1}X_c(\Nd_\ea\eb)^c(\elb-\el) +\le(\EEE(\D,\Nd)\cdot X\ri)_c(\Nd_\ea\eb)^c\ri) \\
    = & -\le(\EEE(\D,\Nd)\cdot X\ri)_c(\Nd_\ea\eb)^c.
    \end{aligned}
  \end{align}

  Thus, we rewrite~\eqref{eq:D2X1} using~\eqref{eq:D2X1annex}, \eqref{eq:D2X1annexbis} and~\eqref{eq:D2X1annexbisbis} as
  \begin{align}\label{eq:D2X2}
    \begin{aligned}
      \D^2_{a,b}X & = \Nd^{2}_{a,b}X + \half r^{-1}\Nd_bX_a\le(\elb-\el\ri) + \half r^{-1}\Nd_aX_b(\elb-\el) \\
      & \quad + \half r^{-1}X_b \D_a(\elb-\el) - \half\chi_{ab}\D_3X - \half \chib_{ab}\D_4X + \EEE^1,
    \end{aligned}
  \end{align}
  where
  \begin{align*}
    \EEE^1 & := \D_a\le(\le(\EEE(\D,\Nd)\cdot X\ri)_b\ri) + \le(\EEE(\D,\Nd)\cdot\Nd_b X\ri)_a -\le(\EEE(\D,\Nd)\cdot X\ri)_c(\Nd_\ea\eb)^c.
  \end{align*}
  We rewrite
  \begin{align}\label{eq:D2X2annex}
    \half r^{-1}X_b \D_a(\elb-\el) & = - r^{-2} X_b \ea + \EEE^2,
  \end{align}
  and
  \begin{align}\label{eq:D2X2annexbis}
    - \half\chi_{ab}\D_3X - \half \chib_{ab}\D_4X & = \half r^{-1}\gd_{ab} (\D_4-\D_3)X + \EEE^3,
  \end{align}
  where
  \begin{align*}
    \EEE^2 & := \half r^{-1}X_b (\D_a\elb +r^{-1}\ea) -\half r^{-1}X_b (\D_a\el-r^{-1}\ea),\\
    \EEE^3 & := -\half (\chi_{ab}-r^{-1}\gd_{ab})\D_3X - \half (\chib_{ab}+r^{-1}\gd_{ab})\D_4X.
  \end{align*}
  Thus, we obtain the following rewriting of~\eqref{eq:D2X2} using~\eqref{eq:D2X2annex} and~\eqref{eq:D2X2annexbis}
  \begin{align}\label{eq:D2X3}
    \begin{aligned}
      \D^2X_{a,b} & = \Nd^{2}_{a,b}X + \half r^{-1}(\Nd_bX_a+\Nd_aX_b)\le(\elb-\el\ri) - r^{-2} X_b \ea \\
      & \quad + \half r^{-1}\gd_{ab} (\D_4-\D_3)X + \EEE^4,
    \end{aligned}
  \end{align}
  with
  \begin{align*}
    \EEE^4 & := \EEE^1+\EEE^2+\EEE^3.
  \end{align*}

  Using formula~\eqref{eq:DNdX} and~\eqref{eq:Riccirel}, we simplify the terms composing $\EEE^4$ as
  \begin{align*}
    \EEE^1 & = \half X_c\Nd_a\le(\chi-r^{-1}\gd\ri)_{bc} \elb + \half X_c\Nd_a\le(\chib+r^{-1}\gd\ri)_{bc}\el \\
           & \quad + \half X_c\le(\chi_{cb}-r^{-1}\gd_{cb}\ri) \le(\chib_{ad}e_d + \ze_a\elb\ri) + \half X_c\le(\chib_{cb}+r^{-1}\gd_{cb}\ri)\le(\chi_{ad}e_d \ze_a\el\ri) \\
           & \quad + \half \Nd_bX_c\le(\chi_{ca}-r^{-1}\gd_{ca}\ri)\elb + \half \Nd_bX_c\le(\chib_{ca}+r^{-1}\gd_{ca}\ri)\el, 
  \end{align*}
  and
  \begin{align*}
    \EEE^2 & = \half r^{-1}X_b\le((\chib_{ac}+r^{-1}\gd_{ac})\ec - (\chi_{ac}-r^{-1}\gd_{ac})\ec + \ze_a\elb + \ze_a\el\ri),
  \end{align*}
  and the formula~\eqref{eq:D2Nd2X} follows. This finishes the proof of the lemma.  
\end{proof}

\section{The last cones geodesic foliation}\label{sec:deflastconesfoliation}
For all $11/4 \leq \ub' \leq \uba$, we define $\CCb_{\ub'}$ to be the incoming null cones with vertex $\o(\ub')$ and we note $\ub'$ its associated optical function. We define $\elb'$ by $\elb' := -(\D\ub')^\sharp$. We define $u'$ to be the associated normalised affine parameter on the cones $\CCb_{\ub'}$, \emph{i.e.}
\begin{align*}
  \elb'(u') & = 2,\\
  u'|_{\o(\ub')} & = \ub'.
\end{align*}
This yields a foliation by $2$-spheres that we note
\begin{align*}
  S' & : = \le(S'_{u',\ub'}\ri)_{11/4\leq \ub' \leq \uba,~ 5/4 \leq u' \leq \uba}
\end{align*}
and that we call the \emph{last cones geodesic foliation}. See Figure~\ref{fig:lastconesfoliation} for a graphic representation of the domain
\begin{align*}
  \MM' & := \bigcup_{11/4\leq \ub' \leq \uba,~ 5/4 \leq u' \leq \uba} S'_{u',\ub'}
\end{align*}
covered by $S'$.\\

We define $\el'$ such that $(\elb',\el')$ forms a null pair orthogonal to $S'$. With respect to this null pair, we let $\chi',\xi',\eta',\ze',\om'$ and $\chib',\xib',\etab',\ze',\omb'$ denote the associated null connection coefficients and $\alpha',\beta',\rho'$, $\sigma',\betab',\alphab'$ the associated null curvature components (see Section~\ref{sec:nulldecomp}).\\

We denote $Y' := -(\D u')^\sharp$, and $\yy'$ the optical defect
\begin{align*}
  \g(Y',Y') & =: -2\yy'.
\end{align*}

From the above definitions, we have
\begin{align}\label{eq:elut}
  \begin{aligned}
    \elb'({u'}) & = 2, & \el'({u'}) & = \yy',\\
    \elb'({\ub'}) & = 0, & \el'({\ub'}) & = 2.
  \end{aligned}
\end{align}
and
\begin{align}\label{eq:Yteltelbt}
  Y' & = \el' + \half\yy'\elb'.
\end{align}

From the geodesic equation $\D_{\elb'}\elb'$ and relations~\eqref{eq:Riccirel}, we have
\begin{align}\label{eq:rellastconegeod1}
  \xib' & = 0, & \omb' & = 0.
\end{align}
From relations~\eqref{eq:elut} and~\eqref{eq:Yteltelbt} we have (see the analogous derivations of the relations of Lemma~\ref{lem:relgeodnull})
\begin{align}\label{eq:rellastconegeod2}
  \ze' & = \eta' = - \etab'.
\end{align}
From relations~\eqref{eq:elut} and~\eqref{eq:Yteltelbt} we have (see the analogous derivations of relations of Lemma~\ref{lem:relCCba})
\begin{align}\label{eq:Ndyyt}
  \elb'(\yy') & = -4\om', &  \Nd'\yy' & = -2\xi'.
\end{align}

\section{The initial layers $\LLb$ and $\LLc$}\label{sec:definitlayer}
\begin{figure}[!h]
  \centering
  \includegraphics[height=8cm]{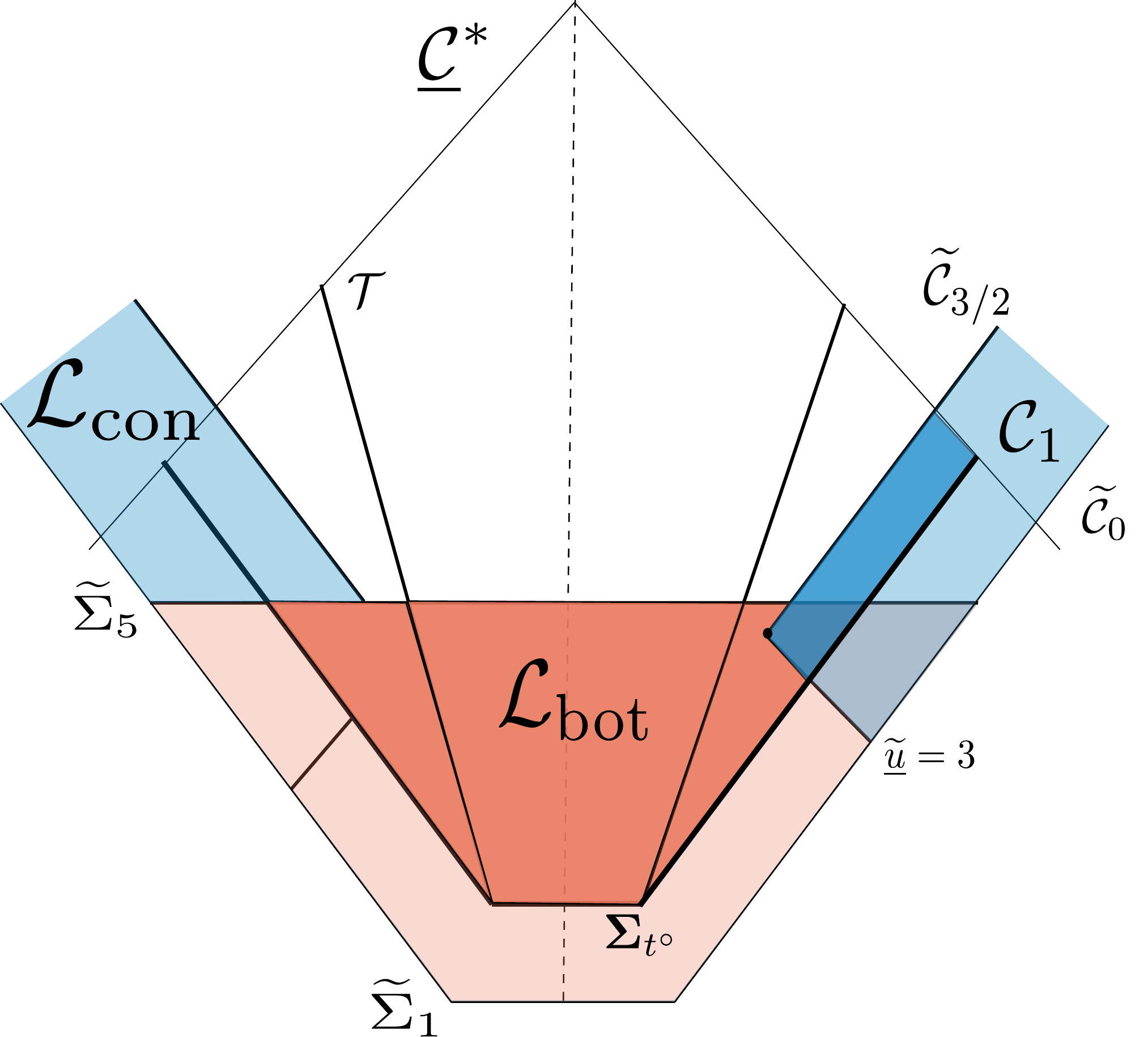}
  \caption{The initial layers $\LLb$ and $\LLc$.}
  \label{fig:initiallayer}
\end{figure}

In this section, we define spacetime regions $\LLb,\LLc \subset \MMt$ which we respectively call the \emph{bottom initial layer} and the \emph{conical initial layer}. The construction and its overlay with the spacetime region $\MM_\uba$ defined in Section~\ref{sec:defcannull} is graphically summarised in Figure~\ref{fig:initiallayer}.

\paragraph{The bottom initial layer $\LLb$}
We define a bottom initial layer $\LLb\subset\MMt$ to be a spacetime region covered by a coordinates system, which we call \emph{bottom initial layer coordinates} $(x^\mu)$ such that in these coordinates $\LLb \simeq \cup_{x^0\in[1,5]} B((x^0,0,0,0),x^0)$, where $B$ denote the coordinate balls. The region $\LLb$ shall be a future neighbourhood of $\Sit_1 = \{x^0=1\}$, with boundary locally coinciding with the null hypersurface $\CCt_0$, \emph{i.e.} $\cup_{x^0\in[1,5]}\pr B((x^0,0,0,0),x^0) \subset \CCt_0$. We denote by $\Sit_{x^0}$ the level sets of the time coordinate $x^0$. We denote by $\Tf^\bott$ the timelike future-pointing unit normal to $\Sit$. We denote by $\Nf^\bott$ the outward-pointing unit normal to the $2$-spheres level sets of $x^0$ and $\sum_{i=1}^3(x^i)^2$ which is orthogonal to $\Tf^\bott$. 

\paragraph{The conical initial layer $\LL_{con}$}
We define a conical initial layer $\LLc\subset\MMt$ to be a spacetime region covered by two optical functions $\ut,\ubt$ such that $\LLc \simeq [0,3/2]_{\ut} \times [3,+\infty)_{\ubt} \times \tilde{S}_{\ut,\ubt}$. Moreover, we require that $\CCt_0 = \le\{\ut=0\ri\}.$\\

We note $\rt$ the area radius of the $2$-spheres $\tilde{S}_{\ut,\ubt}$. We define its null lapse $2(\Omt)^{-2} := -\g(\D \ut,\D\ubt)$ and the associated null pair $(\elbt,\elt)$ by $\elt := - \Omt (\D \ut)^\sharp$ and $\elbt := - \Omt(\D\ubt)^\sharp$.\\ 

We write $\chit,\xit,\etat,\zet,\omt$ and $\chibt,\xibt,\etabt,\zet,\ombt$ the associated null connection coefficients and $\alphat,\betat,\rhot$, $\sigmat,\betabt,\alphabt$ the associated null curvature components (see Section~\ref{sec:nulldecomp}). Note that since the functions $\ut,\ubt$ are assumed to be optical, we have $\xit = \xibt= 0$.\\

We moreover define the following intersections of domains
\begin{align*}
  \LLbint & := \MM^\intr\cap\LLb,\\
  \LLbext & := \MM^\ext\cap\LLb,\\
  \LLcext & := \MM^\ext\cap\LLc.
\end{align*}

\section{General change of null frames}\label{sec:defchgframe}
Let $S = (S_{u,\ub})$ and $S' = (S'_{u',\ub'})$ be two (local) foliations by spacelike $2$-spheres. Let $(\el,\elb)$ and $(\el',\elb')$ be two null pairs associated to the foliations $S$ and $S'$ respectively. We have the following lemma (see~\cite[Lemma 3.1]{Kla.Sze19a} for a proof).
\begin{lemma}\label{lem:transframe}
  There exists two $S$-tangent vectorfields $f$ and $\fb$ and a scalar function $\la$ such that
  \begin{align*}
    \begin{aligned}
      \el' & = \la \left(\el + f + \frac{1}{4}|f|^2\elb\right),\\
      \elb' & = \la^{-1}\left(\left(1+\half f\cdot\fb + \frac{1}{16} |f|^2|\fb|^2\right)\elb +\fb + \frac{1}{4}|\fb|^2f + \frac{1}{4}|\fb|^2\el \right).
    \end{aligned}
  \end{align*}
  Moreover, if $(\ea)_{a=1,2}$ is a (local) orthonormal null frame of $\text{T}S$, the pair of vectorfields $(\ea')_{a=1,2}$ defined for $a=1,2$ by
  \begin{align}\label{eq:chgframeea}
    \begin{aligned}
      \ea' := \ea + \half \fb_a f + \half \fb_a\el + \left(\half f_a + \frac{1}{8}|f|^2\fb_a\right)\elb,
    \end{aligned}
  \end{align}
  is a (local) orthonormal frame of $\text{T}S'$. The triplet $(\la,f,\fb)$ is called the \emph{transition coefficients} of the change of frame.
\end{lemma}

We have the following definition for projections of $S$-tangent tensors to $S'$-tangent tensors.
\begin{definition}
  Let $(\ea)_{a=1,2}$ be an orthonormal frame of $\text{T}S$. For $\phi'$ an $S'$-tangent $r$-tensor, we define its projection $(\phi')^\dg$ to be the $S$-tangent $r$-tensor defined by
  \begin{align*}
    (\phi')^\dg_{a_1\cdots a_r} = (\phi')^\dg\left({e_{a_1}},\cdots,{e_{a_r}}\right) := \phi'\left({e_{a_1}}',\cdots,{e_{a_r}}'\right) = \phi'_{a_1\cdots a_r},
  \end{align*}
  where the frame $(\ea')_{a=1,2}$ is associated to $(\ea)_{a=1,2}$ via formula~\eqref{eq:chgframeea}.\\
  Reciprocally, for $\phi$ an $S$-tangent $r$-tensor, we define its projection $(\phi)^\ddg$ to be the $S'$-tangent $r$-tensor defined by
  \begin{align*}
    \phi^\ddg_{a_1\cdots a_r} = \phi^\ddg\left({e_{a_1}}',\cdots,{e_{a_r}}'\right) := \phi\left({e_{a_1}},\cdots,{e_{a_r}}\right) = \phi_{a_1\cdots a_r}.
  \end{align*}
\end{definition}
\begin{remark}
  In other terms, we have
  \begin{align*}
    (\phi')^\dg := (\phi')^a\ea, \quad \phi^\ddg := \phi^a\ea',
  \end{align*}
  where the frames $\ea'$ and $\ea$ are associated via formula~\eqref{eq:chgframeea}. It can be seen from formula~\eqref{eq:chgframeea} that this definition does not depend on the choice of frame on $S$
\end{remark}

We have the following transition formulas for projected covariant derivatives.
\begin{proposition}\label{prop:chgframederiv}
  Let $\phi'$ be an $S'$-tangent $r$-tensor. We have
  \begin{align*}
    \Nd_4(\phi')^\dg & = \la^{-1}\le(\Nd'_4\phi'\ri)^\dg  - f\cdot\Nd(\phi')^\dg -\frac{1}{4}|f|^2 \Nd_3(\phi')^\dg + (\phi')^\dg\cdot\Err(\Nd_4,\Nd'_4),\\
    \Nd(\phi')^\dg & = \le(\Nd'\phi'\ri)^\dg  -\half \fb  f \cdot\Nd(\phi')^\dg - \half \fb \Nd_4(\phi')^\dg -\le(\half f + \frac{1}{8}|f|^2\fb\ri)\Nd_3(\phi')^\dg + (\phi')^\dg\cdot \Err(\Nd,\Nd'), \\
    \Nd_3(\phi')^\dg & = \la \le(\Nd'_3\phi'\ri)^\dg   - \le(\half f\cdot\fb + \frac{1}{16}|f|^2|\fb|^2 \ri)\Nd_3(\phi')^\dg - \le(\fb + \quar|\fb|^2f \ri)\cdot\Nd (\phi')^\dg \\
  & \quad - \quar|\fb|^2\Nd_4(\phi')^\dg + (\phi')^\dg\cdot\Err(\Nd_3,\Nd_3').
  \end{align*}
  where the $S$-tangent tensors $\Err(\Nd_4,\Nd_4')$, $\Err(\Nd,\Nd')$ and $\Err(\Nd_3,\Nd_3')$ are bilinear (or higher nonlinear) error terms composed of $f$, $\fb$, their first order derivatives and null connection coefficients for the $S$ and $S'$ foliations.\\
  
  Reciprocally, using that for an $S$-tangent $r$-tensor $\phi$ we have $\le(\phi^\ddg\ri)^\dg = \phi$, we obtain symmetric formulas for $\Nd_4'\phi^\ddg$, $\Nd'\phi^\ddg$ and $\Nd_3'\phi^\ddg$. 
\end{proposition}
\begin{proof}
  We compute the tensorial formulas at a point where we choose a normal frame $(\ea)_{a=1,2}$, \emph{i.e.} $\g(\ea,\eb) = \de_{ab}$ and $\D\ea = 0$. We denote by $(\ea')_{a=1,2}$ the frame given by the frame transformation~\eqref{eq:chgframeea}.\footnote{The frame $(\ea')_{a=1,2}$ is not geodesic normal and the error terms $\Err$ will correspond to the projection of its connection coefficients $\D\ea'$.}\\

  We start with the first formula. We shall do the computations for $\phi'$ an $S'$-tangent $1$-form, and the results of Proposition~\ref{prop:chgframederiv} will follow by simple generalisation. With respect to the frames $(\ea)_{a=1,2}$ and $(\ea')_{a=1,2}$, we have
\begin{align*}
  \Nd_4(\phi')^\dg_a & = \el\le(\phi'_a\ri) \\
                     & = \le(\la^{-1}\el' - f -\quar|f|^2\elb\ri)\le(\phi'_a\ri) \\
                     & = \la^{-1}\Nd'_4\phi'_a +\la^{-1}\phi'_b\g(\Nd'_4\ea',\eb') - f\cdot\Nd(\phi')^\dg_a - \quar |f|^2 \Nd_3(\phi')^\dg_a.
\end{align*} 
We deduce that the $S$-tangent tensor $\Err(\Nd_4,\Nd_4')$ can be expressed as
\begin{align*}
  \Err(\Nd_4,\Nd_4')_{ba} = \la^{-1}\g\le(\Nd'_4\ea',\eb'\ri) = \g\le(\la^{-1}\D_{\el'}\ea',\eb'\ri).
\end{align*}
Using the transition formulas from Lemma~\ref{lem:transframe}, we have
\begin{align*}
  \g\le(\la^{-1}\D_{\el'}\ea',\eb'\ri) & = \g\le(\D_{\el+ f}\le(\ea + \half\fb_a \el + \half f_a \elb\ri),\eb'\ri) + \lot \\
                                       & = \g\le(\half\el(\fb_a)\el+\half \el(f_a) \elb,\eb'\ri) + \lot \\
                                       & = \lot,
\end{align*}
where $\lot$ denotes bilinear (or higher) error terms.\\

For the second formula, we have
\begin{align*}
  \Nd_b(\phi')^\dg_a & = \eb (\phi'_a) \\
                     & = \le(\eb' - \half \fb_bf - \half \fb_b\el - \le(\half f_b + \frac{1}{8}|f|^2\fb_b\ri)\elb\ri)(\phi'_a) \\
                     & = \le(\Nd'\phi'\ri)^\dg_{ba} + \phi'_c \g(\D_{\eb'}\ea',\ec') -\half \fb_b f \cdot\Nd(\phi')^\dg_a - \half \fb_b \Nd_4(\phi')^\dg_a -\le( \half f_b + \frac{1}{8}|f|^2\fb_b\ri)\Nd_3(\phi')^\dg_a.
\end{align*}
We therefore have $\Err(\Nd,\Nd')_{cba} = \g(\D_{\eb'}\ea',\ec')$ and we check similarly as before that $\Err(\Nd,\Nd') = \lot$.\\

For the third and last formula, we have
\begin{align*}
  \Nd_3(\phi')^\dg_a & = \elb (\phi'_a) \\
                     & = \le(\la\elb'-\le(\half f\cdot\fb + \frac{1}{16} |f|^2|\fb|^2\ri)\elb -\fb - \frac{1}{4}|\fb|^2f - \frac{1}{4}|\fb|^2\el\ri) \\
                     & = \la \le(\Nd_3'\phi'\ri)^\dg_a + \la \phi'_b \g(\Nd_3'\ea',\eb') - \le(\half f\cdot\fb + \frac{1}{16}|f|^2|\fb|^2 \ri)\Nd_3(\phi')^\dg_a \\
  & \quad - \le(\fb + \quar|\fb|^2f \ri)\cdot\Nd (\phi')^\dg_a - \quar|\fb|^2\Nd_4(\phi')^\dg_a.
\end{align*}
And thus $\Err(\Nd_3,\Nd_3')_{ba} = \g(\D_{\elb'}\ea',\eb') = \lot$ as desired.
\end{proof}

We have the following transformation formulas for the null connection coefficients (see~\cite[Proposition 3.3 and Appendix A]{Kla.Sze19a}).
\begin{proposition}\label{prop:transconn}
  Under the transitions formulas of Lemma~\ref{lem:transframe}, we have
  \begin{subequations}\label{eq:transconn}
    \begin{align}
      \la^{-2}\xi' & = \xi^\ddg +\half \la^{-1}\Nd'_{4}f^\ddg + \frac{1}{4}\trchi f^\ddg + \Err(\xi,\xi'),\label{eq:Nd4f}\\
      \la^2\xib' & = \xib^\ddg + \half\la \Nd'_3\fb^\ddg  +\frac{1}{4}\trchib\fb^\ddg + \Err(\xib',\xib),\label{eq:Nd3fb}\\
      \la^{-1}\chi' & = \chi^\ddg + \Nd' f^\ddg + \Err(\chi,\chi'),\label{eq:Ndf}\\
      \la\chib' & = \chib^\ddg + \Nd'\fb^\ddg + + \Err(\chib,\chib'),\label{eq:Ndfb}\\
      \ze' & = \ze^\ddg -\Nd'(\log\la) -\quar\trchib f^\ddg+\quar\trchi\fb^\ddg +\Err(\ze,\ze'),\label{eq:Ndlogla}\\
      \eta' & = \eta^\ddg +\half\la\Nd'_3f^\ddg +\quar\trchi\fb^\ddg +\Err(\eta,\eta'),\label{eq:Nd3f}\\
      \etab' & = \etab^\ddg + \half \la^{-1}\Nd'_4\fb^\ddg + \quar\trchib f^\ddg +\Err(\etab,\etab'),\label{eq:Nd4fb}\\
      \la^{-1}\om' & = \om -\half\la^{-1}\el'(\log\la)+\Err(\om,\om'),\label{eq:Nd4la}\\
      \la \omb' & = \omb +\half\la\elb'(\log\la)+\Err(\omb,\omb'),\label{eq:Nd3la}
    \end{align} 
    where 
    \begin{align*}
      \Err(\xi,\xi') & := \om f^\ddg + \half f \cdot\chih + \lot,\\
      \Err(\xib,\xib') & := \omb\fb+\half\fb\cdot\chibh + \lot,\\
      \Err(\chi,\chi') & := f^\ddg\otimes \eta^\ddg + f^\ddg\otimes\ze^\ddg + \fb^\ddg\otimes\xi^\ddg + \half (\fb\cdot f) \chi^\ddg- \half (f\cdot\fb) \chi' + \half f^\ddg\cdot\chi'\otimes\fb^\ddg - \quar |f^2| \chib +\lot,\\
      \Err(\chib,\chib') & := \fb^\ddg\otimes \etab^\ddg - \fb^\ddg\otimes\ze^\ddg + f^\ddg\otimes\xib^\ddg + \half (\fb\cdot f) \chib^\ddg- \half (f\cdot\fb) \chib' + \half \fb^\ddg\cdot\chib'\otimes f^\ddg - \quar |\fb^2| \chi +\lot,\\
      \Err(\ze,\ze') & := \om \fb -\omb f  -\half\chibh\cdot f +\quar\la^{-1}\trchi'\fb^\ddg +\half\la^{-1}\fb\cdot\chih' + \lot,\\
      \Err(\eta,\eta') & := -\omb f +\half\fb\cdot\chih' +\lot,\\
      \Err(\etab,\etab') & := \half f\cdot\chibh' +\lot,\\
      \Err(\om,\om') & := \half f\cdot(\ze-\etab)-\frac{1}{8}\trchib|f|^2+\half\la^{-2}\fb\cdot\xi' +\lot,\\
      \Err(\omb,\omb') & := -\half\fb\cdot\ze-\half\fb^\ddg\cdot\eta' +\half f\cdot\xib +\frac{1}{8}(f\cdot\fb)\trchib +\lot,
    \end{align*} 
    where $\lot$ denotes trilinear (or higher nonlinear) error terms composed of a null connection coefficient and transition coefficients.\\

    Using the change of derivatives formulas from Proposition~\ref{prop:chgframederiv}, the same formulas as~\eqref{eq:transconn} hold for derivatives taken with respect to the frame $(\elb,\el,\ea)$, up to additional nonlinear error terms.
  \end{subequations}
\end{proposition}

We have the following transformation formulas for the null curvature components (see~\cite[Proposition 3.3 and Appendix A]{Kla.Sze19a}).
\begin{proposition}\label{prop:transcurv}
  Under the transitions formulas of Lemma~\ref{lem:transframe}, we have
  \begin{subequations}
    \begin{align*}
      \la^{-2}\al' & = \al^\ddg +\Err(\al,\al'), & \la^{-1}\be' & = \be^\ddg +\Err(\be,\be'), &  \rho' & = \rho + \Err(\rho,\rho'),\\
      \la^2\alb' & = \alb^\ddg +\Err(\alb,\alb'), & \la \beb' & = \beb^\ddg + \Err(\beb,\beb'), & \sigma' & = \sigma +\Err(\sigma,\sigma'),   
    \end{align*}  
    where
    \begin{align*}
      \Err(\al,\al') & := (f^\ddg\otimesh\be^\ddg-\dual f^\ddg\otimesh\dual\be^\ddg)  +\lot,\\
      \Err(\be,\be') & := \frac{3}{2}(f^\ddg\rho+\dual f^\ddg\sigma) +\half\al^\ddg\cdot\fb^\ddg +\lot,\\
      \Err(\rho,\rho') & := \fb\cdot\beta-f\cdot\beb+\lot,\\
      \Err(\sigma,\sigma') & := -\fb\cdot\dual\be-f\cdot\dual\beb+\lot,\\
      \Err(\beb,\beb') & := -\frac{3}{2}(\fb^\ddg\rho+\dual\fb^\ddg\sigma) -\half\alb^\ddg\cdot f^\ddg +\lot,\\
      \Err(\alb,\alb') & := -(\fb^\ddg\otimesh\beb^\ddg-\dual\fb^\ddg\otimesh\dual\beb^\ddg)+\lot,
    \end{align*}
    where $\lot$ denotes trilinear error terms composed of a null curvature component and transition coefficients.
  \end{subequations}
\end{proposition}

\begin{remark}
  Although we did not present the full explicit formulas for the nonlinear error terms appearing in Propositions~\ref{prop:chgframederiv}, \ref{prop:transconn}, \ref{prop:transcurv}, we claim that they satisfy a conservation of signature principle,\footnote{See~\cite{Chr.Kla93} for further discussion on signature.} which roughly states that the transition coefficients $f,\la$ and $\fb$ are respectively paired with higher, same and lower signature null curvature components, null connection coefficients, derivatives, etc. Since a drop of signature corresponds at most to a drop of decay by $r^{-1}$, the control of the error terms in Section~\ref{sec:initlayer} can be obtained without decay assumptions on $\la,\fb$, provided that the coefficient $f$ decays as $r^{-1}$. See also~\cite[Remark 4.1.4]{Kla.Sze17}.
\end{remark}

\chapter{Norms, bootstrap assumptions and consequences}\label{sec:normsBA}
\section{Preliminary definitions}
We define the following local frame norms.
\begin{definition}[Frame norms]\label{def:framenorm}
  For an orthonormal frame $(e_\mu)_{\mu=0..3}$, we define the following associated frame norm
  \begin{align*}
    \le|F\ri|^2 & := \sum_{i_1,\cdots i_l = 0\cdots 3}\le|F_{i_{1}\cdots i_{l}}\ri|^2,
  \end{align*}
  for all spacetime tensor $F$ and where $F_{i_1\cdots i_l}$ denotes the evaluation of $F$ on the $l$-uplet $(e_{i_1},\cdots,e_{i_l})$.\\
  In this paper, most frame norms are equivalents. We shall precise with respect to which frame the norms are taken only when it is relevant.
\end{definition}
\begin{remark}\label{rem:defframenorm}
  For an orthonormal frame $(e_\mu)_{\mu=0..3}$, the frame norm from the above definition does not depend on the choice of spacelike orthonormal vectors $e_i$, and we have
  \begin{align*}
    |F|^2 = 2\le|\g(e_0,F)\ri|^2 + \g(F,F).
  \end{align*}
\end{remark}

Using Definition~\ref{def:framenorm}, we define integrals, $L^p$ and $L^pL^q$ norms on the various submanifolds of this paper using the respective (intrinsic) induced metrics. In the case of the null hypersurfaces $\CC$, the integral is defined consistently with the coarea formulas from Lemma~\ref{lem:coarea}. In the case of the null hypersurface $\CCba$, the integral is defined consistently with the foliation by the $2$-spheres of the canonical foliation.

\begin{definition}[$\HHt$ norm]\label{def:Ht}
  Let $(S,\gd)$ be a Riemannian $2$-sphere with area radius $r$.
  Let $F$ be an $S$-tangent tensor. We define the (scaling homogeneous) $\HHt$ norm of $F$ on $(S,\gd)$ to be
  \begin{align*}
    \norm{F}_{\HHt(S)} & := r^{1/2}\norm{F}_{H^{1/2}\le(S,r^{-2}\gd\ri)},
  \end{align*}
  where for a Riemannian $2$-sphere $(S,\ga)$, $H^{1/2}\le(S,\ga\ri)$ denotes the standard fractional Sobolev space on $(S,\ga)$, as defined in~\cite{Sha14}.\footnote{Under the regularity assumptions of this paper (see the Bootstrap Assumptions~\ref{BA:mildsphcoordsast} and \ref{BA:mildsphcoordsext}), these intrinsic fractional norms are also equivalent to coordinate defined norms (see Lemma~\ref{lem:compH12}).} 
\end{definition}
\begin{remark}
  The $\HHt$ norm are schematically of the form
  \begin{align*}
    \norm{F}_{\HHt(S)} & \sim r^{-1/2} \norm{(r\Nd)^{\leq \, \half}F}_{L^2(S)}.
  \end{align*}
\end{remark}
\begin{remark}
  Throughout this paper we will extensively use $L^\infty \HHt(S)$ norms, which should be thought of as an upgraded version of $L^\infty L^4(S)$ norms. Namely:
  \begin{itemize}
  \item They are at the same scaling level, similar Sobolev embeddings, elliptic and transport estimates hold (see Lemmas~\ref{lem:Sobsphere}, \ref{lem:ell}, \ref{lem:evolext}),
  \item $\HHt(S)$ embeds in $L^4(S)$ (see Lemma~\ref{lem:Sobsphere}),
  \item $\HHt(S)$ is a natural trace space for spacelike 3D Dirichlet/Neumann problems with boundary $S$.
  \end{itemize}
  For these reasons, most of the Bootstrap Assumptions (see Section~\ref{sec:BA}) are formulated using $L^\infty\HHt(S)$ norms. For conciseness, we did not write the $L^\infty$ and $L^\infty L^4(S)$ norms obtained via Sobolev embeddings (see Lemma~\ref{lem:Sobsphere}) and shall implicitly assume that these norms are controlled as well.
\end{remark}
\begin{remark}
    Using the $\HHt$ spaces of Definition~\ref{def:Ht} (see also Lemma~\ref{lem:compH12}) forces us to reprove classical (intrinsically obtained) estimates (see the Klainerman-Sobolev $\HHt$ embeddings of Lemmas~\ref{lem:KlSobSitext}, \ref{lem:KlSobast} and their proof in Appendix~\ref{app:KlSobH12}, and the transport estimates of Lemma~\ref{lem:evolext}) relying on coordinate comparison with the Euclidean case. 
\end{remark}




\section{Norms}\label{sec:norms}
In this section, we define the norms of the curvature, connection coefficients and coordinates, upon which the bootstrap argument is constructed (see the Bootstrap Assumptions in Section~\ref{sec:BA}).

\subsection{Norms for the curvature in $\MM^\ext$ and $\protect\CCba$}\label{sec:normnullcurv}
\paragraph{Norms on $\CCba$.}
We define
\begin{align*}
  \qq := \min(r,u).
\end{align*}
\begin{remark}
  In the interior of the cone $\CCba\cap\MM^\intr$, we have
  \begin{align*}
    \qq & \simeq r.
  \end{align*}
  In the exterior of the cone $\CCba\cap\MM^\ext$, we have
  \begin{align*}
    \qq & \simeq u.
  \end{align*}
\end{remark}

We define
\begin{align*}
  \RRast[R] & := \norm{\Ndt^{\leq 2}R}_{L^2(\CCba)},
\end{align*}
where $\Ndt\in\le\{(r\Nd),(r\Nd_4), (\qq\Nd_3)\ri\}$. Using this notation, we have the following definition
\begin{align*}
  \RRast & := \RRast\le[u^2\alb\ri] + \RRast\le[u\ub\beb\ri] + \RRast\le[\ub^2(\rho-\rhoo)\ri] + \RRast\le[\ub^2(\sigma-\sigmao)\ri] \\
         & \quad + \le(\RRast\le[\ub^2\be\ri] + \norm{\ub^2\Ndt^{\leq 1}(r\Nd_3)\be}_{L^2(\CCba)} \ri) \\
         & \quad + \le(\RRast\le[\ub^2\al\ri] + \norm{\ub^2\Ndt^{\leq 1}(r\Nd_3)\al}_{L^2(\CCba)}  + \norm{\ub^2(r\Nd_3)^2\al}_{L^2(\CCba)} - \norm{(r\Nd_4)^2(\ub^2\al)}_{L^2(\CCba)}\ri).
\end{align*}
We define
\begin{align*}
  \RRfast[R] & := \norm{r^{1/2}\qq^{1/2}\Ndt^{\leq 1}R}_{L^\infty_u\HHt(S_{u,\uba})},
\end{align*}
and
\begin{align*}
  \RRfast & := \RRfast\le[u^2\alb\ri] + \RRfast\le[u\ub\beb\ri] + \RRfast\le[\ub^2(\rho-\rhoo)\ri] \\
  & \quad + \RRfast\le[\ub^2(\sigma-\sigmao)\ri] + \RRfast\le[\qq^{-1/2}\ub^{5/2}\be\ri] + \RRfast\le[\qq^{-1/2}\ub^{5/2}\al\ri].
\end{align*}
We define
\begin{align*}
  \RRfoast\le[\rhoo\ri] & := \norm{\ub^3\qq u\Ndt^{\leq 2}\rhoo}_{L^\infty(\CCba)},\\
  \RRfoast\le[\sigmao\ri] & := \norm{\ub^3\qq u \Ndt^{\leq 2}\sigmao}_{L^\infty(\CCba)},
\end{align*}
and
\begin{align*}
  \RRfoast & := \RRfoast\le[\rhoo\ri] + \RRfoast\le[\sigmao\ri].
\end{align*}

\begin{remark}
  Here, as in the rest of this section, we implicitly assume that the $L^\infty\HHt$ bootstrap bounds for a tensor $F$ come together with the corresponding $L^\infty L^4$ bootstrap bounds for $F$ (see the Sobolev embeddings from Lemma~\ref{lem:Sobsphere}), and that the $L^\infty\HHt$ bootstrap bounds for $(r\Nd)^{\leq 1}F$ come together with the corresponding $L^\infty L^\infty$ bound for $F$.
\end{remark}

\begin{remark}
  In the exterior of the cone $\CCba\cap\MM^\ext$ these norms provide the expected optimal decay rates for the spacetime curvature. In the interior of the cone $\CCba\cap\MM^\intr$, these norms provide a suboptimal control in terms of $r$, which is due to the degeneracy of the null decomposition when $r\to 0$. See also Remark~\ref{rem:notoptimalaxisdege}. These norms are only used to estimate the null metric and connection coefficients on $\CCba$. The (optimal) curvature norms used in the treatment of the interior region are presented in Section~\ref{sec:normscurvint}.   
\end{remark}

\paragraph{Norms in $\MM^\ext$.}
We define
\begin{align*}
  \RRext & := \norm{u^{-1/2-\ga} \ub^2\Ndt^{\leq 2}\al}_{L^2(\MM^\ext)} +  \norm{u^{-1/2-\ga} \ub^2\Ndt^{\leq 2}\be}_{L^2(\MM^\ext)} \\
         & \quad +  \norm{\ub^{-1/2-\ga} \ub^2\Ndt^{\leq 2}(\rho-\rhoo)}_{L^2(\MM^\ext)} +  \norm{\ub^{-1/2-\ga} \ub^2\Ndt^{\leq 2}(\sigma-\sigmao)}_{L^2(\MM^\ext)} \\
  & \quad +  \norm{\ub^{-1/2-\ga} \ub u\Ndt^{\leq 2}\beb}_{L^2(\MM^\ext)} +  \norm{\ub^{-1/2-\ga} u^2\Ndt^{\leq 2}\alb}_{L^2(\MM^\ext)},
\end{align*}
for all $\ga>0$ and where here $\Ndt\in\le\{(r\Nd), (\ub\Nd_4), (u\Nd_3)\ri\}$, since $q \simeq u$ in $\MM^\ext$. We define
\begin{align*}
  \RRfext & := \norm{\ub^{3}\Ndt^{\leq 1}\al}_{L^\infty_{u,\ub}\HHt(S_{u,\ub})} + \norm{\ub^{3}\Ndt^{\leq 1}\be}_{L^\infty_{u,\ub}\HHt(S_{u,\ub})} \\
          &\quad + \norm{\ub^{5/2}u^{1/2}\Ndt^{\leq 1}(\rho-\rhoo)}_{L^\infty_{u,\ub}\HHt(S_{u,\ub})} + \norm{\ub^{5/2}u^{1/2}\Ndt^{\leq 1}(\sigma-\sigmao)}_{L^\infty_{u,\ub}\HHt(S_{u,\ub})} \\
          & \quad + \norm{\ub^{3/2}u^{3/2}\Ndt^{\leq 1}\beb}_{L^\infty_{u,\ub}\HHt(S_{u,\ub})} + \norm{\ub^{1/2}u^{5/2}\Ndt^{\leq 1}\alb}_{L^\infty_{u,\ub}\HHt(S_{u,\ub})}.  
\end{align*}
We define
\begin{align*}
  \RRfoext & := \norm{\ub^3u^2\Ndt^{\leq 2}\rhoo}_{L^\infty(\MM^\ext)} + \norm{\ub^3u^2\Ndt^{\leq 2}\sigmao}_{L^\infty(\MM^\ext)}.
\end{align*}

\begin{remark}
  The mean values $\rhoo,\sigmao$ of the null curvature components $\rho,\sigma$ have stronger decay rates than the other curvature components. This is a consequence of the average Bianchi equations~(\ref{eq:Nd4rhoo}), (\ref{eq:Nd3rhoo}) and the vertex limit $r^3\rhoo \to 0$ when $r\to 0$, and of the structure equation~(\ref{eq:Curlze}) for $\sigma$. See Sections~\ref{sec:connestCCba} and~\ref{sec:connest}.
\end{remark}

\subsection{Norms for the curvature in $\MM^\intr$}\label{sec:normscurvint}
We have the following definitions of the curvature boundedness norms in the interior region
\begin{align*}
  \RR^\intr_{\leq 2} & := \sup_{\too \leq t \leq \tast}\norm{t^2 \le|(t\D)^{\leq 2}\R\ri|}_{L^2(\Si_t)},
\end{align*}
where we take the norms in the maximal frame (\emph{i.e.} $e_0=\Tf$). \\

We have the following definitions of the curvature decay norms in the interior region
\begin{align*}
  \RRf^\intr_{\leq 1} & := \sup_{\too \leq t \leq \tast} \le(\norm{t^{7/2}|\R|}_{L^\infty(\Si_t)} + \norm{t^{4}\le|\D\R\ri|}_{L^6(\Si_t)}\ri).
\end{align*}


\begin{remark}
  These norms only cover the bottom interior region $\MM^\intr_\bott$. This is not an issue for the global energy estimates in which these norms are used, since these estimates are performed in that same region $\MM^\intr_\bott$. The top interior region $\MM^\intr_\topp$ is contained in the domain of dependence of the last maximal hypersurface $\Si_{\tast}$ and is treated \emph{via} a local existence argument (see Section~\ref{sec:tipcurvest}). 
\end{remark}

\subsection{Norms for the null connection coefficients on the cone $\protect\CCba$}\label{sec:normnullconnCCba}
For an $S$-tangent tensor $\Ga$, we define the following $L^2(\CCba)$ based norms
\begin{align*}
  \OO^{\ast,\mathfrak{g}}_{\leq\ell}\le[\Ga\ri] & := \norm{r^{-1}\ub^2\Ndt^{\leq \ell}\Ga}_{L^2(\CCba)},\\
  \OO^{\ast,\mathfrak{b}}_{\leq\ell}\le[\Ga\ri] & := \norm{r^{-1}\ub u \Ndt^{\leq \ell}\Ga}_{L^2(\CCba)},
\end{align*}
where $\Ndt \in \le\{(r\Nd), (\qq\Nd_3)\ri\}$ and where $\ell \geq 0$. We also define the following norm for $\omb-\ombo$
\begin{align*}
  \OOastbadbad\le[\omb-\ombo\ri] & := \OOastbad\le[\omb-\ombo\ri] - \norm{r^{-1}\ub u (\qq\Nd_3)^3(\omb-\ombo)}_{L^2(\CCba)}.
\end{align*}
We define
\begin{align*}
  \OOast & := \OOastgood\le[\trchib-\trchibo\ri] + \OOastgood\le[\ze\ri] + \OOastgood\le[\trchi-\trchio\ri] + \OOastgood\le[\chih\ri] \\
  & \quad + \OOastbad\le[\,\chibh\,\ri] + \OOastbadbad\le[\omb-\ombo\ri].
\end{align*}
We define the following norms for the mean values of $\trchi$ and $\trchib$
\begin{align*}
  \OOoast\le[\trchio\ri] & := \norm{u\ub^3(\qq\Nd_3)^{\leq 2}\le(\trchio-\frac{2}{r}\ri)}_{L^\infty(\CCba)}, \\
  \OOoast\le[\trchibo\ri] & := \norm{u^2\ub^2(\qq\Nd_3)^{\leq 2}\le(\trchibo+\frac{2}{r}\ri)}_{L^\infty(\CCba)},
\end{align*}
and
\begin{align*}
  \OOoast := \OOoast\le[\trchio\ri] + \OOoast\le[\trchibo\ri].
\end{align*}
We define the following $L^\infty_u\HHt(S_{u,\uba})$ based norms
\begin{align*}
  \OOfastgood\le[\Ga\ri] & := \norm{r^{-1/2}\qq^{1/2}\ub^2\Ndt^{\leq 2}\Ga}_{L^\infty_u\HHt(S_{u,\uba})},\\
  \OOfastbad\le[\Ga\ri] & := \norm{r^{-1/2}\qq^{1/2}\ub u \Ndt^{\leq 2}\Ga}_{L^\infty_u\HHt(S_{u,\uba})},\\
  \OOfastbadbad\le[\omb-\ombo\ri] & := \OOfastbad\le[\omb-\ombo\ri] - \norm{r^{-1/2}\qq^{1/2}\ub u (\qq\Nd_3)^2(\omb-\ombo)}_{L^\infty_u\HHt(S_{u,\uba})},
\end{align*}  
and
\begin{align*}
  \OOfast & := \OOfastgood\le[\trchib-\trchibo\ri] + \OOfastgood\le[\ze\ri] + \OOfastgood\le[\trchi-\trchio\ri] + \OOfastgood\le[\chih\ri] \\
  & \quad + \OOfastbad\le[\,\chibh\,\ri] + \OOfastbadbad\le[\omb-\ombo\ri].
\end{align*}

\begin{remark}
  Together with the relations between the null connection coefficients from Lemmas~\ref{lem:relgeodnull} and \ref{lem:relCCba}, these norms provide a control for all the null connection coefficients.
\end{remark}

\begin{remark}
  These norms give optimal decay rates in the exterior of the cone $\CCba\cap\MM^\ext$. They give a suboptimal control in $r$ in the region $\CCba\cap\MM^\intr$. This is a consequence of the suboptimal control for the curvature. See Remark~\ref{rem:notoptimalaxisdege}. Except to establish their own control in Section~\ref{sec:connestCCba}, these norms are only used in the exterior region where they are optimal. 
\end{remark}


We also have the following norms for $\OOE$
\begin{align*}
  \OO_{\leq 3}^{\ast,\OOO} & := \OO^{\ast, \mathfrak{g}}_{\leq 2}\le[r^{-1}\Hrot\ri] + \OO^{\ast, \mathfrak{g}}_{\leq 1}\le[\POE\ri],
\end{align*}
where in that case the norms are restricted to the exterior region $\CCba\cap\MM^\ext$ where the exterior rotations are defined. See the definitions from Section~\ref{sec:defrotext}.

\subsection{Norms for the null connection coefficients in $\MM^\ext$}\label{sec:normsnullconn}
We recall that in the exterior region $\MM^\ext$, we have $r\simeq \ub$. For an $S$-tangent tensor $\Ga$, we define the following $\HHt$ norms
\begin{align*}
  \OOfgoodext[\Ga] & := \norm{r^{3/2}u^{1/2}(\Ndt)^{\leq 1}\Ga}_{L^\infty_{u,\ub}\HHt(S_{u,\ub})}, \\
  \OOfbadext[\Ga] & := \norm{r^{1/2}u^{3/2}(\Ndt)^{\leq 1}\Ga}_{L^\infty_{u,\ub}\HHt(S_{u,\ub})},
\end{align*}
where $\Ndt\in\le\{(r\Nd), (u\Nd_3), (\ub\Nd_4)\ri\}$.\\

For an $S$-tangent tensor $\Ga$, we define the following $L^2(\MM^\ext)$ norms
\begin{align*}
  \OOgoodext[\Ga] & := \norm{\ub^{-1/2-\ga}\ub(\Ndt)^{\leq 2}\Ga}_{L^2(\MM^\ext)},\\
  \OObadext[\Ga] & := \norm{\ub^{-1/2-\ga}u(\Ndt)^{\leq 2}\Ga}_{L^2(\MM^\ext)},
\end{align*}
for all $\ga>0$ and where $\Ndt\in\le\{(r\Nd), (u\Nd_3), (\ub\Nd_4)\ri\}$.\\

We define
\begin{align*}
  \OOfexti & := \OOfgoodext\le[\trchi-\trchio\ri] + \OOfgoodext\le[\chih\ri] + \OOfgoodext\le[\trchib-\trchibo\ri] + \OOfgoodext\le[\ze\ri] \\
  & \quad + \OOfbadext\le[\,\chibh\,\ri] + \OOfbadext\le[\omb-\ombo,\ombd\ri] + \OOfbadext\le[\,\xib\,\ri],\\ \\
  \OOfextii & := \OOfexti + \OOfgoodext\le[(r\Nd)(\trchi-\trchio)\ri] + \OOfgoodext\le[(r\Nd)\chih\ri] + \OOfgoodext\le[(r\Nd)\ze\ri] \\
  & \quad + \OOfgoodext\le[(u\Nd_3)\ze\ri] + \OOfbadext\le[(r\Nd)\ombr,(r\Nd)\ombs\ri] + \OOfbadext\le[(r\Nd)\omb, (r\Nd)\ombd\ri],
\end{align*}
where we recall that $\ombd$, $\ombr,\ombs$ are defined in Section~\ref{sec:nulldecompMMext}.\\

We define
\begin{align*}
 \OOexti & := \OOgoodext\le[\trchi-\trchio\ri] + \OOgoodext\le[\chih\ri] + \OOgoodext\le[\trchib-\trchibo\ri] + \OOgoodext\le[\ze\ri] \\
  & \quad + \OObadext\le[\chibh\ri] + \OObadext\le[\omb-\ombo,\ombd\ri] + \OObadext\le[\xib\ri],\\ \\
  \OOextii & := \OOexti + \OOgoodext\le[(r\Nd)(\trchi-\trchio)\ri] + \OOgoodext\le[(r\Nd)\chih\ri] + \OOgoodext\le[(r\Nd)\ze\ri] \\
  & \quad + \OOgoodext\le[(u\Nd_3)\ze\ri] + \OObadext\le[(r\Nd)\ombr,(r\Nd)\ombs\ri] + \OObadext\le[(r\Nd)\omb, (r\Nd)\ombd\ri]. 
\end{align*}

We define
\begin{align*}
  \OOofb[\ombo] & := \norm{\ub^2u^2 \Ndt^{\leq 1}\ombo}_{L^\infty(\MM^\ext)}, \\
  \OOofb[\trchio] & := \norm{\ub^{3}u^{1}\Ndt^{\leq 1}\le(\trchio-\frac{2}{r}\ri)}_{L^\infty(\MM^\ext)},\\
  \OOofb[\trchibo] & := \norm{\ub^2u^2\Ndt^{\leq 1}\le(\trchibo+\frac{2}{r}\ri)}_{L^\infty(\MM^\ext)},
\end{align*}
where $\Ndt\in\le\{(r\Nd), (u\Nd_3), (\ub\Nd_4)\ri\}$, and we note
\begin{align*}
  \OOofb & := \OOofb[\ombo] + \OOofb[\trchio] + \OOofb[\trchibo].
\end{align*}

We define
\begin{align*}
  \OOof[\ombo] & := \norm{\ub^{-1/2-\ga}u^{3/2} \ub \Ndt^{\leq 2}\ombo}_{L^2(\MM^\ext)},\\
  \OOof[\trchio] & := \norm{\ub^{-1/2-\ga} u^{1/2}\ub^2 \Ndt^{\leq 2}\le(\trchio-\frac{2}{r}\ri)}_{L^2(\MM^\ext)},\\
  \OOof[\trchibo] & := \norm{\ub^{-1/2-\ga}u^{3/2} \ub \Ndt^{\leq 2}\le(\trchibo+\frac{2}{r}\ri)}_{L^2(\MM^\ext)},
\end{align*}
for all $\ga>0$ and where $\Ndt\in\le\{(r\Nd), (u\Nd_3), (\ub\Nd_4)\ri\}$.\\

\begin{remark}
  The absence of control for the higher derivatives of the coefficients $\chib$ and $\xib$ is due to the classical loss of regularity for the geodesic foliation. See Section~\ref{sec:connest}.
\end{remark}

We have the following definitions for the norms of the defect $\yy$ in $\MM^\ext$ and on $\TT$
\begin{align*}
  \OOfext[\yy] & := \norm{\ub u^2\Ndt^{\leq 1}\yyo}_{L^\infty(\MM^\ext)} + \norm{r^{-1/2}u^{3/2} \Ndt^{\leq 1}(r\Nd)\yy}_{L^\infty_{u,\ub}\HHt},\\
  \OOext[\yy] & := \norm{\ub^{-1/2-\ga}u^{3/2}\Ndt^{\leq 2}\yyo}_{L^2(\MM^\ext)} + \norm{\ub^{-1/2-\ga}\ub^{-1} u \Ndt^{\leq 2}(r\Nd)\yy}_{L^2(\MM^\ext)},\\
  \OO^\TT_{\leq 2,\ga}[\yy] & := \norm{t^{3/2-\ga}\Ndt^{\leq 2}\yyo}_{L^2(\TT)} + \norm{t^{-\ga}\Ndt^{\leq 2}(r\Nd)\yy}_{L^2(\TT)}, 
\end{align*}
for all $\ga>0$ and where $\Ndt\in\le\{(r\Nd), (u\Nd_3), (\ub\Nd_4)\ri\}$.\\

We define the following norms for $\OOE$ in $\MM^\ext $
\begin{align*}
  \mathfrak{O}^{\ext,\OOO}_{\leq 2} & := \OOfgoodext\le[r^{-1}\Hrot\ri] + \OOfgoodext\le[r^{-1}Y\ri] + \mathfrak{O}^{\ext,\mathfrak{g}}_{\leq 0}\le[\POE\ri],\\
  \OO^{\ext,\OOO}_{\leq 3, \ga} & := \OOgoodext\le[r^{-1}\Hrot\ri] + \OOgoodext\le[r^{-1}Y\ri] + \OO^{\ext,\mathfrak{g}}_{\leq 1,\ga}\le[\POE\ri],
\end{align*}
for all $\ga>0$.\\

\subsection{Norms for the maximal connection coefficients in $\MM^\intr_\bott$}\label{sec:norminteriorconn}

We have the following definitions for the norms for the time lapse in $\MM^\intr_\bott$
\begin{align*}
  \OO^\intr_{\leq 3, \ga}[\nt] & := \norm{(t\nab)^{\leq 3}(\nt-1)}_{L^\infty_tL^2(\Si_t)} + \norm{t(t\nab)^{\leq 2}\Tf(\nt-1)}_{L^\infty_t L^2(\Si_t)} \\
  & \quad + \norm{t^{-1/2-\ga}t^2(t\nab)^{\leq 1}\Tf^2(\nt-1)}_{L^2(\MM^\intr_\bott)},
\end{align*}
for all $\ga>0$.\\

We have the following definitions for the boundedness norms for the plane second fundamental form in $\MM^\intr$
\begin{align*}
  \OO^\intr_{\leq 2}\le[\kt\ri] & := \norm{t\le(t\nab, t\Lieh_\Tf\ri)^{\leq 2}\kt}_{L^\infty_tL^2(\Si_t)}.
\end{align*}

Additionally, we have the following definition for the transition factors $\nut$ at $\TT$
\begin{align*}
  \mathfrak{O}^\TT_{\leq 2}\le[\nut\ri] & := \norm{t\le(t\Nd,tZ\ri)^{\leq 2}\le(\nut-1\ri)}_{L^\infty_t \tilde{H}^{1/2}\le(\pr\Si_t\ri)},
\end{align*}
where $Z$ is the future-pointing unit normal to $\pr\Si_t$ in $\TT$.\\

Moreover, we control the following more regular norms on $\Si_\tast$
\begin{align*}
  \OO^{\Si_\tast}_{\leq 3}[k] & := \norm{t(t\nab)^3k}_{L^2(\Si_\tast)}. 
\end{align*}

\subsection{Norms for the approximate Killing fields in $\MM^\intr_\bott$}\label{sec:normintkill}

We have the following definitions for the control of the derivatives and deformation tensors of $\TI, \XI, \SI,\KI,\OOI$
\begin{align*}
  \OO_{\leq 3, \ga}^\intr\le[\TI\ri] & := \norm{t^{5/2}\D\TI}_{L^\infty(\MM^\intr_\bott)} \\
                                     & \quad + \norm{t^2(t\D)^{\leq 1}\D\TI}_{L^\infty_t L^6(\Si_t)} \\
                                     & \quad + \norm{t^{-1/2-\ga}t(t\D)^{\leq 2}\D\TI}_{L^2(\MM^\intr_\bott)},\\ \\
  \OO_{\leq 3, \ga}^\intr\le[\XI\ri] & := \norm{t^{3/2}(\D\XI-g)}_{L^\infty(\MM^\intr_\bott)} \\
                                     & \quad + \norm{t(t\D)^{\leq 1}(\D\XI-g)}_{L^\infty_tL^6(\Si_t)} \\
                                     & \quad + \norm{t^{-1/2-\ga}(t\D)^{\leq 2}(\D\XI-g)}_{L^2(\MM^\intr_\bott)} \\ \\
  \OO_{\leq 3,\ga}^\intr\le[\SI\ri] & := \norm{t^{3/2}(\D\SI-\g)}_{L^\infty(\MM^\intr_\bott)} \\
                                     & \quad + \norm{t(t\D)^{\leq 1}(\D\SI-\g)}_{L^\infty_tL^6(\Si_t)} \\
                                     & \quad + \norm{t^{-1/2-\ga}(t\D)^{\leq 2}(\D\SI-\g)}_{L^2(\MM^\intr_\bott)} \\ \\
  \OO_{\leq 3,\ga}^\intr\le[\KI\ri] & := \norm{t^{1/2}\le({^{(\KI)}\pi}-4t\g\ri)}_{L^\infty(\MM^\intr_\bott)} \\
                                     & \quad + \norm{(t\D)^{\leq 1}\le({^{(\KI)}\pi}-4t\g\ri)}_{L^\infty_tL^6(\Si_t)} \\
                                     & \quad + \norm{t^{-1/2-\ga}t^{-1}(t\D)^{\leq 2}\le({^{(\KI)}\pi}-4t\g\ri)}_{L^2(\MM^\intr_\bott)},
\end{align*}
and
\begin{align*}
  \OO_{\leq 3, \ga}^\intr\le[\OOI\ri] & := \norm{t^{3/2}{^{(\OOI)}\pi}}_{L^\infty(\MM^\intr_\bott)} \\
                                      & \quad + \norm{t^2 \D^2\OOI}_{L^\infty_tL^6(\Si_t)} \\
                                      & \quad + \norm{t^{-1/2-\ga} t(t\D)^{\leq 1}\D^2\OOI}_{L^2(\MM^\intr_\bott)},
\end{align*}
for all $\ga>0$.


\section{The Bootstrap Assumptions}\label{sec:BA}
In this section, we collect the \emph{Bootstrap Assumptions} which are used throughout this paper. See Section~\ref{sec:proofmainthm} for the description of the associated bootstrap argument.

\subsection{The constants used in this paper}\label{sec:constantsofthepaper}
In this section, we recapitulate the constants used in this paper, and in particular the constants used in the following bootstrap assumptions. All these constants are \emph{independent of the smallness parameter $\varep$}.
\begin{itemize}
\item The constant $\ga_0$ is used in the Bootstrap Assumptions of Sections~\ref{sec:strongBA} and is a fixed constant such that $0 < \ga_0 < 1/4$. The bootstrap assumptions involving $\ga_0$ are improved for all $\ga>0$, and in particular for $\ga=\gao$.
\item The constant $C>0$ is used in the mild Bootstrap Assumptions of Section~\ref{sec:mildBA} and is a fixed large constant.
\item The constant $D>0$ is used in the strong Bootstrap Assumptions of Section~\ref{sec:strongBA} and is a fixed large constant.
\item The constant $0<\wp<1$ is used in the elliptic estimates of Lemma~\ref{lem:ellHodgevar} is a small constant, which depends on $C$.
\item The transition constant $0<\cc_0<1$ is used to determine the timelike transition hypersurface. It is a small constant, which depends on $\wp$ (see Section~\ref{sec:controlRRintr2}). The transition parameter $\cc$ is allowed to range in $[\cco,(1+\cco)/2]$. This freedom is used in the mean value argument of Section~\ref{sec:meanvalue}.
\end{itemize}

\begin{remark}
  The bootstrap assumptions of Sections~\ref{sec:mildBA} and~\ref{sec:strongBA} will be assumed to hold uniformly \emph{for all} transition parameter $\cco \leq \cc \leq (1+\cco)/2$. See Section~\ref{sec:proofmainthm}.
\end{remark}

\subsection{Mild bootstrap assumptions}\label{sec:mildBA}
  


\begin{BA}[Mild bootstrap assumptions for the rotation vectorfield in $\MM^\ext$]\label{BA:mildOOE}
  Let $C>0$ be a (large) numerical constant. We assume that 
  \begin{align*}
    \norm{r^{-1}\OOE}_{L^\infty(\MM^\ext)} + \norm{\Nd\OOE}_{L^\infty(\MM^\ext)} & \leq C.
  \end{align*}
  Moreover, we assume that for all $S$-tangent scalar $f$ and for all $1$-tensor or symmetric traceless $2$-tensor $F$, the following bound holds
  \begin{align*}
    \int_{S_{u,\ub}} \le|(r\Nd) f\ri|^2 & \leq C \sum_{\ell=1}^3\int_{S_{u,\ub}}\le|\OOEi(f)\ri|^2, \\
    \int_{S_{u,\ub}} \le|(r\Nd)^{\leq 1} F\ri|^2 & \leq C \sum_{\ell=1}^3\int_{S_{u,\ub}}\le|\Liedh_\OOEi F\ri|^2,
  \end{align*}
  on all $2$-sphere $S_{u,\ub}\subset \MM^\ext$.
\end{BA}

\begin{BA}[Mild bootstrap assumptions for the maximal hypersurfaces $\Si_t \subset \MM^\intr_\bott$]\label{BA:mildcoordsSit}
  Let $C>0$ be a (large) numerical constant. We assume that on each separate maximal hypersurface $\Si_t$, for $\too\leq t \leq \tast$, there exists (harmonic) global coordinates $(x^i)$ such that
  \begin{align*}
    \le\{\le(x^1\ri)^2 + \le(x^2\ri)^2 + \le(x^3\ri)^2 = \le(\frac{1-\cc}{1+\cc}\ri)^2t^2\ri\} & = \pr\Si_t, \\
    \le\{\le(x^1\ri)^2 + \le(x^2\ri)^2 + \le(x^3\ri)^2 \leq \le(\frac{1-\cc}{1+\cc}\ri)^2 t^2\ri\} & = \Si_t,
  \end{align*}
  and such that we have the following uniform bounds for the metric $g$ in these coordinates
  \begin{align*}
    \le|g_{ij}-\de_{ij}\ri| & < 1/4,\\
    \le|\pr_kg_{ij}\ri| & \leq C,
  \end{align*}
  on each separate maximal hypersurface $\Si_t$.
\end{BA}
\begin{BA}[Mild bootstrap assumptions for the Killing fields in $\MM^\intr_\bott$]\label{BA:mildKillingMMintbot}
  Let $C>0$ be a (large) numerical constant. We assume that the norm of the vectorfields $\TI, \SI$, $\KI, \OOI$ satisfy the following mild bounds in $\MM^\intr_\bott$
  \begin{align}\label{est:mildnormSIKIOOI}
    \begin{aligned}
      |\XI| & \leq C t, & |\SI| & \leq C t, & |\KI| & \leq C t^2, & |\OOI| & \leq C t, & |\D\OOI| & \leq C,  
    \end{aligned}
  \end{align}
  where the norms are taken with respect to the maximal frame. We also assume that $\KI$ is a future-pointing timelike vectorfield in $\MM^\intr_\bott$ and that we have the following mild bounds 
  \begin{align}\label{est:mildtimelikeKI}
    \begin{aligned}
      \g(\KI,\TI) & \leq - C^{-1} t^2, & \le|\KI + \g\le(\KI,\TI\ri)\TI\ri| & \leq \le(1-C^{-1}\ri) |\g(\le(\KI,\TI\ri)|.
    \end{aligned}
  \end{align}
\end{BA}

\subsection{Strong bootstrap assumptions}\label{sec:strongBA}
Let $0 < \gao < 1/4$ be a fixed numerical constant.\footnote{Fixing the constant $\gao$ breaks the scaling in the norms below. However, these norms are only used in the control of the nonlinear error terms of Section~\ref{sec:globener}, which do not require a sharp control. Imposing $\gao<1/4$ is sufficient for this analysis. See Section~\ref{sec:globener}.} We have the following strong bootstrap assumptions in $\MM$.
\begin{BA}[Spacetime curvature in $\CCba$]\label{BA:curvast}
  We assume that on $\CCba$
  \begin{align*}
    \RRast + \RRfast + \RRfoast & \leq D\varep,
  \end{align*}
  and we refer to Section~\ref{sec:normnullcurv} for definitions.
\end{BA}
\begin{BA}[Spacetime curvature in $\MM^\ext$]\label{BA:curvext}
  We assume that on $\MM^\ext$
  \begin{align*}
    \RR^\ext_{\leq 2,\gao} + \mathfrak{R}^\ext_{\leq 1} + \RRfoext & \leq D\varep,
  \end{align*}
  and we refer to Section~\ref{sec:normnullcurv} for definitions.
\end{BA}
\begin{BA}[Spacetime curvature in $\MM^\intr_\bott$]\label{BA:curvint}
  We assume that on $\MM^\intr_\bott$
  \begin{align*}
    \RR^{\intr}_{\leq 2} + \RRf^\intr_{\leq 1} & \leq D\varep,
  \end{align*}
  and we refer to Section~\ref{sec:normscurvint} for definitions.
\end{BA}
\begin{BA}[Null connection in $\CCba$]\label{BA:connCCba}
  We assume that on $\CCba$, we have
  \begin{align*}
    \OO^\ast_{\leq 3} + \mathfrak{O}^\ast_{\leq 2} + \overline{\mathfrak{O}}_{\leq 2} + \OO^{\ast,\OOO}_{\leq 3} & \leq D\varep,
  \end{align*}
  and we refer to Section~\ref{sec:normnullconnCCba} for definitions. Moreover, we assume that we have the following bootstrap bound for $r$ on $\CCba$ 
  \begin{align}
    \label{est:BAarearadiusestimateCCba}
    \le|r-\half(\uba-u)\ri| & \leq D\varep \uba^{-3}r^2 u^{-1}.
  \end{align}
\end{BA}
\begin{BA}[Null connection in $\MM^\ext$]\label{BA:connext}
  We assume that on $\MM^\ext$, we have
  \begin{align*}
    \mathfrak{O}^\ext_{\leq 2} + \OO^\ext_{\leq 3, \gao} + \overline{\OO}^\ext_{2,\gao} & \leq D\varep,\\
    \mathfrak{O}^\ext_{\leq 1}[\yy] + \OO^\ext_{2,\gao}[\yy] + \OO^\TT_{\leq 2,\gao}[\yy] & \leq D\varep,\\
    \mathfrak{O}^{\ext,\OOO}_{\leq 2} + \OO^{\ext,\OOO}_{\leq 3,\gao} & \leq D\varep, 
  \end{align*}
  where we refer to Section~\ref{sec:normsnullconn} for definitions. Moreover, we assume that we have the following bootstrap bound for $r$ on $\MM^\ext$
  \begin{align}
    \label{est:BAarearadiusestimate}
    \le|r-\half(\ub-u)\ri| & \leq D\varep \ub^{-1}u^{-1}.
  \end{align}
\end{BA}
\begin{BA}[Maximal connection in $\MM^\intr_\bott$]\label{BA:connint}
  We assume that on $\MM^\intr_\bott$, we have
  \begin{align*}
    \OO^\intr_{\leq 3,\gao}[\nt] + \OO^\intr_{\leq 2}[k] + \mathfrak{O}^\TT_{\leq 2}[\nu] + \OO^{\Si_\tast}_{\leq 3}[k] & \leq D\varep,
  \end{align*}
  and we refer to Section~\ref{sec:norminteriorconn} for definitions.
\end{BA}

\begin{BA}[Spherical coordinates in $\CCba$]\label{BA:mildsphcoordsast}
  We assume that there exists (two) spherical coordinate systems $(u,\varth,\varphi)$ covering $\CCba\setminus\o(\uba)$ and ranging into 
  \begin{align*}
    [1,\uba]_u\times[\pi/8,7\pi/8]_\varth\times[0,2\pi)_\varphi,
  \end{align*}
  where $u$ is the canonical parameter on $\CCba$ as defined in Section~\ref{sec:definition}.\footnote{\label{foo:sphcoords} The (two) spherical coordinates patch $\varth,\varphi$ correspond to (two) different axis on the sphere. The choice of ranges for $\varth,\varphi$ ensures that two such patches cover the full $2$-spheres $S_{u,\ub}$.} Moreover, we assume that for the induced metric on the $2$-spheres $S_{u,\uba}$ in coordinates $\varth,\varphi$, we have the following bounds
  \begin{align*}
    \le|\pr_a^{\leq 1}\le(r^{-2}\gd_{bc}-(\gd_{\SSS})_{bc}\ri)\ri| & \leq D\varep \uba^{-1}u^{-1}q^{1/2},
  \end{align*}
  where $a,b,c \in\{\varth,\varphi\}$ and
  \begin{align*}
    \gd_{\SSS} & := \d\varth^2 + \sin^2\varth\d\varphi^2.
  \end{align*}
\end{BA}

\begin{BA}[Spherical coordinates in $\MM^\ext$]\label{BA:mildsphcoordsext}
  We assume that there exists (two) spherical coordinate systems $(u,\ub,\varth,\varphi)$ covering $\MM^\ext$ and ranging into 
  \begin{align*}
    [1,\cc\uba]_u\times[\cc^{-1},\uba]_\ub\times[\pi/8,7\pi/8]_\varth\times[0,2\pi)_\varphi,
  \end{align*}
  where $u$ and $\ub$ are respectively the optical and affine parameter functions defined in Section~\ref{sec:definition}, and which coincide with the spherical coordinates system $(u,\varth,\varphi)$ on $\CCba$ from the Bootstrap Assumptions~\ref{BA:mildsphcoordsast}.\footnoteref{foo:sphcoords} Moreover, we assume that for the induced metric on the $2$-spheres $S_{u,\ub}$ in coordinates $\varth,\varphi$, we have the following bounds
  \begin{align*}
    \le|\pr_a^{\leq 1}\le(r^{-2}\gd_{bc}-(\gd_{\SSS})_{bc}\ri)\ri| & \leq D\varep \ub^{-1}u^{-1/2},
  \end{align*}
  where $a,b,c \in\{\varth,\varphi\}$ and
  \begin{align*}
    \gd_{\SSS} & := \d\varth^2 + \sin^2\varth\d\varphi^2.
  \end{align*}
\end{BA}
\begin{BA}[Harmonic coordinates in $\Si_\tast$]\label{BA:harmoSitast}
  We assume that on $\Si_\tast$, we have
  \begin{align*}
      \sum_{i,j=1}^3 \norm{t^{3/2}\le(g(\nab x^i,\nab x^j)-\de_{ij}\ri)}_{L^\infty(\Si_\tast)} + \sum_{i=1}^3 \norm{t(t\nab)^{\leq 3}\nab^2x^i}_{L^2(\Si_\tast)} & \leq D\varep,\\
  \sum_{i=1}^3\norm{t^{1/2}\le(\rast\Nf(x^i)-x^i\ri)}_{L^\infty(S^\ast)} & \leq D\varep,
  \end{align*}
  where we refer to Section~\ref{sec:defUnifandharmo} for the definition of the harmonic coordinates $x^i$.
\end{BA}
\begin{BA}[Interior Killing vectorfields in $\MM^\intr_\bott$]\label{BA:intKill}
  We assume that
  \begin{align}\label{est:BADX}
    \OO^{\intr}_{\leq 3, \gao}\le[\TI,\XI,\SI,\KI,\OOI\ri] & \leq D\varep,
  \end{align}
  and we refer to Section~\ref{sec:normintkill} for definitions. We also assume that
  \begin{align}\label{est:BAXIXI}
    \sup_{\too \leq t \leq \tast}t^{3/2}\le|\sup_{\Si_t}\le(t^{-2}\g\le(\XI,\XI\ri)\ri) - \le(\frac{1-\cc}{1+\cc}\ri)^2\ri| & \leq D\varep,
  \end{align}
  and
  \begin{align}\label{est:BAXITI}
    \sup_{\too \leq t \leq \tast}t^{3/2}\le|(t\D)^{\leq 1}\g(t^{-1}\XI,\TI)\ri| & \leq (D\varep).
  \end{align}
\end{BA}
\begin{BA}[Killing fields at $\TT$]\label{BA:TTKilling}
  We assume that at the interface $\TT$ the following bootstrap bounds hold
  \begin{align}\label{est:XXEXXITT}
    \begin{aligned}
      |\TE-\TI| & \leq D\varep t^{-3/2},\\
      |\SE-\SI| & \leq D\varep t^{-1/2},\\
      |\KE-\KI| & \leq D\varep t^{1/2},\\
      |\OOE-\OOI| & \leq D\varep t^{-1/2}.
    \end{aligned}
  \end{align}
  Moreover, we assume that for the first order derivatives, we have
  \begin{align}\label{est:piXXEXXITT}
    \begin{aligned}
      \norm{t^{-\gao}t^2(t\D)^{\leq 1}\le(\D\TE-\D\TI\ri)}_{L^\infty_tL^4(\pr\Si_t)} & \les \varep, \\
      \norm{t^{-\gao}t(t\D)^{\leq 1}(t\D)^{\leq 1}\le(\D\SE-\D\SI\ri)}_{L^\infty_tL^4(\pr\Si_t)} & \les \varep, \\
      \norm{t^{-\gao}t(t\D)^{\leq 1}\le(\D\OOE-\D\OOI\ri)}_{L^\infty_tL^4(\pr\Si_t)} & \les \varep.
    \end{aligned}
  \end{align}
\end{BA}

\begin{BA}[Last cones geodesic foliation]\label{BA:lastconesfoliation}
  We assume that $\MM' \subset \MM$ and that on the cone $\MM'\cap\CCba$ the following bounds hold
  \begin{align*}
    \le|u-u'\ri| & \leq D\varep \uba^{-3/2}(\uba-u'),
  \end{align*}
  and
  \begin{align*}
    \le|f',\fb',\log{\la'}\ri| & \leq D\varep \uba^{-1}{u'}^{-1/2}, 
  \end{align*}
  and that in the exterior region $\MM'\cap\MM^\ext$ the following bounds hold
  \begin{align*}
    \le|u-{u'}\ri| & \leq D\varep {\ub'}^{-1/2}, & \le|\ub-{\ub'}\ri| & \leq D\varep u'^{-1/2},
  \end{align*}
  and
  \begin{align*}
    \le|f',\log\la'\ri| & \leq D\varep {\ub'}^{-1}{u'}^{-1/2}, & \le|\fb'\ri| & \leq D\varep {u'}^{-3/2},
  \end{align*}
  where $f',\fb',{\la'}$ denote the transition coefficients between the null pairs $(\elb,\el)$ and $(\elb',\el')$, \emph{i.e.} the $S_{u,\ub}$-tangent tensors such that\footnote{To avoid confusion between the transition coefficients for the different change of frames in this paper, we renamed $(f',\fb',\la')$ the transition coefficients $(f,\fb,\la)$ of Section~\ref{sec:defchgframe}.}
  \begin{align*}
    \el' & = \la'\le(\el + f' + \quar |f'|^2\elb\ri),\\
    \elb' & = (\la')^{-1} \le(\le(1+\half f'\cdot\fb' + \frac{1}{16}|f'|^2|\fb'|^2\ri)\elb + \fb' + \quar|\fb'|^2f' + \quar |\fb'|^2\el\ri).
  \end{align*}
  Moreover, we assume that in the interior region $\MM'\cap\MM^\intr_\bott$, the following bounds hold
  \begin{align*}
    \le|t-\half({u'}+{\ub'})\ri| & \leq D\varep {\ub'}^{-1/2},
  \end{align*}
  and
  \begin{align*}
    \le|\g\le(\Tf,\half(\elb'+\el')\ri)+1\ri| & \leq D\varep {\ub'}^{-3/2}.
  \end{align*}
\end{BA}
\begin{BA}[Bottom initial layer]\label{BA:bottom}
  We assume that in the spacetimes regions $\LLbext$ and $\LLbint$, the frames of the respective regions are comparable, \emph{i.e.}
  \begin{align}\label{est:BAcompMMintLLb}
    \le|\g\le(\half(\elb+\el),\Tf^\bott\ri)+1\ri| & \leq D\varep, & \le|\g(\Tf,\Tf^\bott)+1\ri| & \leq D\varep.
  \end{align}
  We moreover assume that in the respective regions $\LLbext$ and $\LLbint$ the following comparisons between time functions hold
  \begin{align}
    \le|x^0 - \half(u+\ub)\ri| & \leq D\varep, & \le|x^0-t\ri| & \leq D\varep.
  \end{align}
\end{BA}
\begin{BA}[Conical initial layer]\label{BA:con}
  We assume that in the spacetime region $\LLcext$, we have
  \begin{align}\label{est:BAuuubub}
    \begin{aligned}
      \le|u-\ut\ri| & \leq D\varep, \\
      \le|\ub-\ubt\ri| & \leq D\varep \ub.
    \end{aligned}
  \end{align}
  Let $(\lat,\ft,\fbt)$ be the transition coefficients associated to the null pairs $(\elbt,\elt)$ and $(\elb,\el)$, \emph{i.e.} the $S_{u,\ub}$-tangent tensors such that
  \begin{align*}
    \elt & = \lat\le(\el + \ft + \quar |\ft|^2\elbt\ri),\\
    \elbt & = (\lat)^{-1} \le(\le(1+\half \ft\cdot\fbt + \frac{1}{16}|\ft|^2|\fbt|^2\ri)\elb + \fbt + \quar|\fbt|^2\ft + \quar |\fbt|^2\el\ri).
  \end{align*}
  We assume that in the region $\LLcext$, the following bounds hold
  \begin{align}\label{est:BAcompnullframe}
    \begin{aligned}
      |\ft| & \leq D\varep \ub^{-1}, & |\log\lat| & \leq D\varep, & |\fbt| & \leq D\varep.
      \end{aligned}
  \end{align}
\end{BA}

\subsection{First consequences of the Bootstrap Assumptions}
In this section, we collect lemmas which follow from the mild and strong Bootstrap Assumptions and which are used throughout Sections~\ref{sec:globener}~--~\ref{sec:initlayer}.

\begin{remark}
  Here and in the following, we write
  \begin{align*}
    f_1 & \les f_2,
  \end{align*}
  if there exists a constant $C'>0$, independent of $\varep$, such that
  \begin{align*}
    f_1 \leq C' f_2.
  \end{align*}
  We write
  \begin{align*}
    f_1 \simeq f_2
  \end{align*}
  if
  \begin{align*}
    f_1 \les f_2  && \text{and} && f_2 \les f_1.
  \end{align*}
\end{remark}

\begin{lemma}[Coarea formulas]\label{lem:coarea}
  In this paper, we use the following coarea formulas
  \begin{align}\label{est:coareaMMint}
    \int_{\MM^\intr_\bott}f & \simeq \int_{\too}^\tast \int_{\Si_t}f \,\d t,
  \end{align}
  and
  \begin{align}\label{est:coareaMMextSitext}
    \begin{aligned}
      \int_{\MM^\ext}f & \simeq \int_{\too}^\tast \int_{\Si^\ext_t}f \,\d t, & \int_{\Si^\ext_t}f & \simeq \int_{2t/(1+\cc^{-1})}^{2t-\uba} \int_{S_{u,\ub = 2t-u}}f \,\d u,
    \end{aligned}
  \end{align}
  and
  \begin{align}\label{est:coareaMMext}
    \begin{aligned}
      \int_{\MM^\ext}f & \simeq \int_{\cc^{-1}}^\uba\int_{1}^{\cc\ub} \int_{S_{u,\ub}} f \,\d u\d\ub & & = \int_{1}^{\cc\uba}\int_{\cc^{-1}u}^\uba \int_{S_{u,\ub}} f \,\d\ub \d u
      =: \int_{1}^{\cc\uba}\int_{\CC_u} f \,\d u,
    \end{aligned}
  \end{align}
  for all scalar function $f$.\footnote{The constants in the above estimates do not depend on $f$.}
\end{lemma}
\begin{proof}
  Formula~\eqref{est:coareaMMint} follows from the expression of $\int_{\MM^\intr_\bott}$ in coordinates
  \begin{align*}
    \int_{\MM^\intr_\bott} f  & = \int_{\too}^\tast \int_{\Si_t} f n \,\d t
  \end{align*}
  and the bound $|n-1|\leq D\varep$ from the Bootstrap Assumptions~\ref{BA:connint}. Formulas~\eqref{est:coareaMMextSitext} follow similarly from the following control of the exterior time lapse $n^\ext$ 
  \begin{align*}
    |n^\ext-1| & \les |\yy| \les D\varep,
  \end{align*}
  and from
  \begin{align*}
    \le|\Nf^\ext(u)+1\ri| & \les |\yy| \les D\varep 
  \end{align*}
  where we used relations~\eqref{eq:timelapseext} and~(\ref{eq:relTfNfextelelb}) and the Bootstrap Assumptions~\ref{BA:connext} for $\yy$.\\

  For the spherical coordinates from the Bootstrap Assumptions~\ref{BA:mildsphcoordsext}, we have the following expression of the coordinate vectorfields $\pr_u,\pr_\ub$
  \begin{align}\label{eq:pruprubbshift}
    \pr_u & = \half \elb-\quar\yy\el + b^a\ea, &   \pr_\ub & = \half\el,
  \end{align}
  which is a consequence of relations~\eqref{eq:elu}. Thus, the spacetime metric $\g$ writes in these coordinates
  \begin{align*}
    \g & = -\half \d u\d\ub - \half\d\ub\d u + \le(\half\yy + |b|^2\ri)\d u^2 +  b^{\flat}\d u +\d u b^{\flat}   + \gd
  \end{align*}
  as an element of $\mathrm{T}^\ast\MM \otimes \mathrm{T}^\ast\MM$, and where $\flat$ is the canonical musical isomorphism.\\
  
  In $\MM^\ext$, using equation~\eqref{eq:pruprubbshift} and relations~\eqref{eq:Riccirel} we have the following formula
  \begin{align*}
    \Nd_4b_a & = b^c\chi_{ca} + 2\ze_a. 
  \end{align*}
  From a straight-forward Gr\"onwall argument, using that $b=0$ on $\CCba$ and the Bootstrap Assumptions~\ref{BA:connext} for $\chi$ and $\ze$ we obtain
  \begin{align}\label{est:bshiftcoords}
    |b| & \les D\varep \ub^{-1}u^{-1/2} 
  \end{align}
  in $\MM^\ext$.\\

  Expressing $\int_{\MM^\ext}$ in coordinates, using the Bootstrap Assumptions~\ref{BA:connext} for $\yy$ and estimate~\eqref{est:bshiftcoords} for $b$ then yields the desired formulas~\eqref{est:coareaMMext}.
\end{proof}

We have the following lemma, which follows from the Bootstrap Assumptions~\ref{BA:mildsphcoordsast} and~\ref{BA:mildsphcoordsext} and the results of~\cite{Sha14}.\footnote{The strong sup-norm bounds from the Bootstrap Assumptions~\ref{BA:mildsphcoordsast} and~\ref{BA:mildsphcoordsext} imply that the assumptions of~\cite{Sha14} are satisfied and that we can build the orthonormal frame of~\cite{Sha14} upon the coordinates vectorfields.}
\begin{lemma}[Coordinates fractional Sobolev spaces]\label{lem:compH12}
  For all $2$-spheres $S_{u,\uba}\subset \CCba$ and all $2$-spheres $S_{u,\ub} \subset \MM^\ext$, and for all $S$-tangent $k$-tensor $F$, we have
  \begin{align*}
    r^{-1/2}\norm{F}_{\HHt(S_{u,\ub})} & \simeq \sum_{(\varth,\varphi)} \sum_{a_1\cdots a_k\in\{\varth,\varphi\}}\norm{r^{-k}F_{a_1\cdot a_k}}_{H^{1/2}_{\varth,\varphi}},
  \end{align*}
  where the sum is taken over the (two) spherical coordinates systems covering $S_{u,\ub}$ given by the Bootstrap Assumptions~\ref{BA:mildsphcoordsast} and \ref{BA:mildsphcoordsext}, and where $H^{1/2}_{\varth,\varphi}$ denotes the standard fractional Sobolev space on $[\pi/8,7\pi/8]_{\varth}\times[0,2\pi)_\varphi$.
\end{lemma}

We have the following Sobolev estimates.
\begin{lemma}[Sobolev estimates on $2$-spheres]\label{lem:Sobsphere}
  For all $2$-spheres $S_{u,\uba}\subset \CCba$ and all $2$-spheres $S_{u,\ub} \subset \MM^\ext$, we have 
  \begin{align*}
    \norm{F}_{L^4(S_{u,\ub})} & \les \norm{F}_{\HHt(S_{u,\ub})}, \\
    \norm{F}_{L^\infty(S_{u,\ub})} & \les r^{-1/2}\norm{(r\Nd)^{\leq 1}F}_{\HHt(S_{u,\ub})},
  \end{align*}
  for all $S$-tangent tensor $F$. Moreover, we have
  \begin{align*}
   r^{-1/2}\norm{F}_{L^2(S_{u,\ub})} & \les \norm{F}_{\HHt(S_{u,\ub})} \les r^{-1/2}\norm{(r\Nd)^{\leq 1}F}_{L^2(S_{u,\ub})}.
  \end{align*}
\end{lemma}
\begin{proof}
  From the Bootstrap Assumptions~\ref{BA:mildsphcoordsast} and~\ref{BA:mildsphcoordsext}, the $2$-spheres $S_{u,\ub}$ admit weakly regular coordinates systems in the sense of~\cite{Sha14}. The proof of the lemma then follows from the results of~\cite{Sha14}.\footnote{It can also be obtained directly from Lemma~\ref{lem:compH12}, consistently with~\cite{Sha14}.}
\end{proof}

We have the following product estimates for $\HHt$ norms (see~\cite[Corollary 3.4]{Sha14}).
\begin{lemma}[Product estimates]\label{lem:prodH12}
  For all $2$-spheres $S_{u,\uba}\subset \CCba$ and all $2$-spheres $S_{u,\ub} \subset \MM^\ext$, we have
  \begin{align*}
    \norm{FG}_{\HHt(S)} & \les \le(\norm{F}_{L^\infty(S)}+r^{-1}\norm{(r\Nd) F}_{L^2(S)}\ri)\norm{G}_{\HHt(S)}, \\
    \norm{FG}_{L^2(S)} & \les \norm{F}_{\HHt(S)}\norm{G}_{\HHt(S)},
  \end{align*}
  for all $S$-tangent tensors $F,G$.
\end{lemma}

We have the following elliptic estimates on $2$-spheres. 
\begin{lemma}[Elliptic estimates on $2$-spheres]\label{lem:ell}
  Let define the following classical Hodge type operators on $2$-spheres
  \begin{align*}
    \Dd_1U & := \le(\Divd U, \Curld U\ri), & \Dd_1^\ast (f,g) & := -\Nd f + \dual\Nd g,\\
    \Dd_2F & := \Divd F, & \Dd_2^\ast U & := -\half\Nd\otimesh U, 
  \end{align*}
  where $U$ is a $S$-tangent $1$-tensor, $f,g$ are scalar functions and $F$ is a $S$-tangent symmetric traceless $2$-tensor.\\

  For all $2$-spheres $S_{u,\uba}\subset \CCba$ and all $2$-spheres $S_{u,\ub} \subset \MM^\ext$ and for all $S$-tangent $1$-tensor $U$, scalar functions $(f,g)$ and $S$-tangent symmetric traceless $2$-tensor $F$, we have for $0 \leq \ell \leq 2$
  \begin{align}\label{est:classellL2}
    \begin{aligned}
    \norm{r^{-1}(r\Nd)^{\leq \ell +1}U}_{L^2(S)} & \les \norm{(r\Nd)^{\leq \ell}\Dd_1 U}_{L^2(S)}, \\
    \norm{r^{-1}(r\Nd)^{\leq \ell+1}\le(f-\overline{f},g-\overline{g}\ri)}_{L^2(S)} & \les \norm{(r\Nd)^{\ell}\Dd_1^\ast (f,g)}_{L^2(S)}, \\
    \norm{r^{-1}(r\Nd)^{\leq \ell+1}F}_{L^2(S)} & \les \norm{(r\Nd)^{\ell}\Dd_2F}_{L^2(S)}.
    \end{aligned}
  \end{align}
  
  For all $2$-spheres $S_{u,\uba}\subset \CCba$ and all $2$-spheres $S_{u,\ub} \subset \MM^\ext$ and for all $S$-tangent $1$-tensor $U$, scalar functions $(f,g)$ and $S$-tangent symmetric traceless $2$-tensor $F$, we have for $0 \leq \ell \leq 1$
  \begin{align}\label{est:classellH12}
    \begin{aligned}
    \norm{r^{-1}(r\Nd)^{\leq \ell +1}U}_{\HHt(S)} & \les \norm{(r\Nd)^{\leq \ell}\Dd_1 U}_{\HHt(S)}, \\
    \norm{r^{-1}(r\Nd)^{\leq \ell+1}\le(f-\overline{f},g-\overline{g}\ri)}_{\HHt(S)} & \les \norm{(r\Nd)^{\ell}\Dd_1^\ast (f,g)}_{\HHt(S)}, \\
    \norm{r^{-1}(r\Nd)^{\leq \ell+1}F}_{\HHt(S)} & \les \norm{(r\Nd)^{\ell}\Dd_2F}_{\HHt(S)}.
    \end{aligned}
  \end{align}

  Moreover, for all $S$-tangent $1$-tensor $U$ satisfying\footnote{Such system of equations are consequences of the canonical foliation choice. See Section~\ref{sec:connestCCba}.}
  \begin{align*}
    \Dd_1U & = \le(\Divd X, \Curld Y\ri), 
  \end{align*}
  we have
  \begin{align}\label{est:ellHodgeXY}
    \norm{U}_{L^2(S)} & \les \norm{X}_{L^2(S)} + \norm{Y}_{L^2(S)}.
  \end{align}
\end{lemma}
\begin{proof}
  From the Bootstrap Assumptions~\ref{BA:curvast}, \ref{BA:curvext}, \ref{BA:connCCba}, \ref{BA:connext} and Gauss equation~\eqref{eq:Gauss}, we have
  \begin{align}\label{est:GaussKL20}
    \le|K-\frac{1}{r^2}\ri| & \les D\varep \uba^{-2} r^{-1} q^{-1/2}, & \le|K-\frac{1}{r^2}\ri| & \les D\varep \ub^{-3}u^{-1/2},
  \end{align}
  in $\CCba$ and $\MM^\ext$ respectively, as well as
  \begin{align}
    \label{est:GaussKL21}
    \norm{(r\Nd)^{\leq 1}\le(K-\frac{1}{r^2}\ri)}_{L^2(S)} & \les D\varep \uba^{-2} q^{-1/2}, & \norm{(r\Nd)^{\leq 1}\le(K-\frac{1}{r^2}\ri)}_{L^2(S)} & \les D\varep \ub^{-3}u^{-1/2},
  \end{align}
  in $\CCba$ and $\MM^\ext$ respectively. Using estimates~\eqref{est:GaussKL20} and~\eqref{est:GaussKL21}, a rescaling in $r$, and standard energy, Bochner, and higher order elliptic estimates, one obtains the desired $L^2(S)$-estimates~\eqref{est:classellL2} for $0\leq \ell \leq 2$.\\

  To obtain the $\HHt$ elliptic estimates~\eqref{est:classellH12}, we have from~\eqref{est:GaussKL20} and since $K(r^{-2}\gd) = r^{2}K(\gd)$ that $(S,r^{-2}\gd)$ is weakly spherical in the sense of~\cite{Sha14}. Thus, the elliptic estimates~\eqref{est:classellH12} for $\ell=0$ are a direct application from~\cite{Sha14} and rescaling in $r$. We claim that the higher order case $\ell=1$ can also be obtained using additionally the $H^1(S)$-control~\eqref{est:GaussKL21} of $K$.\\

  The proof of the additional estimate~\eqref{est:ellHodgeXY} follows from writing $U = \Dd_1^\ast(f,g)$ with $\overline{f}= \overline{g} = 0$, and standard energy estimates for Poisson equation.
\end{proof}
\begin{remark}
  Poincar\'e inequality is a consequence of the estimate~\eqref{est:classellL2} for the elliptic operator $\Dd_1^\ast$. See also~\cite{Sha14}.
\end{remark}

\chapter{Main results}\label{sec:mainresult}
\section{Initial layers $\varep$-close to Minkowski space}\label{sec:Minkowskilayer}
Let $(\MMt,\g)$ be a smooth vacuum spacetime, and let $(\Sit_1,\CCt_0)$ be spacelike-characteristic initial data for $\MMt$. We say that $\MMt$ contains a \emph{bottom initial layer} $\LLb$ and a \emph{conical initial layer} $\LLc$ which are \emph{$\varep$-close to Minkowski space}, if it admits spacetime subregions $\LLb$ and $\LLc$ as defined in Section~\ref{sec:definitlayer} and if the following smallness assumptions holds on $\LLb, \LLc$.

\paragraph{The bottom initial layer $\LLb$}
With respect to the bottom initial layer coordinates $(x^\mu)$ on $\LLb$ defined in Section~\ref{sec:definitlayer}, we have the following spacetime closeness requirement to the Minkowski metric $\etabold$
\begin{subequations}\label{est:botlayass}
  \begin{align}\label{est:botlayassbis}
    \begin{aligned}
      \sum_{\mu,\nu=0}^3\norm{\g_{\mu\nu} - \etabold_{\mu\nu}}_{L^\infty(\LLb)} & \leq \varep,\\
      \sum_{\mu,\nu=0}^3\norm{\pr^{(\leq 3)}\left(\g_{\mu\nu} - \etabold_{\mu\nu}\right)}_{L^2(\LLb)} & \leq \varep,
    \end{aligned}
  \end{align}
  and for all $1\leq x^0 \leq 5$,
  \begin{align}\label{est:initenersourceSit1bis}
    \int_{\Sit_{x^0}}\le|\D^{\leq 2}\R\ri|^2 & \leq \varep^2,
  \end{align}
\end{subequations}
where in the last estimate the norm is taken with respect to the frame associated to the bottom initial layer coordinates.\footnote{With the notations of Definition~\ref{def:framenorm}, this frame is defined such that $e_0=\Tf^\bott$.}

\begin{remark}\label{rem:centreSit1}
  The centre $\o(1)$ of $\Sit_1$ is chosen with respect to bottom initial layer coordinates, \emph{i.e.} such that
  \begin{align*}
    x^i(\o(1)) & = 0,
  \end{align*}
  for all $1\leq i \leq 3$.\\
  
  The closeness to Minkowski assumptions~\eqref{est:botlayass} imply that
  \begin{align*}
    \le|x^i(\o(t))\ri| & \les \varep, & \le|x^0(\o(t))-t\ri| & \les \varep, 
  \end{align*}
  for all $1\leq i \leq 3$, and
  \begin{align*}
    \le|\doto(t)-\pr_{x^0}\ri| & \les \varep,
  \end{align*}
  for all $t$ such that $\o(t)\in\LLb$.
\end{remark}


\paragraph{The conical initial layer $\LLc$}
With respect to the initial layer optical functions on $\LLc$ and the associated null pair (see Section~\ref{sec:definitlayer}), we have the following closeness requirement to Minkowski spacetime 
\begin{subequations}\label{est:conlayass}
  \begin{itemize}
  \item We have the following bounds on the area radius and the null lapse in $\LLc$
    \begin{align}
      \begin{aligned}
        \le|\rt(\ut,\ubt)-\half(\ubt-\ut)\ri| & \leq \varep \rt(\ut,\ubt),\\
        |\Omt-1| & \leq \varep \ubt^{-1},
      \end{aligned}
    \end{align}
  \item We have the following sup-norm bounds on the null connection coefficients in $\LLc$
    \begin{align}\label{est:Minkowskilayerepsconicnullconn}
      \begin{aligned}
        \le|\trchit-\frac{2}{\rt}\ri| & \leq \varep \ubt^{-2}, & \le|\chiht\ri| & \leq \varep \ubt^{-2}, & \le|\etat\ri| & \leq \varep \ubt^{-2} , &     \le|\omt\ri| & \leq \varep \ubt^{-2} , & \le|\ombt\ri| & \leq \varep \ubt^{-1}, \\
        \le|\trchibt+\frac{2}{\rt}\ri| & \leq \varep \ubt^{-2}, & \le|\chibht\ri| & \leq \varep \ubt^{-1},      & \le|\zet\ri| & \leq \varep \ubt^{-2}, &     \le|\etabt\ri| & \leq \varep \ubt^{-2} , 
      \end{aligned}
    \end{align}
    together with $\xit = \xibt = 0$.
  \item We have the following $L^2(\LLc)$ bounds for (derivatives of) the null connection coefficients
    \begin{align}
      \begin{aligned}
        \norm{\ubt^{-1/2-\ga}\ubt\Ndtt^{\leq 2}\le(\trchit-\frac{2}{\rt}\ri)}_{L^2(\LLc)} & \leq_\ga \varep, & \norm{\ubt^{-1/2-\ga}\ubt\Ndtt^{\leq 2}\chiht}_{L^2(\LLc)} & \leq_\ga \varep, \\
        \norm{\ubt^{-1/2-\ga}\ubt\Ndtt^{\leq 2}\le(\trchibt+\frac{2}{\rt}\ri)}_{L^2(\LLc)} & \leq_\ga \varep, & \norm{\ubt^{-1/2-\ga}\Ndtt^{\leq 2}\chibht}_{L^2(\LLc)} & \leq_\ga \varep,\\
        \norm{\ubt^{-1/2-\ga}\ubt\Ndtt^{\leq 2}\omt}_{L^2(\LLc)} & \leq_\ga \varep, & \norm{\ubt^{-1/2-\ga}\Ndtt^{\leq 2}\ombt}_{L^2(\LLc)} & \leq_\ga \varep,\\
         \norm{\ubt^{-1/2-\ga}\ubt\Ndtt^{\leq 2}\etat}_{L^2(\LLc)} & \leq_\ga \varep, & \norm{\ubt^{-1/2-\ga}\ubt\Ndtt^{\leq 2}\etabt}_{L^2(\LLc)} & \leq_\ga \varep,\\
        \norm{\ubt^{-1/2-\ga}\ubt\Ndtt^{\leq 2}\zet}_{L^2(\LLc)} & \leq_\ga \varep,
      \end{aligned}
    \end{align}
    where $\Ndtt \in \le\{\rt\widetilde{\Nd},\rt\widetilde{\Nd}_4,\widetilde{\Nd}_3\ri\}$ and for all $\ga>0$.
  \item We have the following curvature flux bounds on the cones $\CCt_\ut$, for all $0 \leq \ut \leq 3/2$
    \begin{align}
      \label{est:curvfluxCCt0bis}
      \begin{aligned}
        & \int_{3}^{\infty}\int_{\Stt_{\ut,\ubt}}\bigg(\le|\Ndtt^{\leq 2}\betabt\ri|^2 + \le|\ubt\Ndtt^{\leq 2}(\rhot-\rhoot)\ri|^2 + \le|\ubt\Ndtt^{\leq 2}(\sigmat-\sigmaot)\ri|^2 \\
        & \quad\quad\quad\quad\quad\quad + \le|\ubt^2\Ndtt^{\leq 2}\betat\ri|^2 + \le|\ubt^2\Ndtt^{\leq 2}\alphat\ri|^2\bigg) \d\ubt \leq \varep^2,
      \end{aligned}
    \end{align}
    the following curvature flux bounds for all $3\leq \ubt < +\infty$
    \begin{align}
      \label{est:curvfluxCCbt0bis}
      \begin{aligned}
        & \int_{0}^{3/2}\int_{\Stt_{\ut,\ubt}}\bigg(\le|\Ndtt^{\leq 2}\alphabt\ri|^2 + \le|\ubt\Ndtt^{\leq 2}\betabt\ri|^2 + \le|\ubt^2\Ndtt^{\leq 2}(\rhot-\rhoot)\ri|^2 \\
        & \quad\quad\quad\quad\quad\quad+ \le|\ubt^2\Ndtt^{\leq 2}(\sigmat-\sigmaot)\ri|^2+\le|\ubt^2\Ndtt^{\leq 2}\betat\ri|^2\bigg) \d\ut \leq \varep^2,
      \end{aligned}
    \end{align}
    and the following sup-norm estimates in $\LLc$ for the averages $\rhoot$ and $\sigmaot$
    \begin{align}
      \begin{aligned}
        \le|\Ndtt^{\leq 2}\rhoot\ri| & \leq \varep \ubt^{-3}, & \le|\Ndtt^{\leq 2}\sigmaot\ri| & \leq \varep \ubt^{-3},
      \end{aligned}
    \end{align}
    where $\Ndtt \in \le\{\rt\widetilde{\Nd},\rt\widetilde{\Nd}_4,\widetilde{\Nd}_3\ri\}$.
\end{itemize}
\end{subequations}


\paragraph{The bottom and conical initial layer intersection $\LLb\cap\LLc$}
We assume that on $\LLb\cap\LLc$ we have the following frame comparison
\begin{subequations}
  \begin{align}\label{est:LLbotLLconframe}
    \begin{aligned}
      \le|\g\le(\Tf^\bott,\half(\elbt+\elt)\ri)+1\ri| & \leq \varep,\\
      \le|\g\le(\Nf^\bott,\half(\elt-\elbt)\ri)-1\ri| & \les \varep,
    \end{aligned}
  \end{align}
  and the following coordinates comparisons
  \begin{align}\label{est:LLbotLLconoptical}
    \le|x^0 - \half(\ubt+\ut)\ri| & \leq \varep, & \le|\sqrt{\sum_{i=1}^3\le(x^i\ri)^2} - \half(\ubt-\ut)\ri| & \leq \varep.
  \end{align}
\end{subequations}

\section{Main theorem}
The following theorem is the main result of this paper.
\begin{theorem}[Main theorem, version 2]\label{thm:mainv2}
  Let $(\MMt,\g)$ be a smooth vacuum spacetime, and let $(\Sit_1,\CCt_0)$ be smooth spacelike-characteristic initial data for $\MMt$. Assume that $\MMt$ contains a bottom and conical initial layer $\LLb$ and $\LLc$ adapted to $(\Sit_1,\CCt_0)$ which are $\varep$-close to Minkowski space. There exists $\varep_0>0$ such that if $\varep < \varep_0$, the following holds.
  \begin{subequations}\label{est:infty}
  \begin{itemize}
  \item The spacetime $(\MMt,\g)$ is future geodesically complete.
  \item There exists two spacetime regions $\MM^\intr_\infty,\MM^\ext_\infty\subset \MMt$ such that
    \begin{align*}
      \MMt = \MM^\intr_\infty\cup\MM^\ext_\infty\cup\LLb\cup\LLc,
    \end{align*}
    and such that
    \begin{itemize}
    \item $\MM^\intr_\infty$ is foliated by spacelike maximal hypersurfaces $\Si_t$ which are the level sets of a global time function $t$ on $\MM^\intr_\infty$ ranging from $\too = (1+\cc^{-1})/2$ to $+\infty$, where $0<\cc<1$ is a fixed parameter, sufficiently close to $1$. 
    \item $\MM^\ext_\infty$ is foliated by outgoing null hypersurfaces $\CC_u$ which are the level sets of a global optical function $u$ on $\MM^\ext_\infty$ ranging from $1$ to $+\infty$.
    \item There exists a global affine parameter $\ub$ on $\MM^\ext_\infty$ foliating the null hypersurfaces $\CC_u$ ranging from $\cc^{-1}$ to $+\infty$.
    \item The transition hypersurface $\TT$ satisfies
      \begin{align*}
        \TT & := \MM^\intr_\infty\cap\MM^\ext_\infty = \le\{u=\cc\ub\ri\},
      \end{align*}
      and on $\TT$, we have
      \begin{align*}
        t = \half (u+\ub).
      \end{align*}
    \end{itemize}
  \item In $\MM^\intr_\infty$, the following curvature decay holds\footnote{\label{foo:normsSitmainthm} The norms are taken with respect to the frame associated to the spacelike hypersurfaces $\Si_t$.}
    \begin{align}
      |\R| & \les \varep t^{-7/2},
    \end{align}
    we have the following control of the time function $t$
    \begin{align}
      \le|\g(\D t,\D t)+1\ri| & \les \varep t^{-3/2}, & \le|\D^2t\ri| & \les \varep t^{-5/2},
    \end{align}
    and the maximal hypersurfaces $\Si_t$ approach the Euclidean disks in the following (intrinsic) sense\footnote{Using global harmonic coordinates (see Theorem~\ref{thm:globharmonics}), this can be alternatively formulated using coordinates.}
    \begin{align}
      \le|\RRRic\ri| & \les \varep t^{-7/2}, & \le|\trth-\frac{2}{t}\le(\frac{1+\cc}{1-\cc}\ri)\ri| & \les \varep t^{-5/2}, & |\thh| & \les \varep t^{-5/2},
    \end{align}
    where $\RRRic$ is the Ricci curvature tensor of $\Si_t$ and $\th$ is the second fundamental form of the boundaries $\pr\Si_t\subset \TT$.
  \item In $\MM^\ext_\infty$, the following curvature decay holds
    \begin{align}\label{est:inftynullcurv}
      \begin{aligned}
        \le|\al\ri| & \les \varep \ub^{-7/2}, & \le|\be\ri| & \les \varep \ub^{-7/2}, & \le|\rho-\rhoo\ri| & \les \varep \ub^{-3}u^{-1/2},\\
        \le|\alb\ri| & \les \varep \ub^{-1}u^{-5/2}, & \le|\beb\ri| & \les \varep \ub^{-2}u^{-5/2}, & \le|\sigma-\sigmao\ri| & \les \varep \ub^{-3}u^{-1/2},
      \end{aligned}
    \end{align}
    as well as
    \begin{align}\label{est:inftynullcurvaverage}
      |\rhoo| & \les \varep^2 \ub^{-3}u^{-2}, & |\sigmao| & \les \varep^2 \ub^{-3}u^{-2},
    \end{align}
    and we have the following control of the optical function and affine parameter $u$ and $\ub$
    \begin{align}
      \g(\D u,\D u) & = 0, & \g(\D\ub, \D u) & = -2, & \le|\g(\D\ub,\D\ub)\ri| & \les \varep u^{-3/2}, 
    \end{align}
    and
    \begin{align}\label{est:inftynullconn}
      \begin{aligned}
      \le|\trchi-\frac{2}{r}\ri| & \les \varep \ub^{-2}u^{-1/2}, & \le|\chih\ri| & \les \varep \ub^{-2}u^{-1/2}, & \le|\ze\ri| & \les \varep \ub^{-2}u^{-1/2}, & |\omb| & \les \varep \ub^{-1}u^{-3/2}, \\
      \le|\trchib+\frac{2}{r}\ri| & \les \varep \ub^{-2}u^{-1/2}, & \le|\chibh\ri| & \les \varep \ub^{-1}u^{-3/2}, & \le|\xib\ri| & \les \varep \ub^{-1}u^{-3/2},
      \end{aligned}
    \end{align}
    together with $\om = |\xi| = 0$ and $\ze = \eta = -\etab$.\\
    
    Moreover, the $2$-spheres $S_{u,\ub}$ level sets of $u,\ub$ approach the Euclidean $2$-spheres in the following (intrinsic) sense
    \begin{align}
      \le|K-\frac{1}{r^2}\ri| & \les \varep \ub^{-3}u^{-1/2}, & \le|r - \half (\ub-u)\ri| & \les \varep \ub^{-1}u^{-1},   
    \end{align}
    where $r$ denotes the area radius of $S_{u,\ub}$.
  \item The spacetime $\MMt$ admits a future timelike and future null infinity $i^+$ and $\scri^+$, and the future null infinity $\scri^+$ is \emph{future geodesically complete}.\footnote{See~\cite{Wal84} and also~\cite{Chr.Kla93},~\cite{Kla.Nic03} for definitions and discussions.}
    Moreover, $\MMt$ admits well-defined Bondi mass and angular momentum at null infinity, which satisfy respectively a Bondi mass loss formula and an angular momentum evolution equation on $\scri^+$, and which tend to $0$ at timelike infinity $i^+$.\footnote{See Section~\ref{sec:conclusion} for precise definitions.}
  \end{itemize}
  \end{subequations}
\end{theorem}

\paragraph{Remarks on Theorem~\ref{thm:mainv2}}
\begin{enumerate}[ref=\thetheorem\alph*,label=\thetheorem\alph*]
\item Alternatively, the spacetime region $\MM_\infty = \MM_\infty^\intr\cup\MM_\infty^\ext$ can also be foliated by the $2$-spheres $S'_{u',\ub'}$ of the geodesic foliation on the incoming null cones backward emanating from the central axis $\o$ (see the definitions of the optical function $\ub'$ and the geodesic parameter $u'$ of the so-called last cones geodesic foliation in Section~\ref{sec:deflastconesfoliation}). Analogous decay estimates to~(\ref{est:infty}) can be obtained in interior and exterior regions with respect to the null frame adapted to $u'$ and $\ub'$.\\ Since the classical definition of future null infinity and associated asymptotic quantities involve taking limits along the outgoing null cones, and also since the proof of Theorem~\ref{thm:mainv2} foremost relies on the outgoing null cones level sets of $u$, we prefered to state Theorem~\ref{thm:mainv2} using the time, optical and geodesic affine parameter functions $t,u$ and $\ub$.
\item More specific $L^p$-based decay statements, or boundedness statements for $L^2$-fluxes can be obtained for derivatives of the curvature and connection coefficients.
\item\label{item:Spacetimemass} The decay rates of~\eqref{est:infty} match the decay rates obtained in~\cite{Chr.Kla93,Kla.Nic03}. A notable exception to that statement is the strong decay rate~(\ref{est:inftynullcurvaverage}) for the mean value $\rhoo$. This is due to the two different spacetime regions studied in this paper and in~\cite{Chr.Kla93,Kla.Nic03}. In~\cite{Chr.Kla93} and~\cite{Kla.Nic03}, the mean value $\rhoo$ satisfies a weak decay rate of the type $|\rhoo| \les \varep \ub^{-3}$. Namely, $\rhoo$ is not controlled by energy estimates and is only determined by integrating Bianchi equation along $\elb$ from the initial spacelike hypersurface into the spacetime. Since on the initial spacelike hypersurface $\rhoo$ is related to the ADM mass $M$ \emph{via} $\rhoo \sim - 2M/r^{3}$, the decay rate for $\rhoo$ obtained by integration is at most $\rhoo \sim \ub^{-3}$. In the present paper, we determine $\rhoo$ by integrating Bianchi equation \emph{from the central axis} $\o$ backwards in $\MM$. Since $r^3\rhoo \to 0$ when $r\to 0$, the initial value of $\rhoo$ on $\o$ is virtually $0$, thus its decay rate in $\MM$ is only dictated by the nonlinear terms in Bianchi equation~(\ref{eq:Nd3rhoo}), from which we deduce the strong decay rate $\rhoo \sim \ub^{-3}u^{-2}$. This corresponds to obtaining bounds for the Bondi mass in $\scri^+$ integrating Bondi mass loss formula from timelike infinity $i^+$ backwards on $\scri^+$, provided that it is known that the final Bondi mass vanishes at $i^+$.\footnote{All the statements for $\rhoo$ can also be formulated in terms of the Hawking mass/average of mass aspect function.}
\end{enumerate}

\section{Auxiliary theorems}
In this section, we state auxiliary results to the main theorem (Theorems~\ref{thm:globharmonics},~\ref{thm:vertex} and~\ref{thm:canonical}), which are independent and of more general interest. Their respective proofs are given in Appendix~\ref{sec:globharmo},~\ref{app:vertexlimits} and~\ref{app:canlocalex}. We also state an existence of initial layers theorem (see Theorem~\ref{thm:initlayer}), which we claim can be obtained from previous literature results.

\subsection{Global harmonic coordinates}
We first have the following definition.
\begin{definition}[Weakly regular $3$-disk]\label{def:weakreg3DStab}
  Let $\Si$ be a smooth Riemannian manifold diffeomorphic to the unit disk $\DDD$ of $\RRR^3$. Let $C>0$. We say that $\Si$ is a \emph{weakly regular} $3$-disk with constant $C>0$, if the following functional estimates hold.\\

  For all $\Si$-tangent tensors $F$, we have the following Sobolev estimates
  \begin{align}\label{est:sobeucl}
    \begin{aligned}
      \norm{F}_{L^6(\Si)} & \leq C\norm{\nab^{\leq 1}F}_{L^2(\Si)}, \\
      \norm{F}_{L^\infty(\Si)} & \leq C\norm{\nab^{\leq 2}F}_{L^2(\Si)}.
    \end{aligned}
  \end{align}
  
  For all $\Si$-tangent tensor $F$ we have the following $L^{3/2}(\Si)-L^1(\pr\Si)$-Poincar\'e inequality\footnote{This Poincar\'e estimate holds true in the Euclidean disk case, \emph{via} the following sequence of more standard Poincar\'e and Sobolev embeddings
  \begin{align*}
    \norm{F}_{L^1(\Si)} & \les \norm{\nab F}_{L^1(\Si)} + \norm{F}_{L^1(\pr\Si)},\\
    \norm{F}_{L^{3/2}(\Si)} & \les \norm{\nab^{\leq 1} F}_{L^1(\Si)},\\
    \norm{F}_{L^2(\Si)} & \les \norm{\nab^{\leq 1}F}_{L^{3/2}(\Si)}.
  \end{align*}}
  \begin{align}
    \label{est:PoincareL3/2}
    \norm{F}_{L^2(\Si)} & \leq C\le(\norm{\nab F}_{L^{3/2}(\Si)} + \norm{F}_{L^1(\pr\Si)}\ri).
  \end{align}

  For all $\Si$-tangent tensor $F$, we have the following trace estimate from $\Si$ on $\pr\Si$
  \begin{align}\label{est:traceeucl}
    \norm{F}_{H^{1/2}(\pr\Si)} & \leq C\norm{\nab^{\leq 1}F}_{L^2(\Si)}.
  \end{align}
\end{definition}


\begin{theorem}[Global harmonic coordinates on $\Si$]\label{thm:globharmonics}
  Let $\Si$ be a Riemannian manifold diffeomorphic to the unit coordinate disk $\DDD$ of $\RRR^3$. Assume that $\Si$ is a weakly regular $3$-disk with constant $C>0$. 
  Let $\varep>0$ and assume that on $\Si$ the following $L^2$ bounds for the Ricci curvature tensor of $\Si$ and the second fundamental form $\th$ of the boundary $\pr\Si$ hold\footnote{The second fundamental form $\th$ of $\pr\Si\subset\Si$ is defined by
    \begin{align*}
      \th(X,Y) & := g(\nab_X N,Y), 
    \end{align*}
    for $X,Y\in T\pr\Si$, where $N$ is the outward-pointing unit normal to $\pr\Si$ in $\Si$.
}
  \begin{align}\label{est:L2RicTh}
    \begin{aligned}
      \norm{\RRRic}_{L^2(\Si)} & \leq \varep, \\
      \norm{\trth-2}_{H^{1/2}(\pr\Si)} +\Vert\thh\Vert_{H^{1/2}(\pr\Si)} & \leq \varep. 
    \end{aligned}
  \end{align}
  Then, there exists $\varep_0(C)>0$ such that if $\varep < \varep_0$, there exists global harmonic coordinates $(x^i)$ from $\Si$ onto $\DDD$ with the following bounds
  \begin{align}\label{est:nabxi}
    \sum_{i,j=1}^3\norm{\g\le(\nab x^i,\nab x^j\ri) -\de_{ij}}_{L^\infty(\Si)} + \sum_{i=1}^3\norm{\nab^{\leq 1}\nab^2 x^i}_{L^2(\Si)} + \sum_{i=1}^3 \norm{N(x^i)-x^i}_{L^\infty(\pr\Si)} & \les \varep,
  \end{align}
  where $N$ denotes the outward-pointing unit normal to $\pr\Si$. Moreover, for all $k\geq 0$, we have the following higher regularity estimates
  \begin{align}\label{est:nabknabxi}
    \sum_{i=1}^3\norm{\nab^{k}\nab^3x^i}_{L^2(\Si)} & \leq C_k\le(\norm{\nab^{\leq k}\RRRic}_{L^2(\Si)} + \norm{\Nd^{\leq k}(\th-\gd)}_{H^{1/2}(\pr\Si)} + \varep\ri).
  \end{align}
\end{theorem}

\paragraph{Remarks on Theorem~\ref{thm:globharmonics}}
\begin{enumerate}[ref=\thetheorem\alph*,label=\thetheorem\alph*]
\item For the metric components $g_{ij}$ in the harmonic coordinate $(x^i)$, we deduce from~\eqref{est:nabxi} and~\eqref{est:nabknabxi} the following respective $L^2$ estimates\footnote{$L^p$ and sup-norm estimates can be deduced by Sobolev embeddings.}
  \begin{align*}
    \norm{\pr^{\leq 2}\le(g_{ij}-\de_{ij}\ri)}_{L^2(\Si)} & \les \varep,\\
    \norm{\pr^{\leq k+2}\le(g_{ij}-\de_{ij}\ri)}_{L^2(\Si)} & \les C_k\le(\norm{\nab^{\leq k}\RRRic}_{L^2(\Si)} + \norm{\Nd^{\leq k}(\th-\gd)}_{H^{1/2}(\pr\Si)} + \varep \ri).
  \end{align*}
\item Theorem~\ref{thm:globharmonics} can be cast as an existence and control of solutions to the \emph{Dirichlet problem for harmonic maps} result with the Euclidean unit disk $\DDD$ of $\RRR^3$ as target manifold.
\item Theorem~\ref{thm:globharmonics} improves on the results of~\cite[Section 7]{Czi.Gra19a} since it only uses elementary (energy and Bochner) estimates and provides optimal quantitative bounds for the metric components.  
\item In the context of the present paper, we use estimate~\eqref{est:nabxi} and estimate~\eqref{est:nabknabxi} with $k=2$.
\item The functional hypothesis from the weak regularity Definition~\ref{def:weakreg3DStab} are in particular satisfied if there exists weakly regular coordinates on $\Si$. These coordinates do not need to be harmonic. In the present paper, the induced coordinates obtained when applying the existence of maximal hypersurfaces result of~\cite{Cho76} are such coordinates.
\item Using the harmonic coordinates of Theorem~\ref{thm:globharmonics} and standard analysis results, one can directly improve on the constants in the assumed functional estimates of Definition~\ref{def:weakreg3DStab} for $\varep>0$ sufficiently small. This suggests that the (already weak) functional assumptions of Definition~\ref{def:weakreg3DStab} could be removed.
\item The crux of the proof of Theorem~\ref{thm:globharmonics} is the obtention of a Bochner identity for the Dirichlet problem on $\Si$ with coercive boundary terms directly controlling Neumann data. See equation~\eqref{eq:refinedBochner}.
\end{enumerate}

\subsection{Axis limits}
The following theorem provides limits for (all derivatives of) the metric, connection and curvature at the central axis $\o$ for the incoming null cones emanating from $\o$ foliated by geodesic affine parameter. Its proof is given in Appendix~\ref{app:vertexlimits}. 
\begin{theorem}[Axis limits]\label{thm:vertex}
  Let $(\MM,\g)$ be a smooth Lorentzian manifold. Let $\o$ be a timelike geodesic. There exists coordinates $x^\mu$ in a neighbourhood of $\o$, smooth in $\MM \setminus \o$, such that the level sets of
  \begin{align*}
    \ub' := x^0+\sqrt{\sum_{i=1}^3\le(x^i\ri)^2}
  \end{align*}
  are the incoming null cones $\CCb_{\ub'}$ emanating from the axis $\o$, and such that
  \begin{align*}
    u' := x^0 - \sqrt{\sum_{i=1}^3\le(x^i\ri)^2}
  \end{align*}
  is the null geodesic affine parameter on the cones $\CCb_{\ub'}$.\footnote{See also the definition of $u'$ in Section~\ref{sec:deflastconesfoliation}.} We call the coordinates $x^\mu$ \emph{(Cartesian) optical  normal coordinates}. \\

  \begin{subequations}\label{est:thmlimits}
  Moreover, there exists (transported along $\CCb_{\ub'}$) spherical coordinates $\varth',\varphi'$ on the $2$-spheres $S_{u',\ub'}$ such that for the induced metric $\gd'$, we have
  \begin{align}
    \le(\pr^k_{u'},(\pr_{u'}+\pr_{\ub'})^l,\pr^m_{\om'}\ri) \le(\gd - \le(\frac{\ub'-u'}{2}\ri)^2\le(\d(\varth')^2 + \sin^2\varth'\d(\varphi')^2\ri)\ri) = O\le(|x|^{|4-k|}\ri), 
  \end{align}
  when $|x| \to 0$, for $k,l,m\geq 0$, and where $\le(\pr^k_{u'},(\pr_{u'}+\pr_{\ub'})^l,\pr^m_{\om'}\ri)$ denotes all combinations of partial derivatives containing respectively $k,l,m$-derivatives of $u',u'+\ub',\om'=\varth',\varphi'$. In particular, for the area radius $r'$ of $S_{u',\ub'}$, we have
  \begin{align}
    \le(\pr^k_{u'},(\pr_{u'}+\pr_{\ub'})^l\ri)\le(r'(u',\ub') -\frac{\ub'-u'}{2}\ri) & = O(|x|^{|3-k|}).
  \end{align}
  Furthermore, the following limits hold when $|x| \to 0$ and for $k,l,m\geq 0$
  \begin{itemize}
  \item for the optical defect $\yy' := -\half\g(\D u',\D u')$
    \begin{align}\label{est:thmlimitsyy}
      \le|\le((r'\Nd'_3)^k,(\Nd'_3+\Nd'_4)^l,(r'\Nd')^m\ri) \yy'\ri| & = O\le(|x|^{2}\ri),
    \end{align}
  \item for the null connection coefficients associated to the null pair $(\elb',\el')$ (defined such that $\elb'=-\D\ub'$)
    \begin{align}\label{est:thmlimitsGa}
      \le|\le((r'\Nd'_3)^k,(\Nd'_3+\Nd'_4)^l,(r'\Nd')^m\ri) \le(\chi'- \frac{1}{r'}\gd', \chib' + \frac{1}{r'}\gd', \eta', \ze', \etab', \om', \omb', \xi', \xib'\ri)\ri| & = O\le(|x|\ri),
    \end{align}

  \item for the null curvature components
    \begin{align}\label{est:thmlimitscurv}
      \le|\le((r'\Nd'_3)^k,(\Nd'_3+\Nd'_4)^l,(r'\Nd')^m\ri) \le(\al',\be',\rho',\sigma',\beb',\alb'\ri)\ri| & = O\le(1\ri),
    \end{align}
  \end{itemize}
\end{subequations}
\end{theorem}
\paragraph{Remarks on Theorem~\ref{thm:vertex}}
\begin{enumerate}[ref=\thetheorem\alph*,label=\thetheorem\alph*]
\item The vertex limits~\eqref{est:thmlimits} are consequences of the fact that in a classical normal coordinates system, the metric and its first derivatives are trivial at the point $O$. They are obtained by a change of coordinates from classical to Cartesian and subsequently spherical optical normal coordinates and by expressing the null coefficients in terms of the spherical optical normal coordinates. See Sections~\ref{sec:vertexlimitsg} and~\ref{sec:limitsnullconn}. 
\item The vanishing of all $\Nd'_3+\Nd'_4$ derivatives is a consequence of the translation invariance along the axis.  
\item The vertex limits~\eqref{est:thmlimitsyy}--\eqref{est:thmlimitscurv} are sharp in terms of the asymptotic behaviour of the $\Nd'_3+\Nd'_4$ and $r'\Nd'$ derivatives. For the $\Nd'_3$ derivatives, these limits allow for a blow-up of the $k$ derivatives of the type $(\Nd'_3)^kF \sim (r')^{-k}$, which we believe is far from being optimal in most cases.
\item Better limits can hold for a vacuum spacetime. For example, using that $\tr'\alb' = 0$ in that case and integrating the Raychaudhuri equation~\eqref{eq:Nd3trchib}, one can obtain better bounds for $\trchib'$.
\item The coordinates $u',\ub'$ from Theorem~\ref{thm:vertex} are the coordinate of the so-called last cones geodesic foliation defined in Section~\ref{sec:deflastconesfoliation}.
\end{enumerate}

\subsection{Well-posedness of the canonical foliation}
The following theorem ensures that the canonical foliation on $\CCba$ defined in Section~\ref{sec:defcannull} is locally well defined and provides vertex limits for its associated metric and null connection coefficients.
\begin{theorem}[The canonical foliation on $\CCba$]\label{thm:canonical}
  Let $(\MM,\g)$ be a smooth Lorentzian manifold. Let $\CCba\subset\MM$ be a smooth null cone emanating from a point $\o(\uba)$, where we also assume that a unit timelike vector $\doto(\uba)$ is given. There exists a function $u$ in a neighbourhood of $\o(\uba)$ in $\CCba$, smooth on $\CCba\setminus\o(\uba)$, such that its level sets $S_u\subset\CCba$ define a \emph{canonical foliation} of $\CCba$ in the following sense (see also the definitions of Section~\ref{sec:defcannull}):
  \begin{itemize}
  \item For the null connection coefficients associated to the null pair $(\elb,\el)$ orthogonal to $S_u$ (defined such that $\elb(u)=2$), the following elliptic condition is satisfied
    \begin{align}\label{eq:ellcanthm}
      \begin{aligned}
      \Divd\eta + \rho & = \rhoo,\\
      \ombo & = 0,
      \end{aligned}
    \end{align}
    on each $2$-sphere $S_u$.
  \item The following limits hold at the vertex $\o(\uba)$
    \begin{align}\label{eq:initvaluecanthm}
      u|_{\o(\uba)} = \uba, && \text{and} && \g(\elb,\doto(\uba))|_{\o(\uba)} = -1. 
    \end{align}
  \end{itemize}
  Moreover, there exists (transported along $\CCba$) spherical coordinates $\varth,\varphi$ such that the induced metric $\gd$ on $S_u$ in these coordinates satisfy the same limits at the vertex $\o(\uba)$ as in Theorem~\ref{thm:vertex}. The limits for the area radius, null connection coefficients associated to the canonical foliation are also identical to the ones of Theorem~\ref{thm:vertex}.
\end{theorem}

\paragraph{Remarks on Theorem~\ref{thm:canonical}}
\begin{enumerate}[ref=\thetheorem\alph*,label=\thetheorem\alph*]
\item Expressed using geodesic affine parameter, the system~\eqref{eq:ellcanthm} rewrites as a coupled system of elliptic and transport equation, with initial value given by~\eqref{eq:initvaluecanthm}. The proof of Theorem~\ref{thm:canonical} then relies on a standard Banach-Picard iteration similar to the one performed in~\cite[Section 6]{Czi.Gra19}. See Appendix~\ref{app:canlocalex}.
\item\label{item:stabilitycanonical} As a byproduct of the Banach-Picard iteration -- using an implicit function theorem --, one obtains that the solutions to the above mentioned system of elliptic-transport equations are also unique and stable under a perturbation of the background spacetime metric. This justifies the continuity argument of Section~\ref{sec:propregeextension}.
\item The vertex limits for the metric and null connection coefficients are used in Section~\ref{sec:connestCCba} as initial data to obtain global estimates on $\CCba$.
\end{enumerate}


\subsection{Existence and control of initial layers}
We have the following theorem, which enables to relate the result of the main Theorem~\ref{thm:mainv2} proved in this paper, to the most general Theorem~\ref{thm:mainv1} stated in the introduction.
\begin{theorem}[Existence and control of initial layers~\cite{Czi.Gra19a},~\cite{Li.Zhu18}]\label{thm:initlayer}
  Let $\varep >0$. There exists $\varep'>0$ such that if $(\MMt,\g)$ is a smooth vacuum spacetime admitting spacelike-characteristic initial data $(\Sit_1,\CCt_0)$ which are $\varep'$-close to Minkowski space, then $(\MMt,\g)$ contains a \emph{bottom initial layer} $\LLb$ and a \emph{conical initial layer} $\LLc$, which are $\varep$-close to Minkowski space, as defined in Section~\ref{sec:Minkowskilayer}.
\end{theorem}
\paragraph{Remarks on Theorem~\ref{thm:initlayer}}
\begin{enumerate}[ref=\thetheorem\alph*,label=\thetheorem\alph*]
\item Combining the existence and control of initial layer Theorem~\ref{thm:initlayer} and the main Theorem~\ref{thm:mainv2}, we obtain the general Theorem~\ref{thm:mainv1}.
\item The existence and control of the bottom initial layer $\LLb$ follows from a small data local existence result for the spacelike-characteristic Cauchy problem. Such a result has been obtained in~\cite{Czi.Gra19,Czi.Gra19a} under regularity assumptions weaker than in the present paper, and applies in particular in the present case.\footnote{It can also easily be reproved at the present regularity using the elementary techniques of this paper.}
\item The conical initial layer decay rates for the metric, connection and curvature from estimates~(\ref{est:conlayass}) correspond to the decay rates of~\cite{Chr.Kla93,Kla.Nic03}. The required regularity is also provided by~\cite{Chr.Kla93,Kla.Nic03}. In particular, a last layer along an outgoing null cone $\CCt_0$ in~\cite{Kla.Nic03} would be an admissible conical initial layer in the present paper.
\item The existence of the conical initial layer $\LLc$ can be obtained from the existence result of~\cite{Li.Zhu18} for initial data posed on characteristic hypersurfaces.
\item The comparisons in the intersection of the initial layers $\LLb\cap\LLc$ can also be obtained from~\cite{Li.Zhu18}. It could also be obtained in the case of a last outgoing conical initial layer from~\cite{Kla.Nic03}.
\end{enumerate}

\section{Proof of the main theorem}\label{sec:proofmainthm}
The main part of the proof of Theorem~\ref{thm:mainv2} is a bootstrap argument. We define $\underline{U}^\ast$ to be the supremum of all $\uba$ such that $\MMt$ admits subregions $^{(\cc)}\MM_\uba$ -- for all $\cco \leq \cc \leq (1+\cco)/2$ -- satisfying a set of mild and strong bootstrap assumptions which are collected in Sections~\ref{sec:mildBA} and~\ref{sec:strongBA}, and which we further call \emph{the Bootstrap Assumptions}.\footnote{\label{foo:uniftrans}We assume that the Bootstrap Assumptions are \emph{uniformly satisfied} for all transition parameter $\cco \leq \cc \leq (1+\cco)/2$.} We assume by contradiction that $\underline{U}^\ast < \infty$. From the closedness of the Bootstrap Assumptions and by propagation of regularity, $\MMt$ admits smooth subregions $^{(\cc)}\MMt_{\underline{U}^\ast}$ satisfying the Bootstrap Assumptions.\footnoteref{foo:uniftrans}

\begin{remark}
  The fact that $\MMt$ admits first subregions, say $^{(\cc)}\MM_{\uba = 4}$, satisfying the Bootstrap Assumptions follows from local constructions in the (domain of dependence of the) initial layer $\LLb$. See also the arguments of Step 2 in Section~\ref{sec:propregeextension}.
\end{remark}

\begin{remark}
  We do not state nor prove a precise propagation of regularity/smoothness result in this paper, which is used to guarantee the existence of the above smooth spacetime subregion $\MM_{\underline{U}^\ast}$ and to run the extraction argument of Section~\ref{sec:conclusion}. We claim that such a result can be obtained repeating the global energy estimates and subsequent arguments of Sections~\ref{sec:globener}--~\ref{sec:initlayer}. The smooth character of the spacetime is crucially used so that the limits of Theorem~\ref{thm:vertex} for (derivatives of) the null connection coefficients at the vertex $\o(\uba)$ of the last cone $\CCba$ hold.\footnote{The limits of Theorem~\ref{thm:vertex} which are actually used in this paper require to obtain limits for up to three derivatives of the connection coefficients. This requires to control the metric in at least a $C^4$-sense at the vertex. Given the optimal possible regularity $\pr^{\leq 4}g_{\mu\nu} \in L^\infty_t L^2$ for the metric component, consistently with the curvature control $\D^{\leq 2}\R \in L^\infty_tL^2$, such a control can only be obtained under an higher regularity assumption. This is the reason for the smoothness assumption in this paper.} It is also used for convenience in the local existence results applied in the extension procedure of Section~\ref{sec:propregeextension}.   
\end{remark}

In Sections~\ref{sec:globener}--\ref{sec:initlayer}, we show that under the Bootstrap Assumptions in the subregions $^{(\cc)}\MM_\uba$ and the $\varep$-closeness to Minkowski of the initial layers $\LLb$ and $\LLc$, all the Bootstrap Assumptions from Sections~\ref{sec:mildBA} and~\ref{sec:strongBA} can be improved. We give an overview of the results obtained in these Sections in Section~\ref{sec:BAproof}.\\

In Section~\ref{sec:propregeextension}, we show that provided that the Bootstrap Assumptions with improved constants hold in $^{(\cc)}\MM_\uba$, the spacetime region $^{(\cc)}\MM_\uba\cup\LLc$ can be extended to a smooth spacetime $^{(\cc)}\MM_{\uba+\de}\cup\LLc$ for $\de>0$ such that the Bootstrap Assumptions hold. This contradicts the finitness of $\underline{U}^\ast$.\\

In Section~\ref{sec:conclusion}, we deduce from $\underline{U}^\ast = + \infty$ the main conclusions of Theorem~\ref{thm:mainv2}. This finishes the proof of Theorem~\ref{thm:mainv2}.

\subsection{Improvement of the Bootstrap Assumptions}\label{sec:BAproof}
In this section, we give an overview of the improvement of the Bootstrap Assumptions which is performed in Sections~\ref{sec:globener}--\ref{sec:initlayer}.
\begin{itemize}
\item In Section~\ref{sec:initlayerenergy}, we prove that under the Bootstrap Assumptions and the $\varep$-closeness assumption of the initial layers to Minkowski space (see Section~\ref{sec:Minkowskilayer}), improved bounds for the energy fluxes for the curvature through the hypersurfaces $\Si_\too$ and $\CC_1$ hold.
\item In Section~\ref{sec:globener}, we show by performing global energy estimates in $\MM^\ext\cup\MM^\intr_\bott$, using the improved initial energy fluxes, that improved energy bounds for the curvature hold on all hypersurfaces $\Si_t$, $\CC_u$, $\Si_t^\ext$ and $\CCba\cap\MM^\ext$. The global energy estimates of Section~\ref{sec:globener} are performed for \emph{one} transition parameter $\cc$ chosen by a mean value argument.
\item In Section~\ref{sec:curvest}, we show that using the improved energy bounds on $\CCba\cap\MM^\ext$ and $\CC_u$, $\Si_t^\ext$, one obtains improved control of the spacetime curvature tensor in $\MM^\ext$.
\item In Section~\ref{sec:planehypcurvest}, we show that using the improved energy bounds on $\Si_t$ and the improved trace bounds for the curvature on the interface $\TT=\MM^\intr\cap\MM^\ext$ from Section~\ref{sec:curvest}, one obtains improved curvature bounds in $\MM^\intr_\bott$. 
\item In Section~\ref{sec:tipcurvest}, we show that using the curvature control of Section~\ref{sec:planehypcurvest} on $\Si_\tast$, applying a rescaling, extension and local existence results and performing a local energy estimate in $\MM^\intr_\topp$, one obtains improved curvature bounds on $\CCba\cap\MM^\intr$. In Section~\ref{sec:alltransparam}, we show by rescaling and local energy estimates that curvature estimates are improved for all the constructions related to \emph{all} transition parameters $\cc$.
\item In Section~\ref{sec:connestCCba}, we show that using the vertex limits of Theorems~\ref{thm:vertex} and~\ref{thm:canonical}, the Bootstrap Assumptions and the improved curvature bounds on $\CCba$, we obtain improved bounds for the null connection coefficients of the canonical foliation. We moreover show that under these improved bounds, the conformal factor and exterior rotation vectorfields satisfy improved bounds on $S^\ast$. Using the improved bounds for the exterior rotation on $S^\ast$, we improve their bounds on $\CCba\cap\MM^\ext$ in Section~\ref{sec:controlOOECCBA}.
\item In Section~\ref{sec:connest}, we show that using the improved bounds for the connection and rotation coefficients on $\CCba\cap\MM^\ext$, we obtain improved bounds for the connection and rotation coefficients on $\MM^\ext$.
\item In Section~\ref{sec:planehypconnest}, we show using the improved curvature estimates in $\MM^\intr_\bott$ and the improved trace estimates for the null connection coefficients on $\TT$ that the connection coefficients related to the maximal foliation are improved, and thus the control of the time translation approximate Killing field $\TI$ is improved. From an application of Theorem~\ref{thm:globharmonics}, we obtain bounds for the harmonic coordinates of the last slice $\Si_\tast$. We further show using these bounds that the remaining interior Killing fields are controlled in $\MM^\intr_\bott$ by integration from $\Si_\tast$. Last, we show that at the interface $\TT$, the differences of interior and exterior Killing fields is controlled.
\item In Section~\ref{sec:initlayer}, we show using the previous improvement for all connection and curvature coefficients that the comparison of the constructions of $\MM$ to the constructions of the initial layers $\LLb$ and $\LLc$ are improved. This finishes the improvement of the Bootstrap Assumptions.
\end{itemize}

\subsection{Extension of $\MM_{\protect\uba}$}\label{sec:propregeextension}
In this section, we assume that $\MMt$ admits a smooth spacetime subregion $^{(\cc)}\MM_\uba$ satisfying the Bootstrap Assumptions of Sections~\ref{sec:mildBA} and~\ref{sec:strongBA} with improved constants. We show that $\MMt$ admits a smooth spacetime region $^{(\cc)}\MMt_{\uba+\de}$ for $\de>0$ satisfying the Bootstrap Assumptions.

\begin{remark}
  Here and in the rest of this section, the transition parameter $\cc$ is fixed. We thus omit the label $\cc$ in the arguments below. 
\end{remark}

\paragraph{Step 1: Spacetime extension.} We first show that there exists a smooth spacetime $\NN$ containing $\MM_\uba\cup\LLc$ such that $\MM_\uba$ is strictly included in $\NN$.\footnote{Here and in the following, strict inclusion means that $\MM_\uba$ is composed of interior points of the manifold $\NN$.}\\

\begin{figure}[h!]
  \centering
  \includegraphics[height=4cm]{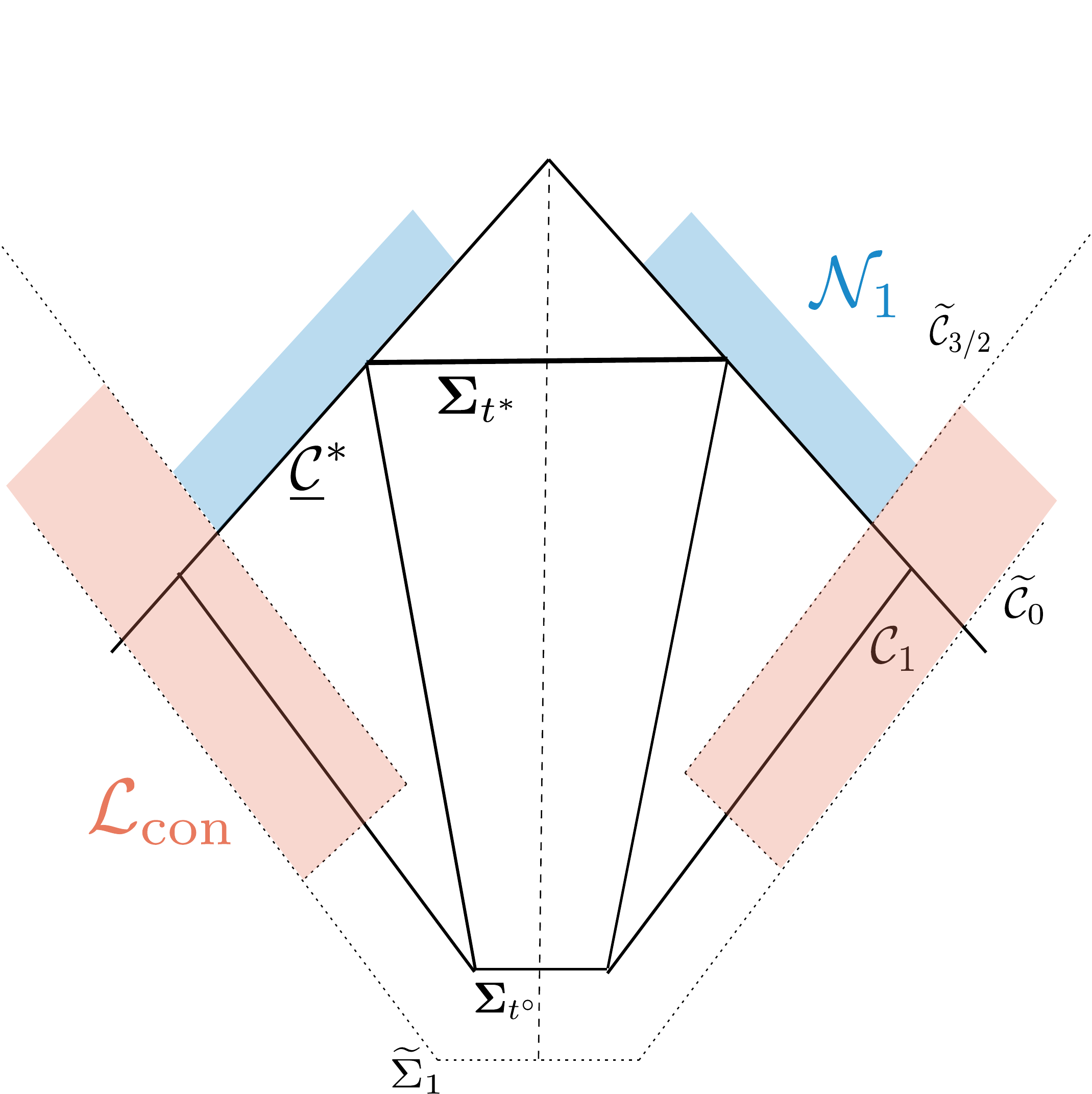}
  \includegraphics[height=4cm]{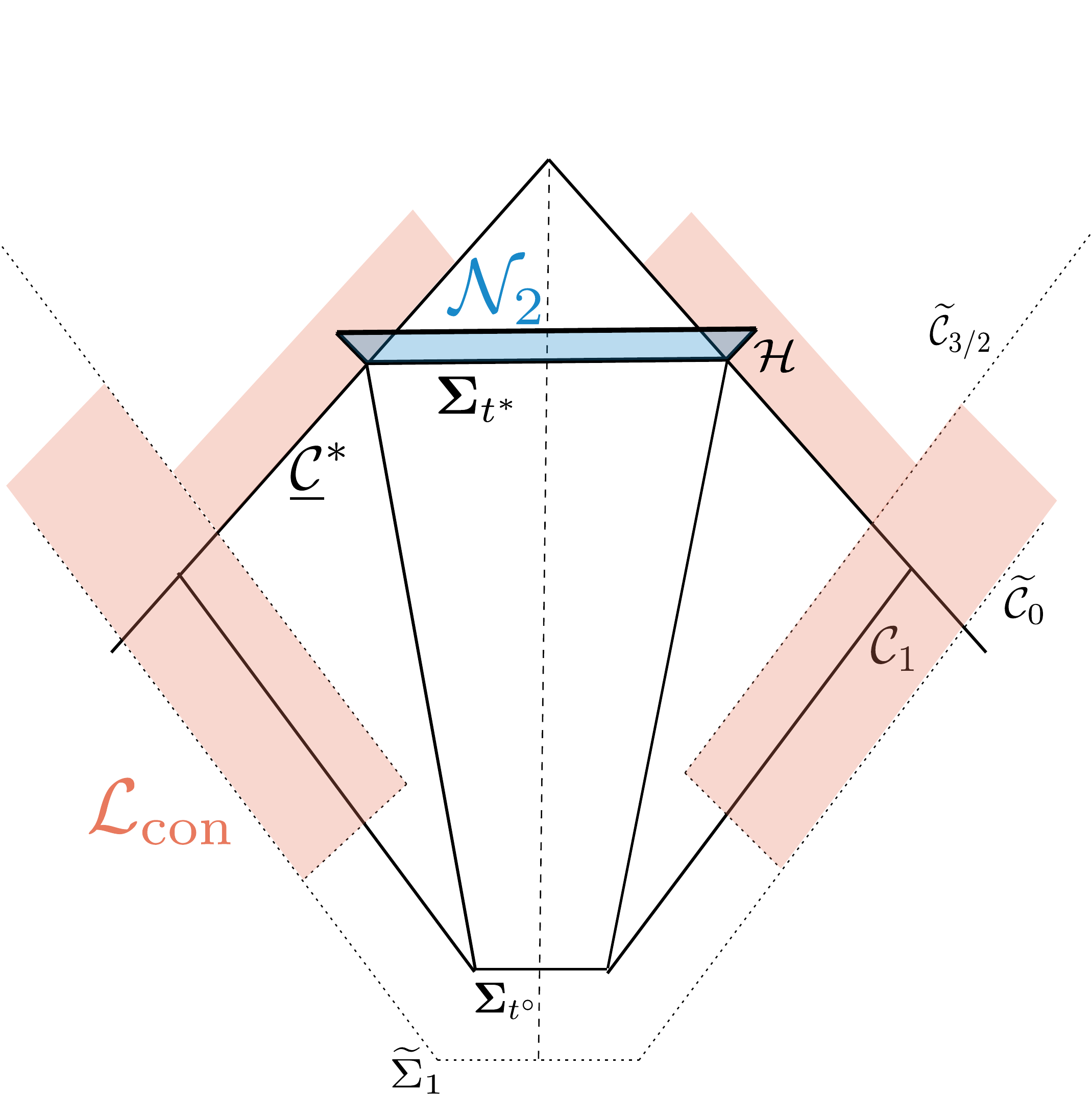}
  \includegraphics[height=4cm]{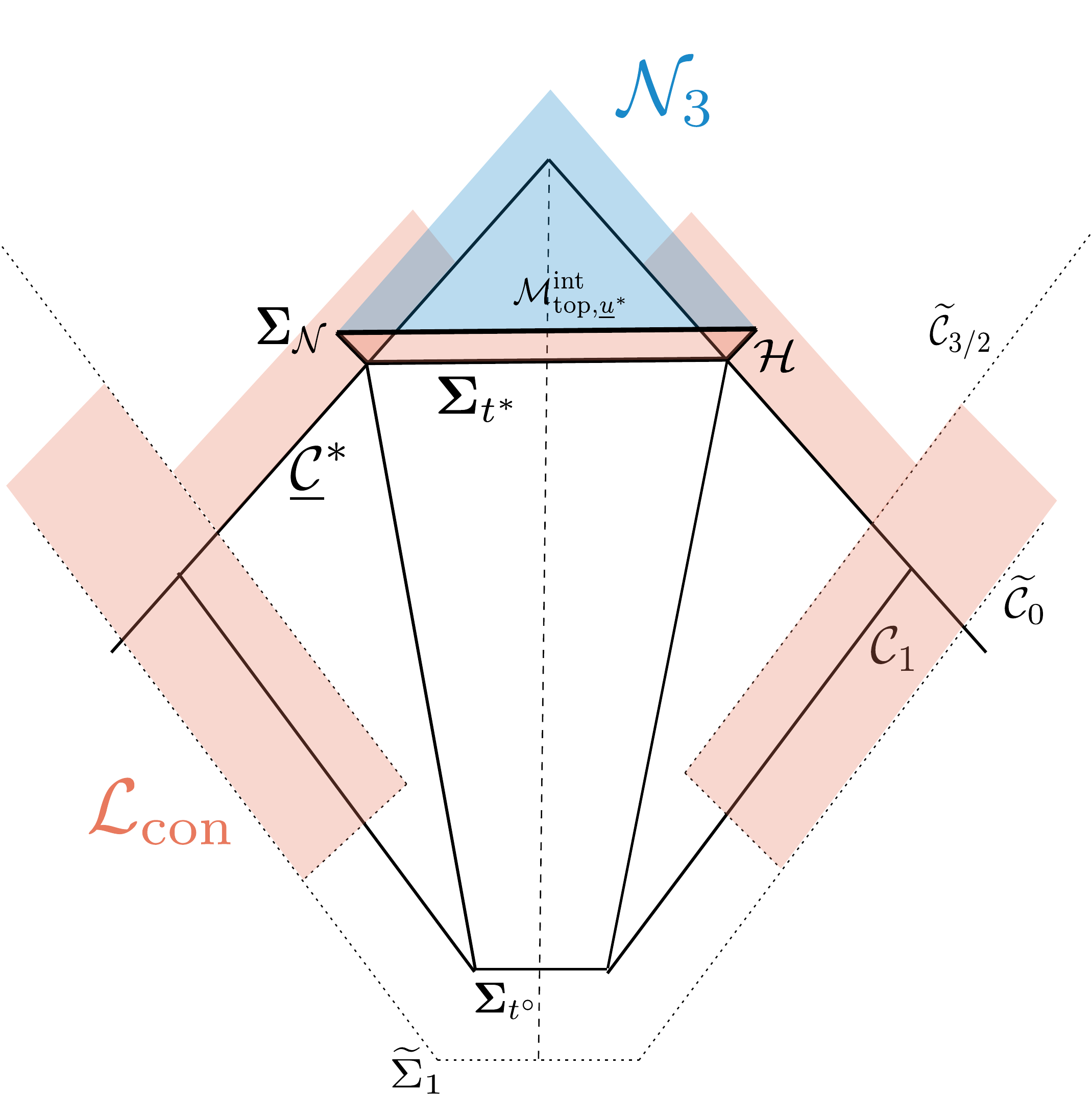}
  \caption{The spacetime $\NN$.}
  \label{fig:extension}
\end{figure}

By a comparison argument, the null hypersurface $\CCba\cap\le\{1\leq u \leq (1+\cc)\uba/2\ri\}$ can be foliated by a smooth non-degenerate geodesic foliation starting from the $2$-sphere $\CCba\cap\CCt_{3/2}$. Similarly, the initial layer null hypersurface $\CCt_{3/2}$ can be foliated by a smooth geodesic foliation starting form $\CCba\cap\CCt_{3/2}$. Applying the local existence result of~\cite{Luk12} to the induced characteristic initial data on the hypersurfaces
\begin{align*}
  \big(\CCba\cap\le\{u\leq (1+\cc)\uba/2\ri\}\big) \bigcup \CCt_{3/2},
\end{align*}
we deduce that there exists a smooth spacetime $\NN_1$ which strictly contains $\CCba\cap\le\{u\leq (1+\cc)\uba/2 \ri\}$ (see also Figure~\ref{fig:extension}).\\

We consider the smooth null hypersurface $\HH$ emanating from $S^\ast$ in $\NN_1$. Applying a spacelike-characteristic local existence result to $\Si_\tast \cup\HH$ (see~\cite[Section 6]{Czi.Gra19a}), we deduce that there exists a spacetime $\NN_2$ which contains $\Si_\tast\cup\HH$ in which the hypersurface $\Si_\tast$ is strictly included.\\ 

In particular, $\NN_2$ contains a maximal spacelike hypersurface $\Sigma_\NN$ strictly in the future of $\Si_\tast$ with boundary included in $\HH$. We define $\NN_3$ to be its smooth maximal globally hyperbolic development of $\Si_\NN$. By continuity, one can assume that on $\Si_\NN$, the harmonic coordinates and second fundamental form bounds from the Bootstrap Assumptions~\ref{BA:connint} and~\ref{BA:harmoSitast} hold. From the $\varep$-smallness of the initial data on $\Si_\NN$, a rescaling to a time-1 situation and a local existence result (see~\cite{Fou52}), one deduces that $\NN_3$ is diffeomorphic and $\varep$-close to a subregion of Minkowski space. Thus $\NN_3$ strictly contains $\MM^\intr_{\topp,\uba} \cap \NN_3$. Patching together $\LLc$, $\MM_\uba$, $\NN_1$, $\NN_2$ and $\NN_3$, we obtain a spacetime $\NN$ strictly extending $\MM_\uba$, as desired.


\paragraph{Step 2: Constructions of optical, affine parameters and time functions.}
Since in particular $\o(\uba) \in \MM_\uba \subset \NN$ is an interior point of $\NN$, the timelike geodesic $\o$ can be continued past $\o(\uba)$ and there exists $\de>0$ sufficiently small such that $\o(\uba+\de) \in \NN$.\\

Let denote by $\CCb^{\ast,\de}$ the incoming null cone backward emanating from $\o(\uba+\de)$, which local existence, smoothness and local foliation by geodesic affine parameter is guaranteed by the (null) geodesic equation and Cauchy-Lipschitz theorem. Since the cone $\CCba$ globally exists and can be foliated by a geodesic foliation going from $s=\uba$ to $s=1/2$, applying Cauchy-Lipschitz again provided that $\de>0$ is sufficiently small, the cone $\CCb^{\ast,\de}$ can be foliated by geodesic affine parameter going from $s=\uba+\de$ to $s=2/3$ and the associated null connection coefficients are close to the null connection coefficients associated to $\CCba$.\\

Applying the well-posedness Theorem~\ref{thm:canonical} for the canonical foliation on $\CCb^{\ast,\de}$ (see in particular Item~\ref{item:stabilitycanonical}), one deduces that for $\de>0$ sufficiently small the canonical foliation exists from $u=\uba+\de$ to $u=1$.\\


We denote by $\CC^{\de}_{u}$ the outgoing null cones backward emanating from the $2$-spheres of the canonical foliation on $\CCb^{\ast,\de}$. From Cauchy-Lipschitz, since the $2$-spheres on $\CCba$ and $\CCb^{\ast,\de}$ are close, the cones globally exist with geodesic affine parameter $\ub$ ranging from $\uba$ to $\cc^{-1}u$, and the associated null connection coefficients satisfy the Bootstrap Assumptions.\\

Since the $2$-spheres of the interfaces $\TT^\de$ and $\TT$ are close, one can apply the implicit function theorem of~\cite{Cho76} (see also a similar application in~\cite[Section 6]{Czi.Gra19a}) and we deduce that there exists a foliation by maximal hypersurfaces $(\Si^\de_{t})_{\too \leq t \leq (1+\cc)(\uba+\de)/2}$ with boundaries $\pr\Si^\de_t = S^\de_{u,\ub} \subset \TT^\de$, which is close to the analogous maximal foliation of $\MM^\intr_{\bott,\uba}$.\\

The existence, uniqueness and control of all centred conformal isomorphism from the new sphere $S^{\ast,\de}$ can be obtained from the Uniformisation Theorem~\cite[Theorem 3.1]{Kla.Sze19}. From induced coordinates on the new last maximal hypersurface $\Si^\de_{(1+\cc)(\uba +\de)/2}$, the hypothesis of Theorem~\ref{thm:globharmonics} are satisfied and there exists global harmonic coordinates on $\Si^\de_{(1+\cc)(\uba +\de)/2}$ adapted to any centred conformal isomorphism on $S^{\ast,\de}$. From the closeness to the constructions of the previous spacetime region, we deduce that the Bootstrap Assumptions are satisfied and we have obtained the desired spacetime region $\MM_{\uba+\de} \subset \NN \subset \MMt$.   

\subsection{Conclusions}\label{sec:conclusion}
In this section, we deduce from $\underline{U}^\ast = + \infty$ the main conclusions of Theorem~\ref{thm:mainv2}.\footnote{In this section, the transition parameter $\cc$ is fixed in $[\cco,(1+\cco)/2]$.}

\paragraph{Global time, optical and affine parameter functions}
First, we infer from $\underline{U}^\ast=+\infty$ that the central axis $\o(t)$ exists for $t=1$ to $+\infty$ and there exists global functions $u',\ub'$ such that its level sets is the last cones geodesic foliation defined in Section~\ref{sec:deflastconesfoliation}.\\

From the Bootstrap Assumptions and propagation of regularity, the functions ${^{(\uba)}\ub}, {^{(\uba)}u}$ and ${^{(\uba)}t}$ and their derivatives are bounded and equicontinuous uniformly in $\uba$ on each compact region $\{\ub'\leq C\}$. Thus, applying Arzel\`a-Ascoli theorem, one can deduce that there exists a sequence $\uba_n\to + \infty$ when $n\to\infty$ and functions $\ub,u,t$ such that
\begin{align*}
  {^{(\uba_n)}\ub} & \to \ub, & {^{(\uba_n)}u} & \to u, & {^{(\uba_n)}t} & \to t,
\end{align*}
in $\mathscr{C}^k$ topology in each compact region $\{\ub' \leq C\}$ and for all $k\geq 0$.\\  

We define the subregions
\begin{align*}
  \MM^\ext_\infty & := \le\{\cc^{-1} \leq \ub < +\infty,~1\leq u\leq\cc\ub\ri\},\\
  \MM^\intr_\infty & := \le\{\too \leq t < +\infty\ri\},\\
  \TT & := \le\{1 \leq u < + \infty,~u=\cc\ub\ri\}.
\end{align*}

Passing to the limits the definitions of $\ub$, $u$ and $t$, one deduces that they are respectively affine parameter, optical functions and maximal time function, and that on the interface $\TT := \{u=\cc\ub\}$, one has $t=\half (u+\ub)$. Moreover, passing to the limit the sup-norm estimates from the Bootstrap Assumptions, one deduces the bounds~(\ref{est:infty}) of Theorem~\ref{thm:mainv2}.\\

Using the bounds~\eqref{est:infty} and taking limits along the null cones $\CC_u$, a proper notion of \emph{future null infinity} $\scri^+$ can be obtained (see~\cite[Section 8]{Kla.Nic03}). From the fact that
\begin{align*}
  \le|\D_\elb\elb\ri| & \les |\omb| + |\xib| \les \varep \ub^{-1}u^{-3/2} \to 0 \quad \text{when}~\ub\to 0,
\end{align*}
that $\elb(u)=2$ and that $u$ ranges from $1$ to $+\infty$, one can deduce that $\scri^+$ is \emph{future geodesically complete}.

\paragraph{Hawking and Bondi mass}
We define the \emph{Hawking mass} $m$ in $\MM^\ext_\uba$ by\footnote{The identity is given by Gauss equation~(\ref{eq:Gauss}) and Gauss-Bonnet formula.}
\begin{align}
  \label{eq:defHawkingmass}
  \begin{aligned}
    m(u,\ub) & := \frac{r}{8\pi}\int_{S_{u,\ub}}\le(-\rho +\half\chih\cdot\chibh\ri) = \frac{r}{2} + \frac{r}{32\pi}\int_{S_{u,\ub}}\trchi\trchib.
  \end{aligned}
\end{align}
From the (improved) Bootstrap Assumptions~\ref{BA:curvext} and~\ref{BA:connext}, we have
\begin{align}\label{est:hawk}
  |m(u,\ub)| & \les \varep^2 u^{-2}.
\end{align}
Using the null structure equations~(\ref{eq:null}) and the Bootstrap Assumptions, one can prove that in $\MM^\ext_\uba$, the following bound holds
\begin{align*}
  \le|\el(m)\ri| & \les \varep \ub^{-3}u^{-3}. 
\end{align*}
Passing to the limit in $\uba$, the estimate still holds in $\MM^\ext_\infty$ and we deduce that $\el(m)$ is integrable along $\el$, and $m(u,\ub)$ admits a limit $M(u)$ when $\ub\to +\infty$. We call $M$ the \emph{Bondi mass}. From~\eqref{est:hawk}, one deduces
\begin{align*}
  |M(u)| & \les \varep^2 u^{-2},
\end{align*}
and in particular
\begin{align*}
  M(u) \to 0 \quad \text{when}~u\to +\infty,
\end{align*}
that is, the \emph{final Bondi mass is $0$}.\\

\paragraph{Bondi mass loss formula}
Using equation~\eqref{eq:Nd3rhoo}~\eqref{eq:Nd3chih} and~\eqref{eq:Nd3chibh} together with the Bootstrap Assumptions, we have
\begin{align*}
  \elb(r^3\rhoo) & = -\half r^3\overline{\chih\cdot\alb}  + O\le(\ub^{-1}\ri),\\
  \Nd_3(r\chih) & = -\half r\trchi \chibh + O\le(\ub^{-2}\ri),\\
  \Nd_3(r^2\chibh) & = -r^2\alb + O\le(1\ri),
\end{align*}
when $\ub\to+\infty$. Thus,
\begin{align*}
  \elb(m) & = -\half \elb\le(r^3\rhoo\ri) + \quar \elb\le(\overline{(r\chih)\cdot(r^2\chibh)}\ri) \\
          & = -\frac{r^3}{8} \overline{\trchi|\chibh|^2} + O\le(\ub^{-1}\ri),
\end{align*}
and passing to the limit when $\uba\to+\infty$ and taking the limit when $\ub\to+\infty$, we infer the following \emph{Bondi mass loss formula}
\begin{align*}
  \frac{\d}{\d u}M(u) & = \lim_{\ub\to+\infty} \le(-\frac{r}{64\pi}\int_{S_{u,\ub}} \trchi|\chibh|^2\ri).
\end{align*}

\paragraph{Angular momentum}
According to~\cite{Riz98}, we define the following local angular momentum
\begin{align}\label{eq:deflocalangular}
  ^{(\ell)}L(u,\ub) & := \frac{1}{8\pi r} \int_{S_{u,\ub}}\ze\cdot\OOEi,
\end{align}
for $\ell=1,2,3$.\footnote{Arguing similarly as before, one can deduce the existence of exterior rotation vectorfields in the limit $\uba\to+\infty$.} Using the Bootstrap Assumptions, we directly obtain
\begin{align}\label{est:Luub}
  \le|^{(\ell)}L(u,\ub)\ri| & \les \varep u^{-1/2}.
\end{align}
Using equations~\eqref{eq:Nd3xiOLD} and~\eqref{eq:Nd4OOE}, we have
\begin{align*}
  r^2\ze & = O(1), & r^{-1}\OOE & = O(1),\\
  \Nd_4(r^2\ze) & = O\le(\ub^{-3/2}\ri), & \Nd_4(r^{-1}\OOE) & = O(\ub^{-2}),
\end{align*}
when $\ub\to+\infty$. Thus, deriving~\eqref{eq:deflocalangular} by $\el$ gives
\begin{align*}
  \el(L) & = \half \el\le(\overline{(r^2\ze)\cdot r^{-1}\OOE}\ri) = O\le(\ub^{-3/2}\ri),
\end{align*}
and the bound still holds when passing to the limit $\uba\to+\infty$. We infer that $\el(L)$ is integrable along $\el$ and we define the \emph{angular momentum at null infinity} $^{(\ell)}L(u)$ to be the limit of $^{(\ell)}L(u,\ub)$ when $\ub \to \infty$. From~\eqref{est:Luub}, we directly obtain
\begin{align*}
  \le|^{(\ell)}L(u)\ri| & \les \varep u^{-1/2},
\end{align*}
and in particular
\begin{align*}
  ^{(\ell)}L(u) & \to 0, \quad \text{when}~u\to +\infty,
\end{align*}
that is, \emph{the final angular momentum is $0$}.\\

Using equations~\eqref{eq:Nd3ze} and~(\ref{eq:relOOEnull}), we have
\begin{align*}
  \Nd_3(r^2\ze) & = -2 r^2\Nd\omb -r^2\beb + O(\ub^{-1}), & \Nd_3(r^{-1}\OOE) & = O(\ub^{-1}),
\end{align*}
when $\ub\to+\infty$. Arguing as previously, we infer the following \emph{angular momentum evolution equation along null infinity}
\begin{align*}
  \frac{\d}{\d u} {^{(\ell)}L}(u) & = \lim_{\ub \to \infty} \le(\frac{1}{16\pi r}\int_{S_{u,\ub}} \le(-2\Nd\omb-\beb\ri)\cdot\OOEi\ri),
\end{align*}
for $\ell=1,2,3$. This finishes the proof of the conclusions of Theorem~\ref{thm:mainv2}.

\chapter{Global energy estimates in $\MM$}\label{sec:globener}
In this section, we perform the global energy estimates in $\MM$, which are used to estimate the curvature in the exterior and interior bottom regions respectively in Sections~\ref{sec:curvest} and~\ref{sec:planehypcurvest}. We first state the following proposition. Its proof is based on local energy estimates in the initial layers and is postponed to Section~\ref{sec:initlayerenergy}.
\begin{proposition}\label{prop:initenercorrec}
  Under the Bootstrap Assumptions and the $\varep$-closeness to Minkowski space assumptions from Section~\ref{sec:Minkowskilayer}, we have the following control of initial energy fluxes
  \begin{align}\label{est:initener}
    \int_{^{(\cc)}\Si_\too}P\cdot \Tf + \int_{^{(\cc)}\CC_1} P\cdot\el & \les \varep^2,
  \end{align}
  where $P$ are the following Bel-Robinson tensors defined respectively in the interior and exterior region $^{(\cc)}\MM^\intr_\bott$ and $^{(\cc)}\MM^\ext$ by
  \begin{align}\label{eq:defBelRob}
    \begin{aligned}
      Q\le(\Lieh_\TX\R\ri)(\KX,\KX,\KX), \quad Q\le(\Lieh_{\OOO}\R\ri)(\KX,\KX,\TX), \\
      Q\le(\Lieh_\OOO\Lieh_\OOO\R\ri)(\KX,\KX,\TX), \quad Q\le(\Lieh_\SX\Lieh_\TX\R\ri)(\KX,\KX,\KX), \quad Q\le(\Lieh_\OOO\Lieh_\TX\R\ri)(\KX,\KX,\TX),
    \end{aligned}
  \end{align}
  where $\TX,\SX,\KX,\OOO$ correspond to the approximate conformal Killing vectorfields defined respectively in the interior region in Section~\ref{sec:defKillingint} and in the exterior region in Section~\ref{sec:defKillingext}.\footnote{The bound~\eqref{est:initener} holds for all transition parameters $\cc$.}
\end{proposition}

The present section is dedicated to the proof of the following proposition.
\begin{proposition}\label{prop:enerestSTAB}
  Recall the results of Proposition~\ref{prop:initenercorrec}. Under these results, the Bootstrap Assumptions and for $\varep>0$ sufficiently small, there exists a transition parameter $\cco\leq\cc\leq(1+\cco)/2$ such that we have
  \begin{align}\label{est:resultglobener}
    \int_{^{(\cc)}\Si_t} P \cdot \Tf +\int_{\CC_u\cap{^{(\cc)}\MM}^\ext} P \cdot \el + \int_{^{(\cc)}\Si_t^\ext}P\cdot \Tf^{\ext} + \int_{\CCba\cap{^{(\cc)}\MM}^\ext} P\cdot\elb \;& \les \varep^2,
  \end{align}
  where $P$ denotes the Bel-Robinson tensor from~\eqref{eq:defBelRob}.
\end{proposition}

\begin{remark}\label{rem:alltransitionhypcontrol}
  The transition parameter $\cc$ is determined by a mean value argument in Section~\ref{sec:meanvalue}. The control for the curvature on the hypersurfaces associated to \emph{all} the transition parameters $\cco \leq \cc \leq (1+\cco)/2$ is obtained in Section~\ref{sec:alltransparam}.
\end{remark}


\begin{figure}[!h]
  \centering
  \includegraphics[height=8cm]{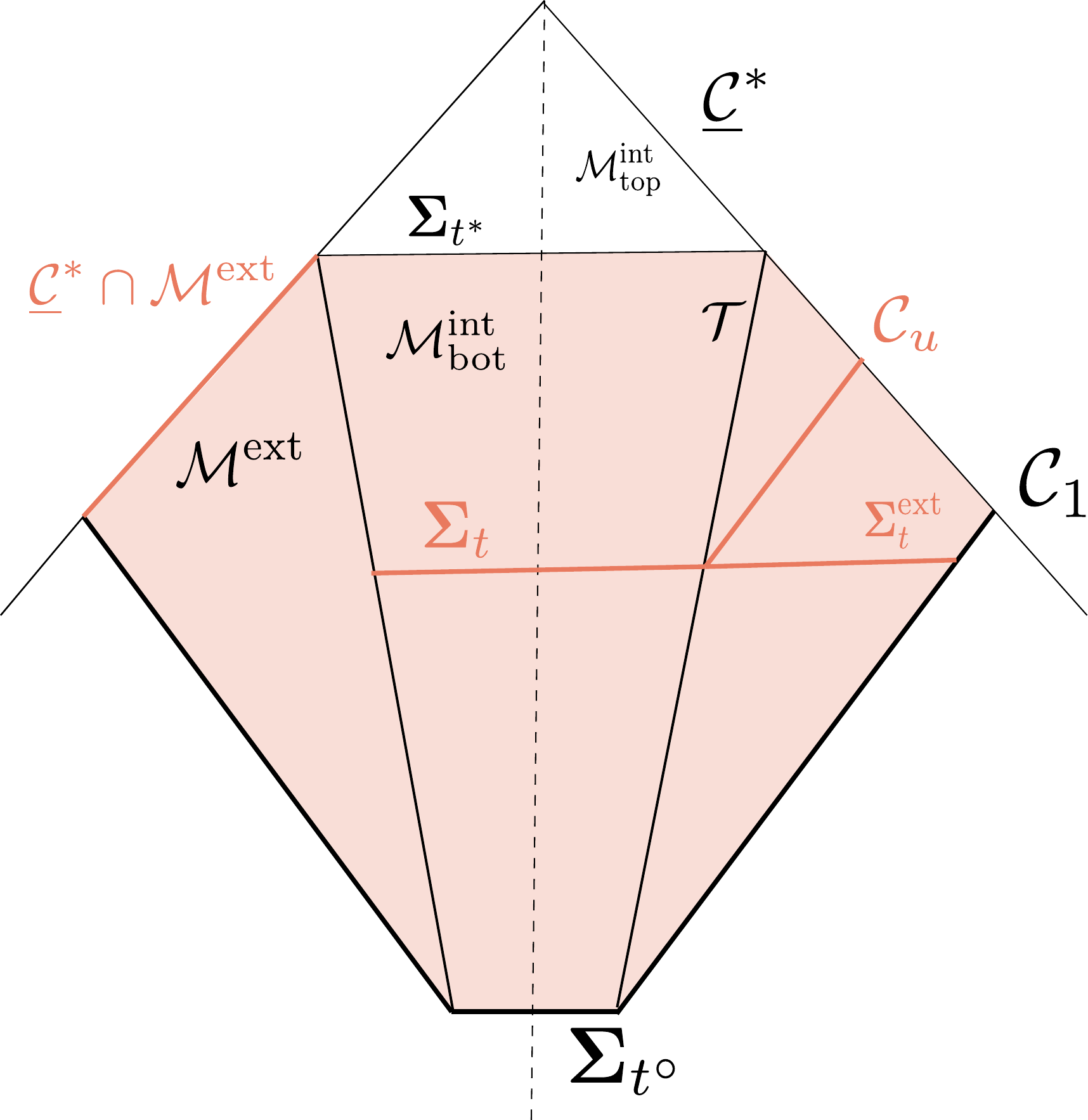}
  \caption{Global energy estimates in $\MM$.}
  \label{fig:globener}
\end{figure}

The proof of Proposition~\ref{prop:enerestSTAB} follows from the following energy estimate
\begin{align*}
  \int_{\Si_t}P\cdot\Tf +\int_{\CC_u\cap\MM^\ext} P \cdot \el+  \int_{\Si_t^\ext} P \cdot \Tf^\ext+ \int_{\CCba\cap\MM^\ext} P\cdot\elb & \les \int_{\Si_{\too}}P\cdot\Tf + \int_{\CC_1}P\cdot\el + \EE^\TT + \EE^\intr + \EE^\ext,
\end{align*}
which holds for all $\too \leq t\leq \tast$ and for all $1 \leq u \leq \cc\uba$, for all contracted and commuted Bel-Robinson tensors $P$ of~\eqref{eq:defBelRob}, and where the nonlinear error terms $\EE$ are defined below.\\ 

The nonlinear error terms are decomposed as follows
\begin{align*}
  \EE^\TT & = \EE^\TT_1 + \EE^\TT_2, & \EE^\intr & = \EE^\intr_1 + \EE^\intr_2, & \EE^\ext & = \EE^\ext_1 + \EE^\ext_2,
\end{align*}
where have the following definitions for each respective factor.
\begin{itemize}
\item We have
  \begin{align*}
    \EE^\TT_1 & := \int_{\TT}\le(P^\ext\cdot N^\TT -P^\intr \cdot N^\TT\ri),
  \end{align*}
  where $N^\TT$ denotes the inward-pointing unit normal to the timelike hypersurface $\TT$ and where
  \begin{align*}
    P \in\le\{Q\le(\Lieh_\TX\R\ri)(\KX,\KX,\KX), \; Q\le(\Lieh_{\OOO}\R\ri)(\KX,\KX,\TX)\ri\}.
  \end{align*}
\item We have
  \begin{align*}
    \EE^\TT_2 & := \int_{\TT}\le(P^\ext\cdot N^\TT -P^\intr \cdot N^\TT\ri),
  \end{align*}
  where
  \begin{align*}
    P \in\le\{Q\le(\Lieh_\OOO\Lieh_\OOO\R\ri)(\KX,\KX,\TX), \; Q\le(\Lieh_\SX\Lieh_\TX\R\ri)(\KX,\KX,\KX), \; Q\le(\Lieh_\OOO\Lieh_\TX\R\ri)(\KX,\KX,\TX)\ri\}.
  \end{align*}
\item We have
  \begin{align*}
    \EE^\intr_1 & := \int_{\MM^\intr_\bott}\DIV\le(P\ri),
  \end{align*}
  where
  \begin{align*}
    P \in\le\{Q\le(\Lieh_\TI\R\ri)(\KI,\KI,\KI), \; Q\le(\Lieh_{\OOI}\R\ri)(\KI,\KI,\TI)\ri\}.
  \end{align*}
\item We have
  \begin{align*}
    \EE^\intr_2 & := \int_{\MM^\intr_\bott}\DIV\le(P\ri),
  \end{align*}
  where
  \begin{align*}
    P \in\bigg\{ & Q\le(\Lieh_\OOI\Lieh_\OOI\R\ri)(\KI,\KI,\TI), \; Q\le(\Lieh_\SI\Lieh_\TI\R\ri)(\KI,\KI,\KI),\\
    & \; Q\le(\Lieh_\OOI\Lieh_\TI\R\ri)(\KI,\KI,\TI)\bigg\}.
  \end{align*}
\item We have
  \begin{align*}
    \EE^\ext_1 & := \int_{\MM^\ext}\DIV(P),
  \end{align*}
  where
  \begin{align*}
    P \in\le\{Q\le(\Lieh_\TE\R\ri)(\KE,\KE,\KE), \; Q\le(\Lieh_{\OOE}\R\ri)(\KE,\KE,\TE)\ri\}.
  \end{align*}
\item We have
  \begin{align*}
    \EE^\ext_2 & := \int_{\MM^\ext}\DIV(P),
  \end{align*}
  where
  \begin{align*}
    P \in\bigg\{ & Q\le(\Lieh_\OOE\Lieh_\OOE\R\ri)(\KE,\KE,\TE), \; Q\le(\Lieh_\SE\Lieh_\TE\R\ri)(\KE,\KE,\KE), \\
    & \; Q\le(\Lieh_\OOE\Lieh_\TE\R\ri)(\KE,\KE,\TE)\bigg\}.
  \end{align*}  
\end{itemize}


This section is dedicated to proving that, under the Bootstrap Assumptions, there exists a transition parameter $\cc$ such that we have
\begin{align}\label{est:controlEEM}
  \EE^\TT_{\leq 2}  + \EE^\intr_{\leq 2} + \EE^\ext_{\leq 2} & \les (D\varep)^3.
\end{align}
The result of Proposition~\ref{prop:enerestSTAB} then directly follows, provided that $\varep>0$ is sufficiently small.\\


To control the interior and exterior error terms in Sections~\ref{sec:errorintr} and~\ref{sec:exterrest}, we make the following additional bootstrap assumption. In view of the above, it will directly be improved by the energy estimate once the control~\eqref{est:controlEEM} of the error terms has been obtained.
\begin{BA}\label{BA:fluxglobener}
  Assume that the above fluxes through the hypersurfaces $\Si_t, \CC_u\cap\MM^\ext$ and $\Si_t^\ext$ satisfy the following additional bootstrap bounds
  \begin{align*}
    \int_{\Si_t} P \cdot \Tf +\int_{\CC_u\cap\MM^\ext} P \cdot \el + \int_{\Si_t^\ext} P \cdot \Tf^\ext & \leq (D\varep)^2,
  \end{align*}
  where $P$ denotes the contracted Bel-Robinson tensors~(\ref{eq:defBelRob}).
\end{BA}


\begin{remark}[Decay in the exterior region]
  From the point of view of the decay, the most difficult error term to treat is the exterior term $\EE^\ext_{\leq 2}$, whose estimate is the crux of~\cite{Chr.Kla93} and~\cite{Kla.Nic03}.
  To this end, one has to exhibit a null structure which pairs curvature/connection coefficients (decomposed in the null directions) with compensating decay rates.
  In this paper, we choose the same commutating and multiplying vectorfields as in~\cite{Chr.Kla93} and~\cite{Kla.Nic03}, and the structure of the error terms in the exterior region is formally identical to the one in these books. Since our vectorfields are constructed upon a geodesic-null foliation which is different from the maximal-null and the double-null foliations of respectively~\cite{Chr.Kla93} and~\cite{Kla.Nic03}, the deformation tensors of our approximate conformal Killing fields satisfy slightly different decay rates.\footnote{This is due to the presence of $\xib$, which vanishes in the double-null setting of~\cite{Kla.Nic03}.} We analyse in detail these differences in Sections~\ref{sec:esterrEE1}, \ref{sec:esterrEE12}, \ref{sec:esterrEE22} and show that the error terms are still integrable in the spacetime region $\MM^\ext$.
\end{remark}

\section{Estimates for the interface error terms $\EE^\TT$}
\subsection{The mean value argument}\label{sec:meanvalue}
In this section, we obtain $L^2(\TT)$ bounds for the spacetime curvature tensor and its first and second derivatives. This is done by a \emph{mean value argument} which selects a particular transition parameter $\cc$.\\

We consider the spacetime region
\begin{align*}
  \DD & := \le\{(u,\ub)\,:\, \cco\ub \leq u \leq (1+\cco)\ub/2\ri\} \cap J^+(\Sit_{3}).
\end{align*}

\begin{figure}[!h]
  \centering
  \includegraphics[height=6cm]{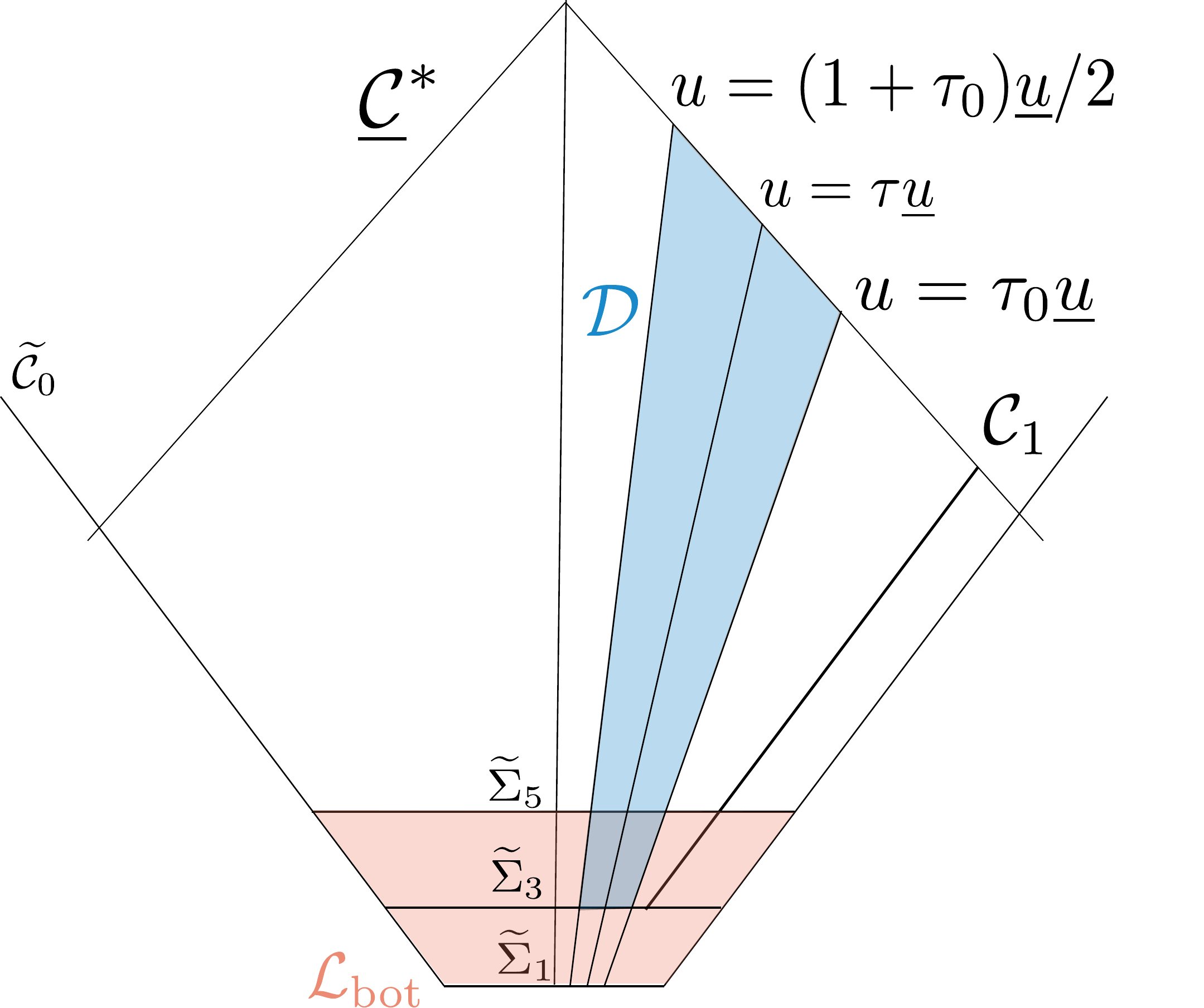}
  \caption{The mean value argument.}
  \label{fig:meanvalueTT}
\end{figure}

From the Bootstrap Assumptions~\ref{BA:curvext} when the transition parameter is set to $\cc=(1+\cco)/2$, we have the following $L^2(\DD)$ control for the spacetime curvature tensor (we recall that in the region $\DD$, we have $u\simeq \ub \simeq r \simeq t$)
\begin{align}\label{est:L2DD}
  \begin{aligned}
    \int_{\DD} \ub^{-1-1/2}\ub^4\le|(\ub\D)^{\leq 2}\R\ri|^2 & \les \int_{\DD}\ub^{-1-1/2}\ub^4\le|\Ndt^{\leq 2}R\ri|^2 + \int_{\DD}\ub^{-1-1/2}\ub^4\le|\Ndt^{\leq 1}R\cdot\ub\Ga\ri|^2 \\
                                                             & \quad + \int_{\DD}\ub^{-1-1/2}\ub^4\le|R \cdot\Ndt(\ub\Ga)\ri|^2 + \int_\DD\ub^{-1-1/2}\ub^4\le|R\cdot\ub\Ga\cdot\ub\Ga\ri|^2 \\  
                                                             & \les \le(\RR^{\ext}_{\leq 2,\gao}+ \overline{\mathfrak{R}}^\ext_{\leq 2}\ri)^2\le(1  + \mathfrak{O}^\ext_{\leq 2} + \overline{\mathfrak{O}}^\ext_{\leq 1} \ri)^4 + \le(\mathfrak{R}^\ext_{\leq 1}+\overline{\mathfrak{R}}^\ext_{\leq 1}\ri)^2\le(1+\OO^\ext_{\leq 3, \gao}\ri)^2\\
    & \les (D\varep)^2,
    \end{aligned}
  \end{align}
  where the norm $|\cdot|$ is taken with respect to the orthonormal frame associated to the null pair $(\elb,\el)$,\footnote{With the notations of Definition~\ref{def:framenorm}, we define $e_0 := \half(\elb+\el)$ and the $e_i$ to be any (local) vectorfields such that $(e_\mu)$ forms an orthonormal frame.} and where
  \begin{align*}
    R & \in \le\{\alb,\beb,\rho,\sigma,\be,\al\ri\},\\
    \Ga & \in \le\{\trchi,\trchib,\chih,\chibh,\ze,\eta,\etab,\omb,\xib\ri\}, \\
    \Ndt & \in \le\{\ub\Nd_4,\ub\Nd_3,\ub\Nd\ri\}.
  \end{align*}
  

The region $\DD$ can be foliated as follows
\begin{align*}
  \DD & = \cup_{\cco \leq \cc \leq (1+\cco)/2} \le\{(u,\ub)~:~ u = \cc\ub \ri\}\cap J^+(\Sit_{3}).
\end{align*}
With respect to this foliation, we have the following coarea formula
\begin{align*}
  & \; \int_{\cco}^{(1+\cco)/2}\le(\int_{\{u=\cc\ub\}\cap J^+(\Sit_3)}\big(\ub^{-1-1/2}\ub^4\le|(\ub\D)^{\leq 2}\R\ri|^2\big)\ub \ri)\,\d \cc \\
  \les & \; \int_{\DD} \ub^{-1-1/2}\ub^4\le|(\ub\D)^{\leq 2}\R\ri|^2 \\
  \les & \; (D\varep)^2. 
\end{align*}
Using the mean value theorem, we deduce that there exists a transition parameter $\cc\in[\cco,(1+\cco)/2]$ such that
\begin{align*}
  \int_{\{u=\cc\ub\}\cap J^+(\Sit_3)}\ub^{-1/2}\ub^4\le|(\ub\D)^{\leq 2}\R\ri|^2  & \les (D\varep)^2.
\end{align*}

This rewrites
\begin{align}\label{est:TTbdy}
 \int_{^{(\cc)}\TT}\ub^{-1/2}\ub^4\le|(\ub\D)^{\leq 2}\R\ri|^2 & \les (D\varep)^2.
\end{align}

\begin{remark}
  The integral bound~\eqref{est:TTbdy} on $\TT\cap J^-(\Sit_3) \subset \LL_{\bott}$ can be obtained using the bottom initial layer assumptions and comparison arguments. Details are left to the reader.
\end{remark}

\begin{remark}\label{rem:optregvsoptdecay}
  One cannot obtain a control on a timelike boundary such as $\TT$ for the curvature components, without either loosing derivatives or decay. For the error terms of $\EE^\TT$, we need to control the highest order terms $\D^2\R$ in $L^2(\TT)$ and are bound to use the above mean value argument which looses $\ub^{1/2}$ decay. However, this loss of decay is (more than) compensated by the decay rates for the difference between interior and exterior Killing fields, see Section~\ref{sec:estEETT12}.
\end{remark}

In the rest of this section, the transition parameter $\cc$ is fixed by the above mean value argument and we drop its label in the notations.

\subsection{Control of $\EE^\TT_1$ and $\EE^\TT_2$}\label{sec:estEETT12}
We have
\begin{align*}
  \EE^\TT_1 & := \int_\TT \le(Q(\Lieh_\TE\R)(\KE,\KE,\KE)-Q(\Lieh_\TI\R)(\KI,\KI,\KI)\ri)\cdot N^\TT \\
            & \quad + \int_\TT\le(Q(\Lieh_\OOE\R)(\KE,\KE,\TE)-Q(\Lieh_\OOI\R)(\KI,\KI,\TI)\ri)\cdot N^\TT,
\end{align*}
where we recall that $N^\TT$ denotes the inward-pointing unit normal to $\TT$.\\

We write schematically 
\begin{align}\label{eq:diffQ}
  \begin{aligned}
  & Q(\Lieh_\TE\R)(\KE,\KE,\KE)-Q(\Lieh_\TI\R)(\KI,\KI,\KI) \\
  = & \;\Lieh_\TE\R \cdot \le(\le(\Lieh_\TE-\Lieh_\TI\ri)\R\ri) \cdot (\KE,\KE,\KE) \\
  & + \le(\le(\Lieh_\TE-\Lieh_\TI\ri)\R\ri) \cdot \Lieh_\TI\R \cdot (\KE,\KE,\KE)\\
  & + Q(\Lieh_\TI\R)(\KE-\KI,\KE,\KE) \\
  & + Q(\Lieh_\TI\R)(\KI,\KE-\KI,\KE) \\
  & + Q(\Lieh_\TI\R)(\KI,\KI,\KE-\KI).
  \end{aligned}
\end{align}

We start with the estimates for the last three terms of~\eqref{eq:diffQ}. From the Bootstrap Assumptions~\ref{BA:TTKilling} at $\TT$, we recall that we have
\begin{align}
  \label{est:DiffKillingTT}
  \begin{aligned}
  \le|\TE-\TI\ri| & \les D\varep t^{-3/2},\\
  \le|\SE-\SI\ri| & \les D\varep t^{-1/2},\\
  \le|\KE-\KI\ri| & \les D\varep t^{1/2}.
  \end{aligned}
\end{align}

Using that from the mild Bootstrap Assumptions~\ref{BA:mildKillingMMintbot} we have $|\KI| \les t^2$ at $\TT$ and using the $L^2(\TT)$ bounds~(\ref{est:TTbdy}) -- where we recall that $t\simeq \ub$ --, we infer
\begin{align*}
  \int_{\TT}\le|Q(\Lieh_\TI\R)(\KI,\KI,\KE-\KI)\ri| & \les \int_\TT \le|Q(\Lieh_\TI\R)\ri|\le(|\KE-\KI|||\KI|^2\ri) \\
            & \les D\varep\int_\TT\le|Q(\Lieh_\TI\R)\ri|t^{9/2} \\
            & \les D\varep\int_\TT t^{-1}t^{4-1/2}\le|(t\D)^{\leq 1} \R\ri|^2 \\
            & \les (D\varep)^3.
\end{align*}

For the first two terms of~\eqref{eq:diffQ}, we have schematically
\begin{align*}
  \le(\Lieh_\TE-\Lieh_\TI\ri)\R & = \D_{\TE-\TI}\R + \le(\D\TE-\D\TI\ri)\cdot\R.
\end{align*}
Using~\eqref{est:DiffKillingTT} and that by the Bootstrap Assumptions~\ref{BA:TTKilling} and~\ref{BA:curvext}, we respectively have
\begin{align*}
  \norm{\D\TE-\D\TI}_{L^4(\pr\Si_t)} & \les D\varep t^{-2+\gao}, \\
  \norm{\R}_{L^4(\pr\Si_t)} & \les D\varep t^{-3},
\end{align*}
we obtain
\begin{align*}
  & \int_{\pr\Si_t}\le|\Lieh_\TE\R \cdot \le(\le(\Lieh_\TE-\Lieh_\TI\ri)\R\ri) \cdot (\KE,\KE,\KE)\ri| \\
  \les & \; \norm{\D\R}_{L^2(\pr\Si_t)}^2 \norm{\TE-\TI}_{L^\infty(\pr\Si_t)}\norm{\KE}_{L^\infty(\pr\Si_t)}^3\\
  & + \norm{\D\R}_{L^2(\pr\Si_t)}\norm{\R}_{L^4(\pr\Si_t)}\norm{\D\TE-\D\TI}_{L^4(\pr\Si_t)}\norm{\KE}_{L^\infty(\pr\Si_t)}^3 \\
  \les & \; (D\varep) t^{6-3/2}\norm{\D\R}_{L^2(\pr\Si_t)}^2 + (D\varep)^2 t^{6-5 +\gao}\norm{\D\R}_{L^2(\pr\Si_t)}.
\end{align*}
Integrating on $\TT$ using a coarea formula, Cauchy-Schwartz, and the estimate~\eqref{est:TTbdy}, one obtains
\begin{align*}
  & \int_\TT \le|\Lieh_\TE\R \cdot \le(\le(\Lieh_\TE-\Lieh_\TI\ri)\R\ri) \cdot (\KE,\KE,\KE)\ri| \\
  \les & \; (D\varep)\int^\tast_{\too} t^{-1} t^{4-1/2}\norm{(t\D)^{\leq 1}\R}_{L^2(\pr\Si_t)}^2 \,\d t + (D\varep)^2\int_{\too}^\tast t^{\gao}\norm{(t\D)^{\leq 1}\R}_{L^2(\pr\Si_t)}\,\d t \\
  \les & \; (D\varep)^3 + (D\varep)^2\le(\int_{\too}^\tast t^2\norm{(t\D)^{\leq 1}\R}^2_{L^2(\pr\Si_t)}\,\d t\ri)^{1/2}\\
  \les & \; (D\varep)^3,
\end{align*}
where from the second to third line, we used that $\gao < 1/4$.\\

Thus, we have obtained
\begin{align*}
  \le|\int_\TT \le(Q(\Lieh_\TE\R)(\KE,\KE,\KE)-Q(\Lieh_\TI\R)(\KI,\KI,\KI)\ri)\cdot N^\TT \ri| & \les (D\varep)^3.
\end{align*}

The estimates for the second term of $\EE^\TT_1$ follow similarly, using that on $\TT$ by the Bootstrap Assumptions~\ref{BA:TTKilling}, we have
\begin{align*}
  \norm{\D\OOE- \D\OOI}_{L^4(\pr\Si_t)} & \les D\varep t^{-1+\gao},
\end{align*}
and this finishes the estimates for $\EE^\TT_1$.\\

The estimates for the highest order error term $\EE^\TT_2$ follow along the same estimates as above, using that by the Bootstrap Assumptions~\ref{BA:TTKilling}, the following bounds hold on $\TT$
\begin{align*}
  \norm{\D\SE-\D\SI}_{L^4(\pr\Si_t)} & \les D\varep t^{-1+\gao},
\end{align*}
and
\begin{align*}
  \norm{(t\D)^{\leq 1}\le(\D\TE-\D\TI\ri)}_{L^4(\pr\Si_t)} & \les D\varep t^{-2+\gao}, \\
  \norm{(t\D)^{\leq 1}\le(\D\OOE-\D\OOI\ri)}_{L^4(\pr\Si_t)} & \les D\varep t^{-1+\gao}.
\end{align*}
This finishes the proof of the desired estimates for $\EE^\TT_2$.
\section{Estimates for the interior error terms $\EE^\intr$}\label{sec:errorintr}
\subsection{Estimates for $\EE^\intr_1$}\label{sec:errorintr1}
We have
\begin{align*}
  \EE^\intr_1 & := \int_{\MM^\intr_\bott} \DIV \le(Q\le(\Lieh_\TI\R\ri)\ri)(\KI,\KI,\KI) \\
        & \quad + \int_{\MM^\intr_\bott} \DIV \le(Q\le(\Lieh_\OOI\R\ri)\ri)(\KI,\KI,\TI) \\
        & \quad +\int_{\MM^\intr_\bott} Q\le(\Lieh_\TI\R\ri)\cdot{^{(\KI)}\pih}\cdot(\KI,\KI) \\
        & \quad + \int_{\MM^\intr_\bott} Q\le(\Lieh_\OOI\R\ri)\cdot{^{(\KI)}\pih}\cdot\KI\cdot\TI \\
        & \quad + \int_{\MM^\intr_\bott} Q\le(\Lieh_\OOI\R\ri)\cdot{^{(\TI)}\pih}\cdot\KI\cdot\KI.
\end{align*}

To treat the last three error terms, we use that from the mild Bootstrap Assumptions~\ref{BA:mildKillingMMintbot} and the Bootstrap Assumptions~\ref{BA:intKill}, we have respectively
\begin{align*}
  \le|\TI\ri| & \les 1, &  \le|\KI\ri| & \les t^2 , \\
  \le|^{(\TI)}\pih\ri| & \les D\varep t^{-5/2}, & \le|^{(\KI)}\pih\ri| & \les D\varep t^{-1/2},
\end{align*}
where the norms are taken with respect to the maximal frame (\emph{i.e.} $e_0 = \Tf$).\\

Using the Bootstrap Assumptions~\ref{BA:fluxglobener} on the flux of the contracted Bel-Robinson tensor, its positivity properties (see~\cite[Chapter 7]{Chr.Kla93}), and a coarea formula, we obtain
\begin{align*}
  & \int_{\MM^\intr_\bott}Q(\Lieh_\OOI\R)\cdot{^{(\TI)}\pih}\cdot\KI\cdot\KI \\
  \les &\; \int_{\too}^\tast \int_{\Si_t}Q(\Lieh_\OOI\R)\le(\TI,\TI,\KI,\KI\ri)\le|^{(\TI)}\pih\ri| \,\d t\\
  \les & \; (D\varep)^3\int_{\too}^\tast t^{-5/2}\,\d t \\
  \les & (D\varep)^3,
\end{align*} 
and the other terms follow similarly.\\

For the first term, we have using Bianchi equations and the formula from~\cite[p. 141]{Chr.Kla93}
\begin{align*}
  \DIV\le(\Lieh_\TI\R\ri) & = {^{(\TI)}\pih}\cdot\D\R + \D{^{(\TI)}\pih}\cdot\R.
\end{align*}

Therefore, using~\cite[p. 137]{Chr.Kla93}, we have
\begin{align*}
  \DIV Q(\Lieh_\TI\R) & = \Lieh_\TI\R \cdot \le(\DIV\le(\Lieh_\TI\R\ri)\ri) \\
                      & = \Lieh_\TI\R \cdot \D\R \cdot {^{(\TI)}\pih} + \Lieh_\TI\R\cdot \R \cdot \D{^{(\TI)}\pih}.
\end{align*} 
The first term is handled as previously. For the second term, we use that, by the Bootstrap Assumptions~\ref{BA:fluxglobener} on the flux of the contracted Bel-Robinson tensors, the positivity properties of the Bel-Robinson tensor, and a coarea formula, one has
\begin{align*}
  \norm{t^{-1/2-\ga}t^3\Lieh_\TI\R}_{L^2(\MM^\intr_\bott)} & \les \le(\int_{\too}^\tast \int_{\Si_t} t^{-1-2\ga} t^6\le|\Lieh_\TI\R\ri|^2\,\d t\ri)^{1/2} \\
                                                           & \les \le(\int_{\too}^\tast \int_{\Si_t} t^{-1-2\ga} Q\le(\Lieh_\TI\R\ri)(\KI,\KI,\KI,\TI) \,\d t\ri)^{1/2} \\
                                                           & \les (D\varep) \le(\int_{\too}^\tast t^{-1-2\ga}\,\d t\ri)^{1/2} \\
  & \les_\ga (D\varep),
\end{align*}
for all $\ga>0$. We also use that by the Bootstrap Assumptions~\ref{BA:curvint} and the Bootstrap Assumptions~\ref{BA:intKill}, we have respectively
\begin{align*}
  \norm{t^{7/2}\R}_{L^\infty(\MM^\intr_\bott)} & \les D\varep,\\
  \norm{t^{-1/2-\gao}t^2\D^{(\TI)}\pih}_{L^2(\MM^\intr_\bott)} & \les D\varep. 
\end{align*}
This gives
\begin{align*}
  & \le|\int_{\MM^\intr_\bott}\Lieh_\TI\R\cdot\R\cdot\D^{(\TI)}\pih\cdot(\KI,\KI,\KI)\ri| \\
   \les & \; (D\varep)\int_{\MM^\intr_\bott} t^{5/2}\le|\Lieh_\TI\R\ri|\le|\D^{(\TI)}\pih\ri| \\
  \les & \; (D\varep)\norm{t^{-1/2-\ga}t^3\Lieh_\TI\R}_{L^2(\MM^\intr_\bott)}\norm{t^{-1/2-\gao}t^2\D{^{(\TI)}\pih}}_{L^2(\MM^\intr_\bott)} \\
  \les & \; (D\varep)^3,
\end{align*}
provided that $\ga$ is such that $\ga< 3/2-\gao$.\\

The remaining term of $\EE^\intr_1$ follows along the same lines, using the Bootstrap Assumptions~\ref{BA:intKill} for $\D^{\leq 1}{^{(\OOI)}\pih}$. This finishes the control of $\EE^\intr_1$.

\subsection{Estimates for $\EE^\intr_2$}\label{sec:estEEint2}
We treat the following error term
\begin{align*}
  \int_{\MM^\intr_\bott} \DIV\le(Q\le(\Lieh_\OOI\Lieh_\OOI\R\ri)\ri)\cdot(\KI,\KI,\TI).
\end{align*}
Using~\cite[p. 141]{Chr.Kla93} and that schematically $\Lieh_\OOI = \D_{\OOI} + \D\OOI \cdot$, we have 
\begin{align*}
  \DIV\le(\Lieh_\OOI\Lieh_\OOI\R\ri) & = \D^2{^{(\OOI)}\pih} \cdot \R \cdot \OOI + \D{^{(\OOI)}\pih} \cdot \D\R \cdot \OOI + {^{(\OOI)}\pih}\cdot\D^2\R\cdot\OOI \\
                                     & \quad + \D{^{(\OOI)}\pih} \cdot \R \cdot \D\OOI + \D^{(\OOI)}\pih \cdot \R \cdot \D\OOI + ^{(\OOI)}\pih\cdot\D\R\cdot\D\OOI \\
  & \quad + ^{(\OOI)}\pih\cdot\R\cdot\D^2\OOI.
\end{align*}
Using the $L^\infty(\MM^\intr_\bott)$ and $L^2(\MM^\intr_\bott)$ of the Bootstrap Assumptions~\ref{BA:curvint} and~\ref{BA:intKill} for $\R$ and $\OOI$ respectively, we check that
\begin{align}\label{est:L2DIVLIEHOOI}
  \norm{t^{4-\gao} \DIV\le(\Lieh_\OOI\Lieh_\OOI\R\ri)}_{L^2(\MM^\intr_\bott)} & \les (D\varep)^2. 
\end{align}
\begin{remark}
  For the second term, we use the $L^\infty_tL^6(\Si_t)$ estimates of the Bootstrap Assumptions~\ref{BA:curvint} and~\ref{BA:intKill} and a coarea formula as follows
  \begin{align*}
    \norm{t^{4-\gao}\D^{(\OOI)}\pih\cdot\D\R\cdot\OOI}_{L^2(\MM^\intr_\bott)} & \les \le(\int_{\too}^\tast \norm{t^{5-\gao}\D^{(\OOI)}\pih\cdot\D\R}^2_{L^2(\Si_t)}\,\d t\ri)^{1/2}\\
                                                                              & \les \le(\int_{\too}^\tast \norm{t^{-1-\gao}}^2_{L^6(\Si_t)}\norm{t^2\D^{(\OOI)}\pih}^2_{L^6(\Si_t)}\norm{t^4\D\R}^2_{L^6(\Si_t)}\,\d t\ri)^{1/2} \\
                                                                              & \les (D\varep)^2 \le(\int_{\too}^\tast t^{-1-2\gao}\,\d t\ri)^{1/2} \\
    & \les (D\varep)^2.
  \end{align*}
\end{remark}

Using~\eqref{est:L2DIVLIEHOOI} and $L^2(\MM^\intr_\bott)$ estimates obtained from the Bootstrap Assumptions~\ref{BA:fluxglobener} on the flux of the contracted Bel-Robinson tensors, we have
\begin{align*}
  & \le|\int_{\MM^\intr_\bott}\DIV \le(Q\le(\Lieh_\OOI\Lieh_\OOI\R\ri)\ri) \cdot\le(\KI,\KI,\TI\ri)\ri| \\
  \les & \; \int_{\MM^\intr_\bott} t^4\le|\Lieh_\OOI\Lieh_\OOI\R\ri|\le|\DIV\le(\Lieh_\OOI\Lieh_\OOI\R\ri)\ri|  \\
  \les & \; \norm{t^{-1/2-\ga}t^2\Lieh_\OOI\Lieh_\OOI\R}_{L^2(\MM^\intr_\bott)}\norm{t^{4-\gao} \DIV\le(\Lieh_\OOI\Lieh_\OOI\R\ri)}_{L^2(\MM^\intr_\bott)} \\
  \les & \; (D\varep)^3,
\end{align*}
provided that $\ga \leq 3/2-\gao$.\\

The estimates for the other error terms of $\EE^\intr_2$ are obtained either arguing as in Section~\ref{sec:errorintr1} or as above, using the estimates for $\D\TI, \D\SI$ from the Bootstrap Assumptions~\ref{BA:intKill}. This finishes the control of $\EE^\intr_2$.


\section{Estimates for the exterior error terms $\EE^\ext$}\label{sec:exterrest}
We denote by $\EE_{1,1}^\ext$, $\EE_{1,2}^\ext$, $\EE^\ext_{2,1}$ and $\EE^\ext_{2,2}$ the exterior error terms produced in the energy estimates, which are defined as follows
\begin{align*}
  \EE^\ext_{1,1} & := \int_{\MM^\ext} \big(Q(\Lieh_\TE\R)(\KE,\KE)\cdot^{(\KE)}\pih + Q(\Lieh_\OOE\R)(\KE,\KE)\cdot^{(\TE)}\pih \\
  & \quad\quad + Q(\Lieh_\OOE\R)(\KE,\TE)\cdot^{(\KE)}\pih\big),\\
  \EE_{1,2}^\ext & := \int_{\MM_{\ext}}\le(\DIV Q(\Lieh_\TE\R)(\KE,\KE,\KE) + \DIV Q(\Lieh_\OOE\R)(\KE,\KE,\TE)\ri),\\
  \EE^\ext_{2,1} & := \int_{\MM^\ext}\big(Q(\Lieh_{\SE}\Lieh_{\TE}\R)(\KE,\KE)\cdot^{(\KE)}\pih + Q(\Lieh^2_{\OOE}\R)(\KE,\TE)\cdot^{(\KE)}\pih \\
             & \quad\quad + Q(\Lieh_{\OOE}^2\R)(\KE,\KE)\cdot^{(\TE)}\pih + Q(\Lieh_{\OOE}\Lieh_{\TE}\R)(\KE,\KE)\cdot^{(\KE)}\pih, \\
  \EE^\ext_{2,2} & := \int_{\MM^\ext} \big(\DIV Q(\Lieh_{\SE}\Lieh_{\TE}\R)(\KE,\KE,\KE) + \DIV Q(\Lieh_{\OOE}^2\R)(\KE,\KE,\TE) \\
  & \quad\quad + \DIV Q(\Lieh_{\OOE}\Lieh_{\TE}\R)(\KE,\KE,\KE) \big).                  
\end{align*} 


\subsection{Preliminary definitions and computational results}\label{sec:exterrpreldef}
\paragraph{Null decompositions of $\pih$.}
We define the following null decompositions of the deformation tensor $\pih$
\begin{align*}
  ^{(X)}\ibf_{ab} & = ^{(X)}\pih_{ab}, & ^{(X)}\jbf & = ^{(X)}\pih_{34},\\
  ^{(X)}\mbf_a & = ^{(X)}\pih_{4a}, & ^{(X)}\mbbf_a & = ^{(X)}\pih_{3a},\\
  ^{(X)}\nbf & = ^{(X)}\pih_{44}, & ^{(X)}\nbbf & = ^{(X)}\pih_{33}.
\end{align*}

We have for the time translation $\TE=\half(\el+\elb)$
\begin{align*}
  \begin{aligned}
    \tr ^{(\TE)}\pi & = -2\omb + \trchi+\trchib,
  \end{aligned}
\end{align*}
and
\begin{align*}
  \begin{aligned}
  ^{(\TE)}\nbf & = 0, &  ^{(\TE)}\nbbf & = -4\omb, \\
  ^{(\TE)}\mbf_a & = -2\ze_a,&  ^{(\TE)}\mbbf_a & = \xib_a + 2\ze_a,
\end{aligned}
\end{align*}
and
\begin{align*}
  \begin{aligned}
    ^{(\TE)}\jbf & = \omb + \half(\trchi+\trchib), \\
    ^{(\TE)}\ibf_{ab} & = \half\omb\gd_{ab}+\quar(\trchi+\trchib)\gd_{ab} +\chih_{ab}+\chibh_{ab}.
  \end{aligned}
\end{align*} 

We have for the scaling vectorfield $\SE=\half(u\elb+\ub\el)$
\begin{align*}
  \tr ^{(\SE)}\pi & = 4 -2u\omb + \ub\trchi+u\trchib,
\end{align*}
and
\begin{align*}
  \begin{aligned}
    ^{(\SE)}\nbf & = 0, & ^{(\SE)}\nbbf & = -2\yy -4\ub\omb,\\
    ^{(\SE)}\mbf_a & = -2u\ze_a, &  ^{(\SE)}\mbbf_a & = u\xib_a+2\ub\ze_a,
  \end{aligned}
\end{align*}
and
\begin{align*}
  \begin{aligned}
    ^{(\SE)}\jbf & = -2 + u\omb +\half(\ub\trchi+u\trchib), \\
    ^{(\SE)}\ibf_{ab} & = -\gd_{ab} + \half u\omb\gd_{ab} + \quar(\ub\trchi +u\trchib)\gd_{ab} +\ub\chih_{ab}+u\chibh_{ab}.
  \end{aligned}
\end{align*}

We have for the conformal Morawetz vectorfield $\KE=\half(u^2\elb+\ub^2\el)$
\begin{align*}
  \tr ^{(\KE)}\pi & = 4(u+\ub) -2u^2\omb + \ub^2\trchi+u^2\trchib,
\end{align*}
and
\begin{align*}
  \begin{aligned}
    ^{(\KE)}\nbf & = 0, & ^{(\KE)}\nbbf & = -4\ub\yy-4\ub^2\omb,\\
    ^{(\KE)}\mbf_a & = -2u^2\ze_a, &  ^{(\KE)}\mbbf_a & = u^2\xib_a + 2\ub^2\ze_a,
  \end{aligned}
\end{align*}
and
\begin{align*}
  \begin{aligned}
    ^{(\KE)}\jbf & = -2(u+\ub) + u^2\omb +\half(\ub^2\trchi+u^2\trchib),\\
    ^{(\KE)}\ibf_{ab} & = -(u+\ub)\gd_{ab} + \half u^2\omb\gd_{ab} + \quar(\ub^2\trchi +u^2\trchib)\gd_{ab} +\ub^2\chih_{ab}+u^2\chibh_{ab}.
  \end{aligned}
\end{align*}

We have for the rotations $\OOE$
\begin{align*}
  \tr ^{(\OOE)}\pi & = 2\ze\cdot\OOE + \tr \HEi,
\end{align*}
and
\begin{align*}
  \begin{aligned}
    ^{(\OOE)}\nbf & = 0, & ^{(\OOE)}\nbbf & = -4\xib\cdot\OOE,\\
    ^{(\OOE)}\mbf_a & = 0, & ^{(\OOE)}\mbbf_a & = \YEi_a,
  \end{aligned}
\end{align*}
and
\begin{align*}
  \begin{aligned}
    ^{(\OOE)}\jbf & = -\ze\cdot\OOE + \half\tr \HEi,\\                                     
    ^{(\OOE)}\ibf_{ab} & = \HEi_{ab} -\quar\le(2\ze\cdot\OOE+\tr \HEi\ri)\gd_{ab}.
  \end{aligned}
\end{align*}

\paragraph{Null decomposition of $\D\pih$}
We have the following definitions for (contractions of) $\D\pih$
\begin{align*}
  ^{(X)}p_\mu &:= \DIV ^{(X)}\pih_\mu,\\
  ^{(X)}q_{\mu\nu\la} & := \D_\nu^{(X)}\pih_{\la\mu}-\D_\la^{(X)}\pih_{\nu\mu} - \frac{1}{3}\le(^{(X)}p_\la\g_{\mu\nu}-^{(X)}p_\nu\g_{\mu\la}\ri).
\end{align*}
We note $p_3,p_4$ and $\psl$ the null decomposition of the spacetime vectorfield $p$.\\

For a spacetime $3$-tensor $F$, we define its null decomposition\footnote{Note that this decomposition is only relevant for $3$-tensors with the same symmetries as $q$ or the currents $J$ defined below.} by (see~\cite[p. 212]{Chr.Kla93})
\begin{align*}
  \Lambda(F) & := \quar F_{434}, & \Lambdab(F) & := \quar F_{343},\\
  K(F) & := \quar \in^{ab}F_{4ab}, & \Kbb(F) & := \quar \in^{ab}F_{3ab},\\
  \Xi(F)_{a} & := \half F_{44a}, & \Xib(F)_a & := \half F_{33a},\\
  I(F)_a & := \half F_{34a}, & \Ib(F)_a & := \half F_{43a},\\
  \Theta(F)_{ab} & := F_{a4b} + F_{b4a}-(\gd^{ab}F_{a4b})\gd_{ab}, & \Thetab(F)_{ab} & := F_{a3b}+F_{b3a}-(\gd^{ab}F_{a3b})\gd_{ab}.
\end{align*}

\paragraph{The currents $J$}
We define the \emph{current} $J$ of a Weyl tensor $W$ to be
\begin{align*}
  J(W)_{\nu\la\ga} & := \D^\mu W_{\mu\nu\la\ga}.
\end{align*}


For a vectorfield $X$, and a Weyl tensor $W$, we note $J(X,W)$ the current of $\Lieh_XW$. We have the following formula (see~\cite[p. 205]{Chr.Kla93})
\begin{align}\label{eq:JdecompoJ1J2J3}
  J(X,W) & = \half \le(J^1(X,W) + J^2(X,W) + J^3(X,W)\ri),
\end{align}
where
\begin{align*}
  J^1(X,W)_{\be\ga\de} & := {^{(X)}\pih}^{\mu\nu}\D_\nu W_{\mu\be\ga\de},\\
  J^{2}(X,W)_{\be\ga\de} & := {^{(X)}p}_\la W^{\la}_{\,\be\ga\de},\\
  J^3(X,W)_{\be\ga\de} & := {^{(X)}q}_{\al\be\la}W^{\al\la}_{\,\,\ga\de} + {^{(X)}q}_{\al\ga\la}W^{\al\, \la}_{\,\be\,\de} + {^{(X)}q}_{\al\de\la}W^{\al\,\,\la}_{\,\be\ga}.
\end{align*}

The divergence of the Bel-Robinson tensors is related to the null decomposition of the associated Weyl field and the null decomposition of its current. When the Weyl field is a modified Lie derivative $\Lieh_XW$, its current $J(X,W)$ can be expressed in terms of the null decomposition of $^{(X)}\pih$, the null decomposition of $\D^{(X)}\pih$ and null decompositions of $W$ and $\D W$. See~\cite[pp. 214--218]{Chr.Kla93} for the computations. We use these formulas in the following Sections~\ref{sec:esterrEE1}--\ref{sec:prooflemdecayesterrtermOOESEW}.    

\subsection{Preliminary sup-norm estimates for the deformation tensors $\pi$}
We define the following decay norms for the null decompositions of the deformation tensors of $\TE,\SE,\KE$ and $\OOE$.
\begin{align*}
  & \DDf^\ext_0[\TE] := \norm{u^{3/2}\ub^{(\TE)}\ibf}_{L^\infty} + \norm{u^{3/2}\ub^{(\TE)}\jbf}_{L^\infty} + \norm{\ub^2u^{1/2} {^{(\TE)}\mbf}}_{L^\infty} + \norm{\ub u^{3/2}{^{(\TE)}\mbbf}}_{L^\infty} \\
  & \quad \quad \quad \quad + \norm{\ub u^{3/2}{^{(\TE)}\nbbf}}_{L^\infty}\\ \\
  & \DDf^\ext_0[\SE] := \norm{\ub u^{1/2}{^{(\SE)}\ibf}}_{L^\infty} + \norm{\ub u^{1/2}{^{(\SE)}\jbf}}_{L^\infty} + \norm{\ub^{2}u^{1/2}{^{(\SE)}\mbf}}_{L^\infty} + \norm{\ub u^{1/2}{^{(\SE)}\mbbf}}_{L^\infty}  \\
  & \quad \quad \quad \quad + \norm{u^{3/2}{^{(\SE)}\nbbf}}_{L^\infty},\\ \\
  & \DDf^\ext_0[\KE] := \norm{u^{1/2}{^{(\KE)}\ibf}}_{L^\infty} + \norm{u^{1/2}{^{(\KE)}\jbf}}_{L^\infty} + \norm{\ub^{2}u^{-3/2}{^{(\KE)}\mbf}}_{L^\infty} + \norm{u^{1/2}{^{(\KE)}\mbbf}}_{L^\infty}  \\
  & \quad \quad \quad \quad + \norm{u^{3/2}\ub^{-1}{^{(\KE)}\nbbf}}_{L^\infty},\\
  & \DDf^\ext_0[\OOE] := \norm{\ub u^{1/2}{^{(\OOE)}\ibf}}_{L^\infty} + \norm{\ub u^{1/2}{^{(\OOE)}\jbf}}_{L^\infty} + \norm{\ub u^{1/2}{^{(\OOE)}\mbbf}}_{L^\infty} + \norm{u^{3/2}{^{(\OOE)}\nbbf}}_{L^\infty},
\end{align*}
together with
\begin{align*}
  ^{(\TE)}\nbf & = ^{(\SE)}\nbf = ^{(\KE)}\nbf = ^{(\OOE)}\nbf = 0, \\
  ^{(\OOE)}\mbf & = 0.
\end{align*}

From the Bootstrap Assumptions~\ref{BA:connext} and the formulas from Section~\ref{sec:exterrpreldef} relating the null decompositions of the deformation tensors of $\TE,\SE,\KE$ and $\OOE$, the null connection coefficients and the rotation coefficients $\YEi$ and $\HEi$, we have
\begin{align}\label{est:DDfext0}
  \DDf^\ext_0[\TE,\SE,\KE,\OOE] & \les D\varep.
\end{align}

\begin{remark}
  We detail the estimate for $^{(\SE)}\jbf$. We have
  \begin{align*}
    ^{(\SE)}\jbf & = -2+\half(\ub\trchi+u\trchib) +u\omb \\
    & = -2 + \frac{\ub-u}{r} +\half\ub\le(\trchi-\frac{2}{r}\ri) + \half u\le(\trchib+\frac{2}{r}\ri) +u\omb \\
     & = \frac{\ub-u-2r}{r} +\half\ub\le(\trchi-\frac{2}{r}\ri) + \half u\le(\trchib+\frac{2}{r}\ri) +u\omb,
  \end{align*}
  and therefore
  \begin{align*}
    \norm{\ub u^{1/2} {^{(\SE)}\jbf}}_{L^\infty(\MM^\ext)} & \les \norm{u^{1/2}\le(r-\half(\ub-u)\ri)}_{L^\infty(\MM^\ext)} + \mathfrak{O}^\ext_{\leq 1} \les D\varep, 
  \end{align*}
  where we used the bootstrap bound~\eqref{est:BAarearadiusestimate} for $r$ and where we refer to Section~\ref{sec:normsnullconn} for the definition of $\mathfrak{O}^\ext_{\leq 1}$.
\end{remark}

\begin{remark}
  The decay rates~\eqref{est:DDfext0} are different from the ones obtained in~\cite[p. 223]{Chr.Kla93}. They are stronger to the ones obtained in~\cite[pp. 251-254]{Kla.Nic03} -- in particular due to the stronger bounds~\eqref{est:BAarearadiusestimate} for the area radius $r$ that we derived in this paper--, except for $^{(\TE)}\mbbf$ which lose $\ub u^{-1}$ and for $^{(\OOE)}\mbbf$ which is non zero in our case. Both these differences are due to the presence of $\xib$ which is vanishing in the double null foliation of~\cite{Kla.Nic03}.
  In the next sections we prove in particular that the control of the error terms is still valid despite these differences.  
\end{remark}

\subsection{Estimates for $\EE^\ext_{1,1}, \EE^\ext_{2,1}$}\label{sec:esterrEE1}
The estimates for the error terms $\EE^\ext_{1,1}$ and $\EE^\ext_{2,1}$ are obtained exactly in the same manner and each split into obtaining a control of the spacetime integrals in $\MM^\ext$ of the following three integrands:
\begin{align}\label{eq:integrandsEE1}
  Q(W)(\KE,\KE)\cdot^{(\KE)}\pih, && Q(W)(\KE,\KE)\cdot^{(\TE)}\pih, && Q(W)(\KE,\TE)\cdot^{(\KE)}\pih,
\end{align}
where $W$ is a Weyl tensor such that energy estimates are performed for the following contracted Bel-Robinson tensors respectively
\begin{align*}
  Q(W)(\KE,\KE,\KE), && Q(W)(\KE,\KE,\TE), && Q(W)(\KE,\KE,\TE).
\end{align*}

Since we have same or stronger decay rates for the sup-norm of $^{(\KE)}\pih$ as in~\cite{Kla.Nic03}, we only treat the second term of~\eqref{eq:integrandsEE1} which involves $^{(\TE)}\pih$, \emph{i.e.}
\begin{align*}
  \int_{\MM^\ext} Q(W)(\KE,\KE)\cdot ^{(\TE)}\pih.
\end{align*}

Arguing as in Section~\ref{sec:errorintr1}, one can obtain from the Bootstrap Assumptions~\ref{BA:fluxglobener} on the boundedness of the energy fluxes through the hypersurfaces $\CC_u$ and $\Si_t^\ext$ the following $L^2(\MM^\ext)$ bounds for the null decomposition of $W$
\begin{align}\label{est:L2MMextW}
  \begin{aligned}
    \norm{u^{-1/2-\ga}\ub^2\al(W)}_{L^2(\MM^\ext)} & \les_\ga D\varep, & \norm{u^{-1/2-\ga}\ub^2\be(W)}_{L^2(\MM^\ext)} & \les_\ga D\varep,\\
    \norm{\ub^{-1/2-\ga}\ub^{2}\rho(W)}_{L^2(\MM^\ext)} & \les_\ga D\varep, & \norm{\ub^{-1/2-\ga}\ub^2\sigma(W)}_{L^2(\MM^\ext)} & \les_\ga D\varep,\\
    \norm{\ub^{-1/2-\ga}\ub u\beb(W)}_{L^2(\MM^\ext)} & \les_\ga D\varep, & \norm{\ub^{-1/2-\ga}u^2\alb(W)}_{L^2(\MM^\ext)} & \les_\ga D\varep,
  \end{aligned}
\end{align}
for all $\ga>0$.\\

We now decompose $Q(W)(\KE,\KE)\cdot{^{(\TE)}\pih}$ in terms of the null decompositions of $W$ and $^{(\TE)}\pih$. We have
\begin{align*}
  Q(W)(\KE,\KE)\cdot{^{(\TE)}\pih} & = \ub^4 Q(W)_{44}^{\mu\nu}{^{(\TE)}\pih}_{\mu\nu} + \ub^2u^2Q(W)_{34}^{\mu\nu}{^{(\TE)}\pih}_{\mu\nu} + u^4Q(W)_{33}^{\mu\nu}{^{(\TE)}\pih}_{\mu\nu}.
\end{align*}
Rewriting schematically the formulas (6.2.45) (6.2.46) and (6.2.47) of~\cite[p. 275]{Kla.Nic03} (see also~\cite[pp. 248--250]{Chr.Kla93}), we have
\begin{align*}
  \begin{split}
  & \ub^4 Q(W)_{44}^{\mu\nu}{^{(\TE)}\pih}_{\mu\nu} \\
  = & \; \ub^4\bigg(|\al|^2\nbbf + (\rho^2 +\sigma^2)\nbf + |\be|^2\jbf + \al\be\mbbf + \rho\be\mbf +\sigma\be\mbf + |\be|^2\ibf + \rho\al\ibf + \sigma\al\ibf\bigg),
\end{split}\\ \\
  \begin{split}
    & \ub^2u^2 Q(W)_{34}^{\mu\nu}{^{(\TE)}\pih}_{\mu\nu} \\
    = & \;\ub^2u^2\bigg(|\be|^2\nbbf + |\beb|^2\nbf + (\rho^2+\sigma^2)\jbf + \rho\be\mbbf + \sigma\be\mbbf + \rho\beb\mbf + \sigma\beb\mbf +(\rho^2+\sigma^2)\ibf + \be\beb\ibf\bigg),
  \end{split}\\ \\
  \begin{split}
    & u^4Q(W)_{33}^{\mu\nu}{^{(\TE)}\pih}_{\mu\nu} \\
    = & \; u^4\bigg((\rho^2+\sigma^2)\nbbf + |\alb|^2\nbf + |\beb|^2\jbf + \alb\beb\mbf + \rho\beb\mbbf + \sigma\beb\mbbf + |\beb|^2\ibf + \rho\alb\ibf + \sigma\alb\ibf\bigg),
  \end{split}
\end{align*}
where $(\al,\be,\rho,\sigma,\beb,\alb) = \nulld(W)$ and $(\ibf, \jbf, \mbf, \mbbf, \nbf, \nbbf) = \le(^{(\TE)}\ibf, ^{(\TE)}\jbf, ^{(\TE)}\mbf, ^{(\TE)}\mbbf, ^{(\TE)}\nbf, ^{(\TE)}\nbbf\ri)$.\\

Using the sup norm estimates~\eqref{est:DDfext0} for the null decomposition $(\ibf, \jbf, \mbf, \mbbf, \nbf, \nbbf)$ of the deformation tensor $^{(\TE)}\pih$, we obtain that 
\begin{align}\label{est:ErrQWTE}
  \begin{split}
  & \le|\ub^{4}Q(W)_{44}^{\mu\nu} {^{(\TE)}\pih_{\mu\nu}}\ri| \\
  \les \; & D\varep \bigg(\ub^{3}u^{-3/2}|\al|^2  + \ub^{3}u^{-3/2}|\be|^2 + \ub^{3}u^{-3/2}|\al||\be| \\
  & + \ub^{2}u^{-1/2}|\rho||\be| +\ub^{2}u^{-1/2}|\sigma||\be| + \ub^{3}u^{-3/2}|\be|^2 + \ub^3u^{-3/2}|\rho||\al|\\
  & + \ub^3u^{-3/2}|\sigma||\al|\bigg),
\end{split} \\ \\
  \begin{split}
    & \le|\ub^2u^2 Q(W)_{34}^{\mu\nu}{^{(\TE)}\pih}_{\mu\nu}\ri| \\
    \les \; & D\varep \bigg( \ub u^{1/2}|\be|^2 + \ub u^{1/2}(\rho^2+\sigma^2) + \ub u^{1/2}|\rho||\be| \\
    & + \ub u^{1/2}|\sigma||\be| + u^{3/2}|\rho||\beb| + u^{3/2}|\sigma||\beb| +\ub u^{1/2}(\rho^2+\sigma^2)\\
    & + \ub u^{1/2}|\be||\beb|\bigg),
  \end{split}\\ \\
  \begin{split}
    & \le|u^4Q(W)_{33}^{\mu\nu}{^{(\TE)}\pih}_{\mu\nu}\ri| \\
    \les & \; D\varep\bigg(\ub^{-1}u^{5/2} (\rho^2+\sigma^2) + \ub^{-1}u^{5/2}|\beb| + \ub^{-2}u^{7/2}|\alb||\beb| \\
    & + \ub^{-1}u^{5/2}|\rho||\beb| + \ub^{-1}u^{5/2}|\sigma||\beb| + \ub^{-1}u^{5/2}|\beb|^2 + \ub^{-1}u^{5/2}|\rho||\alb|\\
    & + \ub^{-1}u^{5/2}|\sigma||\alb|\bigg).
  \end{split}
\end{align}

Applying Cauchy-Schwartz, one checks that the spacetime integral in $\MM^\ext$ of all the terms from~\eqref{est:ErrQWTE} can be controlled by the $L^2(\MM^\ext)$ norms of~(\ref{est:L2MMextW}), which thus gives
\begin{align*}
  \le|\int_{\MM^\ext}Q(W)(\KE,\KE)\cdot{^{(\TE)}\pih}\ri|  & \les (D\varep)^3. 
\end{align*}

This concludes the estimate of the error terms $\EE^\ext_{1,1}$ and $\EE^\ext_{2,1}$.




\subsection{Preliminary $L^\infty L^4(S)$ estimates for $\D\pi$}
We have the following definition of the decay norms for the null decompositions of the tensors $p$ and $q$
\begin{align*}
  \DDf^\ext_{1,p}[\TE] & := \norm{\ub^{1/2}u^{5/2}{^{(\TE)}p_3}}_{L^\infty_{u,\ub}L^4(S_{u,\ub})} + \norm{\ub^{3/2}u^{3/2}{^{(\TE)}p_4}}_{L^\infty_{u,\ub}L^4(S_{u,\ub})} \\
                       & \quad + \norm{\ub^{3/2}u^{3/2}{^{(\TE)}\psl}}_{L^\infty_{u,\ub}L^4(S_{u,\ub})}, \\
  \DDf^\ext_{1,q}[\TE] & := \norm{\ub^{3/2}u^{3/2}\Lambda(^{(\TE)}q)}_{L^\infty_{u,\ub}L^4(S_{u,\ub})} + \norm{\ub^{5/2}u^{1/2}K(^{(\TE)}q)}_{L^\infty_{u,\ub}L^4(S_{u,\ub})} \\
                       & \quad + \norm{\ub^{5/2}u^{1/2}\Xi(^{(\TE)}q)}_{L^\infty_{u,\ub}L^4(S_{u,\ub})} + \norm{\ub^{3/2}u^{3/2}I(^{(\TE)}q)}_{L^\infty_{u,\ub}L^4(S_{u,\ub})} \\
                       & \quad + \norm{\ub^{3/2}u^{3/2}\Theta(^{(\TE)}q)}_{L^\infty_{u,\ub}L^4(S_{u,\ub})} + \norm{\ub^{1/2}u^{5/2}\Lambdab(^{(\TE)}q)}_{L^\infty_{u,\ub}L^4(S_{u,\ub})} \\
                       & \quad + \norm{\ub^{3/2}u^{3/2}\Kbb(^{(\TE)}q)}_{L^\infty_{u,\ub}L^4(S_{u,\ub})} + \norm{\ub^{1/2}u^{5/2}\Xib(^{(\TE)}q)}_{L^\infty_{u,\ub}L^4(S_{u,\ub})} \\
                       & \quad + \norm{\ub^{3/2}u^{3/2}\Ib(^{(\TE)}q)}_{L^\infty_{u,\ub}L^4(S_{u,\ub})} + \norm{\ub^{1/2}u^{5/2}\Thetab(^{(\TE)}q)}_{L^\infty_{u,\ub}L^4(S_{u,\ub})}, 
\end{align*}
and
\begin{align*}
  \DDf^\ext_{1,p}[\SE] & := \norm{\ub^{1/2}u^{3/2}{^{(\SE)}p_3}}_{L^\infty_{u,\ub}L^4(S_{u,\ub})} + \norm{\ub^{3/2}u^{1/2}{^{(\SE)}p_4}}_{L^\infty_{u,\ub}L^4(S_{u,\ub})} \\
                       & \quad + \norm{\ub^{3/2}u^{1/2} {^{(\SE)}\psl}}_{L^\infty_{u,\ub}L^4(S_{u,\ub})}, \\
  \DDf^\ext_{1,q}[\SE] & := \norm{\ub^{3/2}u^{1/2}\Lambda(^{(\SE)}q)}_{L^\infty_{u,\ub}L^4(S_{u,\ub})} + \norm{\ub^{5/2}u^{-1/2}K(^{(\SE)}q)}_{L^\infty_{u,\ub}L^4(S_{u,\ub})} \\
                       & \quad + \norm{\ub^{5/2}u^{-1/2}\Xi(^{(\SE)}q)}_{L^\infty_{u,\ub}L^4(S_{u,\ub})} + \norm{\ub^{3/2}u^{1/2} I(^{(\SE)}q)}_{L^\infty_{u,\ub}L^4(S_{u,\ub})} \\
                       & \quad+ \norm{\ub^{3/2}u^{1/2}\Theta(^{(\SE)}q)}_{L^\infty_{u,\ub}L^4(S_{u,\ub})} + \norm{\ub^{1/2}u^{3/2}\Lambdab(^{(\SE)}q)}_{L^\infty_{u,\ub}L^4(S_{u,\ub})} \\
                       & \quad + \norm{\ub^{3/2}u^{1/2}\Kbb(^{(\SE)}q)}_{L^\infty_{u,\ub}L^4(S_{u,\ub})} + \norm{\ub^{1/2}u^{3/2}\Xib(^{(\SE)}q)}_{L^\infty_{u,\ub}L^4(S_{u,\ub})} \\
                       & \quad + \norm{\ub^{3/2}u^{1/2}\Ib(^{(\SE)}q)}_{L^\infty_{u,\ub}L^4(S_{u,\ub})} + \norm{\ub^{1/2}u^{3/2}\Thetab(^{(\SE)}q)}_{L^\infty_{u,\ub}L^4(S_{u,\ub})}, 
\end{align*}
and
\begin{align*}
  \DDf^\ext_{1,p}[\OOE] & := \norm{\ub^{1/2}u^{3/2}{^{(\OOE)}p_3}}_{L^\infty_{u,\ub}L^4(S_{u,\ub})} + \norm{\ub^{3/2}u^{1/2}{^{(\OOE)}p_4}}_{L^\infty_{u,\ub}L^4(S_{u,\ub})},\\
  & \quad + \norm{\ub^{3/2}u^{1/2}{^{(\OOE)}\psl}}_{L^\infty_{u,\ub}L^4(S_{u,\ub})}, \\
  \DDf^\ext_{1,q}[\OOE] & := \norm{\ub^{3/2}u^{1/2}\Lambda(^{(\OOE)}q)}_{L^\infty_{u,\ub}L^4(S_{u,\ub})} + \norm{\ub^{5/2}u K(^{(\OOE)}q)}_{L^\infty_{u,\ub}L^4(S_{u,\ub})} \\
                       & \quad + \norm{\ub^{3/2}u^{1/2} I(^{(\OOE)}q)}_{L^\infty_{u,\ub}L^4(S_{u,\ub})} \\
                       & \quad+ \norm{\ub^{3/2}u^{1/2}\Theta(^{(\OOE)}q)}_{L^\infty_{u,\ub}L^4(S_{u,\ub})} + \norm{\ub^{1/2}u^{3/2}\Lambdab(^{(\OOE)}q)}_{L^\infty_{u,\ub}L^4(S_{u,\ub})} \\
                       & \quad + \norm{\ub^{3/2}u^{1/2}\Kbb(^{(\OOE)}q)}_{L^\infty_{u,\ub}L^4(S_{u,\ub})} + \norm{\ub^{1/2}u^{3/2}\Xib(^{(\OOE)}q)}_{L^\infty_{u,\ub}L^4(S_{u,\ub})} \\
                       & \quad + \norm{\ub^{3/2}u^{1/2}\Ib(^{(\OOE)}q)}_{L^\infty_{u,\ub}L^4(S_{u,\ub})} + \norm{\ub^{1/2}u^{3/2}\Thetab(^{(\OOE)}q)}_{L^\infty_{u,\ub}L^4(S_{u,\ub})},
\end{align*}
together with $\Xi(^{(\OOE)}q) = 0$.\\


Using formulas from~\cite[pp. 231--232]{Chr.Kla93}, the Bootstrap Assumptions~\ref{BA:connext} for $\HHt$ norms of (one derivative of) the null connection coefficients in $\MM^\ext$, the product estimate from Lemma~\ref{lem:prodH12}, the Sobolev embeddings from Lemma~\ref{lem:Sobsphere}, we have
\begin{align}\label{est:DDf1p1q}
  \DDf^\ext_{1,p} + \DDf^\ext_{1,q} & \les \mathfrak{O}^\ext_{\leq 1} \les D\varep.
\end{align}

\begin{remark}
  The decay rate of each term is easily checked using that $r\Nd, \ub\Nd_4$ and $u\Nd_3$ preserve the decay rates, and that the decay of the $L^4(S)$ norm loses a weight $r^{1/2}$ with respect to the decay of the $L^\infty(S)$ norm.  
\end{remark}


\begin{remark}
  The decay rates in the norms $\DDf^\ext_{1}[\TE]$, $\DDf^\ext_{1}[\SE]$ and $\DDf^\ext_1[\OOE]$ are the same or better as the ones obtained in~\cite[pp. 250-258]{Kla.Nic03}\footnote{By a scaling consideration, we do believe that there is a typo in the estimates (6.1.53) which should be multiplied by $r^{-1}$. Moreover, we do also believe that there is a typo in the estimate for $\Xib(^{(\SE)}q)$ in~\cite{Kla.Nic03} which should be multiplied by $r^{-1}$.} except for $\Xib(^{\TE}q)$, which loses a $\ub u^{-1}$ factor due to the non vanishing of $\xib$ in the present paper.
\end{remark}

\subsection{Estimates for $\EE^\ext_{1,2}$}\label{sec:esterrEE12}
We start with the estimate for
\begin{align*}
  \int_{\MM^\ext}\DIV Q(\Lieh_\TE\R)(\KE,\KE,\KE).
\end{align*}

Arguing as in Section~\ref{sec:errorintr1} using the Bootstrap Assumption~\ref{BA:fluxglobener}, we first record the following $L^2(\MM^\ext)$ bounds for the null decomposition of $\Lieh_\TE\R$
\begin{align}\label{est:L2MMextW1}
  \begin{aligned}
    \norm{u^{-1/2-\ga}\ub^3\al\le(\Lieh_\TE\R\ri)}_{L^2(\MM^\ext)} & \les_\ga D\varep, & \norm{\ub^{-1/2-\ga}\ub^3\be\le(\Lieh_\TE\R\ri)}_{L^2(\MM^\ext)} & \les_\ga D\varep,\\
    \norm{\ub^{-1/2-\ga}\ub^{2}u\rho\le(\Lieh_\TE\R\ri)}_{L^2(\MM^\ext)} & \les_\ga D\varep, & \norm{\ub^{-1/2-\ga}\ub^2u\sigma\le(\Lieh_\TE\R\ri)}_{L^2(\MM^\ext)} & \les_\ga D\varep,\\
    \norm{\ub^{-1/2-\ga}\ub u^2\beb\le(\Lieh_\TE\R\ri)}_{L^2(\MM^\ext)} & \les_\ga D\varep, & \norm{\ub^{-1/2-\ga}u^3\alb\le(\Lieh_\TE\R\ri)}_{L^2(\MM^\ext)} & \les_\ga D\varep,
  \end{aligned}
\end{align}
for all $\ga>0$.\\

We have
\begin{align}\label{eq:DIVQTR}
  \begin{aligned}
    \DIV Q(\Lieh_\TE\R)(\KE,\KE,\KE) & = \ub^6\DIV Q(\Lieh_\TE\R)_{444} + \ub^4u^2\DIV Q(\Lieh_\TE\R)_{443} \\
    & \quad + \ub^2u^4 \DIV Q(\Lieh_\TE\R)_{433} + u^6\DIV Q(\Lieh_\TE\R)_{333}.
  \end{aligned}
\end{align}

From~\cite[p. 213]{Chr.Kla93}, we have
\begin{align}\label{eq:nullDIVLTR}
  \begin{aligned}
    \ub^6\DIV Q(\Lieh_\TE\R)_{444} & = \ub^6\al\le(\Lieh_\TE\R\ri)\Theta(J(\Lieh_\TE\R)) \\
    & \quad + \ub^6\be\le(\Lieh_\TE\R\ri)\Xi(J(\Lieh_\TE\R)),\\
    \ub^4u^2\DIV Q(\Lieh_\TE\R)_{443} & = \ub^4u^2\rho\le(\Lieh_\TE\R\ri)\Lambda(J(\Lieh_\TE\R)) \\
    &\quad + \ub^4u^2 \sigma\le(\Lieh_\TE\R\ri)K(J(\Lieh_\TE\R)) \\
    & \quad + \ub^4u^2\be\le(\Lieh_\TE\R\ri)I(J(\Lieh_\TE\R)),\\
    \ub^2u^4\DIV Q(\Lieh_\TE\R)_{433} & = \ub^2u^4\rho\le(\Lieh_\TE\R\ri)\Lambdab(J(\Lieh_\TE\R)) \\
    & \quad + \ub^2u^4 \sigma\le(\Lieh_\TE\R\ri)\Kb(J(\Lieh_\TE\R)) \\
    & \quad + \ub^2u^4\be\le(\Lieh_\TE\R\ri)\Ib(J(\Lieh_\TE\R)),\\
    u^6\DIV Q(\Lieh_\TE\R)_{333} & = u^6\alb\le(\Lieh_\TE\R\ri)\Thetab(J(\Lieh_\TE\R)) \\
    & \quad + u^6\beb\le(\Lieh_\TE\R\ri)\Xib(J(\Lieh_\TE\R)),
  \end{aligned}
\end{align}

We only examine the term
\begin{align}\label{eq:intbexi}
  \int_{\MM^\ext}\ub^6\be(\Lieh_\TE\R)\cdot\Xi(J).
\end{align}
of~\eqref{eq:nullDIVLTR} which is the hardest to treat due to the pairing of high weights in $\ub$ with slower decaying null components of $\Lieh_\TE\R$. Our goal is to obtain the following estimate
\begin{align*}
  \le|\int_{\MM^\ext}\ub^6\be(\Lieh_\TE\R)\cdot\Xi(J) \ri| & \les (D\varep)^3.
\end{align*}

Using the $L^2(\MM^\ext)$ estimates~(\ref{est:L2MMextW1}), the last estimate is proved provided that we obtain the following $L^2(\MM^\ext)$ estimates for $\Xi(J)$
\begin{align}\label{est:wishL2MMextXiJ}
  \norm{\ub^{3+1/2+\ga}\Xi(J)}_{L^2(\MM^\ext)} & \les (D\varep)^2.
\end{align}

Using the formulas from~\cite[p. 215]{Chr.Kla93} -- where we recall the definition~\eqref{eq:JdecompoJ1J2J3} of the decomposition of $J$ into $J^1,J^2$ and $J^3$ --, the estimates~(\ref{est:DDfext0}) for the sup-norm and the estimates~(\ref{est:DDf1p1q}) for the $L^\infty L^4$ norm of the null decompositions of $\pi$ and $\D\pi$ respectively, we have
\begin{align}\label{est:L2SXiJ1}
  \begin{aligned}
    \norm{\Xi(J^1(\Lieh_\TE\R))}_{L^2(S_{u,\ub})} & \les (D\varep)\bigg(\ub^{-1} u^{-3/2} \norm{\Nd\al(\R)}_{L^2(S_{u,\ub})} \\
    & \quad + \ub^{-2} u^{-1/2} \norm{\Nd_3\al(\R)}_{L^2(S_{u,\ub})} \\
    & \quad + \ub^{-1} u^{-3/2} \norm{\Nd_4\al(\R)}_{L^2(S_{u,\ub})} \\
    & \quad + \ub^{-2} u^{-1/2}  \norm{\Nd\be(\R)}_{L^2(S_{u,\ub})} \\
    & \quad + \ub^{-1} u^{-3/2} \norm{\Nd_4\be(\R)}_{L^2(S_{u,\ub})} \\
    & \quad + \ub^{-2}u^{-3/2}\norm{\al(\R)}_{L^2(S_{u,\ub})} \\
    & \quad + \ub^{-2}u^{-3/2}\norm{\be(\R)}_{L^2(S_{u,\ub})} \\
    & \quad + \ub^{-3}u^{-1/2}\norm{(\rho(\R),\sigma(\R))}_{L^2(S_{u,\ub})}\bigg),
  \end{aligned}
\end{align}
and
\begin{align}\label{est:L2SXiJ2}
  \begin{aligned}
    \norm{\Xi(J^2(\Lieh_\TE\R))}_{L^2(S_{u,\ub})} & \les (D\varep)\bigg(\ub^{-3/2}u^{-3/2} \norm{\al(\R)}_{L^4(S_{u,\ub})} \\
    & \quad + \ub^{-3/2} u^{-3/2} \norm{\be(\R)}_{L^4(S_{u,\ub})}\bigg), \\
    \norm{\Xi(J^3(\Lieh_\TE\R))}_{L^2(S_{u,\ub})}  & \les (D\varep)\bigg(\ub^{-3/2} u^{-3/2} \norm{\al(\R)}_{L^4(S_{u,\ub})} \\
    & \quad + \ub^{-3/2} u^{-3/2} \norm{\be(\R)}_{L^4(S_{u,\ub})} \\
    & \quad + \ub^{-5/2} u^{-1/2} \norm{(\rho(\R),\sigma(\R))}_{L^4(S_{u,\ub})}\bigg).
  \end{aligned}
\end{align}

Using the Sobolev estimates on the spheres $S_{u,\ub}$ from Lemma~\ref{lem:Sobsphere}, we further deduce from~\eqref{est:L2SXiJ2}
\begin{align}\label{est:L2SXiJ2bis}
  \begin{aligned}
    \norm{\Xi(J^2(\Lieh_\TE\R))}_{L^2(S_{u,\ub})} & \les (D\varep)\bigg(\ub^{-2}u^{-3/2} \norm{(r\Nd)^{\leq 1}\al(\R)}_{L^2(S_{u,\ub})} \\
    & \quad + \ub^{-2} u^{-3/2} \norm{(r\Nd)^{\leq 1}\be(\R)}_{L^2(S_{u,\ub})}\bigg), \\
    \norm{\Xi(J^3(\Lieh_\TE\R))}_{L^2(S_{u,\ub})}  & \les (D\varep)\bigg(\ub^{-2} u^{-3/2} \norm{(r\Nd)^{\leq 1}\al(\R)}_{L^2(S_{u,\ub})} \\
    & \quad + \ub^{-2} u^{-3/2} \norm{(r\Nd)^{\leq 1}\be(\R)}_{L^2(S_{u,\ub})} \\
    & \quad + \ub^{-3} u^{-1/2} \norm{(r\Nd)^{\leq 1}(\rho(\R),\sigma(\R))}_{L^2(S_{u,\ub})}\bigg).
  \end{aligned}
\end{align}

From~\eqref{est:L2SXiJ1} and~\eqref{est:L2SXiJ2bis}, we obtain
\begin{align}\label{est:XiJL2MMext}
  \begin{aligned}
    \norm{\ub^{3+1/2+\ga}\Xi\le(J\le(\Lieh_\TE\R\ri)\ri)}_{L^2(\MM^\ext)} & \les (D\varep)\bigg(\norm{\ub^{1+1/2+\ga} u^{-3/2} (r\Nd)^{\leq 1}\al(\R)}_{L^2(\MM^\ext)} \\
    & \quad +  \norm{\ub^{1+1/2+\ga} u^{-3/2}(u\Nd_3)\al(\R)}_{L^2(\MM^\ext)} \\
    & \quad +  \norm{\ub^{1+1/2+\ga} u^{-3/2}(\ub\Nd_4)\al(\R)}_{L^2(\MM^\ext)} \\
    & \quad +  \norm{\ub^{1+1/2+\ga} u^{-1/2} (r\Nd)^{\leq 1}\be(\R)}_{L^2(\MM^\ext)} \\
    & \quad +  \norm{\ub^{1+1/2+\ga} u^{-3/2}(\ub\Nd_4)\be(\R)}_{L^2(\MM^\ext)} \\
    & \quad +\norm{\ub^{1/2+\ga}u^{-1/2}(r\Nd)^{\leq 1}(\rho(\R),\sigma(\R))}_{L^2(\MM^\ext)}\bigg).
  \end{aligned}
\end{align}

We check that for $\ga < 1/4$, all the above $L^2(\MM^\ext)$ norms can be bounded by $\RR^\ext_{\leq 2, \gao}$, where we recall that $\gao < 1/4$. Thus, from the Bootstrap Assumptions~\ref{BA:curvext} on the $L^2(\MM^\ext)$ norms of the null curvature components, we have
\begin{align*}
  \norm{\ub^{3+1/2+\ga}\Xi\le(J\le(\Lieh_\TE\R\ri)\ri)}_{L^2(\MM^\ext)} & \les (D\varep) \RR^\ext_{\leq 2,\gao} \les (D\varep)^2, 
\end{align*}
as desired. This finishes the control of~\eqref{eq:intbexi}.\\ 



We now repeat the procedure to estimate 
\begin{align*}
  \int_{\MM^\ext}\DIV Q(\Lieh_\OOE\R)(\KE,\KE,\TE),
\end{align*}
which is the second error term of $\EE^\ext_{1,2}$.\\

We first record the following $L^2(\MM^\ext)$ bounds for the null decomposition of $\Lieh_\OOE\R$
\begin{align}\label{est:L2MMextW2}
  \begin{aligned}
    \norm{u^{-1/2-\ga}\ub^2\al\le(\Lieh_\OOE\R\ri)}_{L^2(\MM^\ext)} & \les_\ga D\varep, & \norm{u^{-1/2-\ga}\ub^2\be\le(\Lieh_\OOE\R\ri)}_{L^2(\MM^\ext)} & \les_\ga D\varep,\\
    \norm{\ub^{-1/2-\ga}\ub^{2}\rho\le(\Lieh_\OOE\R\ri)}_{L^2(\MM^\ext)} & \les_\ga D\varep, & \norm{\ub^{-1/2-\ga}\ub^2\sigma\le(\Lieh_\OOE\R\ri)}_{L^2(\MM^\ext)} & \les_\ga D\varep,\\
    \norm{\ub^{-1/2-\ga}\ub u\beb\le(\Lieh_\OOE\R\ri)}_{L^2(\MM^\ext)} & \les_\ga D\varep, & \norm{\ub^{-1/2-\ga}u^2\alb\le(\Lieh_\OOE\R\ri)}_{L^2(\MM^\ext)} & \les_\ga D\varep,
  \end{aligned}
\end{align}
for all $\ga>0$.\\

We have
\begin{align}\label{eq:DIVQSER}
  \begin{aligned}
    \DIV Q(\Lieh_\OOE\R)(\KE,\KE,\TE) & \les \ub^4 \DIV Q(\Lieh_\OOE\R)_{444} + \ub^4\DIV Q(\Lieh_\OOE\R)_{443} \\
    & \quad + \ub^2u^2\DIV Q(\Lieh_\OOE\R)_{433} + u^4\DIV Q(\Lieh_\OOE\R)_{333}.
  \end{aligned}
\end{align}

The most critical terms of~\eqref{eq:DIVQSER} are $\ub^4\DIV Q(\Lieh_\OOE\R)_{443}$ since it has the highest $\ub$ weight and the lowest signature.\\

We have
\begin{align}\label{eq:DIVQOOE443}
  \begin{aligned}
    \ub^4\DIV Q(\Lieh_\OOE\R)_{443} & = \ub^4 \rho\le(\Lieh_\OOE\R\ri)\Lambda\le(J(\Lieh_\OOE\R)\ri) \\
    & \quad + \ub^4\sigma\le(\Lieh_\OOE\R\ri)K\le(J(\Lieh_\OOE\R)\ri) \\
    & \quad + \ub^4 \beta\le(\Lieh_\OOE\R\ri)\cdot I(J(\Lieh_\OOE\R)).
  \end{aligned}
\end{align}




We only check the first error term of~\eqref{eq:DIVQOOE443}, the second will follow by duality and the third is easier to treat because of the stronger decay for $\be\le(\Lieh_\OOE\R\ri)$. Our goal is thus to obtain the following estimate
\begin{align*}
  \le|\int_{\MM^\ext}\ub^4\rho\le(\Lieh_\OOE\R\ri)\Lambda\le((J\le(\Lieh_\OOE\R\ri)\ri)) \ri| & \les (D\varep)^3.
\end{align*}

In view of the bounds~\eqref{est:L2MMextW2}, proving this last estimate reduces to obtain the following bound
\begin{align}
  \label{est:wishrholambda}
  \norm{\ub^{2+1/2+\ga}\Lambda\le(J\le(\Lieh_\OOE\R\ri)\ri)}_{L^2(\MM^\ext)} & \les (D\varep)^2.
\end{align}

Using the formulas from~\cite[p. 216]{Chr.Kla93}, the estimates~(\ref{est:DDfext0}) for the sup-norm and the estimates~(\ref{est:DDf1p1q}) for the $L^\infty L^4$ norm of the null decompositions of $\pi$ and $\D\pi$ respectively, we have
\begin{align}\label{est:LambdaJOOER1}
  \begin{aligned}
    \norm{\Lambda(J^1(\Lieh_\OOE\R))}_{L^2(S_{u,\ub})} & \les (D\varep)\bigg(\ub^{-1}u^{-1/2}\norm{\Nd\be(\R)}_{L^2(S_{u,\ub})} \\
    & \quad + \ub^{-1}u^{-1/2}\norm{\Nd_4\be(\R)}_{L^2(S_{u,\ub})} \\
    & \quad + \ub^{-1}u^{-1/2}\norm{(\rho_4(\R),\sigma_4(\R))}_{L^2(S_{u,\ub})} \\
    & \quad + \ub^{-2}u^{-1/2}\norm{\al(\R)}_{L^2(S_{u,\ub})} \\
    & \quad + \ub^{-2}u^{-1/2}\norm{\be(\R)}_{L^2(S_{u,\ub})} \\
    & \quad + \ub^{-2}u^{-1/2}\norm{\le(\rho(\R),\sigma(\R)\ri)}_{L^2(S_{u,\ub})}\bigg),
  \end{aligned}
\end{align}
and
\begin{align}\label{est:LambdaJOOER2}
  \begin{aligned}
    \norm{\Lambda(J^2(\Lieh_\OOE\R))}_{L^2(S_{u,\ub})} & \les (D\varep)\bigg(\ub^{-3/2} u^{-1/2} \norm{\be(\R)}_{L^4(S_{u,\ub})} \\
    & \quad + \ub^{-3/2}u^{-1/2} \norm{(\rho(\R),\sigma(\R))}_{L^4(S_{u,\ub})}\bigg), \\
    \norm{\Lambda(J^3(\Lieh_\OOE\R))}_{L^2(S_{u,\ub})} & \les (D\varep)\bigg(\ub^{-1/2}u^{-3/2}\norm{\al(\R)}_{L^4(S_{u,\ub})} \\
    & \quad + \ub^{-3/2}u^{-1/2}\norm{(\rho(\R),\sigma(\R))}_{L^4(S_{u,\ub})}\bigg).
  \end{aligned}
\end{align}

Arguing as previously, using Sobolev estimates on the $2$-spheres $S_{u,\ub}$, we deduce from~\eqref{est:LambdaJOOER1} and~\eqref{est:LambdaJOOER2}
\begin{align}\label{est:LambdaJL2MMext}
  \begin{aligned}
    \norm{\ub^{2+1/2+\ga}\Lambda\le(J\le(\Lieh_\OOE\R\ri)\ri)}_{L^2(\MM^\ext)} & \les (D\varep)\bigg(\norm{\ub^{1/2+\ga}u^{-1/2}(r\Nd)\be(\R)}_{L^2(\MM^\ext)} \\
    & \quad + \norm{\ub^{1/2+\ga}u^{-1/2}(\ub\Nd_4)\be(\R)}_{L^2(\MM^\ext)} \\
    & \quad + \norm{\ub^{1/2+\ga}u^{-1/2}(\ub\rho_4(\R),\ub\sigma_4(\R))}_{L^2(\MM^\ext)} \\
    & \quad + \norm{\ub^{3/2+\ga}u^{-3/2}(r\Nd)^{\leq 1}\al(\R)}_{L^2(\MM^\ext)} \\
    & \quad + \norm{\ub^{1/2+\ga}u^{-1/2}(r\Nd)^{\leq 1}\be(\R)}_{L^2(\MM^\ext)} \\
    & \quad + \norm{\ub^{1/2+\ga}u^{-3/2}(r\Nd)^{\leq 1}\le(\rho(\R),\sigma(\R)\ri)}_{L^2(\MM^\ext)}\bigg).
  \end{aligned}
\end{align}

We verify that from~\eqref{est:LambdaJL2MMext}, one has
\begin{align*}
  \norm{\ub^{2+1/2+\ga}\Lambda\le(J\le(\Lieh_\OOE\R\ri)\ri)}_{L^2(\MM^\ext)} & \les (D\varep)\RR^\ext_{\leq 2, \gao} \les (D\varep)^2,
\end{align*}
for $0< \ga < 1/4$, and \eqref{est:wishrholambda} is proved. This finishes the control of the error term $\EE^\ext_{1,2}$.



\subsection{Preliminary $L^2\le(\MM^\ext\ri)$ estimates for $\D^2\pi$}


We have the following definitions for the $L^2(\MM^\ext)$ norms for derivatives of $^{(\TE)}p, ^{(\TE)}q$ and $^{(\OOE)}p$, $^{(\OOE)}q$
\begin{align*}
  \DDf^\ext_{2,p, \gao}[\TE] & := \norm{\ub^{-1/2-\gao}u^2\Lieh_X {^{(\TE)}p_3}}_{L^2(\MM^\ext)} + \norm{\ub^{-1/2-\gao}\ub u\Lieh_X {^{(\TE)}p_4}}_{L^2(\MM^\ext)} \\
  & \quad + \norm{\ub^{-1/2-\gao}\ub u \Lieh_X {^{(\TE)}p_a}}_{L^2(\MM^\ext)}\\
  \DDf^\ext_{2,q, \gao}[\TE] & := \norm{\ub^{-1/2-\gao}\ub u\Lieh_X\Lambda(^{(\TE)}q)}_{L^2(\MM^\ext)} + \norm{\ub^{-1/2-\gao}\ub^{2}\Lieh_X K(^{(\TE)}q)}_{L^2(\MM^\ext)} \\
                       & \quad + \norm{\ub^{-1/2-\gao}\ub^{2}\Lieh_X\Xi(^{(\TE)}q)}_{L^2(\MM^\ext)} + \norm{\ub^{-1/2-\gao}\ub u\Lieh_X I(^{(\TE)}q)}_{L^2(\MM^\ext)} \\
                       & \quad + \norm{\ub^{-1/2-\gao}\ub u\Lieh_X\Theta(^{(\TE)}q)}_{L^2(\MM^\ext)} + \norm{\ub^{-1/2-\gao}u^{2}\Lieh_X\Lambdab(^{(\TE)}q)}_{L^2(\MM^\ext)} \\
                       & \quad + \norm{\ub^{-1/2-\gao}\ub u\Lieh_X\Kbb(^{(\TE)}q)}_{L^2(\MM^\ext)} + \norm{\ub^{-1/2-\gao}u^{2}\Lieh_X\Xib(^{(\TE)}q)}_{L^2(\MM^\ext)} \\
                       & \quad + \norm{\ub^{-1/2-\gao}\ub u\Lieh_X\Ib(^{(\TE)}q)}_{L^2(\MM^\ext)} + \norm{\ub^{-1/2-\gao}u^{2}\Lieh_X\Thetab(^{(\TE)}q)}_{L^2(\MM^\ext)},
\end{align*}
where $X \in \{\OOE,\SE\}$, and
\begin{align*}
  \DDf^\ext_{2,p, \gao}[\OOE] & := \norm{\ub^{-1/2-\gao}u\Lieh_\OOE{^{(\OOE)}p_3}}_{L^2(\MM^\ext)} + \norm{\ub^{-1/2-\gao}\ub \Lieh_\OOE{^{(\OOE)}p_4}}_{L^2(\MM^\ext)},\\
  & \quad + \norm{\ub^{-1/2-\gao}\ub \Lieh_\OOE{^{(\OOE)}p_a}}_{L^2(\MM^\ext)},\\
  \DDf^\ext_{2,q, \gao}[\OOE] & := \norm{\ub^{-1/2-\gao}\ub \Lieh_\OOE\Lambda(^{(\OOE)}q)}_{L^2(\MM^\ext))} + \norm{\ub^{-1/2-\gao}\ub^{2}u^{1/2}\Lieh_\OOE K(^{(\OOE)}q)}_{L^2(\MM^\ext)} \\
                       & \quad + \norm{\ub^{-1/2-\gao}\ub \Lieh_\OOE I(^{(\OOE)}q)}_{L^2(\MM^\ext)} \\
                       & \quad+ \norm{\ub^{-1/2-\gao}\ub \Lieh_\OOE\Theta(^{(\OOE)}q)}_{L^2(\MM^\ext)} + \norm{\ub^{-1/2-\gao}u\Lieh_\OOE\Lambdab(^{(\OOE)}q)}_{L^2(\MM^\ext)} \\
                       & \quad + \norm{\ub^{-1/2-\gao}\ub \Lieh_\OOE\Kbb(^{(\OOE)}q)}_{L^2(\MM^\ext)} + \norm{\ub^{-1/2-\gao}u\Lieh_\OOE\Xib(^{(\OOE)}q)}_{L^2(\MM^\ext)} \\
                       & \quad + \norm{\ub^{-1/2-\gao}\ub \Lieh_\OOE\Ib(^{(\OOE)}q)}_{L^2(\MM^\ext)} + \norm{\ub^{-1/2-\gao}u\Lieh_\OOE\Thetab(^{(\OOE)}q)}_{L^2(\MM^\ext)},
\end{align*}
together with $\Lieh_\OOE\Xi(^{(\OOE)}q) = 0$.\\

Differentiating the formulas for $p,q$ from~\cite[pp. 231--232]{Chr.Kla93} and the formulas for the null decompositions of $^{(\TE)}\pih$ and $^{(\OOE)}\pih$ from Sections~\ref{sec:exterrpreldef}, using the Bootstrap Assumptions~\ref{BA:connext} for the $L^2(\MM^\ext)$ norms of (two derivatives) of the null connection coefficients, we obtain
\begin{align}\label{est:DDfext2}
  \DDf^\ext_{2, \gao} & \les \OO^\ext_{\leq 2, \gao} \les D\varep,
\end{align}
for all $\ga>0$.\\

\begin{remark}
  The decay rates are easily checked, using that deriving with respect to $\OOE$ or $\SE$ does not change the asymptotic behaviour of the components.
\end{remark}

\subsection{Estimates for $\EE^\ext_{2,2}$}\label{sec:esterrEE22}
We start with the estimate for
\begin{align}\label{eq:intL2OOE}
  \int_{\MM^\ext}\DIV Q(\Lieh_\OOE^2\R)(\KE,\KE,\TE).
\end{align}

We first record the following $L^2(\MM^\ext)$ bounds for the null decomposition of $\Lieh_\OOE^2\R$ which are consequences of the Bootstrap Assumptions~\ref{BA:fluxglobener}
\begin{align}\label{est:L2MMextW3}
  \begin{aligned}
    \norm{u^{-1/2-\ga}\ub^2\al\le(\Lieh_\OOE^2\R\ri)}_{L^2(\MM^\ext)} & \les_\ga D\varep, & \norm{u^{-1/2-\ga}\ub^2\be\le(\Lieh_\OOE^2\R\ri)}_{L^2(\MM^\ext)} & \les_\ga D\varep,\\
    \norm{\ub^{-1/2-\ga}\ub^{2}\rho\le(\Lieh_\OOE^2\R\ri)}_{L^2(\MM^\ext)} & \les_\ga D\varep, & \norm{\ub^{-1/2-\ga}\ub^2\sigma\le(\Lieh_\OOE^2\R\ri)}_{L^2(\MM^\ext)} & \les_\ga D\varep,\\
    \norm{\ub^{-1/2-\ga}\ub u\beb\le(\Lieh_\OOE^2\R\ri)}_{L^2(\MM^\ext)} & \les_\ga D\varep, & \norm{\ub^{-1/2-\ga}u^2\alb\le(\Lieh_\OOE^2\R\ri)}_{L^2(\MM^\ext)} & \les_\ga D\varep,
  \end{aligned}
\end{align}
for all $\ga>0$.\\

Similarly to Section~\ref{sec:esterrEE12}, the decompositions on the null directions $(\elb,\el)$ gives
\begin{align}\label{eq:DIVQOOEOOER}
  \begin{aligned}
    \DIV Q(\Lieh^2_\OOE\R)(\KE,\KE,\TE) & \les \ub^4\DIV Q(\Lieh^2_\OOE\R)_{444} + \ub^4\DIV Q(\Lieh^2_\OOE\R)_{443} \\
    & \quad + \ub^2u^2 \DIV Q(\Lieh^2_\OOE\R)_{433} + u^4\DIV Q(\Lieh^2_\OOE\R)_{333}.
  \end{aligned}
\end{align}
and
\begin{align}\label{eq:nullDIVLOOEOOER}
  \begin{aligned}
    \ub^4\DIV Q(\Lieh^2_\OOE\R)_{444} & = \ub^4\al\le(\Lieh^2_\OOE\R\ri)\Theta(J(\Lieh^2_\OOE\R)) \\
    & \quad + \ub^4\be\le(\Lieh^2_\OOE\R\ri)\Xi(J(\Lieh^2_\OOE\R)),\\
    \ub^4\DIV Q(\Lieh^2_\OOE\R)_{443} & = \ub^4\rho\le(\Lieh^2_\OOE\R\ri)\Lambda(J(\Lieh^2_\OOE\R)) \\
    &\quad + \ub^4 \sigma\le(\Lieh^2_\OOE\R\ri)K(J(\Lieh^2_\OOE\R)) \\
    & \quad + \ub^4\be\le(\Lieh^2_\OOE\R\ri)I(J(\Lieh^2_\OOE\R)),\\
    \ub^2u^2\DIV Q(\Lieh^2_\OOE\R)_{433} & = \ub^2u^2\rho\le(\Lieh^2_\OOE\R\ri)\Lambdab(J(\Lieh^2_\OOE\R)) \\
    & \quad + \ub^2u^2 \sigma\le(\Lieh^2_\OOE\R\ri)\Kb(J(\Lieh^2_\OOE\R)) \\
    & \quad + \ub^2u^2\be\le(\Lieh^2_\OOE\R\ri)\Ib(J(\Lieh^2_\OOE\R)),\\
    u^4\DIV Q(\Lieh^2_\OOE\R)_{333} & = u^4\alb\le(\Lieh^2_\OOE\R\ri)\Thetab(J(\Lieh^2_\OOE\R)) \\
    & \quad + u^4\beb\le(\Lieh^2_\OOE\R\ri)\Xib(J(\Lieh^2_\OOE\R)).
  \end{aligned}
\end{align}

We have from the formulas~\cite[p. 206]{Chr.Kla93}
\begin{align*}
  J(\Lieh^2_\OOE\R) & = \Lieh_\OOE J(\Lieh_\OOE\R) + J^i(\Lieh_\OOE\Lieh_\OOE\R),
\end{align*}
where the terms $J^i(\Lieh_\OOE\Lieh_\OOE\R)$ are the same as the terms $J^i(\Lieh_\OOE\R)$ from Section~\ref{sec:esterrEE12} with $\R$ replaced by $\Lieh_\OOE\R$. Therefore, the estimates of that section carry over and we only need to treat the term
\begin{align*}
  \Lieh_\OOE J(\Lieh_\OOE\R).
\end{align*}
This term again has the same structure as the terms treated in Section~\ref{sec:esterrEE12}, although differentiated by $\Lieh_\OOE$. This does not change the decay of the components, and the difference with Section~\ref{sec:esterrEE12} is that one does not have a control for the $L^\infty_{u,\ub}L^4(S_{u,\ub})$ norm of the components of the type $\Lieh_\OOE\D\pih$.\\

That case is actually (more) simply handled using the $L^2(\MM^\ext)$ estimates~(\ref{est:L2MMextW3}) for the null decomposition of $\Lieh_\OOE^2\R$, the $L^2(\MM^\ext)$ estimates~(\ref{est:DDfext2}) for $\Lieh_\OOE\D\pih$, the Bootstrap Assumptions~\ref{BA:curvext} for the sup-norm of the curvature $\R$ and Cauchy-Schwartz. As an example, we treat the first term of~\eqref{eq:nullDIVLOOEOOER} containing a term of the type $\D^2\pih$ which is
\begin{align*}
  \ub^4\al\le(\Lieh^2_\OOE\R\ri)\Lieh_\OOE {^{(\OOE)}p_3} \cdot \al(\R).  
\end{align*}

We have
\begin{align*}
  & \le|\int_{\MM^\ext}\ub^4\al\le(\Lieh^2_\OOE\R\ri)\Lieh_\OOE {^{(\OOE)}p_3} \cdot \al(\R)\ri| \\
  \les & \; (D\varep)\int_{\MM^\ext}\ub^{1/2}\le|\al\le(\Lieh^2_\OOE\R\ri)\ri|\le|\Lieh_\OOE{^{(\OOE)}p_3}\ri|\\
  \les & \; (D\varep) \norm{u^{-1/2-\ga}\ub^2\al\le(\Lieh^2_\OOE\R\ri)}_{L^2(\MM^\ext)}\norm{u^{1/2+\ga}\ub^{-3/2}\Lieh_\OOE{^{(\OOE)}p_3}}_{L^2(\MM^\ext)} \\
  \les & \; (D\varep) \norm{u^{-1/2-\ga}\ub^2\al\le(\Lieh^2_\OOE\R\ri)}_{L^2(\MM^\ext)} \norm{\ub^{-1/2-\gao}u\Lieh_\OOE{^{(\OOE)}p_3}}_{L^2(\MM^\ext)} \\
  \les & \; (D\varep)^3,
\end{align*}
provided that $0 < \ga < 1/2$.\\

All the other terms of~\eqref{eq:intL2OOE} follow similarly, using the analysis of the decay rates already performed in Section~\ref{sec:esterrEE12}.\\

To handle the last error terms of $\EE^\ext_{2,2}$, we have the following lemma.
\begin{lemma}\label{lem:decayerrtermOOESEW}
  Let $W$ be a current-free Weyl field such that the following $L^2(\MM^\ext)$ bounds hold
  \begin{align}\label{est:L2MMEXTW}
    \begin{aligned}
      \norm{u^{-1/2-\ga}\ub^3\Ndt^{\leq 1}\al\le(W\ri)}_{L^2(\MM^\ext)} & \les_\ga D\varep, & \norm{\ub^{-1/2-\ga}\ub^3\Ndt^{\leq 1}\be\le(W\ri)}_{L^2(\MM^\ext)} & \les_\ga D\varep,\\
      \norm{\ub^{-1/2-\ga}\ub^{2}u\Ndt^{\leq 1}\rho\le(W\ri)}_{L^2(\MM^\ext)} & \les_\ga D\varep, & \norm{\ub^{-1/2-\ga}\ub^2u\Ndt^{\leq 1}\sigma\le(W\ri)}_{L^2(\MM^\ext)} & \les_\ga D\varep,\\
      \norm{\ub^{-1/2-\ga}\ub u^2\Ndt^{\leq 1}\beb\le(W\ri)}_{L^2(\MM^\ext)} & \les_\ga D\varep, & \norm{\ub^{-1/2-\ga}u^3\Ndt^{\leq 1}\alb\le(W\ri)}_{L^2(\MM^\ext)} & \les_\ga D\varep,
    \end{aligned}
  \end{align}
  for all $\ga>0$ and where $\Ndt \in\le\{r\Nd, u\Nd_3, \ub\Nd_4 \ri\}$. Then, we have
  \begin{align*}
    \int_{\MM^\ext}\DIV Q(\Lieh_\OOE W)(\KE,\KE,\KE) & \les (D\varep)^3,\\
    \int_{\MM^\ext}\DIV Q(\Lieh_\SE W)(\KE,\KE,\KE) & \les (D\varep)^3.
  \end{align*}
\end{lemma}
The proof of Lemma~\ref{lem:decayerrtermOOESEW} is similar to the control of the error terms from Section~\ref{sec:esterrEE12} and consists in checking that the integrands have the appropriate $u,\ub$ decay and are integrable on $\MM^\ext$. It is postponed to Section~\ref{sec:prooflemdecayesterrtermOOESEW}.\\

The estimates for the terms
\begin{align*}
  \int_{\MM^\ext}\DIV Q(\Lieh_\OOE\Lieh_\TE\R)(\KE,\KE,\KE)
\end{align*}
and
\begin{align*}
  \int_{\MM^\ext}\DIV Q(\Lieh_\SE\Lieh_\TE\R)(\KE,\KE,\KE)
\end{align*}
now follow from Lemma~\ref{lem:decayerrtermOOESEW}, arguing as above, and using the decay rates analysis for
\begin{align*}
  \int_{\MM^\ext}\DIV Q(\Lieh_\TE\R)(\KE,\KE,\KE)
\end{align*}
performed in Section~\ref{sec:esterrEE12} to handle the terms from the non-vanishing current $\Lieh_\OOE J\le(\Lieh_\TE\R\ri)$ and $\Lieh_\SE J\le(\Lieh_\TE\R\ri)$. This finishes the control of the error terms $\EE^\ext$.

\subsection{Proof of Lemma~\ref{lem:decayerrtermOOESEW}}\label{sec:prooflemdecayesterrtermOOESEW}
Since from an inspection of~\eqref{est:DDfext0} and~\eqref{est:DDf1p1q}, the deformation tensors of $\OOE$ satisfy better decay estimates than the deformation tensor of $\SE$, it is enough to control
\begin{align*}
  \int_{\MM^\ext}\DIV Q(\Lieh_\SE W)(\KE,\KE,\KE).
\end{align*}

Decomposing the divergence on the null directions $(\elb,\el)$, the most difficult terms to treat are
\begin{align*}
  \int_{\MM^\ext}\ub^6 \be\le(\Lieh_\SE W\ri)\cdot\Xi(J(\Lieh_\SE W)),
\end{align*}
since it has the highest $\ub$ weight and the $\be$ component satisfies a weaker control than $\al$.\\

Using estimates~\eqref{est:DDfext0} and~\eqref{est:DDf1p1q} for the null decomposition of $^{(\SE)}\pih$, the formulas from~\cite[p. 215--217]{Chr.Kla93}, and Sobolev estimates on the $2$-sphere $S_{u,\ub}$ we have 
\begin{align*}
  \begin{aligned}
    \norm{\Xi(J^1(\Lieh_\SE W))}_{L^2(S_{u,\ub})} & \les (D\varep)\bigg(\ub^{-1}u^{-1/2} \norm{\Nd\al(W)}_{L^2(S_{u,\ub})} \\
    & \quad + \ub^{-2}u^{1/2} \norm{\Nd_3\al(W)}_{L^2(S_{u,\ub})} \\
    & \quad + \ub^{-1} u^{-1/2} \norm{\Nd_4\al(W)}_{L^2(S_{u,\ub})} \\
    & \quad + \ub^{-2}u^{1/2} \norm{\Nd\be(W)}_{L^2(S_{u,\ub})} \\
    & \quad + \ub^{-1}u^{-1/2} \norm{\Nd_4\be(W)}_{L^2(S_{u,\ub})} \\
    & \quad + \ub^{-2}u^{-1/2}\norm{\al(W)}_{L^2(S_{u,\ub})} \\
    & \quad + \ub^{-2}u^{-1/2}\norm{\be(W)}_{L^2(S_{u,\ub})} \\
    & \quad + \ub^{-3}u^{1/2}\norm{(\rho(W),\sigma(W))}_{L^2(S_{u,\ub})}\bigg),
  \end{aligned}
\end{align*}
and
\begin{align*}
  \begin{aligned}
    \norm{\Xi(J^2(\Lieh_\SE\R))}_{L^2(S_{u,\ub})} & \les (D\varep)\bigg(\ub^{-2}u^{-1/2} \norm{(r\Nd)^{\leq 1}\al(W)}_{L^2(S_{u,\ub})} \\
    & \quad + \ub^{-2} u^{-1/2}\norm{(r\Nd)^{\leq 1}\be(W)}_{L^2(S_{u,\ub})}\bigg), \\
    \norm{\Xi(J^3(\Lieh_\SE\R))}_{L^2(S_{u,\ub})}  & \les (D\varep)\bigg(\ub^{-2} u^{-1/2} \norm{(r\Nd)^{\leq 1}\al(W)}_{L^2(S_{u,\ub})} \\
    & \quad + \ub^{-2} u^{-1/2}\norm{(r\Nd)^{\leq 1}\be(W)}_{L^2(S_{u,\ub})} \\
    & \quad + \ub^{-3}u^{1/2} \norm{(r\Nd)^{\leq 1}(\rho(W),\sigma(W))}_{L^2(S_{u,\ub})}\bigg).
  \end{aligned}
\end{align*}

Thus,
\begin{align*}
  \begin{aligned}
    \norm{\ub^{3+1/2+\ga}\Xi(J(\Lieh_\SE W))}_{L^2(\MM^\ext)} & \les (D\varep)\bigg(\norm{\ub^{1+1/2+\ga}u^{-1/2}(r\Nd)^{\leq 1}\al(W)}_{L^2(\MM^\ext)} \\
    & \quad + \norm{\ub^{1+1/2+\ga}u^{-1/2} (u\Nd_3)\al(W)}_{L^2(\MM^\ext)} \\
    & \quad +  \norm{\ub^{1+1/2+\ga} u^{-1/2}(\ub\Nd_4)\al(W)}_{L^2(\MM^\ext)} \\
    & \quad +  \norm{\ub^{1+1/2+\ga}u^{-1/2}(r\Nd)^{\leq 1}\be(W)}_{L^2(\MM^\ext)} \\
    & \quad +  \norm{\ub^{1+1/2+\ga}u^{-1/2}(\ub\Nd_4)\be(W)}_{L^2(\MM^\ext)} \\
    & \quad + \norm{\ub^{1/2+\ga}u^{1/2}(r\Nd)^{\leq 1}(\rho(W),\sigma(W))}_{L^2(\MM^\ext)}\bigg),
  \end{aligned}
\end{align*}
and using the $L^2(\MM^\ext)$ bounds~\eqref{est:L2MMEXTW}, we have
\begin{align*}
  \norm{\ub^{3+1/2+\ga}\Xi(J(\Lieh_\SE W))}_{L^2(\MM^\ext)} & \les (D\varep)^2,
\end{align*}
provided that $0< \ga < 1$.
Arguing as previously, using Cauchy-Schwartz and the $L^2(\MM^\ext)$ bounds~\eqref{est:L2MMEXTW} gives the desired control. This finishes the proof of Lemma~\ref{lem:decayerrtermOOESEW}.

\chapter{Null curvature estimates in $\protect\CCba\cap\MM^\ext$ and $\MM^\ext$}\label{sec:curvest}
In this section, we prove the following proposition.

\begin{proposition}\label{prop:curvestSTAB}
  Recall that from Proposition~\ref{prop:enerestSTAB}, we have
  \begin{align}\label{est:enerestasscurvest}
    \int_{^{(\cc)}\Si_t^\ext}P\cdot\Tf^\ext + \int_{\CC_u\cap{^{(\cc)}\MM}^\ext}P\cdot\el + \int_{\CCba\cap{^{(\cc)}\MM}^\ext}P\cdot\elb & \les \varep^2,
  \end{align}
  for a fixed transition parameter $\cc$, for all $1\leq u \leq \cc\uba$ and for all $\too \leq t \leq \tast$, and where $P$ denote the following contracted and commuted Bel-Robinson tensors
  \begin{align*}
    & Q\le(\Lieh_\TE\R\ri)(\KE,\KE,\KE), \quad Q\le(\Lieh_{\OOE}\R\ri)(\KE,\KE,\TE), \\
    & Q\le(\Lieh_\OOE\Lieh_\OOE\R\ri)(\KE,\KE,\TE), \quad Q\le(\Lieh_\SE\Lieh_\TE\R\ri)(\KE,\KE,\KE),\\
    & Q\le(\Lieh_\OOE\Lieh_\TE\R\ri)(\KE,\KE,\TE).
  \end{align*}

  Under the Bootstrap Assumptions, the energy estimates~\eqref{est:enerestasscurvest}, and for $\varep>0$ sufficiently small, we have
  \begin{subequations}\label{est:curvext}
  \begin{itemize}
  \item the following $L^2$ bounds on the exterior cone $\CCba\cap{^{(\cc)}\MM}^\ext$ 
    \begin{align}\label{est:RRCCbaext}
      \RR^{\ast}_{\leq 2} & \les \varep,
    \end{align}
  \item the following $L^2$ bounds on ${^{(\cc)}\MM}^\ext$
    \begin{align}\label{est:RRextMM}
      \RR^{\ext}_{\leq 2,\ga} & \les_\ga \varep,
    \end{align}
    for all $\ga>0$,\footnote{In particular, the bound holds with $\ga = \gao$ where $0<\gao<1/4$ is the fixed numerical constant of Section~\ref{sec:strongBA}. This improves the bound from the Bootstrap Assumption~\ref{BA:curvext}.}
  \item the following $L^\infty_{u}\HHt$ bounds on $\CCba\cap{^{(\cc)}\MM}^\ext$
    \begin{align}\label{est:RRCCbadecaynew}
      \mathfrak{R}^{\ast}_{\leq 1} & \les \varep,
    \end{align}
    where the norms are restricted to $\CCba\cap{^{(\cc)}\MM}^\ext$,
  \item the following $L^\infty_{u, \ub}\HHt$ bounds in ${^{(\cc)}\MM}^\ext$
    \begin{align}\label{est:RRextdecay}
      \mathfrak{R}^\ext_{\leq 1} & \les \varep.
    \end{align}
  \end{itemize}
  We refer the reader to the norm definitions of Section~\ref{sec:normnullcurv}.
\end{subequations}
\end{proposition}

\begin{remark}
  Proposition~\ref{prop:curvestSTAB} does not provide bounds for the mean value $\rhoo,\sigmao$. These are obtained in Sections~\ref{sec:connestCCba} and \ref{sec:connest}. See also Remark~\ref{rem:fasterdecayerrorrhoosigmaoo}.
\end{remark}


The proof of the estimates~\eqref{est:curvext} relies on the following localised control on the $2$-spheres $S_{u,\ub}$ for the null curvature tensors and their derivatives in terms of the (contracted) Bel-Robinson tensors used in Section~\ref{sec:globener}. The proof of Proposition~\ref{prop:curvS} is provided in Section~\ref{sec:proofcurvS}.
\begin{proposition}\label{prop:curvS}
  On each $2$-sphere $S_{u,\ub}$ of the exterior region $\MM^\ext$, the following control holds
  \begin{subequations}\label{est:RRRRb}
    \begin{align}
      \le(\RRb_{\leq 1}(u,\ub)\ri)^2 & \les \QQb_1(u,\ub) +(D\varep)^2\norm{\ub(\rhoo,\sigmao)}_{L^2(S_{u,\ub})}^2, \label{est:RRb1}\\
      \le(\RR_{\leq 1}(u,\ub)\ri)^2 & \les \QQ_{1}(u,\ub) +(D\varep)^2\norm{\ub(\rhoo,\sigmao)}_{L^2(S_{u,\ub})}^2, \label{est:RR1}\\
      \le(\RRb_{\leq 2}(u,\ub)\ri)^2 & \les \QQb_{\leq 2}(u,\ub) +(D\varep)^2\norm{\ub(\rhoo,\sigmao)}_{L^2(S_{u,\ub})}^2,\label{est:RRb2}\\
      \le(\RR_{\leq 2}(u,\ub)\ri)^2 & \les \QQ_{\leq 2}(u,\ub) +(D\varep)^2\norm{\ub(\rhoo,\sigmao)}_{L^2(S_{u,\ub})}^2,\label{est:RR2}
    \end{align}
  \end{subequations}
  where we have the following definitions
  \begin{align*}
    \RRb_0[\alb](u,\ub) & :=  \norm{u^2\alb}_{L^2(S_{u,\ub})} , & \RRb_0[\beb](u,\ub) & := \norm{\ub u \beb}_{L^2(S_{u,\ub})} ,\\
    \RRb_0[\rho](u,\ub) & := \norm{\ub^2(\rho-\rhoo)}_{L^2(S_{u,\ub})} , & \RRb_0[\sigma](u,\ub) & := \norm{\ub^2(\sigma-\sigmao)}_{L^2(S_{u,\ub})} ,\\
    \RRb_0[\be](u,\ub) & := \norm{\ub^2\be}_{L^2(S_{u,\ub})}, & \RRb^{add}_{0}[\al](u,\ub) & := \norm{\ub^2\al}_{L^2(S_{u,\ub})}, 
  \end{align*}
  and
  \begin{align*}
    \RRb_1[\alb](u,\ub) & := \norm{u^2(r\Nd)\alb}_{L^2(S_{u,\ub})} + \norm{u^2(r\Nd_4)\alb}_{L^2(S_{u,\ub})} + \norm{u^2(\qq\Nd_3)\alb}_{L^2(S_{u,\ub})}, \\
    \RRb_1[\beb](u,\ub) & := \norm{u\ub (r\Nd)\beb}_{L^2(S_{u,\ub})} + \norm{u \ub (r\Nd_4)\beb}_{L^2(S_{u,\ub})} + \norm{u \ub (\qq\Nd_3)\beb}_{L^2(S_{u,\ub})}, \\
    \RRb_1[\rho](u,\ub) & := \norm{\ub^2 (r\Nd)\rho}_{L^2(S_{u,\ub})} + \norm{\ub^2 r \Nd_4(\rho-\rhoo)}_{L^2(S_{u,\ub})} + \norm{\ub^2 \qq \Nd_3(\rho-\rhoo)}_{L^2(S_{u,\ub})}, \\
    \RRb_1[\sigma](u,\ub) & := \norm{\ub^2 (r\Nd)\sigma}_{L^2(S_{u,\ub})} + \norm{\ub^2 r\Nd_4(\sigma-\sigmao)}_{L^2(S_{u,\ub})} + \norm{\ub^2 \qq\Nd_3(\sigma-\sigmao)}_{L^2(S_{u,\ub})}, \\
    \RRb_1[\be](u,\ub) & := \norm{\ub^2 (r\Nd)\be}_{L^2(S_{u,\ub})} + \norm{\ub^2 (r\Nd_4)\be}_{L^2(S_{u,\ub})} + \norm{\ub^2 (r\Nd_3)\be}_{L^2(S_{u,\ub})}, \\
    \RRb_1[\al](u,\ub) & := \norm{\ub^2 (r\Nd)\al}_{L^2(S_{u,\ub})} + \norm{\ub^2 (r\Nd_3)\al}_{L^2(S_{u,\ub})}, \\
    \RRb_{1}^{add}[\al](u,\ub) & := \norm{\ub^2(r\Nd_4)\al}_{L^2(S_{u,\ub})},
  \end{align*}
  and
  \begin{align*}
    \RRb_2[\alb](u,\ub) & := \norm{u^2r^2\Nd^2\alb}_{L^2(S_{u,\ub})} + \norm{u^2r^2\Nd\Nd_4\alb}_{L^2(S_{u,\ub})}+ \norm{u^2r^2\Nd_4\Nd_4\alb}_{L^2(S_{u,\ub})} \\
                        & \quad + \norm{u^2r\qq\Nd\Nd_3\alb}_{L^2(S_{u,\ub})} + \norm{u^2\qq^2\Nd_3\Nd_3\alb}_{L^2(S_{u,\ub})} + \norm{u^2(r\Nd_3)(r\Nd_4)\alb}_{L^2(S_{u,\ub})}, \\ \\
    \RRb_2[\beb](u,\ub) & :=  \norm{u\ub r^2\Nd^2\beb}_{L^2(S_{u,\ub})} + \norm{u \ub r^2\Nd\Nd_4\beb}_{L^2(S_{u,\ub})} + \norm{u\ub r^2\Nd_4\Nd_4\beb}_{L^2(S_{u,\ub})}\\
                        & \quad + \norm{u\ub r\qq\Nd\Nd_3\beb}_{L^2(S_{u,\ub})} + \norm{u\ub\qq^2\Nd_3\Nd_3\beb}_{L^2(S_{u,\ub})} + \norm{u \ub(r\Nd_3)(r\Nd_4)\beb}_{L^2(S_{u,\ub})}, \\ \\
    \RRb_2[\rho](u,\ub) & :=  \norm{\ub^2r^2\Nd^2\rho}_{L^2(S_{u,\ub})} + \norm{\ub^2r^2\Nd\Nd_4(\rho-\rhoo)}_{L^2(S_{u,\ub})} + \norm{\ub^2r^2\Nd_4\Nd_4(\rho-\rhoo)}_{L^2(S_{u,\ub})} \\
                        & \quad + \norm{\ub^2r \qq\Nd\Nd_3(\rho-\rhoo)}_{L^2(S_{u,\ub})} + \norm{u^2\qq^2\Nd_3\Nd_3(\rho-\rhoo)}_{L^2(S_{u,\ub})} \\
                        &\quad + \norm{\ub^2(r\Nd_3)r\Nd_4(\rho-\rhoo)}_{L^2(S_{u,\ub})}, \\ \\
    \RRb_2[\sigma](u,\ub) & :=  \norm{\ub^2r^2\Nd^2\sigma}_{L^2(S_{u,\ub})} + \norm{\ub^2r^2\Nd\Nd_4(\sigma-\sigmao)}_{L^2(S_{u,\ub})} + \norm{\ub^2r^2\Nd_4\Nd_4(\sigma-\sigmao)}_{L^2(S_{u,\ub})} \\
                        & \quad + \norm{\ub^2r\qq\Nd\Nd_3(\sigma-\sigmao)}_{L^2(S_{u,\ub})} + \norm{u^2\qq^2\Nd_3\Nd_3(\sigma-\sigmao)}_{L^2(S_{u,\ub})} \\
                        & \quad + \norm{\ub^2(r\Nd_3)r\Nd_4(\sigma-\sigmao)}_{L^2(S_{u,\ub})}, \\ \\
    \RRb_2[\be](u,\ub) & := \norm{\ub^2r^2\Nd^2\be}_{L^2(S_{u,\ub})} + \norm{\ub^2r^2\Nd\Nd_4\be}_{L^2(S_{u,\ub})} + \norm{\ub^2r^2\Nd_4\Nd_4\be}_{L^2(S_{u,\ub})} \\
                        & \quad + \norm{\ub^2r^2\Nd\Nd_3\be}_{L^2(S_{u,\ub})} + \norm{u^2r\qq\Nd_3\Nd_3\be}_{L^2(S_{u,\ub})} + \norm{\ub^2(r\Nd_3)(r\Nd_4)\be}_{L^2(S_{u,\ub})}, \\ \\
    \RRb_2[\al](u,\ub) & :=  \norm{\ub^2r^2\Nd^2\al}_{L^2(S_{u,\ub})} + \norm{\ub^2r^2\Nd\Nd_4\al}_{L^2(S_{u,\ub})} + \norm{\ub^2r^2\Nd\Nd_3\al}_{L^2(S_{u,\ub})} \\
                        & \quad + \norm{\ub^2r^2\Nd_3\Nd_3\al}_{L^2(S_{u,\ub})} + \norm{\ub^2(r\Nd_3)(r\Nd_4)\al}_{L^2(S_{u,\ub})}.\\
  \end{align*}
  We moreover define the dual norms
  \begin{align*}
    \RR_{0}^{add}[\alb](u,\ub) & :=  \norm{u^2\alb}_{L^2(S_{u,\ub})} , & \RR_0[\beb](u,\ub) & := \norm{u^2 \beb}_{L^2(S_{u,\ub})} ,\\
    \RR_0[\rho](u,\ub) & := \norm{u\ub(\rho-\rhoo)}_{L^2(S_{u,\ub})} , & \RR_0[\sigma](u,\ub) & := \norm{u\ub(\sigma-\sigmao)}_{L^2(S_{u,\ub})} ,\\
    \RR_0[\be](u,\ub) & := \norm{\ub^2\be}_{L^2(S_{u,\ub})}, & \RR_0[\al](u,\ub) & := \norm{\ub^2\al}_{L^2(S_{u,\ub})}, 
  \end{align*}
  and
  \begin{align*}
    \RR_1[\alb](u,\ub) & := \norm{u^2(\qq\Nd)\alb}_{L^2(S_{u,\ub})} + \norm{u^2(r\Nd_4)\alb}_{L^2(S_{u,\ub})}, \\
    \RR_{1}^{add}[\alb](u,\ub) & := \norm{u^2(\qq\Nd_3)\alb}_{L^2(S_{u,\ub})},\\
    \RR_1[\beb](u,\ub) & := \norm{u^2(r\Nd)\beb}_{L^2(S_{u,\ub})} + \norm{u^2(r\Nd_4)\beb}_{L^2(S_{u,\ub})} + \norm{u^2(\qq\Nd_3)\beb}_{L^2(S_{u,\ub})}, \\
    \RR_1[\rho](u,\ub) & := \norm{\ub u (r\Nd)\rho}_{L^2(S_{u,\ub})} + \norm{\ub u r\Nd_4(\rho-\rhoo)}_{L^2(S_{u,\ub})} + \norm{\ub u \qq\Nd_3(\rho-\rhoo)}_{L^2(S_{u,\ub})}, \\
    \RR_1[\sigma](u,\ub) & := \norm{\ub u (r\Nd)\sigma}_{L^2(S_{u,\ub})} + \norm{\ub u r\Nd_4(\sigma-\sigmao)}_{L^2(S_{u,\ub})} + \norm{\ub u \qq\Nd_3(\sigma-\sigmao)}_{L^2(S_{u,\ub})}, \\
    \RR_1[\be](u,\ub) & := \norm{\ub^2 (r\Nd)\be}_{L^2(S_{u,\ub})} + \norm{\ub^2 (r\Nd_4)\be}_{L^2(S_{u,\ub})} + \norm{\ub^2 (\qq\Nd_3)\be}_{L^2(S_{u,\ub})}, \\
    \RR_1[\al](u,\ub) & := \norm{\ub^2 (r\Nd)\al}_{L^2(S_{u,\ub})} + \norm{\ub^2 (r\Nd_3)\al}_{L^2(S_{u,\ub})} + \norm{\ub^2(r\Nd_4)\al}_{L^2(S_{u,\ub})},
  \end{align*}
  and
  \begin{align*}
    \RR_2[\alb](u,\ub) & := \norm{u^2\qq r\Nd^2\alb}_{L^2(S_{u,\ub})} + \norm{u^2r^2\Nd\Nd_4\alb}_{L^2(S_{u,\ub})}+ \norm{u^2r^2\Nd_4\Nd_4\alb}_{L^2(S_{u,\ub})} \\
                       & \quad + \norm{u^2\qq^2\Nd\Nd_3\alb}_{L^2(S_{u,\ub})} + \norm{u^2(\qq\Nd_3)(r\Nd_4)\alb}_{L^2(S_{u,\ub})}, \\ \\
    \RR_2[\beb](u,\ub) & :=  \norm{u^2 r^2\Nd^2\beb}_{L^2(S_{u,\ub})} + \norm{u^2 r^2\Nd\Nd_4\beb}_{L^2(S_{u,\ub})} + \norm{u^2 r^2\Nd_4\Nd_4\beb}_{L^2(S_{u,\ub})}\\
                       & \quad + \norm{u^2 r\qq\Nd\Nd_3\beb}_{L^2(S_{u,\ub})} + \norm{u^2\qq^2\Nd_3\Nd_3\beb}_{L^2(S_{u,\ub})} + \norm{u^2(r\Nd_3)(r\Nd_4)\beb}_{L^2(S_{u,\ub})}, \\ \\
    \RR_2[\rho](u,\ub) & :=  \norm{\ub u r^2\Nd^2\rho}_{L^2(S_{u,\ub})} + \norm{\ub u r^2\Nd\Nd_4(\rho-\rhoo)}_{L^2(S_{u,\ub})} + \norm{\ub u r^2\Nd_4\Nd_4(\rho-\rhoo)}_{L^2(S_{u,\ub})} \\
                       & \quad + \norm{\ub u r \qq\Nd\Nd_3(\rho-\rhoo)}_{L^2(S_{u,\ub})} + \norm{\ub u \qq^2\Nd_3\Nd_3(\rho-\rhoo)}_{L^2(S_{u,\ub})} \\
                       & \quad + \norm{\ub u (r\Nd_3)r\Nd_4(\rho-\rhoo)}_{L^2(S_{u,\ub})}, \\ \\
    \RR_2[\sigma](u,\ub) & :=  \norm{\ub u r^2\Nd^2\sigma}_{L^2(S_{u,\ub})} + \norm{\ub u r^2\Nd\Nd_4(\sigma-\sigmao)}_{L^2(S_{u,\ub})} + \norm{\ub u r^2\Nd_4\Nd_4(\sigma-\sigmao)}_{L^2(S_{u,\ub})} \\
                       & \quad + \norm{\ub u r\qq\Nd\Nd_3(\sigma-\sigmao)}_{L^2(S_{u,\ub})} + \norm{\ub u\qq^2\Nd_3\Nd_3(\sigma-\sigmao)}_{L^2(S_{u,\ub})} \\
                       & \quad + \norm{\ub u (r\Nd_3)r\Nd_4(\sigma-\sigmao)}_{L^2(S_{u,\ub})}, \\ \\
    \RR_2[\be](u,\ub) & := \norm{\ub^2r^2\Nd^2\be}_{L^2(S_{u,\ub})} + \norm{\ub^2r^2\Nd\Nd_4\be}_{L^2(S_{u,\ub})} + \norm{\ub^2r^2\Nd_4\Nd_4\be}_{L^2(S_{u,\ub})} \\
               & \quad + \norm{\ub^2r \qq\Nd\Nd_3\be}_{L^2(S_{u,\ub})} + \norm{\ub^2\qq^2\Nd_3\Nd_3\be}_{L^2(S_{u,\ub})} + \norm{\ub^2(r\Nd_3)(r\Nd_4)\be}_{L^2(S_{u,\ub})}, \\ \\
    \RR_2[\al](u,\ub) & :=  \norm{\ub^2r^2\Nd^2\al}_{L^2(S_{u,\ub})} + \norm{\ub^2r^2\Nd\Nd_4\al}_{L^2(S_{u,\ub})} + \norm{\ub^2r^2\Nd\Nd_3\al}_{L^2(S_{u,\ub})} \\
                       & \quad + \norm{\ub^2r^2\Nd_3\Nd_3\al}_{L^2(S_{u,\ub})} + \norm{\ub^2(r\Nd_3)(r\Nd_4)\al}_{L^2(S_{u,\ub})} + \norm{\ub^2(r\Nd_4)^2\al}_{L^2(S_{u,\ub})}. 
  \end{align*}
  We define
  \begin{align*}
    \RRb_{\leq 1}(u,\ub) & := \RRb_{0}(u,\ub) + \RRb_0^{add}(u,\ub) + \RRb_{1}(u,\ub), & \RRb_{\leq 2}(u,\ub) & := \RRb_2(u,\ub) +\RRb_1^{add}(u,\ub) + \RRb_{\leq 1}(u,\ub),\\
    \RR_{\leq 1}(u,\ub) & := \RR_{0}(u,\ub) + \RR_0^{add}(u,\ub) + \RR_{1}(u,\ub), & \RR_{\leq 2}(u,\ub) & := \RR_2(u,\ub) +\RR_1^{add}(u,\ub) + \RR_{\leq 1}(u,\ub).\\
  \end{align*}
  We have the following definitions for integrals of the Bel-Robinson tensors
  \begin{align*}
    \QQb_{1}(u,\ub) & := \int_{S_{u,\ub}}\le(Q(\Lieh_\TE\R)(\KE,\KE,\KE,\elb) + Q(\Lieh_\OOE\R)(\KE,\KE,\TE,\elb)\ri), \\
    \QQ_{1}(u,\ub) & := \int_{S_{u,\ub}}\le(Q(\Lieh_\TE\R)(\KE,\KE,\KE,\el) + Q(\Lieh_\OOE\R)(\KE,\KE,\TE,\el)\ri),\\
    \QQb_{2}(u,\ub) & := \int_{S_{u,\ub}} \big(Q(\Lieh_\OOE^2\R)(\KE,\KE,\TE,\elb) + Q(\Lieh_\OOE\Lieh_\TE\R)(\KE,\KE,\KE,\elb) \\
                    & \quad\quad\quad  + Q(\Lieh_\SE\Lieh_\TE\R)(\KE,\KE,\KE,\elb)\big),\\
    \QQ_2(u,\ub) & := \int_{S_{u,\ub}}\big(Q(\Lieh_\OOE^2\R)(\KE,\KE,\TE,\elb) + Q(\Lieh_\OOE\Lieh_\TE\R)(\KE,\KE,\KE,\elb) \\
                    & \quad\quad\quad  + Q(\Lieh_\SE\Lieh_\TE\R)(\KE,\KE,\KE,\elb)\big),
  \end{align*}
  and
  \begin{align*}
    \QQb_{\leq 2}(u,\ub) & := \QQb_{1}(u,\ub) + \QQb_{2}(u,\ub), & \QQ_{\leq 2}(u,\ub) & := \QQ_1(u,\ub) + \QQ_2(u,\ub).
  \end{align*}
\end{proposition}

\begin{remark}\label{rem:fasterdecayerrorrhoosigmaoo}
  The terms $\rhoo$, $\sigmao$ in the right-hand side of~\eqref{est:RRRRb} are not controlled in terms of the Bel-Robinson tensors.\footnote{They are controlled by integrating the associated Bianchi equation~\eqref{eq:Nd3rho} and by using the null structure~(\ref{eq:Curlze}) respectively. See Sections~\ref{sec:connestCCba} and~\ref{sec:connest}.} However, as it will be shown in the next section, the terms $\norm{\ub(\rhoo,\sigmao)}_{L^2(S_{u,\ub})}$ and $\norm{u(\rhoo,\sigmao)}_{L^2(S_{u,\ub})}$ have sufficient decay to be integrated in $u$ or $\ub$.
\end{remark}

\section{Proof of the $L^2$ bounds~(\ref{est:RRCCbaext}) and~(\ref{est:RRextMM}) on $\protect\CCba\cap\MM^\ext$ and $\MM^\ext$}
In this section, we use the control of Proposition~\ref{prop:curvS} and the energy estimates of Section~\ref{sec:globener} to derive the estimates~\eqref{est:RRCCbaext} and~\eqref{est:RRextMM}.\\

We first have the following control of the curvature fluxes through the hypersurfaces $\CC_u,\Si_t^\ext$ and $\CCba\cap\MM^\ext$.
\begin{lemma}\label{lem:bdedcurvflux}
  We have
  \begin{subequations}\label{est:bdedcurvflux}
    \begin{itemize}
    \item for all $1\leq u \leq \cc\uba$
      \begin{align}\label{est:curvfluxCC}
        \int_{\cc^{-1}u}^\uba \le(\RR_{\leq 2}(u,\ub)\ri)^2 \,\d\ub & \les \varep^2,
      \end{align}
    \item for all $\too \leq t \leq \tast$
      \begin{align}\label{est:curvfluxSitext}
        \int_{2t/(1+\cc^{-1})}^{2t-\uba} \le(\RR_{\leq 2}(u,\ub = 2t-u)\ri)^2 + \le(\RRb_{\leq 2}(u,\ub=2t-u)\ri)^2 \, \d u & \les \varep^2,
      \end{align}
    \item on $\CCba\cap\MM^\ext$
      \begin{align}\label{est:curvfluxCCba}
      \int_{1}^{\cc\uba}\le(\RRb_{\leq 2}(u,\uba)\ri)^2\,\d u & \les \varep^2.
      \end{align}
  \end{itemize}
  \end{subequations}
\end{lemma}
\begin{proof}
  From the bounds~(\ref{est:enerestasscurvest}) through the hypersurface $\CC_u$, the Bootstrap Assumptions~\ref{BA:curvext} on $\rhoo$ and $\sigmao$, and the results of Proposition~\ref{prop:curvS}, we have for all $1\leq u \leq \cc\uba$
  \begin{align*}
    \int_{\cc^{-1}u}^\uba \le(\RR_{\leq 2}(u,\ub)\ri)^2 \,\d\ub & \les \int_{\cc^{-1}u}^\uba \QQ_{\leq 2}(u,\ub) \, \d\ub + \int_{\cc^{-1}u}^\uba (D\varep)^2\norm{\ub(\rhoo,\sigmao)}_{L^2(S_{u,\ub})}^2 \,\d\ub\\
                                                                & \les \varep^2 + (D\varep)^4\int_{\cc^{-1}u}^\uba \ub^4 \ub^{-6}u^{-4}\,\d\ub \\
                                                                & \les \varep^2+(D\varep)^4 \\
                                                                & \les \varep^2,
  \end{align*}
  which proves~\eqref{est:curvfluxCC}. Estimates~\eqref{est:curvfluxSitext} and~\eqref{est:curvfluxCCba} follow similarly from the bounds~(\ref{est:enerestasscurvest}) on the Bel-Robinson tensors through respectively $\Si_t^\ext$ and $\CCba\cap\MM^\ext$. Details are left to the reader.  
\end{proof}

From the result of Lemma~\ref{lem:bdedcurvflux}, and an inspection of the definitions of Section~\ref{sec:normnullcurv}, we have on $\CCba\cap\MM^\ext$
\begin{align*}
  \RR^\ast_{\leq 2} & \simeq \le(\int_1^{\cc\uba}\le(\RRb_{\leq 2}(u,\uba)\ri)^2\,\d u\ri)^{1/2} \les \varep,
\end{align*}
which proves the desired bound~\eqref{est:RRCCbaext}.\\

From an inspection of the norms defined in Section~\ref{sec:normnullcurv} and the results of Lemma~\ref{lem:bdedcurvflux}, we deduce
\begin{align*}
  \le(\RR_{\leq 2,\ga}^\ext\ri)^2 & \les \int_{1}^{\cc\uba}\int_{\cc^{-1}u}^\uba \le(u^{-1-2\ga}\le(\RR_{\leq 2}(u,\ub)\ri)^2 + \ub^{-1-2\ga}\le(\RRb_{\leq 2}(u,\ub)\ri)^2 \ri)\,\d\ub\d u \\
                             & \les \int_{1}^{\cc\uba} u^{-1-2\ga} \int_{\cc^{-1}u}^\uba \le(\RR_{\leq 2}(u,\ub)\ri)^2 \,\d\ub\d u \\
                                  & \quad + \int_{\too}^\tast t^{-1-2\ga}\int_{2t/(1+\cc^{-1})}^{2t-\uba} \le(\RR_{\leq 2}(u,\ub = 2t-u)\ri)^2 + \le(\RRb_{\leq 2}(u,\ub=2t-u)\ri)^2 \, \d u \d t \\
                                  & \les \varep^2 \int_{1}^{\cc\uba} u^{-1-2\ga} \,\d u + \varep^2 \int_{\too}^\tast t^{-1-2\ga}\,\d t \\
                             & \les_\ga \varep^2,
\end{align*}
for all $\ga>0$, and where we used the coarea formulas from Lemma~\ref{lem:coarea} and $t\simeq \ub$ in $\MM^\ext$. This finishes the proof of estimate~\eqref{est:RRextMM}.

\section{Proof of the $L^\infty\HHt$ estimates~(\ref{est:RRCCbadecaynew}) and~(\ref{est:RRextdecay})}\label{sec:decayH12extcurv}

The proof of the estimates~(\ref{est:RRCCbadecaynew}) and~(\ref{est:RRextdecay}) boils down to the following two Klainerman-Sobolev estimates on $\CCba\cap\MM^\ext$ and $\Si_t^\ext$ respectively. Their proof are postponed to Appendix~\ref{app:KlSobH12}.
\begin{lemma}[Klainerman-Sobolev estimates on $\CCba\cap\MM^\ext$]\label{lem:KlSobast}
  For all $S$-tangent tensor $F$ on $\CCba\cap\MM^\ext$, we have the following $L^\infty_{u}\HHt\le(S_{u,\uba}\ri)$ estimates in $\CCba\cap\MM^\ext$
  \begin{align*}
    \norm{r F}_{L^\infty_{1\leq u \leq \cc\uba}\HHt\le(S_{u,\uba}\ri)} & \les \norm{F}_{L^2(\CCba\cap\MM^\ext)} + \norm{r\Nd F}_{L^2(\CCba\cap\MM^\ext)} +\norm{r\Nd_3F}_{L^2(\CCba\cap\MM^\ext)},
  \end{align*}
  and
  \begin{align*}
    \norm{r^{1/2}\qq^{1/2}F}_{L^\infty_{1\leq u\leq \cc\uba}\HHt\le(S_{u,\uba}\ri)} & \les \norm{F}_{L^2(\CCba\cap\MM^\ext)} + \norm{r\Nd F}_{L^2(\CCba\cap\MM^\ext)} + \norm{\qq\Nd_3F}_{L^2(\CCba\cap\MM^\ext)}.
  \end{align*}
  

\end{lemma}


\begin{lemma}[Klainerman-Sobolev estimates on $\Si_t^\ext$]\label{lem:KlSobSitext}
  For all $S$-tangent tensor $F$ we have the following $L^\infty_{u,\ub}\HHt\le(S_{u,\ub}\ri)$ estimates in $\MM^\ext$
  \begin{align*}
    \norm{\ub F}_{L^\infty_{u,\ub}\HHt\le(S_{u,\ub}\ri)} & \les \norm{\uba F}_{L^\infty_{u}\HHt(S_{u,\uba})} + \norm{F}_{L^\infty_tL^2(\Si_t^\ext)} \\
    & \quad + \norm{\ub\Nd F}_{L^\infty_tL^2(\Si_t^\ext)} +\norm{\ub\Nd_3F}_{L^\infty_tL^2(\Si_t^\ext)} +\norm{\ub\Nd_4F}_{L^\infty_tL^2(\Si_t^\ext)},
  \end{align*}
  and
  \begin{align*}
    \norm{\ub^{1/2}u^{1/2}F}_{L^\infty_{u,\ub}\HHt\le(S_{u,\ub}\ri)} & \les \norm{\uba^{1/2}u^{1/2}F}_{L^\infty_u\HHt(S_{u,\uba})} + \norm{F}_{L^\infty_tL^2(\Si_t^\ext)} \\
    & \quad + \norm{\ub\Nd F}_{L^\infty_tL^2(\Si_t^\ext)} + \norm{u\Nd_3F}_{L^\infty_tL^2(\Si^\ext_t)} + \norm{u\Nd_4F}_{L^\infty_tL^2(\Si^\ext_t)}.
  \end{align*}
\end{lemma}

First, we apply the degenerate version of the Klainerman-Sobolev estimates of Lemma~\ref{lem:KlSobast}, with $F$ the following respective tensors
\begin{align*}
  u^2\Ndt^{\leq 1}\alb,~u\ub\Ndt^{\leq 1}\beb,~\ub^2\Ndt^{\leq 1}(\rho-\rhoo),~\ub^2\Ndt^{\leq 1}(\sigma-\sigmao),
\end{align*}
and the non-degenerate version with $F$ respectively
\begin{align*}
  \ub^2\Ndt^{\leq 1}\be,~\ub^2(r\Nd)^{\leq 1}\al,~\ub^2(r\Nd_3)\al,
\end{align*}
where $\Ndt \in\{(r\Nd),(q\Nd_3)\}$. From an inspection of the definitions of Section~\ref{sec:normnullcurv}, we deduce from the bounds~(\ref{est:curvfluxCCba}) that
\begin{align*}
  \mathfrak{R}^\ast_{\leq 1} & \les \RR^\ast_{\leq 2} \les \varep,
\end{align*}
where the norms are restricted to $\CCba\cap\MM^\ext$. This proves~\eqref{est:RRCCbadecaynew}.\\

From an inspection of the definitions of Section~\ref{sec:normnullcurv}, the Klainerman-Sobolev estimates from Lemma~\ref{lem:KlSobSitext}, the coarea formulas from Lemma~\ref{lem:coarea}, the above bound for $\mathfrak{R}^\ast_{\leq 1}$ and the bounds from Lemma~\ref{lem:bdedcurvflux}, we have
\begin{align*}
  \mathfrak{R}^\ext_{\leq 1} & \les \; \mathfrak{R}^{\ast}_{\leq 1} + \sup_{\too \leq t \leq \tast} \le(\int_{2t/(1+\cc^{-1})}^{2t-\uba} \le(\RR_{\leq 2}(u,\ub = 2t-u)\ri)^2 + \le(\RRb_{\leq 2}(u,\ub=2t-u)\ri)^2 \, \d u\ri) \\
  & \les \varep^2. 
\end{align*}
This finishes the proof of~(\ref{est:RRextdecay}).

\section{Proof of Proposition~\ref{prop:curvS}}\label{sec:proofcurvS}

\subsection{Control of $\protect\RRb_{\leq 1}(u,\protect\ub)$ and $\RR_{\leq 1}(u,\protect\ub)$}\label{sec:curvestlin1}
From the definition of $\TE,\KE$, and the decomposition~\cite[p. 150]{Chr.Kla93} of the Bel-Robinson tensors in the null frame $(\elb,\el)$, we have
\begin{align*}
  & \le|u^2\alb\le(\Lieh_\OOE\R\ri)\ri|^2 + \le|u\ub\beb\le(\Lieh_\OOE\R\ri)\ri|^2 + \le|\ub^2\rho\le(\Lieh_\OOE\R\ri)\ri|^2 + \le|\ub^2\sigma\le(\Lieh_\OOE\R\ri)\ri|^2 + \le|\ub^2\be\le(\Lieh_\OOE\R\ri)\ri|^2 \\
  \les & \; Q\le(\Lieh_\OOE\R\ri)\le(\KE,\KE,\TE,\elb\ri), 
\end{align*}
where we denote by $\OOE$ any exterior rotation vectorfield $\OOEi$ for $\ell=1,2,3$.\\

Moreover, we have
\begin{align*}
  & \le|u^2\Liedh_\OOE\alb\ri|^2 + \le|u\ub\Liedh_\OOE \beb\ri|^2 + \le|\ub^2\Liedh_\OOE\rho\ri|^2 + \le|\ub^2\Liedh_\OOE\sigma\ri|^2 + \le|\ub^2\Liedh_\OOE\be\ri|^2\\
  \les & \; \le|u^2\alb\le(\Lieh_\OOE\R\ri)\ri|^2 + \le|u\ub\beb\le(\Lieh_\OOE\R\ri)\ri|^2 + \le|\ub^2\rho\le(\Lieh_\OOE\R\ri)\ri|^2 + \le|\ub^2\sigma\le(\Lieh_\OOE\R\ri)\ri|^2 + \le|\ub^2\be\le(\Lieh_\OOE\R\ri)\ri|^2 \\
  & \quad + \le|\Err_1\ri|^2,
\end{align*}
where $\Err_1$ writes
\begin{align*}
  \Err_1 & := u^2\le(\Liedh_\OOE\alb - \alb\le(\Lieh_\OOE\R\ri)\ri) + u\ub\le(\Liedh_\OOE\beb-\beb\le(\Lieh_\OOE\R\ri)\ri) \\
         & \quad + \ub^2\le(\Liedh_\OOE\rho-\rho\le(\Lieh_\OOE\R\ri)\ri) +\ub^2\le(\Liedh_\OOE\sigma-\sigma\le(\Lieh_\OOE\R\ri)\ri) \\
         & \quad +\ub^2\le(\Liedh_\OOE\be-\be\le(\Lieh_\OOE\R\ri)\ri).
\end{align*}
The treatment of the error term $\Err_1$ is postponed to Section~\eqref{sec:curvesterr}, where it is proved that
\begin{align*}
  \int_{S_{u,\ub}}|\Err_1|^2 & \les (D\varep)^2\le(\RRb_{\leq 1}(u,\ub)\ri)^2 + (D\varep)^2 \int_{S_{u,\ub}}\le|\ub(\rhoo,\sigmao)\ri|^2.
\end{align*}
Thus, combining these estimates and integrating on $S_{u,\ub}$, using the definition of $\QQb_1$, we obtain
\begin{align}
  \label{est:LiehOOEnulld1}
  \begin{aligned}
    & \int_{S_{u,\ub}}\le|u^2\Liedh_\OOE\alb\ri|^2 + \le|u\ub\Liedh_\OOE \beb\ri|^2 + \le|\ub^2\Liedh_\OOE\rho\ri|^2 + \le|\ub^2\Liedh_\OOE\sigma\ri|^2 + \le|\ub^2\Liedh_\OOE\be\ri|^2 \les \FF(u,\ub),
  \end{aligned}
\end{align}
where
\begin{align*}
  \FF(u,\ub) & := \QQb_{1}(u,\ub) +  (D\varep)^2\le(\RRb_{\leq 1}(u,\ub)\ri)^2 + (D\varep)^2\int_{S_{u,\ub}}\le|\ub(\rhoo,\sigmao)\ri|^2.
\end{align*}




Using the estimates of the mild Bootstrap Assumptions~\ref{BA:mildOOE} and the Poincar\'e estimates from Lemma~\ref{lem:ell} we deduce from~\eqref{est:LiehOOEnulld1} the following estimate
\begin{align}
  \label{est:NdRRRb1}
  \begin{aligned}
    & \int_{S_{u,\ub}}\bigg(\le|u^2(r\Nd)^{\leq 1}\alb\ri|^2 + \le|u\ub(r\Nd)^{\leq 1}\beb\ri|^2 + \le|\ub^2(r\Nd)^{\leq 1}(\rho-\rhoo)\ri|^2 \\
    & + \le|\ub^2(r\Nd)^{\leq 1}(\sigma-\sigmao)\ri|^2 + \le|\ub^2(r\Nd)^{\leq 1}\be\ri|^2 \bigg) \\
    & \int_{S_{u,\ub}}\le|u^2\Liedh_\OOE\alb\ri|^2 + \le|u\ub\Liedh_\OOE \beb\ri|^2 + \le|\ub^2\Liedh_\OOE\rho\ri|^2 + \le|\ub^2\Liedh_\OOE\sigma\ri|^2 + \le|\ub^2\Liedh_\OOE\be\ri|^2  \\
    \les & \; \FF(u,\ub).
  \end{aligned}
\end{align}


Using the Bianchi identities~\eqref{eq:bianchi}, we deduce from~\eqref{est:NdRRRb1} the following control for the $\Nd_3,\Nd_4$ derivatives
\begin{align}\label{est:bianchiRRb10}
  \begin{aligned}
    & \int_{S_{u,\ub}} \bigg(\le|u^2(r\Nd_4)\alb\ri|^2 + \le|u\ub(u\Nd_3)\beb\ri|^2 + \le|u\ub (r\Nd_4)\beb\ri|^2 \\
    & + \le|\ub^2u\le(\Nd_3\rho + \frac{3}{2}\trchib\rho\ri)\ri|^2 + \le|\ub^2u\le(\Nd_3\sigma +\frac{3}{2}\trchib\sigma\ri)\ri|^2 \\
    & + \le|\ub^2\ub\le(\Nd_4\rho + \frac{3}{2}\trchi\rho\ri)\ri|^2 + \le|\ub^2\ub\le(\Nd_4\sigma +\frac{3}{2}\trchi\sigma\ri)\ri|^2 + \le|\ub^2\ub\Nd_3\be\ri|^2\bigg)  \\
    \les & \; \FF(u,\ub) + \int_{S_{u,\ub}}|\Err_2|^2 \\
    \les & \; \FF(u,\ub),
  \end{aligned}
\end{align}
where the error term $\Err_2$ writes
\begin{align*}
  \Err_2 & := u^2r\le(\Nd_4\alb +r^{-1}\alb + \Nd\otimesh\beb\ri) + u^2\ub\le(\Nd_3\beb-4r^{-1}\beb+\Divd\alb\ri) \\
         & \quad + u\ub r\le(\Nd_4\beb+2r^{-1}\beb+\Nd\rho-\dual\Nd\sigma\ri) + \ub^2u\le(\Nd_3\rho+\frac{3}{2}\trchib\rho+\Divd\beb\ri) \\
         & \quad +\ub^2u\le(\Nd_3\sigma+\frac{3}{2}\trchib\sigma+\Curld\beb\ri) +\ub^2\ub\le(\Nd_4\rho+\frac{3}{2}\trchi\rho-\Divd\be\ri) \\
         & \quad +\ub^2\ub\le(\Nd_4\sigma+\frac{3}{2}\sigma+\Curld\be\ri) +\ub^2\ub\le(\Nd_3\be-2r^{-1}\be-\Nd\rho-\dual\Nd\sigma\ri).
\end{align*}
The control for the error term $\Err_2$ by $\FF$ is postponed to Section~\ref{sec:curvesterr}. From~\eqref{est:bianchiRRb10} and the control~\eqref{est:NdRRRb1} for $\rho-\rhoo, \sigma-\sigmao$ and the mild bounds $|r\trchi|+|r\trchib| \les 1$, we deduce
\begin{align}\label{est:bianchiRRb1}
  \begin{aligned}
    & \int_{S_{u,\ub}} \bigg(\le|u^2(r\Nd_4)\alb\ri|^2 + \le|u\ub(u\Nd_3)\beb\ri|^2 + \le|u\ub (r\Nd_4)\beb\ri|^2 \\
    & + \le|\ub^2u\Nd_3(\rho-\rhoo)\ri|^2 + \le|\ub^2u\Nd_3(\sigma-\sigmao)\ri|^2 \\
    & + \le|\ub^2\ub\Nd_4(\rho-\rhoo)\ri|^2 + \le|\ub^2\ub\Nd_4(\sigma-\sigmao)\ri|^2 + \le|\ub^2\ub\Nd_3\be\ri|^2\bigg)  \\
    \les & \; \FF(u,\ub).
  \end{aligned}
\end{align}
Using the other contracted Bel-Robinson tensor of $\QQb_1$, we have
\begin{align*}
  \int_{S_{u,\ub}}\le|u^2u\Nd_{\TE}\alb\ri|^2 + \le|\ub^2\ub\Nd_{\TE}\be\ri|^2 & \les \int_{S_{u,\ub}}\le(Q(\Lieh_\TE\R)(\KE,\KE,\KE,\elb)\ri) + \int_{S_{u,\ub}}|\Err_3|^2 \\
  & \les \FF(u,\ub),
\end{align*}
where the error term $\Err_3$ is defined by
\begin{align*}
  \Err_3 & := u^2u\le(\Nd_\TE\alb-\alb\le(\Lieh_\TE\R\ri)\ri) + \ub^2\ub\le(\Nd_\TE\be-\be\le(\Lieh_\TE\R\ri)\ri).
\end{align*}
The control of the error term $\Err_3$ by $\FF$ is postponed to Section~\ref{sec:curvesterr}. Combining with~\eqref{est:bianchiRRb1}, we obtain
\begin{align}\label{est:addRRb1}
  \int_{S_{u,\ub}}\le|u^2u\Nd_3\alb\ri|^2 + \int_{S_{u,\ub}}\le|\ub^3\Nd_4\be\ri|^2 & \les \FF(u,\ub).
\end{align}

Using the estimates~\eqref{est:NdRRRb1} and~\eqref{est:addRRb1} respectively for $\be$ and $\Nd_4\be$, and the elliptic estimates on $2$-spheres from Lemma~\ref{lem:ell}, we obtain
\begin{align}\label{est:correcalR0add}
  \begin{aligned}
  \int_{S_{u,\ub}}\le|\ub^2(r\Nd)^{\leq 1}\al\ri|^2 & \les \int_{S_{u,\ub}}\ub^4\le|(r\Nd_4)^{\leq 1}\be\ri|^2 + \int_{S_{u,\ub}}|\Err_4|^2 \\
  & \les \FF(u,\ub),
  \end{aligned}
\end{align}
where $\Err_4$ is defined by
\begin{align*}
  \Err_4 & := \ub^2r\le(\Nd_4\be+2\trchi\be-\Divd\al\ri).
\end{align*}
The control of the error term $\Err_4$ is obtained in Section~\ref{sec:curvesterr}.
\begin{remark}
  Estimate~\eqref{est:correcalR0add} provides the desired control of the additional $0$-order norms $\RRb_0^{add}[\al]$.
\end{remark}

We further obtain
\begin{align*}
  \int_{S_{u,\ub}}\le|\ub^2r\Nd_3\al\ri|^2 & \les \FF(u,\ub) + \int_{S_{u,\ub}}|\Err_5|^2 \les \FF(u,\ub),
\end{align*}
where the term $\Err_5$ is defined by
\begin{align*}
  \Err_5 & := \ub^2r\le(\Nd_3\al-r^{-1}\al-\Nd\otimesh\be\ri),
\end{align*}
and is controlled in Section~\ref{sec:curvesterr}.\\

Summarising the bounds obtained in this section, we have proved
\begin{align*}
  \le(\RRb_{\leq 1}(u,\ub)\ri)^2 & \les \FF(u,\ub) = \QQb_1(u,\ub) + (D\varep)^2\le(\RRb_{\leq 1}(u,\ub)\ri)^2 + (D\varep)^2\norm{\ub(\rhoo,\sigmao)}^2_{L^2(S_{u,\ub})},
\end{align*}
which after an absorption argument, concludes the desired control of $\RRb_{\leq 1}(u,\ub)$. The control of $\RR_{\leq 1}(u,\ub)$ is obtained along the same lines and is left to the reader.

\subsection{Error term estimates from Section~\ref{sec:curvestlin1}}\label{sec:curvesterr}
In this section, we show that, under the Bootstrap Assumptions~\ref{BA:connext} on the null connection and rotation coefficients in $\MM^\ext$, we have for the error terms $\Err_1,\cdots, \Err_5$ from Section~\ref{sec:curvestlin1}
\begin{align*}
  \norm{\Err}_{L^2(S_{u,\ub})} &  \les (D\varep) \RRb_{\leq 1}(u,\ub) + (D\varep)\norm{\ub(\rhoo,\sigmao)}_{L^2(S_{u,\ub})}.
\end{align*}
The error terms for the dual estimates for $\RR$ are controlled similarly, and their estimate is left to the reader.

\paragraph{Estimates for $\Err_1$ and $\Err_3$}
We first recall from~\cite[pp. 152--153]{Chr.Kla93}, relations~\eqref{eq:Riccirel} and~(\ref{eq:relOOEnull}) that for the vectorfields $\TE$ and $\OOE$, we have
\begin{align*}
  [\TE,\elb] & = -2\ze_a\ea - \omb\el, \\
  [\TE,\el] & = 2\ze_a\ea + \omb\el, \\
  [\TE,\ea] & = \Pi\le([T,\ea]\ri) + \half \xib_a\el,
\end{align*}
and
\begin {align*}
  [\OOE,\elb] & = -(\ze\cdot\OOE)\elb -(\xib\cdot\OOE)\el  -Y,\\
  [\OOE,\el] & = 0,\\
  [\OOE,\ea] & = \Pi\le([\OOE,\ea]\ri), 
\end{align*}
where $\Pi([T,\ea])$ and $\Pi([\OOE,\ea])$ are the projection on $S$ of $[T,\ea]$ and $[\OOE,\ea]$ respectively, and play no role in the argument below.\\

For a general vectorfield $X$, we have the following definition (see~\cite[pp. 152--153]{Chr.Kla93}) for the coefficients $P,\Pb, Q, \Qb, M, \Mb, N, \Nb$ by
\begin{align*}
  [X,\elb] & = {^{(X)}\Pb}_a\ea + {^{(X)}\Mb}\elb + {^{(X)}\Nb}\el, \\
  [X,\el] & = {^{(X)}P}_a\ea + {^{(X)}N}\elb+{^{(X)}M}\el,\\
  [X,\ea] & = \Pi([X,\ea]) + \half {^{(X)}Q}_a\elb + \half {^{(X)}\Qb}_a\el.
\end{align*}

In terms of the above notation, we thus have
\begin{align}\label{eq:ChrLieT}\begin{aligned}
  ^{(\TE)}\Pb & = -2\ze, & ^{(\TE)}P & = 2\ze,\\
  ^{(\TE)}\Mb & = 0, & ^{(\TE)}M & = \omb,\\
  ^{(\TE)}\Nb & = -\omb, & ^{(\TE)}N & = 0,\\
  ^{(\TE)}\Qb & = \xib, & ^{(\TE)}Q & = 0,
  \end{aligned}
\end{align}
and
\begin{align}\label{eq:ChrLieOOE}
  \begin{aligned}
  ^{(\OOE)}\Pb & = -Y, & ^{(\OOE)}P & = 0,\\
  ^{(\OOE)}\Mb & = -\ze\cdot\OOE, & ^{(\OOE)}M & = 0,\\
  ^{(\OOE)}\Nb & = -\xib\cdot\OOE, & ^{(\OOE)}N & = 0,\\
  ^{(\OOE)}\Qb & = 0, & ^{(\OOE)}Q & = 0.
  \end{aligned}
\end{align}

From direct computation or from relations~(\ref{eq:relpiOOEnull}), we also have
\begin{align}\label{eq:trTpi}
  \begin{aligned}
    \tr^{(\TE)}\pi & = -2\omb + \le(\trchi+\trchib\ri),
  \end{aligned}
\end{align}
and
\begin{align}\label{eq:trOOEpi}
  \begin{aligned}
    \tr^{(\OOE)}\pi & = 2\ze\cdot\OOE + \tr \Hrot
  \end{aligned}
\end{align}
and that
\begin{align}\label{eq:pihT}
  \begin{aligned}
    ^{(T)}\pih_{ab} & = \half\omb +\quar(\trchi+\trchib)\gd_{ab} + \chih_{ab} + \chibh_{ab},
  \end{aligned}
\end{align}
\begin{align}\label{eq:pihOOE}
  \begin{aligned}
    ^{(\OOE)}\pih_{ab} & = \Hrot_{ab} -\quar (2\ze\cdot\OOE + \tr \Hrot) \gd_{ab}.
  \end{aligned}
\end{align}

From the formulas of~\cite[pp. 152-153]{Chr.Kla93}, we have
\begin{align}\label{eq:commnullLie}
  \nulld\le(\Lieh_X\R\ri) - \Liedh_X\nulld(\R) & = A \cdot R,
\end{align}
where $A$ are Lie coefficients $P,\Pb, Q, \Qb, M, \Mb, \tr\pi, \pih_{ab}$ from~\eqref{eq:ChrLieT},~\eqref{eq:trTpi} and~\eqref{eq:pihT} for $X=\TE$ and~\eqref{eq:ChrLieOOE},~\eqref{eq:trOOEpi} and~\eqref{eq:pihOOE} for $X=\OOE$ and where $R \in \{\alb,\cdots,\al\}$ with signature $s(R) = s(\nulld) \pm 1 $.\\

Recalling the definition of $\Err_1$ from Section~\ref{sec:curvestlin1}, we want to obtain
\begin{align}\label{obj:LienulldOOE}
  \norm{u^p\ub^q\le(\nulld\le(\Lieh_\OOE\R\ri)-\Liedh_\OOE\nulld(\R)\ri)}_{L^2(S_{u,\ub})} & \les (D\varep) \RRb_{\leq 1}(u,\ub) + (D\varep)\norm{\ub(\rhoo,\sigmao)}_{L^2(S_{u,\ub})},
\end{align}
where $u^p,\ub^q$ are the appropriate powers of $u$ and $\ub$ given in Section~\ref{sec:curvestlin1}.\\

From the expression of the Lie coefficients $A$ of $\OOE$ given in~\eqref{eq:ChrLieOOE},~\eqref{eq:trOOEpi} and~\eqref{eq:pihOOE} and the decay estimates of the Bootstrap Assumptions~\ref{BA:connext}, we have
\begin{align*}
  \norm{\ub A}_{L^\infty(\MM^\ext)} & \les D\varep.
\end{align*}

From formula~\eqref{eq:commnullLie}, we have
\begin{align*}
  \norm{u^p\ub^q \le(\nulld\le(\Lieh_\OOE\R\ri)-\Liedh_\OOE\nulld(\R)\ri)}_{L^2(S_{u,\ub})} & \les \norm{\ub A}_{L^\infty}\norm{u^p\ub^{q-1}R}_{L^2(S_{u,\ub})},
\end{align*}
where $R\in\{\alb,\cdots,\al\}$ with $s(R) = s(\nulld) \pm 1$.\\

We have
\begin{align*}
  \norm{u^p\ub^{q-1}R}_{L^2(S_{u,\ub})} & \les \RRb_{\leq 1}(u,\ub) + \norm{\ub(\rhoo,\sigmao)}_{L^2(S_{u,\ub})},
\end{align*}
and the control of~\eqref{obj:LienulldOOE} and of $\Err_1$ follows.\\

From the expression of the Lie coefficients $A$ of $\TE$ given in~\eqref{eq:ChrLieT},~\eqref{eq:trTpi} and~\eqref{eq:pihT} and the Bootstrap Assumptions~\ref{BA:connext}, we have
\begin{align*}
  \norm{\ub u A}_{L^\infty(\MM^\ext)} & \les D\varep,\\
  \norm{\ub^2\le(P,Q\ri)}_{L^\infty(\MM^\ext)} & \les D\varep.
\end{align*}
and we therefore deduce from an inspection of~\cite[pp. 152--153]{Chr.Kla93}
\begin{align*}
  \norm{u^p\ub^q \le(\nulld\le(\Lieh_\TE\R\ri)-\Liedh_\TE\nulld(\R)\ri)}_{L^2(S_{u,\ub})} & \les (D\varep) \le(\RRb_{\leq 1}(u,\ub) + \norm{\ub(\rhoo,\sigmao)}_{L^2(S_{u,\ub})}\ri),
\end{align*}
where $p,q$ are the adequate powers given in the definition of $\Err_3$ in Section~\ref{sec:curvestlin1}.\\

Moreover, we have schematically for $R\in\le\{\alb,\beb,\be\ri\}$
\begin{align*}
  \Liedh_\TE R - \Nd_\TE R & = \le(\chi+\chib\ri)\cdot R,
\end{align*}
and we therefore deduce, using the Bootstrap Assumptions~\ref{BA:connext}, that
\begin{align*}
  \norm{u^p\ub^q\le(\Liedh_\TE R - \Nd_\TE R\ri)}_{L^2(S_{u,\ub})} & \les \norm{\ub\le(\chi+\chib\ri)}_{L^\infty}\norm{u^p\ub^{q-1}R}_{L^2(S_{u,\ub})} \\
  & \les (D\varep)\RRb_{\leq 1}(u,\ub),
\end{align*}
which concludes the control of $\Err_3$.

\paragraph{Estimates for $\Err_2, \Err_4$ and $\Err_5$}
We shall only treat the term $\Err_5$, for the error terms $\Err_2,\Err_4$ will follow along the same lines.\\

From Bianchi equation~\eqref{eq:Nd3al} and the relations of Lemma~\ref{lem:relgeodnull}, we have
\begin{align*}
  \Err_5 & = \ub^3\le(2\omb\al -3\le(\chih\rho+\dual\chih\sigma\ri) + 5\ze\otimes\be\ri).
\end{align*}
Thus, using the Bootstrap Assumptions~\ref{BA:connext}
\begin{align*}
  \norm{\Err_5}_{L^2(S_{u,\ub})} \les & \; (D\varep)\le(\norm{\ub^2\al}_{L^2(S_{u,\ub})} + \norm{\ub^2\be}_{L^2(S_{u,\ub})} + \norm{\ub^2(\rho-\rhoo,\sigma-\sigmao)}_{L^2(S_{u,\ub})}\ri)  \\
  & + (D\varep)\norm{\ub(\rhoo,\sigmao)}_{L^2(S_{u,\ub})} \\
  \les & \; (D\varep) \le(\RRb_{\leq 1}(u,\ub) + \norm{\ub(\rhoo,\sigmao)}_{L^2(S_{u,\ub})}\ri),
\end{align*}
as desired. This finishes the control of the error terms $\Err_2, \Err_4$ and $\Err_5$.




\subsection{Control of $\protect\RRb_{\leq 2}(u,\protect\ub)$ and $\protect\RR_{\leq 2}(u,\protect\ub)$}\label{sec:curvestlin2}
In this section, we argue as in Section~\ref{sec:curvestlin1} to obtain the control of $\RRb_{\leq 2}(u,\ub)$ and $\RR_{\leq 2}(u,\ub)$. The error terms are dealt with arguing as in Section~\ref{sec:curvesterr}, using $L^\infty(\MM^\ext)$ for the connection and rotation coefficients and $L^\infty_{u,\ub}L^4(S_{u,\ub})$ estimates for derivatives of the connection and rotation coefficients. This treatment of the error terms is left to the reader.
\begin{remark}
  From Bianchi equations~\eqref{eq:Nd3rho},~\eqref{eq:Nd3sigma}, \eqref{eq:Nd4rho} and~\eqref{eq:Nd4sigma} for $\Nd_3,\Nd_4$ derivatives of $\rho,\sigma$, and the Bootstrap Assumptions~\ref{BA:connext} for the null connection coefficients, we have
  \begin{align*}
    \norm{\ub^2\le(\Nd_3,\Nd_4\ri)\le(\rhoo,\sigmao\ri)}_{L^2(S_{u,\ub})} & \les (D\varep)\RRb_{\leq 2}(u,\ub) + \norm{\ub(\rhoo,\sigmao)}_{L^2(S_{u,\ub})}.
  \end{align*}
  Thus, derivatives of $\rhoo,\sigmao$ shall in the sequel directly be replaced using the above estimate. 
\end{remark}
\begin{remark}
  Using the mild Bootstrap Assumptions~\ref{BA:mildOOE} and the Bootstrap Assumptions~\ref{BA:connext} for the rotation vectorfields $\OOE$, one can obtain the following higher order tangential derivatives estimates
  \begin{align}\label{est:NdLied2}
    \int_{S_{u,\ub}}\le|(r\Nd)^2F\ri|^2 & \les \sum_{\ell,\ell'=1}^3\int_{S_{u,\ub}}\le|\Lieh_\OOEi\Lieh_{^{(\ell')}\OOE} F\ri|^2 + \int_{S_{u,\ub}}\le|(r\Nd)^{\leq 1}F\ri|^2,
  \end{align}
  for all $S$-tangent tensor $F$ and for all $2$-sphere $S_{u,\ub}\subset\MM^\ext$.
\end{remark}

Arguing as in Sections~\ref{sec:curvestlin1} and~\ref{sec:curvesterr}, using~\eqref{est:NdLied2} we have
\begin{align*}
  \begin{aligned}
    & \int_{S_{u,\ub}}\bigg(\le|u^2(r\Nd)^{\leq 2}\alb\ri|^2 + \le|u\ub(r\Nd)^{\leq 2}\beb\ri|^2 + \le|\ub^2(r\Nd)^{\leq 2}(\rho-\rhoo)\ri|^2 
    + \le|\ub^2(r\Nd)^{\leq 2}(\sigma-\sigmao)\ri|^2 + \le|\ub^2(r\Nd)^{\leq 2}\be\ri|^2  \bigg)\\
  \les & \; \int_{S_{u,\ub}} Q\le(\Lieh_\OOE\Lieh_\OOE\R\ri)(\KE,\KE,\TE,\elb) + (D\varep)^2\le(\RRb_{\leq 2}(u,\ub)\ri)^2 + (D\varep)^2\norm{\ub(\rhoo,\sigmao)}_{L^2(S_{u,\ub})}^2,
  \end{aligned}
\end{align*}
and
\begin{align*}
  \begin{aligned}
    & \int_{S_{u,\ub}} \bigg(\le|u^3(r\Nd)\Nd_\TE\alb\ri|^2 + \le|\ub^3(r\Nd)\Nd_{\TE}\be\ri|^2 \bigg) \\
    \les & \; \int_{S_{u,\ub}}Q\le(\Lieh_\OOE\Lieh_\TE\R\ri)(\KE,\KE,\KE,\elb) + (D\varep)^2\le(\RRb_{\leq 2}(u,\ub)\ri)^2 + (D\varep)^2\norm{\ub(\rhoo,\sigmao)}_{L^2(S_{u,\ub})}^2. 
  \end{aligned}
\end{align*}
Using these two estimates with Bianchi equations~\eqref{eq:bianchi} commuted with $(r\Nd)$ as in Section~\ref{sec:curvestlin1}, we deduce
\begin{align*}
  & \int_{S_{u,\ub}} \bigg(\le|u^2(r\Nd)(u\Nd_3)\alb\ri|^2 + \le|u^2(r\Nd)(r\Nd_4)\alb\ri|^2 + \le|u\ub(r\Nd)(u\Nd_3)\beb\ri|^2 + \le|u\ub (r\Nd)(r\Nd_4)\beb\ri|^2 \\
  & + \le|\ub^2u(r\Nd)\Nd_3(\rho-\rhoo)\ri|^2 + \le|\ub^2u(r\Nd)\Nd_3(\sigma-\sigmao)\ri|^2+ \le|\ub^2\ub(r\Nd)\Nd_4(\rho-\rhoo)\ri|^2 \\
  & + \le|\ub^2\ub(r\Nd)\Nd_4(\sigma-\sigmao)\ri|^2 + \le|\ub^2\ub(r\Nd)\Nd_3\be\ri|^2 + \le|\ub^2(r\Nd)(\ub\Nd_4)\be\ri|^2\bigg)  \\
  \les & \; \QQb_{\leq 2}(u,\ub) + (D\varep)^2\le(\RRb_{\leq 2}(u,\ub)\ri)^2 + (D\varep)^2\norm{\ub(\rhoo,\sigmao)}_{L^2(S_{u,\ub})}^2.
\end{align*}
Taking $\Nd_3,\Nd_4$ derivative in Bianchi equations~\eqref{eq:bianchi}, we deduce from the above estimates
\begin{align*}
  & \le(\RRb_{\leq 2}(u,\ub)\ri)^2 - \int_{S_{u,\ub}} \le|u^2(u\Nd_3)^2\alb\ri|^2 - \int_{S_{u,\ub}}\le|\ub^2(\ub\Nd_4)^2\be\ri|^2 - \int_{S_{u,\ub}}\le|\ub^4(\Nd,\Nd_3)(\Nd,\Nd_3,\Nd_4)\al\ri|^2 \\
  \les & \; \QQb_{\leq 2}(u,\ub) + (D\varep)^2\le(\RRb_{\leq 2}(u,\ub)\ri)^2 + (D\varep)^2\norm{\ub(\rhoo,\sigmao)}_{L^2(S_{u,\ub})}^2.
\end{align*}
Using the last Bel-Robinson tensor of $\QQb_{\leq 2}$ and the above estimate, we have
\begin{align*}
  & \int_{S_{u,\ub}}\bigg(\le|u^4\Nd_3^2\alb\ri|^2 + \le|\ub^4\Nd_4^2\be\ri|^2\bigg) \\
  \les & \int_{S_{u,\ub}}Q\le(\Lieh_\SE\Lieh_\TE\R\ri)(\KE,\KE,\KE,\elb) + (D\varep)^2\le(\RRb_{\leq 2}(u,\ub)\ri)^2 + (D\varep)^2\norm{\ub(\rhoo,\sigmao)}_{L^2(S_{u,\ub})}^2.
\end{align*}
Using (commuted) Bianchi equations~\eqref{eq:Nd3al}~\eqref{eq:Nd4be} for $\Nd_3\al$ and $\Divd\al$, the elliptic estimates of Lemma~\ref{lem:ell} and the above estimates for (derivatives of) $\be$ we deduce
\begin{align*}
  \int_{S_{u,\ub}}\le|\ub^4(\Nd,\Nd_3)(\Nd,\Nd_3,\Nd_4)\al\ri|^2 & \les \int_{S_{u,\ub}}\le|\ub^2(\ub\Nd_4,\ub\Nd)^{\leq 1}(\ub\Nd,\ub\Nd_3,\ub\Nd_4)^{\leq 1}\be\ri|^2 \\
  & \les \QQb_{\leq 2}(u,\ub) + (D\varep)^2\le(\RRb_{\leq 2}(u,\ub)\ri)^2 + (D\varep)^2\norm{\ub(\rhoo,\sigmao)}_{L^2(S_{u,\ub})}^2.
\end{align*}
Using an absorption argument, this concludes the desired control of $\RRb_{\leq 2}(u,\ub)$. The control of $\RR_{\leq 2}(u,\ub)$ follows along the same lines and this finishes the proof of Proposition~\ref{prop:curvS}.

\chapter{Maximal curvature estimates in $\MM^\intr_\bott$}\label{sec:planehypcurvest}
In this section, we prove the following proposition.
\begin{proposition}\label{prop:planehypcurvestSTAB}
  Recall that from Proposition~\ref{prop:enerestSTAB}, we have in $^{(\cc)}\MM^\intr_\bott$
  \begin{subequations}\label{est:intQint}
    \begin{align}
      & \int_{^{(\cc)}\Si_t} Q(\Lieh_\TI \R)(\KI,\KI,\KI,\Tf) \les \varep^2, \label{est:intQ1intSit}\\
      & \int_{^{(\cc)}\Si_t}Q(\Lieh_\SI \Lieh_\TI\R)(\KI,\KI,\KI,\Tf) \les \varep^2,\label{est:intQ2intSit}
    \end{align}
  \end{subequations}
  for a fixed transition parameter $\cc$ and for all $^{(\cc)}\too \leq t \leq {^{(\cc)}\tast}$. Recall that from Proposition~\ref{prop:curvestSTAB}, we have on $^{(\cc)}\TT$
    \begin{align}\label{est:MMbotbdy}
      \norm{t^3\Ndt^{\leq 1}R}_{L^\infty_t\HHt(\pr\Si_t)} & \les \varep,
    \end{align}
    for $\Ndt\in\le\{(t\Nd,t\Nd_3,t\Nd_4)\ri\}$ and $R\in\le\{\al,\be,\rho,\sigma,\beb,\alb\ri\}$.\\
    
    Under the Bootstrap Assumptions, and the estimates~\eqref{est:intQint} and~\eqref{est:MMbotbdy}, we have the following bound in $^{(\cc)}\MM^\intr_\bott$ for $\varep >0$ sufficiently small
    \begin{align}\label{est:RRintbott}
      \begin{aligned}
        \RR^\intr_{\leq 2} & \les \varep,
      \end{aligned}
    \end{align}
    where we refer the reader to Section~\ref{sec:normscurvint} for definitions.
\end{proposition}
\begin{remark}
  Using Sobolev estimates (see Lemma~\ref{lem:SobSitHHrr}), one directly deduces from~\eqref{est:RRintbott} that
  \begin{align*}
    \RRfb^\intr_{\leq 1} & \les \varep
  \end{align*}
  in $^{(\cc)}\MM^\intr_\bott$. See Section~\ref{sec:normscurvint} for definitions.
\end{remark}

\section{Preliminary results}\label{sec:prelimintcurvest}
We have the following Sobolev embeddings on the maximal hypersurfaces $\Si_t$ of $\MM^\intr_\bott$.
\begin{lemma}[Sobolev estimates on $\Si_t$]\label{lem:SobSitHHrr}
  Under the mild Bootstrap Assumptions~\ref{BA:mildcoordsSit}, we have for all $\too \leq t \leq \tast$ and for all $\Si_t$-tangent tensor $F$
  \begin{align*}
    \norm{F}_{L^6(\Si_t)} & \les \norm{\nab F}_{L^2(\Si_t)}^{1/2}\norm{F}_{L^2(\Si_t)}^{1/2} + t^{-1}\norm{F}_{L^2(\Si_t)},\\
    \norm{F}_{L^\infty(\Si_t)} & \les \norm{\nab^2F}_{L^2(\Si_t)}^{3/4}\norm{F}_{L^2(\Si_t)}^{1/4} + t^{-3/2}\norm{F}_{L^2(\Si_t)}.
  \end{align*}
\end{lemma}
\begin{proof}
  A rescaling in $t$ of the estimates of~\cite[Section 3]{Czi.Gra19a} gives the result. 
\end{proof}

We have the following elliptic estimates for div-curl systems on $\Si_t$.
\begin{lemma}[Elliptic estimates for div-curl systems on $\Si_t$]\label{lem:ellHodge}
  Assume that the mild Bootstrap Assumptions~\ref{BA:mildcoordsSit} are satisfied and that the following bounds hold
  \begin{align*}
    \norm{t^{5/2}\le(\trth -\frac{2}{r}\ri)}_{L^\infty(\TT)} & \les D\varep,\\
    \norm{t^{5/2}\thh}_{L^\infty(\TT)} & \les D\varep.
  \end{align*}
  Then, for $\varep>0$ sufficiently small, we have for all symmetric traceless $\Si_t$-tangent $2$-tensor $F$ 
  \begin{align*}
    \int_{\Si_t}\le|(t\nab)^{\leq 1} F\ri|^2 & \les \int_{\Si_t}\le(t^2|\Div F|^2+ t^2|\Curl F|^2\ri) + t^2\le|\int_{\pr\Si_t}\Fslash_\Nf\cdot\Nd\le(F_{\Nf\Nf}\ri)\ri|,
  \end{align*}
  where $\Fslash_\Nf$ is the $\pr\Si_t$-tangent $1$-tensor defined by $\Fslash_{\Nf a} := F_{\Nf a}$.
\end{lemma}
\begin{proof}
  The result of the lemma is a straight-forward generalisation of the estimates obtained for $k$ in~\cite[Section 4.7]{Czi.Gra19a} and a $t$-rescaling (see in particular estimate (4.32) in that paper). Details are left to the reader. 
\end{proof}

We have the following variation of Lemma~\ref{lem:ellHodge}, which will be used in Section~\ref{sec:controlRRintr2} to control $\Lieh_\Tf E$ and $\Lieh_\Tf H$ (see also~\cite[pp. 104-105]{Chr.Kla93} for an analogous result in the case without boundary).
\begin{lemma}[Elliptic estimates for modified div-curl systems on $\Si_t$]\label{lem:ellHodgevar}
  Assume that the hypothesis of Lemma~\ref{lem:ellHodge} hold. There exists a sufficiently small constant $\wp>0$, such that if $X$ is a $\Si_t$-tangent vectorfield with $|X|\leq \wp$ and $|(t\nab) X| \les 1$ on $\Si_t$, then, for all symmetric traceless $\Si_t$-tangent $2$-tensors $F,G$, we have
  \begin{align*}
    \int_{\Si_t}\le|(t\nab)F\ri|^2 + \le|(t\nab)G\ri|^2 & \les \int_{\Si_t}\big(t^2|\Div F|^2 + t^2|\Div G|^2\big) \\
                                                        & \quad + \int_{\Si_t} t^2\le(\le|\Curl F + \Lieh_XG\ri|^2 + \le|\Curl G - \Lieh_XF\ri|^2\ri) \\
                                                        & \quad + \int_{\Si_t} \le(|F|^2 + |G|^2\ri)\\
         & \quad + t^2\le|\int_{\pr\Si_t}\Fslash_N\cdot\Nd\le(F_{NN}\ri)\ri| + t^2\le|\int_{\pr\Si_t}\Gslash_N\cdot\Nd\le(G_{NN}\ri)\ri|
  \end{align*}
\end{lemma}
\begin{proof}
  The result follows from a rescaling of $\Si_t$ to a disk of radius $1$, a generalisation of the elliptic estimates obtained for $k$ in~\cite[Section 4.7]{Czi.Gra19a} and an absorption argument, provided that $|X|\leq \wp$ is sufficiently small. Details are left to the reader. 
\end{proof}

We have the following higher order elliptic estimates for div-curl systems, which is used to control $\nab^2E$ and $\nab^2H$ in Section~\ref{sec:controlRRintr2}.
\begin{lemma}[Higher order elliptic estimates on $\Si_t$]\label{lem:ellSitHrr}
  Under the Bootstrap Assumptions~\ref{BA:mildcoordsSit}, we have for all traceless symmetric $\Si_t$-tangent 2-tensor $F$
  \begin{align*}
    \int_{\Si_t}|(t\nab)^2F|^2 & \les \int_{\Si_t}t^4|\nab \Div F|^2 + \int_{\Si_t}t^4|\nab \Curl F|^2 + \int_{\Si_t}\le|(t\nab)^{\leq 1}F\ri|^2\\
    & \quad  + t^4\norm{\Nd \Fslash}^2_{\tilde{H}^{1/2}(\pr\Si_t)} +t^4\norm{\Nd \Fslash_{\Nf}}^2_{\tilde{H}^{1/2}(\pr\Si_t)} + t^4\norm{\Nd (F_{\Nf\Nf})}^2_{\tilde{H}^{1/2}(\pr\Si_t)},
  \end{align*}
  where $\Fslash,\Fslash_\Nf$ are respectively the projections of $F$ and $F_{\Nf\cdot}$ as $\pr\Si_t$-tangent tensors. 
\end{lemma}
\begin{proof}
  The proof follows from rescaling in $t$ in the results of~\cite[Lemma A.2]{Czi.Gra19a}, using the definition of the operators $\Div,\Curl$ to express the boundary terms only in terms of tangential derivatives, and using standard $H^{1/2}\times H^{-1/2}$ and trace estimates. Details are left to the reader.
\end{proof}

\section{Control of $\RR_{\leq 1}^\intr$}\label{sec:curvestbot1}
\paragraph{Control of $\Lieh_\Tf E, \Lieh_\Tf H$}
Using~(\ref{est:intQ1intSit}), the fact that $\TI=\Tf$, and that -- by the mild Bootstrap Assumptions~\ref{BA:mildKillingMMintbot} --, $\KI$ is a future-pointing timelike vectorfield satisfying comparison bounds~(\ref{est:mildtimelikeKI}) with $\Tf$, we have
\begin{align*}
  \int_{\Si_t} t^6\le|\Lieh_\Tf\R\ri|^2 & \les \int_{\Si_t}Q(\Lieh_\TI\R)(\KI,\KI,\KI,\Tf) \les \varep^2.
\end{align*}
One has schematically
\begin{align*}
  \Lieh_\Tf\Et,~\Lieh_\Tf H & = \Lieh_\Tf\R + \R\cdot\D\Tf,
\end{align*}
from which, using the Bootstrap Assumptions~\ref{BA:curvint} and~\ref{BA:intKill} on the curvature and on $\D\Tf$ respectively, one deduces that
\begin{align}\label{est:LiehTfEHL2}
  \norm{t^3\Lieh_\Tf \Et}_{L^2(\Si_t)} +\norm{t^3\Lieh_\Tf\Ht}_{L^2(\Si_t)} & \les \varep.
\end{align}

\paragraph{Control of $\nab^{\leq 1} (E,H)$}
From Maxwell equations~\eqref{eq:Maxwellt}, the traceless symmetric $\Si_t$-tangent $2$-tensor $\Et$ satisfies schematically the following div-curl system
\begin{align}\label{eq:HodgeE}
  \begin{aligned}
    \tr E & = 0,\\
    \Div E & = \Err(\Div,E),\\
    \Curl E & = \Lieh_\Tf H + \Err(\Curl,E),
  \end{aligned}
\end{align}
where the error terms $\Err(\Div, E)$, $\Err(\Curl, E)$ are of the following form
\begin{align*}
  \Err & = (\nab n, k)\cdot(E,H).
\end{align*}
Using~\eqref{eq:HodgeE} and the elliptic estimate from Lemma~\ref{lem:ellHodge}, we have
\begin{align*}
  \int_{\Si_t}t^4\le|(t\nab)^{\leq 1} E\ri|^2 & \les \int_{\Si_t}t^6|\Lieh_\Tf H|^2 + t^6\le|\int_{\pr\Si_t}\Eslash_\Nf\cdot\Nd(E_{\Nf\Nf})\ri| + \int_{\Si_t}t^6|\Err|^2\\
  & \les \varep^2 + t^6\le|\int_{\pr\Si_t}\Eslash_\Nf\cdot\Nd(E_{\Nf\Nf})\ri|,
\end{align*}
where the error terms are estimated using the Bootstrap Assumptions~\ref{BA:curvint} and~\ref{BA:intKill}.\\

From relations~\eqref{eq:defnutnur} and~\eqref{eq:NTT} between $\Tf,\Nf$ and $\elb,\el$, one has
\begin{align*}
  E_{\Nf\Nf} & = \R(\Tf,\Nf,\Tf,\Nf) \\
             & = \quar\R(\elb,\el,\elb,\el) \\
             & = \rho.
\end{align*}

Using an $H^{1/2}\times H^{-1/2}$ estimate on $\pr\Si_t$,\footnote{Such an estimate can be obtained from the definitions~\cite[pp. 834--836]{Sha14} and~\cite[Appendix B]{Czi.Gra19}.} we have
\begin{align*}
  t^6\le|\int_{\pr\Si_t}\Eslash_\Nf\cdot\Nd(E_{\Nf\Nf})\ri| & \les t^6\norm{\Eslash_N}_{\tilde{H}^{1/2}(\pr\Si_t)}\norm{\rho-\rhoo}_{\tilde{H}^{1/2}(\pr\Si_t)}. 
\end{align*}

Using the trace estimates~\eqref{est:MMbotbdy} for $\rho-\rhoo$ and an absorption argument for $\Eslash_N$, we therefore conclude
\begin{align*}
  \int_{\Si_t}t^4\le|(t\nab)^{\leq 1}E\ri|^2 & \les \varep^2,
\end{align*}
as desired. The estimates for $\nab^{\leq 1}H$ follow similarly. This concludes the control of $\RR^\intr_{\leq 1}$.

\section{Control of $\RR^\intr_{2}$}\label{sec:controlRRintr2}
\paragraph{Control of $\Lieh_\SI\Lieh_\Tf (\Et,\Ht)$}
Using~\eqref{est:intQ2intSit} and the mild Bootstrap Assumptions~\ref{BA:mildKillingMMintbot}, one has
\begin{align}\label{est:SITfEH}
  \int_{\Si_t}\le(t^6\le|\Lieh_\SI\Lieh_\Tf E\ri|^2 + t^6\le|\Lieh_\SI\Lieh_\Tf H\ri|^2\ri) & \les \varep^2.
\end{align}

\paragraph{Control of $\nab \Lieh_\Tf(E,H)$ and $\Lieh_\Tf^2(E,H)$}
Commuting Maxwell equations~\eqref{eq:Maxwellt} with $\Lieh_\Tf$, one has the following div-curl system for the symmetric traceless $\Si_t$-tangent $2$-tensor $\Lieh_\Tf E$ and $\Lieh_\Tf H$
\begin{align}\label{eq:HodgeLieTfE}
  \begin{aligned}
  \tr \Lieh_\Tf E & = 0,\\
  \Div \Lieh_\Tf E & = \Err(\Div,\Lieh_\Tf E),\\
  \Curl \Lieh_\Tf E & = \Lieh_\Tf\Lieh_\Tf H + \Err(\Curl,\Lieh_\Tf E),
  \end{aligned}
\end{align}
and
\begin{align}\label{eq:HodgeLieTfH}
  \begin{aligned}
  \tr \Lieh_\Tf H & = 0,\\
  \Div \Lieh_\Tf H & = \Err(\Div,\Lieh_\Tf H),\\
  \Curl \Lieh_\Tf H & = -\Lieh_\Tf\Lieh_\Tf E + \Err(\Curl,\Lieh_\Tf H),
  \end{aligned}
\end{align}
where the error terms are of the form
\begin{align*}
  \Err & = \le(\nab^{\leq 2}(n-1) + \nab^{\leq 1} \Lieh_\Tf^{\leq 1}(n-1) + \le(\nab,\Lieh_\Tf\ri)^{\leq 1}k\ri)\cdot(E,H)\\
  & \quad + \le(\nab^{\leq 1} (n-1) + k\ri)\cdot\le((\nab,\Lieh_\Tf)^{\leq 1}(E,H)\ri).
\end{align*}
From the definition of $\SI$ and $\XI$ in Section~\ref{sec:defKillingint}, we have
\begin{align}\label{eq:TfTfSITf}
  \Lieh_\Tf\Lieh_\Tf(E,H) & = \Lieh_{t^{-1}\SI}\Lieh_\Tf(E,H) - \Lieh_{t^{-1}\XI}\Lieh_\Tf(E,H).
\end{align}
Using~(\ref{est:LiehTfEHL2}),~\eqref{est:SITfEH} and the Bootstrap Assumptions~\ref{BA:connint} to control the derivatives of $t$, we have
\begin{align}\label{est:SITfEH2}
  \int_{\Si_t}t^8\le|\Lieh_{t^{-1}\SI}\Lieh_\Tf(E,H)\ri|^2 & \les \varep^2.
\end{align}
We define $X$ to be the projection of $t^{-1}\XI$ on $\Si_t$, \emph{i.e.}
\begin{align*}
  t^{-1}\XI & =: -\g(t^{-1}\XI, \Tf) \Tf +X.
\end{align*}
From the bootstrap bound~\eqref{est:BAXITI}, we have
\begin{align*}
  \le|(t\D)^{\leq 1}\g(t^{-1}\XI,\TI)\ri| & \les (D\varep) t^{-3/2},
\end{align*}
and the $\Tf$ component of $\XI$ in~\eqref{eq:TfTfSITf} can be treated as an error term in the following rewriting of the div-curl systems~(\ref{eq:HodgeLieTfE}) and (\ref{eq:HodgeLieTfH})
\begin{align}\label{eq:HodgeLieTfEbis}
  \begin{aligned}
  \tr \Lieh_\Tf E & = 0,\\
  \Div \Lieh_\Tf E & = \Err,\\
  \Curl \Lieh_\Tf E +\Lieh_X\Lieh_\Tf H & = \Lieh_{t^{-1}\SI}\Lieh_\Tf H + \Err,
  \end{aligned}
\end{align}
and
\begin{align}\label{eq:HodgeLieTfHbis}
  \begin{aligned}
  \tr \Lieh_\Tf H & = 0,\\
  \Div \Lieh_\Tf H & = \Err,\\
  \Curl \Lieh_\Tf H -\Lieh_X\Lieh_\Tf E & = -\Lieh_{t^{-1}\SI}\Lieh_\Tf E + \Err.
  \end{aligned}
\end{align}

From the bootstrap bound~(\ref{est:BAXIXI}), we have
\begin{align*}
  \le|\sup_{\Si_t} |X| - \frac{1-\cc}{1+\cc}\ri| & \les \varep t^{-3/2}.
\end{align*}
Therefore, there exists a numerical constant $0 < \cco < 1$, such that for all transition parameters $\cc$, \emph{i.e.} $\cco \leq \cc \leq (1+\cco)/2$, and for $\varep>0$ sufficiently small, we have
\begin{align*}
  |X| \leq \wp,
\end{align*}
in $\MM^\intr_\bott$, where $\wp$ is the constant from Lemma~\ref{lem:ellHodgevar}. Moreover, from the Bootstrap Assumptions~\ref{BA:intKill} on $\D\XI$, we have
\begin{align*}
  \le|(t\nab)X\ri| & \les 1.
\end{align*}

One can therefore apply the result of Lemma~\ref{lem:ellHodgevar} to the modified div-curl systems~\eqref{eq:HodgeLieTfEbis} and~\eqref{eq:HodgeLieTfHbis} for $\Lieh_\Tf E$ and $\Lieh_\Tf H$, using estimate~\eqref{est:SITfEH2} and the estimates for the lower order terms obtained in Section~\ref{sec:curvestbot1} and the Bootstrap Assumptions~\ref{BA:curvint} and~\ref{BA:connint} to control the nonlinear error terms, and we have
\begin{align*}
  \int_{\Si_t}t^{6} \le|(t\nab)\Lieh_\Tf(E,H)\ri|^2 & \les \varep^2 + \le|\int_{\pr\Si_t}\Lieh_\Tf E_{\Nf} \cdot \Nd\le(\Lieh_\Tf E_{NN}\ri) \ri| + \le|\int_{\pr\Si_t}\Lieh_\Tf H_{\Nf} \cdot \Nd\le(\Lieh_\Tf H_{NN}\ri) \ri|.
\end{align*}

We have on $\pr\Si_t$
\begin{align*}
  \Lieh_\Tf E_{NN} & = \elb(\rho) + \el(\rho) + \Err.
\end{align*}
Therefore using the trace estimates~(\ref{est:MMbotbdy}) and $H^{-1/2}\times H^{1/2}$ estimates as previously gives the control of the integral boundary term for $\Lieh_\Tf E$. Arguing similarly gives the control of the integral boundary term for $\Lieh_\Tf H$, and we finally obtain
\begin{align}\label{est:nabLiehTfEH}
  \int_{\Si_t}t^{6} \le|(t\nab)^{\leq 1}\Lieh_\Tf(E,H)\ri|^2 & \les \varep^2.
\end{align}

Using~\eqref{eq:TfTfSITf} and the above mentioned controls for $\XI$, we also deduce from~\eqref{est:SITfEH2} and~\eqref{est:nabLiehTfEH} that
\begin{align*}
  \int_{\Si_t}t^{8} \le|\Lieh_\Tf^2(E,H)\ri|^2 & \les \varep^2.
\end{align*}

\paragraph{Control of $\nab^{2}(E,H)$}
The control of $\nab^2E$ follows from an application of the higher order elliptic estimates of Lemma~\ref{lem:ellSitHrr} to the div-curl systems~(\ref{eq:HodgeE}), using the estimate~\eqref{est:nabLiehTfEH} for $\nab\Lieh_\Tf H$ and the trace estimates~\eqref{est:MMbotbdy}. The estimate for $\nab^2H$ follow similarly.

\paragraph{Control of $\D^{\leq 2}\R$}
Using the definition~(\ref{eq:defEH}) of $E,H$, the control of $\Tf,\D\Tf,\D^2\Tf$ from the Bootstrap Assumptions~\ref{BA:intKill}, and the Bootstrap Assumptions~\ref{BA:curvint}, one deduces from the bounds obtained for $E,H$ that
\begin{align*}
  \norm{t^2(t\D)^{\leq 2}\R}_{L^2(\Si_t)} & \les \varep,
\end{align*}
for all $\too \leq t \leq \tast$, and where the norms are taken with respect to the maximal frame. This finishes the proof of~(\ref{est:RRintbott}).

\chapter{Remaining curvature estimates}\label{sec:remainingcurvestfinal}
In this section, we prove that curvature estimates hold in the interior of the cone $\CCba\cap^{(\cc)}\MM^\intr$. We also prove that the curvature estimates hold on all hypersurfaces associated to all transition parameter $\cco\leq\cc\leq (1+\cco)/2$. Both arguments go by rescaling and local energy estimates.

\section{Null curvature estimates on $\protect\CCba\cap\MM^\intr$}\label{sec:tipcurvest}
In this section, we prove the following proposition.
\begin{proposition}\label{prop:tipcurvestSTAB}
  Recall that from Proposition~\ref{prop:planehypcurvestSTAB}, we have on $^{(\cc)}\Si_\tast$
  \begin{align}\label{est:enerestlastSi}
  \norm{\tast^2 (\tast\D)^{\leq 2}\R}_{L^2(^{(\cc)}\Si_\tast)} & \les \varep,
  \end{align}
  for a fixed transition parameter $\cc$.\\
  
  Under the Bootstrap Assumptions and~\eqref{est:enerestlastSi} and for $\varep>0$ sufficiently small, we have the following bounds for the null curvature components on the interior of the last cone $\CCba\cap{^{(\cc)}\MM}^\intr$
  \begin{align}\label{est:RRbintr2}
    \RR^\ast_{\leq 2} & \les \varep.
  \end{align}
  Moreover, we also have the following $L^\infty \HHt$ estimates in $\CCba\cap{^{(\cc)}\MM}^\intr$
  \begin{align}\label{est:RRfbintrCCbafinish}
    \mathfrak{R}^\ast_{\leq 1} & \les \varep.
  \end{align}
  We refer the reader to Section~\ref{sec:normnullcurv} for definitions. 
\end{proposition}



\begin{remark}\label{rem:notoptimalaxisdege}
  The $L^2$-norms of~\eqref{est:RRbintr2} (see the definitions of Section~\ref{sec:normnullcurv}) allow for a degeneracy of the type $R \sim r^{-3/2}$ when $r\to 0$ (see the Sobolev embeddings of Lemma~\ref{lem:KlSobast}). This is suboptimal by far and could be easily improved using that we obtain independently the uniform boundedness of the spacetime curvature tensor $|\R| \les \varep$ in the extended spacetime of Section~\ref{sec:tastSitast}. In this section, we could improve the $L^2$-bounds~\eqref{est:RRbintr2} by a more careful analysis of relations between Cartesian and spherical derivatives, see Remark~\ref{rem:rGanotoptimal}. However -- since the null decompositions of $\R$ are multi-valued at the vertex --, there is no hope to obtain $L^2$ norms on $\CCba$ sufficiently regular to get boundedness of $R$ when $r\to 0$ \emph{via} Sobolev embeddings, since continuity would follow as well. As a consequence, we do not seek for rates in $r$ in the $L^2$ norms better than the (easily obtained) ones of~\eqref{est:RRbintr2}. Fortunately, we can still obtain a control of the null connection coefficients consistent with the potential $r^{-3/2}$ singularity of the null curvature components on the cone $\CCba$ (see Remark~\ref{rem:nullconndegeaxis} and~\cite{Wan09}).  
\end{remark}

\subsection{$\rast$-rescaling, extension and local existence from $\Si_\tast$}\label{sec:tastSitast}
We perform a $\rast$-rescaling of the last maximal hypersurface $\Si_\tast$ to a maximal hypersurface $\Sitt_1$ of size $1$, \emph{i.e.} such that $\pr\Sitt_1$ has area radius $1$.\footnote{Recall that $\rast$ is defined to be the area radius of $\pr\Si_\tast=S^\ast$ and that $\rast \simeq \tast$.} The bounds~\eqref{est:enerestlastSi} rewrite as
\begin{align}\label{est:scaleenerestlastSi}
  \norm{\D^{\leq 2}\R}_{L^2(\Sitt_1)} & \les \varep_1,
\end{align}
where $\D,\R$ are associated to the rescaled metric, where the frame norm is adapted to the unit normal to $\Sitt_1$ for the rescaled metric, and where
\begin{align*}
  \varep_1  & := \varep (\tast)^{-3/2}.
\end{align*}
The bounds from the Bootstrap Assumptions~\ref{BA:harmoSitast} and~\ref{BA:connint} rewrite as\footnote{Bounds for harmonic coordinates are equivalent to bounds for the metric in these coordinates. See also Theorem~\ref{thm:globharmonics}.}
\begin{align}\label{est:rescaleBASitast}
  \begin{aligned}
    \norm{\pr^{\leq 4}\le(g_{ij}-\de_{ij}\ri)}_{L^2(\Sitt_1)} & \les D\varep_1,\\
    \norm{\nab^{\leq 3}k}_{L^2(\Sitt_1)} & \les D\varep_1.
  \end{aligned}
\end{align}

Using the extension theorem~\cite[Theorem 3.1]{Czi18}, we obtain an extension $\tilde{g}_{ij},\tilde{k}_{ij}$ of $(g_{ij},k_{ij})$ defined on $\RRR^3$, \emph{i.e.} such that $\tilde {g}_{ij} = g_{ij}$ and $\tilde{k}_{ij} = k_{ij}$ on $\Sitt_1 \simeq \mathbb{D} \subset \RRR^3$, that satisfy the maximal constraint equations, and such that
\begin{align*}
  \norm{\pr^{\leq 4}\le(\tilde{g}_{ij}-\de_{ij}\ri)}_{L^2(\RRR^3)} & \les D\varep_1,\\
  \norm{\pr^{\leq 3}\tilde{k}_{ij}}_{L^2(\RRR^3)} & \les D\varep_1,
\end{align*}
with suitable fall-off rate at infinity (see~\cite{Czi18} for precisions).\\

Using the local existence result for such an initial data set satisfying the maximal constraint equations on $\RRR^3$ (see~\cite[Theorem 10.2.2]{Chr.Kla93}), and provided that $D\varep_1$ is sufficiently small, there exists a spacetime $D\varep_1$-close to $[0,3/2]\times\RRR^3 \subset \RRR^{1+3}$ (see~\cite{Chr.Kla93} for the precise definitions). This spacetime admits a maximal time function $\tilde{t}$ such that $\tilde{t}=0$ on $\Sitt_1$ and we have the following control for the norms of the spacetime curvature tensor in the maximal frame norm
\begin{align}\label{est:curvextension}
  \begin{aligned}
    \norm{\R}_{L^\infty} & \les D\varep_1,\\
    \norm{\D\R}_{L^\infty_{\tilde{t}}L^6(\Si_{\tilde{t}})} & \les D\varep_1,\\
    \norm{\D^2\R}_{L^\infty_{\tilde{t}}L^2(\Si_{\tilde{t}})} & \les D\varep_1,
  \end{aligned}
\end{align}
locally in the domain of dependence of $\Sitt_1$. Moreover, there exists approximate Cartesian Killing vectorfields $\pr_\mu$ for $\mu = 0,1,2,3$, such that
\begin{align}\label{est:Killingextension}
  \begin{aligned}
  \norm{\D\pr_\mu}_{L^\infty} & \les D\varep_1,\\
  \norm{\D^2 \pr_\mu}_{L^\infty_{\tilde{t}}L^2(\Si_{\tilde{t}})} & \les D\varep_1,\\
  \norm{\D^3\pr_\mu}_{L^\infty_{\tilde{t}}L^2(\Si_{\tilde{t}})} & \les D\varep_1,
  \end{aligned}
\end{align}
and where $\pr_0$ is defined as
\begin{align*}
  \pr_0 := \frac{\D\tilde{t}}{(-\g(\D \tilde{t},\D \tilde{t}))^{1/2}}.
\end{align*}

\subsection{Energy estimates in $\MM^\intr_\topp$}\label{sec:enerDD}
We define the spacetime region $\DD_1$ to be the future domain of dependence of $\Sitt_1$.\footnote{$\DD_1$ corresponds to the $\rast$ rescaling of $\MM^\intr_\topp$.} In this section, we apply energy estimates in the region $\DD_1$ to obtain the desired curvature control on the $\rast$-rescaling of the cone $\CCba\cap\MM^\intr$.\\

More precisely, we integrate the following contracted Bel-Robinson tensors in $\DD_1$
\begin{align*}
  Q(\R)(\pr_0,\pr_0,\pr_0), \quad Q(\Lieh_{\pr_\mu}\R)(\pr_0,\pr_0,\pr_0), \quad Q(\Lieh_{\pr_\mu}\Lieh_{\pr_\nu}\R)(\pr_0,\pr_0,\pr_0),
\end{align*}
where the vectorfields $\pr_\mu$ are the approximate Cartesian Killing fields from Section~\ref{sec:tastSitast}.\\

For each estimates, the spacetime error term in $\DD_1$ is controlled using the estimates~\eqref{est:curvextension} and~\eqref{est:Killingextension} (see also the treatment of the error terms in $\MM^\intr_\bott$ and $\MM^\ext$ in Section~\ref{sec:globener}). We therefore obtain the following estimates on the future boundary $\CCbb_1$ of $\DD_1$
\begin{align*}
  \int_{\CCbb_1}Q(\R)(\pr_0,\pr_0,\pr_0,\elb) & \les \int_{\Sitt_1}Q(\R)(\pr_0,\pr_0,\pr_0,\Tf) + (D\varep_1)^3,\\
  \int_{\CCbb_1}Q(\Lieh_{\pr_\mu}\R)(\pr_0,\pr_0,\pr_0,\elb) & \les \int_{\Sitt_1}Q(\Lieh_{\pr_\mu}\R)(\pr_0,\pr_0,\pr_0,\Tf) + (D\varep_1)^3, \\
  \int_{\CCbb_1}Q(\Lieh_{\pr_\mu}\Lieh_{\pr_\nu}\R)(\pr_0,\pr_0,\pr_0,\elb) & \les \int_{\Sitt_1}Q(\Lieh_{\pr_\mu}\Lieh_{\pr_\nu}\R)(\pr_0,\pr_0,\pr_0,\Tf) + (D\varep_1)^3,
\end{align*}
where $\Tf$ is the future-pointing unit normal (with respect to the rescaled metric) to $\Sitt_1$ and $\elb,\el$ is the rescaled null pair on $\CCbb_1$ (\emph{i.e.} such that $\elb,\el$ is a null pair for the rescaled metric).\\

Using the improved curvature estimate~(\ref{est:scaleenerestlastSi}) and the control~(\ref{est:Killingextension}) of the vectorfields $\pr_\mu$ on $\Sitt_1$, this yields
\begin{align}\label{est:fluxCCb1}
  \begin{aligned}
  \int_{\CCbb_1}Q(\R)(\pr_0,\pr_0,\pr_0,\elb) & \les \varep_1^2,\\
  \int_{\CCbb_1}Q(\Lieh_{\pr_\mu}\R)(\pr_0,\pr_0,\pr_0,\elb) & \les \varep_1^2, \\
  \int_{\CCbb_1}Q(\Lieh_{\pr_\mu}\Lieh_{\pr_\nu}\R)(\pr_0,\pr_0,\pr_0,\elb) & \les \varep_1^2,
  \end{aligned}
\end{align}
for $\varep_1$ sufficiently small.\\


We now compare the coordinate vectorfields $\pr_0$ to $T := \half(\elb + \el)$ on $\CCbb_1$. Since $\pr_0 = \Tf$ on $\Sitt_1$, we have on $\pr\Sitt_1$
\begin{align*}
  \g(\pr_0,T) & = \g(\Tf,T) = -\half\le(\nut^{-1}+\nut\ri),
\end{align*}
where here $\nut$ corresponds to a rescaling of the original $\nut$. From the (rescaled) bootstrap bound on $\nu$ from the Bootstrap Assumptions~\ref{BA:connint}, this gives on $\pr\Sitt_1$
\begin{align*}
  \le|\g(\pr_0,T) +1 \ri| & \les D\varep_1.
\end{align*}
 Using the relations~(\ref{eq:Riccirel}), we further have
\begin{align*}
  |\D_{\elb}(\g(\pr_0,T))| & \les |\D\pr_0||\elb| + |\D_\elb \elb| + |\D_\elb\el| \\
                    & \les |\D\pr_0||\elb| + |\eta||\ea| + |\omb|(|\elb|+|\el|),
\end{align*}
where the norms are taken with respect to the maximal frame (\emph{i.e.} $e_0 = \pr_0$ in the notations of Definition~\ref{def:framenorm}) and where the null coefficients correspond to rescaling of the original coefficients.\footnote{Alternatively, they are the null coefficients associated to the rescaled null pair $\elb,\el$.}\\

Using the (rescaled) Bootstrap Assumptions~\ref{BA:connCCba} on the null connection coefficients on $\CCba$ and~(\ref{est:Killingextension}), we obtain
\begin{align*}
  |\D_{\elb}(\g(\pr_0,T))| & \les r^{-1/2}D\varep_1\le(|\elb| + |\el| + |\ea|\ri).
\end{align*}

From Remark~\ref{rem:defframenorm}, using that $e_3,e_4$ are future-pointing, we have
\begin{align*}
  |\elb| + |\el| & = -\sqrt{2}\g(\elb,\pr_0) - \sqrt{2}\g(\el,\pr_0) = 2\sqrt{2}\g(T,\pr_0),\\
  \sum_{a=1,2}|\ea|^2 & = 2+2\sum_{a=1,2}|\g(\pr_0,\ea)|^2 \leq 2+\le(|\pr_0|_{T}^2\ri) \leq 1 + 2|\g(\pr_0,T)|^2,
\end{align*}
where $|\cdot|_T$ is the frame norm with respect to $T$.\\

Thus, we have
\begin{align*}
  \le|\D_\elb(\g(\pr_0,T))\ri| & \les r^{-1/2}D\varep_1\le(1 + \le|\g(\pr_0,T)\ri|\ri),
\end{align*}
and from a Gr\"onwall argument, integrating from $r=1$ to $r=0$, we obtain
\begin{align}\label{est:comppr0T}
  \begin{aligned}
  |\g(\pr_0,T) + 1| & \les D\varep_1
  \end{aligned}
\end{align}
on $\CCbb_1$, \emph{i.e.} the frames adapted to $\pr_0$ and $T$ are comparable.\\

Using~\eqref{est:fluxCCb1} and~\eqref{est:comppr0T} we obtain
\begin{align}\label{est:fluxCCb1bis}
  \begin{aligned}
    \int_{\CCbb_1}|\nulld(\R)|^2 & \les \varep_1^2,\\
    \int_{\CCbb_1}\le|\nulld(\Lieh_{\pr_\mu}\R)\ri|^2 & \les \varep_1^2, \\
    \int_{\CCbb_1}\le|\nulld(\Lieh_{\pr_\mu}\Lieh_{\pr_\nu}\R)\ri|^2 & \les \varep_1^2,
  \end{aligned}
\end{align}
where in that case $\nulld \in\le\{\alb,\cdots,\be\ri\}$. This provides in particular the first desired estimates for the $L^2$ norms of the curvature.\\

Using~\eqref{est:Killingextension} and~\eqref{est:scaleenerestlastSi}, we have
\begin{align*}
  \int_{\CCbb_1}\le|\nulld(\D_{\pr_\mu}\R)\ri|^2 & \les \int_{\CCbb_1}\le|\nulld(\Lieh_{\pr_\mu}\R)\ri|^2  +(D\varep_1)^3, \\
  \int_{\CCbb_1}\le|\nulld(\D^2_{\pr_\mu,\pr_\nu}\R)\ri|^2 & \les \int_{\CCbb_1}\le|\nulld(\Lieh_{\pr_\mu}\Lieh_{\pr_\nu}\R)\ri|^2 + (D\varep_1)^3,
\end{align*}
which, using~\eqref{est:comppr0T} and the resulting comparison between $\pr_\mu$ and $\elb,\el,\ea$ further gives
\begin{align*}
  \int_{\CCbb_1}\le|\nulld(\D_{X}\R)\ri|^2 & \les \varep_1^2,\\
  \int_{\CCbb_1}\le|\nulld(\D^2_{X,Y}\R)\ri|^2 & \les \varep_1^2,
\end{align*}
where $X,Y\in\{\elb,\el,\ea\}$.\\

We have schematically
\begin{align*}
  \Nd_X\nulld(\R) & = \nulld(\D_X\R) + \R\cdot\Ga,
\end{align*}
where\footnote{The null coefficients are here rescaling of the original coefficients. Alternatively they are the null coefficients associated to the rescaled null pair $\elb,\el$ used in this section.}
\begin{align*}
  \Ga\in\{\trchi,\trchib,\chih,\chibh,\eta,\ze,\etab,\omb\}.
\end{align*}
Multiplying by $r$, and since by the Bootstrap Assumptions~\ref{BA:connCCba}, 
\begin{align}\label{est:rGanotoptimal}
  |r\Ga| & \les 1,
\end{align}
we therefore infer that for $X\in\{\elb,\el,\ea\}$ and since $r\les 1$,
\begin{align*}
  \int_{\CCbb_1}\le|r\Nd_X\nulld(\R)\ri|^2 & \les \int_{\CCbb_1} r^2\le|\nulld(\D_{X}\R)\ri|^2 + \int_{\CCbb_1}|r\Ga|^2|\R|^2 \\
  & \les \varep_1^2.
\end{align*}
Similarly, since we have schematically
\begin{align*}
  \Nd_X\Nd_Y\nulld(\R) & = \nulld(\D_{X,Y}\R) + \D\R\cdot\Ga + \R\cdot\Ga\cdot\Ga,
\end{align*}
for $X,Y\in\{\elb,\el,\ea\}$, multiplying by $r^2$ and taking the $L^2$ norm on $\CCbb_1$, we conclude
\begin{align*}
  \int_{\CCbb_1}\le|r^2\Nd^2_{X,Y}\nulld(\R)\ri|^2 & \les \varep_1^2.
\end{align*}
Rescaling back theses estimates on the original cone $\CCba\cap\MM^\intr$ gives the desired bounds for the curvature, and finishes the proof of~\eqref{est:RRbintr2}.

\begin{remark}\label{rem:rGanotoptimal}
  Estimate~\eqref{est:rGanotoptimal} is suboptimal in terms of $r$-weight for most of the connection coefficients. More precise bounds together with Hardy estimates (see Lemma~\ref{lem:transportast}) would lead to improve bounds. See Remark~\ref{rem:notoptimalaxisdege}. 
\end{remark}


\subsection{Proof of the decay estimate~\eqref{est:RRfbintrCCbafinish}}
The proof of~\eqref{est:RRfbintrCCbafinish} boils down to the following Klainerman-Sobolev embeddings. 
\begin{lemma}[Klainerman-Sobolev estimates on $\CCba\cap\MM^\intr$]\label{lem:KlSobastintr}
  For all $S$-tangent tensor $F$ with vertex limit
  \begin{align*}
    r^{3/2}|(r\Nd)^{\leq 1}F|\to 0 \quad \text{when $r\to 0$,} 
  \end{align*}
  we have the following $L^\infty_{u}\HHt\le(S_{u,\uba}\ri)$ estimates in $\CCba\cap\MM^\intr$
  \begin{align*}
    \norm{r F}_{L^\infty_{u\geq \cc\uba}\HHt\le(S_{u,\uba}\ri)} & \les \norm{F}_{L^2(\CCba\cap\MM^\intr)} + \norm{r\Nd F}_{L^2(\CCba\cap\MM^\intr)} +\norm{r\Nd_3F}_{L^2(\CCba\cap\MM^\intr)}.
  \end{align*}
  

\end{lemma}
The proof of Lemma~\ref{lem:KlSobastintr} is postponed to Appendix~\ref{app:KlSobH12}. We apply the Klainerman-Sobolev estimates of Lemma~\ref{lem:KlSobastintr}, with $F$ the following respective tensors
\begin{align*}
  u^2\Ndt^{\leq 1}\alb,~u\ub\Ndt^{\leq 1}\beb,~\ub^2\Ndt^{\leq 1}(\rho-\rhoo),~\ub^2\Ndt^{\leq 1}(\sigma-\sigmao),~\ub^2\Ndt^{\leq 1}\be,~\ub^2(r\Nd)^{\leq 1}\al,~\ub^2(r\Nd_3)\al,
\end{align*}
where $\Ndt \in\{(r\Nd),(q\Nd_3)\}$. From Theorems~\ref{thm:vertex} and~\ref{thm:canonical}, these tensors have the asymptotic behaviour $r^{3/2}(r\Nd)^{\leq 1}R = O\le(r^{3/2}\ri)$ when $r\to 0$ and satisfy in particular the required limits of Lemma~\ref{lem:KlSobastintr}. From an inspection of the definitions of Section~\ref{sec:normnullcurv}, and using that in $\CCba\cap\MM^\intr$ we have $r\simeq q$ and $u\simeq\ub$, we deduce from the bounds~(\ref{est:RRbintr2}) that
\begin{align*}
  \mathfrak{R}^\ast_{\leq 1} & \les \RR^\ast_{\leq 2} \les \varep,
\end{align*}
where the norms are restricted to $\CCba\cap\MM^\intr$. This proves~(\ref{est:RRfbintrCCbafinish}).\\

\section{Curvature estimates for all transition parameters}\label{sec:alltransparam}
In this section, we call $\cc_1:=\cc$ the fixed transition parameter given by the mean value argument of Section~\ref{sec:meanvalue}. We prove the following proposition.
\begin{proposition}\label{prop:alltransparamSTAB}
  Recall that from Propositions~\ref{prop:planehypcurvestSTAB}, we have
  \begin{align}\label{est:sourcecc1intr}
    ^{(\cc_1)}\RR^\intr_{\leq 2} & \les \varep,
  \end{align}
  in $^{(\cc_1)}\MM^\intr_\bott$. Recall that from Lemma~\ref{lem:bdedcurvflux}, we have
  \begin{align}\label{est:sourcecc1CCu}
    \int_{\cc_1^{-1}u}^{\cco^{-1}u}\le(\RR_{\leq 2}(u,\ub)\ri)^2 \,\d\ub & \les \varep^2,
  \end{align}
  for all $1 \leq u \leq \cc_1 \uba$.\\

  Under the Bootstrap Assumptions and estimates~\eqref{est:sourcecc1intr} and~\eqref{est:sourcecc1CCu}, and for $\varep>0$ sufficiently small, we have the following bounds \emph{for all transition parameters $\cco \leq \cc \leq (1+\cco)/2$}
  \begin{align}
    ^{(\cc)}\RR^{\ext}_{\leq 2,\ga} & \les_\ga \varep,\label{est:cc1goal1}\\
    ^{(\cc)}\mathfrak{R}^\ext_{\leq 1} & \les \varep,\label{est:cc1goal2}
  \end{align}
  for all $\ga>0$ and
  \begin{align}
    ^{(\cc)}\RR^{\intr}_{\leq 2} & \les \varep.\label{est:cc1goal3}
  \end{align}
  We refer to Sections~\ref{sec:normnullcurv} and \ref{sec:normscurvint} for the definitions of these norms. Here we indicate by $^{(\cc)}$ the dependency on either just the covered domain (in estimates~\eqref{est:cc1goal1},~\eqref{est:cc1goal2}) or the covered domain and the nature of the hypersurfaces (in estimate~\eqref{est:cc1goal3}). 
\end{proposition}

\begin{remark}\label{rem:alltransrem}
  As a result of Proposition~\ref{prop:alltransparamSTAB}, the curvature estimates of Propositions~\ref{prop:curvestSTAB} and~\ref{prop:planehypcurvestSTAB} obtained for the fixed transition parameter $\cc_1$ given by the mean value argument of Section~\ref{sec:meanvalue}, hold \emph{for all transition parameter $\cco \leq \cc \leq (1+\cco)/2$}. In the remaining Sections~\ref{sec:connestCCba}~--~\ref{sec:initlayer}, we shall thus assume that $\cc$ is any parameter in $[\cco,(1+\cco)/2]$. This will improve the mild and strong Bootstrap Assumptions (see Sections~\ref{sec:mildBA} and \ref{sec:strongBA}) \emph{for all transition parameter $\cc$}.
\end{remark}

The remaining $L^2(\MM^\ext)$ estimates of~\eqref{est:cc1goal1} not already obtained in Proposition~\ref{prop:curvestSTAB} are directly obtained from~\eqref{est:sourcecc1intr} and a comparison of frame argument. The decay estimates~\eqref{est:cc1goal2} are obtained \emph{via} the Klainerman-Sobolev embeddings of Lemma~\ref{lem:KlSobSitext} and $L^2$ bounds on $^{(\cc)}\Si_t^\ext$ (see the proof of the analogous estimates in Section~\ref{sec:decayH12extcurv}). Thus Proposition~\ref{prop:planehypcurvestSTAB} boils down to the following three lemmas.
\begin{lemma}\label{lem:proofproptrans1}
  For all $\cc_1 \leq \cc \leq (1+\cco)/2$, we have
  \begin{align}\label{est:proofproptrans1}
    \int_{^{(\cc)}\Si_t} t^4 |(t\D)^{\leq 2}\R|^2 + \int_{\cc^{-1}u}^{\cc_1^{-1}u}\le(\RR_{\leq 2}(u,\ub)\ri)^2\,\d\ub & \les \varep^2,
  \end{align}
  for all $^{(\cc)}\too \leq t \leq ^{(\cc)}\tast$ and for all $1 \leq u \leq \cc\uba$ and where the frame norm is adapted to $^{(\cc)}\Tf$.
\end{lemma}

\begin{lemma}\label{lem:proofproptrans2}
  For all $\cco \leq \cc \leq \cc_1$, we have
  \begin{align}\label{est:proofproptrans2}
    \int_{^{(\cc)}\Si_t} t^4 |(t\D)^{\leq 2}\R|^2 & \les \varep^2,
  \end{align}
  for all $^{(\cc)}\too \leq t \leq ^{(\cc)}\tast$ and where the frame norm is adapted to $^{(\cc)}\Tf$.
\end{lemma}

\begin{lemma}\label{lem:proofproptrans3}
  For all $\cc_1 \leq \cc \leq (1+\cco)/2$, we have
  \begin{align}\label{est:proofproptrans3}
    \int_{^{(\cc)}\Si^\ext_t \cap {^{(\cc_1)}\MM}^\intr} t^4 |(t\D)^{\leq 2}\R|^2 & \les \varep^2,
  \end{align}
  for all $^{(\cc)}\too \leq t \leq ^{(\cc)}\tast$ and where the frame norm is adapted to $\Tf^{\ext}$.
\end{lemma}

\subsection{Proof of Lemma~\ref{lem:proofproptrans1}}\label{sec:proofproptrans1}
In this section, we assume that $\cc_1 \leq \cc \leq (1+\cco)/2$.\\

Let $^{(\cc)}\too \leq t \leq ^{(\cc)}\tast$. Let
\begin{align*}
  u & := \frac{2}{1+\cc^{-1}}t, & t_1 & := \frac{1+\cc_1^{-1}}{2}u.
\end{align*}

\begin{remark}
  In the case where $t_1\geq ^{(\cc_1)}\tast$, \emph{i.e.} when $^{(\cc)}\Si_t$ is in the domain of dependence of the last slice $^{\cc_1}\Si_{^{(\cc_1)}\tast}$ the bound~\eqref{est:proofproptrans1} follows from the arguments of Section~\ref{sec:tipcurvest} and a comparison of frame. In the following, we restrict to the case $t_1\leq ^{(\cc_1)}\tast$.
\end{remark}

\begin{figure}[h!]
  \centering
  \includegraphics[height=6cm]{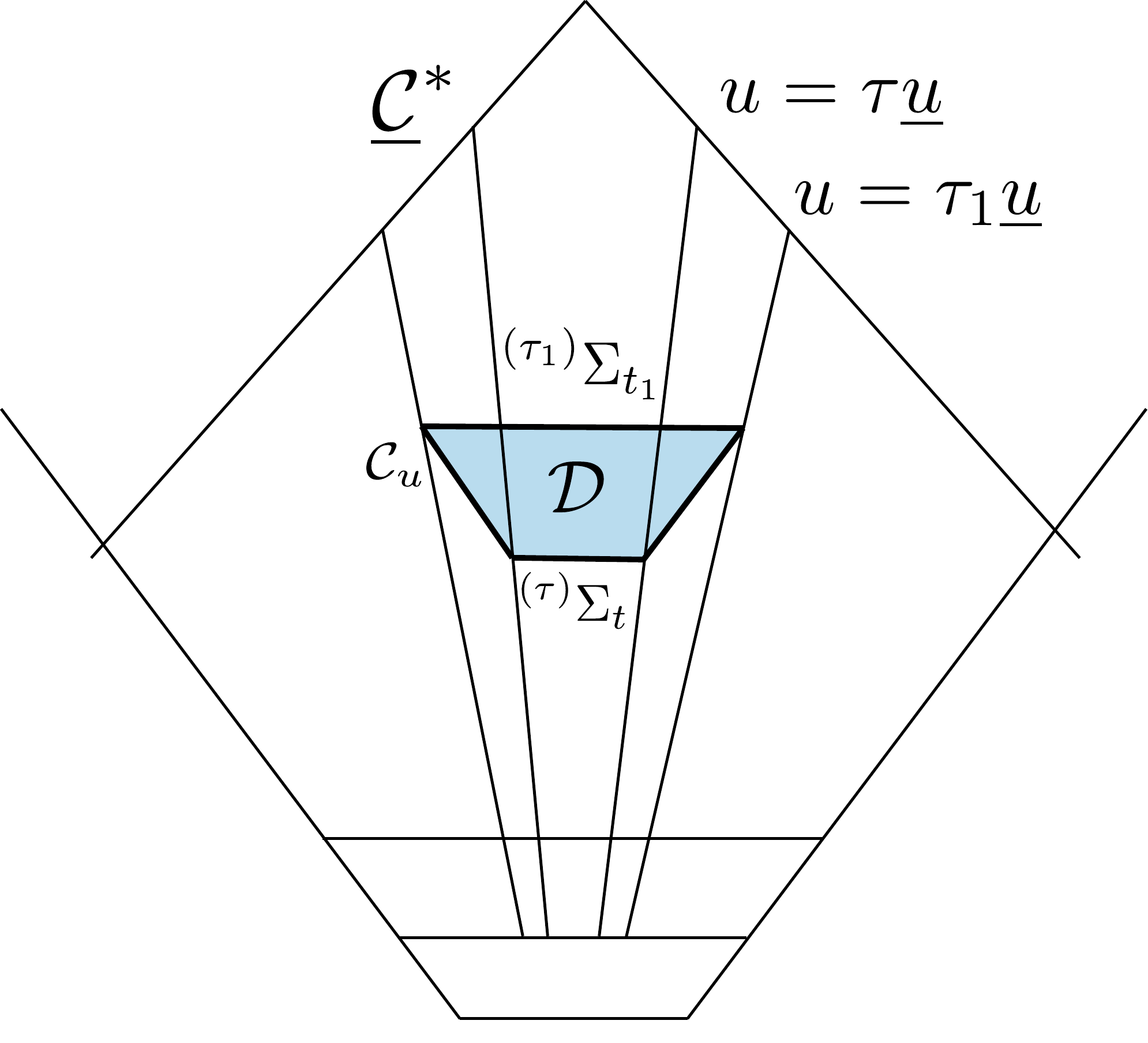}
  \caption{Local energy estimates for $\cc \geq \cc_1$.}
  \label{fig:localenergccggcc1}
\end{figure}

We consider the maximal hypersurface $^{(\cc_1)}\Si_{t_1}$. We perform a $t_1$-rescaling of the region $\DD$ enclosed by $^{(\cc)}\Si_t\cup\CC_u\cup{^{(\cc_1)}\Si}_{t_1}$, which rescales to a size $1$ region (see Figure~\ref{fig:localenergccggcc1}).\footnote{From its definition, we have $t_1 \simeq t$.} We call $\varep_1$ the rescaled smallness parameter $\varep$ (see the similar definition in Section~\ref{sec:tastSitast}).\\

Using that by the Bootstrap Assumptions~\ref{BA:curvint}, \ref{BA:connext}, \ref{BA:connint}, we have
\begin{align*}
  \norm{(\nab)^{\leq 1}\RRRic}_{L^2({^{(\cc_1)}\Si}_{t_1})} & \les D\varep_1,\\
  \norm{\Nd^{\leq 1}(\th-\gd)}_{H^{1/2}(\pr{^{(\cc_1)}\Si}_{t_1})} & \les D\varep_1,
\end{align*}
and applying the results of Theorem~\ref{thm:globharmonics}, there exists harmonic coordinates $x^i$ on $^{(\cc_1)}\Si_{t_1}$ (see also the definition of Section~\ref{sec:defUnifandharmo}) such that the following bounds hold on $^{(\cc_1)}\Si_{t_1}$
\begin{align}\label{est:harmocoordsaux}
  \norm{\nab^{\leq 2}\nab^2x^i}_{L^2(^{(\cc_1)}\Si_{t_1})} & \les D\varep_1,
\end{align}
for $i=1\ldots 3$. Let define the Cartesian vectorfields $X_i$ by parallel transport of $\nab x^i$, \emph{i.e.}
\begin{align}\label{eq:transportXiKill}
  \begin{aligned}
  \D_{^{(\cc_1)}\Tf}(X_i) & = 0,\\
  X_i|_{^{(\cc_1)}\Si_{t_1}} & = \nab x^i,
  \end{aligned}
\end{align}
for $i = 1 \ldots 3$. Commuting equation~\eqref{eq:transportXiKill} by $\D^{\leq 3}$ and integrating in $^{(\cc_1)}\Tf$, using the Bootstrap Assumptions~\ref{BA:curvint} and~\ref{BA:connint} in $^{(\cc_1)}\MM^\intr_{\bott}$, estimate~\eqref{est:harmocoordsaux} on $^{(\cc_1)}\Si_{t_1}$ and the Sobolev estimates of Lemma~\ref{lem:SobSitHHrr}, we have the following bounds in (the rescaled) region $\DD$
\begin{align}\label{est:transportXiKill}
  \norm{\D X_i}_{L^\infty(\DD)} & \les D\varep_1, & \norm{\D^2X_i}_{L^\infty_t L^6} & \les D\varep_1, & \norm{\D^3X_i}_{L^2(\DD)} & \les D\varep_1.
\end{align}
Thus, one can apply energy estimates (see Section~\ref{sec:errorintr} or Section~\ref{sec:tipcurvest}) in the rescaled domain $\DD$ for the following contracted and commuted Bel-Robinson tensors
\begin{align*}
  Q(\Lieh_{X_{\mu}}\R)(X_0,X_0,X_0) , \quad\quad Q(\Lieh_{X_\mu}\Lieh_{X_\nu}\R)(X_0,X_0,X_0)
\end{align*}
for $\mu,\nu = 0 \ldots 3$ and where $X_0 := ^{(\cc_1)}\Tf$, and we obtain
\begin{align}\label{est:enerestXiKill}
  \begin{aligned}
    & \; \int_{^{(\cc)}\Si_t} Q(\Lieh_{X_\mu}^{\leq 1}\Lieh_{X_\nu}\R)(X_0,X_0,X_0,{^{(\cc)}\Tf})  + \int_{\CC_u\cap\DD}Q(\Lieh_{X_\mu}^{\leq 1}\Lieh_{X_\nu}\R)(X_0,X_0,X_0,\el) \\
    \les & \; \int_{^{(\cc_1)}\Si_{t_1}} Q(\Lieh_{X_\mu}^{\leq 1}\Lieh_{X_\nu}\R)(X_0,X_0,X_0,{^{(\cc_1)}\Tf}) + (D\varep_1)^3 \\
  \les  & \; \varep_1,
  \end{aligned}
\end{align}
where we used the curvature bounds~\eqref{est:sourcecc1intr} and~\eqref{est:transportXiKill}.\\

Performing a frame comparison argument as in Section~\ref{sec:tipcurvest} and scaling back in $t_1$, we deduce from~\eqref{est:enerestXiKill} that estimate~\eqref{est:proofproptrans1} holds. This finishes the proof of Lemma~\ref{lem:proofproptrans1}.

\subsection{Proof of Lemma~\ref{lem:proofproptrans2}}
In this section, we assume that $\cco \leq \cc \leq \cc_1$.\\

\begin{figure}[h!]
  \centering
  \includegraphics[height=6cm]{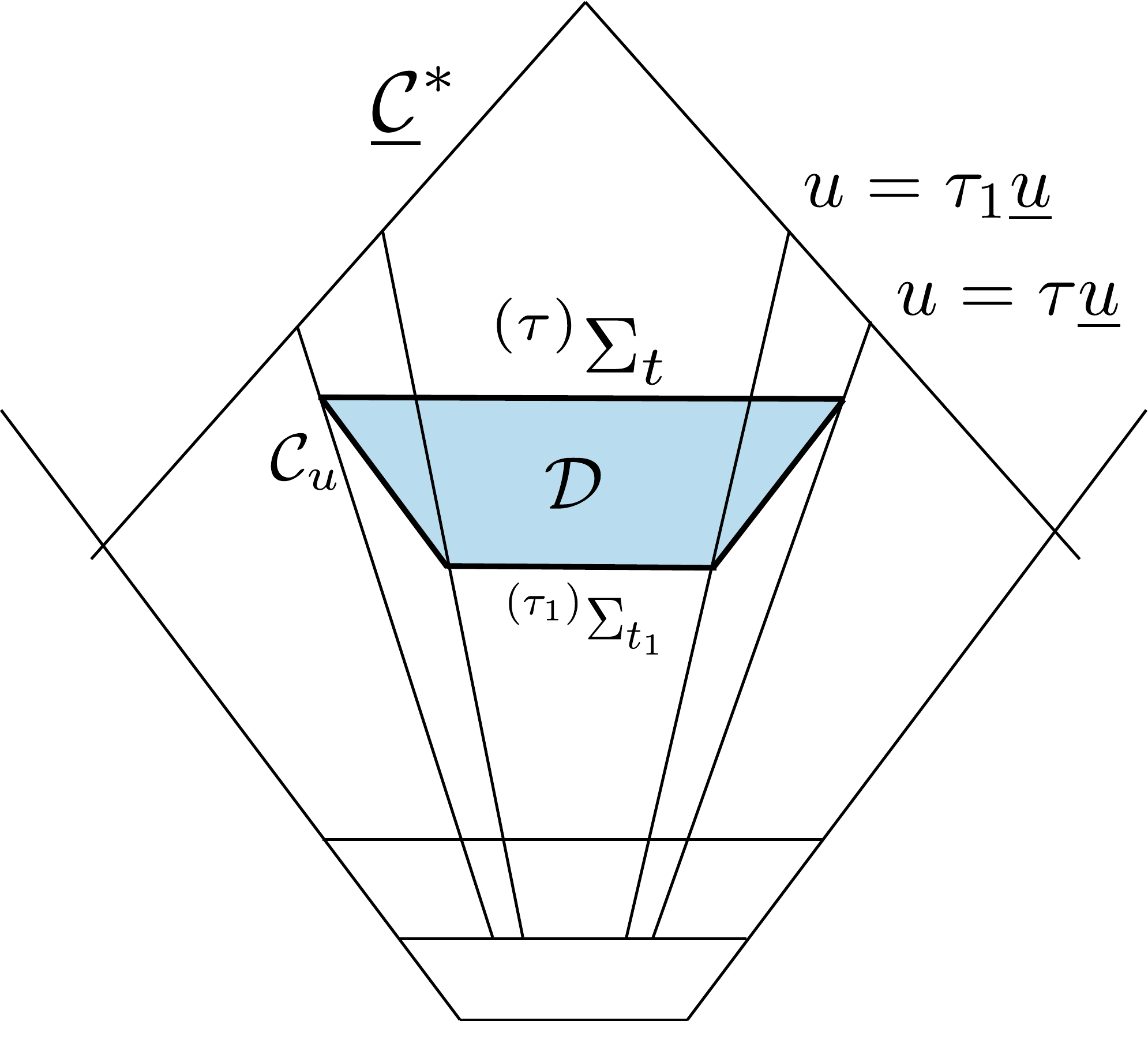}
  \caption{Local energy estimates for $\cc \leq \cc_1$.}
  \label{fig:localenergccllcc1}
\end{figure}

Let $^{(\cc)}\too \leq t \leq ^{(\cc)}\tast$. Let
\begin{align*}
  u & := \frac{2}{1+\cc^{-1}}t, & t_1 & := \frac{1+\cc_1^{-1}}{2}u.
\end{align*}

\begin{remark}
  In the case $t \leq 3\frac{1+\cco^{-1}}{1+\cc_1^{-1}}$, the proof of estimate~\eqref{est:proofproptrans2} follows from energy estimates in a bottom initial layer\footnote{That bottom initial layer is potentially larger than the initial layer used in this paper.} and is left to the reader (see similar estimates in Section~\ref{sec:initlayerenergy}). We shall thus consider that $t_1 \geq 3$ (that is, the hypersurface $^{(\cc_1)}\Si_{t_1}$ is well defined).
\end{remark}

We consider the maximal hypersurfaces $^{(\cc_1)}\Si_{t_1}$ and $^{(\cc)}\Si_{t}$ (see Figure~\ref{fig:localenergccllcc1}). Exchanging their roles in the proof of Lemma~\ref{lem:proofproptrans1} in Section~\ref{sec:proofproptrans1}, defining harmonic coordinates on the hypersurface $^{(\cc)}\Si_t$ and running the same procedure yields estimate~\eqref{est:proofproptrans2}, where we use the curvature bounds~(\ref{est:sourcecc1intr}) through $^{(\cc_1)}\Si_{t_1}$ and the curvature bounds~(\ref{est:sourcecc1CCu}) through $\CC_u$. This finishes the proof of Lemma~\ref{lem:proofproptrans2}.

\subsection{Proof of Lemma~\ref{lem:proofproptrans3}}
In this section, we assume that $\cc_1 \leq \cc \leq (1+\cco)/2$.\\

\begin{figure}[h!]
  \centering
  \includegraphics[height=6cm]{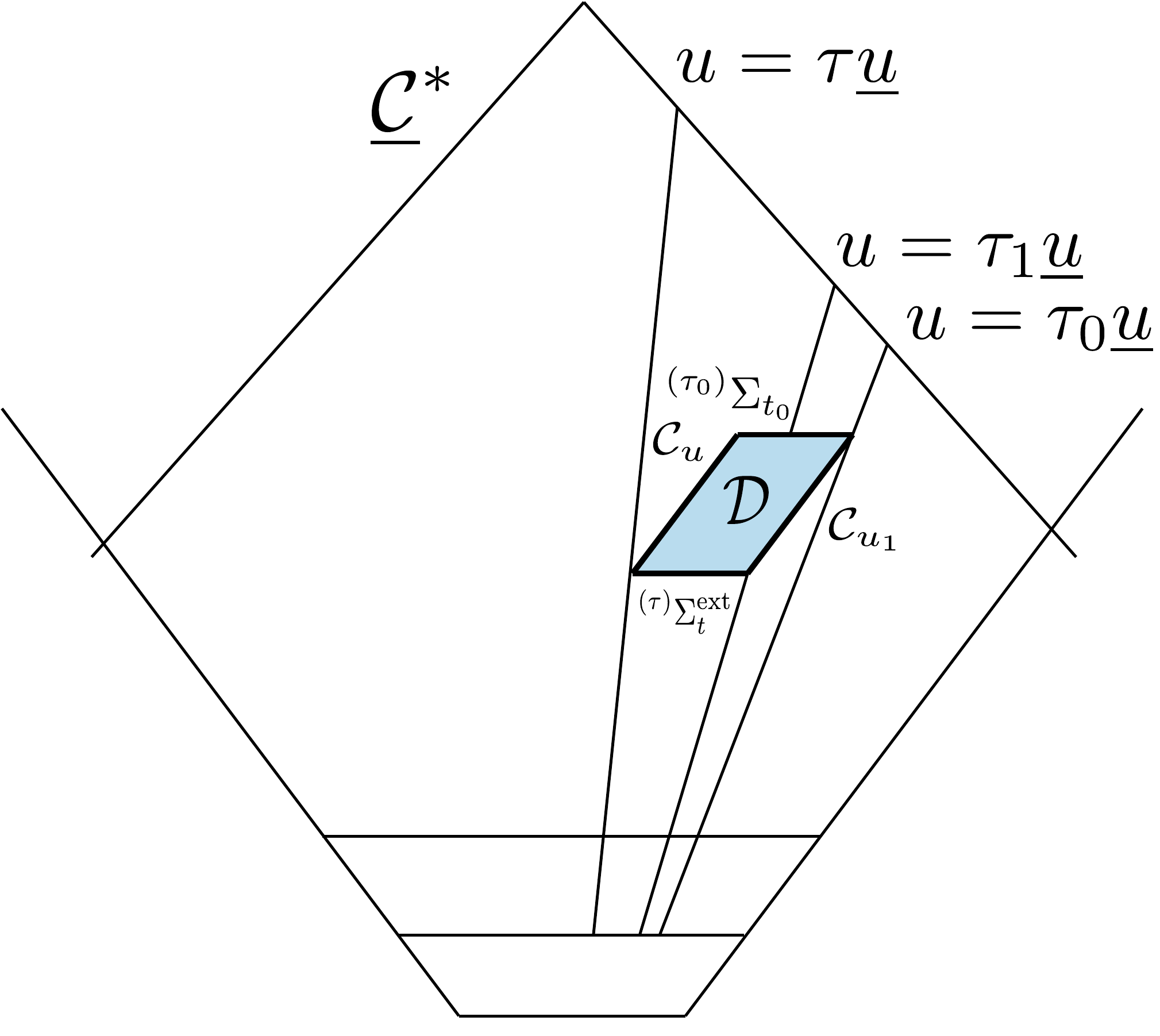}
  \caption{Local energy estimates for $^{(\cc)}\Si_t^{\ext}$.}
  \label{fig:localenergcc1Sitext}
\end{figure}

Let $^{(\cc)}\too \leq t \leq ^{(\cc)}\tast$. Let
\begin{align*}
  u & := \frac{2}{1+\cc}t, & u_1 & := \frac{2}{1+\cc_1}t, & t_0 & := \frac{1+\cco}{2}u_1.
\end{align*}
Let denote by $\DD$ the domain enclosed by $^{(\cc)}\Si_t^\ext$, $\CC_{u_1}$, $^{(\cco)}\Si_{t_0}$ and $\CC_u$ (see Figure~\ref{fig:localenergcc1Sitext}). 

\begin{remark}
  If $t_0\geq ^{(\cco)}\tast$, \emph{i.e.} $^{(\cc)}\Si^\ext_t$ is in the domain of dependence of the last slice $^{(\cco)}\Si_{^{(\cco)}\tast}$, estimates follow along the arguments of Section~\ref{sec:tipcurvest}. In the following we assume that $t_0\leq ^{(\cco)}\tast$.
\end{remark}

Defining harmonic coordinates on the full slice $^{(\cco)}\Si_{t_0}$, and arguing as in Section~\ref{sec:proofproptrans1}, one can perform energy estimates in the (rescaled) region $\DD$, using the commutator vectorfields given by parallel transport of the harmonic vectorfields on $^{(\cco)}\Si_{t_0}$. Estimate~\eqref{est:proofproptrans3} thus follows from the energy bound~(\ref{est:proofproptrans2}) through $^{(\cco)}\Si_{t_0}$ from Lemma~\ref{lem:proofproptrans2} and the energy bound~(\ref{est:sourcecc1CCu}) and (\ref{est:proofproptrans1}) on $\CC_u$, and a comparison argument. This finishes the proof of Lemma~\ref{lem:proofproptrans3} and of Proposition~\ref{prop:alltransparamSTAB}.

\chapter{Null connection estimates on $\protect\CCba$}\label{sec:connestCCba}
In this section, we prove the following proposition.
\begin{proposition}\label{prop:connestCCbaSTAB}
  Recall that from Propositions~\ref{prop:curvestSTAB} and \ref{prop:tipcurvestSTAB}, the following estimates hold for the curvature boundedness norms on $\CCba$ (see the definitions of Section~\ref{sec:normnullcurv}) 
  \begin{align}\label{est:sourcecurvast}
    \RRast + \RRfast & \les \varep.
  \end{align} 
  
  Under the Bootstrap Assumptions, estimates~\eqref{est:sourcecurvast}, the limits for the null connection coefficients at the vertex of the cone $\CCba$ from Theorems~\ref{thm:vertex} and~\ref{thm:canonical}, we have for $\varep>0$ sufficiently small (see the definitions of Sections~\ref{sec:normnullcurv} and~\ref{sec:normnullconnCCba})
  \begin{align}\label{est:OOb3ast}
    \RRfoast + \OOast + \OOoast + \OO_{\leq 3}^{\ast,\OOO} &\les \varep,
  \end{align}
  and
  \begin{align}\label{est:OOfb2ast}
    \OOfast &\les \varep.
  \end{align}
  Moreover, the Bootstrap Assumptions~\ref{BA:mildsphcoordsast} for the spherical coordinates on $\CCba$ (see Section~\ref{sec:mildBAimpast}), the bootstrap bound~(\ref{est:BAarearadiusestimateCCba}) on the area radius (see Section~\ref{sec:arearadiusestimateCCba}), and the mild Bootstrap Assumptions~\ref{BA:mildOOE} for the rotation vectorfields $\OOE$ on $\CCba\cap\MM^\ext$ (see Section~\ref{sec:mildOOE}) are improved.\footnote{This amounts to an improvement of the strong Bootstrap Assumptions~\ref{BA:connCCba}.}
\end{proposition}





\begin{remark}
  The bound~\eqref{est:OOfb2ast} will follow from~\eqref{est:OOb3ast} and the Klainerman-Sobolev estimates of Lemma~\ref{lem:KlSobcomplet}.
\end{remark}

\begin{remark}\label{rem:nullconndegeaxis}
  The norms $\OOfast$ (see the definitions of Section~\ref{sec:normnullconnCCba}) allow for a degeneracy of the null connection coefficients of the type $\Gac \sim r^{-1/2}$ when $r \to 0$ on the cone $\CCba$ , where $\Gac \in\{\trchi-\frac{2}{r},\trchib+\frac{2}{r},\chih,\chibh,\ze,\omb\}$. This is inherited from the degeneracy in $r$ of the null curvature components (see Remark~\ref{rem:notoptimalaxisdege}). These degeneracies in $r$ still allow to close the estimates for the null connection coefficients. This is consistent with the low regularity control for the geodesic foliation on a null cone obtained in~\cite{Wan09}. In that article, it is namely proven that one can close (low regularity) estimates for the null connection coefficients of the geodesic foliation, only assuming an $L^2$ control of the null curvature components on the cone. In the present paper, such an $L^2$ control holds and is consistent with a potential $r^{-3/2}$ singularity for the curvature at the axis, and we therefore expect that the $r$ degeneracy is not an obstruction to close our estimates.\footnote{In view of the nonlinearities in the null structure equations, the critical degeneracy at which one cannot close (degenerate) estimates for the null connection coefficients is given by the (scaling-critical) blow-up rate $R\sim r^{-2}$.}  
\end{remark}


\section{Hardy estimates}
We first state the following preliminary general Hardy estimates, which hold on the whole cone $\CCba$. We recall that $\qq :=  \min(r,u)$. 
\begin{lemma}[Hardy estimates on $\CCba$]\label{lem:transportast}
  For all $\kappa\in\RRR$, the following holds. Assume that $U$ is an $S$-tangent tensor satisfying
  \begin{align*}
    \Nd_3U + \frac{\kappa}{2} \trchib U & = F,
  \end{align*}
  with vertex limit
  \begin{align*}
    r^\kappa |U| \to 0, \quad \text{when $r\to 0$.}
  \end{align*}
  We have the following Hardy inequality on $\CCba$
  \begin{align}\label{est:Hardyast}
    \norm{r^\lambda U}_{L^2(\CCba)} & \les \norm{r^{\lambda}\qq F}_{L^2(\CCba)}.
  \end{align}
  for all $\lambda < \kappa - 3/2$.
\end{lemma}
\begin{remark}
  In the following, we will use the Hardy inequality, with $\lambda = -1$ for the null connection coefficients $\trchi,\trchib$, which satisfy transport equations with $\kappa\geq 1$. See the null structure equations~(\ref{eq:null}).\footnote{The fact that the null connection coefficients intrinsic to the geometry of the cones $\chi,\chib,\ze$ satisfy transport equations with $\kappa \geq 1$ 
    enables the control of the intrinsic geometry of null cones at the $L^2$ level in~\cite{Wan09}.} 
\end{remark}
\begin{proof}
  

  We first prove the following two classical Hardy estimates corresponding respectively to the interior and exterior region
  \begin{align}
    \int_{\cc\uba}^\uba \le(\frac{1}{(\uba-u)^{\al}}\int_u^\uba f \,\d u'\ri)^2 (\uba-u)^2 \,\d u & \les \int_{\cc\uba}^\uba (\uba-u)^{2-2\al} |f|^2 (\uba-u)^2 \, \d u, \label{est:Hardyclassint}\\
        \int_{1}^{u_0}\le(\int_{u}^{u_0}f \,\d u'\ri)^2 \d u & \les \int^{u_0}_1 u^2|f|^2 \,\d u,\label{est:Hardyclassext} 
  \end{align}
  where $f$ is a scalar function of the variable $u$, and where $\al>3/2$. Here $u_0>0$ is a parameter which will be determined in the sequel.\\

  To prove~\eqref{est:Hardyclassint}, using Cauchy-Schwartz we write
  \begin{align*}
     & \; \int_{\cc\uba}^\uba\le(\int_u^\uba f \,\d u'\ri)^2 (\uba-u)^{2-2\al}\,\d u \\
    \leq & \; \int_{\cc\uba}^\uba\bigg(\int_{u}^\uba (\uba-u')^{2\ga}|f|^2\,\d u'\bigg)\bigg( \int_{u}^\uba (\uba-u')^{-2\ga} \,\d u'\bigg) (\uba-u)^{2-2\al} \,\d u,
  \end{align*}
  where $\ga\in\RRR$. Provided that $-2\ga + 1 >0$, we have
  \begin{align*}
    \int_{u}^\uba (\uba-u')^{-2\ga}\,\d u' & \les (\uba-u)^{-2\ga+1}.
  \end{align*}

  Using Fubini theorem, and provided that $-2\ga+4-2\al <0$, we obtain
  \begin{align*}
    & \; \int_{\cc\uba}^\uba\le(\int_u^\uba f \,\d u'\ri)^2 (\uba-u)^{2-2\al} \,\d u \\
    \les & \; \int_{\cc\uba}^\uba(\uba-u')^{2\ga} |f|^2 \le(\int_{\cc\uba}^{u'} (\uba-u)^{-2\ga+1+2-2\al} \,\d u\ri) \,\d u' \\
    \les & \; \int_{\cc\uba}^\uba(\uba-u')^{2\ga} |f|^2 (\uba-u')^{-2\ga+4-2\al} \,\d u' \\
    = & \; \int_{\cc\uba}^\uba(\uba-u')^{2-2\al}|f|^2 (\uba-u')^{2}\,\d u'.
  \end{align*}
  The conditions $-2\ga+1 >0$ and $-2\ga+4-2\al <0$ can both be satisfied if and only if $\al > 3/2$, which finishes the proof of~\eqref{est:Hardyclassint}.\\

  To prove~\eqref{est:Hardyclassext}, we write using Cauchy-Schwartz
  \begin{align*}
    & \; \int_{1}^{u_0}\le(\int_{u}^{u_0}f \, \d u' \ri)^2\,\d u \\
    \leq & \; \int_{1}^{u_0}\bigg(\int_{u}^{u_0}(u')^{2\ga} |f|^2\,\d u'\bigg)\bigg(\int_{u}^{u_0}(u')^{-2\ga}\,\d u'\bigg) \, \d u \\
    \les & \; \int_{1}^{u_0}\bigg(\int_{u}^{u_0}(u')^{2\ga} |f|^2\,\d u'\bigg) u^{-2\ga+1} \, \d u,
  \end{align*}
  provided that $-2\ga+1<0$. Using Fubini theorem, we thus have
  \begin{align*}
    & \; \int_{1}^{u_0}\le(\int_{u}^{u_0}f \, \d u' \ri)^2\,\d u \\
    \les & \; \int_{1}^{u_0}(u')^{2\ga} |f|^2 \le(\int_{1}^{u'} u^{-2\ga+1}\,\d u\ri)\,\d u' \\
    \les & \; \int_{1}^{u_0}(u')^{2\ga} |f|^2 (u')^{-2\ga+2} \, \d u',
  \end{align*}
  provided that $-2\ga+2 >0$. Thus, choosing $\ga = 3/4$, we obtained~\eqref{est:Hardyclassext}.\\

  Using formula~\eqref{eq:elbr}, we have
  \begin{align}\label{eq:transportrkappaUcorrec}
    \Nd_3(r^\kappa U) & = r^\kappa F + \frac{\kappa}{2}\le(\trchib-\trchibo\ri)r^\kappa U.
  \end{align}
  Using the vertex limit $r^\kappa |U| \to 0$, the Bootstrap Assumptions~\ref{BA:connCCba} on $\trchib$ and a Gr\"onwall argument, we deduce
  \begin{align}\label{est:rkappaU}
    \norm{r^{-1}r^{\kappa}|U|}_{L^2(S_{u,\uba})} & \les \int_{u}^\uba \norm{r^{-1}r^\kappa|F|}_{L^2(S_{u',\uba})} \, \d u'.
  \end{align}
  

  Applying the classical Hardy estimate~\eqref{est:Hardyclassint} with $f(u) = \norm{r^{-1}r^{\kappa}F}_{L^2(S_{u,\uba})}$, using the Bootstrap Assumptions~\ref{BA:mildsphcoordsast} on the metric in spherical coordinates and the bootstrap bounds~(\ref{est:BAarearadiusestimateCCba}) on the area radius, we deduce from~\eqref{est:rkappaU} in the interior region
  
  \begin{align*}
    \norm{r^{\lambda}U}_{L^2(\CCba\cap\MM^\intr)} & = \norm{r^{-\alpha}r^{\kappa}U}_{L^2(\CCba\cap\MM^\intr)} \les \norm{r^{1-\al}r^{\kappa}F}_{L^2(\CCba\cap\MM^\intr)} = \norm{r^{\lambda+1}F}_{L^2(\CCba\cap\MM^\intr)},
  \end{align*}
  where $\la := \kappa-\alpha$ and $\alpha>3/2$. Rewriting this estimate, using that $r\simeq \qq$ in the interior region, we thus have  
  \begin{align}\label{est:Hardyint}
    \begin{aligned}
    \norm{r^\lambda U}_{L^2(\CCba\cap\MM^\intr)} & \les \norm{r^{\lambda+1}F}_{L^2(\CCba\cap\MM^\intr)} \les \norm{r^\lambda \qq F}_{L^2(\CCba)},
                                                                     \end{aligned}
  \end{align}
  for all $\lambda < \kappa -3/2$.\\
  
  From a mean value argument, there exists $\cc\uba \leq u_0\leq (1+\cc)\uba/2$ such that
  \begin{align}\label{est:rlambdaUu0}
    \begin{aligned}
    \norm{r^\lambda U}_{L^2(S_{u_0,\uba})} & \les \uba^{-1/2}\norm{r^\lambda U}_{L^2(\CCba\cap\{\cc\uba\leq u_0\leq (1+\cc)\uba/2\})} \les \uba^{-1/2}\norm{r^\lambda U}_{L^2(\CCba\cap\MM^\intr)}.
  \end{aligned}
  \end{align}
  Integrating equation~\eqref{eq:transportrkappaUcorrec} and using a Gr\"onwall argument, we have for $\la < \kappa -3/2$ and for $u \leq u_0$
  \begin{align*}
    \norm{r^{-1}r^\lambda U}_{L^2(S_{u,\uba})} & \les \norm{r^{-1}\le(r^\lambda U\ri)}_{L^2(S_{u_0,\uba})} + \int_{u}^{u_0}\norm{r^{-1}r^\lambda F}_{L^2(S_{u',\uba})}\,\d u'.
  \end{align*}
  Thus,
  \begin{align*}
    \norm{r^\lambda U}_{L^2(\CCba\cap\MM^\ext)} & \les \norm{(r^\lambda U)|_{u=u_0}}_{L^2(\CCba\cap\MM^\ext)} + \norm{\int_u^{u_0}\norm{r^{-1}r^\lambda F}_{L^2(S_{u',\uba})}\,\d u'}_{L^2(\CCba\cap\MM^\ext)}.
  \end{align*}
  
  Using~\eqref{est:rlambdaUu0} and~\eqref{est:Hardyint}, we have for the first term
  \begin{align*}
    \norm{(r^\lambda U)|_{u=u_0}}_{L^2(\CCba\cap\MM^\ext)} & \les \uba^{1/2} \norm{r^\lambda U}_{L^2(S_{u_0,\uba})} \\
                                                                  & \les \norm{r^\lambda U}_{L^2(\CCba\cap\MM^\intr)} \\
                                                                  & \les \norm{r^\lambda \qq F}_{L^2(\CCba)}.
  \end{align*}
  For the second term, applying the classical Hardy estimate~\eqref{est:Hardyclassext} to $f(u) = \norm{r^{-1}r^{\lambda}|F|}_{L^2(S_{u,\uba})}$, using the Bootstrap Assumption~\ref{BA:mildsphcoordsast} on the metric in spherical coordinates and the bootstrap bounds~(\ref{est:BAarearadiusestimateCCba}) on the area radius, gives
  \begin{align*}
    \norm{\int_u^{u_0}\norm{r^{-1}r^\lambda F}_{L^2(S_{u',\uba})}\,\d u'}_{L^2(\CCba\cap\MM^\ext)} & \les \norm{r^\lambda \qq F}_{L^2(\CCba)},
  \end{align*} 
  where we used that $\qq \simeq u$ in the exterior region and $\qq \simeq r$ in the interior region. This finishes the proof of the lemma.
\end{proof}


\section{Klainerman-Sobolev estimates}
We have the following Klainerman-Sobolev estimates, which are used to obtain estimates~(\ref{est:OOfb2ast}). Its proof is postponed to Appendix~\ref{app:KlSobH12}.
\begin{lemma}[Klainerman-Sobolev estimates on $\CCba$]\label{lem:KlSobcomplet}
  For all $S$-tangent tensor $F$ on $\CCba$ with vertex limit
  \begin{align*}
    r^{3/2}|(r\Nd)^{\leq 1}F|\to 0 \quad \text{when $r\to 0$,} 
  \end{align*}
  we have the following $L^\infty_{u}\HHt\le(S_{u,\uba}\ri)$ estimates in $\CCba$
  \begin{align*}
    \norm{r F}_{L^\infty_{u}\HHt\le(S_{u,\uba}\ri)} & \les \norm{F}_{L^2(\CCba)} + \norm{r\Nd F}_{L^2(\CCba)} +\norm{r\Nd_3F}_{L^2(\CCba)},
  \end{align*}
  and
  \begin{align*}
    \norm{r^{1/2}\qq^{1/2}F}_{L^\infty_{u}\HHt\le(S_{u,\uba}\ri)} & \les \norm{F}_{L^2(\CCba)} + \norm{r\Nd F}_{L^2(\CCba)} + \norm{\qq\Nd_3F}_{L^2(\CCba)}.
  \end{align*}
\end{lemma}

\section{Control of $\rhoo$ and $\sigmao$}\label{sec:controlrhoosigamooCCba}
From the vertex limits of Theorems~\ref{thm:vertex} and~\ref{thm:canonical}, we have $r^3\rhoo \to 0$ when $r\to 0$. Integrating equation~\eqref{eq:Nd3rhoo} thus gives 
\begin{align*}
  r^3|\rhoo| & \les \int_{u}^\uba r^3|\Err\le(\Nd_3,\rhoo\ri)|\,\d u'.
\end{align*}
From the Bootstrap Assumptions~\ref{BA:connCCba} and~\eqref{est:sourcecurvast}, we have
\begin{align}\label{est:Err3rhoo}
  \begin{aligned}
  \le|\Err\le(\Nd_3,\rhoo\ri)\ri| & \les \le|\chih\ri|\le|\alb\ri| + \le|\ze\ri|\le|\beb\ri| + \le|\trchib-\trchibo\ri|\le|\rho-\rhoo\ri| \\
  & \les (D\varep)\qq^{-1/2}\ub^{-2} \cdot \varep r^{-1}\qq^{-1/2}u^{-2}.
  \end{aligned}
\end{align}
Thus,
\begin{align*}
  r^3|\rhoo| & \les D\varep^2\uba^{-2}\int_u^\uba r(u',\uba)^2\qq(u',\uba)^{-1}(u')^{-2}\,\d u'. 
\end{align*}
In the interior region $u\geq \cc\uba$, this yields
\begin{align*}
  r^3|\rhoo| & \les D\varep^2\uba^{-4}\int_u^\uba r(u',\uba)\,\d u' \les D\varep^2\uba^{-4}r^2.
\end{align*}
In the exterior region $u\leq \cc\uba$, this yields
\begin{align*}
  r^3|\rhoo| & \les \uba^3|\rhoo||_{S^\ast} + D\varep^2\int_{u}^\uba (u')^{-3}\,\d u'\\
             & \les D\varep^2u^{-2}.
\end{align*}
Thus, we summarise the bounds on $\rhoo$ as
\begin{align}\label{est:rhooLinfcorrec}
  \norm{u\qq\uba^3\rhoo}_{L^\infty(\CCba)} & \les \varep.
\end{align}
Estimating directly the $(\qq\Nd_3)$ derivative using equation~\eqref{eq:Nd3rhoo}, and estimates~\eqref{est:Err3rhoo} and~\eqref{est:rhooLinfcorrec}, we have
\begin{align*}
  \le|\qq\Nd_3\rhoo\ri| & \les |\rhoo| + \qq\le|\Err(\Nd_3,\rhoo)\ri| \les \varep u^{-1}\qq^{-1}\uba^{-3}.
\end{align*}
Deriving by $(\qq\Nd_3)$ in equation~\eqref{eq:Nd3rhoo} and directly estimating, we have
\begin{align*}
  \le|(\qq\Nd_3)^2\rhoo\ri| & \les \le|(\qq\Nd_3)^{\leq 1}\rhoo\ri| + \le|(\qq\Nd_3)^{\leq 1}\Err(\Nd_3,\rhoo)\ri|.
\end{align*}
For the last term, we only treat the term
\begin{align*}
  \le|(\qq\Nd_3) \overline{\chih\cdot\alb}\ri| & \les \le|\overline{(\qq\Nd_3)\chih\cdot\alb}\ri| + \le|\overline{\chih\cdot(\qq\Nd_3)\alb}\ri| + \le|\overline{\chih\cdot\alb \qq(\trchib-\trchibo)}\ri| \\
                                               & \les r^{-2}\norm{(\qq\Nd_3)^{\leq 1}\chih}_{L^2(S_{u,\uba})}\norm{(\qq\Nd_3)^{\leq 1}\alb}_{L^2(S_{u,\uba})} \\
                                               & \les (D\varep)^2u^{-1}\qq^{-1}\uba^{-3}.
\end{align*}
where we used formula~(\ref{eq:commelbov}), Cauchy-Schwartz, and the Bootstrap Assumptions~\ref{BA:curvast} and~\ref{BA:connCCba} together with the Sobolev estimates of Lemma~\ref{lem:Sobsphere}.\\

Summarising the above estimates, we obtained
\begin{align}\label{est:Linfrhooast}
  \RRfoast[\rhoo] & \les \varep.
\end{align}

Based on the analysis performed for $\rhoo$, we make the following general remark.
\begin{remark}[Higher order estimates]\label{rem:DecayandRegularity}
  Deriving/commuting an equation by $\Ndt \in\le\{r\Nd,\qq\Nd_3,\ub\Nd_4\ri\}$ \emph{does not change/worsen the decay} of the terms in that equation. 
  Thus, to obtain higher order estimates, one does only need to check that the terms in the derived equation are controlled \emph{from a regularity point of view}, that is, \emph{without having to explicit the estimates in terms of $u,\ub$ weights.}
\end{remark}

In the sequel, we shall use the above remark to infer control of higher order derivatives only based on a regularity check of the terms involved. To illustrate how this proceeds, we detail the simple case of $\sigmao$.\\

Taking the sup-norm in equation~\eqref{eq:sigmao}, using the Bootstrap Assumptions~\ref{BA:connCCba}, we first have the estimate for the non-derived quantity
\begin{align*}
  \le|\sigmao\ri| & \les \le|\chih\ri| \le|\chibh\ri| \les \le(D\varep \qq^{-1/2}\ub^{-2}\ri)\le(\qq^{-1/2}u^{-1}\ub^{-1}\ri) \les \varep \qq^{-1}u^{-1}\ub^{-3}.
\end{align*}

Commuting equation~\eqref{eq:sigmao} with $\qq\Nd_3$, then gives \emph{directly the same estimate as above} for $\qq\Nd_3\sigmao$ since \emph{from the regularity point of view}, by a simple check of the Bootstrap Assumptions~\ref{BA:connCCba}, we have that the derivatives of $\chih,\chibh$ are controlled in $L^\infty$. We thus have the following estimate (which can be independently checked)
\begin{align*}
  \le|(\qq\Nd_3)\sigmao\ri| & \les \varep \qq^{-1}u^{-1}\ub^{-3}.
\end{align*}

This observation still holds when deriving~\eqref{eq:sigmao} with $(\qq\Nd_3)^2$. In fact, in that case, from the regularity point of view, we only have that two derivatives of $\chih,\chibh$ are controlled in $L^\infty_uL^4(S_{u,\uba})$. This is however still enough to obtain the $L^\infty$ estimate for $(\qq\Nd_3)^2\sigmao$, using Cauchy-Schwartz as in the case of $(\qq\Nd_3)^2\rhoo$. Here again -- as for the estimate for $(\qq\Nd_3)\sigmao$ --, we do not need to check the actual weigths in $u,\ub$, and we have the following estimate (which can be independently checked)
\begin{align*}
  \le|(\qq\Nd_3)^{2}\sigmao\ri| & \les \varep \qq^{-1}u^{-1}\ub^{-3}.
\end{align*}

Thus, we have proved
\begin{align}\label{est:Linfsigmaoast}
  \RRfoast[\sigmao] & \les \varep,
\end{align}
which together with the estimates for $\rhoo$ gives the desired estimate
\begin{align*}
  \RRfoast & \les \varep.
\end{align*}

\section{Control of $\protect\trchibo+\frac{2}{r}$}
Integrating the transport equation~\eqref{eq:Nd3trchibo} for $\trchibo+\frac{2}{r}$, using that from Theorems~\ref{thm:vertex} and~\ref{thm:canonical}, one has the limit $r\le(\trchibo+\frac{2}{r}\ri) \to 0$ when $r\to 0$, we have
\begin{align*}
  \le|r\le(\trchibo+\frac{2}{r}\ri)\ri| & \les \int^\uba_{u}r\le|\Err\le(\Nd_3,\trchibo+\frac{2}{r}\ri)\ri| \, \d u'.
\end{align*}
From the error terms, we only treat the leading term $|\chibh|^2$ for which we have
\begin{align*}
  \int_u^\uba r |\chibh|^2 \,\d u' & \les (D\varep)^2\int_u^\uba r\uba^{-2}u^{-2}\qq^{-1}\,\d u' \\
                                   & \les \varep r \uba^{-2}u^{-2}.
\end{align*}
Thus, we have
\begin{align*}
  \norm{\uba^2u^{2}\le(\trchibo+\frac{2}{r}\ri)}_{L^\infty} & \les \varep.
\end{align*}

Estimating directly equation~\eqref{eq:Nd3trchibo} for $\Nd_3\le(\trchibo+\frac{2}{r}\ri)$, one further infers
\begin{align*}
  \le|\qq\Nd_3\le(\trchibo+\frac{2}{r}\ri)\ri| & \les \le|\trchibo+\frac{2}{r}\ri| + \qq\le|\Err\le(\Nd_3,\trchibo+\frac{2}{r}\ri)\ri| \\
  & \les \varep \uba^{-2}u^{-2}.
\end{align*}

Deriving equation~\eqref{eq:Nd3trchibo} by $\qq\Nd_3$, -- checking from the Bootstrap Assumptions~\ref{BA:connCCba} that one has a control for one derivative of the error term $\Err\le(\Nd_3,\trchibo+\frac{2}{r}\ri)$ -- also gives
\begin{align*}
  \le|(\qq\Nd_3)^2\le(\trchibo+\frac{2}{r}\ri)\ri| & \les \varep \uba^{-2}u^{-2}.
\end{align*}

Thus, we have proved
\begin{align}\label{est:Linftrchiboast}
  \OOoast\le[\trchibo\ri] & \les \varep.
\end{align} 

\section{Control of $\trchio-\frac{2}{r}$}
Integrating the transport equation~\eqref{eq:Nd3trchio} for $\trchio-\frac{2}{r}$, using that from Theorems~\ref{thm:vertex} and~\ref{thm:canonical}, one has the limit $r\le(\trchio-\frac{2}{r}\ri)\to 0$ when $r\to 0$, one has
\begin{align*}
  \le|r\le(\trchio-\frac{2}{r}\ri)\ri| & \les \int_{u}^\uba r\le(|\rhoo|+ \le|\Err\le(\Nd_3,\trchio-\frac{2}{r}\ri)\ri|\ri)\,\d u.
\end{align*}
From the bound~(\ref{est:Linfrhooast}) obtained for $\rhoo$, and an inspection of the error term, we have
\begin{align*}
  r\le(|\rhoo|+ \le|\Err\le(\Nd_3,\trchio-\frac{2}{r}\ri)\ri|\ri) & \les \varep r u^{-1} \qq^{-1} \uba^{-3} + (D\varep)^2 r \uba^{-3}u^{-1}\qq^{-1} \\
  & \les \varep u^{-2}\uba^{-2}.
\end{align*}
Thus we have
\begin{align*}
  \le|\trchio-\frac{2}{r}\ri| & \les \varep r^{-1} \int_u^\uba (u')^{-2}\uba^{-2} \,\d u' \\
                      & \les \varep \uba^{-3}u^{-1},
\end{align*}
that is,
\begin{align*}
  \norm{u\uba^3\le(\trchio-\frac{2}{r}\ri)}_{L^\infty(\CCba)} & \les \varep.
\end{align*}
Estimating directly equation~\eqref{eq:Nd3trchio} for $\Nd_3\le(\trchio-\frac{2}{r}\ri)$, and also taking further a $(\qq\Nd_3)$ derivative in~\eqref{eq:Nd3trchio}, we infer
\begin{align}\label{est:Linftrchioast}
  \OOoast\le[\trchio\ri] & \les \varep.
\end{align}

This finishes the improvement of
\begin{align*}
  \OOoast & \les \varep.
\end{align*}

\section{Control of $\ze$}
Using the definition~(\ref{eq:can}) of the canonical foliation on $\CCba$, $\ze$ satisfies the following elliptic equation
\begin{align}\label{eq:Dd1zeCCba}
  \Dd_1\ze & = \le(-\rho+\rhoo, \sigma-\half\chih\wedge\chibh\ri).
\end{align}

Using the $L^2$ elliptic estimates of Lemma~\ref{lem:ell} with $\ell=0$, the estimates~(\ref{est:sourcecurvast}) for the curvature and the Bootstrap Assumptions~\ref{BA:connCCba}, one has
\begin{align*}
  \norm{r^{-1}(r\Nd)^{\leq 1}\ze}_{L^2(\CCba)} & \les \norm{\rho-\rhoo}_{L^2(\CCba)} + \norm{\sigma-\sigmao}_{L^2(\CCba)} + \norm{u^{-1}\chih}_{L^\infty}\norm{u\chibh}_{L^2(\CCba)} \\
                                               & \les \varep \uba^{-2} + \le(D\varep \uba^{-2}\ri) \le(D\varep\ri)\\
                                                 & \les \varep \uba^{-2}.
\end{align*}

Using the $L^2$ elliptic estimates of Lemma~\ref{lem:ell} with $\ell=2$, we have
\begin{align*}
  \norm{r^{-1}(r\Nd)^{\leq 3}\ze}_{L^2(\CCba)} & \les \norm{(r\Nd)^{\leq 2}(\rho-\rhoo)}_{L^2(\CCba)} + \norm{(r\Nd)^{\leq 2}(\sigma-\sigmao)}_{L^2(\CCba)} + \norm{(r\Nd)^{\leq 2}\chih\wedge\chibh}_{L^2(\CCba)}\\
                                               & \les \varep\uba^{-2},
\end{align*}
where we only verified that from a regularity point of view, using the Bootstrap Assumption~\ref{BA:connCCba}, the nonlinear terms can be handled by an $L^\infty(\CCba)\times L^2(\CCba)$ estimate as before.\\ 

Multiplying equation~\eqref{eq:Dd1zeCCba} with $r$ and commuting with $(\qq\Nd_3)$ gives
\begin{align}\label{eq:Dd1uNd3zeCCba}
  r\Dd_1((\qq\Nd_3)\ze) & = \le((\qq\Nd_3)\le(-r\rho+r\rhoo\ri), (\qq\Nd_3)(r\sigma-\half r\chih\wedge\chibh)\ri) + \Err,
\end{align}
where, using~\eqref{eq:commNd3Nd}, we have schematically
\begin{align*}
  \Err & := \le[(\qq\Nd_3),(r\Dd_1)\ri]\ze \\
       & = \qq r \le((\trchib-\trchibo) \Nd\ze + \chibh\cdot\Nd\ze + \chib\cdot\ze\cdot\ze + \beb\cdot\ze\ri).
\end{align*}

Applying the elliptic estimates of Lemma~\ref{lem:ell} gives
\begin{align*}
  \norm{r^{-1}(r\Nd)^{\leq 2}(\qq\Nd_3)\ze}_{L^2(\CCba)} & \les \norm{r^{-1}(r\Nd)^{\leq 1}(\qq\Nd_3)\le(r(\rho-\rhoo\ri)}_{L^2(\CCba)} + \norm{r^{-1}(r\Nd)^{\leq 1}(\qq\Nd_3)(r(\sigma-\sigmao))}_{L^2(\CCba)} \\
                                                         & \quad + \norm{r^{-1}(r\Nd)^{\leq 1}(\qq\Nd_3)\le(r\chih\wedge\chibh\ri)}_{L^2(\CCba)} + \norm{r^{-1}(r\Nd)^{\leq 1}\Err}_{L^2(\CCba)} \\
                                                         & \les \varep \uba^{-2},
\end{align*}
where, applying the principle of Remark~\ref{rem:DecayandRegularity}, we only check that using the Bootstrap Assumptions~\ref{BA:connCCba}, all nonlinear terms from $(r\Nd)^{\leq 1}(\qq\Nd_3)\le(r\chih\wedge\chibh\ri)$ and $(r\Nd)^{\leq 1}\Err$ can be handled using $L^\infty(\CCba)\times L^2(\CCba)$ estimates.\\

Commuting~\eqref{eq:Dd1uNd3zeCCba} further with $(\qq\Nd_3)$ and arguing along the same lines gives
\begin{align*}
  \norm{r^{-1}(r\Nd)^{\leq 1}(\qq\Nd_3)^2\ze}_{L^2(\CCba)} & \les \varep \uba^{-2}.
\end{align*}

Commuting~\eqref{eq:Dd1uNd3zeCCba} with $(\qq\Nd_3)^2$, using the Bianchi equations~\eqref{eq:Nd3rho}, \eqref{eq:Nd3sigma} for $\Nd_3\rho$ and $\Nd_3\sigma$ respectively
\begin{align*}
  (r\Dd_1)(\qq\Nd_3)^{3}\ze & = \le((\qq\Nd_3)^3(-r\rho+r\rhoo), (\qq\Nd_3)^3\le(r\sigma - \half r \chih\wedge\chibh\ri)\ri) + \le[(\qq\Nd_3)^3,(r\Dd_1)\ri]\ze \\
  & = \bigg(\Divd\le(-r\qq(\qq\Nd_3)^2\beb\ri), \Curld\le(r\qq(\qq\Nd_3)^2\beb\ri)\bigg) + \EEE,
\end{align*}
where $\EEE$ is composed of lower order linear and nonlinear terms and we schematically have from~\eqref{eq:Nd3rho},~\eqref{eq:Nd3sigma},~\eqref{eq:Nd3rhoo},~\eqref{eq:sigmao} and~\eqref{eq:commNd3Nd}
\begin{align*}
  \EEE & := r(\qq\Nd_3)^{\leq 2}\le(\rho-\rhoo\ri) + r(\qq\Nd_3)^{\leq 2}(\sigma-\sigmao) + r(\qq\Nd_3)^{\leq 2}\le(\chih\cdot\alb + \ze\cdot\beb+ \le(\trchib-\trchibo\ri)(\rho-\rhoo)\ri) \\
       & \quad + (\qq\Nd_3)^{\leq 1}\bigg(\qq^2 r \le((\trchib-\trchibo) \Nd\beb + \chibh\cdot\Nd\beb + \chib\cdot\ze\cdot\beb + \beb\cdot\beb\ri)\bigg) \\
       &\quad + (\qq\Nd_3)^{\leq 3}\le(r\chih\wedge\chibh\ri) \\
       & \quad +(\qq\Nd_3)^{\leq 2}\bigg(\qq r \le((\trchib-\trchibo) \Nd\ze + \chibh\cdot\Nd\ze + \chib\cdot\ze\cdot\ze + \beb\cdot\ze\ri)\bigg).
\end{align*}
All the linear and nonlinear terms of $\EEE$ can be controlled in $L^2(\CCba)$, using $L^\infty(\CCba)\times L^2(\CCba)$ estimates and the Bootstrap Assumptions~\ref{BA:connCCba} for the nonlinear error terms.\footnote{The main point here is to verify that the uncontrolled term $\Nd_3^3\omb$ does not appear in $\EEE$.}\\

Using the elliptic estimate~\eqref{est:ellHodgeXY} of Lemma~\ref{lem:ell} with $X = -r\qq(\qq\Nd_3)^2\beb$ and $Y=r\qq(\qq\Nd_3)^2\beb$, we have
\begin{align*}
  \norm{r^{-1}(\qq\Nd_3)^{3}\ze}_{L^2(\CCba)} & \les \norm{r^{-1}\qq(\qq\Nd_3)^2\beb}_{L^2(\CCba)} + \norm{\EEE}_{L^2(\CCba)}\\
  & \les \varep \uba^{-2}.
\end{align*}

We thus have proved
\begin{align}\label{est:NdNd3zeast}
  \OOastgood\le[\ze\ri] & \les \varep.
\end{align}


\section{Control of $\protect\omb-\protect\ombo$}
Using the results of Lemmas~\ref{lem:relgeodnull} and \ref{lem:relCCba}, equation~\eqref{eq:Nd3ze} rewrites
\begin{align*}
  2\Nd\omb & = -\Nd_3\ze - 2\chib\cdot\ze - \beb.
\end{align*}

From the estimates~\eqref{est:NdNd3zeast} for $\ze$, the estimates~(\ref{est:sourcecurvast}) for $\beb$ and the Bootstrap Assumptions~\ref{BA:connCCba}, one directly deduces
\begin{align*}
  \begin{aligned}
    \norm{r^{-1}(r\Nd)^{\leq 2}\qq\Nd\omb}_{L^2(\CCba)} & \les \varep \uba^{-2},\\
    \norm{r^{-1}(r\Nd)^{\leq 1}(\qq\Nd_3)\qq\Nd\omb}_{L^2(\CCba)} & \les \varep \uba^{-2},\\
    \norm{r^{-1}(r\Nd)^{\leq 1}(\qq\Nd_3)^2\qq\Nd\omb}_{L^2(\CCba)} & \les \varep \uba^{-2},
  \end{aligned}
\end{align*}
which, using Poincar\'e estimates (see Lemma~\ref{lem:ell}) further gives
\begin{align}\label{est:NdNd3ombast}
  \OOastbadbad\le[\omb-\ombo\ri] & \les \varep.
\end{align}


\begin{remark}\label{rem:lossNd3logOm}
  It is expected for the canonical foliation that we do not control all $\Nd_3$ derivatives of the coefficient $\omb$, since the foliation is only defined by an elliptic equation on $2$-spheres. See also the estimates for the lapse $\Om$ in~\cite{Czi.Gra19}.
\end{remark}

\section{Control of $\protect\trchib-\protect\trchibo$}
Applying the Hardy estimate of Lemma~\ref{lem:transportast} with $\kappa=2, \la = -1$ to the transport equation~\eqref{eq:Nd3trchibtrchibo} for $\trchib-\trchibo$, using that from Theorems~\ref{thm:vertex} and~\ref{thm:canonical} one has the vertex limit $r^2\le(\trchib+\frac{2}{r}\ri) \to 0$ when $r\to 0$, we have
\begin{align*}
  \norm{r^{-1}\le(\trchib-\trchibo\ri)}_{L^2(\CCba)} & \les \norm{r^{-1}\qq \trchib(\omb-\ombo)}_{L^2(\CCba)} + \norm{r^{-1}\qq\Err\le(\Nd_3,\trchib-\trchibo\ri)}_{L^2(\CCba)}.
\end{align*}
Using the improved bound~\eqref{est:NdNd3ombast} for $\omb$ and the Bootstrap Assumptions~\ref{BA:connCCba}, we have
\begin{align*}
  \norm{r^{-1}\qq \trchib(\omb-\ombo)}_{L^2(\CCba)} & \les \varep \uba^{-2},\\
  \norm{r^{-1}\qq\Err\le(\Nd_3,\trchib-\trchibo\ri)}_{L^2(\CCba)} & \les \norm{r^{-1}\qq u^{-1}\chibh}_{L^\infty(\CCba)} \norm{u\chibh}_{L^2(\CCba)} + \text{faster decaying terms} \\
                                                    & \les \le(D\varep \uba^{-2}\ri) (D\varep) \\
  & \les \varep \uba^{-2}.
\end{align*}
Thus,
\begin{align*}
  \norm{r^{-1}\le(\trchib+\frac{2}{r}\ri)}_{L^2(\CCba)} & \les \varep \uba^{-2}.
\end{align*}
Commuting equation~\eqref{eq:Nd3trchibtrchibo} with $(r\Nd)^{\leq 3}$ or directly estimating its $(\qq\Nd_3)^{\leq 2}$ derivatives, applying the principle of Remark~\ref{rem:DecayandRegularity} and checking that from a regularity point of view all terms can be estimated using $L^\infty(\CCba) \times L^2(\CCba)$ estimates, we further infer
\begin{align}\label{est:trchibtrchiboL2ast}
  \OOastgood\le[\trchib-\trchibo\ri] & \les \varep.
\end{align}

\section{Control of $\protect\chibh$}
Using the elliptic equation~(\ref{eq:Divdchibh}), the elliptic estimates from Lemma~\ref{lem:ell}, estimates~(\ref{est:NdNd3zeast}), (\ref{est:trchibtrchiboL2ast}), the curvature estimates~(\ref{est:sourcecurvast}) and the Bootstrap Assumptions~\ref{BA:connCCba}, we have
\begin{align*}
  \begin{aligned}
    \norm{r^{-1}u(r\Nd)^{\leq 3}\chibh}_{L^2(\CCb_\ub)} & \les \norm{u(r\Nd)^{\leq 2}\beb}_{L^2(\CCb_\ub)} + \norm{u (r\Nd)^{\leq 2} \Nd\trchib}_{L^2(\CCb_\ub)} \\
    & \quad + \norm{u (r\Nd)^{\leq 2}(\trchib\ze)}_{L^2(\CCb_\ub)} + \Err \\
    & \les \varep \ub^{-1}.
  \end{aligned}
\end{align*}

Using the just obtained estimate for $\chibh$ and directly estimating equation~\eqref{eq:Nd3chibh} for $\Nd_3\chibh$, one further has
\begin{align}\label{est:chibhL2ast}
  \OOastbad\le[\chibh\ri] & \les \varep
\end{align}


\section{Control of $\trchi-\trchio$}
Applying the Hardy estimate of Lemma~\ref{lem:transportast} with $\kappa =1$ and $\la = -1$, to the transport equation~\eqref{eq:Nd3trchitrchio} for $\trchi-\trchio$, using that from Theorems~\ref{thm:vertex} and~\ref{thm:canonical}, one has the limit $r\le(\trchi-\trchio\ri) \to 0$ when $r\to 0$, we have
\begin{align*}
  \norm{r^{-1}\le(\trchi-\trchio\ri)}_{L^2(\CCba)} & \les \norm{r^{-1}\qq \trchi(\omb-\ombo)}_{L^2(\CCba)} + \norm{r^{-1}\qq \Err\le(\Nd_3,\trchi-\trchio\ri)}_{L^2(\CCba)}.
\end{align*}
Using estimate~(\ref{est:NdNd3ombast}) for $\omb-\ombo$ and the Bootstrap Assumptions~\ref{BA:connCCba}, one deduces
\begin{align*}
  \norm{r^{-1}\le(\trchi-\trchio\ri)}_{L^2(\CCba)} & \les \varep \uba^{-2}.
\end{align*}

Commuting the transport equation~\eqref{eq:Nd3trchitrchio} further with $(r\Nd)^{\leq 3}$ or directly estimating its $(\qq\Nd_3)^{\leq 2}$ derivatives, we further infer
\begin{align}\label{est:trchitrchioL2ast}
  \OOastgood\le[\trchi-\trchio\ri] & \les \varep
\end{align}

\section{Control of $\chih$}
Using elliptic equation~(\ref{eq:Divdchih}) for $\chih$, the elliptic estimate from Lemma~\ref{lem:ell}, and the bounds (\ref{est:sourcecurvast}), (\ref{est:NdNd3zeast}) and (\ref{est:trchitrchioL2ast}), one has
\begin{align*}
  \begin{aligned}
  \norm{r^{-1} (r\Nd)^{\leq 3}\chih}_{L^2(\CCba)} & \les \norm{(r\Nd)^{\leq 2}\be}_{L^2(\CCba)} + \norm{(r\Nd^{\leq 2})\Nd\trchi}_{L^2(\CCba)} + \text{l.o.t.} \\
  & \les \varep \uba^{-2}.
  \end{aligned}
\end{align*}

Using equation~(\ref{eq:Nd3chih}) for $\Nd_3\chih$ and estimating directly, one further infers
\begin{align}\label{est:chihL2ast}
  \OOastgood\le[\chih\ri] & \les \varep.
\end{align}

This finishes the improvement of
\begin{align*}
  \OOast & \les \varep,
\end{align*}
from which, using the Klainerman-Sobolev estimates of Lemma~\ref{lem:KlSobast} and the vertex limits of Theorems~\ref{thm:vertex} and~\ref{thm:canonical}, one also deduces
\begin{align*}
  \OOfast & \les \varep.
\end{align*}


\begin{remark}\label{rem:optimalregchi}
  Due to the choice of the canonical foliation, we have obtained the optimal tangential regularity for $\chi$, namely an $L^2(\CCb)$ control of $\Nd^3\chi$. In view of the transport equation~(\ref{eq:Nd3trchi}), this would not have been the case for other foliations, such as the geodesic foliation. The need of this optimal regularity in the extension argument\footnote{We recall that optimal regularity of the last slice $\Si_\tast$ can be obtained provided that optimal regularity holds for its boundary $\pr\Si_\tast = S^\ast = S_{\cc\uba,\uba} \subset \CCba$, which is a $2$-sphere of the canonical foliation on $\CCba$.} of Section~\ref{sec:tipcurvest} motivates the choice of the canonical foliation. See also its use in~\cite{Czi.Gra19,Czi.Gra19a,Kla.Nic03} for similar tangential regularity reasons.   
\end{remark}


\section{Uniformisation of $S^\ast$}
We have the following lemma, which is a consequence of the Uniformisation Theorem~\cite[Theorem 3.1]{Kla.Sze19}.
\begin{lemma}[Uniformisation of $S^\ast$]\label{lem:unifSastthm}
  There exists a unique -- up to isomorphisms of $\SSS$ -- centred conformal isomorphism $\Phi~:~S^\ast \to \SSS$ (see the definitions in Section~\ref{sec:defUnifandharmo}). The associated conformal factor $\phi$ on $S^\ast$ satisfies
  \begin{align}\label{est:phiSast}
    \norm{t(t\Nd)^{\leq 3}(\phi-1)}_{\HHt(S^\ast)} & \les \varep.
  \end{align}
\end{lemma}
\begin{proof}
  Using the improved estimate for $\rho$ from Sections~\ref{sec:curvest} and estimate~(\ref{est:Linfrhooast}), and using the improved estimates~(\ref{est:Linftrchiboast}), (\ref{est:Linftrchioast}), (\ref{est:trchibtrchiboL2ast}), (\ref{est:trchitrchioL2ast}) together with the $\HHt$ Klainerman-Sobolev embeddings of Lemma~\ref{lem:KlSobast}, we have on $S^\ast$
  \begin{align*}
    \norm{t^{3}(t\Nd)^{\leq 1}\rho}_{\HHt(S^\ast)} & \les \varep, \\
    \norm{t^2(t\Nd)^{\leq 2}\le(\trchi-\frac{2}{\rast}\ri)}_{\HHt(S^\ast)} & \les \varep,\\
    \norm{t^2 (t\Nd)^{\leq 2}\le(\trchib+\frac{2}{\rast}\ri)}_{\HHt(S^\ast)} & \les \varep.
  \end{align*}
  Using the above estimates, Bootstrap Assumptions~\ref{BA:connCCba} and Gauss equation~\eqref{eq:Gauss}, we have
  \begin{align*}
    \norm{t^{3}(t\Nd)^{\leq 1}\le(K-\frac{1}{\rast^2}\ri)}_{\HHt(S^\ast)} & \les \varep,
  \end{align*}
  from which, applying the Uniformisation Theorem~\cite[Theorem 3.1]{Kla.Sze19} one deduces the existence, uniqueness up to isometry of a centred conformal isomorphism $\Phi$, and the estimate~\eqref{est:phiSast} for $\phi$.
\end{proof}

\section{Control of the rotation vectorfields $\OOE$ on $S^\ast$}\label{sec:rotationSast}
\begin{lemma}[Mild bounds for $\OOE$ on $S^\ast$]\label{lem:mildOOESast}
  We have the following mild bounds on $S^\ast$
  \begin{align}\label{est:mildOOESast}
    \begin{aligned}
      \norm{r^{-1}\OOE}_{L^\infty(S^\ast)} & \les 1,\\
      \norm{r^{-2}(r\Nd)^{\leq 3}\OOE}_{L^2(S^\ast)} & \les 1.
    \end{aligned}
  \end{align}
  Moreover, for all $S$-tangent scalar $f$ and for all $1$-tensor or symmetric traceless $2$-tensor $F$, we have
  \begin{align}\label{est:mildestimateOOE}
    \begin{aligned}
    \norm{(r\Nd)f}_{L^2(S^\ast)} & \les \sum_{\ell=1}^3\norm{\OOEi(f)}_{L^2(S^\ast)},\\
    \norm{(r\Nd)^{\leq 1}F}_{L^2(S^\ast)} & \les \sum_{\ell=1}^3\norm{\Lieh_\OOEi F}_{L^2(S^\ast)}.
    \end{aligned}
  \end{align}
\end{lemma}
\begin{proof}
  By rescaling, we shall assume that $r^\ast = 1$. By definition (see Section~\ref{sec:defUnifandharmo}), we have the following uniform bounds on $S^\ast$ for the Cartesian coordinates $x^i$ and their derivatives with respect to the Euclidean metric $\gd_\SSS$
  \begin{align}\label{est:confxiSast}
    \begin{aligned}
      \le|x^i\ri| & \les 1,\\
      \norm{\le(^{\gd_\SSS}\Nd\ri)^k x^i}_{\gd_{\SSS}} & \les_k 1,
    \end{aligned}
  \end{align}
  for all $k\geq 1$ and all $i=1,2,3$.\\

  We recall the definition of $^{(3)}\OOE$
  \begin{align*}
    ^{(3)}\OOE & = x^1 {^{\gd}\Nd} x^2 - x^2 {^{\gd}\Nd} x^1 = x^1 {^{\gd_{\SSS}}\Nd} x^2 - x^2 {^{\gd_\SSS}\Nd} x^1. 
  \end{align*}
  Thus, by the definition of the conformal factor $\phi$ from Section~\ref{sec:defUnifandharmo}, we deduce the following uniform bound on $S^\ast$
  \begin{align*}
    \norm{{^{(3)}\OOE}}_{\gd} & = \phi^{-1}\norm{{^{(3)}\OOE}}_{\gd_{\SSS}} \les 1,
  \end{align*}
  where we used that from the bounds~\eqref{est:phiSast} and the Sobolev estimates of Lemma~\ref{lem:Sobsphere}, one has
  \begin{align}\label{est:phimild}
    |\phi| + |\phi^{-1}| \les 1.
  \end{align}
  We further have
  \begin{align}\label{eq:Ndgd3OOE}
    \begin{aligned}
      {^{\gd}\Nd} {^{(3)}\OOE} & = {^{\gd}\Nd}x^1 {^{\gd}\Nd} x^2 - {^{\gd}\Nd}x^2 {^{\gd}\Nd} x^1 \\
      & \quad + x^1 {^{\gd}\Nd}{^{\gd}\Nd}x^2 - x^2 {^{\gd}\Nd}{^{\gd}\Nd}x^1.
    \end{aligned}
  \end{align}
  From a standard computation of the Christoffel symbols for conformal metrics, we have for all scalar function $f$
  \begin{align}\label{eq:NdNdgdgdNdNd}
    \begin{aligned}
      {^{\gd}\Nd}_i{^{\gd}\Nd}_j f & = {^{\gd_\SSS}\Nd}_i{^{\gd_\SSS}\Nd}_jf + {^{\gd_\SSS}\Nd}_i(\log\phi){^{\gd_\SSS}\Nd}_jf + {^{\gd_\SSS}\Nd}_j(\log\phi){^{\gd_\SSS}\Nd}_if \\
      & \quad - \le({^{\gd_\SSS}\Nd}(\log\phi)\cdot{^{\gd_\SSS}\Nd}f\ri)(\gd)_{ij}.
      \end{aligned}
  \end{align}
  From the bounds~\eqref{est:phiSast} for $\phi$ and arguing as previously, we have
  \begin{align}\label{est:Ndlogphi}
    \norm{^{\gd_\SSS}\Nd(\log\phi)}_{L^2(S^\ast)} & \les \varep \les 1.
  \end{align}
  Thus, applying formula~\eqref{eq:NdNdgdgdNdNd} to equation~\eqref{eq:Ndgd3OOE} using (\ref{est:confxiSast}), \eqref{est:phimild} and ~(\ref{est:Ndlogphi}), we have
  \begin{align*}
    \norm{{^{\gd}\Nd} {^{(3)}\OOE}}_{L^2(S^\ast)} & \les 1.
  \end{align*}
  Differentiating equation~\eqref{eq:Ndgd3OOE} further by $\Nd^{\leq 2}$, generalising formula~\eqref{eq:NdNdgdgdNdNd}, using that by~\eqref{est:phiSast}, one has
  \begin{align}\label{est:Nd3logphi}
    \norm{\le(^{\gd_\SSS}\Nd\ri)^{\leq 3}(\log\phi)}_{L^2(S^\ast)} & \les \varep \les 1,
  \end{align}
  and arguing as previously, we obtain
  \begin{align*}
    \norm{{^{\gd}\Nd}^{\leq 3} {^{(3)}\OOE}}_{L^2(S^\ast)} & \les 1,
  \end{align*}
  as desired.\\

  We turn to the proof of~\eqref{est:mildestimateOOE}. From an exact computation in the Euclidean case, we have
  \begin{align*}
    \le|r\Nd f\ri|_{\gd_\SSS}^2 & \les \sum_{\ell=1}^3 \le|\OOEi(f)\ri|_{\gd_\SSS}^2,\\
    \le|(r^{\gd_\SSS}\Nd)F\ri|_{\gd_{\SSS}}^2 + \le|F\ri|_{\gd_\SSS}^2& \les \sum_{\ell=1}^3\le|\Lieh_\OOEi F\ri|_{\gd_\SSS}^2,
  \end{align*}
  for all $S$-tangent scalar $f$ and for all $1$-tensor or symmetric traceless $2$-tensor $F$. From the first bound and~\eqref{est:phimild}, we directly deduce~\eqref{est:mildestimateOOE} for scalar functions.\\
  
  From the second bounds, noting that $\Lieh_\OOEi$ is invariant under conformal change and arguing as earlier, using the control~\eqref{est:Ndlogphi} for $\Nd\log\phi$, one has
  \begin{align*}
    \le|(r^{\gd}\Nd)F\ri|_{\gd}^2 + \le|F\ri|_{\gd}^2 & \les \le|r^{\gd_\SSS}\Nd F\ri|_{\gd_\SSS}^2 + (1+\le|\Nd\log\phi\ri|) \le|F\ri|_{\gd_\SSS}^2 \\
                                                      & \les \le|r^{\gd_\SSS}\Nd F\ri|_{\gd_\SSS}^2 + \le|F\ri|_{\gd_\SSS}^2 \\
                                                      & \les \sum_{\ell=1}^3\le|\Lieh_\OOEi F\ri|^2_{\gd_\SSS} \\
    & \les \sum_{\ell=1}^3\le|\Lieh_\OOEi F\ri|^2_{\gd}.
  \end{align*}
  This finishes the proof of the lemma.
\end{proof}

\begin{lemma}[Precise bounds for $\OOE$ on $S^\ast$]\label{lem:preciseOOESast}
  We have the following bounds
  \begin{align}
    \begin{aligned}
      \norm{r^{-1}(r\Nd)^{\leq 2}\Hrot}_{L^2(S^\ast)} & \les \varep,\\
      \norm{(r\Nd)^{\leq 1}\POE}_{L^2(S^\ast)} & \les \varep.
    \end{aligned}
  \end{align}
\end{lemma}
\begin{proof}
  By a rescaling argument, we may assume that $r^\ast=1$. Arguing as in the proof of Lemma~\ref{lem:mildOOESast}, using that from~\eqref{est:phiSast}, one has for the conformal factor $\phi$, $|\phi| + |\phi^{-1}| \les 1$ and $\norm{\Nd^{\leq 3}\log\phi}_{L^2(S^\ast)} \les \varep$, it is enough to prove that in the Euclidean case $\Hrot=0$ and $\POE=0$.\\

  From a direct computation, one has in the Euclidean case
  \begin{align}\label{eq:D2OOEuclidean}
    \begin{aligned}
      \D^2{^{(\ell)}\OOO} & = \D^2\le(\in_{\ell i j} x^i \pr_{x^j}\ri) = 0.
    \end{aligned}
  \end{align}
  Using formula~(\ref{eq:D2Nd2X}) in the Euclidean case where\footnote{We refer to Lemma~\ref{lem:D2Nd2X} for the definition of $\EEE(\D^2,\Nd^2)$.}
  \begin{align*}
    \D_4\OOO & = \OOO,&  \D_3\OOO & = -\OOO, & \EEE(\D^2,\Nd^2) & = 0,
  \end{align*}
  and the definitions~\eqref{eq:defH} and~\eqref{eq:defPOE} of $\Hrot$ and $\POE$, we deduce that
  \begin{align*}
    0 & = \D^2_{a,b}\OOO \\
      & = \Nd^2_{a,b}\OOO - \OOO_b\ea +\gd_{ab}\OOO + \half \le(\Nd_a\OOO_b+\Nd_b\OOO_a\ri) (\elb-\el) \\
      & = \POE_{abc}\ec + \half \Hrot_{ab} (\elb-\el),
  \end{align*}
  from which we deduce that $\POE = 0$ and $\Hrot=0$ in the Euclidean case. This finishes the proof of the lemma.
\end{proof}

\section{Control of $\OOE$ on $\protect\CCba\cap\MM^\ext$}\label{sec:controlOOECCBA}
\subsection{Mild control of $\OOE$}\label{sec:mildOOE}
For convenience, we write $\OOO$ any rotation vectorfield $\OOEi$. We rewrite equation~(\ref{eq:Nd3OOE}) under the following form
\begin{align}\label{eq:Nd3OOENEW}
  \Nd_3\le(r^{-1}\OOO\ri) & = r^{-1}\chibh\cdot \OOO  + \half r^{-1}(\trchib-\trchibo)\OOO.
\end{align}
Integrating equation~\eqref{eq:Nd3OOENEW} from $u=\cc\uba$ --\emph{i.e.} from the sphere $S^\ast$--, using that from the result of Lemma~\ref{lem:mildOOESast}, one has an improved mild bound $\norm{r^{-1}\OOO}_{L^\infty(S^\ast)} \les 1$, we have
\begin{align*}
  |r^{-1}\OOO| & \les 1 + \int_u^{\cc\uba}\le(|r^{-1}\chih\cdot \OOO|+ |r^{-1}(\trchib-\trchibo)\OOO|\ri)\,\d u',
\end{align*}
and from a Gr\"onwall argument, we infer
\begin{align*}
  \norm{r^{-1}\OOO}_{L^\infty(\CCba\cap\MM^\ext)} & \les 1.
\end{align*}

Arguing similarly, using the commuted equation~(\ref{eq:Nd3NdO}), and the improved mild bounds of Lemma~\ref{lem:mildOOESast}, we obtain
\begin{align}\label{est:mildOOE}
  \norm{r^{-1}(r\Nd)^{\leq 1}\OOO}_{L^\infty(\CCba\cap\MM^\ext)} & \les 1.
\end{align}

Let $S_{u,\uba}\subset \CCba$. Let $f$, $F$ be respectively an $S_{u,\uba}$-tangent scalar function and $1$-tensor or symmetric traceless $2$-tensor. We want to improve the bootstrap bounds
\begin{align*}
  \int_{S_{u,\uba}}\le|(r\Nd)f\ri|^2 & \les \sum_{\ell=1}^3\int_{S_{u,\uba}}\le|\OOEi(f)\ri|^2,\\
  \int_{S_{u,\uba}}\le|(r\Nd)^{\leq 1}F\ri|^2 & \les \sum_{\ell=1}^3\int_{S_{u,\uba}}\le|\Lieh_\OOEi F\ri|^2.
\end{align*}
Let extend $f,F$ on $\CCba$ as $S$-tangent tensors by parallel transport, \emph{i.e.}
\begin{align*}
  \elb(f) &  = 0, & \Nd_3 F  & = 0.
\end{align*}
From the definition of $\OOEi$ on $S^\ast$, we obtained in Lemma~\ref{lem:mildOOESast} that there exists $c>0$ such that
\begin{align*}
  \int_{S^\ast}\le|(r\Nd)f\ri|^2 & \leq c\sum_{\ell=1}^3\int_{S^\ast}\le|\OOEi(f)\ri|^2,\\
  \int_{S^\ast}\le|(r\Nd)^{\leq 1}F\ri|^2 & \leq c \sum_{\ell=1}^3\int_{S^\ast}\le|\Lieh_\OOEi F\ri|^2.
\end{align*}
Let define
\begin{align*}
  \FF_f(u) & := \int_{S_{u,\uba}}\le|(r\Nd)f\ri|^2 - c \sum_{\ell=1}^3\int_{S_{u,\uba}}\le|\OOEi(f)\ri|^2 .
\end{align*}
Using formula~\eqref{eq:commelbint} and~\eqref{eq:commNd3Nd}, we have
\begin{align*}
  \elb\le(\int_{S_{u,\uba}}\le|(r\Nd)f\ri|^2\ri) & = \int_{S_{u,\uba}} \trchib|r\Nd f|^2 + \int_{S_{u,\uba}} 2(r\Nd f)\cdot\Nd_3(r\Nd)f \\
                                                & = \int_{S_{u,\uba}}\trchib|r\Nd f|^2 + \Err\le[(r\Nd) f\ri],
\end{align*}
where
\begin{align*}
  \Err\le[(r\Nd)f\ri] & := \int_{S_{u,\uba}} 2(r\Nd f)\cdot\le((\trchibo-\trchib)(r\Nd)f+\chibh\cdot(r\Nd f)\ri).
\end{align*}
Using formula~\eqref{eq:commelbint} and the Lie transport of $\OOEi$~(\ref{eq:Lie3OOE}), we have
\begin{align*}
  \elb\le(\int_{S_{u,\uba}}\le|\OOEi(f)\ri|^2\ri) & = \int_{S_{u,\uba}}\trchib\le|\OOEi (f)\ri|^2.
\end{align*}
Thus, we deduce that
\begin{align}\label{eq:Nd3FFf}
  \elb \FF_f & = \trchibo\FF_f + \Err\le[(r\Nd f)\ri]. 
\end{align}
Integrating~\eqref{eq:Nd3FFf} along $u$, using that $\FF_f(\cc\uba) \leq 0$ and the bounds~(\ref{est:OOfb2ast}), we have
\begin{align*}
  \FF_f(u) & \les \le|r^2\int_{u}^{\cc\uba}r^{-2} \le(\varep r^{-1}(u')^{-3/2}\ri) \int_{S_{u',\uba}}|(r\Nd)f|^2 \,\d u'\ri| \\
           & \les \varep \sup_{u \leq u'\leq \cc\uba}\int_{S_{u',\uba}}|(r\Nd)f|^2 \\
           & \les \varep \int_{S_{u,\uba}}|(r\Nd)f|^2,
\end{align*}
where one obtains the last estimate by using the transport of $f$ along $\elb$. Using the definition of $\FF_f$, we deduce that
\begin{align*}
  \int_{S_{u,\uba}}\le|(r\Nd) f\ri|^2 & \leq c \sum_{\ell=1}^3\int_{S_{u,\uba}}\le|\OOEi(f)\ri|^2 + \varep\int_{S_{u,\uba}}|(r\Nd)f|^2,
\end{align*}
and the desired bound follows by absorption. The bounds when $F$ is a $1$-tensor or symmetric traceless $2$-tensor are obtained similarly and left to the reader. This finishes the improvement of the mild Bootstrap Assumptions~\ref{BA:mildOOE} on $\CCba\cap\MM^\ext$.






\subsection{Control of $\Hrot$}
Applying similar Hardy estimates to the ones of Lemma~\ref{lem:transportast} with $\kappa=0, \lambda = -2$ -- but integrating from $S^\ast$ -- to the transport equation~(\ref{eq:Nd3H}) for $\Hrot$, using that from Lemma~\ref{lem:preciseOOESast} we have $\norm{\uba^{1/2} \Hrot}_{L^2(S^\ast)} \les \varep$, we have
\begin{align*}
  \norm{r^{-2}\Hrot}_{L^2(\CCba\cap\MM^\ext)} & \les \uba^{1/2}\norm{r^{-2}\Hrot}_{L^2(S^\ast)} + \norm{r^{-2}\qq\Nd_3\Hrot}_{L^2(\CCba\cap\MM^\ext)} \\
                                          & \les \varep \uba^{-2} + \norm{r^{-2}u\Nd_3\Hrot}_{L^2(\CCba\cap\MM^\ext)}.
\end{align*}
From an inspection of the terms composing $\Nd_3\Hrot$, we only treat the terms $\Nd\chibh\cdot\OOO$ and $\chibh \cdot\Nd\OOO$ and the remaining terms will follow similarly. We have using estimate (\ref{est:chibhL2ast}) and the mild estimate~\eqref{est:mildOOE} for $\OOO$
\begin{align*}
  \norm{r^{-2}\qq \Nd\chibh \cdot \OOO}_{L^2(\CCba)} & \les \norm{r^{-1}\qq \Nd\chibh}_{L^2(\CCba)}\norm{r^{-1}\OOO}_{L^\infty(\CCba)} \\
                                                     & \les \uba^{-2}\norm{r^{-1}\uba u (r\Nd)\chibh}_{L^2(\CCba)} \\
  & \les \varep \uba^{-2},
\end{align*}
and
\begin{align*}
  \norm{r^{-2}\qq\chibh\cdot\Nd\OOO}^2_{L^2(\CCba)} & = \int_1^\uba \norm{r^{-2}\qq\chibh\cdot\Nd\OOO}_{L^2(S_{u,\uba})}^2\,\d u \\
                                                  & \les \int_1^\uba \norm{r^{-1}\qq\chibh}^2_{L^\infty(S_{u,\uba})}\norm{r^{-2}(r\Nd)\OOO}^2_{L^\infty_uL^2(S_{u,\uba})} \,\d u \\
                                                  & \les \int_1^\uba \le(r^{-1}\norm{r^{-1}\qq(r\Nd)^{\leq 2}\chibh}_{L^2(S_{u,\uba})}\ri)^2\,\d u\\
                                                  & \les \norm{r^{-1}\uba^{-1}u(r\Nd)^{\leq 2}\chibh}^2_{L^2(\CCba)}\\
                                                  & \les \varep^2 \uba^{-4},
\end{align*}
where we used the Sobolev embeddings from Lemma~\ref{lem:Sobsphere} in the third line.
Thus,
\begin{align*}
  \norm{r^{-2}\Hrot}_{L^2(\CCba\cap\MM^\ext)} & \les \varep \uba^{-2}.
\end{align*}

Estimating directly equation~\eqref{eq:Nd3H}, one also obtains
\begin{align*}
  \norm{r^{-2}\qq\Nd_3\Hrot}_{L^2(\CCba\cap\M^\ext)} & \les \varep \uba^{-2}.
\end{align*}
From the definitions~\eqref{eq:defH} and~\eqref{eq:defPOE} of $\Hrot$ and $\POE$, the higher derivative of $\Hrot$ are controlled by $\POE$ and from the results of this section we have
\begin{align}\label{est:HPOE}
  \OO_{\leq 2}^{\ast,\mathfrak{g}}\le[r^{-1}\Hrot\ri] & \les \varep + \OO_{\leq 1}^{\ast,\mathfrak{g}}\le[\POE\ri].
\end{align}
We thus refer to the control of $\POE$ in the next section to obtain the full control of $\Hrot$.


\subsection{Control of $\POE$}\label{sec:POECCba}
The control for $\POE$ follows from applying the same Hardy estimate as in the previous section to the (commuted with $r\Nd$) transport equation~\eqref{eq:Nd3POE} for $r\POE$, using the initial estimates from Lemma~\ref{lem:preciseOOESast} on $S^\ast$. Details are left to the reader and we have
\begin{align}\label{est:POEastL2}
  \OO^{\ast,\mathfrak{g}}_{\leq 1}\le[\POE\ri] & \les \varep.
\end{align}
We moreover deduce from~\eqref{est:HPOE} and~\eqref{est:POEastL2} 
\begin{align*}
  \OO^{\ast,\mathfrak{g}}_{\leq 2}\le[r^{-1}\Hrot\ri] & \les \varep,
\end{align*}
and this finishes the improvement of~(\ref{est:OOb3ast}).


\section{Control of the area radius}\label{sec:arearadiusestimateCCba}
\begin{lemma}
  Under the estimates obtained in the previous sections, we have on $\CCba$
  \begin{align}\label{est:arearadiusestimateCCba}
    \le|r(u,\uba) -\half (\uba-u)\ri| & \les \varep \uba^{-3}r^2 u^{-1}.
  \end{align}
\end{lemma}
\begin{proof}
  Using relations~\eqref{eq:elu} and~\eqref{eq:elbr}, we have on $\CCba$
  \begin{align*}
    \elb\le(r(u,\uba)-\half(\uba-u)\ri) & = \half r\trchibo + 1.
  \end{align*}
  Integrating in $u$, using the improved estimate~(\ref{est:Linftrchiboast}) for $\trchibo+\frac{2}{r}$ and the limit $r(u,\uba)\to 0$ when $u\to\uba$ (see Theorems~\ref{thm:vertex} and~\ref{thm:canonical}), we have
  \begin{align*}
    \le|r(u,\uba)-\half(\uba-u)\ri| & \les \lim_{u\to\uba}\le|r(u,\uba)-\half(\uba-u)\ri| + \int_u^\uba r \le|\trchibo+\frac{2}{r}\ri|\,\d u' \\
                                    & \les \varep \int_u^\uba r \uba^{-2}(u')^{-2}\,\d u' \\
                                    & \les \varep \uba^{-3}r^2u^{-1},
  \end{align*}
  as desired.
\end{proof}

\section{Control of spherical coordinates on $\protect\CCba$}\label{sec:mildBAimpast}
For a spherical coordinate system as described in the Bootstrap Assumption~\ref{BA:mildsphcoordsast}, we have using~(\ref{eq:elbgd}) the following transport equation in the $\elb$ direction
\begin{align*}
  \elb(\gd_{ab}) & = 2\chib_{ab} = \trchibo \gd_{ab} + (\trchib-\trchibo)\gd_{ab} + 2\chibh_{ab},
\end{align*}
where $a,b\in\{\varth,\varphi\}$, which rewrites using~\eqref{eq:elbr} and the notation of the Bootstrap Assumption~\ref{BA:mildsphcoordsast}
\begin{align}\label{eq:transportrenormgdast}
  \elb\le(r^{-2}\gd_{ab}-(\gd_{\SSS})_{ab}\ri) & = (\trchib-\trchibo)r^{-2}\gd_{ab} + 2r^{-2}\chibh_{ab}.
\end{align}
Integrating along $\elb$ from $\o(\uba)$, this yields
\begin{align*}
  \le|r^{-2}\gd_{ab}-(\gd_{\SSS})_{ab}\ri|(u,\varth,\varphi) & \les \lim_{u\to\uba}\le(\le|r^{-2}\gd_{ab}-(\gd_{\SSS})_{ab}\ri|(u,\varth,\varphi)\ri) \\
                                                             & \quad + \int_{u}^\uba \le(\le|(\trchib-\trchibo)r^{-2}\gd_{ab}\ri| + 2r^{-2}\le|\chibh_{ab}\ri|\ri)\,\d u' \\
                                                                 & \les \varep \int_{\ub}^\uba (\uba)^{-1}(u')^{-1}q^{-1/2}\,\d u' \\
                                                                 & \les \varep \uba^{-1}u^{-1} q^{1/2},
\end{align*}
where we used the vertex limits from Theorems~\ref{thm:vertex} and~\ref{thm:canonical}, the improved bounds (\ref{est:trchibtrchiboL2ast}), (\ref{est:chibhL2ast}) for $\trchib$ and $\chibh$ and the Bootstrap Assumptions~\ref{BA:mildsphcoordsast} for the coordinate component of $\gd$.\\

Commuting~\eqref{eq:transportrenormgdast} with $\pr$ and integrating using the sup-norm bounds for $(r\Nd)(\trchib-\trchibo)$ and $(r\Nd)\chibh$, we further have
\begin{align*}
  \le|\pr^{\leq 1}\le(r^{-2}\gd_{ab}-(\gd_{\SSS})_{ab}\ri)\ri| & \les \varep \uba^{-1}u^{-1}q^{1/2}.
\end{align*}
This finishes the improvement of the Bootstrap Assumptions~\ref{BA:mildsphcoordsast}.

\chapter{Null connection estimates in $\protect\MM^\ext$}\label{sec:connest}
In this section, we prove the following proposition.
\begin{proposition}\label{prop:connestSTAB}
  Recall that from Proposition~\ref{prop:curvestSTAB}, the following estimates hold for the curvature norms in $\MM^\ext$ (see the definitions of Section~\ref{sec:normnullcurv}) 
  \begin{align}\label{est:sourcecurvext}
    \begin{aligned}
      \RRfext + \RRext & \les \varep,
  \end{aligned}
  \end{align}
  for all $\ga>0$. Recall that from Proposition~\ref{prop:connestCCbaSTAB} the following estimates hold for the null connection coefficient norms on $\CCba$ (see the definitions of Sections~\ref{sec:normnullcurv} and~\ref{sec:normnullconnCCba})
  \begin{align}\label{est:initconnCCba}
    \RRfoast + \OOast + \OOoast + \OOfast + \OO_{\leq 3}^{\ast,\OOO} & \les \varep,
  \end{align}
  and that the Bootstrap Assumptions~\ref{BA:mildOOE}, \ref{BA:connCCba}, \ref{BA:mildsphcoordsast} are improved on $\CCba\cap\MM^\ext$.\\

  Under the Bootstrap Assumptions, the estimates~\eqref{est:sourcecurvext} and~\eqref{est:initconnCCba}, and for $\varep>0$ sufficiently small, we have the following bounds for the connection coefficients (see the definitions of Section~\ref{sec:normsnullconn})
  \begin{align}\label{est:OO1}
    \begin{aligned}
      \RRfoext + \OOfextii + \OOextii + \OOof & \les_\ga \varep,\\ \\
      \mathfrak{O}^\ext_{\leq 1}[\yy] + \OO^\ext_{2, \ga}[\yy] +\OO^\TT_{\leq 2, \ga}[\yy] & \les_\ga \varep, \\ \\
      \mathfrak{O}^{\ext,\OOO}_{\leq 2} + \OO^{\ext,\OOO}_{\leq 3, \ga} & \les_\ga \varep,
    \end{aligned}
  \end{align}
  for all $\ga>0$. Moreover, the Bootstrap Assumptions~\ref{BA:mildsphcoordsext} for the spherical coordinates (see Section~\ref{sec:mildBAimpext}), the bootstrap bound~(\ref{est:BAarearadiusestimate}) on the area radius (see Section~\ref{sec:arearadiusestimate}), and the mild Bootstrap Assumptions~\ref{BA:mildOOE} for the rotation vectorfields $\OOE$ (see Section~\ref{sec:mildOOEext}) are improved.\footnote{These last bounds together with~\eqref{est:OO1} also amount to an improvement the strong Bootstrap Assumptions~\ref{BA:connext}.}
\end{proposition}

\section{Evolution estimates}
Let us first recall that in the region $\MM^\ext$, we have
\begin{align*}
  r\simeq \ub.
\end{align*}

The following lemma provides estimates for solutions to transport equations in the $\el$ direction.
\begin{lemma}[Transport estimates in $\MM^\ext$]\label{lem:evolext}
  For all $\kappa\in\RRR$, the following holds. Assume that $U$ is an $S$-tangent tensor satisfying
  \begin{align*}
    \Nd_4U + \frac{\kappa}{2} \trchi U & = F.
  \end{align*}
  
  \begin{itemize}
  \item We have the following $L^\infty L^\infty$ estimates
    \begin{align}\label{est:transLinf}
      \norm{r^\la U}_{L^\infty(S_{u,\ub})} & \les \norm{r^\la U}_{L^\infty(S_{u,\uba})} + \int_{\ub}^\uba \norm{r^\la F}_{L^\infty(S_{u,\ub'})}\,\d\ub',
    \end{align}
    for all $S_{u,\ub}\subset\MM^\ext$ and for all $\la\geq \kappa$.\\
    
  \item We have the following $L^\infty H^{1/2}$ estimate
    \begin{align}\label{est:transH12}
      \norm{r^{\la} U}_{\HHt(S_{u,\ub})} & \les \norm{r^\la U}_{\HHt(S_{u,\uba})} + \int_{\ub}^\uba \norm{r^\la F}_{\HHt(S_{\ub',u})}\,\d\ub',
    \end{align}
    for all $S_{u,\ub}\subset\MM^\ext$ and for all $\la \geq \kappa -1/2$.
    
  \item We have the following $L^\infty_\ub L^2_uL^2(S_{u,\ub})$ estimate
    \begin{align}\label{est:transLinfL2}
      \norm{r^\la U}_{L^2_uL^2(S_{u,\ub})} & \les \norm{r^\la U}_{L^2(\CCba)} + \int_\ub^\uba \norm{r^\la F}_{L^2_uL^2(S_{u,\ub'})}\,\d\ub',
    \end{align}
    for all $\la\geq\kappa-1$.
  \end{itemize}
\end{lemma}
\begin{proof}
  We only perform the proof for the $\HHt$ estimate and for a $S$-tangent $1$-tensor $U$. Rewriting the transport equation satisfied by $U$, using~\eqref{eq:elbr}, we have
  \begin{align}\label{eq:Nd4Urenorm}
    \Nd_4(r^{\kappa}U) & = r^{\kappa}F + \frac{\kappa}{2}(\trchio-\trchi)r^{\kappa}U.
  \end{align}
  Using the spherical coordinates in $\MM^\ext$ from the Bootstrap Assumptions~\ref{BA:mildsphcoordsext}, we have
  \begin{align*}
    \Nd_4(r^{\kappa}U)_a & = 2 \pr_\ub\le(r^{\kappa}U_a\ri) - 2 r^{\kappa}U^b \chi_{ab},
  \end{align*}
  where $a,b\in\{\varth,\varphi\}$. Thus, using~\eqref{eq:elbr} and equation~\eqref{eq:Nd4Urenorm}, we infer
  \begin{align*}
    2\pr_\ub\le(r^{-1+\kappa}U_a\ri) & = r^{-1}\Nd_4(r^{\kappa}U)_a + r^{-1+\kappa}U_a(\trchi-\trchio) + 2r^{-1+\kappa}U_b\chih_{a}^b \\
                                     & = r^{-1+\kappa}F_a + (\frac{\kappa}{2}-1)(\trchio-\trchi)r^{-1+\kappa}U_a + 2r^{-1+\kappa}U_b\chih_a^b.
  \end{align*}
  Integrating the above equation, we obtain
  \begin{align*}
    r^{-1+\kappa}U_a(\ub,u,\varth,\varphi) & = r^{-1+\kappa}U_a(\uba,u,\varth,\varphi) \\
                                           & + \half\int_{\ub}^\uba \le(r^{-1+\kappa}F_a + (\frac{\kappa}{2}-1)(\trchio-\trchi)r^{-1+\kappa}U_a + 2r^{-1+\kappa}U_b\chih_a^b\ri)(\ub',u,\varth,\varphi)\,\d\ub'.
  \end{align*}
  Taking the coordinate $H^{1/2}_{\varth,\varphi}$ norm in the above, we obtain
  \begin{align*}
    \norm{r^{-1+\kappa}U_a(\ub,u)}_{H^{1/2}_{\varth,\varphi}} & \les \norm{r^{-1+\kappa}U_a(\uba,u)}_{H^{1/2}_{\varth,\varphi}} \\
                                                              & + \int_{\ub}^\uba \norm{\le(r^{-1+\kappa}F_a + (\frac{\kappa}{2}-1)(\trchio-\trchi)r^{-1+\kappa}U_a + 2r^{-1+\kappa}U_b\chih_a^b\ri)(\ub',u)}_{H^{1/2}_{\varth,\varphi}}\,\d\ub'.
  \end{align*}
  Using the fractional Sobolev space comparison Lemma~\ref{lem:compH12}, we thus infer
  \begin{align*}
    \norm{r^{-1/2+\kappa}U}_{\HHt(S_{u,\ub})} & \les \norm{r^{-1/2+\kappa}U}_{\HHt(S_{u,\uba})} + \int_{\ub}^\uba \norm{r^{-1/2+\kappa}F}_{\HHt(S_{u,\ub'})} \,\d\ub' \\
                                              & + \int_{\ub}^\uba\le(\norm{r^{-1/2+\kappa}(\trchi-\trchio)U}_{\HHt(S_{u,\ub'})} + \norm{r^{-1/2+\kappa}\chih\cdot U}_{\HHt(S_{u,\ub'})}\ri)\,\d\ub'.
  \end{align*}
  Using the $\HHt$ product estimates from Lemma~\ref{lem:prodH12}, and the Bootstrap Assumptions~\ref{BA:connext} for the $L^\infty$ norms of $\chi$ and $\Nd\chi$, we obtain
  \begin{align*}
    \norm{r^{-1/2+\kappa}U}_{\HHt(S_{u,\ub})} & \les \norm{r^{-1/2+\kappa}U}_{\HHt(S_{u,\ub})} + \int_{\ub}^\uba \norm{F}_{\HHt(S_{u,\ub'})}\,\d\ub' \\
                                              & + D\varep \int_{\ub}^\uba (\ub')^{-2}u^{-1/2}\norm{U}_{\HHt(S_{u,\ub'})}\,\d\ub',
  \end{align*}
  and the conclusion follows from a standard Gr\"onwall argument. The result for $\la\geq \kappa-1/2$ follows since $r(u,\ub) \les r(u,\ub')$ for $\ub \leq \ub' \leq\uba$. The remaining estimates are obtained similarly. This finishes the proof of the lemma.
\end{proof}
\begin{remark}
  The conclusions of Lemma~\ref{lem:evolext} also hold with $\trchi$ replaced by $\trchio$ or $\frac{2}{r}$.
\end{remark}

\section{Control of $\rhoo$ and $\sigmao$}
Applying the $L^\infty$ estimate of Lemma~\ref{lem:evolext}, with $\kappa=3, \la=3$ to the transport equation~\eqref{eq:Nd4rhoo}, we have
\begin{align*}
  \norm{r^3\rhoo}_{L^\infty(S_{u,\ub})} & \les \norm{r^3\rhoo}_{L^\infty(S_{u,\uba})} + \int_\ub^\uba \norm{r^3\Err(\Nd_4,\rhoo)}_{L^\infty(S_{u,\ub'})}\,\d\ub'.
\end{align*}
Using the bounds~(\ref{est:initconnCCba}) for $\rhoo$ on $\CCba$, the curvature bounds~\eqref{est:sourcecurvext} and the Bootstrap Assumptions~\ref{BA:connext}, with the expression of $\Err\le(\Nd_4,\rhoo\ri)$ from~\eqref{eq:Nd4rhoo}, we have
\begin{align*}
  \norm{r^3\rhoo}_{L^\infty(S_{u,\uba})} & \les \varep u^{-2},\\
  \int_\ub^\uba \norm{r^3\Err(\Nd_4,\rhoo)}_{L^\infty(S_{u,\ub'})}\,\d\ub' & \les (D\varep)^2\int_\ub^\uba r^3(\ub')^{-1}(\ub')^{-7/2}u^{-3/2}\,\d\ub'\\
  & \les \varep u^{-2}.
\end{align*}

Thus,
\begin{align}\label{est:rhoo}
  \norm{r^3u^2\rhoo}_{L^\infty(\MM^\ext)} & \les \varep.
\end{align}

Estimating directly equation~\eqref{eq:Nd4rhoo} for $\Nd_4\rhoo$, using~\eqref{est:rhoo},~\eqref{est:sourcecurvext} and the Bootstrap Assumptions~\ref{BA:connext}, we have
\begin{align*}
  \norm{r^3u^2(\ub\Nd_4)\rhoo}_{L^\infty} & \les \varep.
\end{align*}

Taking the average in equation~\eqref{eq:Nd3rho} for $\Nd_3\rho$, using the bound~\eqref{est:rhoo} for $\rhoo$, estimates~\eqref{est:sourcecurvext} and the Bootstrap Assumptions~\ref{BA:connext}, one has
\begin{align*}
  \norm{r^3u^2(u\Nd_3)\rhoo}_{L^\infty} & \les \varep. 
\end{align*}

Taking one more derivative in the above equations, applying the principle of Remark~\ref{rem:DecayandRegularity}, checking that the derivatives of the nonlinear error terms can be handled by Cauchy-Schwartz as in the control of $\rhoo$ in Section~\ref{sec:controlrhoosigamooCCba}, we further have
\begin{align*}
  \RRfoext[\rhoo] & \les \varep.
\end{align*}

Estimating directly equation~\eqref{eq:sigmao} for $\sigmao$, using the Bootstrap Assumptions~\ref{BA:connext}, we have
\begin{align*}
  \begin{aligned}
    \norm{r^3u^2\sigmao}_{L^\infty} & \les \norm{r^3u^2\chih\wedge\chibh}_{L^\infty} \\
    & \les (D\varep)^2 \\
    &\les \varep.
  \end{aligned}
\end{align*}

Taking the average and directly estimating equations~\eqref{eq:Nd3sigma} and~\eqref{eq:Nd4sigma} for $\Nd_3\sigma$ and $\Nd_4\sigma$ respectively, one further has
\begin{align*}
  \RRfoext[\sigmao] & \les \varep,
\end{align*}
and we have thus proved
\begin{align}\label{est:rhoosigmao}
  \RRfoext & \les \varep.
\end{align}

\section{Control of $\protect\ombo$}
Applying the $L^\infty$ estimate of Lemma~\ref{lem:evolext}, with $\kappa=0, \la = 1$, to the transport equation~\eqref{eq:Nd4ombo} for $\ombo$, using the bounds~\eqref{est:rhoosigmao} obtained for $\rhoo$ and the Bootstrap Assumptions~\ref{BA:connext}, and using that $\ombo=0$ on $\CCba$, we have

\begin{align*}
  \norm{r\ombo}_{L^\infty(S_{u,\ub})} & \les \int_\ub^\uba \le(r|\rhoo| + \Err\le(\Nd_4,\ombo\ri)\ri)\,\d\ub' \\
                                      & \les \varep\ub^{-1}u^{-2}
\end{align*}


Commuting~\eqref{eq:Nd4ombo} with $\Ndt^{\leq 1}$ and using the principle of Remark~\ref{rem:DecayandRegularity}, we further have
\begin{align}\label{est:omboLinf}
  \OOofb[\ombo] & \les \varep.
\end{align}

Applying the $L^\infty_\ub L^2_uL^2(S_{u,\ub})$ estimate of Lemma~\ref{lem:evolext}, with $\kappa=0, \la = 1-\ga$, to the transport equation~\eqref{eq:Nd4ombo} for $\ombo$ multiplied by $u^{3/2-\ga}$, using the estimates~(\ref{est:rhoosigmao}) for $\rhoo$, from which one deduces $\norm{u^{3/2-\ga}\rhoo}_{L^2_uL^2(S_{u,\ub})} \les_\ga \varep \ub^{-2}$ for all $\ga>0$, and the Bootstrap Assumptions~\ref{BA:connext}, we have
\begin{align*}
  \norm{u^{3/2-\ga}\ub^{-\ga} \ub \ombo}_{L^2_uL^2(S_{u,\ub})} & \les \int_\ub^\uba \le(\norm{u^{3/2-\ga}(\ub')^{1-\ga}\rhoo}_{L^2_uL^2(S_{u,\ub'})}+\norm{u^{-3/2-\ga}(\ub')^{1-\ga}\Err(\Nd_4,\ombo)}_{L^2_uL^2(S_{u,\ub'})}\ri)\,\d\ub' \\
                                                          & \les \varep \int_\ub^\uba (\ub')^{-1-\ga} \,\d\ub'\\
  & \les_\ga \varep.
\end{align*}
Thus, arguing as in Section~\ref{sec:errorintr1}, using that $u\les \ub$ in $\MM^\ext$
\begin{align*}
  \norm{u^{3/2}\ub^{-1/2-\ga}\ub\ombo}_{L^2(\MM^\ext)} & \les_\ga \varep,
\end{align*}
for all $\ga>0$. Commuting equation~\eqref{eq:Nd4ombo} with $\Ndt^{\leq 2}$ and using the principle of Remark~\ref{rem:DecayandRegularity} gives
\begin{align}\label{est:omboL2}
  \OOof[\ombo] & \les \varep,
\end{align}
for all $\ga>0$.

\section{Control of $\trchio-\frac{2}{r}$}
Applying the $L^\infty$ estimate of Lemma~\ref{lem:evolext}, with $\kappa=1, \la = 2$, to the transport equation~\eqref{eq:Nd4trchio} for $\trchio-\frac{2}{r}$, we have
\begin{align*}
  \norm{r^2\le(\trchio-\frac{2}{r}\ri)}_{L^\infty(S_{u,\ub})} & \les \norm{r^2\le(\trchio-\frac{2}{r}\ri)}_{L^\infty(S_{u,\uba})} + \int_\ub^\uba \norm{r^2\Err\le(\Nd_4,\trchio-\frac{2}{r}\ri)}_{L^\infty(S_{u,\ub'})}\,\d\ub'. 
\end{align*}
From the bound~(\ref{est:initconnCCba}) on $\CCba$, and the Bootstrap Assumptions~\ref{BA:connext}, we have
\begin{align*}
  \norm{r^2\le(\trchio-\frac{2}{r}\ri)}_{L^\infty(S_{u,\uba})} & \les \varep \uba^{-1}u^{-1},
\end{align*}
\begin{align*}
  \int_\ub^\uba \norm{r^2\Err\le(\Nd_4,\trchio-\frac{2}{r}\ri)}_{L^\infty(S_{u,\ub'})}\,\d\ub' & \les \int_\ub^\uba \le(r^2\norm{\chih}_{L^\infty}^2 +r^2\norm{\trchi-\trchio}_{L^\infty}^2 \ri)\,\d\ub'\\
                                                               & \les (D\varep)^2\int_\ub^\uba r^2(\ub')^{-4}u^{-1} \,\d\ub' \\
  & \les \varep \ub^{-1}u^{-1},
\end{align*}
and we therefore deduce
\begin{align}\label{est:trchio}
  \norm{u\ub^3\le(\trchio-\frac{2}{r}\ri)}_{L^\infty} & \les \varep.
\end{align}

Estimating directly equation~\eqref{eq:Nd4trchio} or commuting equation~\eqref{eq:Nd4trchio} by $u\Nd_3$, and applying analogously the $L^\infty$ estimates of Lemma~\ref{lem:evolext} -- and the principle of Remark~\ref{rem:DecayandRegularity}, one obtains
\begin{align}\label{est:Nd4Nd3trchio}
  \OOofb[\trchio] & \les \varep.
\end{align}

Commuting with $\Ndt^{\leq 2}$ derivatives and applying the $L^\infty_\ub L^2$ estimates of Lemma~\ref{lem:evolext} as in the previous section further gives
\begin{align}\label{est:Nd2trchioL2}
  \OOof[\trchio] & \les \varep,
\end{align}
for all $\ga>0$.


\section{Control of $\protect\trchibo+\frac{2}{r}$}
Applying the $L^\infty$ estimate of Lemma~\ref{lem:evolext}, with $\kappa=1, \la = 2$, to the transport equation~\eqref{eq:Nd4trchibo} for $\trchibo+\frac{2}{r}$, 
using the bounds~\eqref{est:rhoosigmao} obtained for $\rhoo$ and the Bootstrap Assumptions~\ref{BA:connext}, we have

\begin{align}\label{est:trchibo}
  \OOofb[\trchibo] & \les \varep.
\end{align}

Applying the $L^\infty_\ub L^2$ estimate of Lemma~\ref{lem:evolext} as in the previous sections gives
\begin{align}
  \OOof[\trchio] & \les \varep,
\end{align}
for all $\ga>0$.

\section{Control of $\trchi-\trchio$}\label{sec:controltrchitrchioext}
Applying the $L^\infty\HHt$ estimates of Lemma~\ref{lem:evolext} with $\kappa=2,\la=3/2$ to the transport equation~\eqref{eq:Nd4trchitrchio}, we have
\begin{align*}
  \norm{r^{3/2}(\trchi-\trchio)}_{\HHt(S_{u,\ub})} & \les \norm{r^{3/2}(\trchi-\trchio)}_{\HHt(S_{u,\uba})} \\
  & \quad + \int_\ub^\uba \norm{r^{3/2}\Err\le(\Nd_4,\trchi-\trchio\ri)}_{\HHt(S_{u,\ub'})}\,\d\ub'.
\end{align*}
Using the bounds~(\ref{est:initconnCCba}) on $\CCba$, the $\Ht$ product estimates of Lemma~\ref{lem:prodH12} and the Bootstrap Assumptions~\ref{BA:connext}, we have
\begin{align*}
  \norm{r^{3/2}(\trchi-\trchio)}_{\HHt(S_{u,\uba})} & \les \varep u^{-1/2},
\end{align*}
and
\begin{align*}
  & \int_\ub^\uba \norm{r^{3/2}\Err\le(\Nd_4,\trchi-\trchio\ri)}_{\HHt(S_{u,\ub'})}\,\d\ub' \\
  \les &\; \int_\ub^\uba \le(\norm{r^{3/2}\chih\cdot\chih}_{\HHt(S_{u,\ub'})} + \norm{r^{3/2}(\trchi-\trchio)^2}_{\HHt(S_{u,\ub'})}\ri)\,\d\ub' \\
  \les & \;\int_\ub^\uba \bigg(\le(\norm{\chih}_{L^\infty(S_{u,\ub'})} + r^{-1}\norm{r\Nd\chih}_{L^2(S_{u,\ub'})} \ri) \norm{r^{3/2}\chih}_{\HHt(S_{u,\ub'})} \\
  & \quad +\le(\norm{\trchi-\trchio}_{L^\infty(S_{u,\ub'})}+r^{-1}\norm{r\Nd(\trchi-\trchio)}_{L^2(S_{u,\ub'})}\ri)\norm{r^{3/2}\le(\trchi-\trchio\ri)}_{\HHt(S_{u,\ub'})} \bigg)\, \d\ub' \\
  \les & \, (D\varep)^2\int_{\ub}^\uba \le(r^{-2}u^{-1/2}\ri) u^{-1/2}\,\d\ub' \\
  \les & \;\varep u^{-1/2}.
\end{align*}
Thus,
\begin{align}\label{est:H12trchitrchio}
  \norm{r^{3/2}u^{1/2}\le(\trchi-\trchio\ri)}_{L^\infty_{u,\ub}\HHt} & \les \varep.
\end{align}

Arguing similarly, commuting by $(r\Nd)$, $(u\Nd_3)$ or directly estimating equation~\eqref{eq:Nd4trchitrchio}, we obtain
\begin{align}\label{est:H12Ndtrchitrchio}
  \OOfgoodext[\trchi-\trchio] & \les \varep.
\end{align}


Arguing similarly, commuting first by $(r\Nd)$, 
we obtain
\begin{align}\label{est:H12Nd2trchitrchio}
  \OOfgoodext[(r\Nd)\trchi] & \les \varep.
\end{align}



Applying the $L^\infty_\ub L^2$ estimates of Lemma~\ref{lem:evolext} with $\kappa=2, \la=1$ to the transport equation~\eqref{eq:Nd4trchitrchio}, we obtain
\begin{align*}
  \norm{r\le(\trchi-\trchio\ri)}_{L^2_uL^2(S_{u,\ub})} & \les \norm{r(\trchi-\trchio)}_{L^2(\CCba)} \\
  & \quad + \int_{\ub}^\uba \norm{r\Err\le(\Nd_4,\trchi-\trchio\ri)}_{L^2_uL^2(S_{u,\ub'})}\,\d\ub'.
\end{align*}
Using the estimate~(\ref{est:initconnCCba}) on $\CCba$, we have
\begin{align*}
  \norm{r(\trchi-\trchio)}_{L^2(\CCba)} & \les \varep.
\end{align*}
We only treat the error terms involving $\chih$ and the other terms follow (more easily). Using the Bootstrap Assumptions~\ref{BA:connint} and Cauchy-Schwartz, gives
\begin{align*}
  \int_\ub^\uba \norm{\ub' |\chih|^2}_{L^2_uL^2(S_{u,\ub'})}\, \d\ub' & \les  \int_\ub^\uba (D\varep) (\ub')^{-1} \norm{u^{-1/2}\chih}_{L^2_uL^2(S_{u,\ub'})}\,\d\ub' \\
                                                                         & \les (D\varep) \norm{u^{-1/2}\ub^{-\ga}\chih}_{L^2(\MM^\ext)}\\
                                                                         & \les  (D\varep) \norm{\ub^{-1/2-\ga} \ub \chih}_{L^2(\MM^\ext)}\\
  & \les (D\varep)^2.
\end{align*}
Thus,
\begin{align*}
  \norm{r\le(\trchi-\trchio\ri)}_{L^\infty_{\ub}L^2_uL^2(S_{u,\ub})} & \les \varep,
\end{align*}
and we infer
\begin{align*}
  \norm{\ub^{-1/2-\ga}\ub \le(\trchi-\trchio\ri)}_{L^2(\MM^\ext)} & \les \varep
\end{align*}
for all $\ga>0$.\\

Commuting equation~\eqref{eq:Nd4trchitrchio} by $(r\Nd)^{\leq 1}$ and then by $\Ndt^{\leq 2}$ and arguing similarly then gives
\begin{align}\label{est:H12Nd3trchitrchio}
  \OOgoodext[\trchi-\trchio] & \les \varep, & \OOgoodext[(r\Nd)\trchi] & \les \varep,
\end{align}
for all $\ga>0$.

\section{Control of $\ze$}
\subsection{Control of $\mu-\muo$ and $\Nd^{\leq 1}\ze$}
Arguing as in the control for $\trchi-\trchio$ in Section~\ref{sec:controltrchitrchioext}, using that $\mu-\muo=0$ on $\CCba$ and the transport equation~(\ref{eq:Nd4mumuo}), and using the control obtained for $\trchi-\trchio$, we obtain using the $\HHt$ estimates of Lemma~\ref{lem:evolext} with $\kappa=3, \la=5/2$
\begin{align}\label{est:NdmumuoH12}
  \OOfgoodext[r(\mu-\muo)] & \les \varep,
\end{align}
and, using the $L^\infty_\ub L^2_uL^2(S_{u,\ub})$ estimates of Lemma~\ref{lem:evolext} with $\kappa=3, \la=5/2$ and arguing a previously,
\begin{align}\label{est:Nd2mumuoL2}
  \OOgoodext[r(\mu-\muo)] & \les \varep,
\end{align}
for all $\ga>0$.


Using the elliptic equation~(\ref{eq:Dd1etab}) and the elliptic estimates from Lemma~\ref{lem:ell}, one obtains from the $\HHt$ estimates~\eqref{est:NdmumuoH12} for $\mu-\muo$,
\begin{align}\label{est:NdNdzeH12}
  \begin{aligned}
    \OOfgoodext[\ze] & \les \varep, & \OOfgoodext[(r\Nd)\ze] & \les \varep,
  \end{aligned}
\end{align}
and from the $L^2(\MM^\ext)$ estimates~\eqref{est:Nd2mumuoL2} for $\mu-\muo$
\begin{align}\label{est:Nd2NdzeL2}
  \begin{aligned}
    \OOgoodext[\ze] & \les \varep, & \OOgoodext\le[(r\Nd)\ze\ri] & \les \varep,
  \end{aligned}
\end{align}
for all $\ga>0$.

\subsection{Control of $\io, \protect\ombr, \protect\ombs$}
Applying the $\HHt$ estimates of Lemma~\ref{lem:evolext} with $\kappa=0,\la=-1/2$ to the transport equation~\eqref{eq:Nd4io}, we have
\begin{align*}
  \norm{r^{-1/2}(r^2\io)}_{\HHt(S_{u,\ub})} & \les \norm{r^{-1/2}(r^2\io)}_{\HHt(S_{u,\uba})} + \int_\ub^\uba \norm{r^{-1/2}\Err(\Nd_4,\io)}_{\HHt(S_{u,\ub'})}\,\d\ub'.
\end{align*}
Using that $\io=\beb$ on $\CCba$ and the estimates for the curvature~(\ref{est:sourcecurvext}), we have
\begin{align*}
  \norm{r^{-1/2}(r^2\io)}_{\HHt(S_{u,\uba})} & \les \norm{r^{3/2}\beb}_{\HHt(S_{u,\uba})} \les \varep u^{-3/2}.
\end{align*}

From the error term $\Err(\Nd_4,\io)$, we only treat the term $r^2(\trchi-\trchio)$ and the estimates for the other terms will follow similarly. Using the product estimates from Lemma~\ref{lem:prodH12}, estimates~\eqref{est:sourcecurvext} and the Bootstrap Assumptions~\ref{BA:connext}, we have
\begin{align*}
  & \int_{\ub}^\uba \norm{r^{-1/2}r^2(\trchi-\trchio)\beb}_{\HHt(S_{u,\ub'})}\,\d\ub' \\
  \les & \; \int_\ub^\uba \le(\norm{\trchi-\trchio}_{L^\infty(S_{u,\ub'})}+r^{-1}\norm{(r\Nd)\trchi}_{L^2(S_{u,\ub'})}\ri)\norm{r^{-1/2}r^2\beb}_{\HHt(S_{u,\ub'})} \,\d\ub' \\
  \les & \; \int_\ub^\uba (D\varep)\ub^{-2}u^{-1/2}(D\varep)u^{-3/2}\,\d\ub'\\
  \les & (D\varep)^2u^{-3/2},
\end{align*}
Thus,
\begin{align}\label{est:ioH12}
  \norm{r^{3/2}u^{3/2}\io}_{L^\infty_{u,\ub}\HHt} & \les \varep.
\end{align}

Commuting equation~\eqref{eq:Nd4io} by $(r\Nd)$, $(u\Nd_3)$ and $\ub\Nd_4$, one further obtains
\begin{align}\label{est:NdioH12}
  \OOfbadext[r\io] & \les \varep.
\end{align}

Using the $L^\infty_\ub L^2_uL^2(S_{u,\ub})$ estimates of Lemma~\ref{lem:evolext}, we have
\begin{align*}
  \norm{(\ub^{-1} u)(r^2\io)}_{L^2_uL^2(S_{u,\ub})} & \les \norm{(\ub^{-1} u)(r^2\io)}_{L^2(\CCba)} + \int_\ub^\uba \norm{((\ub')^{-1} u)\Err(\Nd_4,\io)}_{L^2_uL^2(S_{u,\ub'})}\,\d\ub'.
\end{align*}
From~\eqref{est:initconnCCba}, we have
\begin{align*}
  \norm{(\ub^{-1} u)(r^2\io)}_{L^2(\CCba)} & \les \varep.
\end{align*}
From estimates~\eqref{est:sourcecurvext} and the Bootstrap Assumptions~\ref{BA:connext} and using Cauchy-Schwartz, we have
\begin{align*}
 \int_\ub^\uba \norm{((\ub')^{-1} u)r^2\le(\trchi-\trchio\ri)\beb}_{L^2_uL^2(S_{u,\ub'})}\,\d\ub' & \les \int_\ub^\uba \norm{(\ub')^{-1} u r^2\beb(\trchi-\trchio)}_{L^2_uL^2(S_{u,\ub'})}\,\d\ub' \\
  & \les  \int_\ub^\uba (\ub')^{-3/2+\ga}\norm{(\ub')^{-1/2-\ga}\ub' u^{1/2} \beb}_{L^2_uL^2(S_{u,\ub'})}\,\d\ub' \\
  & \les \norm{\ub^{-1/2-\ga} \ub u^{1/2} \beb}_{L^2(\MM^\ext)} \\
  & \les \norm{\ub^{-1/2-\ga} \ub u \beb}_{L^2(\MM^\ext)} \\
     & \les (D\varep)^2. 
\end{align*}
Thus,
\begin{align*}
  \norm{(\ub^{-1} u)(r^2\io)}_{L^\infty_\ub L^2_uL^2(S_{u,\ub})} & \les \varep.
\end{align*}
Commuting with $\Ndt^{\leq 2}$ and arguing as previously, we infer
\begin{align}\label{est:Nd2ioL2}
  \OObadext[r\io] & \les \varep,
\end{align}
for all $\ga>0$.\\

Using the elliptic equation~(\ref{eq:defio}), the elliptic estimates of Lemma~\ref{lem:ell}, the estimates for the curvature~(\ref{est:sourcecurvext}), we have using the $\HHt$ estimates~\eqref{est:NdioH12}
\begin{align}\label{est:NdNdombrsH12}
  \OOfbadext\le[(r\Nd)(\ombr,\ombs)\ri] & \les \varep,
\end{align}
and using the $L^2(\MM^\ext)$ estimates~\eqref{est:Nd2ioL2}
\begin{align}\label{est:Nd2NdombrsL2}
  \begin{aligned}
    \OObadext\le[(r\Nd)(\ombr,\ombs)\ri] & \les \varep,
  \end{aligned}
\end{align}
for all $\ga>0$.

\subsection{Control of $\Pso$ and $\Nd_3\ze$}
Applying the estimates from Lemma~\ref{lem:evolext} with $\kappa=0, \lambda=-1/2$ to the transport equation~(\ref{eq:Nd4Pso}) for $\Pso$, we have
\begin{align*}
  \norm{r^{-1/2}r^2\Pso}_{\HHt(S_{u,\ub})} & \les \norm{r^{-1/2}r^2\Pso}_{\HHt(S_{u,\uba})} \\
                                           & \quad + \int_\ub^\uba \le(\norm{r^{-1/2}r^2\trchib\be}_{\HHt(S_{u,\ub'})} + \norm{r^{-1/2}r^2\trchi\trchib\ze}_{\HHt(S_{u,\ub'})}\ri)\,\d\ub' \\
  & \quad + \int_\ub^\uba \norm{r^{-1/2}\Err(\Nd_4,\Pso)}_{\HHt(S_{u,\ub'})}\,\d\ub'.
\end{align*}
From~(\ref{est:initconnCCba}), we have
\begin{align*}
  \norm{r^{-1/2}r^{2}\Pso}_{\HHt(S_{u,\uba})} & = \norm{r^{-1/2}r^{2}\Nd_3\ze}_{\HHt(S_{u,\uba})} \les \varep u^{-3/2}.
\end{align*}
From the curvature estimates~(\ref{est:sourcecurvext}), the Bootstrap Assumptions~\ref{BA:connext} for $\trchi,\trchib$ and the previously obtained estimates~(\ref{est:NdNdzeH12}) for $\ze$, we have
\begin{align*}
  & \int_\ub^\uba \le(\norm{r^{-1/2}r^2\trchib\be}_{\HHt(S_{u,\ub'})} + \norm{r^{-1/2}r^2\trchi\trchib\ze}_{\HHt(S_{u,\ub'})}\ri)\,\d\ub'\\
  \les & \; \int_{\ub}^\uba \varep(\ub')^{-5/2} + \varep (\ub')^{-2}u^{-1/2} \,\d\ub' \\
  \les & \; \varep u^{-3/2}.
\end{align*}

From the product estimates from Lemma~\ref{lem:prodH12}, estimates~\eqref{est:sourcecurvext} and the Bootstrap Assumptions~\ref{BA:connext} and an inspection of the nonlinear terms composing $\Err(\Nd_4,\Pso)$, we have
\begin{align*}
  & \int_\ub^\uba \norm{r^{-1/2}\Err(\Nd_4,\Pso)}_{\HHt(S_{u,\ub'})}\,\d\ub' \\
  \les & \; \int_\ub^\uba (D\varep)^2(\ub')^{-2}u^{-3/2}\,\d\ub' \\
  \les & \; (D\varep)^2u^{-3/2}.
\end{align*}
Thus,
\begin{align}\label{est:PsoH12}
  \norm{r^{-1/2}u^{3/2}r^2\Pso}_{L^\infty_{u,\ub}\HHt} & \les \varep.
\end{align}

Commuting by $(r\Nd)$, $(u\Nd_3)$ and arguing similarly, or directly estimating $\Nd_4\Pso$, we have
\begin{align}\label{est:NdPsoH12}
  \OOfgoodext[u\Pso] & \les \varep.
\end{align}


Applying the $L^\infty_\ub L^2$ estimates of Lemma~\ref{lem:evolext} with $\kappa=0,\la=-1$ to the transport equation~\eqref{eq:Nd4Pso} multiplied with $u^{1-\ga}$ gives
\begin{align*}
  & \norm{r^{-1} u^{1-\ga} r^2\Pso}_{L^2_uL^2(S_{u,\ub})} \\
  \les & \; \norm{r^{-1}u^{1-\ga} r^2\Pso}_{L^2(\CCba)} \\
  & + \int_\ub^\uba \le(\norm{r^{-1}u^{1-\ga}r^2(\trchib\be)}_{L^2_uL^2(S_{u,\ub'})} + \norm{r^{-1}u^{1-\ga}r^2(\trchi\trchib\ze)}_{L^2_uL^2(S_{u,\ub'})}\ri)\,\d\ub' \\
  & + \int_\ub^\uba\norm{r^{-1}u^{1-\ga}\Err(\Nd_4,\Pso)}_{L^2_uL^2(S_{u,\ub'})}\,\d\ub'.
\end{align*}

From the estimates~(\ref{est:initconnCCba}) for $\Nd_3\ze$ on $\CCba$, we have
\begin{align*}
   \norm{r^{-1}u^{1-\ga} r^2\Pso}_{L^2(\CCba)} & = \norm{u^{1-\ga}r\Nd_3\ze}_{L^2(\CCba)} \les \varep.
\end{align*}

Using Cauchy-Schwartz, the estimates~(\ref{est:Nd2NdzeL2}) for the $L^2(\MM^\ext)$ norm of $\ze$, and the Bootstrap Assumptions~\ref{BA:connext}, we have
\begin{align*}
  \int_{\ub}^\uba \norm{r^{-1}u^{1-\ga}r^2\trchi\trchib\ze}_{L^2_uL^2(S_{u,\ub'})}\,\d\ub'  &\les \int_{\ub}^\uba \norm{r^{-1}u^{1-\ga}\ze}_{L^2_uL^2(S_{u,\ub'})} \, \d\ub' \\
  &\les \int_\ub^\uba \ub^{-1/2-\ga/2} \norm{\ub^{-1/2-\ga/2}\ub^{2} r^{-1}\ze}_{L^2_uL^2(S_{u,\ub'})}\,\d\ub'\\
  &\les \norm{\ub^{-1/2-\ga/2}\ub^{2} r^{-1}\ze}_{L^2(\MM^\ext)}\\
  &\les \varep.
\end{align*}

The linear source term $\be$ and the error term $\Err(\Nd_4,\Pso)$ are estimated similarly, using the curvature estimate~(\ref{est:sourcecurvext}) and the Bootstrap Assumptions~\ref{BA:connext}, and we have
\begin{align*}
\int_{\ub}^\uba \norm{r^{-1}u^{1-\ga}r^2\trchib\be}_{L^2_uL^2(S_{u,\ub'})}\,\d\ub'\les \varep.
\end{align*}
Thus,
\begin{align*}
  \norm{r^{-1}u^{1-\ga}r^2\Pso}_{L^\infty_\ub L^2_uL^2(S_{u,\ub})} & \les \varep,
\end{align*}
for all $\ga >0$, from which we also deduce
\begin{align*}
  \norm{\ub^{-1/2-\ga}\ub(u\Pso)}_{L^2(\MM^\ext)} & \les \varep,
\end{align*}
for all $\ga>0$. Commuting by $\Ndt^{\leq 2}$ and arguing similarly further gives
\begin{align}\label{est:Nd2PsoL2}
  \OOgoodext\le[u\Pso\ri] & \les \varep.
\end{align}

Using equation~(\ref{eq:defPso}) for $\Nd_3\ze$, the $\HHt$ estimates~(\ref{est:NdNdombrsH12}) for $\Nd\ombr,\Nd\ombs$, the $\HHt$ estimates~(\ref{est:PsoH12}), (\ref{est:NdPsoH12}), for $\Pso$, we deduce
\begin{align}\label{est:NdNd3zeH12}
  \OOfgoodext\le[(u\Nd_3)\ze\ri] & \les \varep.
\end{align}


From the $L^2(\MM^\ext)$ estimates~(\ref{est:Nd2NdombrsL2}) for $\Nd\ombr,\Nd\ombs$, and the $L^2(\MM^\ext)$ estimates~(\ref{est:Nd2PsoL2}) for $\Pso$, we have
\begin{align}\label{est:Nd2Nd3zeL2}
  \begin{aligned}
    \OOgoodext\le[(u\Nd_3)\ze\ri] & \les \varep,
  \end{aligned}
\end{align}
for all $\ga>0$.


\section{Control of $\chih$}
Using the elliptic estimates of Lemma~\ref{lem:ell} with the elliptic equation~(\ref{eq:Divdchih}) for $\chih$, together with the improved estimates~(\ref{est:H12Ndtrchitrchio}), (\ref{est:H12Nd2trchitrchio}), (\ref{est:H12Nd3trchitrchio}) for $\Nd\trchi$, the improved estimates~(\ref{est:NdNdzeH12}), (\ref{est:NdNd3zeH12}), for $\ze$, the Bootstrap Assumptions~\ref{BA:connext} and the curvature estimates~\eqref{est:sourcecurvext}, we have
\begin{align}\label{est:NdchihH12}
  \OOfgoodext\le[\chih\ri] & \les \varep, & \OOfgoodext\le[(r\Nd)\chih\ri] & \les \varep,
\end{align}
and using the corresponding $L^2(\MM^\ext)$ estimates,
\begin{align}\label{est:NdchihL2}
  \OOgoodext\le[\chih\ri] & \les \varep, & \OOgoodext\le[(r\Nd)\chih\ri] & \les \varep,
\end{align}
for all $\ga>0$.\\


\section{Control of $\protect\trchib-\protect\trchibo$ and $\protect\chibh$}
Applying the $\HHt$ estimates of Lemma~\ref{lem:evolext} with $\kappa = 1, \lambda = 1/2$ to the (commuted by $\Ndt^{\leq 1}$) transport equation~\eqref{eq:Nd4trchibtrchibo} for $\trchib-\trchibo$, using the improved estimates~(\ref{est:NdNdzeH12}) for $\Nd\ze$ and the estimates~\eqref{est:sourcecurvext} for $\rho-\rhoo$, we have
\begin{align}\label{est:trchibtrchiboH12}
  \OOfgoodext[\trchib-\trchibo] & \les \varep,
\end{align}
Using the corresponding $L^2(\MM^\ext)$ estimates, we also have
\begin{align}\label{est:trchibtrchiboL2}
  \OOgoodext\le[\trchib-\trchibo\ri] & \les \varep,
\end{align}
for all $\ga>0$.\\

Arguing similarly, using the transport equation~\eqref{eq:Nd4chibh} for $\chibh$, we have
\begin{align}\label{est:chibhH12}
  \OOfbadext\le[\chibh\ri] & \les \varep,
\end{align}
and
\begin{align}\label{est:chibhL2}
  \OObadext\le[\chibh\ri] & \les \varep,
\end{align}
for all $\ga>0$.

\section{Control of $\protect\omb-\protect\ombo$}
\subsection{Control of $\protect\omb-\protect\ombo$}
Applying the $\HHt$ estimates of Lemma~\ref{lem:evolext} with $\kappa=0, \lambda=-1/2$ to the transport equation~\eqref{eq:Nd4ombombo} for $\omb-\ombo$, we obtain
\begin{align}\label{est:NdombomboH12}
  \OOfbadext\le[\omb-\ombo\ri] & \les \varep,
\end{align}
and, applying the $L^2(\MM^\ext)$ estimates
\begin{align}\label{est:Nd2ombomboL2}
  \begin{aligned}
    \OObadext\le[\omb-\ombo\ri] & \les \varep,
  \end{aligned}
\end{align}
for all $\ga>0$.

\subsection{Control of $\protect\iob$ and $\Nd(\protect\omb,\protect\ombd)$}
Applying the $\HHt$ estimates of Lemma~\ref{lem:evolext} with $\kappa=0, \lambda=-1/2$ to the transport equation~\eqref{eq:Nd4iob} for $\iob$,
\begin{align}\label{est:iobH12}
  \OOfbadext\le[r\iob\ri] & \les \varep,
\end{align}
and applying the $L^2(\MM^\ext)$ estimates of Lemma~\ref{lem:evolext}, we have
\begin{align}\label{est:iobL2}
  \begin{aligned}
    \OObadext[r\iob] & \les \varep,
  \end{aligned}
\end{align}
for all $\ga>0$.\\

From the elliptic equation~(\ref{eq:defio}), the elliptic estimates of Lemma~\ref{lem:ell}, we deduce from~\eqref{est:iobH12} and~\eqref{est:iobL2} and the estimates~\eqref{est:sourcecurvext} for $\beb$ 
\begin{align}\label{est:NdombombdH12}
  \OOfbadext\le[(r\Nd)\omb\ri] & \les \varep, & \OOfbadext\le[(r\Nd)\ombd\ri] & \les \varep,
\end{align}
and
\begin{align}\label{est:Nd2NdombombdL2}
  \OObadext\le[(r\Nd)\omb\ri] & \les \varep, & \OObadext\le[(r\Nd)\ombd\ri] & \les \varep,
\end{align}
for all $\ga>0$.


\section{Control of $\protect\xib$}
We rewrite equation~\eqref{eq:Nd3ze} together with the results of Lemma~\ref{lem:relgeodnull} as
\begin{align}\label{eq:trchibxib}
  \half\trchi\xib & = \Nd_3\ze+\beb +2\Nd\omb +2\chib\cdot\ze - \chih\cdot\xib.
\end{align}

Estimating directly equation~\eqref{eq:trchibxib} together with the $\HHt$ estimates~(\ref{est:NdNd3zeH12}) for $\Nd_3\ze$, (\ref{est:NdombombdH12}) for $\Nd\omb$, and the estimates~\eqref{est:sourcecurvext} for $\beb$, we have
\begin{align}\label{est:xibH12}
  \begin{aligned}
    \OOfbadext\le[\xib\ri] & \les \varep.
  \end{aligned}
\end{align}

Using the $L^2(\MM^\ext)$ estimates~(\ref{est:Nd2Nd3zeL2}), (\ref{est:Nd2NdombombdL2}) and~\eqref{est:sourcecurvext}, we further have
\begin{align}\label{est:xibL2}
  \begin{aligned}
    \OObadext\le[\xib\ri] & \les \varep,
  \end{aligned}
\end{align}
for all $\ga>0$.

\section{Control of $\yy$}
\subsection{Control of $\yyo$}\label{sec:controlyyo}
Applying the $L^\infty$ estimate of Lemma~\ref{lem:evolext} with $\kappa=0, \la=0$, to the transport equation
\begin{align}\label{eq:Nd4yyo}
  \Nd_4\yyo = -4\ombo +\overline{(\trchi-\trchio)(\yy-\yyo)},
\end{align}
obtained by taking the average in~(\ref{eq:Nd4yy}) and using formula~\eqref{eq:commelbov}, we have
\begin{align*}
  \norm{\yyo}_{L^\infty(S_{u,\ub})} & \les \int_{\ub}^\uba\le(\norm{\ombo}_{L^\infty(S_{u,\ub'})} + \norm{(\trchi-\trchio)(\yy-\yyo)}_{L^\infty(S_{u,\ub'})}\ri) \,\d\ub', 
\end{align*}
where we used that $\yy=0$ on $\CCba$. From the bounds~(\ref{est:omboLinf}) for $\ombo$ and from the Bootstrap Assumptions~\ref{BA:connext}, we have
\begin{align*}
  \int_{\ub}^\uba\norm{\ombo}_{L^\infty(S_{u,\ub'})}\,\d\ub' & \les \varep\int_{\ub}^\uba (\ub')^{-2}u^{-2}\,\d\ub' \les \varep \ub^{-1}u^{-2},\\
  \int_\ub^\uba \norm{(\trchi-\trchio)(\yy-\yyo)}_{L^\infty(S_{u,\ub'})} \,\d\ub' & \les (D\varep)^2\int_{\ub}^\uba (\ub')^{-2}u^{-1/2}u^{-3/2}\,\d\ub' \les \varep \ub^{-1}u^{-2}.
\end{align*}
Thus,
\begin{align*}
  \norm{\ub u^2\yyo}_{L^\infty(\MM^\ext)} \les \varep.
\end{align*}
Commuting equation~\eqref{eq:Nd4yyo} by $\Ndt\in\{r\Nd,u\Nd_3,\ub\Nd_4\}$ and arguing similarly, we obtain
\begin{align}\label{est:NdyyoLinf}
  \norm{\ub u^2\Ndt\yyo}_{L^\infty(\MM^\ext)} & \les \varep.
\end{align}

Applying the $L^\infty_\ub L^2_uL^2(S_{u,\ub})$ estimate of Lemma~\ref{lem:evolext} with $\kappa =0,\la=-\ga$ to equation~\eqref{eq:Nd4yyo} multiplied by $u^{3/2}$, using Cauchy-Schwartz, the $L^2(\MM^\ext)$ estimates~(\ref{est:omboL2}) for $\ombo$ and the Bootstrap Assumptions~\ref{BA:connext}, we have
\begin{align*}
  \norm{u^{3/2}\ub^{-\ga}\yyo}_{L^2_uL^2(S_{u,\ub})} & \les \int^\uba_\ub \le(\norm{u^{3/2}(\ub')^{-\ga}\ombo}_{L^2_uL^2(S_{u,\ub'})} + \norm{u^{3/2}(\ub')^{-\ga}(\trchi-\trchio)(\yy-\yyo)}_{L^2_uL^2(S_{u,\ub'})}\ri)\,\d\ub' \\
                                                 & \les \norm{u^{3/2}\ub^{-\ga/2}\ub^{1/2}\ombo}_{L^2(\MM^\ext)} + (D\varep)\norm{u\ub^{-3/2-\ga/2}(\yy-\yyo)}_{L^2(\MM^\ext)} \\
  & \les_\ga \varep.
\end{align*}
Thus,
\begin{align*}
  \norm{\ub^{-1/2-\ga}u^{3/2}\yyo}_{L^2(\MM^\ext)} & \les \varep,
\end{align*}
for all $\ga>0$. Commuting with $\Ndt^{\leq 2}$ and arguing similarly, we obtain
\begin{align}\label{est:Nd2yyoLinf}
  \norm{\ub^{-1/2-\ga}u^{3/2}\Ndt^{\leq 2}\yyo}_{L^2(\MM^\ext)} & \les \varep.
\end{align}

To prove the estimates on $\TT$, we write the $L^2(\TT)$ norm of $\yyo$ of Section~\ref{sec:normsnullconn} as
\begin{align*}
  \int_1^{\cc\uba} u^{3-2\ga} \yyo^2 u^2\,\d u & \les \int_1^{\cc\uba}u^{3-2\ga}\le(\int_{\cc^{-1}u}^\uba |\Nd_4\yyo|\,\d\ub'\ri)^2u^2\,\d u \\
                                              & \les \int_1^{\cc\uba}\int_{\cc^{-1}u}^\uba u^{3-2\ga} (\ub')^{-2\al}|\Nd_4\yyo|^2 u^{2\al+1}u^2\,\d\ub'\d u,\\
                                              & \les \norm{u^{3-\ga+\al}\ub^{-1-\al}\Nd_4\yyo}_{L^2(\MM^\ext)}^2 \\
                                              & \les \norm{u^{3/2}\ub^{1/2-\ga}\Nd_4\yyo}_{L^2(\MM^\ext)}^2,
\end{align*}
where we used that $\yy=0$ on $\CCba$ and where in the second line we chose $\al$ such that $2\al+1<0$ and in the third line we chose $\al$ such that $\al+\frac{3}{2}-\ga >0$.\\

Using equation~\eqref{eq:Nd4yyo} for $\Nd_4\yyo$, the $L^2(\MM^\ext)$ bounds~(\ref{est:omboL2}) for $\ombo$ and the Bootstrap Assumptions, we have
\begin{align*}
  \norm{u^{3/2}\ub^{1/2-\ga}\Nd_4\yyo}_{L^2(\MM^\ext)} & \les \norm{u^{3/2}\ub^{1/2-\ga}\ombo}_{L^2(\MM^\ext)} + \norm{u^{3/2}\ub^{1/2-\ga}(\trchi-\trchio)(\yy-\yyo)}_{L^2(\MM^\ext)} \\
                                                           & \les \varep + (D\varep)\norm{u\ub^{-3/2-\ga}(\yy-\yyo)}_{L^2(\MM^\ext)} \\
                                                           & \les \varep + (D\varep)^2 \\
                                                           & \les \varep.
\end{align*}
Thus,
\begin{align*}
  \norm{u^{3/2-\ga}\yyo}_{L^2(\TT)} & \les \varep,
\end{align*}
for all $\ga>0$. Commuting by $\Ndt^{\leq 2}$ and arguing similarly further gives the desired estimate on $\TT$
\begin{align}\label{est:yyoL2TT}
  \norm{u^{3/2-\ga}\Ndt^{\leq 2}\yyo}_{L^2(\TT)} & \les \varep.
\end{align}

\subsection{Control of $\Nd\yy$}
From formula~(\ref{eq:Ndyy}), we have $\Nd\yy = -2\xib$. Compiling the bounds for $\yyo$ obtained in Section~\ref{sec:controlyyo} and the bounds~\eqref{est:xibH12} for $\xib$, one directly deduces
\begin{align}\label{est:yyMMext}
  \OOfext[\yy] & \les \varep, & \OOext[\yy] & \les \varep.
\end{align}

Arguing as in Section~\ref{sec:controlyyo}, we have
\begin{align*}
  \norm{u^{1-\ga}\xib}_{L^2(\TT)} & \les \norm{u\ub^{1/2-\ga}\Nd_4\xib}_{L^2(\MM^\ext)},
\end{align*}
where we used that $\xib=0$ on $\CCba$. Rewriting equation~(\ref{eq:Nd3etab}), we have the following equation for $\Nd_4\xib$.
\begin{align}\label{eq:Nd4xibyy}
  \Nd_4\xib & = -\Nd_3\ze-\beb -2\chib\cdot\ze.
\end{align}

Thus,
\begin{align*}
  \norm{u^{1-\ga}\xib}_{L^2(\TT)} & \les \norm{u\ub^{1/2-\ga}\Nd_3\ze}_{L^2(\MM^\ext)} +  \norm{u\ub^{1/2-\ga}\beb}_{L^2(\MM^\ext)} + \norm{u\ub^{-1/2-\ga}\chib\cdot\ze}_{L^2(\MM^\ext)}. 
\end{align*}
Using the $L^2(\MM^\ext)$ bounds~(\ref{est:Nd2Nd3zeL2}) for $\Nd_3\ze$, the curvature bounds (\ref{est:sourcecurvext}) for $\beb$, the bounds~(\ref{est:Nd2NdzeL2}) for $\ze$ and the Bootstrap Assumptions~\ref{BA:connext}, one infers
\begin{align*}
  \norm{u^{1-\ga}\xib}_{L^2(\TT)} & \les \varep.
\end{align*}
Commuting equation~\eqref{eq:Nd4xibyy} by $\Ndt^{\leq 2}$ and arguing similarly, we finally obtain the desired $L^2(\TT)$ bounds for (derivatives of) $\xib = -\half\Nd\yy$, and we have
\begin{align}\label{est:OOTTyy}
  \OO^\TT_{\leq 2,\ga}[\yy] & \les \varep,
\end{align}
for all $\ga>0$, as desired.

\section{Control of $\OOE$}
\subsection{Mild control of $\OOE$}\label{sec:mildOOEext}
We rewrite equation~(\ref{eq:Nd4OOE}) under the following form
\begin{align}\label{eq:Nd4OOENEW}
  \Nd_4\le(r^{-1}\OOO\ri) & = r^{-1}\chih\cdot \OOO  + \half r^{-1}(\trchi-\trchio)\OOO.
\end{align}

Applying the $L^\infty$ estimates of Lemma~\ref{lem:evolext} with $\kappa=0,\la=0$, to the transport equation~\eqref{eq:Nd4OOENEW}, using the improvement of the mild Bootstrap Assumptions for $\OOE$ on $\CCba\cap\MM^\ext$ obtained in Section~\ref{sec:mildOOE}, we obtain the following improvement of the mild Bootstrap Assumptions~\ref{BA:mildOOE} on the full domain $\MM^\ext$
\begin{align}\label{est:mildOOEext}
  \norm{r^{-1}(r\Nd)^{\leq 1}\OOE}_{L^\infty(\MM^\ext)} & \les 1.
\end{align}

Arguing similarly as in Section~\ref{sec:mildOOE}, one also obtains the last estimates of the Bootstrap Assumptions~\ref{BA:mildOOE}. This finishes the improvement of~\ref{BA:mildOOE} on $\MM^\ext$.




\subsection{Control of $\Hrot$ and $\POE$}
Applying the $\HHt$ estimates of Lemma~\ref{lem:evolext} with $\kappa=0, \la =-1/2$ to the transport equation~\eqref{eq:Nd4H} for $\Hrot$, we have
\begin{align*}
  \norm{r^{-1/2}\Hrot}_{\HHt(S_{u,\ub})} & \les \norm{r^{-1/2}\Hrot}_{\HHt(S_{u,\uba})} + \int_\ub^\uba\norm{r^{-1/2}\Nd_4\Hrot}_{\HHt(S_{u,\ub'})} \,\d\ub'.
\end{align*}
From estimates~\eqref{est:initconnCCba}, we have
\begin{align*}
  \norm{r^{-1/2}\Hrot}_{\HHt(S_{u,\uba})} & \les \varep.
\end{align*}
Using the product estimates of Lemma~\ref{lem:prodH12} and the estimates~(\ref{est:NdchihH12}), we check that
\begin{align*}
  \int_\ub^\uba\norm{r^{-1/2}\Nd\chih\cdot \OOO}_{\HHt(S_{u,\ub'})} \,\d\ub' & \les \int_\ub^\uba\norm{r^{-1/2}(r\Nd)\chih}_{\HHt(S_{u,\ub'})} \,\d\ub' \les \varep u^{-1/2}\ub^{-1},
\end{align*}
and
\begin{align*}
  \int_\ub^\uba\norm{r^{-1/2}\chih\cdot\Nd\OOO}_{\HHt(S_{u,\ub'})} \,\d\ub' & \les \varep \int_\ub^\uba (\ub')^{-2}u^{-1/2} \norm{r^{-1/2}\Nd\OOO}_{\HHt(S_{u,\ub'})}\,\d\ub' \\
                                                                            & \les \varep \int_\ub^\uba (\ub')^{-2}u^{-1/2}\le(\norm{r^{-1}\OOO}_{L^\infty(S_{u,\ub'})} + \norm{\Nd\OOO}_{L^\infty(S_{u,\ub'})} + \norm{\POE}_{L^2(S_{u,\ub'})}\ri)\,\d\ub' \\
  & \les \varep \int_\ub^\uba (\ub')^{-2}u^{-1/2}\le(1 + D\varep\ri)\,\d\ub' \\
  & \les \varep u^{-1/2}\ub^{-1},
\end{align*}
and the control of the other terms composing $\Nd_4\Hrot$ follows similarly. Thus, we have
\begin{align*}
  \norm{r^{-1/2}u^{1/2}\ub^2r^{-1}\Hrot}_{L^\infty_{u,\ub} \HHt} & \les \varep.
\end{align*}
Commuting equation~\eqref{eq:Nd4H} with $(\Ndt)^{\leq 2}$, where $\Ndt\in\{(r\Nd),u\Nd_3\}$ and applying the $\HHt$ estimates of Lemma~\ref{lem:evolext}, or directly estimating the $\Nd_4$ derivatives from equation~\eqref{eq:Nd4H} gives similar bounds and we have
\begin{align}\label{est:LinfH12H}
  \OOfgoodext\le[r^{-1}\Hrot\ri] & \les \varep.
\end{align}

Applying the $L^\infty_\ub L^2_uL^2(S_{u,\ub})$ estimates of Lemma~\ref{lem:evolext} with $\kappa =0, \lambda = -\ga$ to the transport equation~\eqref{eq:Nd4H} for $\Hrot$ gives
\begin{align*}
  \norm{\ub^{-\ga}\Hrot}_{L^2_uL^2(S_{u,\ub})} & \les \norm{\ub^{-\ga}\Hrot}_{L^2(\CCba)} + \int^\uba_\ub \le(\norm{\ub^{-\ga}\Nd_4\Hrot}_{L^2_uL^2(S_{u,\ub'})}\ri)\,\d \ub'.
\end{align*}
From estimate~\eqref{est:initconnCCba} on $\CCba$, we have
\begin{align*}
  \norm{\ub^{-\ga}\Hrot}_{L^2(\CCba)} & \les \varep.
\end{align*}
From the $L^2(\MM^\ext)$ estimate~(\ref{est:NdchihL2}) for $\chih$ and Cauchy-Schwartz, we have
\begin{align*}
  \int^\uba_\ub \norm{(\ub')^{-\ga}\Nd\chih\cdot\OOO}_{L^2_uL^2(S_{u,\ub'})}\,\d\ub' & \les \int_{\ub}^\uba \norm{(\ub')^{-\ga}(r\Nd)\chih}_{L^2_uL^2(S_{u,\ub'})}\,\d\ub' \\
                                                                                                 & \les \norm{\ub^{1/2}\ub^{-\ga/2}(r\Nd)\chih}_{L^2(\MM^\ext)} \\
                                                                                                 & \les \varep.
\end{align*}
The other terms follow similarly. Thus, we deduce that
\begin{align*}
  \norm{\ub^{-1/2-\ga}\Hrot}_{L^2(\MM^\ext)} & \les \varep,
\end{align*}
for all $\ga>0$. Commuting by $(r\Nd)^{\leq 2}$ and arguing similarly, we obtain
\begin{align}\label{est:HextL2}
  \OOgoodext\le[r^{-1}\Hrot\ri] & \les \varep,
\end{align}
for all $\ga>0$.\\

The control of $\POE$ follows from similar estimates as the ones performed for $\Hrot$, using the transport equation~\eqref{eq:Nd4POE} for $r\POE$. Thus, we have
\begin{align}\label{est:POEext}
  \mathfrak{O}^{\ext,\mathfrak{g}}_{\leq 0}\le[\POE\ri] & \les \varep, & \OO^{\ext,\mathfrak{g}}_{\leq 1,\ga}\le[\POE\ri] & \les \varep,
\end{align}
for all $\ga>0$.



\subsection{Control of $Y$}
Applying the $\HHt$ estimates of Lemma~\ref{lem:evolext} with $\kappa = -1, \la = -3/2$ to the transport equation~(\ref{eq:Nd4YOOE}) for $Y$, using that $Y=0$ on $\CCba$, we have
\begin{align*}
  \norm{r^{-3/2}Y}_{\HHt(S_{u,\ub})} & \les \int_\ub^\uba \norm{r^{-3/2}\le(\Nd_4Y-\half\trchi Y\ri)}_{\HHt(S_{u,\ub'})} \,\d\ub'.
\end{align*}

From an inspection of the terms in the right-hand side of equation~\eqref{eq:Nd4YOOE} the term with slowest decay is $2(\Divd\ze) \OOO$, and we have
\begin{align*}
  \int_\ub^\uba \norm{r^{-3/2}(\Divd\ze)\OOO}_{\HHt(S_{u,\ub'})} \,\d\ub' & \les \int_\ub^\uba \norm{r^{-1/2}\Nd\ze}_{\HHt(s_{u,\ub'})}\,\d\ub'\\
                                                                           & \les \varep \int_\ub^\uba (\ub')^{-3}u^{-1/2} \,\d\ub' \\
  & \les \varep \ub^{-2}u^{-1/2},
\end{align*}
where we used the product estimates of Lemma~\ref{lem:prodH12}, the estimates~(\ref{est:NdNdzeH12}) for $\ze$ and~(\ref{est:mildOOEext}) for $\OOO$. Commuting with $(r\Nd), (\qq\Nd_3)$ or estimating directly equation~\eqref{eq:Nd4YOOE} and arguing similarly, we obtain
\begin{align*}
  \OOfgoodext\le[r^{-1}Y\ri] & \les \varep.
\end{align*}

Applying the $L^\infty_\ub L^2_uL^2(S_{u,\ub})$ estimates of Lemma~\ref{lem:evolext} with $\kappa = -1, \lambda = -\ga$, using that $Y= 0$ on $\CCba$, we have
\begin{align*}
  \norm{r^{-1-\ga} \ub Y}_{L^2_uL^2(S_{u,\ub})} & \les \int_\ub^\uba \norm{r^{-1-\ga}\ub'\le(\Nd_4Y-\half\trchi Y\ri)}_{L^2_uL^2(S_{u,\ub'})} \,\d\ub'.
\end{align*}
The right-hand side of the previous estimate can be bounded using the transport equation~\eqref{eq:Nd4YOOE}. We only treat the slowest decaying term $2(\Divd\ze)\OOO$ of~\eqref{eq:Nd4YOOE}, and we claim that the other terms from the right-hand side of the previous inequality are treated similarly
\begin{align*}
  \int_\ub^\uba \norm{r^{-1-\ga}\ub'\Divd\ze\OOO}_{L^2_uL^2(S_{u,\ub'})}\,\d\ub' & \les \int_\ub^\uba \norm{r^{-1-\ga}\ub'(r\Nd)\ze}_{L^2_uL^2(S_{u,\ub'})} \,\d\ub' \\
                                 & \les \norm{r^{-1/2-\ga/2}\ub(r\Nd)\ze}_{L^2(\MM^\ext)} \\
                                 & \les \varep.
\end{align*}
Therefore, we obtain
\begin{align*}
  \norm{\ub^{-1/2-\ga} \ub Y}_{L^2(\MM^\ext)} & \les \varep, 
\end{align*}
for all $\ga>0$. Commuting with $\Ndt^{\leq 2}$ and arguing similarly, we have
\begin{align*}
  \OOgoodext\le[r^{-1}Y\ri] & \les \varep,
\end{align*}
for all $\ga>0$. This finishes the proof of (\ref{est:OO1}).

\section{Control of the area radius}\label{sec:arearadiusestimate}
\begin{lemma}
  Under the improved estimates obtained in this section, we have in $\MM^\ext$
  \begin{align*}
    \le|r(u,\ub) -\half (\ub-u)\ri| & \les \varep \ub^{-1} u^{-1}.
  \end{align*}
\end{lemma}
\begin{proof}
  Using relations~\eqref{eq:elu} and~\eqref{eq:elbr}, we have
  \begin{align*}
    \el\le(r(u,\ub)-\half(\ub-u)\ri) & = \half r\trchio - 1.
  \end{align*}
  Integrating in $\ub$, using the improved estimate for $r$ on $\CCba$ from Section~\ref{sec:arearadiusestimateCCba} and using the improved estimate~(\ref{est:trchio}) for $\trchio$, we have
  \begin{align*}
    \le|r(u,\ub)-\half(\ub-u)\ri| & \les \le|r(u,\uba)-\half(\uba-u)\ri| + \int_\ub^\uba r \le|\trchio-\frac{2}{r}\ri|\,\d \ub' \\
                                  & \les \varep \uba^{-1}u^{-1} + \varep \int_\ub^\uba r (\ub')^{-3}u^{-1}\,\d\ub'\\
    & \les \varep \ub^{-1}u^{-1},
  \end{align*}
  as desired.
\end{proof}

\section{Control of the spherical coordinates in $\MM^\ext$}\label{sec:mildBAimpext}
For a spherical coordinate system as described in the Bootstrap Assumption~\ref{BA:mildsphcoordsext}, we have using~(\ref{eq:elbgd}), the following transport equation in the $\el$ direction
\begin{align*}
  \el(\gd_{ab}) & = 2\chi_{ab} = \trchio \gd_{ab} + (\trchi-\trchio)\gd_{ab} + 2\chih_{ab},
\end{align*}
where $a,b\in\{\varth,\varphi\}$, which rewrites using~\eqref{eq:elbr} and the notation of the Bootstrap Assumption~\ref{BA:mildsphcoordsext}
\begin{align}\label{eq:transportrenormgd}
  \el\le(r^{-2}\gd_{ab}-(\gd_{\SSS})_{ab}\ri) & = (\trchi-\trchio)r^{-2}\gd_{ab} + 2r^{-2}\chih_{ab}.
\end{align}
Integrating along $\el$, this yields
\begin{align*}
  \le|r^{-2}\gd_{ab}-(\gd_{\SSS})_{ab}\ri|(u,\ub,\varth,\varphi) & \les \le|r^{-2}\gd_{ab}-(\gd_{\SSS})_{ab}\ri|(u,\uba,\varth,\varphi) + \int_{\ub}^\uba \le(\le|(\trchi-\trchio)r^{-2}\gd_{ab}\ri| + 2r^{-2}\le|\chih_{ab}\ri|\ri)\,\d\ub' \\
                                                                 & \les \varep + \varep \int_{\ub}^\uba (\ub')^{-2}u^{-1/2}\,\d\ub' \\
                                                                 & \les \varep \ub^{-1}u^{-1/2},
\end{align*}
where we used the improved bounds for the metric in coordinates on $\CCba$ from Section~\ref{sec:mildBAimpast}, the improved bounds~(\ref{est:H12Ndtrchitrchio}), (\ref{est:NdchihH12}) for $\trchi$ and $\chih$ and the Bootstrap Assumptions~\ref{BA:mildsphcoordsext} for the coordinate component of $\gd$.\\

Commuting~\eqref{eq:transportrenormgd} with $\pr$ and integrating using the sup-norm bounds for $(r\Nd)(\trchi-\trchio)$ and $(r\Nd)\chih$, we further have
\begin{align*}
  \le|\pr^{\leq 1}\le(r^{-2}\gd_{ab}-(\gd_{\SSS})_{ab}\ri)\ri| & \les \varep \ub^{-1}u^{-1/2}.
\end{align*}
This finishes the improvement of the Bootstrap Assumptions~\ref{BA:mildsphcoordsext}.

\chapter{Maximal connection estimates in $\MM^\intr_\bott$}\label{sec:planehypconnest}
In this section, we prove the following proposition.
\begin{proposition}\label{prop:planehypconnestSTAB}
  Recall that from Proposition~\ref{prop:planehypcurvestSTAB} we have the following curvature control on the maximal hypersurfaces $\Si_t$ (see the definitions of Section~\ref{sec:normscurvint})
  \begin{align}\label{est:curvestsourceplanehyp}
    \RR^\intr_{\leq 2} + \RRfb^\intr_{\leq 1} & \les \varep.
  \end{align}
  Recall that from Proposition~\ref{prop:connestSTAB} we have the following control of the null connection coefficients in $\MM^\ext$ -- and in particular on the interface $\TT$ -- 
  \begin{align}\label{est:bdysourcenullconn}
    \OOfexti + \OOofb + \OOfexti\le[\yy\ri] + \OO^\TT_{2,\ga}\le[\yy\ri] & \les \varep,
  \end{align}
  for all $\ga>0$ (see the definitions from Section~\ref{sec:normsnullconn}). Recall that from Proposition~\ref{prop:connestCCbaSTAB}, we have on $\CCba$ (see the definitions of Section~\ref{sec:normnullconnCCba})
  \begin{align}
    \label{est:sourcechichibSast}
    \mathfrak{O}_{\leq 2}^\ast & \les \varep.
  \end{align}
  Recall that from Proposition~\ref{prop:connestSTAB} we have the following control of the exterior rotation vectorfields in $\MM^\ext$ -- and in particular on the interface $\TT$ -- (see the definitions from Section~\ref{sec:normsnullconn})
  \begin{align}\label{est:sourcerotext}
    \mathfrak{O}^{\ext,\OOO}_{\leq 2} & \les \varep.
  \end{align}

  Under the Bootstrap Assumptions and estimates~\eqref{est:curvestsourceplanehyp}, \eqref{est:bdysourcenullconn}, we have for $\varep>0$ sufficiently small the following control on the connection coefficients of the maximal hypersurfaces (see Section~\ref{sec:norminteriorconn} for definitions)
  \begin{align}\label{est:OOmaxconn}
    \OO^\intr_{\leq 3, \ga}[n] + \OO^\intr_{\leq 2}[k] + \mathfrak{O}^\TT_{\leq 2}[\nu] & \les_{\ga} \varep,
  \end{align}
  for all $\ga>0$. Moreover, under the same hypothesis and~\eqref{est:sourcechichibSast}, we have the following harmonic coordinates control on the last maximal slice $\Si_\tast$
  \begin{align}\label{est:sourceharmo}
    \begin{aligned}
      \sum_{i,j=1}^3 \norm{t^{3/2}\le(g(\nab x^i,\nab x^j)-\de_{ij}\ri)}_{L^\infty(\Si_\tast)} + \sum_{i=1}^3 \norm{(t\nab)^{\leq 3}\nab^2x^i}_{L^2(\Si_\tast)} & \les \varep,\\
      \sum_{i=1}^3\norm{t\le(\rast\Nf(x^i)-x^i\ri)}_{L^\infty(S^\ast)} & \les \varep.
    \end{aligned}
  \end{align}
  Under the same hypothesis, the respectively mild and strong Bootstrap Assumptions~\ref{BA:mildKillingMMintbot} and \ref{BA:intKill} for the interior approximate conformal Killing fields $\TI,\SI,\KI,\OOI$ are improved. Under the same hypothesis and using~\eqref{est:sourcerotext}, the Bootstrap Assumptions~\ref{BA:TTKilling} on the difference of the interior and exterior approximate Killing fields on $\TT$ are improved.
\end{proposition}

\begin{remark}
  We do not control all $3$-derivatives of $\nt$, which is a far-reaching consequence of the lack of regularity for $\xib$. However, we do control all $2$-derivatives of $\nab\nt$, which is enough to control $\D^{(\TI)}\pi$ and $\D^2{^{\TI}\pi}$ in the error term estimates of Section~\ref{sec:globener} and in the curvature estimates of Section~\ref{sec:planehypcurvest}, where it is used.
\end{remark}




\section{Elliptic estimates}
In Section~\ref{sec:estkn}, we will need the following elliptic estimates for the Laplace equation on $\Si_t$.
\begin{lemma}[Elliptic estimate for Laplace equation on $\Si_t$]\label{lem:ellLap}
  Under the Bootstrap Assumptions~\ref{BA:mildcoordsSit}, we have for all $1\leq t \leq \tast$ and for all scalar function $f$
  \begin{align*}
    \norm{(t\nab) f}_{L^2(\Si_t)} + \norm{f}_{L^2(\Si_t)} & \les \norm{t^2\Delta f}_{L^2(\Si_t)} + \norm{t f}_{\HHt(\pr\Si_t)},
  \end{align*}
  and for all $k\geq 2$
  \begin{align*}
    \norm{(t\nab)^{\leq k}f}_{L^2(\Si_t)} & \les \norm{t^2(t\nab)^{k-2}(\Delta f)}_{L^2(\Si_t)} + \norm{t(t\Nd)^{\leq k-1}f}_{\HHt(\pr\Si_t)}. 
  \end{align*}
\end{lemma}
\begin{proof}
  The proof follows from a rescaling in $t$ and the results from~\cite{Czi19} (see also~\cite{Czi.Gra19a}).
\end{proof}

\section{Control of the second fundamental form $k$ and lapse $n$}\label{sec:estkn}
\subsection{Control of $k, \nab k$ on $\Si_t$}\label{sec:knabk}
In this section, we use energy estimates for the elliptic div-curl systems~\eqref{eq:Hodgekt} satisfied by $\kt$ on $\Si_t$. The source term for this elliptic equation is $\Ht$ and the boundary conditions on $\pr\Si_t$ are given by mixed implicit Dirichlet and Neumann conditions for $\kt$ (see Lemma~\ref{lem:transkTT}). This estimate has already been done in~\cite[Sections 4.7 and 4.8]{Czi.Gra19a}, where, rescaling in $t$ these estimates, it holds that
\begin{align*}
  & \norm{t(t\nab)\kt}_{L^2(\Si_t)} + \norm{t\kt}_{L^2(\Si_t)} + t^{3/2}\norm{k}_{L^2(\pr\Si_t)} \\
   \les & \; \norm{t^2\Ht}_{L^2(\Si_t)} + \norm{t^2\le(\trchi-\frac{2}{r}\ri)}_{\HHt(\pr\Si_t)} + \norm{t^2\le(\trchib+\frac{2}{r}\ri)}_{\HHt(\pr\Si_t)} + \norm{t^2\ze}_{\HHt(\pr\Si_t)} +\text{error terms}.
\end{align*}

Using the $L^2(\Si_t)$ bounds~(\ref{est:curvestsourceplanehyp}) for $\Ht$ and the $\HHt(\pr\Si_t)$ bounds~(\ref{est:bdysourcenullconn}) for $\trchi-\frac{2}{r}, \trchib+\frac{2}{r}, \ze$, we thus obtain
\begin{align*}
  \norm{t(t\nab)^{\leq 1}\kt}_{L^2(\Si_t)} & \les \varep . 
\end{align*}

As a byproduct of these estimates, it is shown in~\cite[Section 4]{Czi.Gra19a} that we have the following control for $\nu$
\begin{align*}
  \norm{t(t\Nd)^{\leq 1}(\nut-1)}_{\HHt(\pr\Si_t)} & \les\varep.
\end{align*}
\begin{remark}
  The $\HHt(\pr\Si_t)$ control of $\nu$ uses the first two relations of~\eqref{eq:transkTTt}. See~\cite[Lemma 4.9]{Czi.Gra19a} for the full argument.
\end{remark}

\subsection{Control of $\nab^2k$ on $\Si_t$}\label{sec:nab2knab3k}
To estimate $\nab^2\kt$ on $\Si_t$, one cannot use the higher regularity elliptic estimates from Lemma~\ref{lem:ellSitHrr} since the boundary terms for $\kt$ are only implicitly related to the null connection coefficients on $\pr\Si_t$. We therefore have to show that commuting the elliptic equations~(\ref{eq:Hodgekt}) by a derivative does not change the coercive structure of the boundary terms which appeared in the original energy estimate. This analysis has already been carried out in~\cite[Sections 5.1.3 and 5.2.3]{Czi.Gra19a}, and we refer the reader to that paper.\\

We therefore shall assume that we have the following elliptic estimate on $\Si_t$
\begin{align*}
  & \norm{t^3\nab^2 \kt}_{L^2(\Si_t)} + \norm{t^{5/2}\nab k}_{L^2(\pr\Si_t)} \\
  \les & \; \norm{t^3\nab \Ht}_{L^2(\Si_t)} +\norm{t^2\nab\kt}_{L^2(\Si_t)} + \norm{t\kt}_{L^2(\Si_t)}+ \norm{t^3\Nd\trchi}_{\HHt(\pr\Si_t)} + \norm{t^3\Nd\trchib}_{\HHt(\pr\Si_t)} \\
  & + \norm{t^3\Nd\ze}_{\HHt(\pr\Si_t)} + \norm{t^2\le(\trchi-\frac{2}{r}\ri)}_{\HHt(\pr\Si_t)} + \norm{t^2\le(\trchib+\frac{2}{r}\ri)}_{\HHt(\pr\Si_t)} + \norm{t^2\ze}_{\HHt(\pr\Si_t)}.
\end{align*}

Using the $L^2(\Si_t)$ bounds~(\ref{est:curvestsourceplanehyp}) for $\Ht$ and the $\HHt(\pr\Si_t)$ bounds~(\ref{est:bdysourcenullconn}) for $\trchi-\frac{2}{r}, \trchib+\frac{2}{r}, \ze$, we thus obtain
\begin{align}\label{est:nab2ktL2}
  \norm{t^3\nab^2 \kt}_{L^2(\Si_t)} & \les \varep.
\end{align}

As a byproduct of these estimates, one further has the following improved estimates for the tangential derivatives of $\nut$
\begin{align}\label{est:Nd2nutH12}
  \norm{t(t\Nd)^{\leq 2}(\nut-1)}_{\HHt(\pr\Si_t)} & \les\varep.
\end{align}

\subsection{Optimal control for $\nab^3k$ on the last hypersurface $\Si_\tast$}\label{sec:nab3kSitast}
On the boundary of the last hypersurface $\pr\Si_\tast = S^\ast \subset \CCba$, we have from $\OOfast \les \varep$ (see the definitions in Section~\ref{sec:normnullconnCCba} and the improvement in Section~\ref{sec:connestCCba}) the following additional bounds
\begin{align}\label{est:addbdysourceSitast}
  \norm{t^{2}(t\Nd)^2\le(\trchi-\frac{2}{r}\ri)}_{\HHt(S^\ast)} + \norm{t^{2}(t\Nd)^2\le(\trchib+\frac{2}{r}\ri)}_{\HHt(S^\ast)} + \norm{t^{2}(t\Nd)^2\ze}_{\HHt(S^\ast)} & \les \varep.
\end{align}

From the elliptic estimates of~\cite[Sections 5.1.3 and 5.2.3]{Czi.Gra19a}, we have on the last slice $\Si_\tast$
\begin{align*}
    & \norm{t^4\nab^3 \kt}_{L^2(\Si_\tast)} + \norm{t^{7/2}\nab^2 k}_{L^2(\pr\Si_\tast)} \\
  \les & \; \norm{t^4\nab^2 \Ht}_{L^2(\Si_\tast)} +\norm{t^3\nab^2\kt}_{L^2(\Si_\tast)} +\norm{t^2\nab\kt}_{L^2(\Si_\tast)} + \norm{t\kt}_{L^2(\Si_\tast)}\\
    & + \norm{t^4\Nd^2\trchi}_{\HHt(S^\ast)}  + \norm{t^4\Nd^2\trchib}_{\HHt(S^\ast)} + \norm{t^4\Nd^2\ze}_{\HHt(S^\ast)} \\
    & + \norm{t^3\Nd\trchi}_{\HHt(S^\ast)}  + \norm{t^3\Nd\trchib}_{\HHt(S^\ast)} + \norm{t^3\Nd\ze}_{\HHt(S^\ast)} \\
  & + \norm{t^2\le(\trchi-\frac{2}{r}\ri)}_{\HHt(S^\ast)} + \norm{t^2\le(\trchib+\frac{2}{r}\ri)}_{\HHt(S^\ast)} + \norm{t^2\ze}_{\HHt(S^\ast)}.
\end{align*}

Thus, using the curvature estimates~(\ref{est:curvestsourceplanehyp}) and the additional bounds~\eqref{est:addbdysourceSitast}, we obtain
\begin{align*}
  \norm{t^4\nab^3 \kt}_{L^2(\Si_\tast)} & \les \varep.
\end{align*}

\begin{remark}
  The optimal control in regularity for $k$, $\nab^3k \in L^2$, can only be achieved provided that $\trchi$, $\trchib$ and $\ze$ have optimal regularity on the boundary. The optimal regularity for these coefficients is only obtained on the last cone $\CCba$ and is a consequence of the choice of the canonical foliation (see Section~\ref{sec:connestCCba}). The optimal regularity for $\chib$ in $\MM^\ext\setminus\CCba$ cannot be obtained due to the classical loss of regularity for the geodesic foliation on the cones $\CC_u$. In the present paper, this is not an issue since the optimal bounds for $k$ are only needed on $\Si_\tast$ to carry out the extension procedure of $\Si_\tast$ in Section~\ref{sec:tipcurvest}. 
\end{remark}

\subsection{Control of $\nab^{\leq 3}n$ on $\Si_t$}\label{sec:nab3n}
From relation~\eqref{eq:relntnut}, we have on $\TT$
\begin{align}\label{eq:n-1}
  n-1 & = \frac{\cc}{\cc+1}\le((\nu-1) + \cc^{-1}(\nu^{-1}-1) - \half \yy\nu^{-1}\ri).
\end{align}
From the assumption~(\ref{est:bdysourcenullconn}) (see the definitions of Section~\ref{sec:normsnullconn} and the improvements of Section~\ref{sec:connest}), and from the just obtained~\eqref{est:Nd2nutH12}, we have
\begin{align*}
  \norm{t(t\Nd)^{\leq 2}\yy}_{\HHt(\pr\Si_t)} & \les \OOfexti[\yy] \les \varep,\\
  \norm{t(t\Nd)^{\leq 2}(\nut-1)}_{\HHt(\pr\Si_t)} & \les \varep.
\end{align*}
Thus, using relation~\eqref{eq:n-1} and the product estimates of Lemma~\ref{lem:prodH12}, we have
\begin{align}\label{est:H12n-1}
  \norm{t(t\Nd)^{\leq 2}(n-1)}_{\HHt(\pr\Si_t)} & \les \varep.
\end{align}
Applying the elliptic estimates of Lemma~\ref{lem:ellLap}, using Laplace equation~\eqref{eq:Lapnt} for $\nt$, the Bootstrap Assumptions~\ref{BA:connint} and the boundary estimate~\eqref{est:H12n-1}, we therefore deduce
\begin{align}\label{est:nab3L2n-1}
  \norm{(t\nab)^{\leq 3}(n-1)}_{L^2(\Si_t)} & \les \varep,  
\end{align}
as desired.

\subsection{Control of $\Lieh_\Tf\kt, \nab\Lieh_\Tf\kt$ on $\Si_t$}\label{sec:TkSit}
The improved estimates for $\Lieh_\Tf\kt, \nab\Lieh_\Tf\kt$ on $\Si_t$ are directly obtained using equation~\eqref{eq:LiehTfkt} for $\Lieh_\Tf\kt$ together with the assumed estimates~(\ref{est:curvestsourceplanehyp}) for the curvature and the just improved estimates~\eqref{est:nab3L2n-1} for $\nab^{\leq 3}n$ from Section~\ref{sec:nab3n}, and we have
\begin{align}\label{est:LiehTfktL2}
  \begin{aligned}
    \norm{t^2(t\nab)^{\leq 1}\Lieh_\Tf\kt}_{L^2(\Si_t)} & \les \norm{t^2(t\nab)^{\leq 1}E}_{L^2(\Si_t)} + \norm{t^2(t\nab)^{\leq 1}\nab^2n}_{L^2(\Si_t)} + \text{nonlinear error terms}\\
    & \les \varep.
  \end{aligned}
\end{align}

\subsection{Control of $\nab^{\leq 2}(\Tf(\nt))$ on $\Si_t$}\label{sec:TntSit}
On the boundary $\pr\Si_t$, we have
\begin{align*}
  Z & = \le(a^{1/2}\nut^{-1}+a^{-1/2}\nut\ri) \Tf + \le(a^{1/2}\nut^{-1}-a^{-1/2}\nut\ri)\Ntf,
\end{align*}
where $Z$ is the future-pointing unit normal to $\pr\Si_t$ in $\TT$, defined in the proof of Lemma~\ref{lem:relntnut}, where $a = \cc^{-1}-\half\yy$ and where we recall that $\Ntf$ is the outward-pointing unit normal to $\pr\Si_t$ in $\Si_t$.\\

We therefore have
\begin{align}\label{eq:TnZnNn}
  \le(a^{1/2}\nut^{-1}+a^{-1/2}\nut\ri)\Tf(\nt) & = Z(\nt) - \le(a^{1/2}\nut^{-1}-a^{-1/2}\nut\ri)\Ntf(\nt). 
\end{align}

From a trace estimate, using the improved estimates from Section~\ref{sec:nab3n}, the second term and its tangent derivative are controlled in $H^{1/2}$.\\

Using relation~\eqref{eq:relntnut} on $\TT$, we have 
\begin{align*}
  Z(\nt) & = Z\le(\frac{\cc}{\cc+1}\le(\nu + (\cc^{-1}-\half\yy)\nut^{-1}\ri)\ri),
\end{align*}
and $\Nd^{\leq 1}Z(\nt)$ is schematically composed of 
\begin{align*}
  \Nd^{\leq 1} (\Nd_3,\Nd_4)\yy, &&  \Nd^{\leq 1}Z(\nut).
\end{align*}

From the assumption~(\ref{est:bdysourcenullconn}) (see the definitions of Section~\ref{sec:normsnullconn} and the improvements of Section~\ref{sec:connest}), we have
\begin{align}\label{est:Nd34yyTT}
  \norm{t(t\Nd)^{\leq 1}(t\Nd_3)\yy}_{\HHt(\pr\Si_t)} + \norm{t(t\Nd)^{\leq 1}(t\Nd_4)\yy}_{\HHt(\pr\Si_t)} & \les \OOfexti[\yy] \les \varep.
\end{align}

Using relation~\eqref{eq:transkTTt} for $\tr\kapt$, we have on $\TT$
\begin{align*}
  \nut^{-1}-\nut & = -r\delt +r \le(\trchib+\frac{2}{r}\ri) + r\le(\trchi-\frac{2}{r}\ri) + \Err,
\end{align*}
which derived in $Z$ gives
\begin{align*}
  Z(\nut) & = \half Z\le(r\delt\ri) - \half Z\le(r \le(\trchib+\frac{2}{r}\ri)\ri) - \half Z\le(r\le(\trchi-\frac{2}{r}\ri)\ri) + \Err.
\end{align*}
Using the improved estimates~(\ref{est:nab2ktL2}), (\ref{est:LiehTfktL2}) for respectively $\nab^{\leq 2}\kt$ and $\nab^{\leq 1}\Lieh_\Tf\kt$ and trace estimates, one deduces from the above formula
\begin{align}\label{est:ZnutH121}
  \norm{t^2Z(\nut)}_{\HHt(\pr\Si_t)} & \les \varep.
\end{align}
Using relations~\eqref{eq:transkTTt} for $\ep$, and the improved estimates~(\ref{est:nab2ktL2}),~(\ref{est:LiehTfktL2}) for $\kt$ and trace estimates, we have
\begin{align}\label{est:ZnutH122}
  \norm{t^2(t\Nd)Z(\nut)}_{\HHt(\pr\Si_t)} & \les \varep.
\end{align}

Thus, combining~(\ref{est:Nd34yyTT}), (\ref{est:ZnutH121}), (\ref{est:ZnutH122}) and equation~\eqref{eq:TnZnNn}, we obtain
\begin{align*}
  \norm{t(t\Nd)^{\leq 1}\Tf(\nt)}_{\HHt(\pr\Si_t)} & \les \varep.
\end{align*}
Applying the elliptic estimate of Lemma~\ref{lem:ellLap}, commuting Laplace equation~\eqref{eq:Lapnt} with $\Tf$ and using the above boundary estimate, we obtain
\begin{align}\label{est:nabTfntL2}
  \norm{t(t\nab)^{\leq 2}\Tf(\nt)}_{L^2(\Si_t)} & \les \varep,
\end{align}
as desired.
\subsection{Control of $\Lieh^2_\Tf\kt$ on $\Si_t$}
The improved estimates for $\Lieh^2_\Tf\kt$ are obtained directly using equation~\eqref{eq:LiehTfkt}, the estimates just obtained in Section~\ref{sec:TntSit} for $\nab^{\leq 2}\Tf(\nt)$ together with the assumed estimates~(\ref{est:curvestsourceplanehyp}) for the curvature, and we have
\begin{align}\label{est:T2ktL2Sit}
  \begin{aligned}
  \norm{t^3\Lieh^2_\Tf\kt}_{L^2(\Si_t)} & \les \norm{t^2(t\Lieh_\Tf)^{\leq 1} E}_{L^2(\Si_t)} + \norm{t(t\nab)^{\leq 2} (t \Tf)^{\leq 1}(n)}_{L^2(\Si_t)} + \text{nonlinear error terms} \\
  & \les \varep.
  \end{aligned}
\end{align}

\subsection{Control of $\nab^{\leq 1}\Tf^2(\nt)$ on $\Si_t$}\label{sec:nabTf2ntSit}
Arguing as in Section~\ref{sec:TntSit}, our goal is to obtain boundary estimates for $\Tf^2(\nt)$ on $\pr\Si_t$. As in~\eqref{eq:TnZnNn}, one can obtain schematically obtain a formula
\begin{align}\label{eq:T2nprSit}
  \Tf^2(\nt) & = Z^2(\nt) + \Nf^2(n) + \Nf\Tf(n) + \text{lower order terms}.
\end{align}
The term $Z^2(\nt)$ is schematically composed of
\begin{align*}
  Z^{\leq 2}(\yy) && Z^{\leq 2}(\nu). 
\end{align*}
Using the embeddings of Lemma~\ref{lem:Sobsphere}, we have the following $L^2_t\HHt(\pr\Si_t)$ bounds for $Z^2(\yy)$
\begin{align*}
  \norm{t^{1/2-\ga}t^2 Z^2(\yy)}_{L^2_t \HHt(\pr\Si_t)} & \les \norm{t^{-\ga}t^2 (t\Nd)^{\leq 1}Z^2(\yy)}_{L^2_tL^2(\pr\Si_t)} \les \OO^\TT_{\leq 2, \ga}[\yy] \les \varep,
\end{align*}
for all $\ga>0$.\\

As in Section~\ref{sec:TntSit}, $\HHt$ estimates for $\nu$ can be obtained from trace estimates and the improved estimates~(\ref{est:nab2ktL2}), (\ref{est:LiehTfktL2}), (\ref{est:T2ktL2Sit}) for $\kt$, and we have
\begin{align*}
  \norm{t^3Z^{2}(\nut)}_{L^\infty_t \HHt(\pr\Si_t)} & \les \varep.
\end{align*}
Using the above estimate, we deduce the following $L^2_t\HHt(\pr\Si_t)$ control
\begin{align*}
  \norm{t^{1/2-\ga}t^2Z^2(\nut)}_{L^2_t \HHt(\pr\Si_t)} & \les \norm{t^3Z^{2}(\nut)}_{L^\infty_t \HHt(\pr\Si_t)} \les \varep,
\end{align*}
for all $\ga>0$.\\

For the lower order terms and the terms $\Nf^2(n)$ and $\Nf\Tf(n)$, one can obtain $L^\infty_t\HHt(\pr\Si_t)$ bounds, using trace estimates and the previous bounds~(\ref{est:nab3L2n-1}), (\ref{est:nabTfntL2}), and these bounds can similarly be turned into weaker $L^2_t\HHt(\pr\Si_t)$ bounds. Combining these estimates, we obtain the following boundary control
\begin{align}\label{est:T2nL2tH12}
  \norm{t^{1/2-\ga}t^2\Tf^2(n)}_{L^2_t\HHt(\pr\Si_t)} & \les \varep.
\end{align}

Applying the elliptic estimates of Lemma~\ref{lem:ellLap}, using that by~(\ref{eq:Lapnt}) $\Delta(n)$ is composed of nonlinear lower order error terms, we obtain on each separate slice $\Si_t$
\begin{align*}
  \norm{t^{-1/2-\ga}t^2(t\nab)^{\leq 1}\Tf^2(n)}_{L^2(\Si_t)} & \les \norm{t^{1/2-\ga}t^2\Tf^2(n)}_{\HHt(\pr\Si_t)} + \text{nonlinear error terms}.
\end{align*}
Taking the $L^2_t$ norm in the above estimate and using the boundary bound~\eqref{est:T2nL2tH12}, we obtain
\begin{align*}
  \norm{t^{-1/2-\ga}t^2(t\nab)^{\leq 1}\Tf^2(n)}_{L^2_t L^2(\Si_t)} & \les \norm{t^{1/2-\ga}t^2\Tf^2(n)}_{L^2_t\HHt(\pr\Si_t)} + \text{nonlinear error terms} \les \varep.
\end{align*}
This finishes the proof of~(\ref{est:OOmaxconn}).\\

\subsection{Control of $\D\TI$}
From~\eqref{est:OOmaxconn}, the formula~(\ref{eq:piTI}) and the Sobolev embeddings of Lemma~\ref{lem:SobSitHHrr}, we directly deduce the following control for $\D\TI$ and in particular for the deformation tensor $^{(\TI)}\pih$ of the interior approximate Killing field $\TI$
\begin{align}\label{est:controlTI}
  \begin{aligned}
    \norm{t^{5/2}{\D\TI}}_{L^\infty(\MM^\intr_\bott)} & \les \varep,\\
    \norm{t^2(t\D)^{\leq 1}{\D\TI}}_{L^\infty_t L^6(\Si_t)} + \norm{t^2(t\D)^{\leq 1}{\D\TI}}_{L^\infty_tL^4(\pr\Si_t)} & \les \varep,\\
    \norm{t^{-1/2-\ga}t(t\D)^{\leq 2}{\D\TI}}_{L^2(\MM^\intr_\bott)} & \les \varep,
  \end{aligned}
\end{align}
for all $\ga>0$.

\section{Control of the harmonic Cartesian coordinates on $\Si_\tast$}\label{sec:harmocoordsSitast}
From Gauss equation~(\ref{eq:GaussSit}) on $\Si_\tast$, the curvature bounds~(\ref{est:curvestsourceplanehyp}) and the Bootstrap Assumptions~\ref{BA:connint}, we have
\begin{align}\label{est:RicSitast}
  \norm{t^2(t\Nd)^{\leq 2}\RRRic}_{L^2(\Si_\tast)} & \les \varep.
\end{align}
From relation~(\ref{eq:NTT}) for the unit normal $\Nf$ to $\pr\Si_\tast$, the assumed bounds~(\ref{est:sourcechichibSast}) for $\chi$ and $\chib$ on $S^\ast \subset \CCba$, the improved bounds~(\ref{est:Nd2nutH12}) for $\nut$ on $\pr\Si_\tast$, we have
\begin{align}\label{est:thprSitast}
  \norm{t^2(t\Nd)^{\leq 2}\le(\tr\th-\frac{2}{\rast}\ri)}_{\HHt(S^\ast)} + \norm{t^2(t\Nd)^{\leq 2}\thh}_{\HHt(S^\ast)} & \les \varep.
\end{align}
Thus, applying the (rescaled) results of Theorem~\ref{thm:globharmonics}, we deduce that (for all the centred conformal isomorphisms of $S^\ast$), estimate~(\ref{est:sourceharmo}) holds for the harmonic Cartesian coordinates of $\Si_\tast$.

\section{Control of the interior Killing fields $\TI, \SI, \KI$ and $\OOI$ in $\MM^\intr_\bott$}\label{sec:controlKillingint}
In Section~\ref{sec:mildcontrolTISIKI}, we improve the mild Bootstrap Assumptions~\ref{BA:mildKillingMMintbot} on the interior approximate Killing vectorfields $\TI, \SI, \KI$ and on the vectorfield $\XI$. In Section~\ref{sec:XIXI}, we improve the control~(\ref{est:BAXIXI}) for $\XI$ from the Bootstrap Assumptions~\ref{BA:intKill}. In Section~\ref{sec:controlpihint}, we improve the bounds~(\ref{est:BADX}) and~\eqref{est:BAXITI} on (derivatives of) the deformation tensors, which finishes the improvement of the Bootstrap Assumptions~\ref{BA:intKill}.

\subsection{Control of $\D\XI$, $\D\SI, \D\KI$ and $\D{\OOI}$}\label{sec:controlpihint}
\begin{lemma}\label{lem:DXI}
  The following bounds hold on $\MM^\intr_\bott$ 
  \begin{align}\label{est:DXIg}
    \begin{aligned}
    \norm{t^{3/2}\le(\D\XI - g\ri)}_{L^\infty(\MM^\intr_\bott)} & \les \varep,\\
    \norm{t(t\D)^{\leq 1}\le(\D\XI-g\ri)}_{L^\infty_tL^6(\Si_t)} & \les \varep,\\
    \norm{t^{-1/2-\ga}t(t\D)^{\leq 2}\le(\D\XI-g\ri)}_{L^2(\MM^\intr_\bott)} & \les \varep,
    \end{aligned}
  \end{align}
  for all $\ga>0$. We also have the following bounds on $\pr\Si_t$
  \begin{align}\label{est:DXigTT}
    \norm{t^{-\ga}t(t\D)^{\leq 1}\le(\D\XI-g\ri)}_{L^\infty_tL^4(\pr\Si_t)},
  \end{align}
  for all $\ga>0$. Moreover, we have the following control on $\MM^\intr_\bott$
  \begin{align}
    \label{est:XITIprecise}
    \le|(t\D)^{\leq 1}\g(\XI,\TI)\ri| & \les \varep t^{-1/2}.
  \end{align}
\end{lemma}
\begin{proof}
  From the definition~(\ref{eq:defXISitast}) on $\Si_\tast$, $\XI$ is an $\Si_\tast$ tangent vectorfield, such that
  \begin{align*}
    \nab_k\XI_l & = \sum_{i=1}^3\nab_k x^i \nab_l x^i + \sum_{i=1}^3x^i\nab^2_{k,l}x^i \\
                & = g_{kl} + (\de_{kl}-g_{kl}) + \sum_{i=1}^3 x^i\nab^2_{k,l}x^i. 
  \end{align*}
  and from the bound~(\ref{est:sourceharmo}), we deduce
  \begin{align}\label{est:nabXSitast}
    \begin{aligned}
      \norm{(t\nab)^{\leq 2}\le(\nab \XI -g\ri)}_{L^2(\Si_\tast)} & \les \norm{(t\nab)^{\leq 2}(\de_{kl}-g_{kl})}_{L^2(\Si_\tast)} + \norm{t(t\nab)^{\leq 2}\nab^2x^i}_{L^2(\SI_\tast)}\\
      & \les \varep.
    \end{aligned}
  \end{align}

  From the definition~\eqref{eq:defXI}, we have schematically
  \begin{align*}
    \D \XI & = \nab \XI + \XI \cdot \kt,
  \end{align*}
  thus, we deduce from~\eqref{est:nabXSitast}, the bounds~(\ref{est:OOmaxconn}) for $\kt$ and the mild Bootstrap Assumptions~\ref{BA:mildKillingMMintbot} and~\ref{BA:intKill} for $\XI$ that
  \begin{align}
    \label{est:DXSitast}
    \norm{(t\D)^{\leq 2}\le(\D\XI -g\ri)}_{L^2(\Si_\tast)} & \les \varep.
  \end{align}
  
  Commuting the transport equation~\eqref{eq:defXI} for $\XI$, we have
  \begin{align*}
    \D_\Tf\D_\mu\XI_\nu & = -\D_{\D_\mu\Tf}\XI_\nu + \XI^{\la}\Tf^{\la'}\R_{\la\nu\la'\mu}.
  \end{align*}
  Moreover, we have
  \begin{align*}
    \D_\Tf g & = \D_\Tf(\g + \Tf \otimes \Tf) = \D_\Tf\Tf \otimes \Tf + \Tf\otimes\D_\Tf\Tf,
  \end{align*}
  and therefore $\D\XI-2g$ satisfy the following transport equation in $\MM^\intr_\bott$
  \begin{align}\label{eq:TfDXIg}
    \D_\Tf\le(\D_\mu\XI_\nu-g_{\mu\nu}\ri) & = -\D_{\D_\mu\Tf}\XI_\nu + \XI^{\la}\Tf^{\la'}\R_{\la\nu\la'\mu} - \D_\Tf\Tf_\mu \Tf_\nu + \D_\Tf\Tf_\nu\Tf_\mu.
  \end{align}

  Commuting the transport equation~\eqref{eq:TfDXIg} with $(t\D)^{\leq 2}$, using the estimates~(\ref{est:curvestsourceplanehyp}) for the curvature, the mild Bootstrap Assumptions~\ref{BA:mildKillingMMintbot} for $\TI,\XI$, the bootstrap bounds~\eqref{est:BADX} for $\D\XI$, and the bounds~(\ref{est:controlTI}) for $\D\TI$, we deduce from~\eqref{eq:TfDXIg} that
  \begin{align}\label{est:DTfDXIgLinf}
    \norm{t^{-1/2-\ga}t\D_\Tf(t\D)^{\leq 2}\le(\D\XI-g\ri)}_{L^2(\MM^\intr_\bott)} & \les \varep.
  \end{align}
  Integrating~\eqref{est:DTfDXIgLinf}, using the bound~\eqref{est:DXSitast} on $\Si_\tast$, we obtain
  \begin{align*}
    \norm{t^{-\ga}t(t\D)^{\leq 2}\le(\D\XI-g\ri)}_{L^\infty_t L^2(\Si_t)} & \les \varep,
  \end{align*}
  for all $\ga>0$, from which we also deduce the desired bounds
  \begin{align*}
    \norm{t^{-1/2-\ga}t(t\D)^{\leq 2}\le(\D\XI-g\ri)}_{L^2(\MM^\intr_\bott)} & \les \varep,
  \end{align*}
  and, using a trace estimate on $\pr\Si_t$,
  \begin{align*}
    \norm{t^{-\ga}t(t\D)^{\leq 1}\le(\D\XI-g\ri)}_{L^\infty_tL^4(\pr\Si_t)},
  \end{align*}
  for all $\ga>0$. The $L^\infty(\MM^\intr_\bott)$ and $L^\infty_tL^6(\Si_t)$ estimates of~\eqref{est:DXIg} follow similarly and are left to the reader. \\

  To obtain~(\ref{est:XITIprecise}), we use the definition~(\ref{eq:defXI}) of $\XI$, from which we have
  \begin{align}\label{eq:TfgTIXI}
    \begin{aligned}
      \Tf\le(\g(\TI,\XI)\ri) & = ^{(\TI)}\pi\le(\Tf,\XI\ri), \\
      \g(\TI,\XI)|_{\Si_\tast} & = 0.
    \end{aligned}
  \end{align}
  Thus, integrating~\eqref{eq:TfgTIXI}, using the mild Bootstrap Assumptions~\ref{BA:mildKillingMMintbot} on $\XI$ and the bounds~(\ref{est:controlTI}) on $\D\TI$, we have
  \begin{align}\label{est:gTIXImax}
    |\g(\TI,\XI)| & \les \varep t^{-1/2}.
  \end{align}
  Using the bounds~\eqref{est:controlTI} on $\D\TI$, the bounds~\eqref{est:DXIg} on $\D\XI$ and the mild Bootstrap Assumptions~\ref{BA:mildKillingMMintbot} on $\TI,\XI$, we further have
  \begin{align*}
    \le|t\D\g(\TI,\XI)\ri| & \les \le|t\g(\D\TI,\XI)\ri| + \le|t\g(\TI,\D\XI)\ri| \\
                           & \les \varep t^{-1/2},
  \end{align*}
  as desired. This finishes the proof of the lemma.
\end{proof}

\begin{lemma}
  The following bounds hold on $\MM^\intr_\bott$
  \begin{align}\label{est:DSI}
    \begin{aligned}
    \norm{t^{3/2}\le(\D\SI - \g\ri)}_{L^\infty(\MM^\intr_\bott)} & \les \varep,\\
    \norm{t(t\D)^{\leq 1}\le(\D\SI-\g\ri)}_{L^\infty_tL^6(\Si_t)}  & \les \varep,\\
    \norm{t^{-1/2-\ga}t(t\D)^{\leq 2}\le(\D\SI-\g\ri)}_{L^2(\MM^\intr_\bott)} & \les \varep,\\
    \norm{t^{-\ga}t(t\D)^{\leq 1}{\le(\D\SI-\g\ri)}}_{L^\infty_tL^4(\pr\Si_t)} & \les \varep,
    \end{aligned}
  \end{align}
  and
    \begin{align}\label{est:DpiKI}
    \begin{aligned}
      \norm{t^{1/2}\le(^{(\KI)}\pi - 4t\g\ri)}_{L^\infty(\MM^\intr_\bott)} & \les \varep,\\
      \norm{(t\D)^{\leq 1}\le(^{(\KI)}\pi-4t\g\ri)}_{L^\infty_tL^6(\Si_t)}  & \les \varep,\\
      \norm{t^{-1/2-\ga}(t\D)^{\leq 2}\le(^{(\KI)}\pi-4t\g\ri)}_{L^2(\MM^\intr_\bott)} & \les \varep.
    \end{aligned}
  \end{align}
\end{lemma}
\begin{proof}
  The estimates~\eqref{est:DSI} are a direct consequence of the definition~\eqref{eq:defSI} of $\SI$ and the bounds~(\ref{est:OOmaxconn}) for the maximal connection coefficients, the bounds~(\ref{est:controlTI}) for $\D\TI$ and the bounds~(\ref{est:DXIg}) for $\XI$.\\

  From the definition~\eqref{eq:defKI} of $\KI$, we have
  \begin{align*}
    \D\KI & = 2t\D t \otimes\TI + 2\g(\D\XI,\XI)\otimes\TI +\le(t^2+\g(\XI,\XI)\ri)\D\TI + 2\D t \otimes\XI +2t \D\XI\\
          & = -2n^{-1}t\TI\otimes\TI +2\XI\otimes\TI +2\g(\D\XI-g,\XI)\otimes\TI + \le(t^2+\g(\XI,\XI)\ri)\D\TI  \\
    & \quad - 2 n^{-1}\TI \otimes\XI + 2t g + 2t(\D\XI-g) \\
          & = 2t\g +2\XI\otimes\TI - 2\TI\otimes\XI + \EEE,
  \end{align*}
  where
  \begin{align*}
    \EEE & := - 2t(n^{-1}-1)\TI\otimes\TI +2\g(\D\XI-g,\XI)\TI + \le(t^2+\g(\XI,\XI)\ri)\D\TI \\
    & \quad - 2 (n^{-1}-1)\TI \otimes\XI + 2t(\D\XI-g).
  \end{align*}

  From the bounds~(\ref{est:OOmaxconn}) for the maximal connection coefficients, the bounds~(\ref{est:controlTI}) for $\D\TI$ and the bounds~(\ref{est:DXIg}) for $\XI$, and the mild Bootstrap Assumptions~\ref{BA:mildKillingMMintbot}, we have
  \begin{align*}
    \begin{aligned}
      \norm{t^{1/2}\EEE}_{L^\infty(\MM^\intr_\bott)} & \les \varep,\\
      \norm{(t\D)^{\leq 1}\EEE}_{L^\infty_tL^6(\Si_t)} & \les \varep,\\
      \norm{t^{-1/2-\ga}(t\D)^{\leq 2}\EEE}_{L^2(\MM^\intr_\bott)} & \les \varep,
    \end{aligned}
  \end{align*}
  and the desired bounds~\eqref{est:DpiKI} follow. This finishes the proof of the lemma.  
\end{proof}

\begin{lemma}\label{lem:D2OOIlem}
  The following bounds hold on $\MM^\intr_\bott$
  \begin{align}
    \norm{\D\OOI}_{L^\infty(\MM^\intr_\bott)} & \les 1,\label{est:DOOILinf}\\
    \norm{t^{3/2}{^{(\OOI)}\pi}}_{L^\infty(\MM^\intr_\bott)} & \les \varep,\label{est:piOOILinf}\\
    \norm{t^2\D^2\OOI}_{L^\infty_tL^6(\Si_t)} & \les \varep, \label{est:D2OOIL6L4}\\
    \norm{t^{-1/2-\ga}t^2(t\D)^{\leq 1}\D^2\OOI}_{L^2(\MM^\intr_\bott)} & \les \varep,\label{est:D3OOIL2MMintr}\\
    \norm{t^{-\ga}t^2\D^2\OOI}_{L^\infty_tL^4(\pr\Si_t)} & \les \varep\label{est:D2OOITTL4} 
  \end{align}

  for all $\ga>0$.
\end{lemma}
\begin{proof}
  Commuting definition~\eqref{eq:defOOI} of $\OOI$ in $\MM^\intr_\bott$ by $\D$ and arguing as in the proof of Lemma~\ref{lem:DXI}, we have schematically
  \begin{align}\label{eq:DTfpiOOI}
    \D_\Tf \D\OOI & = \D\Tf\cdot\D\OOI + \R\cdot\OOI\cdot\Tf. 
  \end{align}
  Using the definition~\eqref{eq:defOOISitast} of $\OOI$ on $\Si_\tast$, $\OOI$ is a $\Si_\tast$-tangent vectorfield and we have
  \begin{align}\label{eq:piOOISitast1}
    \begin{aligned}
      \nab_l{^{(\ell)}\OOI}_m & = \in_{\ell i j}\nab_l x^i\nab_m x^j + \in_{\ell i j}x^i \nab^2_{l,m}x^j,\\ 
      \nab_l{^{(\ell)}\OOI}_m + \nab_m{^{(\ell)}\OOI}_l & = 2\in_{\ell i j} x^i \nab^2_{l,m}x^j.
    \end{aligned}
  \end{align}
  Moreover, using~\eqref{eq:defOOI}, we have on $\Si_\tast$
  \begin{align}\label{eq:piOOISitast2}
    \begin{aligned}
      \D_\Tf{^{(\ell)}\OOI}_i & = 0,\\
      \D_i{^{(\ell)}\OOI}_{\Tf} & = - \g(^{(\ell)}\OOI, \D_i\Tf) = k_{ij}{^{(\ell)}\OOI}^j,\\
      \D_\Tf{^{(\ell)}\OOI}_{\Tf} & = 0.
    \end{aligned}
  \end{align}
  Thus, we deduce from~\eqref{eq:piOOISitast1} and~\eqref{eq:piOOISitast2} and the bounds~\eqref{est:sourceharmo} and~(\ref{est:controlTI}) that
  \begin{align*}
    \norm{\D\OOI}_{L^\infty(\Si_\tast)} & \les 1,\\
    \norm{t^{3/2}{^{(\OOI)}\pi}}_{L^\infty(\Si_\tast)} & \les \varep.
  \end{align*}
  Integrating~\eqref{eq:DTfpiOOI}, we obtain the desired bounds~\eqref{est:DOOILinf} and~\eqref{est:piOOILinf}.\\

  Differentiating formula~\eqref{eq:piOOISitast1}, we obtain schematically on $\Si_\tast$
  \begin{align*}
    \nab^2{\OOI} & = \nab x\nab^2x + x\nab^3x,
  \end{align*}
  from which, arguing as previously, using the bounds~\eqref{est:sourceharmo} and~(\ref{est:controlTI}), one deduces that 
  \begin{align*}
    \norm{t^2\D^2\OOI}_{L^6(\Si_\tast)} & \les \varep,\\
    \norm{t(t\D)^{\leq 1}\D^2\OOI}_{L^2(\Si_\tast)} & \les \varep.
  \end{align*}


  Arguing as in the proof of Lemma~\ref{lem:DXI}, integrating the transport equation~\eqref{eq:DTfpiOOI} and applying a trace estimate, estimates (\ref{est:D2OOIL6L4}), (\ref{est:D3OOIL2MMintr}) and (\ref{est:D2OOITTL4}) follow. This finishes the proof of the lemma. 
\end{proof}

  

\subsection{Control of $\XI$}\label{sec:XIXI}
This section is dedicated to the improvement of the bound~\eqref{est:BAXIXI} from the Bootstrap Assumptions~\ref{BA:intKill}. 
\begin{proposition}\label{prop:XIXI}
  We have for all $\too \leq t\leq \tast$
  \begin{align*}
    \le|\sup_{\Si_t} \le(t^{-2}\g(\XI,\XI)\ri) - \le(\frac{1-\cc}{1+\cc}\ri)^2\ri|  & \les \varep t^{-3/2}.
  \end{align*}
\end{proposition}
Before turning to the proof of Proposition~\ref{prop:XIXI}, we have the following three lemmas.
\begin{lemma}\label{lem:DXIXI}
  The following bound holds on $\MM^\intr_\bott$
  \begin{align*}
    \le|\D_\XI\XI - \XI\ri| & \les \varep t^{-1/2},
  \end{align*}
  where the norm is taken with respect to the maximal frame.
\end{lemma}
\begin{proof}
  The proof is a straight-forward adaptation of the estimates of Lemma~\ref{lem:DXI}.
\end{proof}

\begin{lemma}\label{lem:supSitprSit}
  For all $\too \leq t\leq\tast$, we have
  \begin{align}\label{eq:supsupXIXI}
    \sup_{\Si_t} t^{-2}\g(\XI,\XI) & = \sup_{\pr\Si_t}t^{-2}\g(\XI,\XI).
  \end{align}
\end{lemma}
\begin{proof}
  Let $\too \leq t \leq \tast$ and assume that the supremum on $\Si_t$ is reached for $p\in \Si_t \setminus \pr\Si_t$. Since $p$ is not on the boundary, we have in particular, for the derivative in the (projected on $\Si_t$) $\XI$-direction
  \begin{align*}
    0 & = \nab_{\XI + \g(\XI,\TI)\TI} \le(t^{-2}\g(\XI,\XI)\ri). 
  \end{align*}
  Using~\eqref{eq:defXI} and relations~(\ref{eq:relTnt}), this gives
  \begin{align*}
    0 & = t^{-2}\g(\XI,\XI) + t^{-2}\g(\D_\XI\XI-\XI,\XI). 
  \end{align*}
  Using the bound from Lemma~\ref{lem:DXIXI}, we thus deduce
  \begin{align*}
    \le|t^{-2}\g(\XI,\XI)\ri| & \les \varep t^{-3/2} \le(\le|t^{-2}\g(\XI,\XI)\ri|\ri)^{1/2}, 
  \end{align*}
  which we rewrite as
  \begin{align*}
    \le|t^{-2}\g(\XI,\XI)\ri| & \les \varep^2 t^{-3}.
  \end{align*}
  For $\varep>0$ sufficiently small, this contradicts the bootstrap bound~(\ref{est:BAXIXI}), and~\eqref{eq:supsupXIXI} follows.
\end{proof}
\begin{remark}
  From the result of Lemma~\ref{lem:supSitprSit} and the bootstrap bound~(\ref{est:BAXIXI}), we also deduce the following mild control for $\g(\XI,\XI)$ on $\TT$
  \begin{align}\label{est:mildXIXITT}
    t^2 & \les \g(\XI,\XI). 
  \end{align} 
\end{remark}

\begin{lemma}\label{lem:ZZterr}
  We have the following control on $\TT$
  \begin{align}\label{est:ZZterr}
    \le|\ZZt - \le(\Tf +\frac{1-\cc}{1+\cc}\Xf\ri)\ri| & \les \varep t^{-3/2},
  \end{align}
  where $\Xf := \g(\XI,\XI)^{-1/2}\XI$ and $\ZZt$ is the $\TT$-tangent vectorfield normal to $\pr\Si_t$ in $\TT$ and such that $\ZZt(t) = 1$.
\end{lemma}
\begin{proof}
  Let define
  \begin{align*}
    \ZZterr & := \ZZt - \le(\Tf + \frac{1-\cc}{1+\cc}\Xf\ri).
  \end{align*}
  From its definition, we have the following expression for $\ZZt$
  \begin{align}\label{eq:ZZt}
    \ZZt & = \le(\frac{\cc}{1+\cc}\ri)\elb + \le(\frac{2-\cc\yy}{2+2\cc}\ri)\el.
  \end{align}
 Using the Bootstrap Assumptions~\ref{BA:connext} and~\ref{BA:connint} on the exterior and interior connection coefficients, one can therefore deduce the following mild control
  \begin{align}\label{est:mildZZt}
    |\ZZt| & \les 1,
  \end{align}
  and using additionally relations~\eqref{eq:Riccirel},
  \begin{align}\label{est:DZZtZZt}
    \le|\D_\ZZt\ZZt\ri| & \les \varep t^{-5/2},
  \end{align}
  where norms are taken in the maximal frame.\\
  
  Using~\eqref{est:controlTI} and~\eqref{est:mildZZt}, we have
  \begin{align}\label{est:DZZtTf}
    \le|\D_\ZZt\Tf\ri| & \les \varep t^{-5/2}.
  \end{align}
  
  We have
  \begin{align}\label{eq:DZZtXf}
    \D_{\ZZt}\Xf & = \D_{\Tf + \frac{1-\cc}{1+\cc}\Xf}\Xf + \ZZterr\cdot\D\Xf.
  \end{align}

  We have
  \begin{align*}
    \D_\Xf\Xf & = \g(\XI,\XI)^{-1/2}\D_{\XI}\le(\g(\XI,\XI)^{-1/2}\XI\ri) \\
              & = \g(\XI,\XI)^{-1/2}\Xf + \g(\XI,\XI)^{-1}\le(\D_\XI\XI-\XI\ri) + \le(\D_{\XI}\le(\g(\XI,\XI)^{-1/2}\ri)\ri)\Xf \\
              & = \g(\XI,\XI)^{-1/2}\Xf + \g(\XI,\XI)^{-1}\le(\D_\XI\XI-\XI\ri) -\le(\g(\XI,\XI)^{-3/2} \g(\D_{\XI}\XI,\XI)\ri)\Xf\\
              & = \g(\XI,\XI)^{-1}\le(\D_\XI\XI-\XI\ri) - \g(\XI,\XI)^{-3/2} \g(\D_\XI\XI-\XI,\XI)\Xf.
  \end{align*}
  
  Using the result of Lemma~\ref{lem:DXIXI} and the mild control $t \les |\XI| \les t$ from the Bootstrap Assumptions~\ref{BA:mildKillingMMintbot} and~\eqref{est:mildXIXITT} on $\TT$, we therefore deduce
  \begin{align*}
    \le|\D_\Xf\Xf\ri| & \les \varep t^{-5/2},
  \end{align*}
  which using that $\D_\Tf\XI = 0$ and~\eqref{eq:DZZtXf} further gives on $\TT$
  \begin{align}\label{est:DZZtXf}
    \le|\D_\ZZt\Xf\ri| & \les \varep t^{-5/2} + \le|\ZZterr\ri|\le|\D\Xf\ri|.
  \end{align}
  We deduce from the definition of $\ZZterr$ and estimates~\eqref{est:DZZtZZt},~\eqref{est:DZZtTf} and~\eqref{est:DZZtXf} the following estimate on $\TT$
  \begin{align}
    \label{est:DZZtZZterr}
    \le|\D_{\ZZt}\ZZterr\ri| & \les \varep t^{-5/2} + \le|\ZZterr\ri|\le|\D\Xf\ri|.
  \end{align}

  Using the expression~\eqref{eq:ZZt}, that $\yy=0$ on $\CCba$ and the relations~(\ref{eq:defnutnur}), (\ref{eq:NTT}), we have on $S^\ast$
  \begin{align}\label{eq:ZZterrSast}
    \begin{aligned}
      \ZZterr & = \frac{\cc}{1+\cc}\elb+\frac{1}{1+\cc}\el - \Tf - \frac{1-\cc}{1+\cc}\Xf \\
      & = \frac{\cc}{1+\cc}\le(1-\nu^{-1}\ri)\elb + \frac{1}{1+\cc}\le(1-\nu\ri)\el -\frac{1-\cc}{1+\cc}\le(\Xf-\Nf\ri).
    \end{aligned}
  \end{align}
  From the bound~\eqref{est:sourceharmo} on the harmonic coordinates on $\Si_\tast$, we have
  \begin{align*}
    \le|\nab x^i - \pr_{x^i}\ri| & \les \varep \tast^{-3/2}.
  \end{align*}
  Thus, from the definition~\eqref{eq:defXISitast} of $\XI$ on $\Si_\tast$, we have
  \begin{align*}
    \le|\XI-x^i\pr_{x^i}\ri| & \les \varep \tast^{-1/2}.
  \end{align*}
  From the bounds~\eqref{est:sourceharmo} on $S^\ast$, we have\footnote{We recall that $\rast$ is the area radius of $S^\ast$.}
  \begin{align*}
    \le|\rast\Nf(x^i)-x^i\ri| & \les \varep \tast^{-1/2}
  \end{align*}
  from which we deduce on $S^\ast$
  \begin{align*}
    \le|\Xf-\Nf\ri| & \les \tast^{-1} \le|\XI-\rast\Nf\ri| \les \tast^{-1} \le|x^i\pr_{x^i}-\rast\Nf\ri| + \varep \tast^{-3/2} \les \varep \tast^{-3/2}.
  \end{align*}
  Using this bound together with the estimates for $\nu$ from~(\ref{est:OOmaxconn}) in equation~\eqref{eq:ZZterrSast}, we deduce that on $S^\ast$
  \begin{align}\label{est:ZZterrSitast}
    \le|\ZZterr\ri| & \les \varep t^{-3/2}.
  \end{align}
  
  Integrating the estimate~\eqref{est:DZZtZZterr} for $\ZZterr$ along $\ZZt$, using the mild control $|\D\Xf| \les 1$, the estimate~\eqref{est:ZZterrSitast} for $\ZZterr$ on $S^\ast$ and a Gr\"onwall argument, the result of the lemma follows.
\end{proof}

\begin{proof}[Proof of Proposition~\ref{prop:XIXI}]
  We have
  \begin{align}\label{eq:ZZtXIXI}
    \ZZt(\g(\XI,\XI)) = \frac{1-\cc}{1+\cc}\Xf(\g(\XI,\XI)) + \ZZterr(\g(\XI,\XI)).
  \end{align}

  From Lemma~\ref{lem:DXIXI} and the mild control $|\XI| \les t$, we have
  \begin{align*}
    \le|\Xf(\g(\XI,\XI)) -2\sqrt{\g(\XI,\XI)}\ri| & \les \varep t^{-1/2}.
  \end{align*}

  Plugging this in~\eqref{eq:ZZtXIXI} and using estimate~\eqref{est:ZZterr} and the mild control $t\les \le|\XI\ri| \les t$ and $\le|\D\XI\ri| \les 1$, we have
  \begin{align}\label{est:ZZtXIXI}
    \le|\frac{\ZZt(\g(\XI,\XI))}{2\sqrt{\g(\XI,\XI)}}- \frac{1-\cc}{1+\cc}\ri| & \les \varep t^{-3/2}. 
  \end{align}

  Integrating~\eqref{est:ZZtXIXI} from $\tast$ to $t$, we obtain
  \begin{align}\label{est:intZZtXIXI}
    \le|\le[\sqrt{\g(\XI,\XI)}\ri]^{\tast}_{t} - \frac{1-\cc}{1+\cc}(\tast-t)\ri| & \les \varep t^{-1/2}.
  \end{align}

  From the definition of $\XI$ on $\Si_\tast$, and the bounds~\eqref{est:sourceharmo} on the harmonic coordinates, we have on $S^\ast$
  \begin{align*}
    \le|\g(\XI,\XI) -(\rast)^2 \ri| & = \le|\sum_{i,j=1}^3 x^ix^j g^{ij} -(\rast)^2\ri| = \le|\sum_{i,j=1}^3 x^ix^j \le(g^{ij}-\de^{ij}\ri)\ri| \les \varep \tast^{1/2},
  \end{align*}
  from which we deduce
  \begin{align*}
    \le|\sqrt{\g(\XI,\XI)}|_{S^\ast}-\rast\ri| & \les \varep \tast^{-1/2}.
  \end{align*}
  From the definition of $\tast$ and the area radius estimate~(\ref{est:arearadiusestimateCCba}), we have
  \begin{align*}
    \le|\rast - \frac{1-\cc}{1+\cc}\tast\ri| & = \le|\rast - \half (1-\cc)\uba \ri| \les \varep \tast^{-2}.
  \end{align*}
  Thus,
  \begin{align*}
    \le|\sqrt{\g(\XI,\XI)}|_{S^\ast}-\frac{1-\cc}{1+\cc}\tast \ri| & \les \varep \tast^{-1/2}.
  \end{align*}
  Plugging this in~\eqref{est:intZZtXIXI}, we infer
  \begin{align*}
    \le|\sqrt{\g(\XI,\XI)}|_{\pr\Si_t} - \frac{1-\cc}{1+\cc}t\ri| & \les \varep t^{-1/2}.
  \end{align*}
  Using the result of Lemma~\ref{lem:supSitprSit}, the result of Proposition~\ref{prop:XIXI} follows.  
\end{proof}

\subsection{Mild control of $\XI, \SI, \KI$ and $\OOI$}\label{sec:mildcontrolTISIKI}
The following lemmas improve the mild Bootstrap Assumptions~\ref{BA:mildKillingMMintbot}. 
\begin{lemma}
  We have the following mild (improved) control in $\MM^\intr_\bott$
  \begin{align}\label{est:XISIKI}
    \le|\XI\ri| & \les t, & \le|\SI\ri| & \les t, & \le|\KI\ri| & \les t^2,
  \end{align}
  where the norm is taken with respect to the maximal frame.
\end{lemma}
\begin{proof}
  The norm of $\XI$ in the maximal frame can be expressed as
  \begin{align}\label{eq:normXImax}
    |\XI|^2 & = 2|\g(\TI,\XI)|^2 + \g(\XI,\XI). 
  \end{align}

  Using the bootstrap bounds~\eqref{est:BAXIXI} and~\eqref{est:BAXITI} from the Bootstrap Assumptions~\ref{BA:intKill}, we have
  \begin{align}\label{est:gXIXImax}
    \begin{aligned}
      \g(\XI,\TI) & \les D\varep t^{-1/2}, \\
      \g(\XI,\XI) & \les t^2 + (D\varep) t^{1/2} \les t^2,
      \end{aligned}
  \end{align}
  for $\varep>0$ sufficiently small, and the first bound of~\eqref{est:XISIKI} is improved,\emph{i.e.}
  \begin{align}\label{est:XImax}
    |\XI| & \les t.
  \end{align}


  Using the estimate~\eqref{est:XImax} obtained for $\XI$, and the definitions~\eqref{eq:defSI} and~\eqref{eq:defKI} for $\SI$ and $\KI$, we have
  \begin{align*}
    t^{-1}|\SI| & \les |\TI| + t^{-1}|\XI| \les 1,
  \end{align*}
  and
  \begin{align*}
    t^{-2}|\KI| & \les t^{-2}|t^2+\g(\XI,\XI)| + t^{-1}|\XI| \les 1,
  \end{align*}
  as desired.
\end{proof}

\begin{lemma}
  The following (improved) mild control holds on $\MM^\intr_\bott$
  \begin{align*}
    \g(\KI,\TI) & \leq -\half t^2, & \le|\KI + \g(\KI,\TI)\TI\ri| & \leq \half |\g(\KI,\TI)|,
  \end{align*}
  for $1-\cc>0$ and $\varep>0$ sufficiently small.
\end{lemma}
\begin{proof}
  From the definition~\eqref{eq:defKI} of $\KI$ and the estimate~\eqref{est:gTIXImax}, we have 
  \begin{align*}
    \le|\g(\KI,\TI) + t^2 + \g(\XI,\XI)\ri| & \les D\varep t^{1/2}.
  \end{align*}
  
  Using that $\g(\XI,\XI) \geq 0$ (this is a consequence of the definition~\eqref{eq:defXI}), we have
  \begin{align*}
    \g(\KI,\TI) & \les - t^2 - \g(\XI,\XI) + (D\varep)t^{1/2} \\
                & \les -t^2 + (D\varep)t^{1/2} \\
                & \leq -\half t^2,
  \end{align*}
  provided that $\varep>0$ is sufficiently small.\\ 
  
  Moreover, from the definition of $\KI$ and the bound~\eqref{est:BAXIXI}, we have
  \begin{align*}
    \le|\KI+\g(\KI,\TI)\TI\ri| & \leq 2t|\XI| + (D\varep) t^{1/2} \\
                               & \leq 2\le(\frac{1-\cc}{1+\cc}\ri)^2t^2 + (D\varep) t^{1/2} \\
                               & \leq \frac{1}{4} t^2,
  \end{align*}
  provided that $1-\cc>0$ and $\varep>0$ are sufficiently small. This finishes the proof of the lemma.
\end{proof}

\begin{lemma}
  The following (improved) mild control holds on $\MM^\intr_\bott$
  \begin{align}\label{est:OOImax}
    |\OOI| & \les t.
  \end{align}
\end{lemma}
\begin{proof}
  Using the definition~\eqref{eq:defOOI}, we have
  \begin{align}\label{eq:TgTIOOI}
    \begin{aligned}
      \Tf\le(\g(\TI,\OOI)\ri) & = ^{(\TI)}\pih\le(\Tf,\OOI\ri),\\
      \g(\Tf,\OOI)|_{\Si_\tast} & = 0.
    \end{aligned}
  \end{align}
  Integrating~\eqref{eq:TgTIOOI}, using the control~\eqref{est:controlTI} and the mild Bootstrap Assumptions~\ref{BA:mildKillingMMintbot}, we have
  \begin{align}\label{est:TgTIOOI}
    |\g(\TI,\OOI)| & \les (D\varep) t^{-1/2}.
  \end{align}

  Using the definitions~\eqref{eq:defOOI},~\eqref{eq:defOOISitast}, and estimates~\eqref{est:sourceharmo} for the global harmonic coordinates, we have
  \begin{align}\label{est:gOOIOOI}
    \sup_{\Si_t}\g(\OOI,\OOI) \leq \sup_{\Si_\tast}\g(\OOI,\OOI) \les t^2.  
  \end{align}

  Thus, combining~\eqref{est:TgTIOOI} and~\eqref{est:gOOIOOI}, we have
  \begin{align*}
    |\OOI|^2 & = 2|\g(\TI,\OOI)|^2 + \g(\OOI,\OOI) \les t^2, 
  \end{align*}
  which improves~\eqref{est:OOImax} as desired.
\end{proof}

\section{Control of the Killing fields at the interface $\TT$}\label{sec:killingTT}
This section is dedicated to the improvement of the Bootstrap Assumptions~\ref{BA:TTKilling}.
\begin{lemma}
  The following bounds hold on $\TT$
  \begin{align}
    \le|\TE-\TI\ri| & \les \varep t^{-3/2}, \label{est:TETI}\\
    \le|\SE-\SI\ri| & \les \varep t^{-1/2},\label{est:SESI}\\
    \le|\KE-\KI\ri| & \les \varep t^{1/2}.\label{est:KEKI}
  \end{align}
\end{lemma}
\begin{proof}
  From the transition relations~\eqref{eq:defnutnur}, we have
  \begin{align*}
    \TE-\TI & = \half\le(1-\nut^{-1}\ri)\elb + \half\le(1-\nut\ri)\el,
  \end{align*}
  and~\eqref{est:TETI} follows from the bounds~(\ref{est:OOmaxconn}) for $\nut$.\\

  From the definition~(\ref{eq:defTESEKE}) for $\SE$, and from the definitions of Section~\ref{sec:defcannull}, we have 
  \begin{align*}
    \begin{aligned}
      \D\SE & = \half \le(\D u\otimes\elb + u\D\elb + \D\ub \otimes \el + \ub\D\el\ri) \\
      & = -\half \el\otimes\elb - \half \elb\otimes\el -\quar \yy\el\otimes\el+ \half u\D\elb + \half\ub\D\el \\
      & = -\half \el\otimes\elb - \half \elb\otimes\el + \gd +\EEE,
    \end{aligned}
  \end{align*}
  where
  \begin{align*}
    \EEE & := \le(- \half \frac{u}{r}\gd + \half \frac{\ub}{r}\gd -\gd \ri)-\quar \yy\el\otimes\el + \half u\le(\D\elb + \frac{1}{r}\gd\ri) + \half\ub\le(\D\el-\frac{1}{r}\gd\ri).
  \end{align*}
  Using the relations~\eqref{eq:Riccirel} and the bounds~(\ref{est:bdysourcenullconn}) for the null connection coefficients in $\MM^\ext$, we have
  \begin{align*}
    \norm{t^{3/2}\EEE}_{L^\infty(\TT)} & \les \varep.
  \end{align*}
  Thus, combining this bound together with the estimate~(\ref{est:DSI}) for $\D\SI$, we have
  \begin{align}\label{est:DSESI}
    \begin{aligned}
      \norm{t^{3/2}\D\le(\SE-\SI\ri)}_{L^\infty(\TT)} & \les \norm{t^{3/2}\le(\D\SE-\g\ri)}_{L^\infty(\TT)} + \norm{t^{3/2}\le(\D\SI-\g\ri)}_{L^\infty(\TT)} \\
      & \les \varep.
    \end{aligned}
  \end{align}

  Arguing as in the proof of Lemma~\ref{lem:ZZterr}, we have on $S^\ast$
  \begin{align*}
    \le|\XI-\rast\Nf\ri| & \les \varep \tast^{-1/2}.   
  \end{align*}
  Thus, from the definition~\eqref{eq:defSI} of $\SI$, we have
  \begin{align}\label{est:SISast}
    \le|\SI - t^\ast\Tf -\rast\Nf\ri| & \les \varep \tast^{-1/2},
  \end{align}
  on $S^\ast$. From the definition~\eqref{eq:defTESEKE}, relations~\eqref{eq:defnutnur} and~\eqref{eq:NTT} and the bounds~(\ref{est:OOmaxconn}) for $\nu$, we have
  \begin{align}\label{est:SESast}
    \le|\SE-\tast\Tf-r\Nf\ri| & \les \tast |\nu-1| \les \varep \tast^{-1/2}. 
  \end{align}

  Estimate~\eqref{est:SESI} now follows from combining~\eqref{est:SISast}, \eqref{est:SESast} and integrating~\eqref{est:DSESI} along $\TT$.\\

  Estimate~\eqref{est:KEKI} follows from estimates~\eqref{est:TETI} and~\eqref{est:SESI} and the approximate formulas
  \begin{align*}
    \le|\KX-2t\SX + (t^2-r^2)\TX\ri| & \les \varep t^{1/2},
  \end{align*}
  which we leave to the reader. This finishes the proof of the lemma.
\end{proof}

\begin{lemma}\label{lem:OOEOOITT}
  The following bounds hold on $\TT$
  \begin{align}
    \label{est:OOEOOI}
    \le|{^{(\ell)}\OOE}-{^{(\ell)}\OOI}\ri| & \les \varep t^{-1/2},
  \end{align}
  for all $\ell = 1,2,3$.
\end{lemma}
\begin{proof}
  We start by proving the following estimate for the covariant derivative of the interior vectorfields $\OOI$
  \begin{align}
    \label{est:DXfOOI}
    \le|\D_{\XI}{^{(\ell)}\OOI} -{^{(\ell)}\OOI}\ri| & \les \varep t^{-1/2},
  \end{align}
  for all $\ell =1,2,3$.\\
  
  By differentiating in $\Tf$, using the transport equations~\eqref{eq:defXI} and~\eqref{eq:defOOI} for $\XI$ and $\OOI$ and the sup-norm bounds~(\ref{est:controlTI}) for $\D\Tf$ and sup-norm bounds~(\ref{est:curvestsourceplanehyp}) for the curvature, it is enough to prove that~\eqref{est:DXfOOI} holds on the last slice $\Si_\tast$.\\

  From the definitions~\eqref{eq:defXISitast} and~\eqref{eq:defOOISitast} of $\XI$ and $\OOI$ on $\Si_\tast$, we have
  \begin{align*}
    \nab_{\XI}{^{(\ell)}\OOI} & = x^k\nab_{\nab x^k}\le(\in_{\ell i j}x^i\nab x^j\ri) \\
                              & = \in_{\ell i j}x^k \nab_{x^k}x^i \nab x^j + \in_{\ell i j}x^kx^i\nab_{\nab x^k}\nab x^j \\
                              & = \in_{\ell i j}x^i \nab x^j + \EEE,
  \end{align*}
  where
  \begin{align*}
    \EEE & = \in_{\ell i j}x^k \le(\nab_{x^k}x^i-\de_{ki}\ri) \nab x^j + \in_{\ell i j}x^kx^i\nab_{\nab x^k}\nab x^j.
  \end{align*}
  From the bounds~\eqref{est:sourceharmo}, we have
  \begin{align*}
    |\EEE| & \les \varep t^{-1/2}.
  \end{align*}
  Using additionally the control~\eqref{est:controlTI} of $\D\TI$, we deduce that the estimate~\eqref{est:DXfOOI} holds on $\Si_\tast$, and thus~\eqref{est:DXfOOI} is proved.\\

  Using that $\XI = -t\TI + 2\SI$, and the bounds~\eqref{est:SESI} obtained for $\SE-\SI$, the relations~\eqref{eq:defnutnur}, \eqref{eq:NTT} and the bounds for $\nut$, we have
  \begin{align}\label{est:XINf}
    \le|\XI - r\Nf\ri| & \les \varep t^{-1/2},
  \end{align}
  on $\TT$. Thus, we deduce from~\eqref{est:DXfOOI} that at $\TT$
  \begin{align}\label{est:DTfNfOOITT}
    \begin{aligned}
      \le|\D_\Tf\OOI\ri| & \les \varep t^{-3/2}, & \le|\D_\Nf\OOI - \frac{1}{r}\OOI\ri| & \les \varep t^{-3/2}.
    \end{aligned}
  \end{align}

  Using the definitions from Section~\ref{sec:defrotext}, the bounds~(\ref{est:OO1}) on $\YEi$ and on $\chi$ and $\chib$, we have for the exterior rotation vectorfields at $\TT$
  \begin{align*}
    \begin{aligned}
      \le|\D_3\OOE + \frac{1}{r}\OOE\ri| + \le|\D_4\OOE - \frac{1}{r}\OOE\ri| & \les \varep t^{-3/2}, 
    \end{aligned}
  \end{align*}
  from which, using the relations~\eqref{eq:defnutnur} and~\eqref{eq:NTT} and the bounds on $\nu$, we deduce
  \begin{align}\label{est:DTNOOETT}
    \begin{aligned}
      \le|\D_\Nf\OOE - \frac{1}{r}\OOE\ri| & \les \varep t^{-3/2}.
    \end{aligned}
  \end{align}

  Combining~\eqref{est:DTfNfOOITT} and~\eqref{est:DTNOOETT} and estimate~\eqref{est:ZZterr} for $\ZZt$ we obtain
  \begin{align}\label{est:DZOOEOOITT}
    \begin{aligned}
      \le|\D_{\ZZt}\le(\OOE-\OOI\ri) - \frac{1}{r}\le(\frac{1-\cc}{1+\cc}\ri)\le(\OOE-\OOI\ri)\ri| & \les \varep t^{-3/2},
    \end{aligned}
  \end{align}
  where we recall that $\ZZt$ was defined in Section~\ref{sec:XIXI} to be the $\TT$-tangent vectorfield normal to $\pr\Si_t$ in $\TT$ and such that $\ZZt(t) = 1$.\\

  From the definitions~\eqref{eq:defOOISitast} of $\OOI$ on $\Si_\tast$ and in particular on $S^\ast$, we have
  \begin{align}\label{eq:OOIOOESast}
    \begin{aligned}
      ^{(\ell)}\OOI & = \in_{\ell i j}x^i \nab x^j \\
      & = \in_{\ell i j} x^i \Nd x^j + \in_{\ell i j}x^i \Nf(x^j) \Nf \\
      & = \in_{\ell i j}x^i \Nd x^j + \rast^{-1}\in_{\ell i j}x^ix^j\Nf + \in_{\ell i j}x^i\le(\Nf(x^j)-\rast^{-1}x^j\ri)\Nf \\
      & = ^{(\ell)}\OOE + \in_{\ell i j}x^i\le(\Nf(x^j)-\rast^{-1}x^j\ri)\Nf.
    \end{aligned}
  \end{align}

  Thus, using the bounds~\eqref{est:sourceharmo}, we have on $S^\ast$
  \begin{align*}
    \le|\OOEi - {^{(\ell)}\OOI}\ri| & \les \le|\in_{\ell i j}x^i\le(\Nf(x^j)-\rast^{-1}x^j\ri)\ri| \les \varep \tast^{-1/2}.
  \end{align*}
  
  Integrating~\eqref{est:DZOOEOOITT} along $\TT$ then yields~\eqref{est:OOEOOI} and finishes the proof of the lemma.
\end{proof}

\begin{lemma}
  The following bounds hold on $\TT$
  \begin{align}\label{est:DpiXEXI2}
    \begin{aligned}
      \norm{t^{-\ga}t^{2}(t\D)^{\leq 1}\le(\D\TE-\D\TI\ri)}_{L^\infty_t L^4(\pr\Si_t)} & \les_\ga \varep, \\
      \norm{t^{-\ga}t(t\D)^{\leq 1}\le(\D\SE-\D\SI\ri)}_{L^\infty_t L^4(\pr\Si_t)} & \les_\ga \varep,
    \end{aligned}
  \end{align}
  for all $\ga>0$.
\end{lemma}
\begin{proof}
  The proof of the lemma follows directly from the precise control of each separate tensor already obtained in the previous sections and is left to the reader.
\end{proof}

\begin{lemma}\label{lem:DOOEOOITTL4}
  Under the same assumptions as in the previous lemma, we have
  \begin{align}\label{est:DOOEOOITTlem}
    \norm{t^{-\ga}t(t\D)^{\leq 1}\le(\D\OOE-\D\OOI\ri)}_{L^\infty_t L^4(\pr\Si_t)} & \les_\ga \varep,
  \end{align}
  for all $\ga>0$.
\end{lemma}
\begin{proof}
  The main part of the proof is to obtain the following two estimates
  \begin{align}
    \norm{t^{-\ga}t^2\D^2\OOI}_{L^\infty_t L^4(\pr\Si_t)} & \les \varep, \label{est:D2OOITT}\\
    \norm{t^{-\ga}t^2\D^2\OOE}_{L^\infty_t L^4(\pr\Si_t)} & \les \varep.\label{est:D2OOETT}
  \end{align}
  
  We first verify that we can obtain~\eqref{est:DOOEOOITTlem} from~\eqref{est:D2OOITT} and~\eqref{est:D2OOETT}. Using these last estimates and integration along $\ZZt$ as in the proof of Lemma~\ref{lem:OOEOOITT}, it is enough to additionally show that
  \begin{align}\label{est:DOOEOOISast}
    \norm{t(\D\OOE-\D\OOI)}_{L^4(S^\ast)} & \les \varep.
  \end{align}
  Using formula~\eqref{eq:OOIOOESast} and the bounds~\eqref{est:sourceharmo} on the harmonic coordinates, we have
  \begin{align}\label{est:DaOOEOOISast}
    \begin{aligned}
      \norm{t\D_a\le(\OOE-\OOI\ri)}_{L^4(S^\ast)} & \les \norm{t\D_a\le(\le(\in_{\ell i j}x^i(\Nf(x^j)-\rast^{-1}x^j)\ri)\Nf\ri)}_{L^4(S^\ast)} \\
      & \les \norm{t \Nd\le(\in_{\ell i j}x^i(\Nf(x^j)-\rast^{-1}x^j)\ri)}_{L^4(S^\ast)} + \norm{\in_{\ell i j}x^i(\Nf(x^j)-\rast^{-1}x^j) t\D_aN}_{L^4(S^\ast)} \\
      & \les t \sum_{j=1}^3\norm{(t\Nd)^{\leq 1}\le(\Nf(x^j)-\rast^{-1}x^j\ri)}_{L^4(S^\ast)} \\
      & \les \varep.
    \end{aligned}
  \end{align}
  The last line indeed follows from the identity
  \begin{align*}
    \nab_a\le(\Nf(x^i)-\rast^{-1}x^i\ri) & = \nab_a\nab_Nx^i + \le(\th_{ab}-\rast^{-1}\gd_{ab}\ri)\Nd_bx^i, 
  \end{align*}
  together with the bounds~\eqref{est:sourceharmo} and~(\ref{est:bdysourcenullconn}), (\ref{est:OOmaxconn}) and Sobolev estimates on $\Si_t$.\\

  Using the estimates~\eqref{est:DTfNfOOITT} and~\eqref{est:DTNOOETT} already obtained in the proof of Lemma~\ref{lem:OOEOOITT} for respectively $\OOI$ and $\OOE$ on $\TT$, we have in particular on $S^\ast$ and for the $L^4$ norm
  \begin{align*}
    \norm{t\D_\Tf\OOI}_{L^4(S^\ast)} + \norm{t\le(\D_\Nf\OOI-\frac{1}{r}\OOI\ri)}_{L^4(S^\ast)} & \les \varep,\\
    \norm{t\D_\Tf\OOE}_{L^4(S^\ast)} + \norm{t\le(\D_\Nf\OOE -\frac{1}{r}\OOE\ri)}_{L^4(S^\ast)} & \les \varep.
  \end{align*}
  These last estimates combined with~\eqref{est:DaOOEOOISast} yield~\eqref{est:DOOEOOISast} as desired.\\

  We now turn to the proof of~\eqref{est:D2OOITT} and~\eqref{est:D2OOETT}, where we note that estimate~\eqref{est:D2OOITT} was already obtained in Lemma~\ref{lem:D2OOIlem}. To prove~\eqref{est:D2OOETT}, we compute that for each component of $\D^2_{\mu,\nu}\OOE$, we have
  \begin{align}\label{eq:D2OOEmunu}
    \D^2_{\mu,\nu}\OOE & = \EEE^1\cdot\FFF^0 + \EEE^1 \cdot \EEE^0,
  \end{align}
  where
  \begin{align*}
    \EEE^1 && \text{is a linear combination of} && \le(\D^{\leq 1}\Gac, \D^{\leq 1}Y, \D^{\leq 1}\Hrot, \POE , \R\ri), \\
    \EEE^0 && \text{is a (nonlinear) combination of} && \le(\Gac, Y, \Hrot\ri),\\
    \FFF^0 && \text{is a (nonlinear) combination of} && \le(1,\D\OOE\ri),
  \end{align*}
  with
  \begin{align*}
    \Gac & \in\le\{\trchi-\frac{2}{r}, \trchib+\frac{2}{r}, \chih, \chibh, \ze, \omb, \xib\ri\}.
  \end{align*}
  Note that we do not precise the weights in $u,\ub,r,t$ in~\eqref{eq:D2OOEmunu}, which are recovered by a simple scaling consideration on $\TT$ where $u\simeq \ub \simeq r \simeq t$. Using formula~\eqref{eq:D2OOEmunu}, the $L^\infty_t\HHt(\pr\Si_t)$ estimates of~(\ref{est:bdysourcenullconn}) for $\Gac,Y,\Hrot,\POE,\R$ and $\Nd\OOE$, the $\HHt$ product estimates of Lemma~\ref{lem:prodH12} and the Sobolev embeddings of Lemma~\ref{lem:Sobsphere}, we obtain the desired estimate~\eqref{est:D2OOETT}.\\

  The rest of the proof is therefore dedicated to obtaining the formulas~\eqref{eq:D2OOEmunu}. For simplicity, we call $\OOO$ any exterior rotation vectorfield $^{(\ell)}\OOE$. From~(\ref{eq:Lie4OOE}) and relations~\eqref{eq:Riccirel}, we have
  \begin{align}
    \label{eq:D4OOE}
    \D_4\OOO & = \D_\OOO\el = \chi_{ab}\OOO^b\ea - \ze_b\OOO^b\el = r^{-1}\OOO + \EEE^0,  
  \end{align}
  from which we also deduce
  \begin{align}\label{eq:D4rOOE}
    \D_4\le(r^{-1}\OOO\ri) & = \EEE^0.
  \end{align}
  Similarly, using the relations~(\ref{eq:relOOEnull}), we have
  \begin{align}\label{eq:D3OOE}
    \begin{aligned}
    \D_3\OOO & = -r^{-1}\OOO + \EEE^0, \\
    \D_3(r^{-1}\OOO) & = \EEE^0.
    \end{aligned}
  \end{align}

  Deriving equation~\eqref{eq:D4OOE} by $\D_4$ and using~\eqref{eq:D4rOOE} and relations~\eqref{eq:Riccirel}, we obtain
  \begin{align*}
    \D_{4,4}^2\OOO & = \D_4(\D_4\OOO) + \EEE^1\cdot\FFF^0 = \D_4(r^{-1}\OOO) + \D_4\le(\EEE^0\ri) + \EEE^1\cdot\FFF^0 = \EEE^1\cdot\FFF^0 + \EEE^1\cdot\EEE^0.
  \end{align*}
  Arguing similarly using $\D_3,\D_4$ derivatives, we have
  \begin{align*}
    \D_{4,4}^2\OOO,~\D_{3,4}^2\OOO,~\D_{4,3}^2\OOO,~\D^2_{3,3}\OOO & = \EEE^1\cdot\FFF^0 + \EEE^1\cdot\EEE^0.
  \end{align*}

  For $\D^2_{a,4}$, using equation~\eqref{eq:D4OOE}, we have
  \begin{align*}
    \D^2_{a,4}\OOO & = \D_a(\D_4\OOO) - (\D_a\el)^{\mu}\D_{\mu}\OOO \\
                   & = \D_a(\D_4\OOO) - r^{-1}\D_a\OOO + \EEE^1\cdot\FFF^0 \\
                   & = \D_a(r^{-1}\OOO) -r^{-1}\D_a\OOO + \EEE^{1}\cdot\FFF^0 + \EEE^1\cdot\EEE^0 \\
                   & = \EEE^1\cdot\FFF^0 + \EEE^1\cdot\EEE^0.
  \end{align*}
  Using the above, we also obtain for $\D^2_{4,a}$
  \begin{align*}
    \D^2_{4,a}\OOO & = \D^2_{a,4}\OOO + \R = \EEE^1\cdot\FFF^0 + \EEE^1\cdot\EEE^0.
  \end{align*}
  Arguing similarly for the $\elb$ derivatives, we also obtain
  \begin{align*}
    \D^2_{3,a}\OOO,~\D^2_{a,3}\OOO & =  \EEE^1\cdot\FFF^0 + \EEE^1\cdot\EEE^0.
  \end{align*}

  We now turn to $\D^2_{a,b}\OOO$. From the result of Lemma~\ref{lem:D2Nd2X}, using~\eqref{eq:D4OOE} and~\eqref{eq:D3OOE}, and the definitions of $\Hrot$, $Y$ and $\POE$ from Section~\ref{sec:defrotext}, we have
  \begin{align*}
    \D^2_{a,b}\OOO & = \Nd^2_{a,b}\OOO-r^{-2}\OOO_b\ea + \half r^{-1}\gd_{ab}\le(\D_4-\D_3\ri)\OOO \\
                   & \quad + \half r^{-1}\le(\Nd_a\OOO_b+\Nd_b\OOO_a\ri)(\elb-\el) + \le(\EEE\le(\D^2,\Nd^2\ri)\cdot\OOO\ri)_{ab},\\
                   & = \Nd^2_{a,b}\OOO -r^{-2}\OOO_b\ea + r^{-2}\gd_{ab}\OOO + \half r^{-1}\Hrot_{ab}(\elb-\el) + \EEE^1\cdot\FFF^0 \\
                   & = \POE_{abc}\ec + \half r^{-1}\Hrot_{ab}(\elb-\el) + \EEE^1\cdot\FFF^0 \\
                   & = \EEE^1\cdot\FFF^0,
  \end{align*}
  where we used that from the formulas of Lemma~\ref{lem:D2Nd2X}
  \begin{align*}
    \EEE\le(\D^2,\Nd^2\ri)\cdot\OOO & = \EEE^1\cdot\FFF^0.
  \end{align*}
  This finishes the proof of~\eqref{eq:D2OOEmunu} and concludes the proof of the lemma.
\end{proof} 

\chapter{The initial layers}\label{sec:initlayer}
The goal of this section is twofold.\\

In Section~\ref{sec:initlayerenergy} we obtain the initial bounds of Proposition~\ref{prop:initenercorrec} for the energy fluxes of (contracted and commuted) Bel-Robinson tensors through $\Si_\too$ and $\CC_1$, which are used in Section~\ref{sec:globener} to obtain improved energy fluxes bounds in $\MM$. \\

In Sections~\ref{sec:lastconesfoliation} and \ref{sec:controltranscoeff}, we use the improved estimates of Sections~\ref{sec:globener}--\ref{sec:planehypconnest} to improve the Bootstrap Assumptions~\ref{BA:bottom} and~\ref{BA:con} on the comparison between the geometric constructions in $\MM$ and in $\LLb$ and $\LLc$, which were used in Section~\ref{sec:initlayerenergy}. This is done in two steps: in Section~\ref{sec:lastconesfoliation} we control the intermediate foliation by the last cones used to define the region $\MM$ (see its definition in Section~\ref{sec:deflastconesfoliation}). In Section~\ref{sec:controltranscoeff} we use this intermediate control to obtain the desired comparisons in $\LLb$ and $\LLc$.


\section{Initial bounds for energy fluxes through $\CC_1$ and $\Si_\too$}\label{sec:initlayerenergy}
This section is dedicated to the proof of Proposition~\ref{prop:initenercorrec}. We assume that the Bootstrap Assumptions hold and that the initial layers $\LLb$ and $\LLc$ are $\varep$-close to Minkowski (see the definitions of Section~\ref{sec:Minkowskilayer}), and we show that
\begin{align}\label{est:initenerSitoo}
  \int_{\Si_\too} \le|\D^{\leq 2}\R\ri|^2 & \les \varep^2,
\end{align}
and that
\begin{align}\label{est:initenerCC1}
  \begin{aligned}
    \int_{\CC_1} \bigg(\le|\Ndt^{\leq 2}\beb\ri|^2 + \le|\ub\Ndt^{\leq 2}(\rho-\rhoo)\ri|^2 + \le|\ub\Ndt^{\leq 2}(\sigma-\sigmao)\ri|^2 + \le|\ub^2\Ndt^{\leq 2}\be\ri|^2 + \le|\ub^2\Ndt^{\leq 2}\al\ri|^2 \bigg) \les \varep^2
  \end{aligned}
\end{align}
where $\Ndt\in\le\{(r\Nd), (\ub\Nd_4), (u\Nd_3)\ri\}$.

\begin{remark}
Together with the null connection and rotation coefficients Bootstrap Assumptions~\ref{BA:connext} on $\CC_1$, this proves the desired bounds for the Bel-Robinson tensors~\eqref{est:initener} of Proposition~\ref{prop:initenercorrec} used in Section~\ref{sec:globener}. Details are left to the reader. 
\end{remark}

\subsection{Energy fluxes through $\Si_\too$}\label{sec:fluxSitoo}
Under the Bootstrap Assumptions~\ref{BA:bottom}, the hypersurface $\Si_\too$ is included in $\LLb$. One can therefore perform energy estimates in the past region of $\Si_\too$ in $\LLb$ for the following contracted and commuted Bel-Robinson tensors
\begin{align*}
  Q\le(\Lieh_{\pr_\mu}^{\leq 2}\R\ri)\le(\Tf^\bott,\Tf^\bott,\Tf^\bott\ri),
\end{align*}
for $\mu=0,1,2,3$. Arguing as in Section~\ref{sec:errorintr} to estimate the error terms using the bottom initial layer estimates~\eqref{est:botlayass}, one obtains 
\begin{align*}
  \int_{\Si_\too} Q\le(\Lieh^{\leq 2}_{\pr_\mu}\R\ri)\le(\Tf^\bott,\Tf^\bott,\Tf^\bott,\Tf\ri) & \les \int_{\Sit_1} Q\le(\Lieh^{\leq 2}_{\pr_\mu}\R\ri)\le(\Tf^\bott,\Tf^\bott,\Tf^\bott,\Tf^\bott\ri) + \varep^3.
\end{align*}
Using the curvature flux bound~\eqref{est:initenersourceSit1bis} on $\Sit_1$ and the bounds~(\ref{est:botlayassbis}), one has
\begin{align*}
  \int_{\Sit_1}Q\le(\Lieh^{\leq 2}_{\pr_\mu}\R\ri)\le(\Tf^\bott,\Tf^\bott,\Tf^\bott,\Tf^\bott\ri) & \les \varep^2,
\end{align*}
and thus
\begin{align}\label{est:curvfluxSitoo}
  \int_{\Si_\too} Q\le(\Lieh^{\leq 2}_{\pr_\mu}\R\ri)\le(\Tf^\bott,\Tf^\bott,\Tf^\bott,\Tf\ri) & \les \varep^2.
\end{align}
Using the Bootstrap Assumptions~\ref{BA:bottom} for the comparisons of the frames, we deduce from~\eqref{est:curvfluxSitoo}
\begin{align}
  \label{est:curvfluxSitoobisbis}
  \int_{\Si_\too} \le|\Lieh_{\pr_\mu}^2\R\ri|^2 & \les \varep^2,
\end{align}
where the norm is taken with respect to the frame with respect to either $\Tf^\bott$ or $\Tf$.\\

Using the bounds~\eqref{est:botlayass}, one can obtain by integration along $x^0$
\begin{align}\label{est:L2SitooDDpr0}
  \norm{\D^{\leq 1}\D\pr_\mu}_{L^2(\Si_\too)} & \les \norm{\pr^{\leq 2}\le(\g_{\al\be}-\etabold_{\al\be}\ri)}_{L^2(\Si_\too)} \les \norm{\pr^{\leq 3}(\g_{\al\be}-\etabold_{\al\be})}_{L^2(\LLb)} \les \varep. 
\end{align}
Using that we schematically have
\begin{align*}
  \D_{\pr_\mu}^2\R & = \Lieh_{\pr_\mu}^2\R + \D\R\cdot\D\pr_\mu + \R\cdot \D^2\pr_\mu,
\end{align*}
we deduce from~\eqref{est:curvfluxSitoobisbis} and~\eqref{est:L2SitooDDpr0}
\begin{align*}
  \int_{\Si_\too}\le|\D_{\pr_\mu}^{\leq 2}\R\ri|^2 \les \varep^2,
\end{align*}
for $\mu=0,1,2,3$, from which using~\eqref{est:botlayass} again, we deduce
\begin{align*}
  \int_{\Si_\too}\le|\D^{\leq 2}\R\ri|^2 & \les \varep^2,
\end{align*}
in either the frame with respect to $\Tf^\bott$ or $\Tf$. This finishes the proof of~\eqref{est:initenerSitoo}.

\subsection{Energy fluxes through $\CC_1$}
\begin{figure}[h!]
  \centering
  \includegraphics[height=10cm]{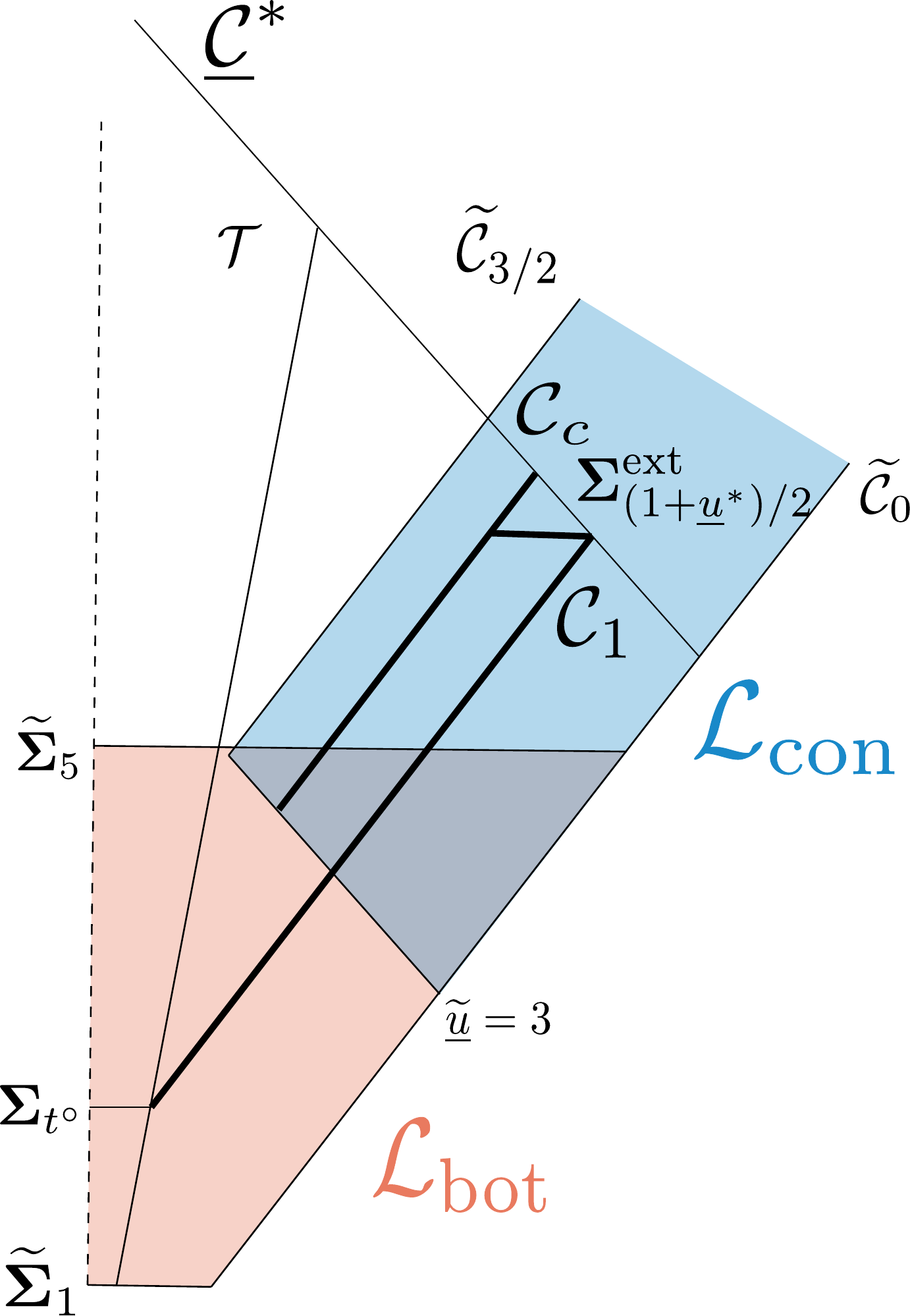}
  \caption{Energy flux through $\CC_1$}
  \label{fig:flux_CC1}
\end{figure}

Under the Bootstrap Assumptions~\ref{BA:bottom} and~\ref{BA:con} (see also Figure~\ref{fig:flux_CC1}), we have
\begin{align*}
  \CC_1 & = \le(\CC_1\cap\LLb\ri) \cup \le(\CC_1\cap\LLc\ri). 
\end{align*}
We first obtain bounds on the bottom initial layer part $\CC_1\cap\LLb$, using the same contracted and commuted Bel-Robinson tensors as in Section~\ref{sec:fluxSitoo} and the bottom initial layer assumptions~\eqref{est:botlayass} to control the error terms, and we have
\begin{align}\label{est:prelimbelrobfluxCC1LLb}
  \int_{\CC_1\cap\LLb} Q\le(\Lieh^{\leq 2}_{\pr_\mu}\R\ri)\le(\Tf^\bott,\Tf^\bott,\Tf^\bott,\el\ri) & \les \varep^2.
\end{align}
Using~\eqref{est:botlayass}, the Bootstrap Assumptions~\ref{BA:bottom} in $\LLbext$ and arguing as in Section~\ref{sec:fluxSitoo}, we deduce from~\eqref{est:prelimbelrobfluxCC1LLb} that the desired bound~(\ref{est:initenerCC1}) holds in $\CC_1\cap\LLb$.\\

Performing a mean value argument between $\CC_1\cap\LLc$ and $\CC_{4/3}\cap\LLc$ using the (integrated) curvature flux bound~\eqref{est:curvfluxCCbt0bis} (see for example Section~\ref{sec:meanvalue}), one can obtain
\begin{align}\label{est:meanvalueCC1}
  \begin{aligned}
    \int_{\CC_c\cap\LLc}\ubt^{-1-\ga}\le|\Ndtt^{\leq 2}\alphabt\ri|^2 & \les_\ga \varep^2,
  \end{aligned}
\end{align}
where $\Ndtt\in\le\{\ubt\widetilde{\Nd},\ubt\widetilde{\Nd}_4,\widetilde{\Nd}_3\ri\}$, for $\ga>0$ and where $1<c<4/3$.\\

Using energy estimates in the spacetime region between $\CCt_0\cup\{\ub'=3\}$ and $\le(\CC_c\cap\LLc\ri)\cup\le(\CCba\cap\LLc\ri)$,\footnote{To obtain no uncontrolled error terms on $\CCba$, we would rather perform en energy estimate replacing $\CCba\cap\LLc$ by a spacelike boundary hypersurface which can be obtained by extension above $\CCba$ in $\LLc$.} for the null Bianchi equations satisfied by the null curvature components with respect to the null pair of the double null foliation of the conical initial layer, commuting with $\ubt\widetilde{\Nd}$ and $\widetilde{\Nd}_3,\widetilde{\Nd}_4$, we claim that one can obtain the following curvature flux control on $\CC_c$ (this follows from an adaptation of~\cite{Li.Zhu18})
\begin{align}\label{est:fluxprimeCC1}
  \begin{aligned}
    &   \int_{\CC_c\cap\LLc}\bigg(\le|\Ndtt^{\leq 2}\betabt\ri|^2+ \le|\ubt\Ndtt^{\leq 2}(\rhot-\rhoot)\ri|^2 + \le|\ubt\Ndtt^{\leq 2}(\sigmat-\sigmaot)\ri|^2 + \le|\ubt^2\Ndtt^{\leq 2}\betat\ri|^2 + \le|\ubt^2\Ndtt^{\leq 2}\alphat\ri|^2\bigg) \\
    & \les \varep^2 + \int_{\CC_c\cap\LLc} (D\varep)^2\ubt^{-2}\bigg(\le|\Ndtt^{\leq 2}\alphabt\ri|^2 + \le|\ubt\Ndtt^{\leq 2}\betabt\ri|^2+ \le|\ubt^2\Ndtt^{\leq 2}(\rhot-\rhoot)\ri|^2 \\
    & \quad \quad \quad \quad + \le|\ubt^2\Ndtt^{\leq 2}(\sigmat-\sigmaot)\ri|^2 + \le|\ubt^2\Ndtt^{\leq 2}\betat\ri|^2\bigg)
\end{aligned}
\end{align}
where $\Ndtt \in\le\{\ubt\widetilde{\Nd},\ubt\widetilde{\Nd}_4,\widetilde{\Nd}_3\ri\}$.
\begin{remark}
  The $(D\varep)^2\ubt^{-2}$ term comes from the fact that using the relation from Lemma~\ref{lem:transframe}, the normal $\el$ to $\CC_c$ can be expressed (at first orders) as $\el = \la^{-1}\elt - \la^{-1}\ft^a\tilde{\ea} -\quar |\ft|^2\elbt$, and thus the difference between $\el$ and $\el'$ yields the presence of $|\ft|^2|\alphabt|^2$, which justifies the claim.
\end{remark}

Thus, we deduce from~\eqref{est:meanvalueCC1} and~\eqref{est:fluxprimeCC1} that
\begin{align}
  \label{est:fluxprimeCC1bis}
  \begin{aligned}
    \int_{\CC_c\cap\LLc}\le|\Ndtt^{\leq 2}\betabt\ri|^2+ \le|\ubt\Ndtt^{\leq 2}(\rhot-\rhoot)\ri|^2 + \le|\ubt\Ndtt^{\leq 2}(\sigmat-\sigmaot)\ri|^2  + \le|\ubt^2\Ndtt^{\leq 2}\betat\ri|^2 + \le|\ubt^2\Ndtt^{\leq 2}\alphat\ri|^2\les \varep^2. 
  \end{aligned}
\end{align}

Using the transition relations from Propositions~\ref{prop:chgframederiv}, \ref{prop:transconn} and \ref{prop:transcurv} and the Bootstrap Assumptions~\ref{BA:con} on the transition coefficients $\ft,\fbt,\lat$, one can deduce from~\eqref{est:fluxprimeCC1bis}
\begin{align}\label{est:fluxnotprimeCCc}
  \begin{aligned}
    \int_{\CC_c\cap\LLc}\le|\Ndt^{\leq 2}\beb\ri|^2 + \le|\ub\Ndt^{\leq 2}(\rho-\rhoo)\ri|^2 + \le|\ub\Ndt^{\leq 2}(\sigma-\sigmao)\ri|^2 + \le|\ub^2\Ndt^{\leq 2}\be\ri|^2 + \le|\ub^2\Ndt^{\leq 2}\al\ri|^2  \les \varep^2
  \end{aligned}
\end{align}
where $\Ndt\in\le\{(r\Nd), (\ub\Nd_4), \Nd_3\ri\}$.
\begin{remark}
  From an inspection of the relations from Propositions~\ref{prop:transconn} and~\ref{prop:transcurv}, the bootstrap bounds $|\ft|\les D\varep \ub^{-1}$, $|\fbt|,|\log\lat| \les D\varep$ are sufficient due to a conservation of signature principle, \emph{i.e.} $\lat,\fbt$ are only paired with lower signature null curvature component hence are not required to decay, while the transition coefficient $\ft$ which is paired with higher signature components satisfies the sufficient decay $\ub^{-1}$. See also~\cite[Remark 4.1.4]{Kla.Sze17} for further discussion.
\end{remark}
\begin{remark}
  We use the relations of Proposition~\ref{prop:transconn} to obtain the control of all derivatives of $\ft,\fbt,\lat$, from the Bootstrap Assumptions~\ref{BA:connext} on the null connection coefficients of $\MM^\ext$, the Bootstrap Assumptions~\ref{BA:con} on $\ft,\fbt,\lat$ and the conical initial layer assumptions~\eqref{est:conlayass} on the connection coefficients.\footnote{$L^2(\CC_c)$ bounds for (one derivative of) the initial conical layer null connection coefficients can be obtained by integration from $\LLc$ to $\CC_c$ of the $L^2(\LLc)$ estimates of~\eqref{est:conlayass}, as in the bottom initial layer case. See Section~\ref{sec:fluxSitoo}.} 
\end{remark}

To obtain the desired bound on $\CC_1\cap\LLc$ from the bound~\eqref{est:fluxnotprimeCCc} on $\CC_c\cap\LLc$, we first obtain an energy bound on the spacelike hypersurface $\Si^\ext_{(1+\uba)/2}\cap\LLc$, arguing similarly as for~\eqref{est:fluxnotprimeCCc}, except that no mean value argument is required to control the error terms. Then performing an energy estimate as in Section~\ref{sec:globener} in the spacetime region comprised between $\le(\CC_1\cap\LLc\ri) \cup \{\ub'=3\}$ and $(\CC_c\cap\LLc)\cup\le(\Si^\ext_{(1+\uba)/2}\cap\LLc\ri)$, one deduces the desired estimate~(\ref{est:initenerCC1}) in $\CC_1\cap\LLc$, which combined with the bound obtained in $\CC_1\cap\LLb$ finishes the proof of~\eqref{est:initenerCC1}.

\section{Control of the last cones geodesic foliation}\label{sec:lastconesfoliation}
This section is dedicated to the independent control of the last cones geodesic foliation. The final goal of these results is the improvement of the Bootstrap Assumptions~\ref{BA:lastconesfoliation}, which will be obtained in Section~\ref{sec:controltranscoeff}.\\

We have the following bounds for the null curvature components of the last cones geodesic foliation.
\begin{lemma}\label{lem:curvgeodsource}
  For all $11/4 \leq {\ub'} \leq \uba,~5/4 \leq {u'} \leq \uba$, the following bounds hold
  \begin{align}
    \begin{aligned}
    |\alpha'| & \les \varep {\ub'}^{-7/2}, & |\beta'| & \les \varep {\ub'}^{-7/2},\\
    |\rho'| & \les \varep {\ub'}^{-3}{u'}^{-1/2}, & |\sigma'| & \les \varep {\ub'}^{-3}{u'}^{-1/2},\\
    |\betab'| & \les \varep {\ub'}^{-2}{u'}^{-3/2}, & |\alphab'| & \les \varep {\ub'}^{-1}{u'}^{-5/2}.
  \end{aligned}
  \end{align}
\end{lemma}
\begin{proof}
  In the exterior region, these bounds follow from the improved curvature bounds $\mathfrak{R}^\ext_{\leq 1} \les \varep$, the Bootstrap Assumptions~\ref{BA:lastconesfoliation} and the transition formulas from Proposition~\ref{prop:transcurv}. In the bottom interior region, these bounds follow from the the improved bounds in the maximal frame $\mathfrak{R}_{\leq 1}^\intr \les \varep$ and the Bootstrap Assumptions~\ref{BA:lastconesfoliation}. In the top interior region, these bounds follow from the following pointwise bound on the spacetime curvature tensor (where the norm is taken with respect to the frame $(\elb',\el')$)
  \begin{align*}
    |\R| & \les \varep \tast^{-7/2}.
  \end{align*}
  This bound holds by a $\tast$-rescaling, extension and local existence as in Section~\ref{sec:curvest}, and comparison between the frames given by the local existence and the frame $(\elb',\el')$. See Section~\ref{sec:tastSitast} and in particular estimate~\eqref{est:curvextension}. This finishes the proof of the lemma. 
\end{proof}

We have the following control for the null connection coefficients of the last cones geodesic foliation. 
\begin{lemma}\label{lem:conngeodcontrol}
  For all $11/4 \leq {\ub'} \leq \uba,~5/4 \leq {u'} \leq \uba$, the following bounds hold
  \begin{align}\label{est:supnormconngeodcontrol}
    \begin{aligned}
      \le|\trchib'+\frac{2}{{r'}}\ri| & \les \varep {\ub'}^{-2}{u'}^{-1/2}, & \le|\chibh'\ri| & \les \varep {\ub'}^{-1}{u'}^{-3/2},\\
      |\ze'| & \les \varep {\ub'}^{-2}{u'}^{-1/2}, & |\om'| & \les \varep {\ub'}^{-5/2},
    \end{aligned}
  \end{align}
  and
  \begin{align}\label{est:L4normconngeodcontrolprelim}
    \begin{aligned}
      \norm{{r'}^{-1/2}(r\Nd)^{\leq 1}\le(\trchib'+\frac{2}{{r'}}\ri)}_{L^\infty L^4(S')} & \les \varep {\ub'}^{-2}{u'}^{-1/2}, & \norm{{r'}^{-1/2}(r\Nd)^{\leq 1}\chibh'}_{L^\infty L^4(S')} & \les \varep {\ub'}^{-1}{u'}^{-3/2},\\
      \norm{{r'}^{-1/2}(r\Nd)^{\leq 1}\ze'}_{L^\infty L^4(S')} & \les \varep {\ub'}^{-2}{u'}^{-1/2}, & \norm{{r'}^{-1/2}(r\Nd)^{\leq 1}\om'}_{L^\infty L^4(S')} & \les \varep {\ub'}^{-5/2},
    \end{aligned}
  \end{align}
  and
  \begin{align}\label{est:L4normconngeodcontrol}
    \begin{aligned}
      \norm{{r'}^{-1/2}\le(\trchi'-\frac{2}{{r'}}\ri)}_{L^\infty L^4(S')} & \les \varep {\ub'}^{-2}{u'}^{-1/2}, & \norm{{r'}^{-1/2}\chih'}_{L^\infty L^4(S')} & \les \varep {\ub'}^{-2}{u'}^{-1/2},\\
      \norm{{r'}^{-1/2}\xi'}_{L^\infty L^4(S')} & \les \varep {\ub'}^{-5/2},
    \end{aligned}
  \end{align}
  together with the additional relations (see~\eqref{eq:rellastconegeod1} and~\eqref{eq:rellastconegeod2})
  \begin{align}\label{eq:relxibtombtzet}
    \begin{aligned}
      \xib' & = 0,&  \omb' & = 0, & \eta' & = \ze' = -\etab',
    \end{aligned}
  \end{align}
  and where we denoted by ${r'}$ the area radius of the $2$-spheres $S'$.\\

  Moreover, we have the following control for the area radius ${r'}$
  \begin{align}\label{est:arearadiuslastconesgeod}
    \le|{r'}-\half({\ub'}-{u'})\ri| & \les \varep {r'} {\ub'}^{-2}{u'}^{-1/2},
  \end{align}
  and for the optical defect $\yy'$
  \begin{align}\label{est:yyt}
  |\yy'| & \les \varep {\ub'}^{-3/2}.
  \end{align}
\end{lemma}

\section{Proof of Lemma~\ref{lem:conngeodcontrol}}
This section is dedicated to the proof of Lemma~\ref{lem:conngeodcontrol}. To ease the notations and since no confusion is possible, we drop the primes in this section.\\

The proof of Lemma~\ref{lem:conngeodcontrol} is based on a continuity argument on the maximal domain starting from the central axis such that a set a bootstrap assumptions holds.\\

As bootstrap assumptions, we assume the conclusions of Lemma~\ref{lem:conngeodcontrol}, \emph{i.e.}
\begin{subequations}\label{est:BAnullgeodconnection}
\begin{align}\label{est:BAsupnormconngeodcontrol}
  \begin{aligned}
    \le|\trchib+\frac{2}{{r}}\ri| & \leq D\varep {\ub}^{-2}{u}^{-1/2}, & \le|\chibh\ri| & \leq D\varep {\ub}^{-1}{u}^{-3/2},\\
    |\ze| & \leq D\varep {\ub}^{-2}{u}^{-1/2}, & |\om| & \leq D\varep {\ub}^{-5/2},
    \end{aligned}
\end{align}
and
\begin{align}\label{est:BAL4normconngeodcontrolprelim}
  \begin{aligned}
    \norm{{r}^{-1/2}(r\Nd)^{\leq 1}\le(\trchib+\frac{2}{{r}}\ri)}_{L^\infty L^4(S)} & \leq D\varep {\ub}^{-2}{u}^{-1/2}, & \norm{{r}^{-1/2}(r\Nd)^{\leq 1}\chibh}_{L^\infty L^4(S)} & \leq D\varep {\ub}^{-1}{u}^{-3/2},\\
    \norm{{r}^{-1/2}(r\Nd)^{\leq 1}\ze}_{L^\infty L^4(S)} & \leq D\varep {\ub}^{-2}{u}^{-1/2}, & \norm{{r}^{-1/2}(r\Nd)^{\leq 1}\om}_{L^\infty L^4(S)} & \leq D\varep {\ub}^{-5/2},
  \end{aligned}
\end{align}
and
\begin{align}\label{est:BAL4normconngeodcontrol}
  \begin{aligned}
    \norm{{r}^{-1/2}\le(\trchi-\frac{2}{{r}}\ri)}_{L^\infty L^4(S)} & \leq D\varep {\ub}^{-2}{u}^{-1/2}, & \norm{{r}^{-1/2}\chih}_{L^\infty L^4(S)} & \leq D\varep {\ub}^{-2}{u}^{-1/2},\\
    \norm{{r}^{-1/2}\xi}_{L^\infty L^4(S)} & \leq D\varep {\ub}^{-5/2}.
  \end{aligned}
\end{align}

Moreover, we assume the following control for the area radius ${r}$
\begin{align}\label{est:BAarearadiuslastconesgeod}
  \le|{r}-\half({\ub}-{u})\ri| & \leq D\varep {r} {\ub}^{-2}{u}^{-1/2},
\end{align}
and for the optical defect $\yy$
\begin{align}\label{est:BAyyt}
  |\yy| & \leq D\varep {\ub}^{-3/2}.
\end{align}

We also assume the following bounds on the metric in transported from the vertex spherical coordinates $\varth,\varphi$ from Theorem~\ref{thm:vertex}
\begin{align}\label{est:BAmetriccoordinatesgeodnull}
  \begin{aligned}
    \le|r^{-2}\gd_{ab}-(\gd_\SSS)_{ab}\ri| & \leq D\varep \ub^{-1}u^{-3/2}\qq,\\
    \norm{r^{-1/2}\pr_c^{\leq 1}\le(r^{-2}\gd_{ab}-(\gd_{\SSS})_{ab}\ri)}_{L^4(S_{u,\ub})} & \leq D\varep \ub^{-1}u^{-3/2}\qq,
  \end{aligned}
\end{align}
where $a,b,c\in\{\varth,\varphi\}$.\\
\end{subequations}

Under these bootstrap assumptions, we have the following three preliminary lemmas.
\begin{lemma}[$L^4(S)$ evolution estimates]\label{lem:transportgeodnull}
  Under the bootstrap bounds~\eqref{est:BAnullgeodconnection}, for all $S$-tangent tensor $U$ satisfying the transport equation
  \begin{align*}
    \Nd_{3}U + \half\kappa\trchib U & = F,
  \end{align*}
  where $\kappa \in \RRR$, and satisfying the vertex limit
  \begin{align*}
    r^\kappa|U| \to 0,
  \end{align*}
  when $r\to 0$, we have
  \begin{align*}
    \norm{r^{\kappa-\frac{1}{2}}(r\Nd)^{\leq \ell}U}_{L^4(S_{u,\ub})} & \les \int_{u}^{\ub}\norm{r^{\kappa-\frac{1}{2}}(r\Nd)^{\leq \ell}F}_{S_{u',\ub}} \,\d u',
  \end{align*}
  for $\ell\in\{0,1\}$.
\end{lemma}
\begin{proof}
  The proof follows from renormalisation and integration along $u$ as in the proof of Lemma~\ref{lem:evolext}. It uses the bootstrap bounds~\eqref{est:BAsupnormconngeodcontrol},~\eqref{est:BAarearadiuslastconesgeod} and~\eqref{est:BAmetriccoordinatesgeodnull}. Details are left to the reader.
\end{proof}
\begin{lemma}[$L^4(S)$ elliptic estimates]\label{lem:ellgeodnull}
  Under the bootstrap bounds~\eqref{est:BAnullgeodconnection} and the estimates of Lemma~\ref{lem:curvgeodsource}, the following holds.\\

  For all $S$-tangent tensor $U$ of appropriate type, we have
  \begin{align*}
    \norm{r^{-1}(r\Nd)^{\leq 1}U}_{L^4(S)} & \les \norm{\Dd_1U}_{L^4(S)},\\
    \norm{r^{-1}(r\Nd)^{\leq 1}(U-\overline{U})}_{L^4(S)} & \les \norm{\Dd_1^\ast U}_{L^4(S)},\\
    \norm{r^{-1}(r\Nd)^{\leq 1}U}_{L^4(S)} & \les \norm{\Dd_2U}_{L^4(S)},
  \end{align*}
  for all $2$-sphere $S$ of the last cones geodesic null foliation. Moreover, for all $S$-tangent traceless symmetric $2$-tensor $U$, satisfying
  \begin{align*}
    \Divd U & = \Nd F + G,
  \end{align*}
  we have
  \begin{align}
    \label{est:ellspecialL4geodnull}
    r^{-1/2}\norm{U}_{L^4(S)} & \les r^{-1/2}\norm{F}_{L^4(S)} + r^{-1/2}\norm{rG}_{L^4(S)}.
  \end{align}
\end{lemma}
\begin{proof}
  Using the bootstrap bounds~\eqref{est:BAnullgeodconnection} and the estimates of Lemma~\ref{lem:curvgeodsource}, we have
  \begin{align*}
    r^{-1/2}\norm{K-\frac{1}{r^2}}_{L^4(S)} & \les D\varep \ub^{-3}u^{-1/2}. 
  \end{align*}
  Using rescaling of these $L^4(S)$ bounds for the Gauss curvature and the results of~\cite[Chapter 8]{Bie.Zip09}, the results of the lemma follow. 
\end{proof}
\begin{lemma}[$L^4(S)$ Sobolev estimates]\label{lem:Sobgeodnull}
  For all $S$-tangent tensor $F$, we have
  \begin{align*}
    \norm{F}_{L^\infty(S_{u,\ub})} & \les r^{-1/2}\norm{(r\Nd)^{\leq 1} F}_{L^4(S_{u,\ub})},
  \end{align*}
  for all $2$-sphere $S$ of the last cones geodesic null foliation.
\end{lemma}

We now turn to the improvement of the bootstrap bounds~\eqref{est:BAnullgeodconnection}.\\

Using relations~\eqref{eq:relxibtombtzet}, equation~\eqref{eq:Nd3trchib} rewrites as
\begin{align}\label{eq:Nd3trchib'}
  \Nd_3\trchib + \half (\trchib)^2 & = -|\chibh|^2.
\end{align}
Taking the average in~\eqref{eq:Nd3trchib'}, using formula~(\ref{eq:commelbov}) -- which still holds in the last cones geodesic null foliation--, we obtain\footnote{For the computations, we refer the reader to the proof of~(\ref{eq:Nd3trchibtrchibo}), (\ref{eq:Nd3trchibo}), which are done in the more general case of the canonical foliation.}
\begin{align}
  \Nd_3(\trchib-\trchibo) + \trchibo (\trchib-\trchibo) & = -|\chibh|^2 + \overline{|\chibh|^2} -\half(\trchib-\trchibo)^2 -\half \overline{(\trchib-\trchibo)^2},\label{eq:Nd3trchibtrchibo'}\\
  \Nd_3\le(\trchibo+\frac{2}{r}\ri) + \half \trchibo\le(\trchibo+\frac{2}{r}\ri) & = - \overline{|\chibh|^2} + \half \overline{(\trchib-\trchibo)^2}.\label{eq:Nd3trchibo'}  
\end{align}

Integrating~\eqref{eq:Nd3trchibo'}, using that from the vertex limits of Theorem~\ref{thm:vertex}, we have $r\le(\trchibo+\frac{2}{r}\ri) \to 0$ when $r\to 0$, and using the bootstrap bounds~(\ref{est:BAsupnormconngeodcontrol}) we have
\begin{align*}
  r\le|\trchibo+\frac{2}{r}\ri| & \les \int_{u}^\ub \le(r|\chibh|^2 + r\le|\trchib-\trchibo\ri|^2 \ri)\,\d u'\\
                                & \les (D\varep)^2 \int_u^\ub r \ub^{-2}(u')^{-3}\,\d u'\\
                                & \les \varep r \ub^{-2}u^{-3}\qq.
\end{align*}
Thus
\begin{align}
  \label{est:trchibo'}
  \le|\trchibo+\frac{2}{r}\ri| & \les \varep \ub^{-2}u^{-3}\qq.
\end{align}
Using Lemma~\ref{lem:ellgeodnull}, the transport equation~\eqref{eq:Nd3trchibtrchibo'} with $\kappa=2$ and $\ell=1$, that from the vertex limits from Theorem~\ref{thm:vertex}, we have $r^2(\trchib-\trchibo) \to 0$ when $r\to 0$, and the bootstrap bounds~(\ref{est:BAL4normconngeodcontrolprelim}) we have
\begin{align*}
  & r^{-1/2}\norm{r^2(r\Nd)^{\leq 1}\le(\trchib-\trchibo\ri)}_{L^4(S_{u,\ub})} \\
  \les & \; \int_u^\ub \le(r^{-1/2}\norm{r^2(r\Nd)^{\leq 1}|\chibh|^2}_{L^4(S_{u',\ub})} + r^{-1/2}\norm{r^2(r\Nd)^{\leq 1}(\trchib-\trchibo)}_{L^4(S_{u',\ub})} \ri)\, \d u' \\
  \les & (D\varep)^2 \int_u^\ub r^2\ub^{-2}(u')^{-3} \,\d u' \\
  \les & \varep r^2\ub^{-2}u^{-3}\qq.
\end{align*}
Thus, combining this bound with~\eqref{est:trchibo'} and using the Sobolev estimates of Lemma~\ref{lem:Sobgeodnull}, we have
\begin{align}\label{est:trchib'improved}
  \le|\trchib+\frac{2}{r}\ri| + r^{-1/2}\norm{(r\Nd)^{\leq 1}\le(\trchib+\frac{2}{r}\ri)}_{L^4(S_{u,\ub})} & \les \varep \ub^{-2}u^{-3}\qq. 
\end{align}

Equation~\eqref{eq:Nd3ze} together with relations~(\ref{eq:relxibtombtzet}) rewrites
\begin{align}\label{eq:Nd3ze'}
  \Nd_3\ze +\trchib \ze & = -\chibh\cdot\ze -\beb.
\end{align}
Using the estimates of Lemma~\ref{lem:transportgeodnull} with $\kappa =2$ and $\ell =0$, using that from the limits of Theorem~\ref{thm:vertex} one has $r^2\ze\to 0$ when $r\to 0$, using the curvature bounds from Lemma~\ref{lem:curvgeodsource}, the bootstrap bounds~(\ref{est:BAsupnormconngeodcontrol}), we have
\begin{align}\label{est:zeL4'}
  \begin{aligned}
    r^{-1/2}\norm{r^2\ze}_{L^4(S_{u,\ub})} & \les \int_u^\ub \le(r^2|\chibh||\ze| + r^2|\beb|\ri)\,\d u'\\
    & \les \varep \ub^{-2}u^{-3/2}r^2\qq.
  \end{aligned}
\end{align}

Commuting equation~\eqref{eq:Nd3ze'} with $\Divd$, using formula~(\ref{eq:commNd3Nd}), and Bianchi equation~\eqref{eq:Nd3rho}, we have
\begin{align}\label{eq:Nd3Divdzerho'}
  \Nd_3(\Divd\ze-\rho) + \frac{3}{2}\trchib (\Divd\ze-\rho) & = \EEE,
\end{align}
where $\EEE$ is composed of lower order terms and we schematically have
\begin{align*}
  \EEE & := \trchib\rho + \chih\cdot\alb + \ze\cdot\beb \\
       & \quad + \Nd\chibh \cdot\ze + \chibh\cdot\Nd\ze \\
       & \quad + \Nd\trchib\cdot\ze + \trchib\cdot\Nd\ze + \chib\cdot\chib\cdot\ze.
\end{align*}
From equation~\eqref{eq:Nd3Divdzerho'}, the estimates of Lemma~\ref{lem:transportgeodnull} with $\kappa=3$ and $\ell=0$, using that from Theorem~\ref{thm:vertex} one has $r^3\Nd\ze,r^3\rho \to 0$ when $r\to 0$, using the estimates from Lemma~\ref{lem:curvgeodsource}, estimates~\eqref{est:zeL4'} and the bootstrap bounds~\eqref{est:BAsupnormconngeodcontrol}, we have
\begin{align*}
  r^{-1/2}\norm{r^3\le(\Divd\ze-\rho\ri)}_{L^4(S_{u,\ub})} & \les \int_u^\ub r^{-1/2}\norm{r^3\EEE}_{L^4(S_{u',\ub})} \,\d u' \\
                                                           & \les \int_{u}^\ub \varep r^2\ub^{-3}(u')^{-1/2}\,\d u' \\
  & \les \varep r^3\ub^{-7/2}.
\end{align*}
Thus, using the estimates of Lemma~\ref{lem:curvgeodsource}, we have
\begin{align}\label{est:Divdze'}
  r^{-1/2}\norm{\Divd\ze}_{L^4(S_{u,\ub})} & \les \varep \ub^{-3}u^{-1/2}.
\end{align}
From~(\ref{eq:Curlze}), the estimates of Lemma~\ref{lem:curvgeodsource} and the bounds~\eqref{est:BAsupnormconngeodcontrol} we moreover have
\begin{align}\label{est:Curldze'}
  \begin{aligned}
    r^{-1/2}\norm{\Curld\ze}_{L^4(S_{u,\ub})} & \les |\sigma| + |\chibh|r^{-1/2}\norm{\chih}_{L^4(S_{u,\ub})} \\
    & \les \varep \ub^{-3}u^{-1/2}.
  \end{aligned}
\end{align}
From the elliptic estimates of Lemma~\ref{lem:ellgeodnull} together with~\eqref{est:Divdze'} and~\eqref{est:Curldze'}, using moreover the Sobolev estimates of Lemma~\ref{lem:Sobgeodnull}, we infer
\begin{align}\label{est:ze'improved}
  |\ze| + r^{-1/2}\norm{(r\Nd)^{\leq 1}\ze}_{L^4(S_{u,\ub})} & \les \varep \ub^{-2}u^{-3/2}\qq.
\end{align}

Applying the elliptic estimates of Lemma~\ref{lem:ellgeodnull} to equation~\eqref{eq:Divdchibh}, using the estimates of Lemma~\ref{lem:curvgeodsource}, estimates~\eqref{est:trchib'improved} and~\eqref{est:ze'improved}, and the bootstrap bounds~\eqref{est:supnormconngeodcontrol}, we have
\begin{align}\label{est:chibh'improved}
  |\chibh| + r^{-1/2}\norm{(r\Nd)^{\leq 1}\chibh}_{L^4(S_{u,\ub})} & \les \varep \ub^{-2}u^{-3/2},
\end{align}
where we also used the Sobolev embeddings of Lemma~\ref{lem:Sobgeodnull}.\\

Equation~\eqref{eq:Nd3om} together with relations~(\ref{eq:relxibtombtzet}) rewrites
\begin{align}\label{eq:Nd3om'}
  \Nd_3\om & = 3|\ze|^2 +\rho.
\end{align}
Directly integrating equation~\eqref{eq:Nd3om'} and using the estimates from Lemma~\ref{lem:curvgeodsource} and the bootstrap bounds~\eqref{est:BAsupnormconngeodcontrol} gives
\begin{align}
  \label{est:om'Linf}
  |\om| + r^{-1/2}\norm{\om}_{L^4(S_{u,\ub})} & \les \varep r\ub^{-7/2}.
\end{align}
where we used that by the limits of Theorem~\ref{thm:vertex} we have $\om \to 0$ when $r\to 0$.\\

Let define $\om_\sigma$ to be the solution of 
\begin{align}
  \label{eq:Nd3omsigma'}
  \begin{aligned}
  \Nd_3\om_\si & = \sigma, \\
  \om_\si & \to 0 \quad \text{when $r\to 0$.}
  \end{aligned}
\end{align}
Integrating equation~\eqref{eq:Nd3omsigma'}, using the estimates of Lemma~\ref{lem:curvgeodsource}, we have
\begin{align}
  \label{est:omsigma'Linf}
  |\om_\si| & \les \varep r\ub^{-7/2}.
\end{align}

Let define
\begin{align*}
  \io & := \Nd\om + \Nds\om_\si.  
\end{align*}
Commuting equations~\eqref{eq:Nd3om'} and~\eqref{eq:Nd3omsigma'}, using Bianchi equation~\eqref{eq:Nd3be}, we have
\begin{align}\label{eq:Nd3io'}
  \Nd_3(\io-\be) + \half \trchi \io & = \EEE,
\end{align}
where we schematically have
\begin{align*}
  \EEE & := \trchi\be + \om\be + \ze\cdot\al + \ze\cdot\Nd\ze + \chih\cdot\Nd\om + \chih\cdot\Nd\om_\sigma.
\end{align*}
Applying the estimates of Lemma~\ref{lem:transportgeodnull} to~\eqref{eq:Nd3io'}, using that from the limits of Theorem~\ref{thm:vertex} and from~\eqref{eq:Nd3omsigma'} one has $r\Nd\om,r\Nd\om_\si \to 0$ when $r\to 0$, and using the estimates from Lemma~\ref{lem:curvgeodsource}, and~\eqref{est:BAsupnormconngeodcontrol}, we have
\begin{align*}
  r^{-1/2}\norm{r(\io-\be)}_{L^4(S_{u,\ub})} & \les \varep \int_{u}^\ub\ub^{-7/2}\,\d u' \les \varep r\ub^{-7/2}. 
\end{align*}
From this bound, the curvature estimates of Lemma~\ref{lem:curvgeodsource}, the elliptic estimates of Lemma~\ref{lem:ellgeodnull}, and estimate~\eqref{est:om'Linf} we infer
\begin{align}
  \label{est:om'improved}
  |\om| + r^{-1/2}\norm{(r\Nd)^{\leq 1}\om}_{L^4(S_{u,\ub})} & \les \varep r\ub^{-7/2}.
\end{align}

Equation~(\ref{eq:Nd3xiOLD}) together with relations~(\ref{eq:relxibtombtzet}) rewrites
\begin{align}
  \label{eq:Nd3xi'}
  \Nd_3\xi + \half\trchib\xi & = 2\Nd\om-\chibh\cdot\xi + 4\om\xi.
\end{align}
Using the estimates of Lemma~\ref{lem:transportgeodnull} with $\kappa =1$ and $\ell=0$, using that from the limits of Theorem~\ref{thm:vertex} we have $r\xi \to 0$ when $r \to 0$, using estimates~\eqref{est:om'improved} and the bootstrap bounds~\eqref{est:BAsupnormconngeodcontrol} and~(\ref{est:BAL4normconngeodcontrol}), we have
\begin{align}\label{est:xi'improved}
  r^{-1/2}\norm{\xi}_{L^4(S_{u,\ub})} & \les \varep r\ub^{-7/2}.
\end{align}

Equation~\eqref{eq:Nd4trchi} together with relations~\eqref{eq:relxibtombtzet} rewrites
\begin{align}\label{eq:Nd3trchi'}
  \Nd_3\trchi + \half \trchib\trchi & = 2\Divd\ze -\chih\cdot\chibh +2|\ze|^2+2\rho.
\end{align}
Taking the average in~\eqref{eq:Nd3trchi'}, we have
\begin{align}\label{eq:Nd3trchitrchio'}
  \begin{aligned}
    \Nd_3(\trchi-\trchio) + \half \trchibo (\trchi-\trchio) & = 2\Divd\ze + 2(\rho -\rhoo) -\chih\cdot\chibh + \overline{\chih\cdot\chibh} \\
    & \quad +2|\ze|^2 -\half \le(\trchib-\trchibo\ri)\le(\trchi-\trchio\ri) \\
    & \quad - \half \overline{\le(\trchib-\trchibo\ri)(\trchi-\trchio)},
  \end{aligned}
\end{align}
and
\begin{align}\label{eq:Nd3trchio'}
  \Nd_3\le(\trchio-\frac{2}{r}\ri) + \half\trchibo\le(\trchio-\frac{2}{r}\ri) & = 2\rhoo -\overline{\chih\cdot\chibh} +2 \overline{|\ze|^2} +\half \overline{(\trchib-\trchibo)(\trchi-\trchibo)}. 
\end{align}

Integrating~\eqref{eq:Nd3trchio'}, using that from the limits of Theorem~\ref{thm:vertex} we have $r\le(\trchio-\frac{2}{r}\ri) \to 0$ when $r\to 0$, using the estimates from Lemma~\ref{lem:curvgeodsource} and the bounds~\eqref{est:BAsupnormconngeodcontrol} and~(\ref{est:BAL4normconngeodcontrol}), we have
\begin{align}\label{est:trchio'Linf}
  \le|\trchio-\frac{2}{r}\ri| & \les \varep r\ub^{-7/2}.
\end{align}
Applying the estimate of Lemma~\ref{lem:transportgeodnull} with $\kappa=1$ and $\ell=0$ to equation~\eqref{eq:Nd3trchitrchio'}, using the vertex limit from Theorem~\ref{thm:vertex} $r(\trchi-\trchio) \to 0$, using the estimates from Lemma~\ref{lem:curvgeodsource}, estimate~\eqref{est:ze'improved} and the bootstrap bounds~\eqref{est:BAsupnormconngeodcontrol} and~\eqref{est:BAL4normconngeodcontrol}, we have   
\begin{align}\label{est:trchitrchio'L4}
  r^{-1/2}\norm{\le(\trchi-\trchio\ri)}_{L^4(S_{u,\ub})} & \les \varep r\ub^{-7/2},
\end{align}
from which, together with~\eqref{est:trchio'Linf} we also deduce
\begin{align*}
  r^{-1/2}\norm{\le(\trchi-\frac{2}{r}\ri)}_{L^4(S_{u,\ub})} & \les \varep r\ub^{-7/2}.
\end{align*}

Applying the elliptic estimate~\eqref{est:ellspecialL4geodnull} of Lemma~\ref{lem:ellgeodnull} with $F=\trchi-\trchio$ and $G=-\be + \trchi\ze -\chi\cdot\ze$, to equation~\eqref{eq:Divdchih}, using the curvature estimate from Lemma~\ref{lem:curvgeodsource}, estimate~\eqref{est:trchitrchio'L4} and the bootstrap bounds~\eqref{est:BAsupnormconngeodcontrol} and~\eqref{est:BAL4normconngeodcontrol}, we have
\begin{align*}
  \begin{aligned}
  r^{-1/2}\norm{\chih}_{L^4(S_{u,\ub})} & \les r^{-1/2}\norm{\trchi-\trchio}_{L^4(S_{u,\ub})} + \norm{r\be}_{L^4(S_{u,\ub})} + r|\ze|r^{-1/2}\norm{\chi}_{L^4(S_{u,\ub})} \\
  & \les \varep r\ub^{-2}u^{-3/2}\qq.
  \end{aligned}
\end{align*}

Arguing as in Section~\ref{sec:arearadiusestimateCCba}, we have
\begin{align*}
  \le|r-\half(\ub-u)\ri| & \les \int_u^\ub r\le|\trchib+\frac{2}{r}\ri| \, \d u'.
\end{align*}
Using~\eqref{est:trchib'improved}, we thus deduce
\begin{align*}
  \le|r - \half (\ub-u)\ri| & \les \varep r \ub^{-2}u^{-3}\qq^2.
\end{align*}

Arguing as in Section~\ref{sec:mildBAimpast}, we have
\begin{align*}
  & r^{-1/2}\norm{\pr^{\leq 1}r^{-2}\gd_{ab}-(\gd_\SSS)_{ab}} \\
  \les & \, \int_u^\ub\le(r^{-1/2}\norm{r^{-2}\pr^{\leq 1}\le((\trchib-\trchibo)\gd_{ab}\ri)}_{L^4(S_{u',\ub})} + r^{-1/2}\norm{r^{-2}\pr^{\leq 1} \chibh_{ab}}_{L^4(S_{u',\ub})}\ri)\,\d u'.
\end{align*}
Using~\eqref{est:trchib'improved} and~\eqref{est:chibh'improved} and the Sobolev estimates of Lemma~\ref{lem:Sobgeodnull}, this gives
\begin{align*}
  \begin{aligned}
    \le|r^{-2}\gd_{ab}-(\gd_\SSS)_{ab}\ri| & \les \varep \ub^{-1}u^{-3/2}\qq,\\
    \norm{r^{-1/2}\pr_c^{\leq 1}\le(r^{-2}\gd_{ab}-(\gd_{\SSS})_{ab}\ri)}_{L^4(S_{u,\ub})} & \les \varep \ub^{-1}u^{-3/2}\qq.
  \end{aligned}
\end{align*}

Integrating equation~(\ref{eq:Ndyyt}), using that from Theorem~\ref{thm:vertex} one has the limit $\yy \to 0$ when $r\to 0$, and using estimate~\eqref{est:om'improved}, we have
\begin{align*}
  |\yy| & \les \int^\ub_u |\om| \, \d u' \les \varep \ub^{-3/2}.
\end{align*}

This finishes the improvement of the bootstrap bounds~(\ref{est:BAnullgeodconnection}). From a standard continuity argument, this finishes the proof of Lemma~\ref{lem:conngeodcontrol}.

\begin{remark}
  We have actually proved stronger estimates than the estimates of Lemma~\ref{lem:conngeodcontrol}. Since these stronger estimates are not needed in the following, we only record the simpler and weaker version of these estimates in Lemma~\ref{lem:conngeodcontrol}.   
\end{remark}

\section{Control of transition coefficients and comparison of foliations}\label{sec:controltranscoeff}
In this section, we recall that the Bootstrap Assumptions hold, that we have the following improved bounds for the null connection coefficients and area radius of the $2$-spheres $S_{u,\ub}$ (see the definitions of Section~\ref{sec:normnullconnCCba})
\begin{align}\label{est:transcoeffnullconnsourceCCba}
  \begin{aligned}
    \mathfrak{O}^\ast_{\leq 2} + \overline{\mathfrak{O}}^\ast_{\leq 2} & \les \varep, \\
    \le|r-\half\le(\uba-u\ri)\ri| & \les \varep \uba^{-3}r^2u^{-1}.
    \end{aligned}
\end{align}
and (see the definitions of Section~\ref{sec:normsnullconn})
\begin{align}\label{est:transcoeffnullconnsource}
  \begin{aligned}
    \mathfrak{O}^\ext_{\leq 1} + \overline{\mathfrak{O}}_{\leq 1}^\ext & \les \varep,\\
    \le|r-\half(\ub-u)\ri| & \les \varep \ub^{-1}u^{-1}.
  \end{aligned}
\end{align}
and for the interior connection coefficients (see the definitions of Section~\ref{sec:norminteriorconn})
\begin{align}
  \label{est:transcoeffintrsource}
  \OO^\intr_{\leq 3,\ga}[n] + \OO^\intr_{\leq 2}[k] + \mathfrak{O}^\TT_{\leq 2}[\nut] & \les \varep,
\end{align}
that the results of Lemma~\ref{lem:conngeodcontrol} hold 
and that the bottom and conical initial layer are $\varep$-close to Minkowski (see the definitions of Section~\ref{sec:Minkowskilayer}).\\

We prove in Section~\ref{sec:imprBAgeod} that under these assumptions, the following bounds hold in $\MM'\cap\CCba$
\begin{align}\label{est:imprffblageodCCba}
  \le|f',\fb',\log{\la'}\ri| & \les \varep \uba^{-1}{u'}^{-1/2},
\end{align}
and
\begin{align}\label{est:compopticalgeodCCba}
  \le|u-{u'}\ri| & \les \varep \uba^{-3/2}(\uba-u'),
\end{align}
the following bounds hold in $\MM' \cap\MM^\ext$
\begin{align}\label{est:imprffblageod}
  \le|f',\log\la'\ri| & \les \varep {\ub'}^{-1}{u'}^{-1/2}, & \le|\fb'\ri| & \les \varep {u'}^{-3/2},
\end{align}
and
\begin{align}
  \label{est:compopticalgeod}
  \le|u-{u'}\ri| & \les \varep {\ub'}^{-1/2}, & \le|\ub-{\ub'}\ri| & \les \varep {u'}^{-1/2},
\end{align}
and in $\MM'\cap\MM^\intr_\bott$
\begin{align}\label{est:imprTfelbtelt}
  \le|\g\le(\Tf,\half(\elb'+\el')\ri)+1\ri| & \les \varep {\ub'}^{-3/2},
\end{align}
and
\begin{align}\label{est:imptimefunctiongeod}
  \le|t-\half({\ub'}+{u'})\ri| & \les \varep {\ub'}^{-1/2}.
\end{align}

We also prove that the following bounds hold in the conical initial layer $\LLcext$
\begin{align}\label{est:imprffbla}
  \le|\ft\ri| & \les \varep \ub^{-1}, & \le|\fbt,\log\lat\ri| & \les \varep,
\end{align}
and
\begin{align}\label{est:compopticalLLcext}
  \le|u-\ut\ri| & \les \varep, & \le|\ub-\ubt\ri| & \les \varep \ub,
\end{align}
and that the following bounds hold in the bottom initial layer:
\begin{align}\label{est:impprx0elbel}
  \le|\g\le(\half(\elb+\el),\Tf^\bott\ri) +1 \ri| & \les \varep,
\end{align}
and
\begin{align}\label{est:comptimeuub}
  \le|x^0-\half (\ub+u)\ri| & \les \varep, 
\end{align}
in $\LLbext$, and
\begin{align}\label{est:impTfprx0}
  \le|\g(\Tf,\Tf^\bott)+1\ri| & \les \varep,
\end{align}
and
\begin{align}\label{est:comptimeLLbint}
  \le|t-x^0\ri| & \les \varep
\end{align}
in $\LLbint$.

\begin{figure}[h!]
  \centering
  \includegraphics[height=10cm]{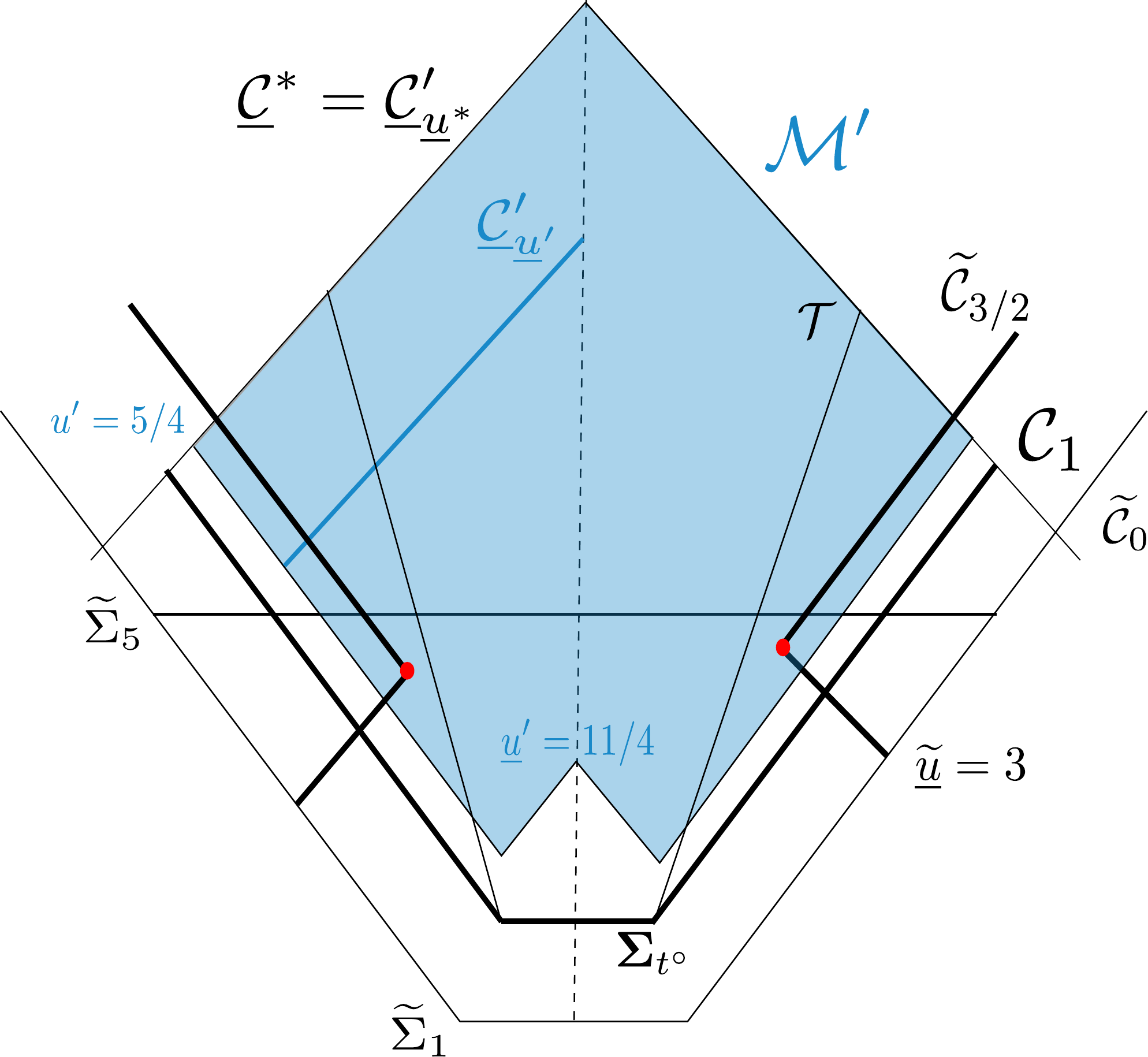}
  \caption{The last cones geodesic foliation and the initial layers.\protect\footnotemark}
  \label{fig:lastconesfoliation}
\end{figure}
\footnotetext{In red the $2$-sphere $\{\ut=3/2,\ubt =3\}$ used in Sections~\ref{sec:conicalcompimpr} and \ref{sec:framecompLLbint}.}

\begin{remark}
  The bounds~\eqref{est:imprffblageodCCba}--\eqref{est:imptimefunctiongeod} improve the Bootstrap Assumptions~\ref{BA:lastconesfoliation} for the comparison to the last cones geodesic foliation. The bounds~\eqref{est:imprffbla} and~\eqref{est:compopticalLLcext} improve the Bootstrap Assumptions~\ref{BA:con} for the comparison to the double null conical initial layer foliation. The bounds~\eqref{est:impprx0elbel}--\eqref{est:comptimeLLbint} improve the Bootstrap Assumptions~\ref{BA:bottom} for the comparison to the bottom initial layer coordinates.
\end{remark}

\subsection{Last cones geodesic foliation comparisons}\label{sec:imprBAgeod}
From the definition of the last cone geodesic foliation, we have the following limits at the tip of the axis $\o(\uba)$
\begin{align}\label{est:ffblaouba}
  |f'|,|\fb'|,\log{\la'} \to 0.
\end{align}
Moreover, since $\CCba = \CCb'_{{\ub'}=\uba}$, we have (see the relations of Lemma~\ref{lem:transframe})\footnote{We recall that for the relations of Lemma~\ref{lem:transframe} as well as for all the relations of Section~\ref{sec:defchgframe}, $f,\fb,\la$ are respectively denoted by $f',\fb',\la'$ in the present section. See also the definition of these coefficients in the Bootstrap Assumptions~\ref{BA:lastconesfoliation}.}
\begin{align}\label{eq:fbtCCba}
  \fb' = 0
\end{align}
on $\CCba$.\\

Using the limits~\eqref{est:ffblaouba} and by integration along $\elb'$, we have
\begin{align}\label{est:r'Nd'3f'v1}
  \begin{aligned}
    \norm{{r'}^{-1/2}(f')^\ddg}_{L^\infty_{{u'}}L^4(S'_{u',\uba})} & \les \int_{5/4}^\uba\norm{{r'}^{-1/2}\Nd'_3(f')^\ddg}_{L^4(S'_{u',\uba})} \, \d u'.
    \end{aligned}
\end{align}
Using equation~\eqref{eq:Nd3f}, we have
\begin{align}\label{est:r'Nd'3f'v2}
  \begin{aligned}
    \norm{{r'}^{-1/2}\Nd'_3(f')^\ddg}_{L^4(S'_{u',\uba})} & \les \norm{{r'}^{-1/2}\eta^\ddg}_{L^4(S'_{u',\uba})} + \norm{{r'}^{-1/2}\eta'}_{L^4(S'_{u',\uba})} \\
    & \quad + \norm{{r'}^{-1/2}\Err\le(\eta,\eta'\ri)}_{L^4(S'_{u',\uba})}.
  \end{aligned}
\end{align}
Using the sup-norm estimates~(\ref{est:transcoeffnullconnsource}) for $\eta=\ze$, we have
\begin{align}\label{est:r'Nd'3f'v2bis}
  \begin{aligned}
    \norm{{r'}^{-1/2}\eta^\ddg}_{L^4(S'_{u',\uba})} \les \norm{(|\ze|)}_{L^\infty(S'_{u',\uba})} \les \varep \norm{\max\le(r^{-1/2},u^{-1/2}\ri)}_{L^\infty(S'_{u',\uba})}\uba^{-2},
  \end{aligned}
\end{align}
where $r$ denotes the area radius of $S_{u,\ub}$. Using~(\ref{est:BAarearadiusestimateCCba}) and the Bootstrap Assumptions~\ref{BA:lastconesfoliation}, we have
\begin{align}\label{est:r'Nd'3f'v3}
  \begin{aligned}
    \norm{r^{-1/2}}_{L^\infty(S'_{u',\uba})} & \les \norm{(\uba-u)^{-1/2}}_{L^\infty(S'_{u',\uba})} \\
    & \les \norm{\le((\uba-u') - |u'-u|\ri)^{-1/2}}_{L^\infty(S'_{u',\uba})} \\
    & \les \le((\uba-u')(1-D\varep)\ri)^{-1/2},
    \end{aligned}
\end{align}
and from the Bootstrap Assumptions~\ref{BA:lastconesfoliation}, we obtain similarly
\begin{align}\label{est:r'Nd'3f'v4}
  \norm{u^{-1/2}}_{L^\infty(S'_{u',\uba})} & \les (u')^{-1/2}.
\end{align}
Combining~\eqref{est:r'Nd'3f'v2bis},~\eqref{est:r'Nd'3f'v3} and~\eqref{est:r'Nd'3f'v4} we have
\begin{align}
  \label{est:r'Nd'3f'v5}
  \begin{aligned}
    \norm{{r'}^{-1/2}\eta^\ddg}_{L^4(S'_{u',\uba})} \les \varep \norm{\max\le(r^{-1/2},u^{-1/2}\ri)}_{L^\infty(S'_{u',\uba})}\uba^{-2} \les \varep \uba^{-2}(\qq')^{-1/2},
  \end{aligned}
\end{align}
where $\qq':=\max(\ub'-u,u')$. Using the the sup-norm estimates~\eqref{est:transcoeffnullconnsource} for $\omb$, the Bootstrap Assumptions~\ref{BA:lastconesfoliation}, relation~(\ref{eq:fbtCCba}) and the estimates~\eqref{est:r'Nd'3f'v3},~\eqref{est:r'Nd'3f'v4}, we have\footnote{We do not treat the lower order error terms of $\Err(\eta,\eta')$. See also Proposition~\ref{prop:transconn}.}
\begin{align}
  \label{est:r'Nd'3f'v2bisbis}
  \begin{aligned}
    \norm{{r'}^{-1/2}\Err\le(\eta,\eta'\ri)}_{L^4(S'_{u',\uba})} \les \norm{|f'||\omb|}_{L^\infty(S'_{u',\uba})} \les D\varep^2\uba^{-2}(u')^{-3/2}(\qq')^{-1/2}.
    \end{aligned}
\end{align}
From~\eqref{est:r'Nd'3f'v1},~\eqref{est:r'Nd'3f'v2},~\eqref{est:r'Nd'3f'v5} and~\eqref{est:r'Nd'3f'v2bisbis}, and using the $L^4$ estimates from Lemma~\ref{lem:conngeodcontrol} for $\eta'=\ze'$, we have
\begin{align}\label{est:ftCCba}
  \begin{aligned}
    \norm{{r'}^{-1/2}f'}_{L^\infty_{{u'}}L^4(S'_{u',\uba})} & \les \int_{5/4}^{\uba}\varep (\qq')^{-1/2}{\uba}^{-2} \,\d u'\\
    & \les \varep\uba^{-3/2}.
  \end{aligned}
\end{align}

Using equation~\eqref{eq:Nd3la} for $\Nd_3'\log\la'$, estimates~\eqref{est:transcoeffnullconnsourceCCba} and the Bootstrap Assumptions~\ref{BA:lastconesfoliation}, we have
\begin{align}\label{est:Nd'3'la'}
  \begin{aligned}
    \norm{{r'}^{-1/2}\Nd'_3\log\la'}_{L^4(S'_{u',\uba})} & \les \norm{(\la')^{-1}\omb}_{L^\infty(S'_{u',\uba})} \les \varep \uba^{-1}(\qq')^{-1/2}(u')^{-1},
  \end{aligned}
\end{align}
where we used that from the relations of Lemmas~\ref{lem:relgeodnull} and~\ref{lem:relCCba} and relation~\eqref{eq:fbtCCba} we have $\Err(\omb,\omb') = 0$ on $\MM'\cap\CCba$. Integrating~\eqref{est:Nd'3'la'}, using~\eqref{est:ffblaouba} and~\eqref{est:Nd'3'la'}, we obtain on $\MM'\cap\CCba$
\begin{align}\label{est:laCCba}
  \begin{aligned}
    \norm{{r'}^{-1/2}\log\la'}_{L^4(S'_{u',\uba})} & \les \int_{u'}^\uba \norm{{r'}^{-1/2}\Nd'_3\log\la'}_{L^4(S'_{v',\uba})} \,\d v' \\
    & \les \varep \uba^{-1}\int_{u'}^\uba (\qq')^{-1/2}(v')^{-1}\,\d v' \\
    & \les \varep {\uba}^{-1}{u'}^{-1/2}.
  \end{aligned}
\end{align}
Using~\eqref{eq:elu}, the relations of Lemma~\ref{lem:transframe} and~\eqref{est:laCCba}, we have
\begin{align*}
  \norm{{r'}^{-1/2}\elb'(u'-u)}_{L^4(S'_{u',\uba})} & = \norm{{r'}^{-1/2}\le(2 - 2(\la')^{-1}\ri)}_{L^4(S'_{u',\uba})} \les \varep \uba^{-1} (u')^{-1/2}.
\end{align*}
Thus by integration in $\elb'$, using that $u'=\uba=u$ at $\o(\uba)$, we deduce
\begin{align}\label{est:L4u'u}
  \norm{{r'}^{-1/2}\le(u'-u\ri)}_{L^4(S'_{u',\uba})} & \les \varep \uba^{-3/2}(\uba-u').
\end{align}

Using equation~\eqref{eq:Ndf}, we have
\begin{align}\label{est:Nd'f'L4v1}
  \begin{aligned}
    \norm{{r'}^{-1/2}r'\Nd'(f')^\ddg}_{L^4(S'_{u',\uba})} & \les \norm{{r'}^{-1/2}r'\le(\chi'-\chi^{\ddg}\ri)}_{L^\infty(S'_{u',\uba})} \\
    & \quad + \norm{{r'}^{-1/2}r'\le((\la')^{-1}-1\ri)\chi'}_{L^4(S'_{u',\uba})} \\
    & \quad + \norm{{r'}^{-1/2}r'\Err\le(\chi,\chi'\ri)}_{L^4(S'_{u',\uba})}.
  \end{aligned}
\end{align}
We have
\begin{align}\label{est:Nd'f'L4v2}
  \begin{aligned}
    \le(\chi'-\chi^\ddg\ri)_{ab} & = \frac{1}{r'}\gd'_{ab} + \half\le(\trchi'-\frac{2}{r'}\ri)\gd'_{ab} + \chih'_{ab} - \le(\frac{1}{r}\gd'_{ab} +\half \le(\trchi-\frac{2}{r}\ri)\gd'_{ab} + \chih_{ab}\ri) \\
    & = \le(\frac{2}{\ub'-u'} - \frac{2}{\ub-u}\ri)\gd'_{ab}  + \le(\frac{1}{r'} - \frac{2}{\ub'-u'}\ri)\gd'_{ab} + \le(\frac{2}{\ub-u}-\frac{1}{r}\ri)\gd'_{ab} \\
    & \quad + \half\le(\trchi'-\frac{2}{r'}\ri)\gd'_{ab} + \chih'_{ab} - \half \le(\trchi-\frac{2}{r}\ri)\gd'_{ab} - \chih_{ab}
  \end{aligned}
\end{align}

Using the improved $L^4$ estimate~\eqref{est:L4u'u} for $u-u'$, we have
\begin{align}\label{est:Nd'f'L4v2bis1}
  \begin{aligned}
    \norm{{r'}^{-1/2}r'\le(\frac{2}{\uba-u'}-\frac{2}{\uba-u}\ri)}_{L^4(S'_{u',\uba})} & \les \norm{{r'}^{-1/2}{r'}^{-1}\le(u-u'\ri)}_{L^4(S'_{u',\uba})} \\
    & \les \varep \uba^{-3/2}.
  \end{aligned}
\end{align}
Using the sup-norm estimates for the respective area radii~\eqref{est:arearadiuslastconesgeod} and~\eqref{est:transcoeffnullconnsourceCCba}, we have
\begin{align}\label{est:Nd'f'L4v2bis2}
  \begin{aligned}
    \norm{{r'}^{-1/2}r'\le(\frac{1}{r'}-\frac{2}{\uba-u'}\ri)}_{L^4(S'_{u',\uba})} & \les \varep \uba^{-2}(u')^{-1/2},\\
    \norm{{r'}^{-1/2}r'\le(\frac{1}{r}-\frac{2}{\uba-u}\ri)}_{L^4(S'_{u',\uba})} & \les \varep \uba^{-3}r'(u')^{-1}
  \end{aligned}
\end{align}
Thus, using equation~\eqref{est:Nd'f'L4v2}, the estimates~\eqref{est:Nd'f'L4v2bis1} and~\eqref{est:Nd'f'L4v2bis2} and the sup-norm and $L^4$ estimates for $\chi$ and $\chi'$ from~\eqref{est:transcoeffnullconnsourceCCba} and Lemma~\ref{lem:conngeodcontrol}, we deduce
\begin{align}\label{est:L4chi'chiddg}
  \begin{aligned}
    \norm{{r'}^{-1/2}r'\le(\chi'-\chi^{\ddg}\ri)}_{L^4(S_{u',\uba})} & \les \varep \uba^{-1}(u')^{-1/2}.
  \end{aligned}
\end{align}
Using the $L^4$ estimate~\eqref{est:laCCba} for $\la'$, the $L^4$ estimates $\chi'$ from Lemma~\ref{lem:conngeodcontrol} and the sup-norm estimates for $\la'$ from the Bootstrap Assumptions~\ref{BA:lastconesfoliation}, we have
\begin{align}\label{est:la'1chi'L4}
  \begin{aligned}
    \norm{{r'}^{-1/2}r'((\la')^{-1}-1)\chi'}_{L^4(S'_{u',\uba})} & \les \norm{{r'}^{-1/2}((\la')^{-1}-1)}_{L^4(S'_{u',\uba})} \\
    & \quad + \norm{{r'}^{-1/2}r'((\la')^{-1}-1)\le(\trchi'-\frac{2}{r'}\ri)}_{L^4(S'_{u',\uba})} \\
    & \quad + \norm{{r'}^{-1/2}r'((\la')^{-1}-1)\chih'}_{L^4(S'_{u',\uba})} \\
    & \les \varep \uba^{-1}(u')^{-1/2} \\
    & \quad + r'\norm{\la'-1}_{L^\infty(S'_{u',\uba})}\norm{{r'}^{-1/2}\le(\trchi'-\frac{2}{r'}\ri)}_{L^4(S'_{u',\uba})} \\
    & \quad + r'\norm{\la'-1}_{L^\infty(S'_{u',\uba})}\norm{{r'}^{-1/2}\chih'}_{L^4(S'_{u',\uba})} \\
    & \les \varep \uba^{-1}(u')^{-1/2}.
  \end{aligned}
\end{align}
From the expression of $\Err(\chi,\chi')$ in Proposition~\ref{prop:transconn}, using the Bootstrap Assumptions~\ref{BA:lastconesfoliation} and the estimates~\eqref{est:transcoeffnullconnsourceCCba} and from Lemma~\ref{lem:conngeodcontrol}, we check that
\begin{align}
  \label{est:Errchichi'L4}
  \norm{{r'}^{-1/2}r'\Err\le(\chi,\chi'\ri)}_{L^4(S'_{u',\uba})} & \les (D\varep)^2\uba^{-1}(u')^{-1/2}.
\end{align}
Thus, plugging~\eqref{est:L4chi'chiddg},~\eqref{est:la'1chi'L4} and~\eqref{est:Errchichi'L4} into~\eqref{est:Nd'f'L4v1}, we obtain
\begin{align*}
  \norm{{r'}^{-1/2}r'\Nd'(f')^\ddg}_{L^4(S'_{u',\uba})} & \les \varep \uba^{-1}(u')^{-1/2},
\end{align*}
from which, using also~\eqref{est:ftCCba}, together with the Sobolev estimates from Lemma~\ref{lem:Sobgeodnull}, we deduce on $\MM'\cap\CCba$
\begin{align*}
  \le|f'\ri| & \les \varep \uba^{-1}(u')^{-1/2}.
\end{align*}
Using equation~\eqref{eq:Ndlogla}, the improved $L^4$ estimate~\eqref{est:ftCCba}, the sup-norm and $L^4$ estimates~\eqref{est:transcoeffnullconnsourceCCba} and of Lemma~\ref{lem:conngeodcontrol}, and the Bootstrap Assumptions~\ref{BA:lastconesfoliation}, we check that
\begin{align*}
  \norm{{r'}^{-1/2}r'\Nd'\log\la'}_{L^4(S'_{u',\uba})} & \les \varep \uba^{-1}(u')^{-1/2},
\end{align*}
which combined with~\eqref{est:laCCba} and the Sobolev estimates of Lemma~\ref{lem:Sobgeodnull} gives on $\MM'\cap\CCba$
\begin{align*}
  \le|\log\la'\ri| & \les \varep \uba^{-1}(u')^{-1/2}.
\end{align*}
By integration along $\elb$, arguing as in the proof of~\eqref{est:L4u'u}, we deduce from that estimate that on $\MM'\cap\CCba$
\begin{align}\label{est:uutCCba}
  \le|u-u'\ri| & \les \varep \uba^{-3/2}(\uba-u').
\end{align}
This finishes the proof of estimates~(\ref{est:imprffblageodCCba}) and (\ref{est:compopticalgeodCCba}) in $\MM'\cap\CCba$.\\

We now turn to the proof of the estimate in $\MM'\cap\MM^\ext$. We first define the following extension of $\MM'\cap\MM^\ext$ and its timelike boundary
\begin{align*}
  \widetilde{\MM^\ext} & := \le\{u'\leq (1+\tau) \ub'/2 \ri\},\\
  \widetilde{{\TT}} & := \le\{u'= (1+\tau)\ub'/2\ri\}.
\end{align*}
\begin{remark}\label{rem:slightextensiongeodnull}
  The geodesic null foliation $(S_{u,\ub})$ can be extended past its original boundary $\TT$ to cover the whole region $\widetilde{\MM^\ext}$ using a standard continuity argument and the same estimates as in Section~\ref{sec:connest}, using the estimates for the curvature in the interior region obtained in Sections~\ref{sec:planehypcurvest} and~\ref{sec:remainingcurvestfinal}. In the following, we will thus assume that the null geodesic foliation covers $\widetilde{\MM^\ext}$ and satisfies consistent extensions of the estimates~\eqref{est:transcoeffnullconnsource}. We further assume that consistent extensions of the bounds from the Bootstrap Assumption~\ref{BA:lastconesfoliation} for the transition coefficients between the geodesic null and the last cones geodesic foliation hold in $\widetilde{\MM^\ext}$.
\end{remark}

By integration along $\el'$, using equation~\eqref{eq:Nd4f}, estimate~\eqref{est:ftCCba}, the estimates~\eqref{est:transcoeffnullconnsource}, the estimates from Lemma~\ref{lem:conngeodcontrol} and the Bootstrap Assumptions~\ref{BA:lastconesfoliation}, we have for all $2$-sphere $S'_{u',\ub'} \subset \widetilde{\MM^\ext}$
\begin{align*}
  \begin{aligned}
    \norm{{r'}^{-1/2}r'f'}_{L^4(S'_{u',\ub'})} & \les \norm{{r'}^{-1/2}r'f'}_{L^4(S'_{u',\uba})} + \int_{\ub'}^{\uba}\norm{{r'}^{-1/2}r'\le(\Nd'_4(f')^\ddg+\half\trchi (f')^\ddg\ri)}_{L^4(S'_{u',\vb'})}\,\d\vb' \\
    & \les \varep \uba^{-1/2} + \int^{\uba}_{\ub'}\le(\norm{{r'}^{-1/2}r'\xi'}_{L^4(S'_{u',\vb'})} + \norm{{r'}^{-1/2}r'\Err\le(\xi,\xi'\ri)}_{L^4(S'_{u',\vb'})}\ri)\,\d\vb' \\
    & \les \varep \uba^{-1/2} + \int_{\ub'}^\uba\le(\varep (\vb')^{-3/2} + (D\varep)^2(\vb')^{-2}\ri)\,\d\vb'\\
    & \les \varep (\ub')^{-1/2}.
  \end{aligned}
\end{align*}
Thus, we deduce that for all $2$-sphere $S'_{u',\ub'} \subset \widetilde{\MM^\ext}$
\begin{align}\label{est:ftMMext}
  \norm{{r'}^{-1/2}f'}_{L^4(S'_{u',\ub'})} & \les \varep {\ub'}^{-3/2}.
\end{align}


Arguing similarly, by integration of~\eqref{eq:Nd4fb} and~\eqref{eq:Nd4la}, using~\eqref{est:ftMMext},~\eqref{eq:fbtCCba},~\eqref{est:laCCba},~\eqref{est:transcoeffnullconnsourceCCba}, the estimates of Lemma~\ref{lem:conngeodcontrol} and the Bootstrap Assumptions~\ref{BA:lastconesfoliation}, we have
\begin{align}
  \norm{{r'}^{-1/2}\fb'}_{L^4(S'_{u',\ub'})} & \les \varep {\ub'}^{-1}{u'}^{-1/2},\label{est:fbtMMext}\\
  \norm{{r'}^{-1/2}\log{\la'}}_{L^4(S'_{u',\ub'})} & \les \varep{\ub'}^{-1}{u'}^{-1/2},\label{est:laMMext}
\end{align}
for all $2$-sphere $S'_{u',\ub'}\subset \widetilde{\MM^\ext}$.\\

From relations~\eqref{eq:elu}, \eqref{eq:elut} and the relations of Lemma~\ref{lem:transframe}, we have
\begin{align*}
  \el'({u'}-u) & = \yy' - \el'(u) = \yy' -\half \la' |f'|^2.
\end{align*}
Integrating the above equation, using the sup-norm estimate~\eqref{est:uutCCba} for $u-u'$ on $\MM'\cap\CCba$, the sup-norm estimate for $\yy'$ from Lemma~\ref{lem:conngeodcontrol} and the sup-norm estimates for $\la',f'$ from the Bootstrap Assumptions~\ref{BA:lastconesfoliation} we deduce that in $\widetilde{\MM^\ext}$
\begin{align}\label{est:supnormu'uMMext}
  |{u'}-u| & \les \varep {\ub'}^{-1/2}.
\end{align}

Let $Z'$ be defined by 
\begin{align*}
  Z' & := \quar(1+\cc-\yy')\elb' + \half \el'.
\end{align*}
Using relations~\eqref{eq:elu},~\eqref{eq:elut} and the relations of Lemma~\ref{lem:transframe}, we have that $Z'$ is tangent to $\widetilde{\TT}$, and 
\begin{align*}
  Z'(\ub-\ub') & = (\la'-1) + \la'|f'|^2\yy + \quar(\la')^{-1}\le(1+\cc-\yy'\ri)\le(\yy+\half f'\cdot\fb'\yy + \frac{1}{16}|f'|^2|\fb'|^2\yy+\half |\fb'|^2\ri). 
\end{align*}
Thus, using the $L^4$ bound~\eqref{est:laMMext} for $\la'$, the sup-norm estimates~\eqref{est:transcoeffnullconnsource} the sup-norm estimates of Lemma~\ref{lem:conngeodcontrol} and the sup-norm estimates from the Bootstrap Assumptions~\ref{BA:lastconesfoliation} we obtain for all $2$-sphere $S'_{u',\ub'} \subset \widetilde{\TT}$
\begin{align*}
  \norm{{r'}^{-1/2}Z'({\ub'}-\ub)}_{L^4(S'_{u',\ub'})} & \les \varep \ub'^{-3/2}.
\end{align*}
Thus, we infer by integration from $\widetilde{\TT}\cap\CCba$ along $Z$, using that $\ub'=\uba=\ub$ on $\CCba$
\begin{align}\label{est:ubtubTTL4}
  \norm{{r'}^{-1/2}({\ub'}-\ub)}_{L^4(S'_{u',\ub'})} & \les \varep {\ub'}^{-1/2},
\end{align}
for all $2$-sphere $S'_{u',\ub'} \subset \widetilde{\TT}$. Using~\eqref{eq:elu}, \eqref{eq:elut} and the relations of Lemma~\ref{lem:transframe}, we have
\begin{align*}
  \elb'(\ub-\ub') & =  \la'^{-1}\le(\le(1+\half f'\cdot\fb' + \frac{1}{16}|f'|^2|\fb'|^2\ri)\yy + \half |\fb'|^2\ri),
\end{align*}
from which, using estimates~\eqref{est:transcoeffnullconnsource} for $\yy$ and the sup-norm estimates of the Bootstrap Assumptions~\ref{BA:lastconesfoliation}, we have the sup-norm estimate
\begin{align*}
  \le|\elb'(\ub-\ub')\ri| & \les \varep u'^{-3/2}.
\end{align*}
Thus integrating the above estimate along $\elb'$ onto $\widetilde{\MM^\ext}$, using the $L^4$ estimate on $\widetilde{\TT}$, we obtain
\begin{align}\label{est:compubtubMMext}
  \norm{r'^{-1/2}(\ub-\ub')}_{L^4(S'_{u',\ub'})} & \les \varep u'^{-1/2},
\end{align}
for all $2$-sphere $S'_{u',\ub'} \subset \widetilde{\MM^\ext}$.\\


Using equation~\eqref{eq:Ndf}, we have
\begin{align}\label{est:Nd'f'L4v1}
  \begin{aligned}
    \norm{{r'}^{-1/2}r'\Nd'(f')^\ddg}_{L^4(S'_{u',\ub'})} & \les \norm{{r'}^{-1/2}r'\le(\chi'-\chi^{\ddg}\ri)}_{L^\infty(S'_{u',\ub'})} \\
    & \quad + \norm{{r'}^{-1/2}r'\le((\la')^{-1}-1\ri)\chi'}_{L^4(S'_{u',\ub'})} \\
    & \quad + \norm{{r'}^{-1/2}r'\Err\le(\chi,\chi'\ri)}_{L^4(S'_{u',\ub'})}.
  \end{aligned}
\end{align}
Arguing as for the estimates on $\CCba$, using the expression for $\Err(\chi,\chi')$ from Proposition~\ref{prop:transconn}, the $L^4$ estimate~\eqref{est:laMMext} for $\la'$, estimates~\eqref{est:transcoeffnullconnsource}, the estimates of Lemma~\ref{lem:conngeodcontrol} and the Bootstrap Assumptions~\ref{BA:lastconesfoliation}, we have
\begin{align}\label{est:Ndf'MMext1}
  \begin{aligned}
  \norm{{r'}^{-1/2}r'\le((\la')^{-1}-1\ri)\chi'}_{L^4(S'_{u',\ub'})} & \les \varep \ub'^{-1}u'^{-1/2},\\
  \norm{{r'}^{-1/2}r'\Err\le(\chi,\chi'\ri)}_{L^4(S'_{u',\ub'})} & \les (D\varep)^2\ub'^{-1}u'^{-1/2}.
  \end{aligned}
\end{align}
Using estimates~\eqref{est:supnormu'uMMext} and~\eqref{est:compubtubMMext}, we have
\begin{align}\label{est:Ndf'MMext2}
  \begin{aligned}
    \norm{r'^{-1/2}r'\le(\frac{2}{\ub'-u'}-\frac{2}{\ub-u}\ri)}_{L^4(S'_{u',\ub'})} & \les \norm{r'^{-1/2}r'^{-1}(\ub'-\ub)}_{L^4(S'_{u',\ub'})} \\
    & \quad + \norm{r'^{-1/2}r'^{-1}(u'-u)}_{L^4(S'_{u',\ub'})} \\
    & \les \varep \ub'^{-1}u'^{-1/2}.
  \end{aligned}
\end{align}
Using estimates~\eqref{est:arearadiuslastconesgeod} and~\eqref{est:transcoeffnullconnsource}, we have
\begin{align}\label{est:Ndf'MMext3}
  \begin{aligned}
    \norm{r'^{-1/2}r'\le(\frac{1}{r'}-\frac{2}{\ub'-u'}\ri)}_{L^4(S'_{u',\ub'})} & \les \varep \ub'^{-2}u'^{-1/2},\\
    \norm{r'^{-1/2}r'\le(\frac{1}{r}-\frac{2}{\ub-u}\ri)}_{L^4(S'_{u',\ub'})} & \les \varep \ub'^{-2}u'^{-1}.
  \end{aligned}
\end{align}
Thus, using formula~\eqref{est:Nd'f'L4v2}, the estimates~\eqref{est:Ndf'MMext1}~\eqref{est:Ndf'MMext2}~\eqref{est:Ndf'MMext3}, the estimates~\eqref{est:transcoeffnullconnsource} for $\chi$, the estimates of Lemma~\ref{lem:conngeodcontrol} for $\chi'$, we obtain
\begin{align}\label{est:Ndf'MMextfinal}
  \begin{aligned}
    \norm{r'^{-1/2}r'\Nd'(f')^\ddg}_{L^4(S'_{u',\ub'})} & \les \varep \ub'^{-1}u'^{-1/2},
  \end{aligned}
\end{align}
for all $2$-sphere $S'_{u',\ub'}\subset \widetilde{\MM^\ext}$. Using the Sobolev estimates from Lemma~\ref{lem:Sobgeodnull} and estimates~\eqref{est:ftMMext},~\eqref{est:Ndf'MMextfinal}, we obtain
\begin{align}\label{est:supnormf'ext}
  |f'| & \les \varep \ub'^{-1}u'^{-1/2}
\end{align}
in $\widetilde{\MM^\ext}$.\\

Arguing along the same lines, using instead equation~\eqref{eq:Ndfb} and replacing the estimates for $\chi,\chi'$ by the estimates for $\chib,\chib'$ from~\eqref{est:transcoeffnullconnsource} and from Lemma~\ref{lem:conngeodcontrol}, we obtain
\begin{align*}
  \norm{r'^{-1/2}r'\Nd'(\fb')^\ddg}_{L^4(S'_{u',\ub'})} & \les \varep u'^{-3/2}.
\end{align*}
Thus, by the Sobolev estimates from Lemma~\ref{lem:Sobgeodnull} and using additionally~\eqref{est:fbtMMext}, we obtain
\begin{align}\label{est:supnormfb'ext}
  |\fb'| & \les \varep u'^{-3/2}.
\end{align}

Taking the $L^4(S'_{u',\ub'})$ in equation~\eqref{eq:Ndlogla} using the estimates~\eqref{est:ftMMext},~\eqref{est:fbtMMext}, the estimates~\eqref{est:transcoeffnullconnsource}, the estimates of Lemma~\ref{lem:conngeodcontrol} and the Bootstrap Assumptions~\ref{BA:lastconesfoliation}, we have
\begin{align*}
  \norm{r'^{-1/2}r'\Nd'\log\la'}_{L^4(S'_{u',\ub'})} & \les \varep \ub'^{-1}u'^{-1/2},
\end{align*}
from which we deduce by the Sobolev estimates from Lemma~\ref{lem:Sobgeodnull} and using additionally~\eqref{est:laMMext}
\begin{align}\label{est:supnormla'ext}
  |\log\la'| & \les \varep \ub'^{-1}u'^{-1/2}.
\end{align}
This finishes the proof of estimates~(\ref{est:imprffblageod}).\\

Arguing along the same lines as previously, using the sup-norm estimate~\eqref{est:supnormla'ext} instead of the $L^4$ estimate~\eqref{est:laMMext}, we obtain
\begin{align*}
  \le|\ub-\ub'\ri| & \les \varep u'^{-1/2}.
\end{align*}
This estimate together with ther previously obtained~\eqref{est:supnormu'uMMext} concludes the proof of~\eqref{est:compopticalgeod}.\\


Let define
\begin{align*}
  \phi & := \g\le(\Tf,\half(\elb'+\el')\ri).
\end{align*}
Using the sup-norm estimates~\eqref{est:supnormf'ext},~\eqref{est:supnormfb'ext},~\eqref{est:supnormla'ext} for $f',\fb',\la'$, relations~(\ref{eq:defnutnur}) and estimates~\eqref{est:transcoeffintrsource} for $\nut$, we obtain on $\TT$
\begin{align}
  \le|\phi+1\ri| & \les \varep {\ub'}^{-3/2}.
\end{align}
Differentiating $\phi$ by $\elb'$ and integrating up to the central axis, using the bounds~\eqref{est:transcoeffintrsource} and the estimates from Lemma~\ref{lem:conngeodcontrol}, one obtains $L^\infty L^4(S')$ estimates in $\MM^\intr_\bott$ for $\phi$. Differentiating by $\Nd'$, one further obtains $L^\infty L^4(S')$ estimates for $\Nd'\phi$, which together with Sobolev embeddings on the $2$-spheres $S'$ gives
\begin{align}\label{est:imprTfMMintrbott}
  \le|\phi+1\ri| & \les \varep {\ub'}^{-3/2}
\end{align}
in $\MM'\cap\MM^\intr_\bott$. This finishes the proof of estimate~(\ref{est:imprTfelbtelt}).\\


We now turn to the proof of~(\ref{est:imptimefunctiongeod}). From the definition of $t$ on $\TT$ and the bounds obtained for ${u'}-u$ and ${\ub'}-\ub$ in $\MM'\cap\MM^\ext$ and on $\MM'\cap\TT$, we have on $\MM'\cap\TT$
\begin{align}\label{est:comptubtutTT}
  \le|t-\half({\ub'}+{u'})\ri| & \les \varep {\ub'}^{-1/2},
\end{align}
On $\MM^\intr_\bott$, using relations~(\ref{eq:relTnt}) and (\ref{eq:Yteltelbt}), we have
\begin{align*}
  -\D\le(t-\half({\ub'}+{u'})\ri) & = n^{-1}\Tf-\half(\elb'+Y') = n^{-1}\Tf-\half(\elb'+\el') -\quar\yy'\elb'. 
\end{align*}
Using the bounds $|\yy'| \les \varep {\ub'}^{-3/2}$, (\ref{est:transcoeffnullconnsource}), and the bounds~(\ref{est:imprTfMMintrbott}) obtained for $\phi$, we deduce
\begin{align*}
  \le|\D\le(t-\half({\ub'}+{u'})\ri)\ri| & \les \varep {\ub'}^{-3/2}, 
\end{align*}
which by integration from the boundary $\TT$ onto $\MM^\intr_\bott$, using estimate~\eqref{est:comptubtutTT} on $\TT$ gives~\eqref{est:imptimefunctiongeod} as desired.

\subsection{Conical initial layer comparisons}\label{sec:conicalcompimpr}
We first have the following lemma, which is a consequence of the $\varep$-closeness of the bottom initial layer to Minkowski (see the definitions from Section~\ref{sec:Minkowskilayer}) and the (local) definition of the last cones geodesic foliation. 
\begin{lemma}\label{lem:compgeodinitlayer}
  In $\MM'\cap\LLb$, the following comparison bounds hold between the last cones geodesic foliation and the bottom initial layer spacetime coordinates
  \begin{align}\label{est:compgeodbottomx0}
    \begin{aligned}
      \le|x^0 - \half ({\ub'}+{u'})\ri| & \les \varep,\\
      \le|\sqrt{\sum_{i=1}^3(x^i)^2} - \half({\ub'}-{u'})\ri| & \les \varep.
    \end{aligned}
  \end{align}
  Moreover, we have in $\MM'\cap\LLb$
  \begin{align}\label{est:compgeodbottomTf}
    \begin{aligned}
      \le|\g\le(\Tf^\bott,\half(\elb'+\el')\ri)+1\ri| & \les \varep,
    \end{aligned}
  \end{align}
  and in $\MM'\cap\LLbext$
  \begin{align}\label{est:compgeodbottomNf}
    \le|\g\le(\Nf^\bott,\half(\el'-\elb')\ri)-1\ri| & \les \varep.
  \end{align}
\end{lemma}
\begin{proof}
  The centre $\o(1)$ of $\Sit_1$ is chosen such that $x^i(\o(1)) = 0$ (see also Remark~\ref{rem:centreSit1}). Moreover, from the assumption that $\doto(1)$ is the unit normal to $\Sit_1$ at $\o(1)$, we have
  \begin{align*}
    |\doto(1)^0-1| & \les \varep, & \le|\doto(1)^i\ri| & \les \varep.
  \end{align*}
  Projecting the geodesic equation $\D_{\doto}\doto = 0$ in the bottom initial layer coordinates, we obtain along $\o$ that
  \begin{align*}
    \le|\doto\le(\doto^\mu\ri)\ri| & \les |\Ga| \les \varep,
  \end{align*}
  where $\Ga$ are the Christoffel symbols in coordinates $x^\mu$, which are estimated using~\eqref{est:botlayassbis}. Integrating this estimate, we obtain
  \begin{align*}
    \le|\doto^0-1\ri| & \les \varep, & \le|\doto^i\ri| & \les \varep.
  \end{align*}
  We have
  \begin{align*}
    \le|\doto\le(x^0-t\ri)\ri| & = \le|\doto^0-1\ri| \les \varep, \\
    \le|\doto\le(x^i\ri)\ri| & = \le|\doto^i\ri| \les \varep,
  \end{align*}
  where we recall that $t$ denotes the affine parameter of $\o$. Integrating these estimates along $\o$, we thus have for all $t\geq 1$
  \begin{align*}
    \le|x^0(\o(t)) - t\ri| & \les \varep, & \le|x^i(\o(t))\ri| & \les \varep.
  \end{align*}
  From the definition of the last cones geodesic foliation, this implies that on the central axis $\o \cap \MM'$
  \begin{align*}
    \le|x^0-\half(u'+\ub')\ri| & \les \varep, & \le|\sqrt{\sum_{i=1}^3(x^i)^2}- \half(\ub'-u')\ri| & \les \varep.  
  \end{align*}
  To prove that the same estimates as above hold in the whole region $\LLb\cap\MM'$, we argue along the same lines, replacing the geodesic $\o$ by the geodesic along $\elb'$ emanating from the central axis $\o$ and using the estimates obtained at the central axis $\o$.\footnote{The same argument can be run in the region of $\MM'\cap\LLb$ which is not covered by the incoming null cones backwards emanating from $\o\cap\LLb$ (see Figure~\ref{fig:lastconesfoliation}). It indeed suffices to continue the initial layer up to time $x^0=100$.} Both~\eqref{est:compgeodbottomx0} and~\eqref{est:compgeodbottomTf},~\eqref{est:compgeodbottomNf} are obtained in this manner and this finishes the proof of the lemma.
\end{proof}

From the estimates~\eqref{est:compopticalgeod} for ${u'},{\ub'}$, one deduces that for $\varep$ sufficiently small, we have the following inclusion (see also Figure~\ref{fig:lastconesfoliation})
\begin{align*}
  \le\{\ut=3/2,\ubt=3\ri\} \subset \LLb\cap\LLc\cap\MM^\ext\cap\MM'.
\end{align*}

From the sup-norm estimates~\eqref{est:supnormf'ext},~\eqref{est:supnormfb'ext},~\eqref{est:supnormla'ext} for the transition coefficients $f',\fb',\la'$ between the null frames $(\elb',\el')$ and $(\elb,\el)$, and from the initial bounds~\eqref{est:compgeodbottomTf} and~\eqref{est:compgeodbottomNf}, we deduce
\begin{align}\label{est:compelbelprx0}
  \begin{aligned}
  \le|\g\le(\Tf^\bott,\half(\elb+\el)\ri)+1\ri| & \les \varep,\\
  \le|\g\le(\Nf^\bott,\half(\el-\elb)\ri)-1\ri| & \les \varep,
  \end{aligned}
\end{align}
on $\{\ut=3/2,\ubt=3\}$. Combining~\eqref{est:compelbelprx0} and the compatibility assumptions~\eqref{est:LLbotLLconframe} between the bottom initial layer and the conical initial layer frames, we deduce the following bounds for the transition coefficients $(\ft,\fbt,\lat)$ between the null frames $(\elbt,\elt)$ and $(\elb,\el)$ on $\{\ut=3/2,\ubt=3\}$
\begin{align}\label{est:ffbla'S}
  \le|\ft,\fbt,\log\lat\ri| & \les \varep.
\end{align}
Integrating successively equations~(\ref{eq:Nd3fb}), (\ref{eq:Nd3f}), (\ref{eq:Nd3la}) along $\elbt$ from $\{\ut=3/2,\ubt=3\}$ to $\CC_1\cap\{\ubt=3\}$, using the bounds~(\ref{est:transcoeffnullconnsource}) and~\eqref{est:Minkowskilayerepsconicnullconn} for the null connection coefficients of each foliation, and the bounds~\eqref{est:ffbla'S} on $\{\ut = 3/2,\ubt=3\}$, we infer that the bound
\begin{align}\label{est:ffbla'MMext}
  \le|\ft,\fbt,\log\lat\ri| & \les \varep
\end{align}
holds on $\MM^\ext \cap \{\ubt=3\}$. We rewrite~\eqref{eq:Nd4f} schematically as
\begin{align*}
  \Nd_4\ft + \half \trchi \ft & = \Err,
\end{align*}
which integrated along $\el$, from $\{\ubt=3\}\cap\MM^\ext$ into the region $\LLcext$, using the Bootstrap Assumptions~\ref{BA:con}, the bounds~(\ref{est:transcoeffnullconnsource}) and~\eqref{est:Minkowskilayerepsconicnullconn} for the null connection coefficients in $\LLcext$ and the bounds~\eqref{est:ffbla'MMext} on $\{\ubt=3\}$ gives
\begin{align}\label{est:f'LLcext}
  \ub |\ft| & \les \varep,
\end{align}
on $\LLcext$. Integrating equations~\eqref{eq:Nd4fb} and~\eqref{eq:Nd4la} from $\{\ubt=3\}\cap\MM^\ext$ into $\LLcext$ yields
\begin{align}\label{est:fbla'LLcext}
  \le|\fbt,\log\lat\ri| & \les \varep,
\end{align}
as desired. This finishes the proof of~\eqref{est:imprffbla}.\\

We turn to the proof of~\eqref{est:compopticalLLcext}. We deduce from the bounds~(\ref{est:compopticalgeod}), (\ref{est:compgeodbottomx0}) and the compatibility assumption~\eqref{est:LLbotLLconoptical} between the bottom and conical initial layers that
\begin{align*}
  \le|\ub-\ubt\ri| & \les \varep, & \le|u-\ut\ri| & \les \varep,
\end{align*}
at the $2$-sphere $\{\ut=3/2, \ubt =3\}$. Differentiating by $\elbt$ and integrating from $\{\ut=3/2,\ubt=3\}$ to $\CC_1\cap\{\ubt=3\}$, we further obtain
\begin{align}\label{est:uutub'3}
  \le|\ub-\ubt\ri| & \les \varep, & \le|u-\ut\ri| & \les \varep,
\end{align}
on $\{\ubt=3\}\cap\MM^\ext$. Using relations~\eqref{eq:elu} and the definitions from Section~\ref{sec:definitlayer} we have the following equations
\begin{align*}
  \begin{aligned}
    \el\le(u-\ut\ri) & = -\quar |\ft|^2 \elbt(\ut) = -\half |\ft|^2\Omt^{-1},\\
    \el\le(\ub-\ubt\ri) & = 2(1-\lat).  
  \end{aligned}
\end{align*}
Integrating the above equations along $\el$ from $\{\ubt=3\}\cap\MM^\ext$ into $\LLcext$, using the bounds~\eqref{est:uutub'3} on $\{\ubt=3\}$ and~\eqref{est:f'LLcext},~\eqref{est:fbla'LLcext},~\eqref{est:Minkowskilayerepsconicnullconn} we deduce
\begin{align*}
  \le|u-\ut\ri| & \les \varep, & \le|\ub-\ubt\ri| & \les \varep \ub,
\end{align*}
as desired. This finishes the proof of~\eqref{est:compopticalLLcext}.

\subsection{Bottom initial layer comparisons}\label{sec:framecompLLbint}
In this section, we prove the comparison estimates~(\ref{est:impprx0elbel}), (\ref{est:comptimeuub}), (\ref{est:impTfprx0}), (\ref{est:comptimeLLbint}) for the bottom initial layer. Following similar arguments as previously, one has
\begin{align*}
  \le|\D \le(\g\le(\Tf^\bott,\half(\elb+\el)\ri) \ri)\ri| & \les \varep
\end{align*}
in $\LLb\cap\MM^\ext$, and we thus deduce from the bound~(\ref{est:compelbelprx0}) on the $2$-sphere $\{\ut=3/2,\ubt=3\} \subset \LLbext$ by integration in $\LLbext$ that~(\ref{est:impprx0elbel}) holds in $\LLbext$. Further integration in $\LLbint$ also gives~\eqref{est:impTfprx0}. The bounds~(\ref{est:comptimeuub}) and~\eqref{est:comptimeLLbint} are then also deduced by integration, using the just obtained bounds~\eqref{est:compelbelprx0} and~\eqref{est:impTfprx0} and the bounds~(\ref{est:compopticalgeod}), (\ref{est:compgeodbottomx0}) on $\{\ut=3/2,\ubt=3\}$.

\begin{appendices}
\chapter{Global harmonic coordinates}\label{sec:globharmo}
This section is dedicated to the proof of Theorem~\ref{thm:globharmonics}.

\section{Definitions and preliminary results}\label{sec:defprelimglobharmo}
In this section, we recall definitions from Section~\ref{sec:definition} and state preliminary (known) results.

\subsection{Definitions and identities on $\pr\Si$}
\begin{definition}[Conformal isomorphism and $x^i$ on $\pr\Si$]\label{def:confiso}
  A \emph{conformal isomorphism} $\Phi$ is a diffeomorphism from $\pr\Si$ onto the Euclidean unit $2$-sphere $\SSS$, such that there exists a \emph{conformal factor} $\phi>0$ satisfying
  \begin{align*}
    \Phi_{\sharp}\gd_{\pr\Si} & = \phi^2\gd_{\SSS},
  \end{align*}
  where $\Phi_\sharp$ denotes the push-forward by $\Phi$.\\

  To a fixed conformal isomorphism, we associate a triplet of function $(x^i)_{i=1,2,3}$ on $\pr\Si$ to be the pull-back by $\Phi$ of the standard Cartesian coordinates restricted on $\SSS$.\\

  We moreover say that $\Phi$ is \emph{centred} if the following conditions hold for the functions $x^i$
  \begin{align*}
    \int_{\pr\Si}x^i & = 0, & i = 1,\cdots, 3.
  \end{align*}
\end{definition}

We have the following identities on $\pr\Si$.
\begin{lemma}[Conformal identities on $\pr\Si$]
  On $\pr\Si$, we have the following identities
  \begin{align}\label{est:sumxiprSi}
    \sum_{i=1}^3 (x^i)^2 & = 1,
  \end{align}
  and
  \begin{align}\label{eq:confLapxi}
    \Ld_{\gd_{\pr\Si}} x^i + 2 x^i & = 2x^i \le(1-\phi^{-2}\ri). 
  \end{align}
\end{lemma}
\begin{proof}
  The first identity~\eqref{est:sumxiprSi} is a direct consequence of Definition~\ref{def:confiso}. For the second identity~\eqref{eq:confLapxi}, we have using the Euclidean Laplacian in spherical coordinates that
  \begin{align*}
    \Ld_{\gd_{\SSS}}x^i + 2 x^i & = 0. 
  \end{align*}
  From the conformal invariance of the Laplacian, we have
  \begin{align*}
    \Ld_{\gd_{\SSS}}x^i = \phi^2\Ld_\gd x^i,
  \end{align*}
  and the desired formula follows.
\end{proof}

\subsection{Uniformisation theorem on $\pr\Si$}
\begin{lemma}[Uniformisation Theorem on $\pr\Si$]\label{lem:unifthm}
  Under the assumptions of Theorem~\ref{thm:globharmonics}, there exists a unique -- up to isomorphisms of $\SSS$ -- centred conformal isomorphism of the boundary $(\pr\Si,\gd)$ to the Euclidean sphere $\SSS$ and we have the following quantitative bounds for the conformal factor $\phi$
  \begin{align}\label{est:assconf}
    \norm{\Nd^{\leq 1}\le(\phi-1\ri)}_{H^{1/2}(\pr\Si)} & \les \varep.
  \end{align}
  Moreover, for all $k\geq 0$, we have the following higher regularity estimates
  \begin{align}
    \norm{\Nd^{\leq k+1}\le(\phi-1\ri)}_{H^{1/2}(\pr\Si)} & \les C_k\le(\norm{\nab^{\leq k}\RRRic}_{L^2(\Si)} + \norm{\Nd^{\leq k}(\th-\gd)}_{H^{1/2}(\pr\Si)} +\varep \ri).
  \end{align}
\end{lemma}
\begin{proof}
  Arguing as in~\cite{Czi19}, one can obtain from the $L^2(\Si)$ and $H^{1/2}(\pr\Si)$ bounds~(\ref{est:L2RicTh}) that 
  \begin{align*}
    \norm{K-1}_{H^{-1/2}(\pr\Si)} & \les \varep,
  \end{align*}
  where $H^{-1/2}(\pr\Si)$ is an (appropriately defined) fractional Sobolev space on $\pr\Si$ (see~\cite{Czi19}). The lemma then follows from an adaptation of Theorem~\cite[Theorem 3.1]{Kla.Sze19} (or from adaptations of the argument of~\cite[Section 6]{Sha14}).
\end{proof}





\subsection{Definitions and identities on $\Si$}
We now have the following definition for $x^i$ on $\Si$.
\begin{definition}[$x^i$ on $\Si$] \label{def:coordsxi}
  We define the $x^i$ on $\Si$ to be the solutions of the following Dirichlet problem
  \begin{align}\label{eq:defcoordsxi}
    \begin{aligned}
      \Delta_g x^i & = 0, \\
      x^i|_{\pr\Si} & = x^i.
    \end{aligned}
  \end{align}
\end{definition}

\begin{remark}
  From the maximum principle, one has on $\Si$
  \begin{align}\label{est:maxppl}
    \le|x^i\ri| & \leq 1.
  \end{align}
\end{remark}

We have the following energy and Bochner identities in $\Si$.
\begin{lemma}[Energy and Bochner identities in $\Si$]\label{lem:enerBochid}
  We have the following energy identity on $\Si$
  \begin{align}\label{est:enercoordsxi0}
    \begin{aligned}
      \norm{\nab x^i}^2_{L^2(\Si)} & = \int_{\pr\Si} x^i N(x^i),
    \end{aligned}
  \end{align}
  and the following Bochner identity on $\Si$
  \begin{align}\label{est:enercoordsxi1}
  \begin{aligned}
    \norm{\nab^2x^i}^2_{L^2(\Si)} & = -\int_{\Si}\RRRic\cdot\nab x^i\cdot\nab x^i + \int_{\pr\Si} \nab x^i \cdot \nab_N \nab x^i,
  \end{aligned}
  \end{align}
  for all $i=1,2,3$ and where here and in the following, $N$ denotes the outward-pointing unit normal to $\pr\Si$.
\end{lemma}
\begin{proof}
  The energy identity~\eqref{est:enercoordsxi0} is obtained by multiplying Laplace equation~\eqref{eq:defcoordsxi} by $x^i$ and integrating by part.\\

  Bochner identity~\eqref{est:enercoordsxi1} is obtained by using the following commutation formula
  \begin{align*}
    \Delta \nab x^i & = \RRRic \cdot \nab x^i, 
  \end{align*}
  contracted with $\nab x^i$ and integrating by part.
\end{proof}

\section{A refined Bochner identity}
This section is dedicated to the proof of the following refined Bochner identity, which is at the centre of the proof of Theorem~\ref{thm:globharmonics}
\begin{align}\label{eq:refinedBochner}
  \begin{aligned}
    \sum_{i=1}^3\norm{\nab^2 x^i}^2_{L^2(\Si)} + 2 \sum_{i=1}^3 \int_{\pr\Si} \le(N(x^i)-x^i\ri)^2 & = -3\sum_{i=1}^3\int_{\Si}\le(\RRRic-\frac{1}{3}Rg\ri)\cdot\nab x^i\cdot\nab x^i + \GGG,
  \end{aligned}
\end{align}
where
\begin{align*}
  \GGG & := \sum_{i=1}^3 \int_{\pr\Si}\bigg((\th_{ab}-\gd_{ab})x^i(\Nd_a\Nd_bx^i +2x^i\gd_{ab})\\
       & \quad\quad\quad -(N(x^i)-x^i)\le(2(\Ld x^i + 2x^i) -(\trth-2) (N(x^i)+x^i)\ri)\bigg)\\
       & + \half\int_{\pr\Si}(\trth-2)^2 - \int_{\pr\Si}|\thh|^2 + 2\sum_{i=1}^3 \int_{\pr\Si}  \le(\RRRic_{NN}-\half R\ri)x^i (N(x^i)-x^i) \\
    & \quad - 2\sum_{i=1}^3\int_{\Si}x^i \le(\RRRic-\half Rg\ri)\cdot \nab^2x^i.
\end{align*}

Identity~\eqref{eq:refinedBochner} is obtained by further computation starting from the original Bochner identity~\eqref{est:enercoordsxi1}. Its proof is postponed to the end of this section and uses the following two lemmas.
\begin{lemma}\label{lem:structureintprSi}
  We have
  \begin{align}\label{eq:structureintprSi}
    \begin{aligned}
      \sum_{i=1}^3 \int_{\pr\Si} \nab x^i \cdot \nab_N \nab x^i & = 2\int_{\pr\Si}(\phi^{-2}-1) - 2\int_{\pr\Si}(\trth-2) -2 \sum_{i=1}^3\int_{\pr\Si}(N(x^i)-x^i)^2 + \EEE,
    \end{aligned}
  \end{align}
  where
  \begin{align*}
    \EEE & : = \sum_{i=1}^3 \int_{\pr\Si}\bigg((\th_{ab}-\gd_{ab})x^i(\Nd_a\Nd_bx^i +2x^i\gd_{ab})\\
           & \quad\quad\quad -(N(x^i)-x^i)\le(2(\Ld x^i + 2x^i) -(\trth-2) (N(x^i)+x^i)\ri)\bigg).
  \end{align*}
\end{lemma}
\begin{proof}
  We decompose the integrand in two parts as follows
  \begin{align}\label{eq:decomp2int}
    \nab x^i \cdot \nab_N\nab x^i & = \nab_ax^i\nab_N\nab_ax^i + \nab_Nx^i\nab_N\nab_N x^i.
  \end{align}

  For the first term of~\eqref{eq:decomp2int}, we have\footnote{We recall that $\th_{ab} = \nab_aN_b$ (see the definition in Theorem~\ref{thm:globharmonics}).}
  \begin{align*}
    \begin{aligned}
      \nab_ax^i\nab_N\nab_ax^i & = \Nd_ax^i\nab_a\nab_Nx^i \\
      & = \Nd_ax^i \le(\Nd_aN(x^i) - \th_{ab}\Nd_bx^i\ri) \\
      & = \Nd_ax^i\le(\Nd_a\le(N(x^i)-x^i\ri) - (\th_{ab}-\gd_{ab})\Nd_bx^i \ri).
    \end{aligned}
  \end{align*} 

  Integrating that relation on $\pr\Si$ and integrating by part gives
  \begin{align*}
    \begin{aligned}
    & \int_{\pr\Si} \nab_ax^i\nab_N\nab_ax^i \\
    = & \; \int_{\pr\Si} \Nd_ax^i\le(\Nd_a\le(N(x^i)-x^i\ri) - (\th_{ab}-\gd_{ab})\Nd_bx^i\ri) \\
    = & \; \int_{\pr\Si}(-2x^i-\Ld x^i)\le(N(x^i)-x^i\ri)  - \int_{\pr\Si} (\th_{ab}-\gd_{ab})\Nd_a x^i\Nd_b x^i \\
    & + \int_{\pr\Si} 2x^i\le(N(x^i)-x^i\ri).
    \end{aligned}
  \end{align*}

  Integrating by part and using relation~(\ref{est:sumxiprSi}), we have
  \begin{align*}
    \begin{aligned}
      & -\sum_{i=1}^3 \int_{\pr\Si} (\th_{ab}-\gd_{ab})\Nd_a x^i\Nd_b x^i  \\
      = & \; \sum_{i=1}^3 \le(\int_{\pr\Si} \half \Nd_a\th_{ab}\Nd_b(x^i)^2 + \int_{\pr\Si}(\th_{ab}-\gd_{ab}) x^i \Nd_a\Nd_bx^i\ri)\\
      = & -\int_{\pr\Si} (\trth -2) + \sum_{i=1}^3 \int_{\pr\Si}(\th_{ab}-\gd_{ab})x^i(\Nd_a\Nd_bx^i +x^i\gd_{ab}).
    \end{aligned}
  \end{align*}

  Summing for $i$ from $1$ to $3$, we therefore infer
  \begin{align}\label{eq:sumiint1}
    \begin{aligned}
      & \sum_{i=1}^3 \int_{\pr\Si} \nab_ax^i\nab_N\nab_ax^i \\
      = & \;  \sum_{i=1}^3 \int_{\pr\Si}(-2x^i-\Ld x^i)\le(N(x^i)-x^i\ri) + \sum_{i=1}^3 \int_{\pr\Si}(\th_{ab}-\gd_{ab})x^i(\Nd_a\Nd_bx^i +x^i\gd_{ab}) \\
      & + \sum_{i=1}^3 \int_{\pr\Si} 2x^i\le(N(x^i)-x^i\ri) - \int_{\pr\Si} (\trth -2).
    \end{aligned}
  \end{align}

  Let turn to the second term of~\eqref{eq:decomp2int}. We first recall that Laplace equation~\eqref{eq:defcoordsxi} can be rewritten as
  \begin{align*}
    \nab_N\nab_Nx^i + \Ld_{\gd}x^i + \trth N(x^i) & = 0.
  \end{align*}

  Using this, we compute
  \begin{align*}
    \nab_Nx^i\nab_N\nab_N x^i & = - N(x^i)\le(\Ld x^i + \trth N(x^i)\ri) \\
                              & = - N(x^i)\le((\Ld x^i + 2x^i) + 2(N(x^i)-x^i) + (\trth-2) N(x^i)\ri) \\
                              & = -x^i\le((\Ld x^i + 2x^i) + 2(N(x^i)-x^i) + (\trth-2) N(x^i)\ri) \\
                              & \quad - (N(x^i)-x^i)\le((\Ld x^i + 2x^i) + 2(N(x^i)-x^i) + (\trth-2) N(x^i)\ri) \\
                              & = -x^i(\Ld x^i + 2x^i) -2x^i(N(x^i)-x^i) - (x^i)^2(\trth-2) \\
                              & \quad - (N(x^i)-x^i)\le((\Ld x^i + 2x^i) + 2(N(x^i)-x^i) + (\trth-2) (N(x^i)+x^i)\ri)
  \end{align*}

  Integrating this relation on $\pr\Si$, gives
  \begin{align*}
    & \int_{\pr\Si} \nab_Nx^i\nab_N\nab_N x^i  \\
    = & \; \int_{\pr\Si}-x^i(\Ld x^i + 2x^i) - \int_{\pr\Si}(x^i)^2(\trth-2) + \int_{\pr\Si}-2x^i(N(x^i)-x^i)\\
    & - \int_{\pr\Si} (N(x^i)-x^i)\le((\Ld x^i + 2x^i) + 2(N(x^i)-x^i) + (\trth-2) (N(x^i)+x^i)\ri).
  \end{align*}

  Using relations~(\ref{eq:confLapxi}), (\ref{est:sumxiprSi}), we have
  \begin{align*}
    \begin{aligned}
      \sum_{i=1}^3-\int_{\pr\Si}x^i(\Ld x^i + 2x^i) & = -\sum_{i=1}^32\int_{\pr\Si}(x^i)^2\le(1-\phi^{-2}\ri) = 2\int_{\pr\Si}(\phi^{-2}-1).
    \end{aligned}
  \end{align*}

  Summing over $i$ and using this relation and~\eqref{est:sumxiprSi}, we infer
  \begin{align}\label{eq:sumiint2}
    \begin{aligned}
      & \sum_{i=1}^3\int_{\pr\Si} \nab_Nx^i\nab_N\nab_N x^i \\
      = & \; 2\int_{\pr\Si}(\phi^{-2}-1) -\int_{\pr\Si}(\trth-2) -2 \sum_{i=1}^3 \int_{\pr\Si}x^i(N(x^i)-x^i) \\
      & -\sum_{i=1}^3 \int_{\pr\Si} (N(x^i)-x^i)\le((\Ld x^i + 2x^i) + 2(N(x^i)-x^i) + (\trth-2) (N(x^i)+x^i)\ri).
    \end{aligned}
  \end{align}

  Summing~\eqref{eq:sumiint1} and~\eqref{eq:sumiint2}, we have
  \begin{align*}
    \begin{aligned}
      \sum_{i=1}^3 \int_{\pr\Si} \nab x^i \cdot \nab_N \nab x^i & =  2\int_{\pr\Si}(\phi^{-2}-1) - 2\int_{\pr\Si}(\trth-2) -2 \sum_{i=1}^3\int_{\pr\Si}(N(x^i)-x^i)^2  \\
      & \quad + \sum_{i=1}^3 \int_{\pr\Si}(\th_{AB}-\gd_{AB})x^i(\Nd_A\Nd_Bx^i +2x^i\gd_{AB}) \\
      & -\sum_{i=1}^3 \int_{\pr\Si} (N(x^i)-x^i)\le(2(\Ld x^i + 2x^i) + (\trth-2) (N(x^i)+x^i)\ri),      
    \end{aligned}
  \end{align*}
  which finishes the proof of the lemma.
\end{proof}

We have the following lemma, which treats the linear boundary term of Lemma~\ref{lem:structureintprSi}.
\begin{lemma}
  \label{lem:Einsteinintegral}
  We have
  \begin{align}
    \label{eq:Einsteinintegral}
    2\int_{\pr\Si}(\phi^{-2}-1) - 2\int_{\pr\Si}(\trth-2) & = -2\sum_{i=1}^3\int_{\Si}\le(\RRRic-\half R g\ri)\cdot\nab x^i\cdot\nab x^i + \FFF,
  \end{align}
  where
  \begin{align*}
    \FFF & := \half\int_{\pr\Si}(\trth-2)^2 - \int_{\pr\Si}|\thh|^2 + 2\sum_{i=1}^3 \int_{\pr\Si}  \le(\RRRic_{NN}-\half R\ri)x^i (N(x^i)-x^i) \\
    & \quad - 2\sum_{i=1}^3\int_{\Si}x^i \le(\RRRic-\half Rg\ri)\cdot \nab^2x^i.
  \end{align*}
\end{lemma}
\begin{proof}
  We define the following Einstein tensor on $\Si$
  \begin{align*}
    G & := \RRRic -\half R g.
  \end{align*}
  Using Gauss equation (see~\cite[equation 5.0.5d]{Chr.Kla93}), we obtain
  \begin{align*}
    G_{NN} & = -K +\quar \trth^2 -\half |\thh|^2,
  \end{align*}
  on $\pr\Si$. Integrating the above on $\pr\Si$, using Gauss-Bonnet formula, we infer
  \begin{align*}
    \int_{\pr\Si}G_{NN} & = -4\pi + \int_{\pr\Si}1 +\int_{\pr\Si}(\trth-2)  + \quar\int_{\pr\Si}(\trth-2)^2 - \half\int_{\pr\Si}|\thh|^2,\\
                         & = -\int_{\pr\Si}(\phi^{-2}-1) +\int_{\pr\Si}(\trth-2)  + \quar\int_{\pr\Si}(\trth-2)^2 - \half\int_{\pr\Si}|\thh|^2.
  \end{align*}
  We thus deduce that
  \begin{align}\label{eq:Gaussflux1}
    2\int_{\pr\Si}(\phi^{-2}-1)-2\int_{\pr\Si}(\trth-2) & = -2 \int_{\pr\Si}G_{NN} + \FFF^1,
  \end{align}
  where
  \begin{align*}
    \FFF^1 & :=  \half\int_{\pr\Si}(\trth-2)^2 - \int_{\pr\Si}|\thh|^2.
  \end{align*}

  Using Bianchi identity on $\Si$
  \begin{align*}
    \Div G & = 0, 
  \end{align*}
  and applying Stokes formula, we have
  \begin{align}\label{eq:StokesG1}
    \begin{aligned}
    \int_{\Si}G\cdot\nab x^i\cdot\nab x^i & = \int_{\pr\Si} G_{lN} x^i \nab_l x^i -\int_{\Si}x^i G\cdot \nab^2x^i.
  \end{aligned}
  \end{align}
  Summing over $i$ and using~\eqref{est:sumxiprSi}, we have 
  \begin{align}\label{eq:StokesG2}
    \begin{aligned}
      \sum_{i=1}^3 \int_{\pr\Si} G_{Nl}x^i \nab_lx^i & = \sum_{i=1}^3 \int_{\pr\Si} G_{NN}x^i N(x^i) + \half\int_{\pr\Si}G_{Na}\Nd_a\le(\sum_{i=1}^3(x^i)^2\ri) \\
      & = \int_{\pr\Si}G_{NN} + \sum_{i=1}^3 \int_{\pr\Si} G_{NN}x^i (N(x^i)-x^i).
    \end{aligned}
  \end{align}
  Combining~\eqref{eq:StokesG1} and~\eqref{eq:StokesG2}, we infer
  \begin{align}\label{eq:StokesG3}
    \begin{aligned}
      -2 \int_{\pr\Si}G_{NN} & = -2\sum_{i=1}^3\int_{\Si}G\cdot\nab x^i\cdot\nab x^i + \FFF^2,
    \end{aligned}
  \end{align}
  where
  \begin{align*}
    \FFF^2 & := 2\sum_{i=1}^3 \int_{\pr\Si} G_{NN}x^i (N(x^i)-x^i) - 2\sum_{i=1}^3\int_{\Si}x^i G\cdot \nab^2x^i.
  \end{align*}
  Combining~\eqref{eq:Gaussflux1} and~\eqref{eq:StokesG3} then gives the desired result.  
\end{proof}
\begin{proof}[Proof of~\eqref{eq:refinedBochner}]
Using the results of Lemmas~\ref{lem:structureintprSi} and \ref{lem:Einsteinintegral}, and Bochner identity~\eqref{est:enercoordsxi1}, we have
\begin{align*}
  \begin{aligned}
    \sum_{i=1}^3\norm{\nab^2 x^i}^2_{L^2(\Si)} & = -\sum_{i=1}^3\int_{\Si}\RRRic\cdot\nab x^i\cdot\nab x^i + \sum_{i=1}^3\int_{\pr\Si} \nab x^i \cdot \nab_N \nab x^i \\
    & = -\sum_{i=1}^3\int_{\Si}\RRRic\cdot\nab x^i\cdot\nab x^i -2 \sum_{i=1}^3\int_{\Si}\le(\RRRic-\half R g\ri)\cdot\nab x^i \cdot \nab x^i \\
    & \quad -2 \sum_{i=1}^3\int_{\pr\Si}(N(x^i)-x^i)^2 + \EEE + \FFF \\
    & = -3\sum_{i=1}^3\int_{\Si}\le(\RRRic-\frac{1}{3} R g\ri)\cdot\nab x^i \cdot \nab x^i -2 \sum_{i=1}^3\int_{\pr\Si}(N(x^i)-x^i)^2+\EEE+\FFF,
  \end{aligned}
\end{align*}
and using the expressions of $\EEE,\FFF$ this concludes the proof of~(\ref{eq:refinedBochner}).
\end{proof}

\section{Refined Bochner estimate}\label{sec:enerBoch}
This section is dedicated to show that, using the refined Bochner identity~\eqref{eq:refinedBochner} together with the assumptions~\eqref{est:L2RicTh} of Theorem~\ref{thm:globharmonics} one has the following estimate
\begin{align}\label{est:refinedBochner}
  \sum_{i=1}^3\norm{\nab^2x^i}^2_{L^2(\Si)} + \sum_{i=1}^3\int_{\pr\Si}(N(x^i)-x^i)^2 + \norm{\BBB}^2_{L^2(\Si)}  & \les \varep^2,
\end{align}
where $B$ is the following $2$-tensor on $\Si$
\begin{align*}
  \BBB & := \sum_{i=1}^3\nab x^i \otimes \nab x^i - g,
\end{align*}
where $\otimes$ denotes the standard tensorial product on $T\Si$. As a consequence, we will also obtain the following estimates
\begin{align}
  \label{est:consrefinedBochner}
  \norm{\nab x^i}_{L^2(\Si)} + \norm{\nab x^i}_{L^6(\Si)} & \les 1.
\end{align}

The proof of~\eqref{est:refinedBochner} and~\eqref{est:consrefinedBochner} is postponed to the end of this section. It relies on the following three lemmas.
\begin{lemma}\label{lem:Bochnerrefined1}
  We have
  \begin{align}
    \begin{aligned}
      & \sum_{i=1}^3\le(\norm{\nab^2x^i}^2_{L^2(\Si)} + \norm{N(x^i)-x^i}^2_{L^2(\pr\Si)}\ri) \\
      \les \; & \varep\le(\norm{\nab^2x^i}_{L^2(\Si)} + \norm{N(x^i)-x^i}_{L^2(\pr\Si)}\ri)\le(1+\norm{\nab^{\leq 2}x^i}_{L^2(\Si)}\ri) + \varep \norm{\BBB}_{L^2(\Si)} + \varep^2.
    \end{aligned}
  \end{align}
\end{lemma}
\begin{proof}
  Using the assumptions~\eqref{est:L2RicTh}, we have
  \begin{align}\label{est:RHSnab2xiB}
    \begin{aligned}
      -3\sum_{i=1}^3\int_{\Si}\le(\RRRic-\frac{1}{3} R g\ri)\cdot\nab x^i \cdot \nab x^i & = -3\int_{\Si}\le(\RRRic-\frac{1}{3} R g\ri)\cdot\BBB \\
      & \les \norm{\RRRic}_{L^2(\Si)}\norm{\BBB}_{L^2(\Si)}\\
      & \les \varep \norm{\BBB}_{L^2(\Si)}.
    \end{aligned}
  \end{align}
  Using a trace estimate and Bianchi identity for Einstein tensor $\Div G = 0$, we claim that the following bilinear estimate holds for all scalar functions $f$ on $\pr\Si$\footnote{This can be obtained by locally defining a geodesic foliation by $2$-spheres extending the boundary $\pr\Si$, extending $f$, integrating along the normal direction using Bianchi identity and integration by part. See~\cite[Proposition 2.9]{Czi19} for the proof of a similar estimate.}
  \begin{align}\label{est:traceEinstein}
    \le|\int_{\pr\Si}G_{NN}f \ri| & \les \norm{G}_{L^2(\Si)}\norm{f}_{H^{1/2}(\pr\Si)}.
  \end{align}
  
  Using~\eqref{est:traceEinstein} and the trace estimate of Definition~\ref{def:weakreg3DStab}, Laplace equation~(\ref{eq:confLapxi}) for $x^i$ on $\pr\Si$, and~\eqref{est:maxppl}, we have the following estimates for the terms composing $\GGG$\footnote{Lines 4 and 5 are obtained by replacing $^{\gd}\Nd$ by their analogue on the round sphere $^{\gd_{\SSS^2}}\Nd$. See also the proof of Lemma~\ref{lem:mildOOESast}.}
  \begin{align*}
    \begin{aligned}
      \le|\int_{\pr\Si}  \le(\RRRic_{NN}-\half R\ri)x^i (N(x^i)-x^i)\ri| & \les \norm{\RRRic}_{L^2(\Si)}\norm{\nab^2x^i}_{L^2(\Si)}, \\
      \le|\int_{\Si}x^i \le(\RRRic-\half Rg\ri)\cdot \nab^2x^i\ri| & \les \norm{\RRRic}_{L^2(\Si)}\norm{\nab^2x^i}_{L^2(\Si)}, \\
      \int_{\pr\Si}(\trth-2)^2 + \int_{\pr\Si}|\thh|^2 & \les \norm{\th-\gd}^2_{L^4(\pr\Si)},\\
      \le|\int_{\pr\Si}(\th_{ab}-\gd_{ab})x^i(\Nd_a\Nd_bx^i +2x^i\gd_{ab})\ri| & \les \norm{\th-\gd}_{L^4(\pr\Si)}\norm{\Nd^{\leq 1}(\phi-1)}_{H^{1/2}(\SSS)},\\
      \le|\int_{\pr\Si}(N(x^i)-x^i)(\Ld x^i + 2x^i)\ri| & \les \norm{N(x^i)-x^i}_{L^2(\pr\Si)}\norm{\phi^{-2}-1}_{L^2(\pr\Si)} \\
      \le|\int_{\pr\Si}(N(x^i)-x^i)(\trth-2) (N(x^i)+x^i)\ri| & \les \norm{N(x^i)-x^i}_{L^2(\pr\Si)} \norm{\trth-2}_{L^4(\pr\Si)}\norm{\nab^{\leq 1}x^i}_{L^4(\pr\Si)} \\
      & \les \norm{N(x^i)-x^i}_{L^2(\pr\Si)} \norm{\trth-2}_{L^4(\pr\Si)}\norm{\nab^{\leq 2}x^i}_{L^2(\Si)}. 
    \end{aligned}
  \end{align*}
  Thus, using estimates~\eqref{est:L2RicTh},~(\ref{est:assconf}) and~(\ref{est:maxppl}), we obtain that
  \begin{align}\label{est:RHSnab2xiG}
    |\GGG| & \les \varep\le(\norm{\nab^2x^i}_{L^2(\Si)} + \norm{N(x^i)-x^i}_{L^2(\pr\Si)}\ri)\le(1+\norm{\nab^{\leq 2}x^i}_{L^2(\Si)}\ri) + \varep^2.
  \end{align}
  Combining~\eqref{est:RHSnab2xiB} and~\eqref{est:RHSnab2xiG} and the refined Bochner identity~\eqref{eq:refinedBochner}, we have
  \begin{align*}
    \begin{aligned}
      & \sum_{i=1}^3\le(\norm{\nab^2x^i}^2_{L^2(\Si)} + \norm{N(x^i)-x^i}^2_{L^2(\pr\Si)}\ri) \\
      \les \; & \varep\le(\norm{\nab^2x^i}_{L^2(\Si)} + \norm{N(x^i)-x^i}_{L^2(\pr\Si)}\ri)\le(1+\norm{\nab^{\leq 2}x^i}_{L^2(\Si)}\ri) + \varep \norm{\BBB}_{L^2(\Si)} + \varep^2,
    \end{aligned}
  \end{align*}
  as desired
\end{proof}

\begin{lemma}\label{lem:energyestnabxi}
  We have
  \begin{align*}
    \norm{\nab x^i}_{L^2(\Si)} & \les \le((\int_{\pr\Si}1) + \norm{N(x^i)-x^i}_{L^2(\pr\Si)}\ri)^{1/2},
  \end{align*}
  and
  \begin{align*}
    \norm{\nab x^i}_{L^6(\Si)} & \les \le((\int_{\pr\Si}1)+ \norm{N(x^i)-x^i}_{L^2(\pr\Si)}\ri)^{1/2} + \norm{\nab^2x^i}_{L^2(\Si)}.
  \end{align*}
\end{lemma}
\begin{proof}
  The first estimate follows from the energy identity~(\ref{est:enercoordsxi0}) and~(\ref{est:maxppl}). The second follows from the first and the Sobolev embeddings of Definition~\ref{def:weakreg3DStab}.
\end{proof}

\begin{lemma}\label{lem:estB}
  We have
  \begin{align}
    \begin{aligned}
    \norm{\BBB}_{L^2(\Si)} & \les \le(\le((\int_{\pr\Si}1)+ \norm{N(x^i)-x^i}_{L^2(\pr\Si)}\ri)^{1/2} + \norm{\nab^2x^i}_{L^2(\Si)} \ri)\norm{\nab^2x^i}_{L^2(\Si)} \\
    & \quad + \norm{N(x^i)-x^i}_{L^2(\pr\Si)} + \norm{N(x^i)-x^i}_{L^2(\pr\Si)}^2 +\varep.
    \end{aligned}
  \end{align}
\end{lemma}
\begin{proof}
  From H\"older estimate and the estimates of Lemma~\ref{lem:energyestnabxi}, we first have
  \begin{align}\label{est:L3/2BBBSi}
    \begin{aligned}
    \norm{\nab \BBB}_{L^{3/2}(\Si)} & \les \norm{\nab x^i}_{L^6(\Si)}\norm{\nab^2x^i}_{L^2(\Si)} \\
    & \les \le(\le((\int_{\pr\Si}1)+ \norm{N(x^i)-x^i}_{L^2(\pr\Si)}\ri)^{1/2} + \norm{\nab^2x^i}_{L^2(\Si)} \ri)\norm{\nab^2x^i}_{L^2(\Si)}.
    \end{aligned}
  \end{align}

  Using~(\ref{est:sumxiprSi}), we have at the boundary $\pr\Si$
  \begin{align*}
    \BBB_{NN} & = \sum_{i=1}^3N(x^i)N(x^i) -1 = 2 \sum_{i=1}^3 x^i(N(x^i)-x^i) + \sum_{i=1}^3(N(x^i)-x^i)^2,\\
    \BBB_{Na} & = \sum_{i=1}^3N(x^i)\Nd_a(x^i) = \sum_{i=1}^3(N(x^i)-x^i)\Nd_a(x^i),\\
    \BBB_{ab} & = \sum_{i=1}^3 \Nd_a(x^i)\Nd_b(x^i) -\gd_{ab} = \le(1-\phi^{2}\ri)(\gd_{\SSS})_{ab}.
  \end{align*}
  Thus, using estimates~(\ref{est:assconf}),~(\ref{est:maxppl}), we deduce
  \begin{align}\label{est:BBBL1prSi}
    \norm{\BBB}_{L^1(\pr\Si)} &  \les \norm{N(x^i)-x^i}_{L^2(\pr\Si)} + \norm{N(x^i)-x^i}_{L^2(\pr\Si)}^2 + \varep.
  \end{align}
  Combining~\eqref{est:L3/2BBBSi} and~\eqref{est:BBBL1prSi} and the Poincar\'e-type estimate~\eqref{est:PoincareL3/2}, the desired result follows.
\end{proof}

\begin{proof}[Proof of~\eqref{est:refinedBochner} and~\eqref{est:consrefinedBochner}]
  Combining the results of Lemmas~\ref{lem:Bochnerrefined1}, \ref{lem:energyestnabxi}, \ref{lem:estB}, we obtain the following estimate
  \begin{align*}
    & \sum_{i=1}^3\norm{\nab^2x^i}^2_{L^2(\Si)} + \sum_{i=1}^3\norm{N(x^i)-x^i}^2_{L^2(\pr\Si)} \\
    \les \; & \varep\le(\norm{\nab^2x^i}_{L^2(\Si)} + \norm{N(x^i)-x^i}_{L^2(\pr\Si)}\ri)\le(1+\norm{\nab^{\leq 2}x^i}_{L^2(\Si)}\ri) + \varep \norm{\BBB}_{L^2(\Si)} + \varep^2 \\
    \les \; & \varep \le(\norm{\nab^2x^i}_{L^2(\Si)} + \norm{N(x^i)-x^i}_{L^2(\pr\Si)}\ri)\le(1+\le((\int_{\pr\Si}1) + \norm{N(x^i)-x^i}_{L^2(\pr\Si)}\ri)^{1/2} + \norm{\nab^{2}x^i}_{L^2(\Si)}\ri) \\
    & + \varep \le(\le((\int_{\pr\Si}1)+ \norm{N(x^i)-x^i}_{L^2(\pr\Si)}\ri)^{1/2} + \norm{\nab^2x^i}_{L^2(\Si)} \ri)\norm{\nab^2x^i}_{L^2(\Si)} \\
    & + \varep\norm{N(x^i)-x^i}_{L^2(\pr\Si)} + \varep\norm{N(x^i)-x^i}_{L^2(\pr\Si)}^2 + \varep^2.
  \end{align*}
  Using Young's inequality and absorption or direct absorption for $\varep>0$ sufficiently small, we obtain from the above estimate that
  \begin{align*}
    \sum_{i=1}^3\norm{\nab^2x^i}^2_{L^2(\Si)} + \sum_{i=1}^3\norm{N(x^i)-x^i}^2_{L^2(\pr\Si)}  & \les \varep^2.
  \end{align*}
  Using this estimate and Lemmas~\ref{lem:energyestnabxi}, \ref{lem:estB}, we further obtain
  \begin{align*}
    \norm{\BBB}_{L^2(\Si)} & \les \varep, & \norm{\nab x^i}_{L^2(\Si)} + \norm{\nab x^i}_{L^6(\Si)} & \les 1. 
  \end{align*}
  This finishes the proof of~\eqref{est:refinedBochner} and~\eqref{est:consrefinedBochner}.
\end{proof}

\section{Higher order estimates}\label{sec:nab3xiSi}
In this section, we prove the following higher order estimate for $\nab^3x^i$. Estimates for higher derivatives will follow similarly (see Lemma~\ref{lem:ellSixi}) and are left to the reader.
\begin{align}\label{est:enercoordsxi2}
  \begin{aligned}
    \norm{\nab^3x^i}_{L^2(\Si)} & \les \norm{\RRRic}_{L^2(\Si)} + \norm{\th-\gd}_{H^{1/2}(\pr\Si)} + \norm{\Nd^{\leq 1}(\phi-1)}_{H^{1/2}(\pr\Si)} + \varep.
  \end{aligned}
\end{align}
As a consequence, using the $L^2$ assumptions~\eqref{est:L2RicTh}, we will have 
\begin{align}\label{est:nab3xiSi}
  \norm{\nab^3x^i}_{L^2(\Si)} & \les \varep.
\end{align}

We turn to the proof of~\eqref{est:enercoordsxi2}. First, the commuted Laplace equation for $\nab^2x^i$ takes the following schematic form
\begin{align}\label{eq:Deltanab2xi}
  \Delta \nab^2x^i & = \nab \RRRic \cdot \nab x^i + \RRRic\cdot \nab^2x^i.
\end{align}
for all $i=1,2,3$.\\

Since we have obtained $L^2$ smallness estimates for $\nab^2x^i$ in the previous section, estimates for derivatives of $\nab^2x^i$ will follow from~\eqref{eq:Deltanab2xi} and from the following standard higher order elliptic estimates. 
\begin{lemma}\label{lem:ellSixi}
  For any $\Si$-tangent tensor $F$ satisfying
  \begin{align*}
    \Delta F = \nab G + H,
  \end{align*}
  we have
  \begin{align*}
    \norm{\nab F}_{L^2(\Si)} & \les \norm{G}_{L^2(\Si)} + \norm{H}_{L^{4/3}(\Si)} + \norm{F}_{L^2(\Si)} + \norm{\Fslash}_{H^{1/2}(\pr\Si)},
  \end{align*}
  where $\Fslash$ denotes the projection of $F$ as a $\pr\Si$-tangent tensor. Moreover, we have
  \begin{align*}
    \norm{\nab^2F}_{L^2(\Si)} & \les \norm{\Delta F}_{L^2(\Si)} +\norm{\nab^{\leq 1}F}_{L^2(\Si)} + \norm{\Nd^{\leq 1}\Fslash}_{H^{1/2}(\pr\Si)}.
  \end{align*}
\end{lemma}
\begin{remark}
  The elliptic estimates from Lemma~\ref{lem:ellSixi} can either be taken as assumptions additional to the functional assumptions of Definition~\ref{def:weakreg3DStab}, or proved using the $L^2$ bounds~\eqref{est:L2RicTh}, Stokes formula, the trace estimates from Definition~\ref{def:weakreg3DStab} and $H^{-1/2}\times H^{1/2}$ estimates on the boundary $\pr\Si$.
\end{remark}

Using the results of Lemma~\ref{lem:ellSixi}, Sobolev embeddings from Definition~\ref{def:weakreg3DStab}, we have
\begin{align}\label{est:nabxi31}
  \begin{aligned}
  \norm{\nab^3x^i}_{L^2(\Si)} & \les \norm{\RRRic\cdot \nab x^i}_{L^2(\Si)} + \norm{\RRRic \cdot\nab^2x^i}_{L^{4/3}} + \norm{\nab^2_{a,b} x^i}_{H^{1/2}(\pr\Si)} \\
  & \les \norm{\RRRic}_{L^2(\Si)}\le(\norm{\nab^{\leq 3}x^i}_{L^2(\Si)}\ri) + \norm{\nab_{a,b}^2x^i}_{H^{1/2}(\pr\Si)}.
  \end{aligned}
\end{align}

We have
\begin{align*}
  \nab^2_{a,b}x^i & = \Nd^2_{a,b}x^i + \th_{ab}N(x^i) = \le(\Nd_{a,b}^2x^i + x^i\gd_{ab}\ri) + x^i(\th_{ab}-\gd_{ab}) + (N(x^i)-x^i)\th_{ab}. 
\end{align*}
Thus, we deduce using standard $H^{1/2}$ product estimates (see~\cite{Sha14})
\begin{align*}
  \norm{\nab^2_{a,b}x^i}_{H^{1/2}(\pr\Si)} & \les \norm{\Nd^{\leq 1}(\phi-1)}_{H^{1/2}(\pr\Si)} + \norm{\th-\gd}_{H^{1/2}(\pr\Si)} \\
  & \quad + \le(\norm{N(x^i)-x^i}_{L^\infty(\pr\Si)} + \norm{\Nd(N(x^i)-x^i)}_{L^2(\pr\Si)}\ri)\norm{\th}_{H^{1/2}(\pr\Si)}.
\end{align*}
By Sobolev embeddings on $\pr\Si$, we have for the last terms
\begin{align*}
  & \norm{N(x^i)-x^i}_{L^\infty(\pr\Si)} + \norm{\Nd(N(x^i)-x^i)}_{L^2(\pr\Si)} \\
  \les & \; \norm{\Nd(N(x^i)-x^i)}^{1/2}_{H^{1/2}(\pr\Si)}\norm{N(x^i)-x^i}^{1/2}_{H^{1/2}(\pr\Si)} + \norm{N(x^i)-x^i}_{H^{1/2}(\pr\Si)},
\end{align*}
thus from Young's inequality and the trace estimates of Definition~\ref{def:weakreg3DStab}, we finally obtain
\begin{align}\label{est:nab2abxi}
  \begin{aligned}
    & \norm{\nab^2_{a,b}x^i}_{H^{1/2}(\pr\Si)} -c\norm{\Nd(N(x^i)-x^i)}_{H^{1/2}(\pr\Si)} \\
   \les \; & \norm{\Nd^{\leq 1}(\phi-1)}_{H^{1/2}(\pr\Si)} + \norm{\th-\gd}_{H^{1/2}(\pr\Si)}+ \norm{N(x^i)-x^i}_{H^{1/2}(\pr\Si)},
    \end{aligned}
\end{align}
where $c>0$ is a small constant. Combining~\eqref{est:nabxi31} and~\eqref{est:nab2abxi} with an absorption argument (provided that $c>0$ is sufficiently small) yields~\eqref{est:enercoordsxi2} as desired.

\section{The $x^i$ are local coordinates on $\Si$}\label{sec:localcoordsxi}
In this section, we prove that the following estimate holds
\begin{align}
  \label{est:unifgijdeij}
  \sum_{i,j=1}^3\le|\g(\nab x^i,\nab x^j) -\de_{ij}\ri| & \les \varep,
\end{align}
uniformly on $\Si$. As a consequence, using the local inverse theorem we deduce that the $x^i$ form a local coordinate system on $\Si$.\\

We turn to the proof of~\eqref{est:unifgijdeij}. Using the Sobolev embeddings from Definition~\ref{def:weakreg3DStab}, the estimates~\eqref{est:nab3xiSi},~\eqref{est:refinedBochner} and \eqref{est:consrefinedBochner}, we obtain for the tensor $\BBB = \sum_{i=1}^3\nab x^i\otimes\nab x^i-g$ defined in Section~\ref{sec:enerBoch}
\begin{align}
  \label{est:supnormBSi}
  \norm{\BBB}_{L^\infty(\Si)} & \les \varep.
\end{align}
We now have the following geometric lemma, which achieves the proof of~\eqref{est:unifgijdeij}.
\begin{lemma}
  Assume that $|\BBB| \les \varep$. Then, for $\varep>0$ sufficiently small, we have
  \begin{align*}
    \sum_{i,j=1}^3\le|\g(\nab x^i,\nab x^j) -\de_{ij}\ri| & \les \varep.
  \end{align*}
\end{lemma}
\begin{proof}
  Taking the trace in $\BBB$, we first have
  \begin{align}\label{est:traceB}
    \le|\sum_{i=1}^3|\nab x^i|^2-3\ri| & \les \varep.
  \end{align}
  Assume that $\nab x^1$ has the maximal norm among the $\nab x^i$. Contracting $\BBB$ with $\nab x^1$ gives
  \begin{align*}
    \le||\nab x^1|^4 + |\nab x^1\cdot\nab x^2|^2 + |\nab x^1 \cdot\nab x^3|^2 - |\nab x^1|^2\ri| & \les \varep |\nab x^1|^2.
  \end{align*}
  Dividing by $|\nab x^1|^2>0$, we deduce
  \begin{align}\label{est:nabx1BBB}
     \le||\nab x^1|^2 + |\nab x^1|^{-2}|\nab x^1\cdot\nab x^2|^2 +|\nab x^1|^{-2}|\nab x^1 \cdot\nab x^3|^2 - 1\ri| & \les \varep.
  \end{align}
  We first infer from~\eqref{est:nabx1BBB} that
  \begin{align*}
    |\nab x^1|^2 -1 \les \varep,
  \end{align*}
  which -- since $\nab x^1$ has been chosen to have maximal norm among the $\nab x^i$ -- using~\eqref{est:traceB} further impose that
  \begin{align}\label{est:giideiigood}
    \le||\nab x^i|^2 - 1\ri| & \les \varep, \quad \text{and} \quad \le||\nab x^i|-1\ri| \les \varep,
  \end{align}
  for all $i=1,2,3$. Injecting the above bounds in~\eqref{est:nabx1BBB}, we obtain 
  \begin{align}\label{est:poorgijdeij}
    |\nab x^1\cdot\nab x^2|^2 + |\nab x^1 \cdot\nab x^3|^2 & \les \varep.
  \end{align}
  Contracting $\BBB$ with $\nab x^2\otimes\nab x^3$ gives
  \begin{align*}
    \le|\nab x^2 \cdot\nab x^3 -\sum_{i=1}^3(\nab x^2\cdot\nab x^i)(\nab x^3\cdot\nab x^i)\ri| & \les \varep
  \end{align*}
  which using~\eqref{est:giideiigood} and~\eqref{est:poorgijdeij} gives
  \begin{align*}
    \le|\nab x^2 \cdot\nab x^3 - 2(\nab x^2\cdot\nab x^3)\ri| & \les \varep,
  \end{align*}
  and the result of the lemma follows.
\end{proof}

\section{The $x^i$ are global coordinates from $\Si$ onto $\mathbb{D}$}
We want to improve the result of Section~\ref{sec:localcoordsxi} and show that the map
\begin{align*}
  \Phi~:~x \mapsto (x^i(x))
\end{align*}
is a global diffeomorphism from $\Si$ onto $\mathbb{D}$. This follows from the following two lemmas.

\begin{lemma}\label{lem:glob1}
  We have
\begin{align*}
  \Phi\left(\Si\right) = \mathbb{D}.
\end{align*}
\end{lemma}
\begin{lemma}\label{lem:glob2}
  We have
\begin{align*}
  \forall p\in\Phi\left(\Si\right),~\#\left\{\Phi^{-1}(p)\right\} = 1.
\end{align*}
\end{lemma}

\begin{proof}[Proof of Lemma~\ref{lem:glob1}]
  Let first prove that $\Phi(\Si) \subset  \mathbb{D}$. We argue by contradiction and suppose that there exists $p\in\Si$ such that $|\Phi(p)|>1$.
  Let define the function $X$ on $\Si$ by
  \begin{align*}
    X := \frac{1}{|\Phi(p)|}\sum_{i=1}^n x^i(p) x^i.
  \end{align*}
  From the definitions of the $x^i$, the function $X$ is harmonic, $|X|\leq 1$ on $\pr\Si$ and $|X(p)| = |\Phi(p)|>1$, which contradicts the maximum principle.\\

  We then show that $\Phi(\Si) = \mathbb{D}$. Since $\Phi(\Si)$ is closed in $\mathbb{D}$ and $\mathbb{D}$ is connected, the result will follow provided that we can prove that $\Phi(\Si)$ is an open subset of $\mathbb{D}$. Applying the maximum principle, we can obtain that $\Phi^{-1}\left(\pr\mathbb{D}\right) = \pr\Si$. Applying the local inverse theorem -- since $\Phi$ is a local diffeomorphism -- (see also Remark~\ref{rem:localinverseboundary}) then ensures that $\Phi\left(\Si\right)$ is open in $\mathbb{D}$. This finishes the proof of Lemma~\ref{lem:glob1}.
\end{proof}

\begin{proof}[Proof of Lemma~\ref{lem:glob2}]
  By the local inverse theorem applied at all points $p\in\Phi\left(\Si\right)$, one has that
  \begin{align*}
    C(p):=\#\left\{\Phi^{-1}(p)\right\}
  \end{align*}
  is locally constant at each point $p\in\Phi\left(\Si\right)$ (see~\cite{Mil97} and Remark~\ref{rem:mil}). Since $\Phi\left(\Si\right)$ is connected, this implies that $C$ is constant on $\Phi\left(\Si\right)$.\\

  It therefore suffices to compute $C$ at $(1,0,\cdots,0)\in\Phi\left(\Si\right)$. By the maximum principle, one easily obtains that $\Phi^{-1}(1,0,\cdots,0) = \{(1,0,\cdots,0)\}$, which finishes the proof of Lemma~\ref{lem:glob2}.
\end{proof}
\begin{remark}\label{rem:localinverseboundary}
  If one is concerned with the applicability of the local inverse theorem at the boundary, one can just \emph{extend} the smooth function by Borel's lemma and then apply the classical local inverse theorem.
\end{remark}
\begin{remark}\label{rem:mil}
  The proof in~\cite{Mil97} is for manifold without boundaries, but extends directly to the case with boundary by considering the map $\Phi$ corestricted to its image (which is a manifold with boundary). Note that this is the reason why we apply this result only for points $p\in\Phi\left(\Si\right)$.
\end{remark}

\chapter{Axis limits}\label{app:vertexlimits}
This section is dedicated to the proof of Theorem~\ref{thm:vertex}. In Section~\ref{sec:nullnormalcoords} we define the Cartesian optical normal coordinates, and show that they form a local (smooth outside of the axis) coordinate system. We then perform a Cartesian to spherical coordinate change of coordinates and derive limits at the axis for the metric components and its derivatives. These limits are the key to obtain the limits for the metric, null connection, curvature tensors in Section~\ref{sec:limitsnullconn} using the expression of these tensors in spherical coordinates.

\begin{remark}
  For simplicity -- and since no confusion is possible--, we drop all the primes from the notations used in Theorem~\ref{thm:vertex} in this Appendix.
\end{remark}

\section{Optical normal coordinates}\label{sec:nullnormalcoords}
Let $(\MM,\g)$ be a smooth Lorentzian manifold. Let $\o$ be a timelike geodesic of $\MM$. We assume that $\o$ is parametrised by its geodesic affine parameter, \emph{i.e.} $\g(\doto,\doto) = -1$ and that a parallel transported orthonormal frame $(e_0:=\doto,e_1,e_2,e_3)$ along $\o$ is given. We assume that the frame $e_\mu$ is smoothly extended outside of $\o$, and that its extension satisfies $\D (e_\mu)|_{\o} = 0$.\footnote{Using normal coordinates vectorfields, it is always possible to obtain such an extension.} We define the map $\Psi~:~\MM \times \RRR^4 \to \MM$, by
\begin{align*}
  \Psi\le(p,\pdot^\nu\ri) & := \exp_{p}\le(\pdot^\nu e_\nu(p)\ri).
\end{align*}

Let $O\in\o$. There exists a neighbourhood of $O$ covered by smooth normal coordinates which we call $(z^\mu)$, such that $O = \{z=0\}$ and $\pr_{z^\mu}|_{O} = e_{\mu}(O)$. In this coordinate system the metric and its first derivatives are trivial at $O$ (see~\cite{Gal.Hul.Laf90}), \emph{i.e.}
\begin{align*}
  \g_{\mu\nu} (O) & = \etabold_{\mu\nu}, \\
  \pr\g_{\mu\nu} (O) & = 0,
\end{align*}
where $\etabold_{\mu\nu}$ denotes the Cartesian coordinates components of the Minkowskian metric.\\


We define the map $\Psibar:=\Psibar_O~:~\RRR^4\times\RRR^4 \to \RRR^4$ by identifying $\MM$ to (a neighbourhood of $0$ in) $\RRR^4$ using the normal coordinates $(z^\mu)$ at $O$.

\begin{remark}
  In the Minkowskian case $\MM = \RRR^4$, we have $\Psibar^{\mathrm{Mink}}(z,\zdot) = z + \zdot$.
\end{remark}

The map $\Psibar$ satisfies the following properties.
\begin{lemma}\label{lem:propPsi}
  The map $\Psibar$ is a smooth map from $\RRR^4\times\RRR^4\to\RRR^4$. Moreover, the following identities hold 
  \begin{align}\label{eq:propPsi1}
    \begin{aligned}
    \Psibar|_{\RRR^4\times\{0\}} & = \Id_{\RRR^4},\\
    \Psibar|_{\{0\}\times\RRR^4} & = \Id_{\RRR^4},
    \end{aligned}
  \end{align}
  and
  \begin{align}\label{eq:propPsi2}
    \begin{aligned}
    \d\Psibar(z,0)\cdot (v,0) & = v,\\
    \d\Psibar(0,\zdot)\cdot (0,\vdot) & = \vdot,
    \end{aligned}
  \end{align}
  and
  \begin{align}\label{eq:propPsi2bis}
    \begin{aligned}
      \d\Psibar(z,0) \cdot (0,\vdot) & = \vdot^\nu e_\nu(z). 
    \end{aligned}
  \end{align}
  and
  \begin{align}\label{eq:propPsi3}
    \begin{aligned}
      \d^2\Psibar(0,0) & = 0,
    \end{aligned}
  \end{align}
  as a bilinear map from $\le(\RRR^4\times\RRR^4\ri)^2 \to \RRR^4$.
\end{lemma}
\begin{proof}
  The smoothness of $\Psibar$ is a consequence of the smoothness of the exponential map (see~\cite[p. 86]{Gal.Hul.Laf90}) and of the vectorfields $e_\mu$.\\

  Identities~\eqref{eq:propPsi1} are straight-forward consequences of the definition of normal coordinates and the exponential map. Identities~\eqref{eq:propPsi2} then follow from identities~\eqref{eq:propPsi1}.\\

  Identity~\eqref{eq:propPsi2bis} is a consequence of the definition of the exponential map.\\

  To derive the last identity~\eqref{eq:propPsi3} we first rewrite \eqref{eq:propPsi2bis} as
  \begin{align*}
    \pr_{\zdot^\nu}\Psibar^\la(z,0) & = (e_{\nu}(z))^\la,
  \end{align*}
  which deriving along $z^\mu$ and evaluating at $0$ gives
  \begin{align*}
    \pr_{z^\mu}\pr_{\zdot^\nu}\Psibar^\la(0,0) & = \le(\pr_{z^\mu}(e_\nu(z))^\la\ri)|_{z=0} \\
                                               & = \D_{\pr_{z^\mu}}(e_\nu)^\la|_{z=0} + \Ga^\la_{\mu\de} (e_\nu(z))^\de|_{z=0} \\
                                               & = 0,
  \end{align*}
  by definition of the extension of $e_\nu$ and since the Christoffel symbols $\Ga$ for the normal coordinates vanish at $z=0$.\\

  Deriving the identities of~\eqref{eq:propPsi2} also gives
  \begin{align*}
    \pr_{z^\mu}\pr_{z^\nu}\Psibar^\la(0,0) & = 0,\\
    \pr_{\zdot^\mu}\pr_{\zdot^\nu} \Psibar^\la(0,0) & = 0.
  \end{align*}

  Using Schwarz theorem, since $\Psibar$ is smooth, all second order partial derivatives vanish at $(0,0)$ and identity~\eqref{eq:propPsi3} follows. This finishes the proof of the lemma. 
\end{proof}

We define the map $\Phi~:~\RRR^4\to\RRR^4\times\RRR^4$ by
\begin{align*}
  \Phi\le(x^0,x^1,x^2,x^3\ri) & := \le[\le(x^0+\sqrt{\sum_{i=1}^3\le(x^i\ri)^2},0,0,0\ri),\le(-\sqrt{\sum_{i=1}^3\le(x^i\ri)^2},x^1,x^2,x^3\ri)\ri],
\end{align*}
and we define the map $\Th~:~\RRR^4\to\RRR^4$ by
\begin{align*}
  \Th & := \Psibar \circ \Phi.
\end{align*}

\begin{remark}
  In the Minkowskian case, we have
  \begin{align*}
    \Th^{\mathrm{Mink}} = \Id_{\RRR^4}. 
  \end{align*}
  However, for a general Lorentzian manifold, the map $\Th$ is not smooth at the axis $\{x^i=0\}$. See Lemma~\ref{lem:Thglobdiffeo}.
\end{remark}

\begin{lemma}\label{lem:Thglobdiffeo}
  The map $\Th$ is a $\mathscr{C}^{1,1}$ diffeomorphism from a neighbourhood $\UU$ of $0\in\RRR^4$ onto a neighbourhood $\VV$ of $0\in\RRR^4$. Moreover $\Th$ is smooth on $\UU\setminus\{x^i=0\}$. We call \emph{(Cartesian) optical normal coordinates} the induced local coordinates $(x^{\mu})$ by $\Th$ in $\MM$.
\end{lemma}
\begin{proof}
  From the explicit definition of $\Phi$ and the smoothness of $\Psibar$ obtained in Lemma~\ref{lem:propPsi}, the map $\Th$ is smooth on $\RRR^4\setminus \{x^i=0\}$.\\

  To obtain that $\Th$ is a $\mathscr{C}^{1,1}$ map at $\{x^i=0\}$, and also that $\d\Th((0,x^i)) = \Id_{\RRR^4}$, it is enough to show that
  \begin{align*}
    \norm{\d\Th(x)-\Id_{\RRR^4}} & =  O\le(\sqrt{\sum_{i=1}^3\le(x^i\ri)^2}\ri),
  \end{align*}
  when $\sqrt{\sum_{i=1}^3\le(x^i\ri)^2} \to 0$.\\

  From the chain rule we have
  \begin{align*}
    \d\Th(x)-\Id_{\RRR^4} & = \d\Psibar(\Phi(x))\cdot\d\Phi(x)-\Id_{\RRR^4} \\
                         & = \le(\d\Psibar(\Phi(x)) - \d\Psibar\le[\le(x^0,0,0,0\ri),0\ri]\ri)\cdot\d\Phi(x) \\
                         & + \quad \d\Psibar\le[\le(x^0,0,0,0\ri),0\ri]\cdot\d\Phi(x)- \Id_{\RRR^4}. 
  \end{align*}

  From properties (\ref{eq:propPsi2}) and (\ref{eq:propPsi2bis}) of Lemma~\ref{lem:propPsi}, and the explicit expression of $\Phi$, we have
  \begin{align*}
    & \d\Psibar\le[\le(x^0,0,0,0\ri),0\ri]\cdot\d\Phi(x)\cdot h \\
    = & \; \d\Psibar\le[\le(x^0,0,0,0\ri),0\ri] \cdot \le[\le(h^0 + \le(\sum_{i=1}^3(x^i)^2\ri)^{-1/2} \sum_{j=1}^3x^jh^j,0,0,0\ri),\le(- \le(\sum_{i=1}^3(x^i)^2\ri)^{-1/2} \sum_{j=1}^3x^jh^j,h^1,h^2,h^3\ri)\ri] \\
    = & \; \Id_{\RRR^4}\cdot h.
  \end{align*}

  From the explicit expression of $\Phi$, we have
  \begin{align*}
    \Phi(x)-\le[\le(x^0,0,0,0\ri),0\ri] & = O\le(\sqrt{\sum_{i=1}^3(x^i)^2}\ri), \\
    \norm{\d\Phi(x)} & = O\le(1\ri),
  \end{align*}
  when $\sqrt{\sum_{i=1}^3(x^i)^2} \to 0$. From the smoothness of $\Psibar$, we thus deduce
  \begin{align*}
    \norm{\le(\d\Psibar(\Phi(x)) - \d\Psibar\le[(x^0,0,0,0),0\ri]\ri)\cdot\d\Phi(x)} & = O\le(\sqrt{\sum_{i=1}^3\le(x^i\ri)^2}\ri),
  \end{align*}
  as desired.\\
  
  We thus have obtained that $\Th$ is a $\mathscr{C}^{1,1}$ map, and that $\d\Th(x)$ is invertible in a neighbourhood $\UU$ of $0\in\RRR^4$ (and more generally of $\{x^i=0\}$). The result of the lemma then follows from an application of the local inverse theorem.
\end{proof}

The following lemma is a direct consequence of the definition of $\Th$.
\begin{lemma}\label{lem:opticalnormalnullcone}
  The level sets of the function
  \begin{align*}
    x^0 + \sqrt{\sum_{i=1}^3\le(x^i\ri)^2}
  \end{align*}
  are the ingoing null cones $\CCb$ emanating from $\o$. The level sets of the function
  \begin{align*}
    x^0 - \sqrt{\sum_{i=1}^3\le(x^i\ri)^2}
  \end{align*}
  intersect $\CCb$ at the $2$-spheres of the geodesic foliation on $\CCb$.
\end{lemma}

\begin{lemma}
  The coordinate functions $x^0+t$ and $x^i$ are independent of the choice of point $O = \o(t)$. In the following, we redefine $x^0:=x^0+t$.
\end{lemma}
\begin{proof}
  Let $p = \le(x^0_1,x^i_1\ri)_{O(t_1)} = \le(x^0_2,x^i_2\ri)_{O(t_2)}$. Since $p$ belongs to the same cone emanating from $\o$, we have
  \begin{align*}
    \o\le(t_1+x^0_1 -\sqrt{\sum_{j=1}^3\le(x_1^j\ri)^2}\ri) & = \o\le(t_2+x^0_2 -\sqrt{\sum_{j=1}^3\le(x_2^j\ri)^2}\ri).
  \end{align*}
  By injectivity of the exponential map and since the $e_\mu$ are independent, we further have
  \begin{align*}
    -\sqrt{\sum_{j=1}^3\le(x_1^j\ri)^2} & = -\sqrt{\sum_{j=1}^3\le(x_2^j\ri)^2},\\
    x^i_1 & = x^i_2,
  \end{align*}
  for all $i=1,2,3$, from which we also infer $t_1+x^0_1 = t_2 + x^0_2$.
\end{proof}


We define the \emph{spherical} optical normal coordinates $\le(\ub,u,\varth,\varphi\ri)$ on $\MM$, to be
\begin{align*}
  \le(\ub,u,\varth,\varphi\ri) & := \Xi^{-1}\le(x^0,x^1,x^2,x^3\ri),
\end{align*}
where $\Xi:\RRR_\ub\times\RRR_u\times(0,\pi)_\varth\times(0,2\pi)_\varphi\to\RRR^4$ is the following spherical-to-Cartesian coordinates map
\begin{align*}
  \Xi(\ub,u,\varth,\varphi) & := \le(\frac{\ub+u}{2},\frac{\ub-u}{2}\sin\varth\cos\varphi, \frac{\ub-u}{2}\sin\varth\sin\varphi, \frac{\ub-u}{2}\cos\varth\ri).
\end{align*}

\section{Axis limits for the metric $\g$ in optical normal coordinates}\label{sec:vertexlimitsg}
Let $O\in\o$. By invariance in $t$, we shall assume that $O = \o(t=0)$. In this section, we use the classical normal coordinates $z$ at $O$ defined in Section~\ref{sec:nullnormalcoords} and the change of coordinates $\Th$ from the classical $z$ to the optical $x$ normal coordinates to derive limits at $O$ for (derivatives of) the optical normal coordinates components of the metric $\g$ (see Lemma~\ref{lem:gxgzCart}). Since it will be shown that the metric and its first derivatives are trivial at $O$ and therefore everywhere on the axis $\o$, we will additionally infer strong limits for derivatives of the metric in the axis direction (see Lemma~\ref{lem:limitsstrongaxis}).\\

We have the following limits when $|x|\to 0$ for the $(\ub,u,\varth,\varphi)$-derivatives of $\Th$.
\begin{lemma}\label{lem:limitsThXi}
  We have
  \begin{align}
    \le(\pr_\ub^k,\pr_u^l,\pr_a^m\ri)\Th\circ\Xi(\ub,u,\varth,\varphi) & = O\le(|x|^{|1-k-l|}\ri), \label{est:limitsThXi} \\
    \le(\pr_\ub^k,\pr_u^l,\pr_a^m\ri) \le((\d\Th)\circ\Xi(\ub,u,\varth,\varphi)-\Id_{\RRR^4}\ri) & = O\le(|x|^{|2-k-l|}\ri) \label{est:limitsdThXi}
  \end{align}
  when $|x| \to 0$, for all $k,l,m\geq 0$, and where here and in the following
  \begin{align*}
    \le(\pr_\ub^k,\pr_u^l,\pr^m_a\ri)
  \end{align*}
  denotes all combinations of partial derivatives containing respectively $k,l,m$-derivatives of $\ub,u$ and $a=\varth,\varphi$.
\end{lemma}
\begin{proof}
  We have
  \begin{align*}
    \Th(\Xi(\ub,u,\varth,\varphi)) & = \Psibar\le[\le(\ub,0,0,0\ri),\le(-\frac{\ub -u}{2}, \frac{\ub-u}{2}\sin\varth\cos\varphi, \frac{\ub-u}{2}\sin\varth\sin\varphi, \frac{\ub-u}{2}\cos\varth\ri)\ri],
  \end{align*}
  and the limits~\eqref{est:limitsThXi} follow from $\Psibar(z,\zdot) = O(z,\zdot)$ and the smoothness of $\Psibar$.\\

  We have
  \begin{align*}
    & \;\le(\d\Th(\Xi(\ub,u,\varth,\varphi))-\Id_{\RRR^4}\ri)\cdot(h^\mu) \\
    = & \; \le(\d\Psibar(\Th(\Xi(\ub,u,\varth,\varphi)))-\d\Psibar(0,0)\ri)\cdot\big[\le(h^0 + \sin\varth\cos\varphi h^1 + \sin\varth\sin\varphi h^2 + \cos\varth h^3,0,0,0\ri),\\
    & \quad\quad\quad\quad\quad\quad\quad\quad\quad\quad\quad\quad\quad\quad\quad\quad\quad\quad \le(-\sin\varth\cos\varphi h^1 -\sin\varth\sin\varphi h^2 -\cos\varth h^3, h^i\ri)\big],
  \end{align*}
  and the second limits follow from $\norm{\d\Psibar(z,\zdot)-\d\Psibar(0,0)} = O((z,\zdot)^2)$, $\norm{\d^2\Psibar(z,\zdot)} = O((z,\zdot))$ which are consequences of Lemma~\ref{lem:propPsi}, the smoothness of $\Psibar$, and the previously obtained~\eqref{est:limitsThXi}. This finishes the proof of the lemma.
\end{proof}

As a consequence of Lemma~\ref{lem:limitsThXi}, we have the following limits for the metric $\g$ in the Cartesian optical normal coordinates.
\begin{lemma}\label{lem:gxgzCart}
  We have
  \begin{align}\label{est:gxgzCart}
    \le(\pr_\ub^k,\pr_u^l,\pr_a^m\ri)\le(\g_{\mu\nu} - \etabold_{\mu\nu}\ri) & = O\le(|x|^{|2-k-l|}\ri),
  \end{align}
  when $|x| \to 0$, for all $k,l,m\geq 0$, and where $\g_{\mu\nu}$ denotes the components of the metric $\g$ in the coordinates $x^\mu$.
\end{lemma}
\begin{proof}
  We have the formula
  \begin{align*}
    \g^x_{\mu\nu}(x) & = \g^z_{\mu'\nu'}\le(\Th(x)\ri)\pr_{x^\mu}\Th^{\mu'}\pr_{x^\nu}\Th^{\nu'},
  \end{align*}
  which we rewrite as
  \begin{align}\label{eq:chgvargxgz}
    \begin{aligned}
      \g^x_{\mu\nu}(x) -\etabold_{\mu\nu} & = \le(\g^{z}_{\mu\nu}(\Th(x))-\etabold_{\mu\nu}\ri) \\
      & \quad + \g^z_{\mu'\nu'}(\Th(x))\le(\pr_{x^\mu}\Th^{\mu'}-\de_\mu^{\mu'}\ri)\pr_{x^\nu}\Th^{\nu'} \\
      & \quad + \g^z_{\mu'\nu'}(\Th(x))\de_\mu^{\mu'}\le(\pr_{x^\nu}\Th^{\nu'}-\de_\nu^{\nu'}\ri).
    \end{aligned}
  \end{align}
  The metric components of the normal coordinates $\g^{z}_{\mu\nu}(z)$ are smooth and we have
  \begin{align}\label{est:limitsnormalcoords}
    \begin{aligned}
    \g^{z}_{\mu\nu}(z)-\etabold_{\mu\nu} & = O(|z|^2),\\
    \pr_z\g^{z}_{\mu\nu}(z) & = O(|z|),\\
    \pr^{\geq 2}_z\g^z_{\mu\nu}(z) & = O(1),
    \end{aligned}
  \end{align}
  when $z \to 0$. The proof of~\eqref{est:gxgzCart} then follows from formula~\eqref{eq:chgvargxgz} and the limits given by~\eqref{est:limitsThXi}, \eqref{est:limitsdThXi} and \eqref{est:limitsnormalcoords}. This finishes the proof of the lemma.
\end{proof}

We have the following stronger limits derivatives along the axis.
\begin{lemma}\label{lem:limitsstrongaxis}
  We have
  \begin{align}\label{est:limitsstrongaxis}
    \le(\pr_u+\pr_\ub\ri)^{n} \le(\pr_\ub^k,\pr_u^l,\pr^m_a\ri)\le(\g_{\mu\nu}-\etabold_{\mu\nu}\ri) & = O\le(|x|^{|2-k-l|}\ri),
  \end{align}
  when $|x|\to 0$ and for all $n\geq 0$ and all $k,l,m\geq 0$.
\end{lemma}
\begin{proof}
  The lemma follows from the smoothness in $(\ub,u,\varth,\varphi)$ and the fact that the limits~\eqref{est:gxgzCart} hold uniformly along the axis $\{u=\ub\}$.
\end{proof}

In the following we combine~\eqref{est:gxgzCart} and~\eqref{est:limitsstrongaxis} into the following equivalent limits
\begin{align}\label{est:limitscombined}
  \le(\pr_u^k,\le(\pr_u+\pr_\ub\ri)^l,\pr_a^m\ri)\le(\g_{\mu\nu} - \etabold_{\mu\nu}\ri) & = O\le(|x|^{|2-k|}\ri),
\end{align}
when $|x|\to 0$ and for all $k,l,m\geq 0$.\\

Since our goal is to obtain vertex limits for the null connection coefficients and rotation vectorfields which are naturally expressed in terms of the metric components in spherical coordinates, we infer the following equivalent of~\eqref{est:limitscombined} for the metric in spherical coordinates. 
\begin{lemma}\label{lem:gtCart}
  We have
  \begin{align}\label{est:gtCart}
    \begin{aligned}
      \le(\pr_u^k,(\pr_u+\pr_\ub)^l,\pr_a^m\ri)\le(\g_{\mubt\nubt}-\etabold_{\mubt\nubt}\ri) & = O\le(|x|^{|2-k|}\ri),
    \end{aligned}
  \end{align}
  when $|x|\to 0$, for all $k,l,m\geq 0$, and where from now on indices $\mu,\nu$ (resp. $a,b$) refer to evaluation of a tensor with respect to the coordinate vectorfields $(\pr_\ub,\pr_u,\pr_\varth,\pr_\varphi)$ (resp. $(\pr_\varth,\pr_\varphi)$), whereas tilde indices $\mubt\nubt$ (resp. $\abt,\bbt$) refer to evaluation with respect to the normalised coordinate vectorfields $(\prt_\ub := \pr_\ub, \prt_u := \pr_u,\prt_\varth := \frac{2}{\ub-u}\pr_{\varth}, \prt_\varphi := \frac{2}{\ub-u}\pr_{\varphi})$ (resp. $(\prt_\varth,\prt_\varphi)$). The tensor $\etabold$ denotes the Minkowski metric, which writes in $(\ub,u,\varth,\varphi)$ coordinates
  \begin{align*}
    \etabold = -\d \ub \d u + \le(\frac{\ub-u}{2}\ri)^2\d\varth^2 + \le(\frac{\ub-u}{2}\ri)^2\sin^2\varth\d\varphi^2. 
  \end{align*}

  We also have
  \begin{align*}
    \le(\pr_u^k,(\pr_u+\pr_\ub)^l,\pr_a^m\ri)\le(\g^{\mubt\nubt}-\etabold^{\mubt\nubt}\ri) & = O\le(|x|^{|2-k|}\ri),
  \end{align*}
  where $\g^{\mubt\nubt} := (\g_{\mubt\nubt})^{-1}$ and also coincides with the renormalisation of $(\g_{\mu\nu})^{-1}$.
\end{lemma}
\begin{proof}
  We have
  \begin{align*}
    \g_{\mubt\nubt} & = (J_1)^{\mu'}_{\mu}\g^x_{\mu'\nu'}(J_1)^{\nu'}_{\nu},
  \end{align*}
  where $J_1$ is the normalised Jacobian matrix of the spherical to Cartesian coordinates change
  \begin{align*}
    J_1 = \le(
    \begin{array}{cccc}
      1/2 & 1/2 & 0 & 0 \\
      1/2\sin\varth\cos\varphi & -1/2\sin\varth\cos\varphi & \cos\varth\cos\varphi & -\sin\varth\cos\varphi \\
      1/2\sin\varth\sin\varphi & - 1/2\sin\varth\sin\varphi & \cos\varth\sin\varphi & \sin\varth\cos\varphi \\
      1/2 \cos\varth & -1/2\cos\varth & -\sin\varth & 0
    \end{array}
          \ri)
  \end{align*}
  and~\eqref{est:gtCart} follows from~\eqref{est:limitscombined}. The limits for the inverse matrix follows by inverting the same formula.
\end{proof}
\begin{corollary}\label{lem:limitssphere}
  We have
  \begin{align*}
    \le(\pr_u^k,(\pr_u+\pr_\ub)^l,\pr_\om^m\ri)\le(\gd_{ab}-\etadt_{ab}\ri) & = O\le(|x|^{|4-k|}\ri),
  \end{align*}
  when $|x| \to 0$, and for all $k,l,m\geq 0$. As a consequence, we have for the area radius $r$ of $S_{u,\ub}$
  \begin{align}\label{est:limitsr}
    \le(\pr_u^k,(\pr_u+\pr_\ub)^l\ri)\le(r(u,\ub)-\frac{\ub-u}{2}\ri) & = O\le(|x|^{|3-k|}\ri).
  \end{align}
\end{corollary}

\begin{remark}
  The limits for the null connection coefficients and the rotation vectorfields are performed in the coordinate chart $(\ub,u,\varth,\varphi)$, in a domain where these chart remains regular, \emph{e.g.} $\pi/10 < \varth < 9\pi/10$. By changing the axis of the spherical coordinates and performing the same computations, one obtains the same limits in the regions $0 \leq \varth \leq \pi/10$, $9\pi/10 \leq \varth \leq \pi$.   
\end{remark}

\begin{lemma}\label{lem:defOScorrec}
  Let $F$ be a $k$-tensor. We say that $F=\OS\le(|x|^\ga\ri)$ if for all $k,l,m\geq 0$, we have
  \begin{align*}
    \sum_{(\mu_1,\cdots,\mu_k)\in\{\ub,u,\varth,\varphi\}}\le|\le(\le((\ub-u)\pr_u\ri)^k,\le(\pr_u+\pr_\ub\ri)^l,\pr_a^m\ri) \le(F_{\mubt_1\cdots \mubt_k}\ri)\ri| &  = O\le(|x|^{\ga}\ri),
  \end{align*}
  when $|x|\to 0$. 
  With this definition, we have
  \begin{subequations}\label{est:limitsOS}
  \begin{align}
      \g -\etabold & = \OS\le(|x|^2\ri),\label{est:limitsOS1}\\
      \Ga(\g)-\Ga(\etabold) & = \OS\le(|x|\ri),\label{est:limitsOS2}\\
      \R & = \OS\le(1\ri),\label{est:limitsOS3}
  \end{align}
  \end{subequations}
  when $|x|\to 0$, where we consider the Christoffel symbols $\Ga$ as $3$-tensors.
\end{lemma}
\begin{proof}
  The limits~\eqref{est:limitsOS1} are direct consequences of~\eqref{est:gtCart}. The limits~\eqref{est:limitsOS2} follow from the definition of the Christoffel symbols and limits~\eqref{est:limitsOS1}. The limits~\eqref{est:limitsOS3} follow from the expression of $\R$ in terms of the (derivatives of) the Christoffel symbols and the limits~\eqref{est:limitsOS2}.
\end{proof}


\begin{remark}
  The limits~\eqref{est:gtCart} are stronger than~\eqref{est:limitsOS} for $\g$. But necessary for $\Ga(\g)$ and $\R$.
\end{remark}

\section{Axis limits for the null connection coefficients}\label{sec:limitsnullconn}
We recall that the null pair $(\elb,\el)$ for the optical normal coordinates is defined by $\elb:=-\D\ub$. By construction of the optical normal spherical coordinates in Section~\ref{sec:nullnormalcoords}, the null pair $(\elb,\el)$ writes
\begin{align}\label{eq:elbelea}
  \elb & = 2\pr_u, & \el & = - \g^{u\nu}\pr_\nu - \half \g^{uu}\pr_u,
\end{align}
which we also rewrite as
\begin{align}\label{eq:elerror}
  \el & = 2\pr_\ub -\le(\g^{u\nu}-\etabold^{u\nu}\ri)\pr_\nu - \half \le(\g^{uu}-\etabold^{uu}\ri)\pr_u.
\end{align}


From the above formulas, we deduce the following expressions of the null connection coefficients associated to the null pair $(\elb,\el)$ in terms of $\g_{\mubt\nubt}$ and $\pr\g_{\mubt\nubt}$.
\begin{lemma}\label{lem:nullconngmunu}
  We have
  \begin{align*}
    \chib_{ab} & = -\frac{2}{\ub-u}\g_{ab} + \quar(\ub-u)^2\pr_u(\g_{\abt\bbt}-\etabold_{\abt\bbt}), \\ \\
    \ze_a & = \pr_a\g_{\ub\ut} + \half\pr_u((\ub-u)\g_{\abt\bbt}) -\half (\g^{\ut\ubt}-\etabold^{\ut\ubt})\le(\pr_a\g_{\ubt\ut} + \half\pr_u((\ub-u)\g_{\abt\ubt})\ri) \\
               & \quad -(\ub-u)^{-1} (\g^{\ut\bbt}-\etabold^{\ubt\bbt}) \le(\quar\pr_u((\ub-u)^2\g_{\abt\bbt})\ri),\\ \\
    \yy & = - \half \le(\g^{\ut\ut}-\etabold^{\ut\ut}\ri), \\ \\
    \chi_{ab} & = \frac{2}{\ub-u} \g_{ab} + \quar(\ub-u)^2\pr_\ub\le(\g_{\abt\bbt}-\etabold_{\abt\bbt}\ri) \\
               & \quad -\quar(\ub-u)\g_{\nubt\bbt}\pr_a\le(\g^{\ut\nubt} - \etabold^{\ut\nubt}\ri) -\quar(\ub-u)\g_{\nubt\abt}\pr_b\le(\g^{\ut\nubt} - \etabold^{\ut\nubt}\ri) \\
               & \quad  - \frac{1}{8} (\g^{\ut\ut}-\etabold^{\ut\ut})\pr_u((\ub-u)^2\g_{\abt\bbt}) \\
               & \quad -\half (\g^{\ut\ubt}-\etabold^{\ut\ubt}) \le(\half\pr_{a}((\ub-u)\g_{\ubt \bbt}) + \quar\pr_\ub((\ub-u)^2\g_{\abt\bbt}) - \half\pr_b((\ub-u)\g_{\abt\ubt})\ri) \\
               & \quad - \frac{1}{4} (\ub-u) (\g^{\ut\cbt}-\etabold^{\ut\cbt})\le(\pr_a\g_{\cbt\bbt} + \pr_c\g_{\abt\bbt} - \pr_b\g_{\abt\cbt}\ri),
  \end{align*}
  together with $\eta = -\etab = \ze$, $\omb = |\xib| = 0$ and $\om = -\quar\elb(\yy), \xi = -\half\Nd\yy$.
\end{lemma}
\begin{proof}
  Using~\eqref{eq:elbelea}, we have
  \begin{align*}
    \chib_{ab} = \chib(\pr_a,\pr_b) = \g(\D_{\pr_a}2\pr_u,\pr_b) = 2\g(\D_{\pr_u}\pr_a,\pr_b) = \g(\D_{\pr_u}\pr_a,\pr_b) + \g(\D_{\pr_u}\pr_b,\pr_a) = \pr_u(\g_{ab}),
  \end{align*}
  where we used the symmetry of $\chib$. We further have
  \begin{align*}
    \pr_u(\g_{ab}) & = \pr_u\le(\le(\frac{\ub-u}{2}\ri)^2\g_{\abt\bbt}\ri),
  \end{align*}
  and the desired formula follows.\\

  Using~\eqref{eq:elbelea}, we have
  \begin{align*}
    \ze_a & = \g\le(\D_{\pr_a}\pr_u, -\g^{u\nu}\pr_\nu-\half \g^{uu}\pr_u\ri) \\
          & = -\g^{u\nu} \g_{\mu\nu} \Ga_{au}^\mu(\g) \\
          & = -\Ga^{u}_{au}(\g).
  \end{align*}
  Using that $\g_{uu}= \g_{au} = 0$, we further have
  \begin{align*}
    -\Ga^u_{au}(\g) & = -\half \g^{u\mu}\le(\pr_a\g_{\mu u} + \pr_u\g_{a\mu} - \pr_{\mu}\g_{a u}\ri) \\
                    & = -\half \g^{\ut\ubt}\le(\pr_a\g_{\ubt\ut} + \half\pr_u((\ub-u)\g_{\abt\ubt})\ri) \\
                    & \quad -(\ub-u)^{-1} \g^{\ut\bbt} \le(\quar\pr_u((\ub-u)^2\g_{\abt\bbt})\ri) \\
                    & = \pr_a\g_{\ub\ut} + \half\pr_u((\ub-u)\g_{\abt\bbt}) -\half (\g^{\ut\ubt}-\etabold^{\ut\ubt})\le(\pr_a\g_{\ubt\ut} + \half\pr_u((\ub-u)\g_{\abt\ubt})\ri) \\
                    & \quad -(\ub-u)^{-1} (\g^{\ut\bbt}-\etabold^{\ubt\bbt}) \le(\quar\pr_u((\ub-u)^2\g_{\abt\bbt})\ri),
  \end{align*}
  as desired.\\

  By definition of $\yy$, we have
  \begin{align*}
    \yy & = -\half \g^{\mu\nu}\pr_{\mu}u\pr_{\nu}u = -\half \g^{uu}= -\half (\g^{\ut\ut}-\etabold^{\ut\ut}),
  \end{align*}
  as desired.\\
  
  Using~\eqref{eq:elbelea} and~\eqref{eq:elerror}, we have
  \begin{align*}
    \chi_{ab} & = 2\g(\D_{\pr_a}\pr_\ub,\pr_b) + \g\le(\D_{\pr_a}\le(-(\g^{u\nu}-\etabold^{u\nu})\pr_\nu -\half \le(\g^{uu}-\etabold^{uu}\ri)\pr_u\ri),\pr_b\ri)\\
              & = 2\g(\D_{\pr_a}\pr_\ub,\pr_b) - \g_{\nu b}\pr_a\le((\g-\etabold)^{u\nu}\ri) - \g_{\mu b}(\g^{u\nu}-\etabold^{u\nu})\Ga_{a\nu}^\mu(\g) - \quar \le(\g^{uu}-\etabold^{uu}\ri)\chib_{ab},
  \end{align*}
  which by symmetrising in $a,b$ gives
  \begin{align*}
    \chi_{ab} & = \g(\D_{\pr_a}\pr_\ub,\pr_b) + \g(\D_{\pr_b}\pr_\ub,\pr_a) \\
              & \quad - \half\g_{\nu b}\pr_a\le((\g-\etabold)^{u\nu}\ri) -\half\g_{\nu a}\pr_b\le((\g-\etabold)^{u\nu}\ri) \\
              & \quad - \half \g_{\mu b}(\g^{u\nu}-\etabold^{u\nu})\Ga_{a\nu}^\mu(\g) - \half \g_{\mu a}(\g^{u\nu}-\etabold^{u\nu})\Ga_{b\nu}^\mu(\g) \\
              & \quad - \quar \le(\g^{uu}-\etabold^{uu}\ri)\chib_{ab}.
  \end{align*}
  We further have
  \begin{align*}
    \pr_\ub(\g_{ab}) & = \frac{2}{\ub-u}\g_{ab} + \quar(\ub-u)^2\pr_\ub(\g_{\abt\bbt}-\etabold_{\abt\bbt}),\\
    - \half\g_{\nu b}\pr_a\le((\g-\etabold)^{u\nu}\ri) & = -\quar(\ub-u)\g_{\nubt\bbt}\pr_a\le(\g^{\ut\nubt} - \etabold^{\ut\nubt}\ri),\\
    - \half \g_{\mu b}(\g^{u\nu}-\etabold^{u\nu})\Ga_{a\nu}^\mu(\g) & = - \quar (\g^{u\nu}-\etabold^{u\nu}) \le(\pr_{a}\g_{\nu b} + \pr_\nu\g_{ab} - \pr_b\g_{a\nu}\ri) \\
                                                                & = - \quar (\g^{\ut\ut}-\etabold^{\ut\ut})\pr_u(\g_{ab}) \\
                                                                & \quad -\quar (\g^{\ut\ubt}-\etabold^{\ut\ubt}) \le(\half\pr_{a}((\ub-u)\g_{\ubt \bbt}) + \quar\pr_\ub((\ub-u)^2\g_{\abt\bbt}) - \half\pr_b((\ub-u)\g_{\abt\ubt})\ri) \\
                                                                & \quad - \frac{1}{8} (\ub-u) (\g^{\ut\cbt}-\etabold^{\ut\cbt})\le(\pr_a\g_{\cbt\bbt} + \pr_c\g_{\abt\bbt} - \pr_b\g_{\abt\cbt}\ri),
  \end{align*}
  and by symmetrising in $a,b$, the desired formula follows. This finishes the proof of the lemma. 
\end{proof}



The following lemma directly follows from the limits from Lemma~\ref{lem:gtCart} and from the formulas from Lemma~\ref{lem:nullconngmunu}.
\begin{lemma}
  We have
  \begin{align}\label{est:limitprchib}
    \begin{aligned}
      \chib + \frac{2}{\ub-u}\gd,~\ze,~\chi-\frac{2}{\ub-u} \gd & = \OS\le(|x|\ri), \\
      \yy & = \OS\le(|x|^{2}\ri),
    \end{aligned}
  \end{align}
  when $|x|\to 0$, where we refer to Lemma~\ref{lem:defOScorrec} for notations.
\end{lemma}

\begin{lemma}\label{lem:FOSNdF}
  Let $F$ be an $S_{u,\ub}$-tangent tensor such that
  \begin{align*}
    F = \OS(|x|^\ga),
  \end{align*}
  when $|x|\to 0$. Then, the following limits hold
  \begin{align*}
    \le|\le((\ub-u)\Nd_3\ri)^k,\le(\Nd_3+\Nd_4\ri)^l,\le((\ub-u)\Nd\ri)^m F\ri| & = O(|x|^\ga),
  \end{align*}
  when $|x|\to 0$ and for all $k,l,m\geq 0$.
\end{lemma}
\begin{proof}
We have the following formulas relating covariant to partial derivatives
\begin{align*}
  \Nd_aF_b & = \pr_a\le(F_b\ri) - \Ga_{ab}^c(\gd) F_c,\\
  (\ub-u)\Nd_3 F_{a} & = 2 (\ub-u) \pr_u\le(F_{a}\ri) -(\ub-u)\chib_{a}^bF_b,\\
  \le(\Nd_3 + \Nd_4\ri)F_a & = 2(\pr_u+\pr_\ub)(F_a) -(\chi+\chib)_a^bF_b\\  
                    & \quad -\le(\g^{u\nu}-\etabold^{u\nu}\ri)\pr_\nu\le(F_{a}\ri) - \half \le(\g^{\ut\ut}-\etabold^{\ut\ut}\ri)\pr_u\le(F_{a}\ri).
\end{align*}
Using~(\ref{est:limitsOS}), (\ref{est:limitprchib}), we obtain that
\begin{align*}
  (\ub-u)\Nd F,~(\ub-u)\Nd_3F,~(\Nd_3+\Nd_4)F & = \OS(|x|^\ga),
\end{align*}
when $|x|\to 0$, and the result follows by iteration.
\end{proof}

From the limits~\eqref{est:limitprchib} and Lemma~\ref{lem:FOSNdF}, we deduce the following lemma. 
\begin{lemma}\label{lem:limitsNdGa}
  We have
  \begin{align}\label{est:limitsNdGa1}
    \begin{aligned}
      \le|\le((\ub-u)\Nd_3\ri)^k,\le(\Nd_3+\Nd_4\ri)^l,\le((\ub-u)\Nd\ri)^m \le(\chib+\frac{2}{\ub-u}\gd,~\ze,~\chi-\frac{2}{\ub-u}\gd\ri)\ri| & = O\le(|x|\ri), \\
      \le|\le((\ub-u)\Nd_3\ri)^k,\le(\Nd_3+\Nd_4\ri)^l,\le((\ub-u)\Nd\ri)^m \yy\ri| & = O\le(|x|^2\ri),
    \end{aligned}
  \end{align}
  and
  \begin{align}\label{est:limitsNdGa2}
    \le|\le((\ub-u)\Nd_3\ri)^k,\le(\Nd_3+\Nd_4\ri)^l,\le((\ub-u)\Nd\ri)^m \le(\om,~\xi\ri)\ri| & = O\le(|x|\ri),
  \end{align}
  when $|x|\to 0$, and for all $k,l,m\geq 0$. 
\end{lemma}
\begin{proof}
  The limits~\eqref{est:limitsNdGa1} follow directly from~\eqref{est:limitprchib} and Lemma~\ref{lem:FOSNdF}. The limits~\eqref{est:limitsNdGa2} are obtained from the limits for $\yy$ and the relations of Lemma~\ref{lem:nullconngmunu}. 
\end{proof}

We also have the following limits for the null decomposition of the spacetime curvature tensor.
\begin{lemma}\label{lem:limitsnullcurv}
  We have
  \begin{align}\label{est:limitsR}
    \le|\le(\le((\ub-u)\Nd_3\ri)^k,\le(\Nd_3+\Nd_4\ri)^l,((\ub-u)\Nd)^m\ri) \le(\al,\be,\rho,\sigma,\beb,\alb\ri)\ri| & = O\le(1\ri),
  \end{align}
  when $|x|\to 0$, for all $k,l,m\geq 0$.
\end{lemma}
\begin{proof}
  By definition of the null decomposition of $\R$ and using the expressions~(\ref{eq:elbelea}) for $\el,\elb$ in terms of the coordinate vectorfields, we have
  \begin{align*}
    \al_{ab} & = \R(\el,\pr_a,\el,\pr_b) \\
             & = \g^{u\nu} \g^{u\nu'} R_{\nu a \nu' b} + \half \g^{u\nu} R_{\nu a u b} + \half \g^{u\nu'} R_{u a \nu' b} + \quar \g^{uu} R_{uaub}. 
  \end{align*}
  We therefore deduce from the limits~\eqref{est:limitsOS} for $\g$ and $\R$ that
  \begin{align*}
    \al & = \OS(1),
  \end{align*}
  and using the result of Lemma~\ref{lem:FOSNdF} the result follows. The other components of the null decomposition of $\R$ are treated similarly.
\end{proof}

\begin{remark}
  Using~\eqref{est:limitsr}, one can replace $(\ub-u)/2$ by the area radius $r$ of the $2$-spheres $S_{u,\ub}$ in the limits~\eqref{est:limitsNdGa1}, \eqref{est:limitsNdGa2} and \eqref{est:limitsR}.
\end{remark}

\chapter{The canonical foliation}\label{app:canlocalex}
In this section we prove Theorem~\ref{thm:canonical}. We follow the scheme of the proof of the local existence result for the canonical foliation in~\cite[Section 6]{Czi.Gra19}. That is: we first associate the canonical foliation choice with respect to the geodesic foliation to a system of coupled transport-elliptic equations, we then perform a Banach-Picard iteration to obtain local existence for this system. In the present case, initial data for the system of transport-elliptic equations are prescribed at the vertex of a cone. Closing the Banach-Picard iteration scheme in this degenerate situation requires a strongly different analysis as in the non-degenerate case of~\cite{Czi.Gra19}.

\section{The system of transport-elliptic equations}\label{sec:transelldefgeom}
By translation invariance, we assume that $\uba = 0$ and we note $z^\mu$ the normal coordinates at $\o(0)$ such that $z^\mu|_{\o(0)} = 0$ for all $\mu=0,1,2,3$ and $\pr_{z^0} = \doto(0)$. The normal coordinates are well defined in a coordinate ball of radius $\varep>0$ and for all fixed integer $k\geq 5$, we have
\begin{align}\label{est:smallnessnormalCCba}
  \le|\pr^{\leq k+1}\le(\g_{\mu\nu}-\etabold_{\mu\nu}\ri)\ri| & \les 1.
\end{align}
We denote by $s := z^0 - \sqrt{\sum_{i=1}^3(z^i)^2}$ the geodesic affine parameter on $\CCba = \le\{ z^0 + \sqrt{\sum_{i=1}^3(z^i)^2} = 0\ri\}$ and by $(\elb',\el')$ its associated null pair (such that $\elb'(s) = 2$).\\

Our goal is to obtain a function $(u,\om) \in [0,\de]\times \SSS^2 \mapsto (s(u,\om),\om) \in \CCba$ such that $\{s=s(u)\}$ defines the level sets of the canonical foliation of Theorem~\ref{thm:canonical} on $\CCba$ (see the detailed geometric set up of~\cite[Section 6]{Czi.Gra19}).\\

For $u$ a function on $\CCba$, we define the null lapse $\Om$ as 
\begin{align*}
  \Om^{-1} & := \pr_su = \half\elb'(u),
\end{align*}
Assume that $u$ defines a regular foliation on $\CCba$ and let $(\la,f,\fb)$ be the transition coefficients between the null frames $(\elb,\el)$ and $(\elb',\el')$ as defined in Section~\ref{sec:defchgframe}. From the relations of Lemma~\ref{lem:transframe}, we have
\begin{align}\label{eq:reltransOm}
  \la & = \Om, & f & = -\Om^{-1}\Nd s, &  \fb & = 0,
\end{align}
where we recall that $\la,f,\fb$ are $S_u$-tangent tensors. From a simple exact recomputation of the relations of Proposition~\ref{prop:transconn} under the relations~\eqref{eq:reltransOm}, we have 
\begin{align}\label{eq:zeze'}
  \begin{aligned}
  \ze_a & = \half \g(\D_\ea\el,\elb) \\
      & = \half \Om \g(\D_{\ea}(\Om^{-1}\el'-f - \quar|f|^2\elb) ,\elb')\\
      & = \ze'_a - \Om \ea(\Om^{-1}) -\half\Om\g(\D_\ea f,\elb') \\
      & = (\ze')^{\dg}_a +\Nd_a(\log\Om) + \half \Om f^\ddg\cdot\chib'.
    \end{aligned}
\end{align}

From a simple exact recomputation of the relations of Proposition~\ref{prop:chgframederiv} under the relations~\eqref{eq:reltransOm}, using a normal frame $(\ea)_{a=1,2}$ at the point of computation as in the proof of Proposition~\ref{prop:chgframederiv}, we have
\begin{align}\label{eq:recomputzeze'}
  \begin{aligned}
    \Divd(\ze')^\dg & = \gd^{ab}\Nd_a(\ze')^\dg_b \\
    & = \gd^{ab}\ea\le(\ze'_b\ri) \\
    & = \gd^{ab}\le(\ea'-\half f_a \Om\elb'\ri)\le(\ze'_b\ri) \\
    & = \Divd'\ze' -\half \Om f\cdot \Nd'_3\ze' + \gd^{ab}\ze'_c\g(\D_{\ea'}\eb',\ec') -\half \Om \gd^{ab} \ze'_cf_a \g(\D_{\elb'}\eb',\ec') \\
    & = \Divd'\ze'-\half\Om f^\ddg\cdot\Nd'_3\ze' + \half \Om \gd^{ab}\ze'_c f_b\chib'_{ac} \\
    & = \Divd'\ze' -\half \Om f^\ddg \cdot\Nd'_3\ze' + \half \Om \ze'\cdot f^\ddg \cdot \chib'.
  \end{aligned}
\end{align}

Using the same relations, we have
\begin{align}\label{eq:Divdchibchib'}
  \begin{aligned}
    \Divd(\chib')^\dg_a & = \gd^{bc}\Nd_b(\chib'_{ca}) \\
    & = \Divd'\chib' -\half \gd^{bc}\Om f_b \Nd'_3\chib'_{ca} + \gd^{bc}\chib'_{da}\g(\Nd'_b\ec',e'_d) + \gd^{bc}\chib'_{cd}\g(\Nd'_b\ea',e_d') \\
    & = \Divd'\chib' - \half \Om f^\ddg \cdot \Nd_3'\chib' + \half \gd^{bc}\chib'_{da}\Om f_c\chib'_{bd} + \half \gd^{bc}\chib'_{cd}\Om f_a \chib'_{bd} \\
    & = \Divd'\chib' -\half \Om f^\ddg\cdot\Nd_3'\chib' + \half \Om (f^\ddg\cdot\chib'\cdot\chib')_a + \half \Om f_a |\chib'|^2.
  \end{aligned}
\end{align}

Again, from a simple exact recomputation of the relation of Proposition~\ref{prop:transcurv} using relations~\eqref{eq:reltransOm}, we have
\begin{align}\label{eq:recomputrhorho'}
  \begin{aligned}
    \rho & = \quar\R(\el,\elb,\el,\elb) \\
    & = \quar \Om^{2}\R\le(\Om^{-1}\el'-f-\quar|f|^2\elb,\elb',\Om^{-1}\el'-f-\quar |f|^2\elb,\elb'\ri) \\
    & = \rho' - \Om f^\ddg \cdot \beb' + \quar \Om^2 f^\ddg\cdot f^\ddg \cdot\alb'.
  \end{aligned}
\end{align}

Thus, using that $\Om f = -\Nd s$, we deduce from equations~\eqref{eq:zeze'},~\eqref{eq:recomputzeze'},~\eqref{eq:Divdchibchib'}, and~\eqref{eq:recomputrhorho'} that
\begin{align}\label{eq:DivdzerhoFi}
  \Divd\ze + \rho & = \Ld(\log\Om) - s^{-1}\Ld(s) - F_1(s)  -F_2(s)\cdot\Nd s - F_3(s)\cdot \Nd s \cdot \Nd s - F_4(s)\cdot\Nd^2s,
\end{align}
where
\begin{align*}
  F_1(s) & := -\Divd'\ze' - \rho',\\
  F_2(s) & := -\half \Nd'_3\ze' + \half \ze'\cdot\chib' + \half \Divd'\chib' -\be', \\
  F_3(s) & := \quar \Nd_3'\chib' - \quar \chib'\cdot\chib' - \quar |\chib'|^2\gd - \quar \alb', \\
  F_4(s) & := \half \chib' -s^{-1}\gd = \quar \le(\trchib'-\frac{4}{s}\ri)\gd + \half\chibh'.
\end{align*}

From a simple exact recomputation of the relations of Proposition~\ref{prop:transconn}, we additionally have 
\begin{align}\label{eq:ombomb'}
  \omb & = -2\elb(\log\Om).
\end{align}

Using the initial conditions~(\ref{eq:initvaluecanthm}) for the canonical foliation, and the definition of the geodesic foliation $s$, we have the initial conditions at the vertex $\o(0)$
\begin{align*}
  \Om & = 1, & s & = u = 0.
\end{align*}
Thus, using~\eqref{eq:DivdzerhoFi} and~\eqref{eq:ombomb'} and formula~\eqref{eq:commelbov}, we rewrite the equations~\eqref{eq:ellcanthm} which define the canonical foliation as the following system of quasilinear elliptic-transport equations
\begin{align}\label{eq:systemcancanbis}
  \begin{aligned}
    s & = \int_{0}^u\Om\,\d u',\\
    \Ld_{\gd(s)}(\log\Om) & = s^{-1}\Ld(s) + F\le(s,\Nd s,\Nd^2s\ri)- \overline{s^{-1}\Ld s + F\le(s,\Nd s,\Nd^2s\ri)}, \\
    \overline{\log\Om}(u) & = \int_0^u \overline{(\trchib-\trchibo)\log\Om}(u')\,\d u',
  \end{aligned}
\end{align}
where
\begin{align*}
  F(s,\Nd s,\Nd^2s) & := F_1(s) + F_2(s)\cdot\Nd s + F_3(s)\cdot\Nd s\cdot\Nd s + F_4(s)\cdot\Nd^2s.
\end{align*}

Using~\eqref{est:smallnessnormalCCba} and the vertex limits of Theorem~\ref{thm:vertex}, there exists constants $A_k>0$ such that
\begin{align}\label{est:Fisk}
  \begin{aligned}
    |(s\Nd)^{\leq k}F_1(s)| & \leq A_k, & \le|(s\Nd)^{\leq k}F_2(s)\ri| & \leq A_k,\\
    \le|(s\Nd)^{\leq k}F_3(s)\ri| & \leq A_ks^{-2}, & \le|(s\Nd)^{\leq k}F_4(s)\ri| & \leq A_k s,
    \end{aligned}
\end{align}
where, more specifically for the last estimate, we used that we have for the area radius of the geodesic foliation $r\sim -\half s$ and that $\trchib+\frac{2}{r} = O(s)$. Using~\eqref{est:smallnessnormalCCba} and the vertex limits of Theorem~\ref{thm:vertex}, we moreover have
\begin{align}
  \label{est:trchibtrchiboiterate}
  \le|\trchib'(s)-\trchibo'(s)\ri| & \leq A_k s.
\end{align}



\section{The Banach-Picard iteration}
To prove local existence for the system of transport-elliptic equations~\eqref{eq:systemcancanbis}, we apply a Banach-Picard fixed point argument as in~\cite{Czi.Gra19}. We define
\begin{align*}
  \begin{aligned}
    s_0 & = u, & \Om_0 & = 1,
  \end{aligned}
\end{align*}
and the iteration
\begin{align}\label{eq:syscancaniterate}
  \begin{aligned}
    s_{n+1} & := \int_{0}^u \Om_n \,\d u',\\
    {^{n+1}\Ld}(\log\Om_{n+1}) & := s_{n+1}^{-1}{^{n+1}\Ld} s_{n+1} + F(s_{n+1},{^{n+1}\Nd} s_{n+1}, {^{n+1}\Nd}^2s_{n+1}) \\
    & \quad -\overline{s_{n+1}^{-1}{^{n+1}\Ld} s_{n+1} + F(s_{n+1},{^{n+1}\Nd} s_{n+1}, {^{n+1}\Nd}^2s_{n+1})}, \\
    \overline{\log\Om_{n+1}} & := \int_0^u \overline{(\trchib_n-\trchibo_n)\log\Om_n}(u')\,\d u',
  \end{aligned}
\end{align}
where $^{n+1}\Nd$ and $^{n+1}\Ld$ are the covariant derivative and Laplace operator on the $2$-spheres $S_u$ associated to the metric $\gd'\le(s_{n+1}(u)\ri)$, and where in the last line the average is taken with respect to the metric $\gd'\le(s_{n+1}(u)\ri)$.\\

In the rest of this iteration argument, we fix an integer $k \geq 5$ such that the estimates~\eqref{est:Fisk} hold. We simply note $A$ the constant $A_k$.\\

Let define
\begin{align*}
  M^{\de,s}_n & := \norm{u^{-4} (u\Nd)^{\leq k}(s_n-u)}_{L^\infty_u([0,\de])L^2(S_u)}, \\
  M^{\de,\Om}_n & := \norm{u^{-3}(u\Nd)^{\leq k}\log\Om_{n}}_{L^\infty_u([0,\de])L^2(S_u)},
\end{align*}
where the norms are taken with respect to the metric $\gd_\SSS$ and $\Nd = {^{\gd_{\SSS}}\Nd}$ with
\begin{align*}
  \gd_\SSS & := \quar u^2\le(\d\varth^2 + \sin^2\varth\d\varphi^2\ri),
\end{align*}
where $\varth,\varphi$ are the transported from the vertex spherical coordinates.\\

\begin{remark}
  From standard Sobolev estimates, we have
  \begin{align*}
    \le|s_n-u\ri| & \les u^3 M^{\de,s}_n, & \le|\log\Om_n\ri| & \les u^2 M^{\de,\Om}_n,
  \end{align*}
  which gives expected limits when $u\to 0$ in view of the results of Theorem~\ref{thm:vertex} and Appendix~\ref{app:vertexlimits}.\footnote{Roughly speaking, $\Om$ is at the level of a metric component $\g_{\mu\nu}$ and thus is expected to decay in $u^2$ towards its Minkowskian value. The parameter $s$ (or reciprocally $u$) is an integral of $\Om$ and is thus expected to decay in $u^3$ towards its Minkowskian value.}
\end{remark}

We have the following boundedness result for the Banach-Picard iteration.
\begin{lemma}\label{lem:boundednessiteratecancan}
  Let $\varep>0$ be the constant defined in Section~\ref{sec:transelldefgeom}. Let $\ga$ be a numerical constant such that
  \begin{align*}
    \frac{1}{3} < \ga < \frac{1}{2}.
  \end{align*}
  There exists a constant $\de_0(A,\ga,\varep)>0$ and a constant $C(A,\ga)>0$ such that for $0<\de<\de_0(A,\ga,\varep)$, we have 
  \begin{align}\label{est:devarepsncan}
    \forall u\in(0,\de),~ 0 < s_n< \varep,
  \end{align}
  and
  \begin{align}\label{est:snulogOmncan}
    M^{\de,s}_n & \leq \ga C(A,\ga), & M^{\de,\Om}_n & \leq C(A,\ga),
  \end{align}
  for all integers $n\geq 0$.
\end{lemma}

\begin{remark}
  The iteration assumption~\eqref{est:devarepsncan} ensures that the iterate $s_n$ remains in the local domain of existence of the geodesic foliation on $\CCba$ of Section~\ref{sec:transelldefgeom}.
\end{remark}

\begin{remark} 
  To close the boundedness result of the iteration it is crucial that the constant $\ga$ in estimate~\eqref{est:snulogOmncan} is comprised in the range $1/3 < \ga < 1/2$. See conditions~(\ref{est:conditiongadeiterate}) and~(\ref{eq:C(A)defiterate}).
\end{remark}

\begin{proof}[Proof of Lemma~\ref{lem:boundednessiteratecancan}]
  Let first remark that the estimates of Lemma~\ref{lem:boundednessiteratecancan} trivially hold true at the initial step $n=0$. Let $n\geq 0$ be a fixed integer and assume that the bounds~\eqref{est:devarepsncan} and~\eqref{est:snulogOmncan} hold true. Let show that these bounds also hold at step $n+1$, provided that the constants $\de,C$ have been appropriately chosen.\\
  
  From~\eqref{eq:syscancaniterate}, we first have
  \begin{align*}
    (u\Nd)^{\leq k}\le(s_{n+1}-u\ri) & = \int_0^u(u\Nd)^{\leq k}\le(\Om_n-1)\ri) \,\d u', \\
                                  & = \int_0^u (u\Nd)^{\leq k}\log\Om_n\, \d u' + \int_0^u (u\Nd)^{\leq k}\le(\Om_n-1-\log\Om_n\ri)\,\d u',
  \end{align*}
  where we recall that $\Nd = {^{\gd_\SSS}\Nd}$. Taking the $L^2(S_u)$ norm (with respect to the metric $\gd_\SSS$), we infer
  \begin{align*}
    & u^{-4}\norm{(u\Nd)^{\leq k}\le(s_{n+1}-u\ri)}_{L^2(S_u)} \\
    \leq & \; u^{-3}\int_{0}^u (u')^{-1}\norm{(u\Nd)^{\leq k}(\log\Om_n)}_{L^2(S_{u'})} \, \d u' + u^{-3} \int_0^u (u')^{-1} \norm{(u\Nd)^{\leq k}\le(\Om_n-1-\log\Om_n\ri)}_{L^2(S_{u'})} \, \d u'\\
    \leq & \; u^{-3}\int_0^u (u')^2 M^{\de,\Om}_n \, \d u' + C' u^{-3}\int_0^u (u')^{-1}\norm{(u\Nd)^{\leq k}\log\Om_n}^2_{L^2(S_{u'})} \, \d u' \\
    \leq & \; \frac{1}{3} M^{\de,\Om}_n + C' u^{-3} \int_0^u (u')^{5} \le((u')^{-3}\norm{(u\Nd)^{\leq k}\log\Om_n}_{L^2(S_{u'})}\ri)^2 \, \d u'\\
    \leq & \; \frac{1}{3} C(A,\ga) + C' u^3 C(A,\ga)^2 \\
    \leq & \; \ga C(A,\ga),
  \end{align*}
  where $C'>0$ is a numerical constant, and where the last line holds provided that $\de>0$ is chosen such that 
  \begin{align}\label{est:conditiongadeiterate}
    \frac{1}{3} + C'\de^3C(A,\ga) < \ga.
  \end{align}
  
  For this choice of $\de$, we thus have
  \begin{align}\label{est:sn+1iteratecancan}
    M_{n+1}^{\de,s} & \leq \ga C(A,\ga).
  \end{align}

  From a standard Sobolev estimate, we infer in particular from~\eqref{est:sn+1iteratecancan} that estimate~(\ref{est:devarepsncan}) holds for $s_{n+1}$.\\
  
  From~\eqref{eq:syscancaniterate}, Sobolev estimates for $\log\Om_n$, and estimate~\eqref{est:trchibtrchiboiterate}, we have
  \begin{align*}
    |\overline{\log\Om_{n+1}}| & \les \int_0^u \le|\trchib_n(s)-\trchibo_n(s)\ri| |\log\Om_n| \, \d u' \\
                               & \les \int_0^u \le(A u'\ri)  \le(M_n^{\de,\Om} (u')^2\ri)\, \d u' \\
                               & \les u^4 A C(A,\ga).
  \end{align*}
  Thus, for $u\leq \de$, we have
  \begin{align}\label{est:controliterateoverlineOm}
    u^{-2}|\overline{\log\Om_{n+1}}| & \les \de^2 A C(A,\ga). 
  \end{align}
  
  We rewrite the Laplace equation of~\eqref{eq:syscancaniterate} under the form
  \begin{align}\label{eq:Laplacemodifiediterate}
    \begin{aligned}
      {^{n+1}\Ld}\le(\log\Om_{n+1}-u^{-1}s_{n+1}\ri) & = \le(s_{n+1}^{-1}-u^{-1}\ri){^{n+1}\Ld}s_{n+1} + F\le(s_{n+1},{^{n+1}\Nd} s_{n+1}, {^{n+1}\Nd}^2s_{n+1}\ri) \\
      & \quad - \overline{\le(s_{n+1}^{-1}-u^{-1}\ri){^{n+1}\Ld}s_{n+1} + F\le(s_{n+1},{^{n+1}\Nd} s_{n+1}, {^{n+1}\Nd}^2s_{n+1}\ri)}.
    \end{aligned}
  \end{align}
  
  Applying standard elliptic estimates to~\eqref{eq:Laplacemodifiediterate}, we have the following estimate
  \begin{align*}
    & u^{-3}\norm{(u\Nd)^{\leq 2}\le(\log\Om_{n+1}-u^{-1}(s_{n+1}-\overline{s_{n+1}})\ri)}_{L^2(S_u)} \\
    \les & \; u^{-2}\le|\overline{\log\Om_{n+1}}\ri| + u\norm{u^{-1}-s_{n+1}^{-1}}_{L^\infty(S_u)} \le(u^{-4}\norm{(u\Nd)^2s_{n+1}}_{L^2(S_u)}\ri) \\
    & + u^{-1}\norm{F_1(s_{n+1})}_{L^2(S_u)} + u^2 \norm{F_2(s_{n+1})}_{L^\infty(S_u)} \le(u^{-4}\norm{(u\Nd) s_{n+1}}_{L^2(S_u)}\ri) \\
    & + u^{-3}\norm{F_3(s_{n+1})}_{L^\infty(S_u)} \norm{(u\Nd) s_{n+1}}_{L^4(S_u)} \norm{(u\Nd)s_{n+1}}_{L^4(S_u)} \\
    & + u \norm{F_4(s_{n+1})}_{L^\infty(S_u)} \le(u^{-4}\norm{(u\Nd)^2s_{n+1}}_{L^2(S_u)}\ri) \\
    & + u^{-1}\norm{\le(\Ld-{^{n+1}\Ld}\ri)\le(\log\Om_{n+1}-u^{-1}(s_{n+1}-\overline{s_{n+1}})\ri)}_{L^2(S_u)}.
  \end{align*}
  From~\eqref{est:Fisk},~\eqref{est:sn+1iteratecancan},~\eqref{est:controliterateoverlineOm} and standard Sobolev estimates, we have the following estimates
  \begin{align*}
    u^{-2}\le|\overline{\log\Om_{n+1}}\ri| & \les \de^2 A C(A,\ga), \\
    u\le|u^{-1}-s_{n+1}^{-1}\ri| \le(u^{-4}\norm{(u\Nd)^2s_{n+1}}_{L^2(S_u)}\ri) & \les u^{-1} \le|s_{n+1}-u\ri| C(A,\ga) \\
                                           & \les \de^2 C(A,\ga)^2, \\
    u^{-1}\norm{F_1(s_{n+1})}_{L^2(S_u)} & \les \norm{F_1(s_{n+1})}_{L^\infty(S_u)} \\
                                           & \les A, \\
    u^2 \norm{F_2(s_{n+1})}_{L^\infty(S_u)} \le(u^{-4}\norm{(u\Nd) s_{n+1}}_{L^2(S_u)}\ri) & \les \de^2 A C(A,\ga) \\
    u^{-3}\norm{F_3(s_{n+1})}_{L^\infty(S_u)} \le(\norm{(u\Nd) s_{n+1}}_{L^4(S_u)}\ri)^2 & \les A u^{-6}\le( \norm{(u\Nd)^{\leq 2}(s_{n+1}-u)}_{L^2(S_u)}\ri)^2 \\
                                           & \les A u^2 \le(u^{-4}\norm{(u\Nd)^{\leq 2}(s_{n+1}-u)}_{L^2(S_u)}\ri)^2 \\
                                           & \les \de^2 A C(A,\ga) \\
    u \norm{F_4(s_{n+1})}_{L^\infty(S_u)} \le(u^{-4}\norm{(u\Nd)^2s_{n+1}}_{L^2(S_u)}\ri) & \les \de^2AC(A,\ga),\\
    u^{-1}\norm{\le(\Ld-{^{n+1}\Ld}\ri)\le(\log\Om_{n+1}-u^{-1}(s_{n+1}-\overline{s_{n+1}})\ri)}_{L^2(S_u)} & \les u^{-1}A\norm{s_n-u}_{L^\infty(S_u)} \\
                                           & \quad \times \bigg(u^{-3}\norm{(u\Nd)^{\leq 2}\log\Om_{n+1}}_{L^2(S_u)} \\
                                           & \quad+ u^{-4}\norm{(u\Nd)^{\leq 2}(s_{n+1}-u)}_{L^2(S_u)}\bigg) \\
    & \les A\de^2\le(\norm{(u\Nd)^{\leq 2}\log\Om_{n+1}}_{L^2(S_u)} + C(A)\ri).
  \end{align*}
  Thus, collecting these estimates, we infer
  \begin{align*}
    u^{-3}\norm{(u\Nd)^{\leq 2}\le(\log\Om_{n+1}-u^{-1}(s_{n+1}-\overline{s_{n+1}})\ri)}_{L^2(S_u)} & \les A + A\de^2\norm{(u\Nd)^{\leq 2}\log\Om_{n+1}}_{L^2(S_u)} \\
    & \quad + \de^2\le(AC(A,\ga) +C(A,\ga)^2\ri). 
  \end{align*}
  Commuting equation~\eqref{eq:Laplacemodifiediterate}, we similarly obtain
  \begin{align}\label{est:Laplacemodifiediterate}
    \begin{aligned}
      & u^{-3}\norm{(u\Nd)^{\leq k}\le(\log\Om_{n+1}-u^{-1}(s_{n+1}-\overline{s_{n+1}})\ri)}_{L^2(S_u)}\\
      \les & \;  A +A\de^2M^{\de,\Om}_{n+1} + \de^2\le(AC(A,\ga) +C(A,\ga)^2\ri).
    \end{aligned}
  \end{align}
  
  Using~\eqref{est:sn+1iteratecancan} and~\eqref{est:Laplacemodifiediterate}, we infer
  \begin{align*}
    \begin{aligned}
      & u^{-3}\norm{(u\Nd)^{\leq k}\le(\log\Om_{n+1}\ri)}_{L^2(S_u)} \\
      \leq & \; u^{-4}\norm{(u\Nd)^{\leq k}(s_{n+1}-\overline{s_{n+1}})}_{L^2(S_u)} +  u^{-3}\norm{(u\Nd)^{\leq k}\le(\log\Om_{n+1}-u^{-1}(s_{n+1}-\overline{s_{n+1}})\ri)}_{L^2(S_u)} \\
      \leq & \; 2 u^{-4}\norm{(u\Nd)^{\leq k}(s_{n+1}-u)}_{L^2(S_u)} +  u^{-3}\norm{(u\Nd)^{\leq k}\le(\log\Om_{n+1}-u^{-1}(s_{n+1}-\overline{s_{n+1}})\ri)}_{L^2(S_u)} \\
      \leq & \; 2 \ga C(A,\ga) +C'A\de^2M^{\de,\Om}_{n+1} + C'\le(A + \de^2\le(AC(A,\ga) +C(A,\ga)^2\ri)\ri),
    \end{aligned}
  \end{align*}
  where $C'>0$ is a numerical constant. We rewrite this estimate as
  \begin{align}\label{est:logOmn+1iteratecan}
    \begin{aligned}
      M_{n+1}^{\de,\Om} & \leq \le(1-C'A\de^2\ri)^{-1}\le(2\ga C(A,\ga) + C'\le(A + \de^2\le(AC(A,\ga) +C(A,\ga)^2\ri)\ri) \ri),
    \end{aligned}
  \end{align}
  provided that
  \begin{align}\label{eq:desmallnessiteratebis}
    1-C'A\de^2>0.
  \end{align}

  Let now define $C(A,\ga)$ by
  \begin{align}\label{eq:C(A)defiterate}
    C(A,\ga) & := \frac{2C'A}{1-2\ga}.
  \end{align}
  Let further impose that $\de$ is chosen such that 
  \begin{align}\label{est:conditiongadeiteratebis}
    \de^2 < \frac{(1-2\ga)}{4C'A+C'C(A,\ga)}.
  \end{align}
  
  Provided that conditions~\eqref{eq:desmallnessiteratebis}~\eqref{eq:C(A)defiterate} and~\eqref{est:conditiongadeiteratebis} hold, we infer from~\eqref{est:logOmn+1iteratecan} and a direct computation that
  \begin{align*}
    M_{n+1}^{\Om,\de} & \leq C(A,\ga),
  \end{align*}
  as desired. This finishes the proof of the boundedness Lemma~\ref{lem:boundednessiteratecancan}.
\end{proof}
  
Let now define
\begin{align*}
  N^\de_n & := \max \bigg(\norm{u^{-4} (u\Nd)^{\leq k}(s_{n+1}-s_n)}_{L^\infty_u([0,\de])L^2(S_u)},\norm{u^{-3}(u\Nd)^{\leq k}(\log\Om_{n+1}-\log\Om_n)}_{L^\infty_u([0,\de])L^2(S_u)}\bigg),
\end{align*}
where the norms and covariant derivative are taken with respect to the metric $\gd_\SSS$.\\

We have the following contraction result for the Banach-Picard iteration.
\begin{lemma}
  \label{lem:contractioncancan}
  Let $\varep>0$ be the constant defined in Section~\ref{sec:transelldefgeom}. Let $\gamma$ be a constant such that
  \begin{align*}
    \frac{1}{3} < \ga < \frac{1}{2}.
  \end{align*}
  There exists $\de_1(A,\ga,\varep)>0$ such that for all $0< \de < \de_1(A,\ga,\varep)$, we have
  \begin{align}
    N^\de_{n+1} < 2\ga N^\de_n,
  \end{align}
  for all integers $n\geq 0$.
\end{lemma}
\begin{proof}
  We only sketch the proof, which goes along the same lines as the proof of Lemma~\ref{lem:boundednessiteratecancan}. Arguing as in Lemma~\ref{lem:boundednessiteratecancan}, we first obtain from the relation
  \begin{align*}
    s_{n+1}-s_n & = \int_0^u \le(\Om_{n}-\Om_{n-1}\ri) \, \d u',
  \end{align*}
  that
  \begin{align}\label{est:contractionsniterate}
    u^{-4}\norm{(u\Nd)^{\leq k}(s_{n+1}-s_n)}_{L^2(S_u)} & \leq \le(\frac{1}{3} + O(\de^2)\ri) N^\de_n.
  \end{align}
  Arguing as in Lemma~\ref{lem:boundednessiteratecancan}, using the rewritten Laplace equation~(\ref{eq:Laplacemodifiediterate}), we have
  \begin{align*}
    ^{n+1}\Ld\le(\log\Om_{n+1}-\log\Om_n - u^{-1}\le(s_{n+1}-s_n\ri)\ri) & = F\le(s_{n+1},{^{n+1}\Nd}s_{n+1},{^{n+1}\Nd}^2s_{n+1}\ri) - F\le(s_n,{^{n}\Nd}s_n,{^{n}\Nd}^2s_n\ri), \\
                                                                         & = G\bigg(s_{n},s_{n+1},\Nd s_n, \Nd s_{n+1},\Nd^2s_n,\Nd^2s_{n+1},\le(s_{n+1}-s_n\ri)\bigg)
  \end{align*}
  from which, we deduce by standard elliptic estimate
  \begin{align*}
    & u^{-3}\norm{(u\Nd)^{\leq k}\le(\log\Om_{n+1}-\log\Om_n\ri) -u^{-1}(s_{n+1}-s_n) + u^{-1}\overline{s_{n+1}-s_n}}_{L^2(S_u)} \\
    \les & \;  O(\de^2) \le(u^{-4}\norm{(u\Nd)^{\leq k}(s_{n+1}-s_n)}_{L^2(S_u)}\ri). 
  \end{align*}
  Using~\eqref{est:contractionsniterate} as in Lemma~\ref{lem:boundednessiteratecancan}, it follows that
  \begin{align*}
    u^{-3}\norm{(u\Nd)^{\leq k}\le(\log\Om_{n+1}-\log\Om_n\ri)}_{L^2(S_u)} & \leq \le(2+O(\de^2)\ri) u^{-4}\norm{(u\Nd)^{\leq k}\le(s_{n+1}-s_n\ri)}_{L^2(S_u)}.
  \end{align*}
  From that last bound and~\eqref{est:contractionsniterate}, we deduce
  \begin{align*}
    N_{n+1}^\de & \leq \le(\frac{2}{3} + O(\de^2)\ri) N_n^\de,
  \end{align*}
  from which the result of the lemma follows.
\end{proof}

\section{Conclusions}
Using the results of Lemmas~\ref{lem:boundednessiteratecancan} and \ref{lem:contractioncancan}, the Banach-Picard iteration converges. Moreover, its limit defines the desired canonical foliation (see also the argument in~\cite[Section 6.3]{Czi.Gra19}), and we have
\begin{align}\label{est:finalsu}
  \norm{u^{-4} (u\Nd)^{\leq k}(s-u)}_{L^\infty_u([0,\de])L^2(S_u)} + \norm{u^{-3}(u\Nd)^{\leq k}\log\Om}_{L^\infty_u([0,\de])L^2(S_u)} & < \infty.
\end{align}
For all $k \geq 5$, estimates~\eqref{est:finalsu} hold true at the vertex, and thus from Sobolev embeddings we deduce that
\begin{align*}
  \begin{aligned}
    \le|(u\Nd)^{k} (s-u)\ri| & = O(u^3), & \le|(u\Nd)^{k}\log\Om\ri| = O(u^2),
  \end{aligned}
\end{align*}
when $u\to 0$ and for all $k\geq 0$. Commuting equation~\eqref{eq:systemcancanbis} by $\Nd_3$, one can further obtain
\begin{align}\label{est:limitssOm}
  \begin{aligned}
    \le|(u\Nd)^{\leq k}(u\Nd_3)^{l} (s-u)\ri| & = O(u^3), & \le|(u\Nd)^{\leq k}(u\Nd_3)^l\log\Om\ri| = O(u^2),
  \end{aligned}
\end{align}
for all $k,l\geq 0$. Using the relations of Lemma~\ref{prop:transconn} between null connection coefficients, relations~(\ref{eq:reltransOm}) and the limits~\eqref{est:limitssOm}, we deduce that the limits of Theorem~\ref{thm:vertex} also hold for the null connection coefficients of the canonical foliation. The spherical coordinates $(s,\varth,\varphi)$ for the geodesic foliation from Theorem~\ref{thm:vertex} also generate spherical coordinates $(u,\varth,\varphi)$ for the canonical foliation. Moreover, since the coordinate vectorfields for $\varth,\varphi$ only differ by null vectors, we have $\gd'_{ab} = \gd_{ab}$, where $a,b\in\{\varth,\varphi\}$. Thus the limits for the induced metric from Theorem~\ref{thm:vertex} also holds for the canonical foliation. This finishes the proof of Theorem~\ref{thm:canonical}.

\chapter{Klainerman-Sobolev estimates}\label{app:KlSobH12}
This section is dedicated to the proof of the Klainerman-Sobolev $\HHt$ estimates of Lemmas~\ref{lem:KlSobast}, \ref{lem:KlSobSitext}, \ref{lem:KlSobastintr} and~\ref{lem:KlSobcomplet}. In all cases we only treat the degenerate versions involving respectively $u$ and $\qq$ weights. The non-degenerate cases follow similarly, or by rescaling of standard $H^1$ to $H^{1/2}$ trace estimate.\\

We first prove the following Klainerman-Sobolev $H^{1/2}$ estimates  in the Euclidean case. The proofs of Lemmas~\ref{lem:KlSobast}, \ref{lem:KlSobSitext}, \ref{lem:KlSobastintr} and~\ref{lem:KlSobcomplet} are deduced from Lemma~\ref{lem:KlSobH12Eucl} at the end of this section.
\begin{lemma}[Euclidean Klainerman-Sobolev $H^{1/2}$ estimate]\label{lem:KlSobH12Eucl}
  Let $S_r$ denote the Euclidean sphere of radius $r$, and let $\HHt(S_r)$ denote the scaling homogeneous fractional Sobolev space on $S_r$, Euclidean analogous to the definition of $\HHt$ from Definition~\ref{def:Ht}. Let $\qq(r)$ be a scalar function of $r$ such that $|\qq| \les r$ and $|\pr_r\qq|\les 1$. For all $r,r'>0$ and for all scalar function $f$, we have
  \begin{align*}
    \norm{r_2^{1/2}\qq(r_2)^{1/2}f}_{\HHt(S_{r_2})} & \les \norm{(r_1)^{1/2}\qq(r_1)^{1/2}f}_{\HHt(S_{r_1})} + \norm{\qq(r)\pr_rf}_{L^2(A_{r_1,r_2})} + \norm{(r\Nd)^{\leq 1}f}_{L^2(A_{r_1,r_2})},
  \end{align*}
  where $A_{r_1,r_2}$ denotes the annulus region comprised between $S_{r_1}$ and $S_{r_2}$.
\end{lemma}
\begin{proof}
  For simplicity, we use the following characterisation of the fractional Sobolev space $\HHt(S_r)$
  \begin{align}\label{eq:Sobcharac}
    \le(r^{-1/2}\norm{f}_{\HHt(S_r)}\ri)^2 & \simeq \int \le|(\Id-\Ld_{\varth,\varphi})^{1/4}f\ri|^2 \,\d\varth\d\varphi.
  \end{align}
  Deriving~\eqref{eq:Sobcharac} in $r$ and integrating from $r_1$ to $r_2$, we infer
  \begin{align*}
    \le(r_2^{-1/2}\norm{r_2\qq(r_2)^{1/2}f}_{\HHt(S_{r_2})}\ri)^2 & \les \le((r_1)^{-1/2}\norm{r_1\qq(r_1)^{1/2}f}_{\HHt(S_{r_1})}\ri)^2 \\
    & + \le|\int_{r_1}^{r_2}\int (\Id-\Ld_{\varth,\varphi})^{1/4}(r\qq^{1/2}f) (\Id-\Ld_{\varth,\varphi})^{1/4}(\pr_r(r\qq^{1/2}f)) \,\d r\d\varth\d\varphi\ri|.
  \end{align*}
  Using the assumptions $|\qq| \les r$ and $|\pr_r\qq| \les 1$, integrating by part and using Cauchy-Schwartz, we have
  \begin{align*}
    & \le|\int_{r_1}^{r_2}\int \le|(\Id-\Ld_{\varth,\varphi})^{1/4}f\ri|^2(r\qq^{1/2}\pr_r(r\qq^{1/2})) \,\d r\d\varth\d\varphi\ri| \\
    \les \; & \le|\int_{r_1}^{r_2}\int \le|(\Id-\Ld_{\varth,\varphi})^{1/4}f\ri|^2 r^2 \,\d r\d\varth\d\varphi\ri| \\
    \les \; & \norm{(r\Nd)^{\leq 1}f}_{L^2(A_{r_1,r_2})}\norm{f}_{L^2(A_{r_1,r_2})}.
  \end{align*}
  Integrating by part and applying Cauchy-Schwartz, we also have
  \begin{align*}
    & \le|\int_{r_1}^{r_2}\int (\Id-\Ld_{\varth,\varphi})^{1/4}(r\qq^{1/2}f) (\Id-\Ld_{\varth,\varphi})^{1/4}(r\qq^{1/2}\pr_rf) \,\d r\d\varth\d\varphi\ri| \\
    \les  \; & \le|\int_{r_1}^{r_2}\int \le((\Id-\Ld_{\varth,\varphi})^{1/2}f\ri) (\qq\pr_rf) r^2\,\d r\d\varth\d\varphi\ri| \\
    \les \; & \norm{(r\Nd)^{\leq 1}f}_{L^2(A_{r_1,r_2})}\norm{\qq\pr_rf}_{L^2(A_{r_1,r_2})}.
  \end{align*}
  Combining the above estimates finishes the proof of the lemma.
\end{proof}

\begin{proof}[Proof of Lemma~\ref{lem:KlSobSitext}]
  We only treat the degenerate version involving $u$ weights which is the most involved estimate of Lemma~\ref{lem:KlSobSitext}. The non-degenerate version follows similarly.\\
  
  Let $F$ be an $S$-tangent $1$-tensor in $\MM^\ext$. Using the comparison between fractional Sobolev spaces of Lemma~\ref{lem:compH12}, we have for all $2$-spheres $S_{u,\ub}\subset \Si_t^\ext$
  \begin{align*}
    \norm{\ub^{1/2}u^{1/2}F}_{\HHt(S_{u,\ub})} & \simeq \sum_{(\varth,\varphi)} \sum_{a\in\{\varth,\varphi\}} \norm{r u^{1/2} (r^{-1}F)_a}_{H^{1/2}_{\varth,\varphi}} \simeq \sum_{a\in\{\varth,\varphi\}} \norm{r^{1/2} \qq(r)^{1/2}(r^{-1}F)_a}_{\HHt(S_{r=(\ub-u)/2})}, 
  \end{align*}
  where here $\qq(r) = u$. Using the Klainerman-Sobolev Euclidean estimates of Lemma~\ref{lem:KlSobH12Eucl} between the $2$-spheres $S_{u,\ub}$ and $S_{u+\ub-\uba,\uba}$ in coordinates, we thus infer
  \begin{align}\label{est:KlSobsemi}
    \begin{aligned}
      \norm{\ub^{1/2}u^{1/2}F}_{\HHt(S_{u,\ub})} & \les \norm{\uba^{1/2}u^{1/2}F}_{\HHt(S_{u+\ub-\uba,\uba})} \\
      & + \sum_{a\in\varth,\varphi}\le(\norm{u(\pr_\ub-\pr_u)(r^{-1}F_a)}_{L^2(\Si_t^\ext)} + \norm{(r\Nd)^{\leq 1}(r^{-1}F_a)}_{L^2(\Si_t^\ext)}\ri).
    \end{aligned}
  \end{align}
  We have $\pr_\ub-\pr_u = \half (1+\half\yy)\el -\half\elb$. Therefore, using the bootstrap bound $|\yy| \les D\varep$, we infer
  \begin{align}\label{est:prubpruFa}
    \begin{aligned}
      |(\pr_\ub-\pr_u)(F_a)| & \les |\Nd_4F_a| + |\Nd_3F_a| + |F\cdot\chi_a| + |F\cdot\chib_a| \\
      & \les r\le(|\Nd_4F| +|\Nd_3F| + r^{-1}|F|\ri).
    \end{aligned}
  \end{align}
  Similarly, we have
  \begin{align*}
    |(r\Nd)(F_a)| & \les |(r\Nd)F_a| + | F\cdot \Ga_a|,
  \end{align*}
  where $\Ga$ are the Christoffel symbols for the metric $\gd$ in coordinates $\varth,\varphi$. Using the Bootstrap Assumptions~\ref{BA:mildsphcoordsext}, we have $|\Ga| \les 1$, thus
  \begin{align}\label{est:rNdFa}
    |(r\Nd)(F_a)| & \les r|(r\Nd)^{\leq 1}F|.
  \end{align}
  Combining~\eqref{est:KlSobsemi},~\eqref{est:prubpruFa} and~\eqref{est:rNdFa}, and that $r^{-1}u \les 1$ in $\MM^\ext$, this concludes the proof of Lemma~\ref{lem:KlSobSitext}.
\end{proof}

\begin{proof}[Proof of Lemmas~\ref{lem:KlSobastintr} and~\ref{lem:KlSobcomplet}]
  Arguing similarly as in the proof of Lemma~\ref{lem:KlSobSitext}, we obtain for all $2$-spheres $S_{u_1,\uba}$ and $S_{u_2,\uba}$ of $\CCba$
  \begin{align}\label{est:KlSobsemiCCba}
    \begin{aligned}
      \norm{r^{1/2}\qq^{1/2}F}_{\HHt(S_{u_2,\uba})} & \les \norm{r^{1/2}\qq^{1/2}F}_{\HHt(S_{u_1,\uba})} + \norm{r^{-1}\qq F}_{L^2(\CCba)} \\
      & \quad +\norm{(r\Nd)F}_{L^2(\CCba)} + \norm{\qq\Nd_3F}_{L^2(\CCba)}, 
    \end{aligned}
  \end{align}
  for $q=\min(u,r)$. Using the estimates of Lemma~\ref{lem:Sobsphere} we have
  \begin{align*}
    \norm{r^{1/2}\qq^{1/2}F}_{\HHt(S_{u_1,\uba})} & \les r^{-1/2}\norm{r^{1/2}\qq^{1/2}(r\Nd)^{\leq 1}F}_{L^2(S_{u_1,\uba})} \\
                                                  & \les r\qq^{1/2} \norm{(r\Nd)^{\leq 1}F}_{L^\infty(S_{u_1,\uba})},
  \end{align*}
  which by assumption tends to $0$ when $u_1\to\uba$. Passing to the limit in~\eqref{est:KlSobsemiCCba}, using additionally that $r^{-1}\qq \les 1$ on $\CCba$, this finishes the proof of the degenerate estimate of Lemma~\ref{lem:KlSobcomplet}. The non-degenerate version follows along the same lines. The proof of Lemma~\ref{lem:KlSobastintr} is obtained along the same lines, restricting the argument to the top region $\CCba\cap\MM^\intr_\bott$.\footnote{Lemma~\ref{lem:KlSobastintr} is simply a restriction of Lemma~\ref{lem:KlSobcomplet} to the interior region. The spliting into two statements is done for the simplicity of the exposition in Section~\ref{sec:remainingcurvestfinal}.}
\end{proof}
\begin{proof}[Proof of Lemma~\ref{lem:KlSobast}]
  Arguing as in the proof of Lemma~\ref{lem:KlSobSitext}, we obtain for all $2$-spheres $S_{u_1,\uba}$ and $S_{u_2,\uba}$ of $\CCba\cap\MM^\ext$
  \begin{align}\label{est:KlSobsemiCCbabis}
    \begin{aligned}
      \norm{r^{1/2}u^{1/2}F}_{\HHt(S_{u_2,\uba})} & \les \norm{r^{1/2}u^{1/2}F}_{\HHt(S_{u_1,\uba})} + \norm{r^{-1}u F}_{L^2(\CCba\cap\MM^\ext)} \\
      & \quad +\norm{(r\Nd)F}_{L^2(\CCba\cap\MM^\ext)} + \norm{u\Nd_3F}_{L^2(\CCba\cap\MM^\ext)},
    \end{aligned}
  \end{align}
  for all $S$-tangent tensor $F$.\\

  Let $\phi:[0,1]\to[0,1]$ be a standard cut-off function such that $\phi|_{[0,1/4]}=0$ and $\phi|_{[3/4,1]}=1$ and let $\tilde{\phi}(u) := \phi\le(\frac{u-1}{\cc\uba-1}\ri)$. Applying~\eqref{est:KlSobsemiCCbabis} to $\tilde{F}$, where $\tilde{F} := \tilde{\phi} F$, with $u_2 = \cc\uba$ and $u_1 = 1$, we obtain
  \begin{align}\label{est:KlSobsemiCCbabisbis}
    \norm{r^{1/2}u^{1/2}F}_{\HHt(S_{\cc\uba,\uba})} & \les \norm{r^{-1}u F}_{L^2(\CCba\cap\MM^\ext)} +\norm{(r\Nd)F}_{L^2(\CCba\cap\MM^\ext)} + \norm{u\Nd_3F}_{L^2(\CCba\cap\MM^\ext)}
  \end{align}
  where we used that $|\tilde{\phi}|\les 1$ and $|u\Nd_3\tilde{\phi}| \les 1$.\footnote{The same argument holds for the non-degenerate estimate since we also have $|r\Nd_3\tilde{\phi}| \les 1$.} Using~\eqref{est:KlSobsemiCCbabis} with $u_2=u$ and $u_1=\cc\uba$ and using~\eqref{est:KlSobsemiCCbabisbis}, we obtain the desired result. The non-degenerate estimates of Lemma~\ref{lem:KlSobast} follow along the same lines. This finishes the proof of Lemma~\ref{lem:KlSobast}.
\end{proof}
\end{appendices}


\end{document}